\documentclass[11pt,a4paper]{amsart}

\usepackage[latin1]{inputenc}

\usepackage{txfonts}

\usepackage{enumerate}

\usepackage{float}
\restylefloat{table}

\usepackage{graphicx}

\usepackage{hyperref}

\usepackage{ulem}

\usepackage{amsmath}

\usepackage{amssymb}

\usepackage{amsthm}

\usepackage{tikz-cd}

\tikzstyle{vertex}=[auto=left,circle,fill=black,minimum size=4pt,inner sep=0pt]

\usepackage{a4wide}

\usepackage[toc, page]{appendix}

\allowdisplaybreaks

\newtheoremstyle{plain}
  {6pt}   
  {6pt}   
  {\itshape}  
  {0pt}       
  {\bfseries} 
  {.}         
  {5pt plus 1pt minus 1pt} 
  {}          

\newtheoremstyle{definition}
  {6pt}   
  {6pt}   
  {\normalfont}  
  {0pt}       
  {\bfseries} 
  {.}         
  {5pt plus 1pt minus 1pt} 
  {}          

\theoremstyle{plain}
\newtheorem*{thm*}{Theorem}
\newtheorem{thm}{Theorem}[section]
\newtheorem{prop}[thm]{Proposition}
\newtheorem{cor}[thm]{Corollary}
\newtheorem{lem}[thm]{Lemma}

\newtheorem{thmx}{Theorem}

\theoremstyle{definition}
\newtheorem{defn}[thm]{Definition}

\newtheorem{ex}[thm]{Example}
\newtheorem{rmk}[thm]{Remark}
\newtheorem*{q}{Question}

\numberwithin{equation}{thm}

\newcommand{\emphbf}[1]{\textbf{#1}}

\DeclareMathAlphabet{\mathpzc}{OT1}{pzc}{m}{it}

\usetikzlibrary{decorations.pathreplacing}

\DeclareMathOperator{\radoperator}{rad}

\DeclareMathOperator{\id}{id}
\DeclareMathOperator{\Hom}{Hom}
\DeclareMathOperator{\Ext}{Ext}

\DeclareMathOperator{\Mod}{Mod}
\DeclareMathOperator{\modu}{mod}

\DeclareMathOperator{\im}{Im}

\DeclareMathOperator{\End}{End}

\DeclareMathOperator{\image}{Im}

\DeclareMathOperator{\op}{op}

\DeclareMathOperator{\tria}{tria}

\DeclareMathOperator{\shift}{\mathbf{s}}

\DeclareMathOperator{\twmod}{twmod}
\DeclareMathOperator{\twcom}{twcom}

\DeclareMathOperator{\mult}{mult}

\DeclareMathOperator{\+}{\boldsymbol{+}}
\renewcommand{\-}{\boldsymbol{-}}

\begin{document}

\title{Uniqueness of exact Borel subalgebras and bocses}
\author{Julian K\"ulshammer}
\author{Vanessa Miemietz}
\date{\today}

\address{Julian K\"ulshammer\\
Departement of Mathematics,
Uppsala University \\Box 480\\ 75106 Uppsala,
Sweden} \email{julian.kuelshammer@math.uu.se}

\address{Vanessa Miemietz\\
School of Mathematics, University of East Anglia\\
Norwich NR4 7TJ\\
United Kingdom} \email{v.miemietz@uea.ac.uk}

\begin{abstract}
In \cite{KKO14}, together with Koenig and Ovsienko, the first author showed that every quasi-hereditary algebra can be obtained as the (left or right) dual of a directed bocs. In this monograph, we prove that if one additionally assumes that the bocs is basic, a notion we define, then this bocs is unique up to isomorphism. This should be seen as a generalisation of the statement that the basic algebra of an arbitrary associative algebra is unique up to isomorphism. The proof associates to a given presentation of the bocs an $A_\infty$-structure on the $\Ext$-algebra of the standard modules of the corresponding quasi-hereditary algebra. Uniqueness then follows from an application of Kadeishvili's theorem. 
\end{abstract}

\maketitle

\tableofcontents

\section{Introduction}

Throughout, let $\Bbbk$ be an algebraically closed field. Quasi-hereditary algebras (and the related notion of a highest weight category) were defined by Cline, Parshall, and Scott (see \cite{Sco87, CPS88}). Important examples include blocks of BGG category $\mathcal{O}$, Schur algebras, as well as algebras of global dimension at most two. Quasi-hereditary algebras have also been used as tools, e.g. in Iyama's proof of finiteness of representation dimension \cite{Iya03} and more recently by Orlov to show geometric realisability of finite dimensional algebras by projective schemes \cite{Orl18}. One of the various ways of defining quasi-hereditary algebras is by the existence of an exceptional collection \cite{Bon89} in their module category such that the projective modules admit a filtration with subquotients isomorphic to the modules in the exceptional collection. A collection of modules $\Delta(\mathtt{1}),\dots,\Delta(\mathtt{n})$ is called \emphbf{exceptional} if 
\begin{itemize}
\item $\End(\Delta(\mathtt{i}))\cong \Bbbk$,
\item $\Hom(\Delta(\mathtt{i}), \Delta(\mathtt{j}))\neq 0\Rightarrow \mathtt{i}\leq \mathtt{j}$,
\item $\Ext^1(\Delta(\mathtt{j}),\Delta(\mathtt{j}))\neq 0\Rightarrow \mathtt{i}<\mathtt{j}$.
\end{itemize}
Exceptional collections (mostly in triangulated categories) also appear in algebraic and symplectic geometry. Their significance there is that, under certain conditions, they provide an equivalence between a geometric and an algebraic category, the first example being the equivalence between the bounded derived category of coherent sheaves on projective space and the bounded derived category of finite dimensional modules over the Be{\u\i}linson algebra, see \cite{Bei78, Bei84}. 

In 1995, inspired by the prime example of BGG category $\mathcal{O}$, Steffen Koenig coined the notion of an exact Borel subalgebra (see \cite{Koe95}). This is a directed subalgebra of a quasi-hereditary algebra with the same number of isomorphism classes of simple modules such that the induction functor is exact and sends simples to standard modules. Here, directedness is meant as an analogue of solvability for the Borel subalgebra of a semisimple Lie algebra and is defined as being quasi-hereditary with simple standard modules or, equivalently, having an acyclic Gabriel quiver with arrows only going in increasing direction. 

Already in the first paper, Koenig posed the question of existence and uniqueness of exact Borel subalgebras and provided examples that both don't hold in general. However, he was able to prove existence in the case of blocks of BGG category $\mathcal{O}$. Despite various efforts, existence in general was an open question, even in the important example of Schur algebras of symmetric groups. Uniqueness was only proved in a very special case by Yuehui Zhang in \cite{Zha00}. 

Together with Steffen Koenig and Sergiy Ovsienko, in \cite{KKO14}, the first-named author gave a new characterisation of quasi-hereditary algebras up to Morita equivalence. Expressed in the language of \cite{BKK20}, it can be formulated as an algebra $R$ having a directed subalgebra $B$ with the same number of simples such that the induction functor is exact and the induced maps 
\[ \Ext^k_B(L(\mathtt{i}),L(\mathtt{j}))\to \Ext^k_R(\Lambda\otimes L(\mathtt{i}), \Lambda\otimes L(j)) \]
are epimorphisms for $k\geq 1$ and isomorphisms for $k>1$. In particular, every quasi-hereditary algebra has an exact Borel subalgebra in a potentially different Morita representative, answering the question of existence. The work of the current article, which started in 2015, aims to shed light on uniqueness of exact Borel subalgebras. The key notion is that of a regular exact Borel subalgebra, meaning that the above induced maps are even isomorphisms for $k\geq 1$. It follows from the results in \cite{KKO14} that such subalgebras always exist up to Morita equivalence. Here we prove the following:

\begin{thmx}
Let $R$ and $S$ be quasi-hereditary algebras together with regular exact Borel subalgebras $A\subseteq R$ and $B\subseteq S$ such that $A$ and $B$ are basic as algebras. Assume that $R$ and $S$ are Morita equivalent as quasi-hereditary algebras. Then there exists an algebra isomorphism $f\colon R\to S$ which restricts to an isomorphism between $A$ and $B$. 
\end{thmx}

We originally proved uniqueness of the subalgebra up to isomorphism in 2015, uniqueness of $R$ up to isomorphism in 2016. We review the easier approach which led to these two uniqueness results in Section \ref{sec:uniquenessofexactborelsubalgebras}. Proving uniqueness of the embedding is much harder and occupies most of this paper. A different, more elementary proof of uniqueness of $R$ can be found in \cite{Con21}. Comparing the situation to that of a Borel subalgebra of a semisimple Lie algebra, it seems plausible to suspect that such Borel subalgebras are even unique up to inner automorphism. A proof of this will, however, require different techniques.

Our proof crucially uses the theory of $A_\infty$-algebras as defined by Stasheff \cite{Sta63}. To explain the idea of our proof, it is convenient to go into the details of \cite{KKO14} and to also compare them to the situation of using simple modules instead of standard modules. To this end, let $A$ be an arbitrary finite dimensional algebra. It was already proved by Morita in \cite{Mor58} that the basic representative of the Morita equivalence class of $A$ is unique up to isomorphism. An overly complicated way of constructing it would be via $A_\infty$-Koszul duality as follows: Let $P^\bullet$ be a projective resolution of the direct sum $\mathbb{L}$ of the simple $A$-modules. Its dg endomorphism ring $\Hom_A(P^\bullet, P^\bullet)$ carries the structure of a dg algebra whose homology is isomorphic to $\Ext^\bullet_A(\mathbb{L}, \mathbb{L})$. By Kadeishvili's theorem \cite{Kad80, Kad82}, the latter  carries the structure of an $A_\infty$-algebra, that is a graded vector space $\mathcal{E}$ together with a collection of graded linear maps $m_n\colon \mathcal{E}^{\otimes n}\to \mathcal{E}$ of degree $2-n$ for $n\geq 1$ such that for all $n$ the identity $\sum_{r+s+t=n} (-1)^{r+st} m_{r+1+t}(\id^{\otimes r} \otimes m_s\otimes \id^{\otimes s})=0$ holds. It turns out that the dual bar construction of this $A_\infty$-algebra, sometimes called $A_\infty$-Koszul dual, is a dg algebra quasi-isomorphic to the basic representative of the Morita equivalence class of $A$. More explicitly, the $m_n$ induce dual maps $d_n\colon \Ext^2_A(\mathbb{L},\mathbb{L})^\#\to (\Ext^1_A(\mathbb{L},\mathbb{L})^\#)^{\otimes n}$ such that the basic representative is given by the following quotient of a tensor algebra:
\[ T_{\mathbb{L}}(\Ext^1_A(\mathbb{L},\mathbb{L})^\#)/\image (\sum d_n).\]

A complicated proof of uniqueness of the basic representative up to isomorphis can be given using a reverse construction. Given a presentation of the basic representative as a quiver with relations $A=\Bbbk Q/I$, and denoting the augmentation ideal of $\Bbbk Q$ by $Q_+$, it is well-known (\cite{Bon83}, see also \cite{Gov73}) that 
\[ \Bbbk Q_1\cong \Ext^1_A(\mathbb{L},\mathbb{L})^\# \text{ and } I/(Q_+ I+IQ_+)\cong \Ext^2_A(\mathbb{L},\mathbb{L})^\#. \]

A result by Keller (stated without proof in \cite{K01, K02}, see also \cite{LPWZ09, Seg08}) states that a splitting of the projection $I\to I/(Q_+I+I Q_+)$ can be chosen such that its dual is the restriction of an $A_\infty$-structure on $\Ext^\bullet_A(\mathbb{L},\mathbb{L})$, which is $A_\infty$-quasi-isomorphic to $\Hom_A(P^\bullet, P^\bullet)$. Two different presentations thus give rise to two different $A_\infty$-structures on $\Ext_A^\bullet(\mathbb{L}, \mathbb{L})$. However, Kadeishvili's theorem implies that these are $A_\infty$-isomorphic. Taking the dual of this isomorphism yields an isomorphism between the two basic representatives. 

Our strategy is to adapt this proof, replacing the $\Ext$-algebra of simple modules with the $\Ext$-algebra of standard modules over a quasi-hereditary algebra. The key difference, making the argument more involved, is that there is further information in the homomorphisms between standard modules, which, by Schur's lemma, is not present in the case of simple modules. 

A further ingredient used is the theory of corings or bocses  (an acronym for bimodule over category with coalgebra structure). The corings with surjective counits we consider are (one-sided) dual to the ring extensions $A\subseteq R$. The notion of regularity of exact Borel subalgebras is in fact motivated by work of Kleiner and Ro\u\i ter on bocses and is a key notion in Drozd's proof of the tame-wild dichotomy. Ro\u\i ter's equivalence between the category of (normal) corings and the category of semifree dg algebras provides a bridge to the theory of $A_\infty$-algebras. Let $\Lambda$ be an arbitrary quasi-hereditary algebra. The semifree dg algebra constructed in \cite{KKO14} is the dual bar construction of the ($A_\infty$-structure on the) $\Ext$-algebra of standard modules over $\Lambda$. More explicitly, the degree zero part of the semifree dg algebra is the exact Borel subalgebra $A$ constructed as the quotient
\[T(\Ext^1_\Lambda(\Delta,\Delta)^\#)/\image\left(\sum d_n\right)\]
where the $d_n$ are dual to the structure maps 
\[m_n\colon \Ext^1_\Lambda(\Delta,\Delta)^{\otimes n}\to \Ext^2_\Lambda(\Delta,\Delta).\]
In some sense, this treats the standard modules as simple modules by forgetting the homomorphisms between them. The semifree dg algebra is determined by the degree one part which is the projective $A$-bimodule generated by $\radoperator_\Lambda(\Delta,\Delta)^\#$, and the differential given by the dual of the structure maps
\[m_n\colon \bigoplus_{i+j=n-1}(\Ext^1_\Lambda(\Delta,\Delta))^{\otimes i}\otimes \radoperator_\Lambda(\Delta,\Delta)\otimes (\Ext^1_\Lambda(\Delta,\Delta))^{\otimes j}\to \Ext^1_\Lambda(\Delta,\Delta)\]
in degree zero and the dual of the structure maps 
\[\bigoplus_{i+j+k=n-2} (\Ext^1_\Lambda(\Delta,\Delta ))^{\otimes i}\otimes \radoperator_\Lambda(\Delta,\Delta)\otimes (\Ext^1_\Lambda(\Delta,\Delta))^{\otimes j}\otimes \radoperator_\Lambda(\Delta,\Delta)\otimes (\Ext^1_\Lambda(\Delta,\Delta))^{\otimes k}\to \radoperator_\Lambda(\Delta,\Delta)\]
in degree one. Our proof of uniqueness, similarly to the classical case, goes via a reverse construction. Suppose that $A=\Bbbk Q/I$ is a basic algebra and that $V$ is an $A$-coring such that $(A,V)$ is a regular directed bocs (i.e. $A$ is a regular exact Borel subalgebra in the one-sided dual algebra $R=\Hom_A(V,A)$). Regularity implies that, similarly to the classical case,
\[Q_1\cong \Ext^1_R(\Delta,\Delta)^\# \text{ and } I/( Q_+ I + I Q_+)\cong \Ext^2_R(\Delta,\Delta)^\#\]
and, in addition, that for the kernel of the counit $\overline{V}$ we have 
\[\overline{V}/(J\overline{V}+\overline{V}J)\cong \radoperator_R(\Delta,\Delta)\]
where $J$ denotes the Jacobson radical of $A$. Then there are splittings $I/(IQ_++Q_+I)\to \Bbbk Q$ and $A\to \Bbbk Q$ and an isomorphism $\overline{V}\cong A\otimes \Phi\otimes A$ where $\Phi$ denotes a generating set of the projective bimodule $\overline{V}$ such that, ignoring the grading for the purpose the introduction, the following technical second main result holds:

\begin{thmx}\label{mainthm}
Let $(A=\Bbbk Q/I,V)$ be a regular directed bocs. Let $R$ be its right algebra with standard modules $\Delta$. 
Then there is an $A_\infty$-structure on $\Ext^\bullet_R(\Delta,\Delta)$ such that 
\begin{enumerate}[(i)]
\item the morphism 
\[m_n\colon \Ext^1_R(\Delta,\Delta)^{\otimes n}\to \Ext^2_R(\Delta,\Delta)\] 
can be identified with the dual of the map
\[I/(IQ_++Q_+I)\to \Bbbk Q_+\twoheadrightarrow Q_1^{\otimes n};\] 
\item the morphism 
\[m_n\colon (\Ext^1_R(\Delta,\Delta))^{\otimes i}\otimes \radoperator_R(\Delta,\Delta)\otimes (\Ext^1_R(\Delta,\Delta))^{\otimes j}\to \Ext^1_R(\Delta,\Delta)\]
can be identified with the dual of the map 
\[Q_1\hookrightarrow A\stackrel{\partial_0}{\to} \overline{V}\cong A\otimes \Phi \otimes A\to \Bbbk Q\otimes \Phi\otimes \Bbbk Q\twoheadrightarrow Q_1^{\otimes i}\otimes \Phi\otimes Q_1^{\otimes j};\] 
where $\partial_0$ denotes the degree $0$ part of the semifree dg algebra corresponding to $(A,V)$ under Ro\u\i ter's equivalence. 
\item the morphism 
\[m_n\colon (\Ext^1_R(\Delta,\Delta))^{\otimes i}\otimes \radoperator_{R}(\Delta,\Delta)\otimes (\Ext^1_R(\Delta,\Delta))^{\otimes j}\otimes \radoperator_{R}(\Delta,\Delta)\otimes (\Ext^1_R(\Delta,\Delta))^{\otimes k}\to \radoperator_R(\Delta,\Delta)\] 
can be identifed with the dual of the map 
\[\Phi\hookrightarrow \overline{V}\stackrel{\partial_1}{\to} \overline{V}\otimes_A \overline{V}\cong A\otimes \Phi\otimes A\otimes \Phi\otimes A\to \Bbbk Q\otimes \Phi\otimes \Bbbk Q\otimes \Phi\otimes \Bbbk Q\to Q_1^{\otimes i}\otimes \Phi\otimes Q_1^{\otimes j}\otimes \Phi\otimes Q_1^{\otimes k},\]
where  $\partial_1$ denotes the degree $1$ part of the semifree dg algebra corresponding to $(A,V)$  under Ro\u\i ter's equivalence. 
\end{enumerate} 
\end{thmx}

Because of the striking similarity of a regular bocs over a basic algebra with the classical case of a quiver and relations, we suggest to call such a bocs basic.

Our paper is structured as follows. In Section \ref{sec:basics}, we recall basics on $A_\infty$-algebras and coalgebras, in Section \ref{sec:duality} we introduce our conventions for graded duality, before recalling bar and cobar constructions and their compatibility in Section \ref{sec:barcobar}. We then present Kadeishvili's theorem in Section \ref{sec:Kadeishvili1} and give Merkulov's construction of an $A_\infty$-structure on homology from a differential graded algebra and its generalisation by Markl, as suggested by Kontsevich and So\u{\i}bel'man in Section \ref{sec:Kadeishvili2}. For convenience of the reader we have included many proofs to make sure that the sign conventions we use are compatible with each other. In Section \ref{sec:bocses}, we recall some general results on bocses and their module categories. In particular, we explain the equivalence between the category of bocses and that of semifree dg algebras. We recall the notion of regularity for bocses and introduce the new notion of being basic. Furthermore, we introduce quasi-hereditary algebras and recall the main result of \cite{KKO14} characterising quasi-hereditary algebras in terms of bocses. As a new result, we show that the construction of \cite{KKO14} is functorial. In Section \ref{sec:extensions}, we show that different notions to define $\Ext$ in the category of modules over a bocs are all equivalent. In Section \ref{sec:projectiveresolution}, we show how, given a bocs $(A,V)$ and a projective bimodule resolution $\mathcal{P}$ of $V$, one constructs an $A_\infty$-coalgebra structure on $\mathcal{P}$ and give an alternative characterisation of regularity for bocses, see Lemma \ref{regularbocsesequivalence}. In the small Section \ref{sec:uniquenessofexactborelsubalgebras}, we give an independent proof of the corollary of the main result that exact Borel subalgebras of quasi-hereditary algebras arising from bocses are unique up to isomorphism. In Section \ref{sec:explicitcalc}, we explicitly compute the relevant terms of an $A_\infty$-coalgebra structure on $\mathcal{P}$, partially relegating the proofs to Appendices \ref{appendixA}--\ref{appendixE}, and then relate them to the differential on the semifree dg algebra associated to $(A,V)$ in Section \ref{sec:muanddel}. After constructing some well-behaved splitting of projections involved in the presentation of $A$ by quiver and relations in Section \ref{splittings}, we finally prove our main theorems in Section \ref{sec:mainthm}.

{\bf Acknowledgements.} V.M. has been supported by grants EP/K011782/1 and  EP/S017216/1 from the EPSRC. Part of this research was conducted while V.M. was visiting Max-Planck Institute, Bonn, whose hospitality and financial support is gratefully acknowledged. Both authors would like to thank the University of East Anglia, University of Stuttgart and Uppsala University for their hospitality during mutual research visits of both authors.

\section{$A_\infty$-algebras, Kadeishvili's theorem and Merkulov's construction}\label{sec:ainfinity}

In this section, we recall basic definitions and results on $A_\infty$-algebras and coalgebras. In particular, we define linear duality and bar and cobar constructions (Sections \ref{sec:duality} and \ref{sec:barcobar}). Furthermore we recall Kadeishvili's theorem that the homology of an $A_\infty$-algebra carries an $A_\infty$-algebra structure (Section \ref{sec:Kadeishvili1}) and Merkulov's explicit construction which was generalised by Markl, based on ideas of Kontsevich and So\u\i belman  (Section \ref{sec:Kadeishvili2}). We analyse the transfer of counitality of an $A_\infty$-coalgebra object via structure transport (Section \ref{sec:counitality}) and in the end briefly recall the notions of twisted modules, twisted complexes and pretriangulated $A_\infty$ categories (Section \ref{sec:twistedcomplexes}).

For further reading we recommend Keller's survey articles \cite{K01, K02, K06}, the series of articles by Lu, Palmieri, Wu, and Zhang \cite{LPWZ04, LPWZ08, LPWZ09}, Herscovich's  article \cite{Her18}, and the PhD theses of Lef\`evre-Hasegawa \cite{L-H03} and Conner \cite{Con11}. 

\subsection{Basic definitions}\label{sec:basics}
Let $\mathsf{C}$ be a monoidal $\Bbbk$-linear abelian complete and cocomplete category. We denote the unit object by $e$ and the tensor product by $\otimes$. The choices we consider in this paper are $\mathsf{C}=\Mod \mathbb{L}$ where $\mathbb{L}=\prod_{\mathtt{i}=\mathtt{1}}^\mathtt{n} \Bbbk$ is a direct product of fields with tensor product being the tensor product $\otimes_\mathbb{L}$ over $\mathbb{L}$, and $\mathsf{C}=\Mod A\otimes A^{\op}$ with tensor product $\otimes_A$ over $A$.  We will often consider $\mathbb{Z}$-graded objects in $\mathsf{C}$. For a homogeneous element $a$, resp. a homogeneous linear map $f$, we denote by $|a|$, resp. $|f|$, its degree. For homogeneous graded maps we use the Koszul sign convention, i.e. 
\[(f\otimes g)(a\otimes b)=(-1)^{|g|\cdot |a|}f(a)\otimes g(b),\]
where $f$ and $g$ are homogeneous morphisms and $a$ and $b$ are homogeneous elements of the corresponding modules. 
Note that this also implies that $(f\otimes f')(g\otimes g')=(-1)^{|f'|\cdot |g|}(fg\otimes f'g')$ for homogeneous morphisms $f,f',g,g'$. Throughout, we use $a\equiv b$ for $a$ equivalent to $b$ modulo $2$.

An \emphbf{$A_\infty$-algebra object} in $\mathsf{C}$ is a $\mathbb{Z}$-graded object 
\[\mathcal{A}=\bigoplus_{j\in \mathbb{Z}} \mathcal{A}^j\]
in $\mathsf{C}$ endowed with homogeneous morphisms
\[m_n\colon \mathcal{A}^{\otimes n}\to \mathcal{A} \quad (n\geq 1)\]
of degree $2-n$ satisfying the relation
\[ \sum_{r+s+t=n} (-1)^{r+st} \,m_{r+1+t} \,(\id^{\otimes r} \otimes m_s \otimes \id^{\otimes t}) = 0.\]
An $A_\infty$-algebra object in $\Mod \mathbb{L}$ is usually called an \emphbf{$A_\infty$-category} (with finitely many objects). An $A_\infty$-category $\mathcal{A}$ is called \emphbf{locally finite} if $\mathcal{A}^j$ is finite dimensional for every $j\in \mathbb{Z}$. We call an $A_\infty$-algebra object a \emphbf{differential graded algebra object} if $m_n=0$ for $n\geq 3$. 

Every algebra object $\mathcal{A}$ in $\mathsf{C}$ becomes an $A_\infty$-algebra object by placing it in degree $0$, letting $m_2$ be the multiplication map, and setting $m_n=0$ for $n\neq 2$. In particular, the unit object $e$ in $\mathsf{C}$ is an $A_\infty$-algebra object in $\mathsf{C}$. 

A \emphbf{morphism of $A_\infty$-algebra objects} $f\colon \mathcal{A}\to \mathcal{B}$ is a family
\[f_n\colon \mathcal{A}^{\otimes n} \to \mathcal{B}\]
of homogeneous morphisms of degree $1-n$ such that, for $n\geq 1$, we have
\[
\sum_{r+s+t=n} (-1)^{r+st}\, f_{r+1+t}\,(\id^{\otimes r}\otimes m_s \otimes \id^{\otimes t}) 
=
\sum_{q=1}^n \sum_{n=j_1+ \cdots +j_q} (-1)^{\sum_{i=1}^q (q-i)(j_i-1)}\, m_q\, (f_{j_1}\otimes f_{j_2}\otimes \ldots \otimes f_{j_q}).
\]
A morphism of $A_\infty$-algebra objects $f\colon \mathcal{A}\to \mathcal{B}$ is called {\emphbf strict} if $f_n=0$ for $n>1$. 
The \emphbf{composition} of two morphisms $f\colon \mathcal{A}\to \mathcal{A}', g\colon \mathcal{A}'\to \mathcal{A}''$ between $A_\infty$-algebra objects is given by 
\[(f\circ g)_n=\sum_{q=1}^n \sum_{n=j_1+ \cdots +j_q} (-1)^{\sum_{i=1}^r (q-i)(j_i-1)} f_q(g_{j_1}\otimes \dots\otimes g_{j_q}).\]
 A morphism $(f_n)_{n\in \mathbb{N}}$ of $A_\infty$-algebra objects such that $f_1$ is a quasi-isomorphism of complexes is called an \emphbf{$A_\infty$-quasi-isomorphism}. 

An $A_\infty$-algebra object $\mathcal{A}$ in $\mathsf{C}$ is called \emphbf{strictly unital} provided there is a morphism of $A_\infty$-algebras $\eta\colon e\to \mathcal{A}$ such that $m_1\eta=0$, $m_2(\id\otimes \eta)=\id=m_2(\eta\otimes \id)$, and $m_n(\id^{\otimes r}\otimes \eta\otimes \id^{\otimes n-r-1})=0$ for all $n\geq 3, r=1,\dots,n-1$, where by slight abuse of notation we denote the canonical isomorphisms $e\otimes \mathcal{A}\cong \mathcal{A}\cong \mathcal{A}\otimes e$ by equalities. A strictly unital $A_\infty$-algebra is called \emphbf{augmented} if there is a morphism of $A_\infty$-algebras $\varepsilon\colon \mathcal{A}\to e$ such that $\varepsilon\eta=\id_e$.

For later use, we record the following result which is stated in \cite[Lemma 4.2]{LPWZ04} without proof.

\begin{lem}\label{changesigns}
Let $(\mathcal{A}, m_n)$ be an $A_\infty$-algebra object in $\mathsf{C}$. Define $\overline{m}_n:=(-1)^n m_n$ for all $n\in \mathbb{N}$, $f_1(a)=(-1)^{|a|} a$, and $f_n=0$ for all $n>1$. Then $(\mathcal{A}, \overline{m}_n)$ is an $A_\infty$-algebra and $f$ defines a strict $A_\infty$-morphism $(\mathcal{A}, m_n)\to (\mathcal{A}, \overline{m}_n)$. In addition $(\mathcal{A}, m_n)$ is strictly unital (respectively augmented) if and only if $(\mathcal{A}, \overline{m}_n)$ is strictly unital (respectively augmented).
\end{lem}

\begin{proof}
We first check that $(\mathcal{A}, \overline{m}_n)$ forms an $A_\infty$-algebra:
\begin{align*}
\sum_{r+s+t=n} (-1)^{r+st} \overline{m}_{r+1+t} (\id^{\otimes r} \otimes \overline{m}_s \otimes \id^{\otimes t})&=\sum_{r+s+t=n} (-1)^{r+st+r+1+t+s} m_{r+t+1} (\id^{\otimes r}\otimes m_s\otimes \id^{\otimes t})\\
&=(-1)^{n+1} \sum_{r+s+t=n} (-1)^{r+st} m_{r+1+t} (\id^{\otimes r}\otimes m_s \otimes \id^{\otimes t})=0.
\end{align*}
To verify that $f$ defines a $A_\infty$-morphism, we have to check that $f_1m_n=\overline{m}_n(f_1\otimes \cdots \otimes f_1)$:
\begin{align*}
\overline{m}_n(f_1\otimes \dots\otimes f_1)(a_1\otimes \dots\otimes a_n)&=\overline{m}_n(f_1(a_1)\otimes \dots \otimes f_n(a_n))=(-1)^{\sum |a_i|} \overline{m}_n (a_1\otimes \dots \otimes a_n)\\
&=(-1)^{n+\sum |a_i|} m_n(a_1\otimes \dots\otimes a_n)=(-1)^{2-n+\sum |a_i|} m_n(a_1\otimes \dots\otimes a_n)\\
&=f_1m_n(a_1\otimes \dots\otimes a_n). 
\end{align*}
To see that $f$ is in fact a strict isomorphism, one can either do a similar calculation for its inverse or one uses the fact that an $A_\infty$-morphism is an isomorphism if and only if $f_1$ is an isomorphism (cf. Lemma \ref{structuretransport}). 

From $\overline{m}_2=m_2$ and $m_n(a_1\otimes \dots\otimes a_n)=0 \Leftrightarrow \overline{m}_n(a_1\otimes \dots \otimes a_n)=0$, it follows immediately that $(\mathcal{A}, m_n)$ is strictly unital if and only if $(\mathcal{A}, \overline{m}_n)$ is. The statement about augmentation is also clear.
\end{proof}

We now recall the dual notions. An \emphbf{$A_\infty$-coalgebra object} in $\mathsf{C}$ is a $\mathbb{Z}$-graded object 
\[\mathcal{C}=\bigoplus_{j\in \mathbb{Z}} \mathcal{C}^j\]
in $\mathsf{C}$ with homogeneous morphisms 
\[\mu_n\colon \mathcal{C}\to \mathcal{C}^{\otimes n} \quad (n\geq 1)\]
of degree $2-n$ such that the morphism 
\[(\mu_n)_{n\in \mathbb{Z}}\colon \mathcal{C}\to \prod \mathcal{C}^{\otimes n}\]
factors through the canonical inclusion $\bigoplus \mathcal{C}^{\otimes n}\to \prod \mathcal{C}^{\otimes n}$ and such that 
\[\sum_{r+s+t=n} (-1)^{rs+t} (\id^{\otimes r}\otimes \mu_s\otimes \id^{\otimes t})\mu_{r+1+t}=0,\]
We call an $A_\infty$-coalgebra object a \emphbf{differential graded coalgebra object} if $\mu_n=0$ for $n\geq 3$. 

Again, every coalgebra object in $\mathsf{C}$ defines an $A_\infty$-coalgebra object in $\mathsf{C}$ by placing it in degree $0$, letting $\mu_2$ be the comultiplication, and setting $\mu_n=0$ for $n\neq 2$. In particular, the unit object $e$ is an $A_\infty$-coalgebra object in $\mathsf{C}$. 

A \emphbf{morphism of $A_\infty$-coalgebra objects} $f\colon \mathcal{C}\to \mathcal{D}$ is a family
\[f_n\colon \mathcal{C}\to \mathcal{D}^{\otimes n}\]
of homogeneous morphisms of degree $1-n$ such that the morphism
\[(f_n)_{n\in \mathbb{Z}}\colon \mathcal{C}\to \prod \mathcal{D}^{\otimes n}\]
factors through the canonical inclusion $\bigoplus \mathcal{D}^{\otimes n}\to \prod \mathcal{D}^{\otimes n}$ and such that 
\[\sum_{r+s+t=n} (-1)^{rs+t}(\id^{\otimes r}\otimes \mu_s\otimes \id^{\otimes t})f_{r+1+t}=\sum_{q=1}^n\sum_{n=j_1+\cdots +j_q} (-1)^{\sum_{i=1}^q(i-1)(j_i+1)}(f_{j_1}\otimes f_{j_2}\otimes \dots\otimes f_{j_q})\mu_q.\] 
A morphism of $A_\infty$-coalgebra objects is called {\emphbf strict} if $f_n=0$ for all $n>1$. 
Furthermore, the composition of two morphisms of $A_\infty$-coalgebra objects $f\colon \mathcal{C}\to \mathcal{C}'$ and $g\colon \mathcal{C}'\to \mathcal{C}''$ is given by
\[(g\circ f)_n:=\sum_{q=1}^n\sum_{n=j_1+\cdots +j_q} (-1)^{\sum_{i=1}^q(i-1)(j_i+1)} (g_{j_1}\otimes \dots\otimes g_{j_r})f_r.\]

A morphism $(f_n)_{n\in \mathbb{N}}$ of $A_\infty$-coalgebra objects such that $f_1$ is a quasi-isomorphism is called an \emphbf{$A_\infty$-quasi-isomorphism}.

An $A_\infty$-coalgebra object $\mathcal{C}$ is called \emphbf{strictly counital} provided there exists a morphism of $A_\infty$-algebras $\tau\colon \mathcal{C}\to e$ such that $\tau \mu_1=0$, $(\id\otimes \tau)\mu_2=\id=(\tau\otimes \id)\mu_2$, and $(\id^{\otimes r}\otimes \tau\otimes \id^{\otimes n-r-1})\mu_n=0$ for all $n>3, r=1,\dots,n-1$. A strictly counital $A_\infty$-coalgebra object $\mathcal{C}$ is called \emphbf{coaugmented} if there is a morphism of $A_\infty$-coalgebra objects $\eta\colon e\to \mathcal{C}$ such that $\tau\eta=\id_e$. 

For future use, we record a technical lemma, which generalises the defining equation for two-step $A_\infty$-comultiplications $\mathcal{C}\to \mathcal{C}^{\otimes n}$ to a similar formula for two-step  comultiplications $\mathcal{C}^{\otimes m}\to \mathcal{C}^{\otimes n}$.

\begin{lem}\label{generalised}
Let $\mathcal{C}$ be an $A_\infty$-coalgebra object in some category $\mathsf{C}$. Then, for each pair of natural numbers $m\leq n$,
\[\sum_{p+q+r=n}\sum_{\substack{a+b+c=p+r+1\\a+1+c=m}}(-1)^{pq+r+ab+c}(\id^{\otimes p}\otimes \mu_q\otimes \id^{\otimes r})(\id^{\otimes a}\otimes \mu_b\otimes \id^{\otimes c})=0\]
\end{lem}

\begin{proof}
We write
\begin{align*}
&\sum_{p+q+r=n}\sum_{\substack{a+b+c=p+r+1\\a+1+c=m}}(-1)^{pq+r+ab+c}(\id^{\otimes p}\otimes \mu_q\otimes \id^{\otimes r})(\id^{\otimes a}\otimes \mu_b\otimes \id^{\otimes c})\\
&=\sum_{p+q+r=n}\sum_{\substack{a+b+c=p+r+1\\a+1+c=m\\a\leq p\\c\leq r}}(-1)^{pq+r+ab+c}\id^{\otimes a}\otimes ((\id^{\otimes (p-a)}\otimes \mu_q\otimes \id^{\otimes (r-c)})\mu_b)\otimes \id^{\otimes c}
\\&+\sum_{p+q+r=n}\sum_{\substack{a+b+c=p+r+1\\a+1+c=m\\\text{(either) } a>p \text{ or } c>r}}(-1)^{pq+r+ab+c}(\id^{\otimes p}\otimes \mu_q\otimes \id^{\otimes r})(1^{\otimes a}\otimes \mu_b\otimes \id^{\otimes c})
\end{align*}
and consider the two summands separately. For the first summand, we substitute $x=p-a, y=r-c$. This yields $b=x+y+1$ and $p+q+r=x+a+q+y+c$ and thus $x+q+y=n-a-c$. Therefore, 
\begin{align*}
&\sum_{p+q+r=n}\sum_{\substack{a+b+c=p+r+1\\a+1+c=m\\a\leq p\\c\leq r}}(-1)^{pq+r+ab+c}\id^{\otimes a}\otimes ((\id^{\otimes (p-a)}\otimes \mu_q\otimes \id^{\otimes (r-c)})\mu_b)\otimes \id^{\otimes c}\\
&=\sum_{\substack{a,c\geq 0\\a+1+c=m}} (-1)^{a(n-a-c+1)}\id^{\otimes a}\otimes \left( \sum_{x+q+y=n-a-c} (-1)^{xq+y} \id^{\otimes x}\otimes \mu_q\otimes \id^{\otimes y} \right)\otimes \id^{\otimes c}=0,
\end{align*}
where the vanishing follows from the $A_\infty$-relations for $\mu$. For the second summand, note that the cases $a>p$ and $c>r$ are mutually exclusive, as otherwise $b$ would have be negative due to the equality $a+b+c=p+r+1$.  
Considering only the sum over the terms with $a>p$ we obtain
\begin{align*}
&\sum_{p+q+r=n}\sum_{\substack{a+b+c=p+r+1\\a+1+c=m\\ a>p}}(-1)^{pq+r+ab+c}(\id^{\otimes p}\otimes \mu_q\otimes \id^{\otimes r})(\id^{\otimes a}\otimes \mu_b\otimes \id^{\otimes c})\\
&=\sum_{p+q+r=n}\sum_{\substack{a+b+c=p+r+1\\a+1+c=m\\ a>p}}(-1)^{pq+r+ab+c+bq}(\id^{\otimes (a+q-1)}\otimes \mu_b\otimes \id^{\otimes c})(\id^{\otimes p}\otimes \mu_q\otimes \id^{\otimes (m-p-1)})
\end{align*}
Using the substitution $p'=a+q-1, q'=b, r'=c, a'=p, b'=q, c'=m-p-1$ we obtain $p'+q'+r'=a+b+c+q-1=n$, $a'+b'+c'=p+q+(m-p-1)=m+q-1=a+c+q=p'+r'+1$, and $a'+c'+1=m$. Furthermore, $a>p$ if and only if $m-p-1>c$ if and only if $c'>r'$. Therefore, the above sum can be rewritten as
\[
\sum_{p'+q'+r'=n}\sum_{\substack{a'+b'+c'=p'+r'+1\\a'+1+c'=m\\c'>r'}} (-1)^v (\id^{\otimes p'}\otimes \mu_{q'}\otimes \id^{\otimes r'})(\id^{\otimes a'}\otimes \mu_{b'}\otimes \id^{\otimes c'})
\]
where we claim that $v\equiv p'q'+r'+a'b'+c'+1$. Indeed,
\begin{align*}
p'q'+r'+a'b'+c'-1&=(a+q-1)b+c+pq+(m-p-1)-1\\
&=pq+ab+c+bq-b+m-p-2\\
&=pq+ab+c+bq-b+(a+1+c)-p-2\\
&\equiv pq+ab+c+bq+r.
\end{align*}
Once can now see that the terms with $a>p$ cancel with the terms with $c>r$. The claim follows. \qedhere
\end{proof}

\subsection{Duality}\label{sec:duality}
We now show that for $\mathsf{C}=\Mod \mathbb{L}$, under suitable finiteness assumptions, one can switch between $A_\infty$-categories and $A_\infty$-cocategories using linear duality. 
We use the simple-preserving duality given by the identity on $\mathbb{L}$, which is an antiautomorphism by commutativity, to obtain a covariant duality $\#$ on the category of $\mathbb{L}$-$\mathbb{L}$-bimodules.
More precisely, 
\[e_\mathtt{j}(M^\#)_se_\mathtt{i} = \Hom_\mathbb{L}(e_\mathtt{j}M_{-s}e_\mathtt{i}, \Bbbk).\]
Let $\iota_{M_1,\dots,M_{n}}\colon M_1^\#\otimes \dots\otimes M_n^\#\to (M_1\otimes \dots\otimes M_n)^\#$ be the map defined by $\phi_1\otimes \dots\otimes \phi_n\mapsto \phi_1\otimes \dots\otimes \phi_n$ where for the latter the Koszul sign rule is taken into account, i.e. $\phi_1\otimes \dots\otimes \phi_n$ is defined by 
\[(\phi_1\otimes\dots\otimes \phi_n)(x_1\otimes \dots\otimes x_n)=(-1)^{\sum_{j<i} |\phi_i|\cdot |x_j|} \phi_1(x_1)\otimes \dots\otimes \phi_n(x_n).\]
If $M_1=\cdots=M_n$ and $M_1$ is clear from the context,  we simply write $\iota_n$ for $\iota_{M_1,\dots,M_{n}}$. Note that the maps $\iota_{M_1,\dots,M_n}$ are monomorphisms and that they are isomorphisms if and only if the vector spaces in question are finite dimensional. 

Notice that $\Hom_{\mathbb{L}}(M,\mathbb{L})$ also carries a natural $\mathbb{L}$-$\mathbb{L}$-bimodule structure given by $e_\mathtt{j}\Hom_{\mathbb{L}}(M,\mathbb{L})e_\mathtt{i} = \Hom_{\mathbb{L}}(Me_\mathtt{j},\mathbb{L}e_\mathtt{i}).$ We can again twist this by the antiautomorphism given by the identity, to define a contravariant $\mathbb{L}$-duality $\flat$ by  $e_\mathtt{j}(M^\flat)_se_\mathtt{i} =  \Hom_{\mathbb{L}}(M_{-s}e_\mathtt{i},\mathbb{L}e_\mathtt{j})$. 

\begin{lem}\label{flathash}
We claim that there is an isomorphism of $\mathbb{L}$-$\mathbb{L}$-bimodules $M^\# \cong M^\flat$ , compatible with tensor product over $\mathbb{L}$.
\end{lem}

\begin{proof}
The  isomorphism of $\mathbb{L}$-$\mathbb{L}$-bimodules is given by  
\[e_\mathtt{j}M^\flat  e_\mathtt{i} =  \Hom_{\mathbb{L}}(Me_\mathtt{i},\mathbb{L}e_\mathtt{j}) =  \Hom_{\mathbb{L}}(e_\mathtt{j}Me_\mathtt{i},\mathbb{L}e_\mathtt{j}) \cong \Hom_{\Bbbk}(e_\mathtt{j}Me_\mathtt{i},\Bbbk) = e_\mathtt{j}M^\# e_\mathtt{i}.\]

This is compatible with tensor product over $\mathbb{L}$ as follows. Let $M,N$ be two $\mathbb{L}$-$\mathbb{L}$-bimodules, $\phi\in M^\flat$,  $\psi\in N^\flat$,$x\in M$ and $y\in N$. Then we can identify $M^\flat \otimes_\mathbb{L} N^\flat$ with $(M\otimes_\mathbb{L} N)^\flat$ by defining
$(\phi\otimes \psi)(x\otimes y) = (-1)^{|\psi||x|}\phi(x\psi(y))$.
Identifying $\phi, \psi$ with elements in $M^\#,N^\#$ respectively, translates this into $(\phi\otimes \psi)(x\otimes y) = (-1)^{|\psi||x|}\phi(x)\psi(y)$ since it considers $\psi(y)$ as a scalar, and both expressions again coincide under the identification of $\flat$ and $\#$.
\end{proof}

For a homogeneous map $f\colon M\to N$, $f^\#\colon N^\#\to M^\#$ is defined as $f^\#(\phi):=(-1)^{|f|\cdot |\phi|}\phi\circ f$ for every homogeneous element $\phi\in N^\#$. Note that with this definition, $(gf)^\#=(-1)^{|f|\cdot |g|} f^\#g^\#$ as 
\begin{align*}
(gf)^\#(\phi)&=(-1)^{|gf|\cdot |\phi|} \phi\circ (gf)=(-1)^{(|g|+|f|)|\phi|} \phi gf\\
&=(-1)^{|g|\cdot |f|}(-1)^{|g|\cdot |\phi|+|f|\cdot |\phi g|} (\phi g)\circ f=(-1)^{|g|\cdot |f|}(-1)^{|g|\cdot |\phi|} f^\#(\phi g)\\
&=(-1)^{|g|\cdot |f|}f^\#g^\#(\phi).
\end{align*}

The following notion, which we adopt from \cite[Definition 4.1.3]{Boo21}, is convenient to give suitable conditions for duality statements. 

\begin{defn}
A graded $\mathbb{L}$-module is called \emphbf{raylike} if it is locally finite dimensional and either bounded above or bounded below. 
\end{defn}

As laid out in \cite{Boo21}, the main property of such modules is that the tensor product of raylike modules is still locally finite dimensional which is not true in general if the assumption is dropped. The following lemma providing a duality between certain $A_\infty$-categories and $A_\infty$-cocategories is well known, see e.g. \cite[Section 2.3]{Her18} where the statement about $A_\infty$-categories is stated without proof.

\begin{lem}\label{dualcategoriescocategories}
\begin{enumerate}[(i)]
\item\label{lem2.1.1}  For an $A_\infty$-cocategory $\mathcal{C}$,  we obtain an $A_\infty$-category $\mathcal{C}^{\#}$ with multiplications given by 
\[m_n= (-1)^n (\mu_n)^\# \iota_n.\]
For a morphism $f\colon \mathcal{C}\to \mathcal{D}$ of $A_\infty$-cocategories, $f^\#$ given by 
\[(f^\# )_n= (f_n)^\# \iota_n\]
defines a morphism of $A_\infty$-categories.
\item\label{lem2.1.2} 
Dually, for a raylike $A_\infty$-category $\mathcal{A}$ such that only finitely many $m_n\colon \mathcal{A}^{\otimes n}\to \mathcal{A}$ are non-zero, we obtain an $A_\infty$-cocategory $\mathcal{A}^{\#}$ with comultiplications given by 
\[\mu_n=(-1)^n\iota_n^{-1}\circ m_n^\#\] 
for $i\geq 1$.
For a morphism $g\colon \mathcal{A}\to \mathcal{B}$ of $A_\infty$-categories such that only finitely many $g_n\colon \mathcal{A}^{\otimes n}\to \mathcal{B}$ are non-zero, $g^\#$ defined by
\[(g^\#)_n=\iota_n^{-1}(g_n)^\#\] 
is a morphism of $A_\infty$-cocategories.
\end{enumerate}
\end{lem}

\begin{proof}
\eqref{lem2.1.1}
Let $\mathcal{C}$ be an $A_\infty$-cocategory and let $m_n$ on $\mathcal{C}^{\#}$ be as defined above. 

We first claim that the diagram 
\begin{equation}\label{eqn211}
\begin{tikzcd}
(\mathcal{C}^{\#})^{\otimes r}\otimes (\mathcal{C}^{\otimes s})^{\#}\otimes  (\mathcal{C}^{\#})^{\otimes t}\arrow{rr}{\iota_r\otimes \id\otimes \iota_t}\ar{dd}{\id^{\otimes r}\otimes \mu_s^{\#}\otimes \id^{\otimes r}}
&&(\mathcal{C}^{\otimes r})^{\#}\otimes (\mathcal{C}^{\otimes s})^{\#}\otimes (\mathcal{C}^{\otimes t})^{\#}\arrow{d}{\iota_{\mathcal{C}^{\otimes r}, \mathcal{C}^{\otimes s}, \mathcal{C}^{\otimes t}}}\\
&&(\mathcal{C}^{\otimes r}\otimes\mathcal{C}\otimes\mathcal{C}^{\otimes t})^{\#}\arrow{d}{(\id^{\otimes r}\otimes \mu_s\otimes id^{\otimes t})^{\#}}\\
(\mathcal{C}^{\#})^{\otimes r}\otimes \mathcal{C}^{\#}\otimes  (\mathcal{C}^{\#})^{\otimes t}\arrow{rr}{\iota_{r+1+t}}&& (\mathcal{C}^{\otimes r}\otimes\mathcal{C}^{\otimes s}\otimes\mathcal{C}^{\otimes t})^{\#}
\end{tikzcd}
\end{equation}
commutes. Indeed, we compute, on the one hand,
\begin{equation*}\begin{split}
\iota_{r+1+t}&(\id^{\otimes r}\otimes \mu_s^{\#}\otimes \id^{\otimes t})(\phi_1\otimes\cdots \otimes \phi_n)\\
&=(-1)^{s\sum_{j=1}^r|\phi_j|}\iota_{r+1+t}(\phi_1\otimes\cdots \otimes \phi_r\otimes \mu_s^{\#}(\phi_{r+1}\otimes\cdots \otimes \phi_{r+s})\otimes \phi_{r+s+1}\otimes\cdots \otimes \phi_n)\\
&=(-1)^{s\sum_{j=1}^{r+s}|\phi_j|}\iota_{r+1+t}(\phi_1\otimes\cdots \otimes \phi_r\otimes (\phi_{r+1}\otimes\cdots \otimes \phi_{r+s})\mu_s\otimes \phi_{r+s+1}\otimes\cdots \otimes \phi_n)\\
&=(-1)^{s\sum_{j=1}^{n}|\phi_j|}\iota_{r+1+t}(\phi_1\otimes\cdots \otimes \phi_n)(\id^{\otimes r}\otimes \mu_s\otimes \id^{\otimes t})\\
&=(-1)^{s\sum_{j=1}^{n}|\phi_j|}(\phi_1\otimes\cdots \otimes \phi_n)(\id^{\otimes r}\otimes \mu_s\otimes \id^{\otimes t})
\end{split}\end{equation*}
and, on the other hand,
\begin{equation*}\begin{split}
(\id^{\otimes r}\otimes &\mu_s\otimes \id^{\otimes t})^{\#} \iota_{\mathcal{C}^{\otimes r}, \mathcal{C}^{\otimes s}, \mathcal{C}^{\otimes t}}(\iota_r\otimes \id\otimes \iota_t) (\phi_1\otimes\cdots \otimes \phi_n)\\
&=(-1)^{s\sum_{j=1}^{n}|\phi_j|} \iota_{\mathcal{C}^{\otimes r}, \mathcal{C}^{\otimes s}, \mathcal{C}^{\otimes t}}(\iota_r\otimes \id\otimes \iota_t) (\phi_1\otimes\cdots \otimes \phi_n)(\id^{\otimes r}\otimes \mu_s\otimes \id^{\otimes t})\\
&=(-1)^{s\sum_{j=1}^{n}|\phi_j|} (\phi_1\otimes\cdots \otimes \phi_n)(\id^{\otimes r}\otimes \mu_s\otimes \id^{\otimes t})
\end{split}\end{equation*}
which proves the claim.

We then check
\begin{equation*}\begin{split}
\sum (-1)^{r+st}m_{r+1+t}&(\id^{\otimes r}\otimes m_s\otimes \id^{\otimes t})= \sum (-1)^{r+st+n+1}\mu_{r+1+t}^{\#}\iota_{r+1+t}(\id^{\otimes r}\otimes \mu_s^{\#}\iota_s\otimes \id^{\otimes t})\\
&= \sum (-1)^{r+st+n+1}\mu_{r+1+t}^{\#}\iota_{r+1+t}(\id^{\otimes r}\otimes \mu_s^{\#}\otimes \id^{\otimes t})(\id^{\otimes r}\otimes \iota_s\otimes \id^{\otimes t})\\
&\stackrel{(\star)}=\sum (-1)^{r+st+n+1}\mu_{r+1+t}^{\#}(\id^{\otimes r}\otimes \mu_s\otimes \id^{\otimes t})^{\#}\iota_{\mathcal{C}^{\otimes r}, \mathcal{C}^{\otimes s}, \mathcal{C}^{\otimes t}}(\iota_r\otimes \id\otimes \iota_t)(\id^{\otimes r}\otimes \iota_s\otimes \id^{\otimes t})\\
&=\sum (-1)^{r+st+n+1+(r+1+t)s}((\id^{\otimes r}\otimes \mu_s\otimes \id^{\otimes t})\mu_{r+1+t})^{\#}\iota_n\\
&=\sum (-1)^{r+st+n+1+rs+s+st}((\id^{\otimes r}\otimes \mu_s\otimes \id^{\otimes t})\mu_{r+1+t})^{\#}\iota_n\\
&=-\sum (-1)^{rs+t}((\id^{\otimes r}\otimes \mu_s\otimes \id^{\otimes t})\mu_{r+1+t})^{\#}\iota_n\\&=0
\end{split}\end{equation*}
where for $(\star)$ we use the commutativity of \eqref{eqn211}. Therefore the $m_n$ define a structure of an $A_\infty$-category on $\mathcal{A}^\#$.

Now suppose $f$ is a morphism of $A_\infty$-cocategories. We claim that the diagram
\begin{equation}\label{eqn212}
\begin{tikzcd}
(\mathcal{C}^{\otimes j_1})^{\#}\otimes \cdots \otimes (\mathcal{C}^{\otimes j_q})^{\#} \arrow{rr}{\iota_{\mathcal{C}^{\otimes j_1},\cdots,\mathcal{C}^{\otimes j_q}}}\arrow{d}{(f_{j_1})^{\#}\otimes \cdots \otimes (f_{j_q})^{\#}}
&& (\mathcal{C}^{\otimes j_1}\otimes \cdots \otimes \mathcal{C}^{\otimes j_q})^{\#} \arrow{d}{(f_{j_1}\otimes \cdots \otimes f_{j_q})^{\#}}\\
(\mathcal{C}^{\#})^{\otimes q}\arrow{rr}{\iota_q}&&(\mathcal{C}^{\otimes q})^{\#}
\end{tikzcd}
\end{equation}
commutes. Indeed, we compute, on the one hand,
\begin{equation*}\begin{split}
\iota_q((f_{j_1})^{\#}&\otimes \cdots \otimes (f_{j_q})^{\#})(\phi_1\otimes \cdots\otimes \phi_n)\\
&=(-1)^{\sum_{k=1}^q (j_k-1)\sum_{l=1}^{\sum_{s=1}^{k-1}j_s}|\phi_l|}\iota_q((f_{j_1})^{\#}(\phi_1\otimes\cdots\otimes \phi_{j_1})\otimes \cdots\otimes   (f_{j_q})^{\#}(\phi_{\sum_{r=1}^{q-1}j_r +1}\otimes\cdots\otimes \phi_{n}) )\\
&=(-1)^{\sum_{k=1}^q (j_k-1)\sum_{l=1}^{\sum_{s=1}^{k}j_s}|\phi_l|}\iota_q((\phi_1\otimes\cdots\otimes \phi_{j_1})f_{j_1}\otimes \cdots\otimes   (\phi_{\sum_{r=1}^{q-1}j_r +1}\otimes\cdots\otimes \phi_{n})f_{j_q})\\
&=(-1)^{\sum_{k=1}^q (j_k-1)\sum_{l=1}^{n}|\phi_l|}\iota_q(\phi_1\otimes \cdots\otimes \phi_{n})(f_{j_1}\otimes \cdots \otimes f_{j_q})\\
&=(-1)^{(n-q)\sum_{l=1}^{n}|\phi_l|}(\phi_1\otimes \cdots\otimes \phi_{n})(f_{j_1}\otimes \cdots \otimes f_{j_q})
\end{split}\end{equation*}
and, on the other hand,
\begin{equation*}\begin{split}
(f_{j_1}\otimes  &\cdots \otimes f_{j_q})^{\#}\iota_{\mathcal{C}^{\otimes j_1},\cdots, \mathcal{C}^{\otimes j_q}}(\phi_1\otimes \cdots\otimes \phi_n)\\
&=(-1)^{(n-q)\sum_{l=1}^{n}|\phi_l|}\iota_{\mathcal{C}^{\otimes j_1},\cdots,  \mathcal{C}^{\otimes j_q}}(\phi_1\otimes \cdots\otimes \phi_n)(f_{j_1}\otimes \cdots \otimes f_{j_q})\\
&=(-1)^{(n-q)\sum_{l=1}^{n}|\phi_l|}(\phi_1\otimes \cdots\otimes \phi_n)(f_{j_1}\otimes \cdots \otimes f_{j_q}),
\end{split}\end{equation*}
which proves the claim.

Then we compute, on the one hand,
\begin{equation*}\begin{split}
\sum &(-1)^{r+st} (f^{\#})_{r+1+t}(\id^{\otimes r}\otimes m_s\otimes \id^{\otimes t})
=\sum (-1)^{r+st+s} (f_{r+1+t})^{\#}\iota_{r+1+t}(\id^{\otimes r}\otimes \mu_s^{\#}\iota_s\otimes \id^{\otimes t})\\
&=\sum (-1)^{r+st+s} (f_{r+1+t})^{\#}\iota_{r+1+t}(\id^{\otimes r}\otimes \mu_s^{\#}\otimes \id^{\otimes t})(\id^{\otimes r}\otimes\iota_s\otimes \id^{\otimes t})\\
&\stackrel{(\star)}{=}\sum (-1)^{r+st+s} (f_{r+1+t})^{\#}(\id^{\otimes r}\otimes \mu_s\otimes \id^{\otimes t})^{\#}\iota_{\mathcal{C}^{\otimes r}, \mathcal{C}^{\otimes s}, \mathcal{C}^{\otimes t}}(\iota_r\otimes \id\otimes \iota_t)(\id^{\otimes r}\otimes\iota_s\otimes \id^{\otimes t})\\
&=\sum (-1)^{r+st+s+(r+t)s)} ((\id^{\otimes r}\otimes \mu_s\otimes \id^{\otimes t})f_{r+1+t})^{\#}\iota_n\\
&=\sum (-1)^{r+s+rs} ((\id^{\otimes r}\otimes \mu_s\otimes \id^{\otimes t})f_{r+1+t})^{\#}\iota_n\\
&=(-1)^n\sum (-1)^{rs+t} ((\id^{\otimes r}\otimes \mu_s\otimes \id^{\otimes t})f_{r+1+t})^{\#}\iota_n\\
&=(-1)^n\sum (-1)^{\sum_{i=1}^q(i-1)(j_i+1)} ((f_{j_1}\otimes \cdots \otimes f_{j_q})\mu_q)^{\#}\iota_n,
\end{split}\end{equation*}
where at $(\star)$ we use again the commutativity of \eqref{eqn211}; and, on the other hand,
\begin{equation*}\begin{split}
\sum (-1)^{\sum_{i=1}^q(q-i)(j_i-1)}
&m_q((f^{\#})_{j_1}\otimes \cdots\otimes (f^{\#})_{j_q})\\
&= \sum (-1)^{q+\sum_{i=1}^q(q-i)(j_i-1)}\mu_q^{\#}\iota_q((f_{j_1})^{\#}\iota_{j_1}\otimes \cdots \otimes(f_{j_q})^{\#}\iota_{j_q})\\
&\stackrel{(\ddagger)}{=} \sum (-1)^{q+\sum_{i=1}^q(q-i)(j_i-1)}\mu_q^{\#}\iota_q((f_{j_1})^{\#}\otimes \cdots\otimes (f_{j_q})^{\#})(\iota_{j_1}\otimes \cdots \otimes\iota_{j_q}))\\
&=\sum (-1)^{q+\sum_{i=1}^q(q-i)(j_i-1)}\mu_q^{\#}(f_{j_1}\otimes \cdots \otimes f_{j_q})^{\#}\iota_{\mathcal{C}^{\otimes j_1},\cdots, \mathcal{C}^{\otimes j_q}}(\iota_{j_1}\otimes \cdots \otimes\iota_{j_q})\\
&=\sum (-1)^{q+\sum_{i=1}^q(q-i)(j_i-1) +(n-q)q}((f_{j_1}\otimes \cdots \otimes f_{j_q})\mu_q)^{\#}\iota_n
\end{split}\end{equation*}
where for $(\ddagger)$ we use the commutativity of \eqref{eqn212}. Comparing signs, we obtain, on the one hand,
\[(-1)^{n+\sum_{i=1}^q(i-1)(j_i+1)}= (-1)^{n+(\sum_{i=1}^q ij_i) + (\sum_{i=1}^q i) -n-q} = (-1)^{(\sum_{i=1}^q ij_i) + (\sum_{i=1}^q i) -q} \]
and, on the other hand,
\[(-1)^{q+\sum_{i=1}^q(q-i)(j_i-1) +(n-q)q} = (-1)^{q+ nq- (\sum_{i=1}^q ij_i) -q^2+(\sum_{i=1}^q i) +nq-q^2} =(-1)^{q- (\sum_{i=1}^q ij_i) +(\sum_{i=1}^q i) } \]
so the signs agree, and $f^{\#}$ is indeed a morphism of $A_\infty$-categories.

To check that this assignment on morphisms respects composition of morphisms, we compute, on the one hand,
\begin{equation*}\begin{split}
((f\circ g)^{\#})_n& = ((f\circ g)_n)^{\#}\iota_n = \sum (-1)^{\sum_{i=1}^q(i-1)(j_i+1)}((f_{j_1}\otimes \cdots \otimes f_{j_q})g_q)^{\#}\iota_n
\end{split}\end{equation*}
while, on the other hand,
\begin{equation*}\begin{split}
(g^{\#}\circ f^{\#})_n &= \sum(-1)^{\sum_{i=1}^q(q-i)(j_i-1)}(g^{\#})_q((f^{\#})_{j_1}\otimes \cdots \otimes(f^{\#})_{j_q})\\
&=\sum(-1)^{\sum_{i=1}^q(q-i)(j_i-1)}(g_q)^{\#}\iota_q((f_{j_1})^{\#}\iota_{j_1}\otimes \cdots \otimes(f_{j_q})^{\#}\iota_{j_q})\\
&=\sum(-1)^{\sum_{i=1}^q(q-i)(j_i-1)}(g_q)^{\#}\iota_q((f_{j_1})^{\#}\otimes \cdots \otimes(f_{j_q})^{\#})(\iota_{j_1}\otimes \cdots \otimes\iota_{j_q})\\
&\stackrel{(\ddagger)}{=}\sum(-1)^{\sum_{i=1}^q(q-i)(j_i-1)}(g_q)^{\#}(f_{j_1}\otimes \cdots \otimes f_{j_q})^{\#}\iota_{\mathcal{C}^{\otimes j_1},\cdots, \mathcal{C}^{\otimes j_q}}(\iota_{j_1}\otimes \cdots \otimes\iota_{j_q})\\
&=\sum(-1)^{\sum_{i=1}^q(q-i)(j_i-1) +(n-q)(q-1)}((f_{j_1}\otimes \cdots \otimes f_{j_q})g_q)^{\#}\iota_n
\end{split}\end{equation*}
where for $(\ddagger)$ we use the commutativity of \eqref{eqn212}. Comparing signs, we obtain, on the one hand,
\[(-1)^{\sum_{i=1}^q(i-1)(j_i+1)} = (-1)^{ (\sum_{i=1}^qij_i)+  (\sum_{i=1}^qi) -n-q }\]
and, on the other hand,
\[(-1)^{\sum_{i=1}^q(q-i)(j_i-1) +(n-q)(q-1)} = (-1)^{nq-(\sum_{i=1}^qij_i) -q^2+  (\sum_{i=1}^qi) +nq-q^2-n+q} =  (-1)^{-(\sum_{i=1}^qij_i) +  (\sum_{i=1}^qi) -n+q} \]
so the two expressions agree, as desired.

\eqref{lem2.1.2}
Let $\mathcal{A}$ be a locally finite $A_\infty$-category such that only finitely many $m_n\colon \mathcal{A}^{\otimes n}\to \mathcal{A}$ are non-zero. As a first step, we show that the diagram 
\begin{equation}\label{eqn213}
\begin{tikzcd}[column sep=4cm]
(\mathcal{A}^\#)^{\otimes r}\otimes \mathcal{A}^\#\otimes (\mathcal{A}^\#)^{\otimes t}\arrow{d}{\id^{\otimes r}\otimes m_s^\#\otimes \id^{\otimes t}}\arrow{r}{\iota_{r+1+t}} &(\mathcal{A}^{\otimes r}\otimes \mathcal{A}\otimes \mathcal{A}^{\otimes t})^\#\arrow{dd}{(\id^{\otimes r}\otimes m_s\otimes \id^{\otimes t})^{\#}}\\
(\mathcal{A}^{\#})^{\otimes r}\otimes (\mathcal{A}^{\otimes s})^\#\otimes (\mathcal{A}^\#)^{\otimes t}\arrow{d}{\iota_r\otimes \id\otimes \iota_t}\\
(\mathcal{A}^{\otimes r})^\#\otimes (\mathcal{A}^{\otimes s})^\#\otimes (\mathcal{A}^{\otimes t})^\#\arrow{r}{\iota_{\mathcal{A}^{\otimes r},\mathcal{A}^{\otimes s},\mathcal{A}^{\otimes t}}} &(\mathcal{A}^{\otimes r}\otimes \mathcal{A}^{\otimes s}\otimes \mathcal{A}^{\otimes t})^\#
\end{tikzcd}
\end{equation}
commutes.
Indeed,
\begin{align*}
&(\id^{\otimes r}\otimes m_s\otimes \id^{\otimes t})^\# \iota_{r+1+t}(\phi_1\otimes\dots\otimes \phi_{r+1+t})\\
=&(\id^{\otimes r}\otimes m_s\otimes \id^{\otimes t})^\#(\phi_1\otimes \dots\otimes \phi_{r+1+t})\\
=&(-1)^{(2-s)\sum_{j=1}^{r+1+t}{|\phi_j|}} (\phi_1\otimes \dots\otimes \phi_{r+1+t})(\id^{\otimes r}\otimes m_s\otimes \id^{\otimes t}),
\end{align*}
while, on the other hand,
\begin{align*}
&\iota_{\mathcal{A}^{\otimes r},\mathcal{A}^{\otimes s},\mathcal{A}^{\otimes t}}(\iota_r\otimes \id\otimes \iota_t)(\id^{\otimes r}\otimes m_s^\#\otimes \id^{\otimes t})(\phi_1\otimes \dots\otimes \phi_{r+1+t})\\
=&(-1)^{(2-s)\sum_{j=1}^r |\phi_j|}\iota_{\mathcal{A}^{\otimes r},\mathcal{A}^{\otimes s}, \mathcal{A}^{\otimes t}}(\iota_r\otimes \id\otimes \iota_t)(\phi_1\otimes \dots\otimes \phi_r\otimes m_s^\#\phi_{r+1}\otimes \phi_{r+2}\otimes \dots\otimes \phi_{r+1+t})\\
=&(-1)^{(2-s)\sum_{j=1}^{r+1}|\phi_j|}\iota_{\mathcal{A}^{\otimes r},\mathcal{A}^{\otimes s},\mathcal{A}^{\otimes t}}(\iota_r\otimes \id\otimes \iota_t)(\phi_1\otimes \dots\otimes \phi_r\otimes \phi_{r+1}m_s\otimes \phi_{r+2}\otimes \dots\otimes \phi_{r+1+t})\\
=&(-1)^{(2-s)\sum_{j=1}^{r+1+t}|\phi_j|}(\phi_1\otimes \dots\otimes \phi_{r+1+t})(\id^{\otimes r}\otimes m_s\otimes \id^{\otimes t}).
\end{align*}

Thus, with $\mu_n=(-1)^n \iota_n^{-1}m_n^\#$, it follows that 
\begin{align*}
&\sum_{r+s+t=n}(-1)^{rs+t}(\id^{\otimes r}\otimes \mu_s\otimes \id^{\otimes t})\mu_{r+1+t}\\
=&\sum_{r+s+t=n}(-1)^{rs+t}(-1)^{s+r+1+t} (\id^{\otimes r}\otimes \iota_s^{-1} m_s^\#\otimes \id^{\otimes t})\iota_{r+1+t}^{-1}m_{r+1+t}^\#\\
=&\sum_{r+s+t=n}(-1)^{rs+s+r+1}(\id^{\otimes r}\otimes \iota_s^{-1}\otimes \id^{\otimes t})(\id^{\otimes r}\otimes m_s^\#\otimes \id^{\otimes t})\iota_{r+1+t}^{-1}m_{r+1+t}^\#\\
\stackrel{(*)}{=}&\sum_{r+s+t=n}(-1)^{rs+s+r+1} (\id^{\otimes r}\otimes \iota_s^{-1}\otimes \id^{\otimes t})(\iota_r\otimes \id\otimes \iota_t)^{-1}\iota_{A^{\otimes r},A^{\otimes s},A^{\otimes t}}^{-1}(\id^{\otimes r}\otimes m_s\otimes \id^{\otimes t})^\#m_{r+1+t}^\#\\
=&\iota_n^{-1}\left(\sum_{r+s+t=n} (-1)^{rs+s+r+1}(-1)^{s(r+1+t)}m_{r+1+t}(\id^{\otimes r}\otimes m_s\otimes \id^{\otimes t})\right)^\#\\
=&-\iota_n^{-1}\left(\sum_{r+s+t=n}(-1)^{r+st}m_{r+1+t}(\id^{\otimes r}\otimes m_s\otimes \id^{\otimes t})\right)^\#\\
=&0,
\end{align*}
where for the equality marked with $(*)$ we use the commutativity of \eqref{eqn213}. Finally the assumption that only finitely many $m_n\colon \mathcal{A}^{\otimes n}\to \mathcal{A}$ are non-zero guarantees that the map $(\mu_n)_{n\in \mathbb{N}}\colon \mathcal{A}^{\#}\to \prod_{n\in \mathbb{N}}(\mathcal{A}^{\#})^{\otimes n}$ factors through the direct sum. 

Next, consider a morphism of $A_\infty$-categories $g_n\colon \mathcal{A}^{\otimes n}\to \mathcal{B}$ and recall that $(g_n)^\#\colon \mathcal{B}^\#\to (\mathcal{A}^{\otimes n})^\#$ is given by $(g_n)^\#(\phi)=(-1)^{(n-1)|\phi|}\phi\circ g_n$. We first claim that the diagram
\begin{equation}\label{dualmor}
\begin{tikzcd}
(\mathcal{B}^\#)^{\otimes q}\arrow{rrr}{\iota_q}\arrow{d}{(g_{j_1})^\#\otimes \cdots \otimes (g_{j_q})^\#} &&& (\mathcal{B}^{\otimes q})^\# \arrow{d}{(g_{j_1}\otimes \cdots \otimes g_{j_q})^\#}\\
(\mathcal{A}^{\otimes j_1})^\#\otimes \cdots \otimes (\mathcal{A}^{\otimes j_q})^\# \arrow{rrr}{\iota_{\mathcal{A}^{\otimes j_1},\dots,\mathcal{A}^{\otimes j_q}}} &&&(\mathcal{A}^{\otimes j_1}\otimes \cdots \otimes \mathcal{A}^{\otimes j_q})^\# 
\end{tikzcd}
\end{equation}
commutes. Indeed, evaluating on $\phi_1\otimes \cdots \otimes \phi_q\in (\mathcal{A}^\#)^{\otimes q}$  we obtain, on the one hand,

\begin{equation*}
\begin{split}
(g_{j_1}\otimes&\cdots\otimes g_{j_q})^\# \iota_q(\phi_1\otimes \dots\otimes \phi_q)\\
&=(g_{j_1}\otimes \dots\otimes g_{j_q})^\#(\phi_1\otimes \dots\otimes \phi_q)\\
&=(-1)^{\sum_{i=1}^q(1-j_i)\sum_{r=1}^q|\phi_r|}(\phi_1\otimes \dots\otimes \phi_q)\circ (g_{j_1}\otimes \dots\otimes g_{j_q})\\
&=(-1)^{\sum_{i=1}^q(1-j_i)\sum_{r=1}^q |\phi_r|+\sum_{i<r}(1-j_i)|\phi_r|} (\phi_1g_{j_1}\otimes \dots\otimes \phi_qg_{j_q})\\
&=(-1)^{\sum_{i\geq r}(1-j_i)|\phi_r|}(\phi_1g_{j_1}\otimes \dots\otimes \phi_qg_{j_q}).
\end{split}
\end{equation*}

On the other hand, we compute
\begin{equation*}
\begin{split}
\iota_{\mathcal{A}^{\otimes j_1},\dots,\mathcal{A}^{\otimes j_q}}&(((g_{j_1})^\#\otimes \cdots \otimes (g_{j_q})^\#)(\phi_1\otimes \cdots \otimes \phi_q))\\
&=(-1)^{\sum_{i>r}(1-j_i)|\phi_r|}\iota_{\mathcal{A}^{\otimes j_1},\dots,\mathcal{A}^{\otimes j_q}}(g_{j_1}^\#(\phi_1)\otimes \cdots \otimes g_{j_q}^\#(\phi_q))\\
&=(-1)^{\sum_{i\geq r}(1-j_i)|\phi_r|}\phi_1g_{j_1}\otimes \cdots\otimes \phi_qg_{j_q}.
\end{split}
\end{equation*}

Now define $f_n:=\iota_n^{-1}g_n^\#$. We verify that this is a morphism of $A_\infty$-coalgebras by computing, on the one hand,
\begin{align*}
&\sum_{r+s+t=n} (-1)^{rs+t} (\id^{\otimes r}\otimes \mu_s\otimes \id^{\otimes t})f_{1+r+t} \\
=&\sum_{r+s+t=n} (-1)^{rs+t+s}(\id^{\otimes r}\otimes \iota^{-1}_sm_s^\#\otimes \id^{\otimes t})\iota_{1+r+t}^{-1}g_{1+r+t}^\#\\
=&\sum_{r+s+t=n}(-1)^{rs+t+s} (\id^{\otimes r}\otimes \iota^{-1}_s\otimes \id^{\otimes t})(\id^{\otimes r}\otimes m_s^\#\otimes \id^{\otimes t})\iota_{1+r+t}^{-1}g_{1+r+t}^\#\\
\stackrel{(\ast)}{=}&\sum_{r+s+t=n} (-1)^{rs+t+s}(\id^{\otimes r}\otimes \iota^{-1}_s\otimes \id^{\otimes t})(\iota_r\otimes \id\otimes \iota_t)^{-1}\iota_{\mathcal{A}^{\otimes r},\mathcal{A}^{\otimes s},\mathcal{A}^{\otimes t}}^{-1}(\id^{\otimes r}\otimes m_s\otimes \id^{\otimes t})^\# g_{1+r+t}^\#\\
=&\sum_{r+s+t=n}(-1)^{rs+t+s}\iota_n^{-1}(\id^{\otimes r}\otimes m_s\otimes \id^{\otimes t})^\#g_{1+r+t}^\#\\
=&\iota_{n}^{-1}\sum_{r+s+t=n} (-1)^{n+r+st}g_{1+r+t}\circ (\id^{\otimes r}\otimes m_s\otimes \id^{\otimes t})^\#\\
=&(-1)^n \iota_n^{-1}\left(\sum_{r+s+t=n} (-1)^{r+st} g_{1+r+t} (\id^{\otimes r}\otimes m_s\otimes \id^{\otimes t}) \right)^\#,
\end{align*}
where for $(\ast)$ we use the commutativity of \eqref{eqn213}; and, on the other hand,
\begin{align*}
&\sum (-1)^{\sum (i-1)(j_i+1)} (f_{j_1}\otimes \dots \otimes f_{j_q})\mu_q\\
=&\sum (-1)^{q+\sum (i-1)(j_i+1)} (\iota_{j_1}^{-1} g_{j_1}^\#\otimes \dots\otimes \iota_{j_q}^{-1} g_{j_q}^\#)\iota_q^{-1}m_q^\#\\
=&\sum (-1)^{q+\sum (i-1)(j_i+1)} (\iota_{j_1}^{-1}\otimes \dots\otimes \iota_{j_q}^{-1})(g_{j_1}^\#\otimes \dots\otimes g_{j_q}^\#)\iota_q^{-1}m_q^\#\\
\stackrel{(\dagger)}{=}&\sum (-1)^{q+\sum (i-1)(j_i+1)} (\iota_{j_1}^{-1}\otimes \dots\otimes \iota_{j_q}^{-1})\iota^{-1}_{\mathcal{A}^{\otimes j_1},\dots,\mathcal{A}^{\otimes j_q}} (g_{j_1}\otimes \dots\otimes g_{j_q})^\# m_q^\#\\
=&\iota_n^{-1}\sum (-1)^{q+\sum (i-1)(j_i+1)+\sum_{i=1}^q q(j_i+1)} (m_q\circ (g_{j_1}\otimes \dots\otimes g_{j_q}))^\#\\
=&(-1)^n\iota_n^{-1}\left( \sum (-1)^{\sum (q-i)(j_i+1)} (m_q\circ (g_{j_1}\otimes \dots \otimes g_{j_q}))^\#\right)
\end{align*}
where the equality $(\dagger)$ follows from commutativity of \eqref{dualmor}.

Finally, we verify that this respects composition, i.e. that for two composable morphisms $f,g$ of $A_\infty$-categories $(f\circ g)^\# = g^\#\circ f^\#$. To see this, we compute, on the one hand,
\begin{equation*}\begin{split}
((f\circ g)^\#)_n & = \iota_n^{-1}((f\circ g)_n)^\#\\
&= \iota_n^{-1}\sum (-1)^{\sum_{i=1}^r(r-i)(j_i-1)}(f_r(g_{j_1}\otimes \cdots \otimes g_{j_r}))^\#\\
&= \iota_n^{-1}\sum (-1)^{\sum_{i=1}^r(r-i)(j_i-1)+(1-r)(r-n)}(g_{j_1}\otimes \cdots \otimes g_{j_r})^\#(f_r)^\#
\end{split}\end{equation*}
and on the other hand
\begin{equation*}\begin{split}
(g^\#\circ f^\#)_n &= \sum(-1)^{\sum_{i=1}^r(i-1)(j_i+1)}((g^\#)_{j_1}\otimes \cdots \otimes(g^\#)_{j_r})(f_r^\#)\\
&=  \sum(-1)^{\sum_{i=1}^r(i-1)(j_i+1)}(\iota_{j_1}^{-1}(g_{j_1})^\#\otimes \cdots \otimes \iota_{j_r}^{-1}(g_{j_r})^\#)\iota_r^{-1}(f_r)^\#\\
&=  \sum(-1)^{\sum_{i=1}^r(i-1)(j_i+1)}(\iota_{j_1}^{-1}\otimes \cdots \otimes \iota_{j_r}^{-1})((g_{j_1})^\#\otimes \cdots \otimes(g_{j_r})^\#)\iota_r^{-1}(f_r)^\#\\
&\stackrel{(\dagger)}{=}\sum (-1)^{\sum_{i=1}^r(i-1)(j_i+1)} (\iota_{j_1}^{-1}\otimes \dots\otimes \iota_{j_r}^{-1})\iota^{-1}_{\mathcal{A}^{\otimes j_1},\dots,\mathcal{A}^{\otimes j_r}} (g_{j_1}\otimes \dots\otimes g_{j_r})^\#(f_r)^\#\\
&= \sum (-1)^{\sum_{i=1}^r(i-1)(j_i+1)} \iota_n^{-1} (g_{j_1}\otimes \dots\otimes g_{j_r})^\#(f_r)^\#
\end{split}\end{equation*}
where again the equation labelled by $(\dagger)$ follows from \eqref{dualmor}.
Comparing signs, on the one hand, the sign is
\[(-1)^{\sum_{i=1}^r(r-i)(j_i-1)+(1-r)(r-n)}=(-1)^{rn-(\sum_{i=1}^rij_i) -r^2+   (\sum_{i=1}^ri )+   r-r^2-n+rn} = (-1)^{(\sum_{i=1}^rij_i )+   (\sum_{i=1}^ri) +   r-n}\]
while, on the other hand, the sign is
\[ (-1)^{\sum_{i=1}^r(i-1)(j_i+1)} = (-1)^{(\sum_{i=1}^r ij_i )-n +(\sum_{i=1}^r i)   -r}.\] These coincide, which completes the proof.
\end{proof}

\begin{ex}
The additional assumption in the preceding lemma, part \eqref{lem2.1.2} is necessary. For example, consider the two-dimensional $A_\infty$-algebra with basis $x,y$ and multiplications defined by $m_n(y,\dots,y)=x$ and $m_n(\dots,x,\dots)=0$. Then the dual map does not factor through the direct sum. These situations do not occur in our setup  since the $A_\infty$-categories we study are directed in a certain sense, implying that there are only finitely many nonzero $m_n$. 
\end{ex}

\begin{lem}\label{dualtocoalgebraobject}
Let $A$ be an associative algebra. Let $\mathcal{C}$ be an $A_\infty$-coalgebra object in $\modu A\otimes A^{\op}$.   Then $\Hom_{A\otimes A^{\op}}(\mathcal{C},\Hom_{\Bbbk}(\mathbb{L},\mathbb{L}))$ is an $A_\infty$-category with multiplications given by 
\[m_n=(-1)^n\Hom_{A\otimes A^{\op}}(\mu_n,\Hom_{\Bbbk}(\mathbb{L},\mathbb{L}))\iota_n'\]
where, similarly to before, $(\iota_n'(\phi_1\otimes \cdots \phi_n))(c_1\otimes \cdots \otimes c_n) = (-1)^{\sum_{j<i} |\phi_i|\cdot |c_j|} \phi_1(c_1)\circ \dots\circ \phi_n(c_n)$.
\end{lem}

\begin{proof}
The proof is similar to the proof of the preceding lemma.
\end{proof}

We remark that this is indeed an $A_\infty$-category, since $\Hom_{\Bbbk}(\mathbb{L},\mathbb{L})$ carries an $\mathbb{L}$-$\mathbb{L}$-bimodule structure given by $e_\mathtt{j}\Hom_{\Bbbk}(\mathbb{L},\mathbb{L})e_\mathtt{i} = \Hom_{\Bbbk}(\mathbb{L}e_\mathtt{i},\mathbb{L}e_\mathtt{j})$, which translates composition $\circ$ into $\otimes_\mathbb{L}$.

\subsection{Bar and cobar constructions}\label{sec:barcobar}
The next thing we recall is the bar construction associating to an $A_\infty$-algebra object a differential graded coalgebra object. Since all $A_\infty$-algebra objects we will consider are augmented, we only describe it in the augmented case. For this, we give the definition of the shift functor on graded objects in $\mathsf{C}$. For a graded object $M$ in $\mathsf{C}$, the \emphbf{shift} $\shift M$ of $M$ is defined to be the graded object with $(\shift M)_i=M_{i+1}$. If $(M,d)$ is differential graded, then we define $d_{\shift M}=-d_M$. We denote the inverse shift by $\shift^{-1}$. If $f\colon M\to N$ is a homogeneous map of degree $p$, then $\shift f$ is defined to be $(-1)^p f$. Note that, by  the Koszul sign rule, we have $(\shift^{\otimes n})^{-1}=(-1)^{\frac{n(n+1)}{2}} (\shift^{-1})^{\otimes n}$. 

Let $\mathcal{A}$ be an augmented $A_\infty$-algebra object in $\mathsf{C}$. Let $\mathfrak{m}$ be the kernel of the augmentation. The \emphbf{bar construction} $\mathbf{B}\mathcal{A}$ of $\mathcal{A}$ is the differential graded coalgebra object with underlying graded object 
\[T^c(\shift\mathfrak{m})=e\oplus \shift\mathfrak{m}\oplus (\shift\mathfrak{m})^{\otimes 2}\oplus \dots.\]
In order to distinguish between the different tensor products, we sometimes write $[\shift a_1|\shift a_2|\dots|\shift a_j]$ for an element of $(\shift \mathfrak{m})^{\otimes j}$. The coproduct on $T^c(\shift \mathfrak{m})$ is given by
\[\mu([\shift a_1|\shift a_2|\dots|\shift a_j])=\sum_{i=0}^j[\shift a_1|\dots|\shift a_i]\otimes [\shift a_{i+1}|\dots |\shift a_j].\]
Furthermore, since $T^c(\shift \mathfrak{m})$ is a cofree object in the category of graded coalgebra objects, a differential on $T^c(\shift \mathfrak{m})$ is determined by maps $b_n\colon (\shift \mathfrak{m})^{\otimes n}\to (\shift \mathfrak{m})$ of degree $1$. We define these maps via $b_n=\shift m_n(\shift ^{\otimes n})^{-1}$. 

\begin{lem}
The following are equivalent:
\begin{enumerate}
\item The $m_n$ satisfy the $A_\infty$-relations.
\item The $b_n$ satisfy the equation
\[\sum_{r+s+t=n} b_{r+1+t}(\id^{\otimes r}\otimes b_s\otimes \id^{\otimes t})=0\]
\item The $b_n$ induce a differential $b$ on $T^c(\shift \mathfrak{m})$, i.e. $b^2=0$. 
\end{enumerate}
\end{lem}

\begin{proof}
To prove (1)$\Rightarrow$(2) we compute
\begin{align*}
m_{r+1+t}(\id^{\otimes r}\otimes m_s\otimes \id^{\otimes t})&=\shift ^{-1}b_{r+1+t}\shift^{\otimes (r+1+t)}(\id^{\otimes r}\otimes \shift^{-1}b_s \shift^{\otimes s}\otimes \id^{\otimes t})\\
&=(-1)^{t(2-s)}\shift^{-1}b_{r+1+t}(\shift^{\otimes r}\otimes b_s\shift^{\otimes s}\otimes \shift^{\otimes t})\\
&=(-1)^{t(2-s)+r}\shift^{-1}b_{r+t+1}(\id^{\otimes r}\otimes b_s\otimes \id^{\otimes t})\shift^{\otimes n}\\
&=(-1)^{r+st}\shift^{-1}b_{r+t+1}(\id^{\otimes r}\otimes b_s\otimes \id^{\otimes t})\shift^{\otimes n}
\end{align*}
This shows that (1) is equivalent to (2). As noted above, the equivalence of (2) and (3) follows from cofreeness of $T^c(\shift \mathfrak{m})$. 
\end{proof}

The bar construction can be extended to a fully faithful functor. Let $\mathcal{A}$ and $\mathcal{B}$ be $A_\infty$-algebra objects with augmentation ideals $\mathfrak{m}$ and $\mathfrak{n}$, respectively. For a morphism $f\colon \mathcal{A}\to \mathcal{B}$ of $A_\infty$-algebra objects we obtain a morphism $\mathbf{B}f$ determined by components $\shift f(\shift^{\otimes j})^{-1}\colon (\shift\mathfrak{m})^{\otimes j}\to \shift\mathfrak{n}$, using cofreeness. 

Occasionally we will also use  the non-unital version of this construction by dropping the direct summand $e$ and denote it by 
\[\overline{T}^c(\shift \mathfrak{m})=\shift \mathfrak{m}\oplus (\shift \mathfrak{m})^{\otimes 2}\oplus \dots.\]

Dually, one defines the {\bf cobar construction}. For a coaugmented $A_\infty$-coalgebra object $\mathcal{C}$ in $\mathsf{C}$ with kernel of the counit $\mathfrak{m}$ let $\mathbf{\Omega} \mathcal{C}$ be the differential graded algebra object with underlying graded object
\[T(\shift^{-1}\mathfrak{m})=e\oplus \shift^{-1}\mathfrak{m}\oplus (\shift^{-1}\mathfrak{m})^{\otimes 2}\oplus \dots.\]
Again we write $[\shift^{-1}a_1|\shift^{-1}a_2|\dots|\shift^{-1}a_j]$ for an element of $(\shift^{-1}\mathfrak{m})^{\otimes j}$ in order to avoid confusion with the different levels of tensor products. The multiplication on $T(\shift^{-1}\mathfrak{m})$ is defined by
\[m([\shift^{-1}a_1|\dots |\shift^{-1}a_i]\otimes [\shift^{-1}a_{i+1}|\dots |\shift^{-1}a_{i+j}])=[\shift^{-1}a_1|\shift^{-1}a_2|\dots|\shift^{-1}a_{i+j}].\]
The differential is, by freeness of $T(\shift^{-1}\mathfrak{m})$, determined uniquely by the collection of $d_n\colon \shift^{-1}\mathfrak{m}\to (\shift^{-1}\mathfrak{n})^{\otimes n}$ given by
\[d_n=(\shift^{\otimes n})^{-1} \mu_n \shift.\]

\begin{lem}
The following are equivalent:
\begin{enumerate}[(1)]
\item The $\mu_n$ satisfy the $A_\infty$-coalgebra relations.
\item The $d_n$ satisfy the relation
\[\sum_{r+s+t=n} (\id^{\otimes r}\otimes d_s\otimes \id^{\otimes t})d_{r+1+t}=0.\]
\item The unique extension of the $d_n$ to $T(\shift^{-1}\mathfrak{m})$ is a differential, i.e. $d^2=0$. 
\end{enumerate}
\end{lem}

\begin{proof}
To prove that (1)$\Leftrightarrow$ (2) we calculate that
\begin{align*}
(\id^{\otimes r}\otimes \mu_s\otimes \id^{\otimes t})\mu_{r+1+t}&=(\id^{\otimes r}\otimes (\shift^{\otimes s}d_s\shift^{-1}) \otimes \id^{\otimes t})\shift^{\otimes (1+r+t)}d_{1+r+t}\shift^{-1}\\
&=(-1)^{r(s-2)}(\shift^{\otimes r}\otimes \shift^{\otimes s}d_s\otimes \shift^{\otimes t})d_{1+r+t}\shift^{-1}\\
&=(-1)^{rs+t}\shift^{\otimes n}(\id^{\otimes r}\otimes d_s\otimes \id^{\otimes t})d_{1+r+t}\shift^{-1}
\end{align*}
As noted above, the equivalence of (2) and (3) follows from freeness of $T(\shift^{-1}\mathfrak{m})$. 
\end{proof}

Also, the cobar construction can be extended to a fully-faithful functor. For a morphism of $A_\infty$-coalgebra objects $f\colon \mathcal{C}\to \mathcal{D}$ with components $f_i\colon \mathcal{C}\to \mathcal{D}^{\otimes i}$, $\mathbf{\Omega} f$ is the unique extension of the map with components 
$(\shift^{\otimes i})^{-1}f_i\shift$, using freeness.

The following lemma gives a compatibility result for the bar/cobar constructions and duality. This is well-known and is, for example stated without proof in \cite[Lemma 10]{EL17}. The condition that $\mathbf{B}\mathcal{E}$ is finite dimensional is quite strong but is satisfied in our setting where $\mathcal{E}$ is the $\Ext$-algebra of the standard modules of a quasi-hereditary algebra since these form an exceptional collection. For a more general statement one has to pass to the completion as done in \cite[Lemma 11]{EL17}.

\begin{lem}\label{barcobariso}
\begin{enumerate}[(i)]          
\item\label{duality1} 
Let $\mathcal{E}$ be an augmented $A_\infty$-category such that $\mathbf{B}\mathcal{E}$ is finite dimensional. Then there is a functorial isomorphism of differential graded categories
\[\mathbf{\Omega}(\mathcal{E}^{\#})\cong (\mathbf{B}\mathcal{E})^{\#}.\]
\item\label{duality2}  Dually, if $\mathcal{C}$ is a coaugmented $A_\infty$-cocategory such that $\mathbf{\Omega} \mathcal{C}$ is finite dimensional, then there is a functorial isomorphism of differential graded cocategories
\[\mathbf{B}(\mathcal{C}^{\#})\cong (\mathbf{\Omega}\mathcal{C})^{\#}.\]
\end{enumerate}
\end{lem}

\begin{proof}

Define $\chi\colon \shift^{-1}\mathcal{E}^{\#} \to(\shift\mathcal{E})^{\#}$ to be $\chi:=-(\shift^{\#})^{-1}\shift=(\shift^{-1})^\#\shift$. This implies that the diagram
\[
\begin{tikzcd}
\shift^{-1}\mathcal{E}^{\#} \arrow{dr}[swap]{-\shift}\arrow{r}{\chi}&(\shift\mathcal{E})^{\#}\arrow{d}{\shift^{\#} }\\
&\mathcal{E}^{\#}
\end{tikzcd}\]
commutes. 
As graded vector spaces, $ (\mathbf{B}\mathcal{E})^{\#}=\bigoplus_{k\geq 0}((\shift\mathfrak{m})^{\otimes k})^\#$ and $\mathbf{\Omega}(\mathcal{E}^{\#}) = \bigoplus_{k\geq 0} (\shift^{-1}(\mathfrak{m}^\#))^{\otimes k}$.
We define an isomorphism of graded vector spaces component-wise by the isomorphism
\[\iota_k\circ\chi^{\otimes k}\colon(\shift^{-1}(\mathfrak{m}^\#))^{\otimes k}\to 
((\shift\mathfrak{m})^\#)^{\otimes k}\to((\shift\mathfrak{m})^{\otimes k})^\#\]
and check that this is multiplicative and preserves the differential.
Using that $\chi$ is of degree $0$ we compute 
\begin{equation}\label{barcobar3}
\begin{split}
\iota_k\circ\chi^{\otimes k}\left(m_{\mathbf{\Omega}(\mathcal{E}^{\#})}([\shift^{-1}a_1^\#|\dots|\shift^{-1}a_t^\#]\otimes[\shift^{-1}a_{t+1}^\#|\dots|\shift^{-1}a_k^\#])\right)&=
\iota_k\circ\chi^{\otimes k}([\shift^{-1}a_1^\#|\dots|\shift^{-1}a_k^\#])\\
&= (-1)^k
\iota_k([(\shift a_1)^\#|\dots|(\shift a_k)^\#])
\end{split}
\end{equation}
and
\begin{equation}\label{barcobar4}
\begin{split}
m_{\mathbf{B}(\mathcal{E})^{\#}}&\left(((\iota_t\circ\chi^{\otimes t})\otimes (\iota_{k-t}\circ\chi^{\otimes k-t}))([\shift^{-1}a_1^\#|\dots|\shift^{-1}a_t^\#]\otimes[\shift^{-1}a_{t+1}^\#|\dots|\shift^{-1}a_k^\#])\right)\\
&= (-1)^k m_{\mathbf{B}(\mathcal{E})^{\#}}(\iota_t\otimes \iota_{k-t})([(\shift a_1)^\#|\dots|(\shift a_t)^\#]\otimes [(\shift a_{t+1})^\#|\dots|(\shift a_k)^\#])\\
\end{split}
\end{equation}

so we need to show that the evaluation of, on the one hand, $\iota_k([(\shift a_1)^\#|\dots|(\shift a_k)^\#])$ and, on the other hand, $m_{B(\mathcal{E})^{\#}}(\iota_t\otimes \iota_{k-t})([(\shift a_1)^\#|\dots|(\shift a_t)^\#]\otimes [(\shift a_{t+1})^\#|\dots|(\shift a_k)^\#])$ on some $[x_1|\dots |x_k]$ coincide.

Indeed,

\[\iota_k([(\shift a_1)^\#|\dots|(\shift a_k)^\#])([x_1|\dots |x_k]) = (-1)^{\sum_{i=1}^k\sum_{j=1}^{i-1}|(\shift a_i)^\#||x_j|}((\shift a_1)^\#(x_1))\cdots((\shift a_k)^\#(x_k))\]
while
\begin{equation}
\begin{split}
m&_{B(\mathcal{E})^{\#}}(\iota_t\otimes \iota_{k-t})([(\shift a_1)^\#|\dots|(sa_t)^\#]\otimes [(\shift a_{t+1})^\#|\dots|(\shift a_k)^\#])([x_1|\dots |x_k]) \\
& =(\iota_t[(\shift a_1)^\#|\dots|(\shift a_t)^\#]\otimes \iota_{k-t}[(\shift a_{t+1})^\#|\dots|(\shift a_k)^\#])([x_1|\dots |x_t]\otimes [x_{t+1}|\dots |x_k] )\\
&=(-1)^{(\sum_{i=t+1}^k|(\shift a_i)^\#|)(\sum_{j=1}^t|x_j|)}(\iota_t[(\shift a_1)^\#|\dots|(\shift a_t)^\#])([x_1|\dots |x_t])(\iota_{k-t}[(\shift a_{t+1})^\#|\dots|(\shift a_k)^\#])( [x_{t+1}|\dots |x_k])\\
&= (-1)^{(\sum_{i=t+1}^k|(\shift a_i)^\#|)(\sum_{j=1}^t|x_j|) +\sum_{i=1}^t\sum_{j=1}^{i-1}|(\shift a_i)^\#||x_j| +\sum_{i=t+1}^k\sum_{j=t+1}^{i-1}|(\shift a_i)^\#||x_j|}((\shift a_1)^\#(x_1))\cdots((\shift a_k)^\#(x_k))\\
&=(-1)^{\sum_{i=1}^k\sum_{j=1}^{i-1}|(\shift a_i)^\#||x_j|}((\shift a_1)^\#(x_1))\cdots((\shift a_k)^\#(x_k)).
\end{split}
\end{equation}

Thus we have an isomorphism of graded algebras. In order to confirm that the differential is preserved, taking into account  $\iota_1=\id$ and the sign of the differential of the shifted complex, it suffices, by freeness, to check that 
the diagram
\[
\begin{tikzcd}
\shift ^{-1}(\mathcal{E}^{\#}))\arrow{rr}{\chi}\arrow{d}{d_n}&&(\shift \mathcal{E})^{\#}\arrow{d}{-b_n^{\#}}\\
(\shift ^{-1}\mathcal{E}^{\#})^{\otimes n} \arrow{r}{\chi^{\otimes n}}&((\shift\mathcal{E})^{\#})^{\otimes n} \arrow{r}{\iota_n} & ((\shift\mathcal{E})^{\otimes n})^{\#}\\
\end{tikzcd}
\]
commutes. In other words, we claim that
\begin{equation}\label{bndneq}
\iota_n \chi^{\otimes n} d_n= -b_n^{\#} \chi. 
\end{equation}
Inserting the definition of $d_n$, the left hand side becomes
\begin{equation}\label{bocscoreq}
\iota_n \chi^{\otimes n} d_n=\iota_n \chi^{\otimes n} (\shift^{\otimes n})^{-1} \mu_n\shift=(-1)^n\iota_n \chi^{\otimes n} (\shift^{\otimes n})^{-1} \iota_n^{-1}m_n^{\#}\shift,
\end{equation}
while, using the definition of $b_n$, the right hand side yields
\begin{equation*}
-b_n^{\#} \chi = -(\shift m_n(\shift^{\otimes n})^{-1})^{\#}\chi=-(-1)^n((\shift^{\otimes n})^{-1})^{\#}m_n^{\#}\shift^{\#}\chi=(-1)^n((\shift^{\otimes n})^{-1})^{\#}m_n^{\#}\shift.
\end{equation*}
Using $(\shift^{\otimes n})^{-1} = (-1)^{\frac{n(n-1)}{2}}(\shift^{-1})^{\otimes n}$ on both sides and cancelling signs, in order to prove \eqref{bndneq} it suffices to prove
\begin{equation}\label{iota_equation}
\iota_n \chi^{\otimes n} (\shift^{-1})^{\otimes n} \iota_n^{-1}=((\shift^{-1})^{\otimes n})^{\#}
\end{equation}
or equivalently, using 
\[ \chi^{\otimes n} (\shift^{-1})^{\otimes n}  =  (\chi \shift^{-1})^{\otimes n}  =( (\shift^{-1})^{\#})^{\otimes n},\]
(where the left equality follows since $\chi$ is of degree $0$) it suffices to  show
\[\iota_n ( (\shift^{-1})^{\#})^{\otimes n}=((\shift^{-1})^{\otimes n})^{\#}\iota_n. \]

We evaluate both sides and obtain
\begin{equation*}\begin{split}
\iota_n ( (\shift^{-1})^{\#})^{\otimes n}&(\phi_1\otimes \cdots \otimes \phi_n)(\shift m_1\otimes \cdots \otimes \shift m_n)\\
&=(-1)^{\sum_{i=1}^n(n-i+1) |\phi_i|}\iota_n (\phi_1\shift^{-1}\otimes \cdots \otimes \phi_n\shift^{-1})(\shift m_1\otimes \cdots \otimes \shift m_n)\\
&=(-1)^{\sum_{i=1}^n(n-i+1) |\phi_i| +\sum_{j>k}(|\phi_j|+1)(|m_k|+1)} \phi_1(m_1)\cdots \phi_n(m_n)
\end{split}\end{equation*}
and, on the other hand,
\begin{equation*}\begin{split}
((\shift^{-1})^{\otimes n})^{\#}\iota_n&(\phi_1\otimes \cdots \otimes \phi_n)(\shift m_1\otimes \cdots \otimes \shift m_n)\\
&=(-1)^{\sum_{i=1}^nn |\phi_i|}\iota_n(\phi_1\otimes \cdots \otimes \phi_n)(\shift^{-1})^{\otimes n}(\shift m_1\otimes \cdots \otimes \shift m_n)\\
&=(-1)^{\sum_{i=1}^nn |\phi_i| +\sum_{k=1}^n(n-k)(|m_k|+1)}\iota_n(\phi_1\otimes \cdots \otimes \phi_n)(m_1\otimes \cdots \otimes m_n)\\
&=(-1)^{\sum_{i=1}^nn |\phi_i| +\sum_{k=1}^n(n-k)(|m_k|+1)+\sum_{j>l}|\phi_j||m_l|}\phi_1(m_1)\cdots \phi_n(m_n).
\end{split}\end{equation*}
We now check that the signs agree. Indeed, the sign on the left hand side is
\begin{equation*}\begin{split}
(-1)^{\sum_{i=1}^n(n-i+1) |\phi_i| +\sum_{j>k}(|\phi_j|+1)(|m_k|+1)} & = (-1)^{\sum_{i=1}^n(n-i+1) |\phi_i| +\sum_{j>k}|\phi_j||m_k| +\sum_{j>k}|\phi_j|+\sum_{j>k}|m_k|+\sum_{j>k}1} \\
& = (-1)^{\sum_{i=1}^n(n-i+1) |\phi_i| +\sum_{j>k}|\phi_j||m_k| +\sum_{j=1}^n(j-1)|\phi_j|+\sum_{k=1}^n(n-k)|m_k|+\frac{n(n-1)}{2}}\\
&= (-1)^{\sum_{i=1}^n n|\phi_i| +\sum_{j>k}|\phi_j||m_k| +\sum_{k=1}^n(n-k)|m_k|+\frac{n(n-1)}{2}}
\end{split}\end{equation*}
while the sign on the right hand side is
\begin{equation*}\begin{split}
(-1)^{\sum_{i=1}^nn |\phi_i| +\sum_{k=1}^n(n-k)(|m_k|+1)+\sum_{j>l}|\phi_j||m_l|}& = (-1)^{\sum_{i=1}^nn |\phi_i| +\sum_{k=1}^n(n-k)|m_k|+\sum_{k=1}^n(n-k)+\sum_{j>l}|\phi_j||m_l|} \\&=  (-1)^{\sum_{i=1}^nn |\phi_i| +\sum_{k=1}^n(n-k)|m_k|+\frac{n(n-1)}{2}+\sum_{j>l}|\phi_j||m_l|}, \end{split}\end{equation*}
and hence they agree, implying that we have an isomorphism of differential graded algebras.

We now check that this isomorphism is functorial. To this end let $f\colon \mathcal{E} \to \mathcal{E}'$ be a morphism of augmented $A_\infty$-algebras with augmentation ideals $\mathfrak{m}$ and $\mathfrak{n}$ respectively. We need to check that our isomorphism translates $\mathbf{\Omega} (f^{\#})$ into $(\mathbf{B}f)^{\#}$. The former has components
\[(\mathbf{\Omega} (f^{\#}) )_n = (\shift^{\otimes n})^{-1}\iota_n^{-1}f_n^{\#}\shift  \colon \shift^{-1}(\mathfrak{n}^{\#}) \to (\shift^{-1}\mathfrak{m}^{\#})^{\otimes n} \]
while the latter has components
\[((\mathbf{B}f)^{\#})_n =\iota_n^{-1}(\shift f_n(\shift^{\otimes n})^{-1})^{\#}\colon (\shift\mathfrak{n})^{\#} \to  ((\shift\mathfrak{m})^{\#})^{\otimes n}.\]
We claim that the diagram
\[
\begin{tikzcd}
\shift^{-1}(\mathfrak{n}^{\#})\arrow{rr}{\chi}\arrow{d}{(\shift^{\otimes n})^{-1}\iota_n^{-1}f_n^{\#}\shift}&&(\shift\mathfrak{n})^{\#}\arrow{d}{(\shift f_n(\shift^{\otimes n})^{-1})^{\#}}\\
(\shift^{-1}\mathfrak{m}^{\#})^{\otimes n} \arrow{r}{\chi^{\otimes n}}&((\shift\mathfrak{m})^{\#})^{\otimes n} \arrow{r}{\iota_n}& ((\shift\mathfrak{m})^{\otimes n})^{\#}\\
\end{tikzcd}
\]
commutes. That is, we claim that
\[ \iota_n \chi^{\otimes n}(\shift^{\otimes n})^{-1}\iota_n^{-1}f_n^{\#}\shift = (\shift f_n(\shift^{\otimes n})^{-1})^{\#} \chi.\]
The right hand side yields
\[(\shift f_n(\shift^{\otimes n})^{-1})^{\#} \chi = (-1)^{n^2+n-1} ((\shift^{\otimes n})^{-1})^{\#} f_n^{\#}\shift^{\#}\chi = ((\shift^{\otimes n})^{-1})^{\#} f_n^{\#}\shift\]
so again, the problem reduces to 
\[ \iota_n \chi^{\otimes n}(\shift^{\otimes n})^{-1}\iota_n^{-1} =  ((\shift^{\otimes n})^{-1})^{\#},\]
which is \eqref{iota_equation}.

This completes the proof of part \eqref{duality1}. The proof of part \eqref{duality2} is similar and will be omitted since it is not needed for the proof of the main theorem.
\end{proof}

\subsection{Kadeishvili's theorem}\label{sec:Kadeishvili1}

An important advantage of $A_\infty$-categories in contrast to differential graded categories is that the structure can be transferred along homotopy equivalences. This was first observed by Kadeishvili in \cite{Kad80, Kad82} and later different techniques were used to find explicit instead of recursive formulas, see e.g. \cite{Mer99, KS01, Mar06, Kop17}. Here we follow mostly \cite{Kop17}, which provides the most details of all sources known to us,  adapting the sign conventions to our setup.

\begin{thm}\label{kadeishvilitheorem}
Let $\mathcal{A}$ be an $A_\infty$-category with multiplications $m_n$, $n\geq 1$. Let $\mathcal{E}$ be an $\mathbb{L}$-module equipped with a differential $\tilde{m}_1$ such that $\mathcal{E}$ is homotopy equivalent to $\mathcal{A}$ via maps $\mathbf{i}\colon \mathcal{E}\to \mathcal{A}$ and $\mathbf{p}\colon \mathcal{A}\to \mathcal{E}$ compatible with the differential such that $\mathbf{i}\mathbf{p}-\id_{\mathcal{A}}=m_1h+hm_1$ for some $h\colon \mathcal{A}\to \mathcal{A}$ of degree $-1$. Then there exists an $A_\infty$-structure $\tilde{m}_n$, $n\geq 1$ on $\mathcal{E}$ together with an $A_\infty$-quasi-isomorphism $(\mathbf{i}_n)_{n\geq 1}\colon \mathcal{E}\to \mathcal{A}$ such that $\mathbf{i}_1=\mathbf{i}$. 
\end{thm}

The remainder of the section is devoted to a  proof of this theorem and we start by providing the setup.

\begin{defn}\label{lambdakerneldef}
In the setup of the above theorem define the \emphbf{$\lambda$-kernel} as the linear maps $\lambda_n\colon \mathcal{A}^{\otimes n}\to \mathcal{A}$ of degree $2-n$ recursively via $\lambda_2=m_2$ and 
\[\lambda_n=\sum_{\substack{\ell\neq 1\\j_1+\dots+j_\ell=n}}(-1)^{\sum_{i=1}^{\ell-1}(\ell-i)(j_i-1)}m_\ell(h\lambda_{j_1}\otimes \dots\otimes h\lambda_{j_\ell})\]
for $n\geq 3$ where we interpret $h\lambda_1=\id_\mathcal{A}$. 
\end{defn}

\begin{lem}\label{lambdakernellemma}
Let $(\lambda_n)_{n\geq 2}$ be a $\lambda$-kernel. Then 
\[m_1\lambda_n+\sum_{r+t=n-1} (-1)^{n-1} \lambda_n (\id^{\otimes r}\otimes m_1\otimes \id^{\otimes t})+\sum_{\substack{r+s+t=n\\1<s<n}} (-1)^{r+st} \lambda_{r+1+t}(\id^{\otimes r}\otimes \mathbf{i}\mathbf{p}\lambda_s\otimes \id^{\otimes t})=0.\]
\end{lem}

\begin{proof}
In order to eliminate signs we multiply the claim of the lemma from the left with $\shift$ and from the right with $(\shift^{\otimes n})^{-1}$. We thus obtain that the claim of the lemma is equivalent to
\begin{equation}\label{shiftedlambdakernel}
\begin{split}
\underbrace{\shift m_1\shift^{-1}\shift \lambda_n(\shift^{\otimes n})^{-1}}_{A}+\underbrace{\sum_{r+t=n-1} (-1)^{n-1}\shift \lambda_n(\shift^{\otimes n})^{-1}(\shift^{\otimes n})(\id^{\otimes r}\otimes m_1\otimes\id^{\otimes t})(\shift^{\otimes n})^{-1}}_{B}\\
+\underbrace{\sum_{\substack{r+s+t=n\\1<s<n}}(-1)^{r+st}\shift\lambda_{r+1+t}(\shift^{\otimes (r+1+t)})^{-1}(\shift^{\otimes (r+1+t)})(\id^{\otimes r}\otimes i\shift^{-1}\shift p\shift^{-1}\shift\lambda_s\otimes \id^{\otimes t})(\shift^{\otimes n})^{-1}}_{C}=0. 
\end{split}
\end{equation}
Recall that $b_j=\shift m_j(\shift^{\otimes j})^{-1}$ and define similarly $\hat{\lambda}_j:=\shift \lambda_j(\shift^{\otimes j})^{-1}$, $\hat{\mathbf{p}}:=\shift \mathbf{p}\shift^{-1}$, $\hat{\boldsymbol\imath}=\shift \mathbf{i}\shift^{-1}$. 
With these definitions it is immediate that part $A$ of \eqref{shiftedlambdakernel} is equal to $b_1\hat{\lambda}_n$. For part $B$ we compute the expression
\begin{align*}
(-1)^{n-1}\shift^{\otimes n}(\id^{\otimes r}\otimes m_1\otimes \id^{\otimes t})(\shift^{\otimes n})^{-1}
&=(-1)^{n-1+t+\frac{n(n-1)}{2}}(\shift^{\otimes r}\otimes \shift m_1\otimes \shift^{\otimes t})(\shift^{-1})^{\otimes n}\\
&=(-1)^{n-1+t+\frac{n(n-1)}{2}+r+\frac{n(n-1)}{2}}\id^{\otimes r}\otimes \shift m_1\shift^{-1}\otimes \id^{\otimes t}\\
&=\id^{\otimes r}\otimes b_1\otimes \id^{\otimes t}
\end{align*}
Plugging this equality into $B$ of \eqref{shiftedlambdakernel} this part equals
\[\sum_{r+t=n-1}\hat{\lambda}_n(\id^{\otimes r}\otimes b_1\otimes \id^{\otimes t}).\]
Finally for part $C$ of \eqref{shiftedlambdakernel} we compute 
\begin{equation*}
\begin{split}
(-1)^{r+st}&\shift^{\otimes (r+1+t)}(\id^{\otimes r}\otimes \mathbf{i}\shift^{-1}\shift \mathbf{p}\shift^{-1}\shift \lambda_s\otimes \id^{\otimes t})(\shift^{\otimes n})^{-1}\\
&=(-1)^{r+st+(2-s)t+\frac{n(n-1)}{2}}(\shift^{\otimes r}\otimes \hat{\boldsymbol\imath}\hat{\mathbf{p}}\shift \lambda_s\otimes \shift^{\otimes t})(\shift^{-1})^{\otimes n}\\
&=(-1)^{r+\frac{n(n-1)}{2}+rt+r(1-s)+\frac{r(r-1)}{2}+st+\frac{t(t-1)}{2}}(\id^{\otimes r}\otimes \hat{\boldsymbol\imath}\hat{\mathbf{p}}\shift \lambda_s (\shift^{-1})^{\otimes s}\otimes \id^{\otimes t})\\
&=(-1)^{r+\frac{n(n-1)}{2}+rt+r(1-s)+\frac{r(r-1)}{2}+st+\frac{t(t-1)}{2}+\frac{s(s-1)}{2}}(\id^{\otimes r}\otimes \hat{\boldsymbol\imath}\hat{\mathbf{p}}\hat{\lambda}_s\otimes \id^{\otimes t})\\
&=(\id^{\otimes r}\otimes \hat{\boldsymbol\imath}\hat{\mathbf{p}}\hat{\lambda}_s\otimes \id^{\otimes t}),
\end{split}
\end{equation*}
where the last identity follows from the following identity of triangular numbers for $n=r+s+t$:
\[\frac{n(n-1)}{2}=\frac{r(r-1)}{2}+\frac{s(s-1)}{2}+\frac{t(t-1)}{2}+rs+st+rt.\]
Combining $A$, $B$ and $C$ we have reformulated our claim to the sign-free equation
\begin{equation}
\label{inductionhypothesiskernel}
b_1\hat{\lambda}_n+\sum_{r+t=n-1}\hat{\lambda}_n(\id^{\otimes r}\otimes b_1\otimes \id^{\otimes t})+\sum_{\substack{r+s+t=n\\1<s<n}}\hat{\lambda}_{r+1+t}(\id^{\otimes r}\otimes \hat{\boldsymbol\imath}\hat{\mathbf{p}}\hat{\lambda}_s\otimes \id^{\otimes t})=0.
\end{equation}
Our next step is to rewrite the inductive definition 
\[\lambda_n=\sum_{\substack{\ell\neq 1\\j_1+\dots+j_\ell=n}}(-1)^{\sum_{i=1}^{\ell-1}(\ell-i)(j_i-1)}m_\ell(h\lambda_{j_1}\otimes \dots\otimes h\lambda_{j_\ell})\]
of $\lambda_n$ in terms of the $\hat{\lambda}_j$ and $\hat{h}=\shift h\shift^{-1}$. We claim that 
\begin{equation}
\label{signfreedefinitionkernel}
\hat{\lambda}_n=\sum_{\substack{\ell\neq 1\\j_1+\dots+j_\ell=n}}b_\ell(\hat{h}\hat{\lambda}_{j_1}\otimes \dots\otimes \hat{h}\hat{\lambda_{j_\ell}}).
\end{equation}
Indeed,
\begin{align*}
\begin{split}
\hat{\lambda}_n&=\shift\lambda_n(\shift^{\otimes n})^{-1}\\
&=\sum_{\substack{\ell\neq 1\\j_1+\dots+j_\ell=n}}(-1)^{\sum_{i=1}^{\ell-1}(\ell-i)(j_i-1)}\shift m_\ell(\shift^{\otimes \ell})^{-1}\shift^{\otimes \ell}(h\lambda_{j_1}\otimes \dots\otimes h\lambda_{j_\ell})(\shift^{\otimes n})^{-1}\\
&=\sum_{\substack{\ell\neq 1\\j_1+\dots+j_\ell=n}}(-1)^{\sum_{i=1}^{\ell-1}(\ell-i)(j_i-1)+\frac{n(n-1)}{2}+\sum_{i=1}^{\ell-1}(\ell-i)(1-j_i)}b_\ell(\shift h\lambda_{j_1}\otimes\dots\otimes \shift h\lambda_{j_\ell})(\shift^{-1})^{\otimes n}\\
&=\sum_{\substack{\ell\neq 1\\j_1+\dots+j_\ell=n}}(-1)^{\frac{n(n-1)}{2}+\sum_{i<k}j_ij_k}b_\ell(\hat{h}\shift \lambda_{j_1}(\shift^{-1})^{\otimes j_1}\otimes \dots\otimes \hat{h}\shift \lambda_{j_\ell}(\shift^{-1})^{\otimes j_\ell})\\
&=\sum_{\substack{\ell\neq 1\\j_1+\dots+j_\ell=n}}b_\ell(\hat{h}\hat{\lambda}_{j_1}\otimes \dots\otimes \hat{h}\hat{\lambda_{j_\ell}}),
\end{split}
\end{align*}
where for the last equality we again use the equality of triangular numbers (in case $n=j_1+\dots+j_\ell$)
\[\frac{n(n-1)}{2}=\sum_{i=1}^{\ell} \frac{j_i(j_i-1)}{2} +\sum_{i<k}j_ij_k.\]

Now, using sign-free versions of all the equations involved, we prove equation \eqref{inductionhypothesiskernel} for $n\geq 2$ by induction. For the start of the induction notice that the last term vanishes since there is no $s$ with $1<s<2$ and that $\hat{\lambda}_2=b_2$. Therefore, for $n=2$ equation \eqref{inductionhypothesiskernel} reads as 
\[b_1b_2+b_2(b_1\otimes \id+\id\otimes b_1)=0,\]
which is exactly the $A_\infty$-relation for $\mathcal{A}$ for $n=2$. Now suppose that \eqref{inductionhypothesiskernel} holds for all numbers smaller than $n$. We claim that it is then also true for $n$. 

Applying $b_1$ to \eqref{signfreedefinitionkernel} we obtain that
\begin{equation}
\begin{split}\label{finalsteplambdakernel}
b_1\hat{\lambda}_n
&=\sum_{\substack{\ell\neq 1\\j_1+\dots+j_\ell=n}}b_1b_\ell(\hat{h}\hat{\lambda}_{j_1}\otimes \dots\otimes \hat{h}\hat{\lambda}_{j_\ell})\\
&=-\sum_{\substack{\ell\neq 1\\j_1+\dots+j_\ell=n}}\left(\underbrace{\sum_{r+t=\ell-1}b_\ell(\id^{\otimes r}\otimes b_1\otimes \id^{\otimes t})}_{D}+\underbrace{\sum_{\substack{r+s+t=\ell\\1<s<\ell}}b_{r+1+t}(\id^{\otimes r}\otimes b_s\otimes \id^{\otimes t})}_{E}\right)(\hat{h}\hat{\lambda}_{j_1}\otimes \dots\otimes \hat{h}\hat{\lambda}_{j_\ell})\\
\end{split}
\end{equation}
using the $A_\infty$-relations for the $b_k$. 
Considering the second summand, corresponding to $E$, and taking into account that $\hat{h}\hat{\lambda}_j$ is of degree $0$, we claim that
\begin{equation}
\begin{split}\label{kernelE}
\sum_{\substack{\ell\neq 1\\j_1+\dots+j_\ell=\ell}}&\sum_{\substack{r+s+t=\ell\\1<s<\ell}}b_{r+1+t}(\id^{\otimes r}\otimes b_s\otimes \id^{\otimes t})(\hat{h}\hat{\lambda}_{j_1}\otimes \dots\otimes \hat{h}\hat{\lambda}_{j_\ell})\\
&=\sum_{\substack{\ell'\neq 1\\k_1+\dots+k_{\ell'}=n\\k_{r+1}>1}}\sum_{r+t=\ell'-1}b_{\ell'}(\hat{h}\hat{\lambda}_{k_1}\otimes \dots\otimes \hat{h}\hat{\lambda}_{k_r}\otimes \hat{\lambda}_{k_{r+1}}\otimes \hat{h}\hat{\lambda}_{k_{r+2}}\otimes \dots\otimes \hat{h}\hat{\lambda}_{k_{\ell'}}).
\end{split}
\end{equation}
Indeed,
\begin{equation*}
\begin{split}
\sum_{\substack{\ell\neq 1\\j_1+\dots+j_\ell=\ell}}&\sum_{\substack{r+s+t=\ell\\1<s<\ell}}b_{r+1+t}(\id^{\otimes r}\otimes b_s\otimes \id^{\otimes t})(\hat{h}\hat{\lambda}_{j_1}\otimes \dots\otimes \hat{h}\hat{\lambda}_{j_\ell})\\
&=\sum_{\substack{\ell\neq 1\\j_1+\dots+j_\ell=n}}\sum_{\substack{r+s+t=\ell\\1<s<\ell}}b_{r+1+t}(\hat{h}\hat{\lambda}_{j_1}\otimes \dots\otimes \hat{h}\hat{\lambda}_{j_r}\otimes b_s(\hat{h}\hat{\lambda}_{j_{r+1}}\otimes \dots\otimes \hat{h}\hat{\lambda}_{j_{r+s}})\otimes \hat{h}\hat{\lambda}_{j_{r+s+1}}\otimes \dots\otimes \hat{h}\hat{\lambda}_{j_\ell})\\
&\stackrel{(\diamond)}{=}\sum_{\substack{\ell'\neq 1\\k_1+\dots+k_{\ell'}=n}}\sum_{r+t=\ell'-1}b_{r+1+t}(\hat{h}\hat{\lambda}_{k_1}\otimes \dots\otimes \hat{h}\hat{\lambda}_{k_r}\otimes \hat{\lambda}_{k_{r+1}}\otimes \hat{h}\hat{\lambda}_{k_{r+2}}\otimes \dots\otimes \hat{h}\hat{\lambda}_{k_{\ell'}})\\
&=\sum_{\substack{\ell'\neq 1\\k_1+\dots+k_{\ell'}=n}}\sum_{r+t=\ell'-1}b_{\ell'}(\hat{h}\hat{\lambda}_{k_1}\otimes \dots\otimes \hat{h}\hat{\lambda}_{k_r}\otimes \hat{\lambda}_{k_{r+1}}\otimes \hat{h}\hat{\lambda}_{k_{r+2}}\otimes \dots\otimes \hat{h}\hat{\lambda}_{k_{\ell'}}),
\end{split}
\end{equation*}
where for $(\diamond)$ we use \eqref{signfreedefinitionkernel} together with the substitutions $\ell'=\ell-s+1$ and $k_i=j_i$ for $i=1,\dots,r$, $k_{r+1}=j_{r+1}+\dots+j_{r+s}$ and $k_{r+1+i}=j_{r+s+i}$ for $i=1,\dots,t$. 
On the other hand, considering $D$ we obtain
\begin{equation}
\begin{split}\label{kernelD}
\sum_{\substack{\ell\neq 1\\j_1+\dots+j_\ell=n}}\sum_{r+t=\ell-1}&b_\ell(\id^{\otimes r}\otimes b_1\otimes \id^{\otimes t})(\hat{h}\hat{\lambda}_{j_1}\otimes \dots \otimes \hat{h}\hat{\lambda}_{j_\ell})\\
&=\sum_{\substack{\ell\neq 1\\j_1+\dots+j_\ell=n}}\sum_{r+t=\ell-1}b_\ell(\hat{h}\hat{\lambda}_{j_1}\otimes \dots\otimes \hat{\lambda}_{j_r}\otimes b_1\hat{h}\hat{\lambda}_{j_{r+1}}\otimes \hat{h}\hat{\lambda}_{j_{r+2}}\otimes \dots\otimes \hat{h}\hat{\lambda}_{j_\ell}).
\end{split}
\end{equation}
Substituting \eqref{kernelE} and \eqref{kernelD} back into \eqref{finalsteplambdakernel} we obtain
\begin{equation*}
\begin{split}
b_1\hat{\lambda}_n=&-\sum_{\substack{\ell\neq 1\\j_1+\dots+j_\ell=n\\j_{r+1}>1}}b_{\ell}(\hat{h}\hat{\lambda}_{j_1}\otimes \dots\otimes \hat{h}\hat{\lambda}_{j_r}\otimes (b_1\hat{h}+\id)\lambda_{j_{r+1}}\otimes \hat{h}\hat{\lambda}_{j_{r+2}}\otimes \dots\otimes \hat{h}\hat{\lambda}_{j_\ell})\\
&-\sum_{\substack{\ell\neq 1\\j_1+\dots+j_\ell=n\\j_{r+1}=1}}\sum_{r+t=\ell-1}b_{\ell}(\hat{h}\hat{\lambda}_{j_1}\otimes \dots\otimes \hat{\lambda}_{j_r}\otimes b_1\hat{h}\hat{\lambda}_{j_{r+1}}\otimes \hat{h}\hat{\lambda}_{j_{r+2}}\otimes \dots\otimes \hat{h}\hat{\lambda}_{j_\ell}).
\end{split}
\end{equation*}
Note that by definition $\mathbf{i}\mathbf{p}-\id=m_1h+hm_1$ and multiplying by $\shift$ from the left and $\shift^{-1}$ from the right we obtain $\hat{\boldsymbol\imath}\hat{\mathbf{p}}-\id=b_1\hat{h}+\hat{h}b_1$ and thus, 
\begin{equation*}
\begin{split}
b_1\hat{\lambda}_n=&\sum_{\substack{\ell\neq 1\\j_1+\dots+j_\ell=n\\j_{r+1}>1}}b_{\ell}(\hat{h}\hat{\lambda}_{j_1}\otimes \dots\otimes \hat{h}\hat{\lambda}_{j_r}\otimes (\underbrace{\hat{h}b_1}_{G}-\underbrace{\hat{\boldsymbol\imath}\hat{\mathbf{p}}}_{H})\lambda_{j_{r+1}}\otimes \hat{h}\hat{\lambda}_{j_{r+2}}\otimes \dots\otimes \hat{h}\hat{\lambda}_{j_\ell})
\\
&\underbrace{-\sum_{\substack{\ell\neq 1\\j_1+\dots+j_\ell=n\\j_{r+1}=1}}\sum_{r+t=\ell-1}b_{\ell}(\hat{h}\hat{\lambda}_{j_1}\otimes \dots\otimes \hat{\lambda}_{j_r}\otimes b_1\hat{h}\hat{\lambda}_{j_{r+1}}\otimes \hat{h}\hat{\lambda}_{j_{r+2}}\otimes \dots\otimes \hat{h}\hat{\lambda}_{j_\ell})}_{I}.
\end{split}
\end{equation*}
We now analyse the terms corresponding to $G$, $H$, $I$ separately. For $G$ we obtain
\begin{equation*}
\begin{split}
\sum_{\substack{\ell\neq 1\\j_1+\dots+j_\ell=n\\j_{r+1}>1}}&b_{\ell}(\hat{h}\hat{\lambda}_{j_1}\otimes \dots\otimes \hat{h}\hat{\lambda}_{j_r}\otimes \hat{h}b_1\hat{\lambda}_{j_{r+1}}\otimes \hat{h}\hat{\lambda}_{j_{r+2}}\otimes \dots\otimes \hat{h}\hat{\lambda}_{j_\ell})\\
=&\underbrace{-\sum_{\substack{\ell\neq 1\\j_1+\dots+j_\ell=n\\j_{r+1}>1}}b_{\ell}(\hat{h}\hat{\lambda}_{j_1}\otimes \dots\otimes \hat{h}\hat{\lambda}_{j_r}\otimes \sum_{u+w=j_{r+1}-1}\hat{h}\hat{\lambda}_{j_{r+1}}(\id^{\otimes u}\otimes b_1\otimes \id^{\otimes w})\otimes \hat{h}\hat{\lambda}_{j_{r+2}}\otimes \dots\otimes \hat{h}\hat{\lambda}_{j_\ell})}_{G1}\\
&\underbrace{-\sum_{\substack{\ell\neq 1\\j_1+\dots+j_\ell=n\\j_{r+1}>1}}b_{\ell}(\hat{h}\hat{\lambda}_{j_1}\otimes \dots\otimes \hat{h}\hat{\lambda}_{j_r}\otimes \sum_{u+v+w=j_{r+1}}\hat{h}\hat{\lambda}_{u+1+w}(\id^{\otimes u}\otimes \hat{\boldsymbol\imath}\hat{\mathbf{p}}\hat{\lambda}_v\otimes \id^{\otimes w})\otimes \hat{h}\hat{\lambda}_{j_{r+2}}\otimes \dots\otimes \hat{h}\hat{\lambda}_{j_\ell})}_{G2}
\end{split}
\end{equation*}
using the induction hypothesis. For $H$ we obtain
\begin{equation*}
\begin{split}
-\sum_{\substack{\ell\neq 1\\j_1+\dots+j_\ell=n\\j_{r+1}>1}}&b_{\ell}(\hat{h}\hat{\lambda}_{j_1}\otimes \dots\otimes \hat{h}\hat{\lambda}_{j_r}\otimes \hat{\boldsymbol\imath}\hat{\mathbf{p}}\hat{\lambda}_{j_{r+1}}\otimes \hat{h}\hat{\lambda}_{j_{r+2}}\otimes \dots\otimes \hat{h}\hat{\lambda}_{j_\ell})\\
&=-\sum_{\substack{\ell\neq 1\\j_1+\dots+j_\ell=n\\j_{r+1}>1}}b_{\ell}(\hat{h}\hat{\lambda}_{j_1}\otimes \dots\otimes \hat{h}\hat{\lambda}_{j_r}\otimes \hat{h}\hat{\lambda}_1 \otimes \hat{h}\hat{\lambda}_{j_{r+2}}\otimes \dots\otimes \hat{h}\hat{\lambda}_{j_\ell})(\id^{\otimes r}\otimes \hat{\imath}\hat{p}\hat{\lambda}_{j_{r+1}}\otimes \id^{\otimes t})
\end{split}
\end{equation*}
using that $\hat{h}\hat{\lambda}_1=\id$. Finally for $I$ we obtain
\begin{equation*}
\begin{split}
&-\sum_{\substack{\ell\neq 1\\j_1+\dots+j_\ell=n\\j_{r+1}=1}}\sum_{r+t=\ell-1}b_{\ell}(\hat{h}\hat{\lambda}_{j_1}\otimes \dots\otimes \hat{h}\hat{\lambda}_{j_r}\otimes b_1\hat{h}\hat{\lambda}_{j_{r+1}}\otimes \hat{h}\hat{\lambda}_{j_{r+2}}\otimes \dots\otimes \hat{h}\hat{\lambda}_{j_\ell})\\
=&-\sum_{\substack{\ell\neq 1\\j_1+\dots+j_\ell=n\\j_{r+1}=1}}\sum_{r+t=\ell-1}b_{\ell}(\hat{h}\hat{\lambda}_{j_1}\otimes \dots\otimes \hat{h}\hat{\lambda}_{j_r}\otimes \hat{h}\hat{\lambda}_{j_{r+1}}\otimes \hat{h}\hat{\lambda}_{j_{r+2}}\otimes \dots\otimes \hat{h}\hat{\lambda}_{j_\ell})(\id^{\otimes r}\otimes b_1\otimes \id^{\otimes t})
\end{split}
\end{equation*}
again using that $\hat{h}\hat{\lambda}_1=\id$. Using the definition of $\hat{\lambda}_n$ we see that the sum of $G1$ and $I$ is equal to
\[
-\sum_{r+t=n-1}\hat{\lambda}_n(\id^{\otimes r}\otimes b_1\otimes \id^{\otimes t})
\]
while the sum of $G2$ and $H$ is equal to
\[
-\sum_{r+s+t=n}\hat{\lambda}_{r+1+t}(\id^{\otimes r}\otimes \hat{\boldsymbol\imath}\hat{\mathbf{p}}\hat{\lambda}_s\otimes \id^{\otimes t}).
\]
The claim follows. 
\end{proof}

The proof of Kadeishvili's Theorem is now completed by the following lemma.

\begin{lem}\label{lem:lambdakernel}
Let $(\lambda_n)_{n\geq 2}$ be a $\lambda$-kernel. Then the $\tilde{m}_n$ defined by $\tilde{m}_n:=\mathbf{p}\lambda_n\mathbf{i}^{\otimes n}$ for $n\geq 2$ together with the original $\tilde{m}_1$ define an $A_\infty$-structure on $\mathcal{E}$ such that the map $(\mathbf{i}_n)_{n\geq 1}\colon \mathcal{E}\to \mathcal{A}$ defined by $\mathbf{i}_1=\mathbf{i}$, $\mathbf{i}_n=h\lambda_n\mathbf{i}^{\otimes n}$ is an $A_\infty$-quasi-isomorphism. 
\end{lem}

\begin{proof}
By Lemma \ref{lambdakernellemma} we have 
\begin{equation}\label{lambdakernellemma:eqn}
m_1\lambda_n+\sum_{r+t=n-1} (-1)^{n-1} \lambda_n (\id^{\otimes r}\otimes m_1\otimes \id^{\otimes t})+\sum_{r+s+t=n} (-1)^{r+st} \lambda_{r+1+t}(\id^{\otimes r}\otimes \mathbf{i}\mathbf{p}\lambda_s\otimes \id^{\otimes t})=0.\end{equation}
Applying $\mathbf{p}$ from the left and $\mathbf{i}^{\otimes n}$ from the right we obtain
\[\mathbf{p}m_1\lambda_n\mathbf{i}^{\otimes n}+\sum_{r+t=n-1} (-1)^{n-1}\mathbf{p}\lambda_n(\id^{\otimes r}\otimes m_1\otimes \id^{\otimes t})\mathbf{i}^{\otimes n}+\sum_{r+s+t=n} (-1)^{r+st} \mathbf{p}\lambda_{r+1+t}(\id^{\otimes r}\otimes \mathbf{i}\mathbf{p}\lambda_s\otimes \id^{\otimes t})\mathbf{i}^{\otimes n}=0.\]
Using that $\mathbf{p}$ and $\mathbf{i}$ commute with the differential as well as that they are of degree $0$ we obtain
\[\tilde{m}_1\mathbf{p}\lambda_n\mathbf{i}^{\otimes n}+\sum_{r+t=n-1}(-1)^{n-1}\mathbf{p}\lambda_n \mathbf{i}^{\otimes n} (\id^{\otimes r}\otimes \tilde{m}_1\otimes \id^{\otimes t})+\sum_{r+s+t=n} (-1)^{r+st} \mathbf{p}\lambda_{r+1+t}\mathbf{i}^{\otimes (r+1+t)}(\id^{\otimes r}\otimes \mathbf{p}\lambda_s\mathbf{i}^{\otimes s}\otimes \id^{\otimes t})=0.\]
Using the definition of $\tilde{m}_n$ for $n\geq 2$ this translates precisely to the $A_\infty$-algebra conditions on $\mathcal{E}$. 

On the other hand, applying $h$ from the left and $\mathbf{i}^{\otimes n}$ to \eqref{lambdakernellemma:eqn} yields the equation
\[
hm_1\lambda_n\mathbf{i}^{\otimes n}+\sum_{r+t=n-1}(-1)^{n-1} h\lambda_n(\id^{\otimes r}\otimes m_1\otimes \id^{\otimes t})\mathbf{i}^{\otimes n}+\sum_{\substack{r+s+t=n\\1<s<n}}(-1)^{r+st}h\lambda_{r+1+t}(\id^{\otimes r}\otimes \mathbf{i}\mathbf{p}\lambda_s\otimes \id^{\otimes t})\mathbf{i}^{\otimes n}=0.
\]
Using that $hm_1=\mathbf{i}\mathbf{p}-(\id+m_1h)$ and putting the terms coming from $\id+m_1h$ on the other side we obtain
\[
\mathbf{i}\mathbf{p}\lambda_n\mathbf{i}^{\otimes n}+\sum_{r+t=n-1} (-1)^{n-1} h\lambda_n(\id^{\otimes r}\otimes m_1\otimes \id^{\otimes t})\mathbf{i}^{\otimes n}+\sum_{\substack{r+s+t=n\\1<s<n}} (-1)^{r+st}h\lambda_{r+1+t}(\id^{\otimes r}\otimes \mathbf{i}\mathbf{p}\lambda_s\otimes \id^{\otimes t})\mathbf{i}^{\otimes n}=\lambda_n\mathbf{i}^{\otimes n}+m_1h\lambda_n\mathbf{i}^{\otimes n}.
\]
Using the inductive definition of $\lambda_n$ on the right hand side, and the fact that $\mathbf{i}$ is of degree $0$ and  compatible with the differential on the left hand side we obtain 
\begin{equation*}
\begin{split}
\mathbf{i}\mathbf{p}\lambda_n\mathbf{i}^{\otimes n}+\sum_{r+t=n-1} (-1)^{n-1}h\lambda_n\mathbf{i}^{\otimes n}(\id^{\otimes r}\otimes \tilde{m}_1\otimes \id^{\otimes t})+\sum_{\substack{r+s+t=n\\1<s<n}}(-1)^{r+st}h\lambda_{r+1+t}i^{\otimes (r+1+t)}(\id^{\otimes r}\otimes \mathbf{p}\lambda_s\mathbf{i}^{\otimes s}\otimes \id^{\otimes t})\\
=\left(\sum_{\substack{\ell\neq 1\\j_1+\dots+j_\ell=n}}(-1)^{\sum_{i=1}^{\ell-1}(\ell-i)(j_i-1)}m_\ell(h\lambda_{j_1}\otimes \dots\otimes h\lambda_{j_\ell})\right)\mathbf{i}^{\otimes n}+m_1h\lambda_n\mathbf{i}^{\otimes n}
\end{split}
\end{equation*}
Using the definitions of $\mathbf{i}_k$ and $\tilde{m}_k$ on the left hand side we obtain
\[
\sum_{r+s+t=n}(-1)^{r+st}\mathbf{i}_{r+1+t}(\id^{\otimes r}\otimes \tilde{m}_s\otimes \id^{\otimes t})=\sum_{\ell=1}^n \sum_{n=j_1+ \cdots +j_\ell} (-1)^{\sum_{i=1}^\ell (\ell-i)(j_i-1)}m_\ell(\mathbf{i}_{j_1}\otimes \dots\otimes \mathbf{i}_{j_\ell})
\]
which is precisely the condition for $(\mathbf{i}_n)_{n\geq 1}$ to be an $A_\infty$-morphism. 

The map $(\mathbf{i}_n)_{n\geq 1}$ is an $A_\infty$-quasi-isomorphisms since $\mathbf{i}_1=\mathbf{i}$ is a homotopy equivalence and hence a quasi-isomorphism of complexes. 
\end{proof}

\subsection{Special cases of Kadeishvili's Theorem}\label{sec:Kadeishvili2}

For later reference we first record some special cases of  Theorem \ref{kadeishvilitheorem}. 
We call an $A_\infty$-category \emphbf{minimal} if $m_1=0$. The following statement is well known, see e.g. \cite[Exercise 5]{Val14} and  \cite[Theorem 10.4.1]{LV12}. For a proof (with slightly different sign conventions), see \cite[Proposition 5.11]{Wal21}. 

\begin{lem}\label{minimalunique}
Let $(f_n)_{n\in \mathbb{N}}\colon \mathcal{A}\to \mathcal{B}$ be a morphism of $A_\infty$-categories. Then $(f_n)_{n\in \mathbb{N}}$ is an isomorphism if and only if $f_1$ is an isomorphism of $\mathbb{L}$-modules. In particular, two minimal $A_\infty$-categories are quasi-isomorphic if and only if they are isomorphic as $A_\infty$-categories.
\end{lem}

\begin{proof}
The forward direction is trivial. If $(f_n)_{n\in \mathbb{N}}$ is an isomorphism of $A_\infty$-categories with inverse $(g_n)_{n\in \mathbb{N}}$ then in particular $f_1\circ g_1=(f\circ g)_1=(\id)_1=\id$ and $g_1\circ f_1=(g\circ f)_1=(\id)_1=\id$. Therefore $f_1$ is an isomorphism of $\mathbb{L}$-modules. 

For the backward direction define the family of maps $(g_n)_{n\in \mathbb{N}}$ by 
\begin{align*} 
g_1&=f_1^{-1}\\
g_n&=\sum_{\ell=1}^{n-1} (-1)^{1+\sum_{i=1}^\ell(\ell-i)(j_i-1)}g_\ell(f_{j_1}\otimes \dots\otimes f_{j_\ell})(f_1^{-1})^{\otimes n}.
\end{align*}
Denote by $G$ the morphism of non-counital coalgebras $\overline{T}^c(\shift\mathcal{B})\to \overline{T}^c(\shift\mathcal{A})$ induced by  $\shift g_n (\shift^{\otimes n})^{-1}$ using the universal property of the tensor coalgebra. The $g_n$ were defined in such a way that $(g\circ f)=\id$ where composition is defined like composition of morphisms of $A_\infty$-categories and therefore we deduce $G\mathbf{B}(f)=\id$. Similarly, we can define the family of maps 
\begin{align*}
\tilde{g}_1&=f_1^{-1}\\
\tilde{g}_n&=\sum_{\ell=2}^{n} (-1)^{1+\sum_{i=1}^{\ell} (\ell-i)(j_i-1)} f_1^{-1} f_\ell(\tilde{g}_{j_1}\otimes \dots\otimes \tilde{g}_{j_\ell}).
\end{align*}   
Denote by $\tilde{G}$ the corresponding morphism of non-counital coalgebras. Analogously, we obtain $\mathbf{B}(f)\tilde{G}=\id$. Thus, $\mathbf{B}(f)$ is left and right invertible and therefore invertible with $G=\tilde{G}$, whence $g=\tilde{g}$. It remains to check that $G$ is a morphism of (non-counital) dg coalgebras, i.e. that it commutes with the codifferential, which is now straightforward:
\[bG=G\mathbf{B}(f)bG=Gb\mathbf{B}(f)G=Gb.\hfill\qedhere\]
\end{proof}

This allows us to conclude uniqueness of the $A_\infty$-structure provided by Theorem \ref{kadeishvilitheorem} in the case when $\mathcal{E}$ is the homology of $\mathcal{A}$. 

\begin{cor}\label{kadeishvilitheoremhomology}
Let $\mathcal{A}$ be an $A_\infty$-category with multiplications $m_n$, $n\geq 1$. Let $\mathcal{E}=H^{\bullet}(\mathcal{A})$ be its cohomology. Choose $\mathbf{i}\colon \mathcal{E}\to \mathcal{A}$ and $\mathbf{p}\colon \mathcal{A}\to \mathcal{E}$ compatible with the differential such that $\mathbf{i}\mathbf{p}-\id_{\mathcal{A}}=m_1h+hm_1$ for some $h\colon \mathcal{A}\to \mathcal{A}$ of degree $-1$ (note that this is possible since we are working over a field). Then there exists an $A_\infty$-structure $\tilde{m}_n$, $n\geq 1$ on $\mathcal{E}$ with $\tilde{m}_1=0$ together with an $A_\infty$-quasi-isomorphism $(\mathbf{i}_n)_{n\geq 1}\colon \mathcal{E}\to \mathcal{A}$ such that $\mathbf{i}_1=\mathbf{i}$. Moreover, such an $A_\infty$-structure is unique up to (non-unique) isomorphism of $A_\infty$-categories.  
\end{cor}

\begin{proof}
Existence of $\tilde{m}_n$ and $\mathbf{i}_n$ follow immediately from the proof of the above theorem. Uniqueness follows from Lemma \ref{minimalunique} as $H^{\bullet}(\mathcal{A})$ is a minimal $A_\infty$-category.
\end{proof}

Another consequence of Theorem \ref{kadeishvilitheorem} is that every $A_\infty$-category is quasi-isomorphic to a minimal $A_\infty$-category. Specialising even further, we recover Merkulov's construction on the $A_\infty$-structure on the cohomology of a dg category.

\begin{thm}\label{merkulovsconstruction}
Let $\mathcal{A}$ be a dg category with differential $d$ and multiplication $m$. Let $\mathcal{E}=H^{\bullet}(\mathcal{A})$ be its cohomology. Choose a complement $\mathcal{L}^{\bullet}$ of the cycles $\mathcal{Z}^{\bullet}=\ker d$ in $\mathcal{A}$ and a complement $\mathcal{H}^{\bullet}$ of the boundaries $\mathcal{B}^{\bullet}=\im d$ in $Z^{\bullet}$. Choose an isomorphism $\mathbf{i}\colon \mathcal{E}\to \mathcal{H}^{\bullet}$. Then there exists an $A_\infty$-structure $\tilde{m}_n$, $n\geq 1$ on $\mathcal{E}$ with $\tilde{m}_1=0$ and $\tilde{m}_2$ induced by the multiplication $m$ together with an $A_\infty$-quasi-isomorphism $(\mathbf{i}_n)_{n\geq 1}\colon \mathcal{E}\to \mathcal{A}$ such that $\mathbf{i}_1=\mathbf{i}$. Moreover, such an $A_\infty$-structure is unique up to (non-unique) isomorphism of $A_\infty$-categories. 
\end{thm}

\begin{proof}
This theorem follows immediately from Corollary \ref{kadeishvilitheoremhomology} once we construct the map $h\colon \mathcal{A}\to \mathcal{A}$ of degree $-1$ such that $\mathbf{i}\mathbf{p}-\id_{\mathcal{A}}=dh+hd$. Given the direct sum decomposition $\mathcal{A}=\mathcal{B}^{\bullet}\oplus \mathcal{H}^{\bullet}\oplus \mathcal{L}^{\bullet}$ we see that $d$ induces an isomorphism $\mathcal{L}^{\bullet}\to \mathcal{B}^{\bullet}$ (of degree $1$) which by slight abuse of notation we also denote by $d$. We define $h\colon \mathcal{A}\to \mathcal{A}$ by $h|_{\mathcal{H}^{\bullet}\oplus \mathcal{L}^{\bullet}}=0$ and $h|_{\mathcal{B}^{\bullet}}=-d^{-1}$. We check the claim $\mathbf{i}\mathbf{p}-\id=dh+hd$ individually on the direct summands $\mathcal{L}^{\bullet}$, $\mathcal{B}^{\bullet}$, and $\mathcal{H}^{\bullet}$. For $\mathcal{H}^{\bullet}$ we obtain that $\mathbf{i}\mathbf{p}|_{\mathcal{H}^{\bullet}}=\id$ and $(dh+hd)_{\mathcal{H}^{\bullet}}=0$; for $L^{\bullet}$ we obtain that $\mathbf{i}\mathbf{p}|_{L^{\bullet}}=0$, $dh|_{\mathcal{L}^{\bullet}}=0$, but $hd|_{\mathcal{L}^{\bullet}}=-\id_{\mathcal{L}^{\bullet}}$; and for $\mathcal{B}^{\bullet}$ we obtain that $\mathbf{i}\mathbf{p}|_{\mathcal{B}^{\bullet}}=0$, $hd|_{\mathcal{B}^{\bullet}}=0$, and $dh|_{\mathcal{B}^{\bullet}}=-\id_{B^{\bullet}}$. Thus, $\mathbf{i}\mathbf{p}-\id=dh+hd$ holds on all three summands and therefore we can use Corollary \ref{kadeishvilitheoremhomology} to establish the claim. 
\end{proof}

\subsection{Counitality}\label{sec:counitality}
In this section, we give conditions under which an $A_\infty$-coalgebra object whose cohomology is counital is itself strictly counital. This statement is dual to a corresponding statement for $A_\infty$-categories in \cite[Lemma 2.1]{Sei08}. The proof heavily uses the following instance of structure transport (which Seidel calls `formal diffeomorphism' in \cite[(1c)]{Sei08}):

\begin{lem}\label{structuretransport}
Suppose that $(\mathcal{C}, \mu_n)$ is an $A_\infty$-coalgebra object and let $\tilde{\mathcal{C}}$ be a graded object in $\mathsf{C}$, which is isomorphic to $\mathcal{C}$ as a graded object. Let $(f_n\colon \mathcal{C}\to \tilde{\mathcal{C}}^{\otimes n})_{n\geq 0}$ be a sequence of maps such that $f_n$ is of degree $1-n$ such that $f_1$ is an isomorphism in the category of graded objects in $\mathsf{C}$. Then there exists a unique $A_\infty$-coalgebra object structure $(\tilde{\mathcal{C}}, \tilde{\mu}_n)$ on $\tilde{\mathcal{C}}$ such that $f=(f_n)_{n\geq 0}\colon \mathcal{C}\to \tilde{\mathcal{C}}$ is an $A_\infty$-quasi-isomorphism. 
\end{lem}

\begin{proof}
Let $b$ be the codifferential on the cobar construction $\mathbf{\Omega}(\mathcal{C})$. Furthermore, denote by $F$ the unique morphism of (non-unital) algebras $\mathbf{\Omega}(\mathcal{C})\to \overline{T}(\shift^{-1}\tilde{\mathcal{C}})$ induced by  $(\shift^{\otimes n})^{-1}f_n\shift$ using the universal property of the tensor algebra. To define on $\tilde{C}$ an $A_\infty$-coalgebra object structure such that $f$ is an $A_\infty$-quasi-isomorphism is the same as defining on $\overline{T}(\shift^{-1}\tilde{\mathcal{C}})$ the structure of a (non-unital) dg algebra such that $F$ is a morphism of dg algebras. Such a structure is uniquely defined by the requirement that $F$ is a morphism of dg algebras, i.e. $\tilde{b}F=Fb$, or (as $F$ is invertible as a morphism of algebras) $\tilde{b}=FbF^{-1}$. We check that defining $\tilde{b}$ in this way gives indeed a dg algebra structure, i.e. that it is a coderivation which squares to zero. Denote by $m$, resp. $\tilde{m}$, the multiplication on $\mathbf{\Omega}(\mathcal{C})$, resp   $\overline{T}(\shift^{-1}\tilde{\mathcal{C}})$. Then, using that $F$ and $F^{-1}$ are morphisms of algebras and that $b$ is a coderivation, we obtain 
\begin{align*}
\tilde{b}\tilde{m}&=FbF^{-1}\tilde{m}=Fbm(F^{-1}\otimes F^{-1})=Fm(b\otimes \id+\id\otimes b)(F^{-1}\otimes F^{-1})\\
&=m(F\otimes F)(b\otimes \id +\id\otimes b)(F^{-1}\otimes F^{-1})=m(\tilde{b}\otimes \id+\id\otimes \tilde{b}).
\end{align*}
That $\tilde{b}$ squares to zero can be seen from
\[\tilde{b}^{2}=(FbF^{-1})^2=Fb^2F^{-1}=0.\]
Thus $\tilde{b}$ uniquely defines an $A_\infty$-coalgebra object structure on $\tilde{\mathcal{C}}$ whose cobar construction has codifferential $\tilde{b}$. 
\end{proof}

We use this lemma to prove the condition on when a cohomologically counital coalgebra is strictly counital. The additional assumption on the existence of a lift of $\mu_2$, when compared to \cite[Lemma 2.1]{Sei08}, is needed since we work in a not necessarily semisimple monoidal category.

\begin{prop}\label{counitality}
Let $(\mathcal{C},\mu_i)$ be an $A_\infty$-coalgebra object in $\mathsf{C}$ and assume that $H^{\bullet}(\mathcal{C})$ admits a strictly counital $A_\infty$-coalgebra structure with a strict $A_\infty$-quasi-isomorphism $q\colon \mathcal{C}\to H^{\bullet}(\mathcal{C})$. Define $\tau=\varepsilon q$, where $\varepsilon$ denotes the counit $H^{\bullet}(\mathcal{C})\to e$.

Assume that there exists a lift $\tilde{\mu}_2\colon \mathcal{C}\to\mathcal{C}\otimes\mathcal{C}$ of the comultiplication in $H^{\bullet}(\mathcal{C})$, homotopic to $\mu_2$, and satisfying $(\id\otimes \tau)\tilde{\mu}_2=\id=(\tau\otimes \id)\tilde{\mu}_2$. Then there exists a strictly counital  $A_\infty$-coalgebra structure $(\tilde{\mu}_i)$ on $\mathcal{C}$ and a morphism of  $A_\infty$-coalgebra objects $f\colon (\mathcal{C},\tilde{\mu}_i)\to (\mathcal{C},\mu_i)$ such that $f_1$ is an isomorphism.
\end{prop}

\begin{proof}
Recall that $(\mathcal{C}, \tilde{\mu}_i)$ is strictly counital if  $\tau \tilde{\mu}_1=0$, $(\id\otimes \tau)\tilde{\mu}_2=\id=(\tau\otimes \id)\tilde{\mu}_2$ and $(\id^{\otimes r}\otimes \tau\otimes \id^{\otimes (i-r-1)})\tilde{\mu}_i=0$ for all $i>2, r=1,\dots,i-1$.

Using Lemmma \ref{structuretransport} we will inductively construct a sequence of $A_\infty$-coalgebra objects $(\mathcal{C},\mu_i^{(n,d)})$ together with $A_\infty$-quasi-isomorphisms $f^{(n,d)}\colon(\mathcal{C},\mu_i^{(n,d+1)})\to(\mathcal{C},\mu_i^{(n,d)})$, for natural numbers $n\geq 3$ and $d=0,\dots,n$, and such that $
(\mathcal{C},\mu_i^{(n,n)}) =(\mathcal{C},\mu_i^{(n+1,0)})$, where $(\mathcal{C},\mu_i^{(n,d)})$ satisfies
\begin{equation}\label{taumu1}
 \tau{\mu}^{(n,d)}_1=0,
\end{equation}
\begin{equation} \label{taumu2}
 (\id\otimes \tau){\mu}^{(n,d)}_2=\id=(\tau\otimes \id){\mu}^{(n,d)}_2,
\end{equation}
and 
\begin{equation*}(I_{n,d}) \qquad\qquad
(\id^{\otimes r}\otimes \tau\otimes \id^{\otimes (i-r-1)}){\mu}^{(n,d)}_i=0 \quad \text{ for all } 2<i<n, r=1,\dots,i-1  \text{ and for } i=n  \text{ and } r<d.
\end{equation*}

First note that $\tau{\mu}_1= \varepsilon q\mu_1=0$ as $q_1\mu_1= \mu_1^{H^{\bullet}(\mathcal{C})} q_1=0$ since $q_1$ is a morphism of complexes.

We start by constructing an $A_\infty$-coalgebra object $(\mathcal{C},\mu_i^{(3,0)})$, together with an $A_\infty$-quasi-isomorphism $f^{(3,0)}\colon (\mathcal{C},\mu_i^{(3,0)})\to (\mathcal{C},\mu_i)$ which satisfies only the conditions. 
 \eqref{taumu1} and \eqref{taumu2}. To achieve this, we set $\mu_1^{(3,0)} = \mu_1$ and choose $\mu_2^{(3,0)}$ as $\tilde{\mu}_2$ in the statement of the proposition. For  $f^{(3,0)}\colon (\mathcal{C},\mu_i^{(3,0)})\to (\mathcal{C},\mu_i)$, we define $f^{(3,0)}_1$ to be the identity and $f^{(3,0)}_2$ to be a chain homotopy from ${\mu}_2^{(3,0)}$ to $\mu_2$. Then \eqref{taumu1} is satisfied since $\tau\mu_1^{(3,0)} = \tau\mu_1=0$ and \eqref{taumu2} is satisfied by assumption.

We now assume inductively that we have $(\mathcal{C},\mu_i^{(n,d)})$ and wish to construct $(\mathcal{C},\mu_i^{(n,d+1)})$. For notational simplicity, we will write $(\mathcal{C},\mu_i)$ for $(\mathcal{C},\mu_i^{(n,d)})$ and similarly $(\mathcal{C},\tilde{\mu}_i)$ for $(\mathcal{C},\mu_i^{(n,d+1)})$. Furthermore, we will denote the quasi-isomorphism $f^{(n,d)}$ simply by $f$.
 
We choose $f_1$ to be the identity, $f_2, \cdots, f_{n-2}$ to be zero, 
\[f_{n-1}=(-1)^{n-d} (\id^{\otimes d}\otimes \tau\otimes \id^{\otimes (n-d-1)})\mu_n\]
\[f_{n}=(-1)^{n-d+1} (\id^{\otimes d}\otimes \tau\otimes \id^{\otimes (n-d)})\mu_{n+1}. \]
By Lemma \ref{structuretransport} there is a unique $A_{\infty}$-structure $(\tilde{\mu}_i)_{i\geq 1}$ on $\mathcal{C}$ turning this $f$ into a morphism of $A_\infty$-coalgebra objects satisfying our conditions.

From the construction given in the proof of Lemma \ref{structuretransport} is obvious that $\tilde{\mu}_i=\mu_i$ for $i=1, \dots, n-2$, so the first two conditions are satisfied, since they were already satisfied for $(\mathcal{C},\mu_i) = (\mathcal{C},\mu_i^{(n,d)}) $. We next prove $\tilde{\mu}_{n-1}=\mu_{n-1}$. By the assumption that $f$ is an $A_\infty$-morphism, 
\begin{equation}\label{eq55}
\underbrace{\sum_{r+s+t=n-1} (-1)^{rs+t}(\id^{\otimes r}\otimes \mu_s\otimes \id^{\otimes t})f_{r+1+t}}_{\mathbf{L}}=\underbrace{\sum_{j_1+\cdots+ j_q= n-1} (-1)^{\sum_{i=1}^q(i-1)(j_i+1)}(f_{j_1}\otimes f_{j_2}\otimes \dots\otimes f_{j_q})\tilde{\mu}_q}_{\mathbf{R}}.
\end{equation}

We first consider the left hand side $\mathbf{L}$ of the equation, and, using the definition of $f$ (in particular the fact that $f_i=0$ for $i=2,\dots,n-2$), obtain
\begin{equation*}
\begin{split}
\mathbf{L}&= \mu_{n-1}f_1 + \sum_{r=0}^{n-2} (-1)^{n-2}(\id^{\otimes r}\otimes \mu_1\otimes \id^{\otimes (n-2-r)})f_{n-1}\\
&=  \mu_{n-1}+ \underbrace{\sum_{r=0}^{n-2} (-1)^{d} (\id^{\otimes r}\otimes \mu_1\otimes \id^{\otimes (n-2-r)})(\id^{\otimes d}\otimes \tau\otimes \id^{\otimes (n-d-1)})\mu_n}_{(A)}.
\end{split}
\end{equation*}
Using 
\[(\id^{\otimes r}\otimes \mu_1\otimes \id^{\otimes (n-2-r)})(\id^{\otimes d}\otimes \tau\otimes \id^{\otimes (n-d-1)}) = 
\begin{cases}
 (\id^{\otimes d}\otimes \tau\otimes \id^{\otimes (n-d-1)})(\id^{\otimes r}\otimes \mu_1\otimes \id^{\otimes (n-1-r)})& \text{if } r<d \\                                                                                                              
 (\id^{\otimes d}\otimes \tau\otimes \id^{\otimes (n-d-1)})(\id^{\otimes (r+1)}\otimes \mu_1\otimes \id^{\otimes (n-2-r)}) & \text{if } r\geq d                                                                                                \end{cases}
 \]
and $(\id^{\otimes d}\otimes \tau\otimes \id^{\otimes (n-d-1)})(\id^{\otimes d}\otimes \mu_1\otimes \id^{\otimes (n-1-r)}) =0$
we see that
\begin{equation}
\label{eq65}
(A)= \sum_{r=0}^{n-1} (-1)^{d}  (\id^{\otimes d}\otimes \tau\otimes \id^{\otimes (n-d-1)})(\id^{\otimes r}\otimes \mu_1\otimes \id^{\otimes (n-1-r)})\mu_n.
\end{equation}

We now consider the right hand side $\mathbf{R}$ of \eqref{eq55}, which, using the definition of $f$, is
\begin{equation*}
\begin{split}\mathbf{R} &= (f_1\otimes \cdots \otimes f_1)\tilde{\mu}_{n-1} + f_{n-1}\tilde{\mu}_1\\
&= \tilde{\mu}_{n-1} + (-1)^{n-d} (\id^{\otimes d}\otimes \tau\otimes \id^{\otimes (n-d-1)})\mu_n\mu_1.
\end{split}
\end{equation*}
Replacing 
\[\mu_n\mu_1 = \sum_{\substack{r+s+t=n\\s<n}} (-1)^{rs+t+1}(\id^{\otimes r}\otimes \mu_s\otimes \id^{\otimes t})\mu_{r+1+t},\]
we rewrite
\begin{equation*}
\begin{split}
\mathbf{R} = \tilde{\mu}_{n-1} & + \underbrace{\sum_{r=0}^{n-1} (-1)^{d} (\id^{\otimes d}\otimes \tau\otimes \id^{\otimes (n-d-1)})(\id^{\otimes r}\otimes \mu_1\otimes \id^{\otimes (n-r-1)})\mu_{n}}_{(\tilde{A})}\\
&+ \underbrace{\sum_{r=0}^{n-2} (-1)^{d-r+1} (\id^{\otimes d}\otimes \tau\otimes \id^{\otimes (n-d-1)})(\id^{\otimes r}\otimes \mu_2\otimes \id^{\otimes (n-r-2)})\mu_{n-1}}_{(B)}\\
&+\underbrace{\sum_{\substack{r+s+t=n\\2<s<n-1}} (-1)^{n-d+rs+t+1}(\id^{\otimes d}\otimes \tau\otimes \id^{\otimes (n-d-1)})(\id^{\otimes r}\otimes \mu_s\otimes \id^{\otimes t})\mu_{r+1+t}}_{(C)}\\
&+ \underbrace{(-1)^{n-d}(\id^{\otimes d}\otimes \tau\otimes \id^{\otimes (n-d-1)})(\mu_{n-1}\otimes \id)\mu_2}_{(D)}
+\underbrace{ (-1)^{d}(\id^{\otimes d}\otimes \tau\otimes \id^{\otimes (n-d-1)})(\id\otimes \mu_{n-1})\mu_2}_{(E)}.
\end{split}
\end{equation*}
Every summand in $(C)$ is zero by the assumption that $(\mathcal{C},\mu_i)$ satisfies $I_{n,d}$ since $\tau$ hits either $\mu_s$ or $\mu_{r+1+t}$ and both indices are strictly smaller than $n$. Furthermore, $(A)=(\tilde A)$, so to prove that $\tilde{\mu}_{n-1}=\mu_{n-1}$ it suffices to verify $(B)+(D)+(E)=0.$
For $0<d<n-1$, both $(D)$ and $(E)$ are zero, and all terms in $(B)$ are zero, except for those indexed by $r=d-1$ and $r=d$, which equal $\mu_{n-1}$ with opposite signs and hence cancel.
For $d=0$, $(E)=\mu_{n-1}$, while $(D)=0$ and the only nonzero term in $(B)$ is $-\mu_{n-1}$ for $r=0$. Similarly, for $d=n-1$, $(D) = -\mu_{n-1}$, while $(E)=0$ and the only nonzero term in $(B)$ is $\mu_{n-1}$ for $r=n-2$. Hence $\tilde{\mu}_{n-1}=\mu_{n-1}$ and the conditions in $I_{(d+1,n)}$ for $i<n$ hold.

In order to complete the proof, we need to show
$(\id^{\otimes k}\otimes \tau\otimes \id^{\otimes n-k-1})\tilde{\mu}_{n} =0 $ for $k=0\dots, d$. To achieve this, we first determine $\tilde{\mu}_{n}$. The requirement that $f$ is an $A_\infty$-morphism yields
\begin{equation}\label{eq56}
\underbrace{\sum_{r+s+t=n} (-1)^{rs+t}(\id^{\otimes r}\otimes \mu_s\otimes \id^{\otimes t})f_{r+1+t}}_{\mathbf{L'}}=\underbrace{\sum_{j_1+\cdots +j_q= n} (-1)^{\sum_{i=1}^q(i-1)(j_i+1)}(f_{j_1}\otimes f_{j_2}\otimes \dots\otimes f_{j_q})\tilde{\mu}_q}_{\mathbf{R'}}.
\end{equation}

Expanding the left hand side $\mathbf{L'}$, using that $f_2,\dots, f_{n-2}$ are zero,  produces
\begin{equation*}
\begin{split}
\mathbf{L'}&= \mu_nf_1+ \sum_{r=0}^{n-2} (-1)^{n-r}(\id^{\otimes r} \otimes \mu_2\otimes \id^{\otimes (n-2-r)} )f_{n-1} + \sum_{r=0}^{n-1} (-1)^{n-1}(\id^{\otimes r} \otimes \mu_1\otimes \id^{\otimes (n-1-r)} )f_{n}\\
&= \mu_n +\sum_{r=0}^{n-2} (-1)^{d-r}(\id^{\otimes r} \otimes \mu_2\otimes \id^{\otimes (n-2-r)} )(\id^{\otimes d}\otimes \tau\otimes \id^{\otimes (n-d-1)})\mu_n \\ &\qquad \,
+ \sum_{r=0}^{n-1} (-1)^{d}(\id^{\otimes r} \otimes \mu_1\otimes \id^{\otimes (n-1-r)} )(\id^{\otimes d}\otimes \tau\otimes \id^{\otimes (n-d)})\mu_{n+1}
\end{split}
\end{equation*}
and, again (similarly to \eqref{eq65}) using
\[\sum_{r=0}^{n-1} (-1)^{d}(\id^{\otimes r} \otimes \mu_1\otimes \id^{\otimes (n-1-r)} )(\id^{\otimes d}\otimes \tau\otimes \id^{\otimes (n-d)}) = \sum_{r=0}^{n} (-1)^{d}(\id^{\otimes d}\otimes \tau\otimes \id^{\otimes (n-d)})(\id^{\otimes r} \otimes \mu_1\otimes \id^{\otimes (n-r)} ), \]
and
\[(\id^{\otimes r} \otimes \mu_2\otimes \id^{\otimes (n-2-r)} )(\id^{\otimes d}\otimes \tau\otimes \id^{\otimes (n-d-1)})=
\begin{cases}
(\id^{\otimes (d+1)}\otimes \tau\otimes \id^{\otimes (n-d-1)})(\id^{\otimes r} \otimes \mu_2\otimes \id^{\otimes (n-1-r)} ) & \text{ if } 0\leq r\leq d-1\\
(\id^{\otimes d}\otimes \tau\otimes \id^{\otimes (n-d)})(\id^{\otimes (r+1)} \otimes \mu_2\otimes \id^{\otimes (n-r-2)} ) & \text{ if } d\leq r\leq n-2,
\end{cases}
\]
can be rewritten to
\begin{equation*}
\begin{split}
\mathbf{L'}
&= \mu_n 
+\underbrace{\sum_{r=0}^{d-1} (-1)^{d-r} (\id^{\otimes (d+1)}\otimes \tau\otimes \id^{\otimes (n-d-1)})(\id^{\otimes r} \otimes \mu_2\otimes \id^{\otimes (n-1-r)} )\mu_n}_{A}\\&\qquad \,
+\underbrace{\sum_{r=d+1}^{n-1} (-1)^{d-r+1}(\id^{\otimes d}\otimes \tau\otimes \id^{\otimes (n-d)})(\id^{\otimes r} \otimes \mu_2\otimes \id^{\otimes (n-r-1)} )\mu_n}_{B}
\\ &\qquad \,+ \underbrace{\sum_{r=0}^{n} (-1)^{d}(\id^{\otimes d}\otimes \tau\otimes \id^{\otimes (n-d)})(\id^{\otimes r} \otimes \mu_1\otimes \id^{\otimes (n-r)} )\mu_{n+1}}_{C}.
\end{split}
\end{equation*}

On the other hand, expanding the right hand side $\mathbf{R'}$, using the definition of $f$, yields
\begin{equation*}
\begin{split}
\mathbf{R'}&=(f_1\otimes \dots\otimes f_1)\tilde{\mu}_n +(f_{n-1}\otimes f_1)\tilde{\mu}_2 +(-1)^n(f_1\otimes f_{n-1})\tilde{\mu}_2 +f_n\tilde{\mu}_1\\
&= \tilde{\mu}_n +(-1)^{n-d}((\id^{\otimes d}\otimes \tau\otimes \id^{\otimes (n-d-1)})\mu_n \otimes \id)\mu_2 
+(-1)^{d}(\id\otimes (\id^{\otimes d}\otimes \tau\otimes \id^{\otimes (n-d-1)})\mu_n)\mu_2\\&\quad
+(-1)^{n-d+1}(\id^{\otimes d}\otimes \tau\otimes \id^{\otimes (n-d)})\mu_{n+1}\mu_1\\
&= \tilde{\mu}_n +\underbrace{(-1)^{n-d}(\id^{\otimes d}\otimes \tau\otimes \id^{\otimes (n-d)})(\mu_n \otimes \id)\mu_2}_{D} 
+\underbrace{(-1)^{d}(\id^{\otimes (d+1)}\otimes \tau\otimes \id^{\otimes (n-d-1)})(\id\otimes\mu_n)\mu_2}_{E}\\&\quad
+\underbrace{(-1)^{n-d+1}(\id^{\otimes d}\otimes \tau\otimes \id^{\otimes (n-d)})\mu_{n+1}\mu_1}_{F}.
\end{split}
\end{equation*}
Thus
\begin{equation*}
\begin{split}
\tilde{\mu}_n &= \mu_n
+\underbrace{\sum_{r=0}^{d-1} (-1)^{d-r} (\id^{\otimes (d+1)}\otimes \tau\otimes \id^{\otimes (n-d-1)})(\id^{\otimes r} \otimes \mu_2\otimes \id^{\otimes (n-1-r)} )\mu_n}_{A}\\&\qquad \,
+\underbrace{\sum_{r=d+1}^{n-1} (-1)^{d-r+1}(\id^{\otimes d}\otimes \tau\otimes \id^{\otimes (n-d)})(\id^{\otimes r} \otimes \mu_2\otimes \id^{\otimes (n-r-1)} )\mu_n}_{B}\\ &\qquad \,
+ \underbrace{\sum_{r=0}^{n} (-1)^{d}(\id^{\otimes d}\otimes \tau\otimes \id^{\otimes (n-d)})(\id^{\otimes r} \otimes \mu_1\otimes \id^{\otimes (n-r)} )\mu_{n+1}}_{C}\\ &\qquad \,
+\underbrace{(-1)^{n-d+1}(\id^{\otimes d}\otimes \tau\otimes \id^{\otimes (n-d)})(\mu_n \otimes \id)\mu_2}_{D} 
+\underbrace{(-1)^{d+1}(\id^{\otimes (d+1)}\otimes \tau\otimes \id^{\otimes (n-d-1)})(\id\otimes\mu_n)\mu_2}_{E}\\&\qquad\,
+\underbrace{(-1)^{n-d}(\id^{\otimes d}\otimes \tau\otimes \id^{\otimes (n-d)})\mu_{n+1}\mu_1}_{F}.
\end{split}
\end{equation*}
We next observe that 
\begin{equation*}
\begin{split}
C&+F=(-1)^{n-d}(\id^{\otimes d}\otimes \tau\otimes \id^{\otimes (n-d)})\left[\mu_{n+1}\mu_1 +\sum_{r=0}^{n} (-1)^{n}(\id^{\otimes r} \otimes \mu_1\otimes \id^{\otimes (n-r)} )\mu_{n+1}\right]\\
&=(-1)^{n-d+1}(\id^{\otimes d}\otimes \tau\otimes \id^{\otimes (n-d)})\left[ \sum_{\substack{r+s+t=n+1\\1<s<n+1}} (-1)^{rs+t}(\id^{\otimes r}\otimes \mu_s\otimes \id^{\otimes t})\mu_{r+1+t}\right]\\
&=(-1)^{n-d+1}(\id^{\otimes d}\otimes \tau\otimes \id^{\otimes (n-d)})\left[ \sum_{r=0}^{n-1} (-1)^{n-r-1}(\id^{\otimes r}\otimes \mu_2\otimes \id^{\otimes (n-r-1)})\mu_{n} - \underbrace{(\mu_n\otimes\id)\mu_2}_{D'} +\underbrace{(-1)^{n}(\id\otimes \mu_n)\mu_2}_{E'}\right],
\end{split}
\end{equation*}
where the last equality uses that all other summands are zero thanks to the induction hypothesis $I_{n,d}$.

Now,  $D$ and $D'$ cancel, and 
\[E'= (-1)^{d+1}(\id\otimes(\id^{\otimes (d-1)}\otimes \tau\otimes \id^{\otimes (n-d)})\mu_n)\mu_2=0\]
by the induction hypothesis, so 
\begin{equation*}
\begin{split}
C+D+F&= (-1)^{n-d+1}(\id^{\otimes d}\otimes \tau\otimes \id^{\otimes (n-d)}) \sum_{r=0}^{n-1} (-1)^{n-r-1}(\id^{\otimes r}\otimes \mu_2\otimes \id^{\otimes (n-r-1)})\mu_{n} \\
&=\sum_{r=0}^{n-1} (-1)^{r-d}(\id^{\otimes d}\otimes \tau\otimes \id^{\otimes (n-d)}) (\id^{\otimes r}\otimes \mu_2\otimes \id^{\otimes (n-r-1)})\mu_{n}.
\end{split}
\end{equation*}
Noting that 
\begin{equation*}
\begin{split}
B'&:=B+C+D+F\\
&= \sum_{r=d+1}^{n-1} (-1)^{d-r+1}(\id^{\otimes d}\otimes \tau\otimes \id^{\otimes (n-d)})(\id^{\otimes r} \otimes \mu_2\otimes \id^{\otimes (n-r-1)} )\mu_n\\
&\qquad +\sum_{r=0}^{n-1} (-1)^{r-d}(\id^{\otimes d}\otimes \tau\otimes \id^{\otimes (n-d)}) (\id^{\otimes r}\otimes \mu_2\otimes \id^{\otimes (n-r-1)})\mu_{n}\\
&=\sum_{r=0}^{d}(-1)^{d-r}(\id^{\otimes d}\otimes \tau\otimes \id^{\otimes (n-d)})(\id^{\otimes r}\otimes \mu_2\otimes \id^{\otimes (n-r-1)})\mu_{n},
\end{split}
\end{equation*}
we are left with
\begin{equation*}
\begin{split}
\tilde{\mu}_n &= \mu_n
+\underbrace{\sum_{r=0}^{d-1} (-1)^{d-r} (\id^{\otimes (d+1)}\otimes \tau\otimes \id^{\otimes (n-d-1)})(\id^{\otimes r} \otimes \mu_2\otimes \id^{\otimes (n-1-r)} )\mu_n}_{A}\\&\qquad \,
+ \underbrace{\sum_{r=0}^{d}(-1)^{d-r}(\id^{\otimes d}\otimes \tau\otimes \id^{\otimes (n-d)})(\id^{\otimes r}\otimes \mu_2\otimes \id^{\otimes (n-r-1)})\mu_{n}}_{B'}\\&\qquad \,+\underbrace{(-1)^{d+1}(\id^{\otimes (d+1)}\otimes \tau\otimes \id^{\otimes (n-d-1)})(\id\otimes\mu_n)\mu_2}_{E}.
\end{split}
\end{equation*}

We now apply 
$T_k:=(\id^{\otimes k}\otimes \tau\otimes \id^{\otimes (n-k-1)})$, at first for $k< d$, to this equation. Notice that in this case $T_k\mu_n=0$ by the inductive assumption. Furthermore
\[(\id^{\otimes k}\otimes \tau\otimes \id^{\otimes (n-k-1)})(\id^{\otimes d}\otimes \tau\otimes \id^{\otimes (n-d)})=(\id^{\otimes (d-1)}\otimes \tau\otimes \id^{\otimes (n-d)})(\id^{\otimes k}\otimes \tau\otimes \id^{\otimes (n-k)})\]
and thus $T_kB'=0$ by the inductive assumption.

Similarly to the first part of the proof, if $k=0$, the only nonzero summand of $T_0A$ appears for $r=0$, which cancels with $T_0E$. For $1\leq k<d$, the only nonzero summands in $T_kA$ are those for $r=k-1,k$ which cancel, and $T_k$ annihilates $E$ again by the induction hypothesis. Thus $T_k\tilde{\mu}_n=0$ for $0\leq k <d$.

Lastly, we consider the case $k=d$. In this case, 
\begin{equation}\label{eqdd}(\id^{\otimes d}\otimes \tau\otimes \id^{\otimes (n-d-1)})(\id^{\otimes (d+1)}\otimes \tau\otimes \id^{\otimes (n-d-1)})=(\id^{\otimes d}\otimes \tau\otimes \id^{\otimes (n-d-1)})(\id^{\otimes d}\otimes \tau\otimes \id^{\otimes (n-d)}),\end{equation}
thus 
\begin{equation*}
\begin{split}
T_dE&=(-1)^{d+1}(\id^{\otimes d}\otimes \tau\otimes \id^{\otimes (n-d-1)})(\id^{\otimes (d+1)}\otimes \tau\otimes \id^{\otimes (n-d-1)})(\id\otimes\mu_n)\mu_2\\
&=(-1)^{d+1}(\id^{\otimes d}\otimes \tau\otimes \id^{\otimes (n-d-1)})(\id^{\otimes d}\otimes \tau\otimes \id^{\otimes (n-d)})(\id\otimes\mu_n)\mu_2\\
&=0,
\end{split}
\end{equation*}
again by the induction hypothesis.
Furthermore, 
\begin{equation*}
\begin{split}
T_dA &= \sum_{r=0}^{d-1} (-1)^{d-r} (\id^{\otimes d}\otimes \tau\otimes \id^{\otimes (n-d-1)})(\id^{\otimes (d+1)}\otimes \tau\otimes \id^{\otimes (n-d-1)})(\id^{\otimes r} \otimes \mu_2\otimes \id^{\otimes (n-1-r)} )\mu_n\\
&=\sum_{r=0}^{d-1} (-1)^{d-r} (\id^{\otimes d}\otimes \tau\otimes \id^{\otimes (n-d-1)})(\id^{\otimes d}\otimes \tau\otimes \id^{\otimes (n-d)})(\id^{\otimes r} \otimes \mu_2\otimes \id^{\otimes (n-1-r)} )\mu_n\\
\end{split}
\end{equation*}
using \eqref{eqdd}, which equals the sum of the first $d-1$ terms of $T_dB'$.
Therefore,
\begin{equation*}
\begin{split}
T_d\tilde{\mu}_n &= (\id^{\otimes d}\otimes \tau\otimes \id^{\otimes (n-d-1)})\mu_n 
+ \underbrace{(\id^{\otimes d}\otimes \tau\otimes \id^{\otimes (n-d-1)})(\id^{\otimes d}\otimes \tau\otimes \id^{\otimes (n-d)})(\id^{\otimes d}\otimes \mu_2\otimes \id^{\otimes (n-d-1)})\mu_{n}}_{G}\\&\qquad \,
+\underbrace{2\sum_{r=0}^{d-1} (-1)^{d-r}(\id^{\otimes d}\otimes \tau\otimes \id^{\otimes (n-d-1)}) (\id^{\otimes d}\otimes \tau\otimes \id^{\otimes (n-d)})(\id^{\otimes r} \otimes \mu_2\otimes \id^{\otimes (n-1-r)} )\mu_n}_{H}
\end{split}
\end{equation*}
where the second summand is the $r=d$ term of $T_dB'$, and the third term is the sum of the $r\leq d-1$ terms of $T_dA$ and $T_dB'$.
For $0\leq r<d-1$, we obtain 
\[
(\id^{\otimes d}\otimes \tau\otimes \id^{\otimes (n-d)})(\id^{\otimes r} \otimes \mu_2\otimes \id^{\otimes (n-1-r)} )\mu_n = (\id^{\otimes r} \otimes \mu_2\otimes \id^{\otimes (n-2-r)} )(\id^{\otimes (d-1)}\otimes \tau\otimes \id^{\otimes (n-d)})\mu_n =0
\]
by the induction hypothesis, thus the summands indexed by $r=0,\dots, d-2$ in $H$ are zero. On the other hand
\[(\id^{\otimes d}\otimes \tau\otimes \id^{\otimes (n-d)})(\id^{\otimes d}\otimes \mu_2\otimes \id^{\otimes (n-d-1)})= \id^{\otimes n} = (\id^{\otimes d}\otimes \tau\otimes \id^{\otimes (n-d)})(\id^{\otimes (d-1)}\otimes \mu_2\otimes \id^{\otimes (n-d)}),\]
hence 
\[G= (\id^{\otimes d}\otimes \tau\otimes \id^{\otimes (n-d-1)})\mu_{n}, \qquad H=-2(\id^{\otimes d}\otimes \tau\otimes \id^{\otimes (n-d-1)})\mu_n\]
and
\begin{equation*}
\begin{split}
T_d\tilde{\mu}_n &= (\id^{\otimes d}\otimes \tau\otimes \id^{\otimes (n-d-1)})\mu_n + (\id^{\otimes d}\otimes \tau\otimes \id^{\otimes (n-d-1)})\mu_{n}\\&\qquad \,
-2(\id^{\otimes d}\otimes \tau\otimes \id^{\otimes (n-d-1)})\mu_n\\
&=0
\end{split}
\end{equation*}
as desired. 
\end{proof}

\subsection{Pretriangulated $A_\infty$-categories and twisted complexes}\label{sec:twistedcomplexes}

In this subsection, we recall the basics of twisted modules, twisted complexes and pretriangulated $A_\infty$-categories and state Keller--Lef\`evre-Hasegawa's theorem about reconstruction of the extension closure, resp. triangulated closure, of modules from the $A_\infty$-structure on the $\Ext$-algebra of the modules. For further details, we refer the reader to \cite{L-H03, K02, BLM08}. 

\begin{defn}
Let $\mathcal{A}$ be an $A_\infty$-category. 
\begin{enumerate}[(i)]
\item The category of \emphbf{twisted modules} over $\mathcal{A}$, denoted $\twmod \mathcal{A}$ is the $A_\infty$-category with objects $(B,\alpha)$ where $B$ is a sequence $(\mathtt{i}_1,\dots,\mathtt{i}_k)$ of objects in $\mathcal{A}$ and $\alpha=(\alpha_{ij})_{i,j}$ is a strictly upper triangular matrix (i.e. $\alpha_{ij}=0$ for $i\geq j$) of size $k\times k$ with $\alpha_{ij}\in \mathcal{A}(\mathtt{i}_j,\mathtt{i}_i)^1$ such that the Maurer--Cartan equation
\[\sum_{t=1}^\infty (-1)^{\frac{t(t-1)}{2}}m_t(\alpha,\dots,\alpha)=0\]
is satisfied. The morphisms are defined by 
\[(\twmod\mathcal{A})((B,\alpha), (B',\alpha'))=\bigoplus_{i,j}\mathcal{A}(\mathtt{i}_i,\mathtt{i}_j)\]
and 
\[m_n^{\twmod \mathcal{A}}=\sum_{t=0}^\infty \sum \pm m_{n+t}^{\mathcal{A}}(\id^{\otimes i_1}\otimes \alpha^{\otimes j_1}\otimes \id^{\otimes i_2}\otimes \dots\otimes \alpha^{\otimes j_{r-1}}\otimes \id^{\otimes i_r}\otimes \alpha^{\otimes j_r})\]
where the terms of the sum are in bijection with the non-commutative monomials $X^{i_1}Y^{j_1}\otimes \dots\otimes X^{i_r}Y^{j_r}$ of degree $n$ in $X$ and degree $t$ in $Y$ and the sign is given by
\[(\shift X)^{i_1}(\shift Y)^{j_1}\cdots (\shift X)^{i_r}(\shift Y)^{j_r}=\pm \shift^{n+t}X^{i_1}Y^{j_1}\otimes \dots\otimes X^{i_r}Y^{j_r}\]
in the algebra $\Bbbk\langle X,Y, \shift\rangle (\shift X-X\shift, \shift Y+Y\shift)$. 
\item If $\mathcal{A}$ is an $A_\infty$-category, then $\mathbb{Z}\mathcal{A}$ is the $A_\infty$-category with objects $(\mathtt{i},n)$ where $\mathtt{i}\in \mathcal{A}$ and $n\in \mathbb{Z}$ and morphisms $\mathbb{Z}\mathcal{A}((\mathtt{i},n),(\mathtt{i}',n'))=\shift^{n'-n}\mathcal{A}(\mathtt{i},\mathtt{i}')$
and the natural compositions.
\item The $A_\infty$-category of \emphbf{twisted complexes over $\mathcal{A}$} is the $A_\infty$-category $\twcom \mathcal{A}=\twmod \mathbb{Z}\mathcal{A}$. 
\end{enumerate}
\end{defn}

\begin{rmk}
The $A_\infty$-category $\twcom \mathcal{A}$ is triangulated in a canonical way. The shift functor is given on objects by sending $(B,\alpha)$ to $(B,\alpha)[1]:=(B[1],-\alpha)$ where for $B=((\mathtt{i}_1,m_1),\dots,(\mathtt{i}_k,m_k))$ we set $B[1]:=((\mathtt{i}_1,m_1-1),\dots,(\mathtt{i}_k,m_k-1))$. Also, for a morphism $f\colon (B,\alpha)\to (C,\beta)$, its cone is defined to be the object $((B[1],C),\begin{pmatrix}-\alpha&0\\f&\beta\end{pmatrix})$. 
\end{rmk}

\begin{thm}[{\cite[Corollaire 7.6.0.7]{L-H03}, see also \cite[p. 7]{K02},  \cite[Theorem B.3.1]{Mad02}}] 
Let $\Lambda$ be an algebra. Let $M_1,\dots,M_r$ be finite dimensional $\Lambda$-modules. Let $\mathcal{E}$ be $\Ext^{\bullet}_\Lambda(M,M)$ regarded as an $A_\infty$-algebra where $M=\bigoplus M_i$. Then the triangulated hull of $M_1,\dots,M_r$ in $D^b(\modu A)$ is equivalent as a triangulated category to $H^0(\twcom \mathcal{E})$ while the category $\mathcal{F}(M_1,\dots,M_r)$ of $M$-filtered modules is equivalent to $H^0(\twmod \mathcal{E})$.
\end{thm}

\section{Quasi-hereditary algebras and bocses}\label{sec:qhalg_bocses}

In this section, we recall important definitions and background concerning quasi-hereditary algebras and bocses (Section \ref{sec:bocses}). We analyse various descriptions of extensions of modules over a bocs (Section \ref{sec:extensions}), give a specific $A$-bimodule resolution of the coring $V$ involved in the definitions of a bocs, and show that the comultiplication on $V$ can be lifted to an $A_\infty$-coalgebra structure on this bimodule resolution (Section \ref{sec:projectiveresolution}). We finally give a simple proof of uniqueness of exact Borel subalgebras up to isomorphism (Section \ref{sec:uniquenessofexactborelsubalgebras}).

\subsection{Definitions and recollections}\label{sec:bocses}
The main object of study in this article is the class of quasi-hereditary algebras defined by Scott \cite{Sco87}, see also \cite{CPS88}. For further reading we refer to the survey articles \cite{DR92} and \cite{KK99}.

\begin{defn}
A finite dimensional algebra $\Lambda$ with $n$ simples up to isomorphism is called \emphbf{quasi-hereditary} if there exist $\Lambda$-modules $\Delta(\mathtt{1}),\dots, \Delta(\mathtt{n})$ with $\End_\Lambda(\Delta(\mathtt{i}))\cong \Bbbk$ for all $\mathtt{i}$ and $\Ext^s_\Lambda(\Delta(\mathtt{i}),\Delta(\mathtt{j}))=0$ for $\mathtt{i}>\mathtt{j}$ and $s=0,1$ and such that 
\[A\in \mathcal{F}(\Delta)=\{X\in \modu \Lambda\,|\,\exists 0=X_0\subset X_1\subset \dots\subset X_t=X\text{ with } X_j/X_{j-1}\cong \Delta(\mathtt{i}_j)\}.\]
\end{defn}

\begin{rmk}
It follows from dimension shifting that $\Ext_\Lambda^s(\Delta(\mathtt{i}),\Delta(\mathtt{j}))=0$ for all $\mathtt{i}>\mathtt{j}$ and all $s\in \mathbb{N}_0$, i.e. the $\Delta(\mathtt{i})$ form an exceptional collection in $D^b(\modu \Lambda)$. 
\end{rmk}

Important examples of such algebras are Schur algebras, blocks of BGG category $\mathcal{O}$ for a semisimple complex Lie algebra $\mathfrak{g}$ and algebras of global dimension smaller or equal to two. Motivated by the example of blocks of BGG category $\mathcal{O}$, Koenig \cite{Koe95} introduced the notion of an exact Borel subalgebra for a quasi-hereditary algebra.

\begin{defn}
Let $\Lambda$ be a quasi-hereditary algebra with $n$ simples. A subalgebra $A\subseteq \Lambda$ with $n$ simples is called an \emphbf{exact Borel subalgebra} if
\begin{enumerate}[{(B}1{)}]
\item $A$ is \emphbf{directed}, i.e. whenever there is an arrow $\mathtt{i}\to \mathtt{j}$ in the Gabriel quiver of $A$, then $\mathtt{i}\leq \mathtt{j}$.
\item $\Lambda\otimes_A -$ is exact, i.e. $\Lambda$ is projective as a right $A$-module.
\item $\Lambda\otimes_A L_A(\mathtt{i})\cong \Delta_\Lambda(\mathtt{i})$ as $\Lambda$-modules.
\end{enumerate}
\end{defn}

In joint work with Koenig and Ovsienko \cite{KKO14}, the first author proved that exact Borel subalgebras always exist up to Morita equivalence and that moreover the class of quasi-hereditary algebras can be described using bocses, i.e. \underline{b}imodules \underline{o}ver \underline{c}ategories with \underline{c}oalgebra structures. Bocses were introduced by Ro\u\i ter in \cite{Roi79} and most prominently used in the proof of tame-wild dichotomy by Drozd in \cite{D80}. We begin by recalling the definitions and the main result needed from this paper. For further reading on bocses and corings we refer the reader to \cite{BSZ09, BB91, Bur05, BW03, Kul17}. 

\begin{defn}\label{directedbocs}
  \begin{enumerate}[(i)]
  \item A \emphbf{bocs} $\mathfrak{A}=(A,V)$ is a pair consisting of a finite dimensional $\Bbbk$-algebra $A$ and an $A$-\emphbf{coring} $V$, i.e. an $A$-$A$-bimodule $V$ with morphisms of $A$-$A$-bimodules $\mu\colon V\to V\otimes_A V$ and $\varepsilon\colon V\to A$ satisfying the usual axioms for $\mu$ to be coassociative and $\varepsilon$ to be a counit.  In other words, $(\mu\otimes 1)\mu=(1\otimes \mu)\mu$ and $(\varepsilon\otimes 1)\mu$ and $(1\otimes \varepsilon)\mu$ are both the respective canonical isomorphisms. In this paper we will additionally assume that the algebra $A$ is basic.
  \item \label{morofbocsmod} Given a bocs $\mathfrak{A}=(A,V)$ the category of finite dimensional $\mathfrak{A}$-modules $\modu\mathfrak{A}$ has as objects all finite dimensional $A$-modules and as morphisms
    \[\Hom_{\mathfrak{A}}(M,N)=\Hom_{A\otimes A^{\op}}(V,\Hom_{\Bbbk}(M,N))\]
    with composition given by $V\to V\otimes_A V\to \Hom_{\Bbbk}(M,N)\otimes_A \Hom_{\Bbbk}(L,M)\to \Hom_{\Bbbk}(L,N)$.
  \item Given a bocs $\mathfrak{A}=(A,V)$ its \emphbf{right algebra} is defined to be $\End_{\mathfrak{A}}(A)^{\op}$. Dually its \emphbf{left algebra} is defined to be $\End_{\mathfrak{A}^{\op}}(A)$ where $\mathfrak{A}^{\op}=(A^{\op},V,\mu,\varepsilon)$.
\item A bocs $\mathfrak{A}=(A,V)$ is said to have \emphbf{projective kernel} if $\varepsilon$ is surjective and $\overline{V}:=\ker\varepsilon$ is a projective $A$-$A$-bimodule.
\item A bocs $\mathfrak{A}=(A,V)$ is said to be \emphbf{normal} if there is an element $\omega\in V$ with $\mu(\omega)=\omega\otimes_A \omega$ and $\varepsilon(\omega)=1$. In this case, $\omega$ is called \emphbf{grouplike}.
\item Given a bocs with projective kernel, its \emphbf{biquiver} is defined to be the quiver with two types of arrows, solid and dashed: The solid arrows are given by the quiver of $A$, the number of dashed arrows is given by the number of summands isomorphic to $Ae_{\mathtt{j}}\otimes_\Bbbk e_{\mathtt{i}}A$ in $\overline{V}$.
  \item A bocs with projective kernel is called \emphbf{directed} if its biquiver is acyclic.
  \end{enumerate}
\end{defn}

\begin{rmk}\label{equivalentdefinitions}
The standard tensor-hom adjunctions 
\[\Hom_{\mathfrak{A}}(M,N)=\Hom_{A\otimes A^{\op}}(V,\Hom_{\Bbbk}(M,N))\cong \Hom_A(V\otimes_A M,N)\cong \Hom_A(M,\Hom_A(V,N)) ,\]
yield different possible, but equivalent, definitions of the category of representations of a bocs, see \cite[Section 7.1]{KKO14}, \cite{BB91}.
\end{rmk}

\begin{rmk}\label{compflat}
For $\mathbb{L}$ the direct sum of simple $A$-modules, consider the chain of isomorphisms
\begin{equation*}\begin{split}
e_{\mathtt{j}}\End_{\mathfrak{A}}(\mathbb{L}) e_{\mathtt{i}}&= \Hom_{A\otimes A^{\op}}(V,\Hom_{\Bbbk}(\mathbb{L}e_{\mathtt{i}},\mathbb{L}e_{\mathtt{j}}))
\cong \Hom_{A}(V\otimes_A \mathbb{L}e_{\mathtt{i}},\mathbb{L}e_{\mathtt{j}}) \\
&\cong\Hom_{\mathbb{L}}(  \mathbb{L}e_{\mathtt{j}}\otimes_AV\otimes_A\mathbb{L}e_{\mathtt{i}},\mathbb{L})
 = e_{\mathtt{j}} (V^\mathbb{L})^\flat e_{\mathtt{i}},
\end{split}\end{equation*}
where we set $M^\mathbb{L}=\mathbb{L}\otimes_AM\otimes_A \mathbb{L}$ for any $A$-$A$-bimodule $M$. Similarly, setting $f^{\mathbb{L}}= \mathbb{L}\otimes_Af \otimes_A \mathbb{L}$ for any morphism of $A$-$A$-bimodules, this 
translates composition into tensor product over $\mathbb{L}$. Indeed,  composition of $\phi$ and $\psi$ as described in Definition \ref{directedbocs}\eqref{morofbocsmod}, when pushed through the isomorphisms yields
\[V^\mathbb{L} \overset{\mu^\mathbb{L}}{\to} (V\otimes_A V)^\mathbb{L}  \overset{(\id \otimes_A\phi)^\mathbb{L} }{\longrightarrow} V^\mathbb{L} \overset{\psi^\mathbb{L} }{\to} \mathbb{L},\]
which coincides with the tensor product in Lemma \ref{flathash}, keeping in mind that there are no gradings here. Using Lemma \ref{flathash}, we can also further identify $\End_{\mathfrak{A}}(\mathbb{L}) $ with $(V^\mathbb{L})^\#$ and composition again translates to tensor product.
\end{rmk}

\begin{defn}
Let $\mathfrak{A}=(A,V)$ and $\mathfrak{B}=(B,W)$ be bocses. A \emphbf{morphism} $\mathfrak{A}\to\mathfrak{B}$ of bocses is a pair of maps $(f_0, f_1)$ with $f_0\colon A\to B$ and $f_1\colon V\to W$ such that 
\begin{itemize}
\item $f_0$ is an algebra morphism;
\item $f_1$ is a morphism of $A$-$A$-bimodules, where $W$ is given the structure of an $A$-module via $f_0$;
\item $\varepsilon_W \circ f_1 = f_0\circ \varepsilon_V$;
\item $T(f_0, f_1)\circ \mu_V = \mu_W \circ f_1$, where $T(f_0, f_1)(v_1\otimes v_2)= f_1(v_1)\otimes f_1(v_2)$, where $f_0$ is included in the notation as the tensor product changes from being over $A$ to being over $B$ using $f_0$. 
\end{itemize}
\end{defn}

\begin{thm}[{\cite[Theorem 11.3]{KKO14}}]\label{KKOmain}
  An algebra $\Lambda$ is quasi-hereditary if and only if it is Morita equivalent to the right algebra $R$ of a normal directed bocs $\mathfrak{A}=(A,V)$. In this case, $A$ is an exact Borel subalgebra of $R$ and $\modu\mathfrak{A}$ is equivalent to the category $\mathcal{F}(\Delta)$, where $\Delta(\mathtt{i}):=R\otimes_A L_A(\mathtt{i})$ for the simple $A$-module $L_A(\mathtt{i})$ corresponding to the vertex $\mathtt{i}$.
\end{thm}

We briefly recall the construction (in fact the construction given here differs by the original construction by a `global' minus-sign) and hereby show that it is in fact functorial. The construction starts with the $A_\infty$-category structure on the $\Ext$-algebra $\mathcal{E}$ of the standard modules (its restriction to the degrees less than or equal to two suffices). Firstly, one takes the linear dual $\mathcal{E}^{\#}$ as described in Section \ref{sec:duality} to obtain an $A_\infty$-cocategory. To this one applies the cobar construction to obtain the differential graded category $\mathbf{\Omega}(\mathcal{E}^{\#})$. These two steps are functorial as shown in Section \ref{sec:duality}. The remaining two steps, which we describe next, are truncation by the negative degree part and an equivalence between semifree differential graded algebras and bocses by Ro\u\i ter. 

Let $\mathcal{D}=\bigoplus_{i\in \mathbb{Z}} \mathcal{D}_i$ be a differential graded algebra. Define the non-negative truncation $\sigma_{\geq 0}\mathcal{D}$ of $\mathcal{D}$ to be the quotient of $\mathcal{D}$ by the differential ideal generated by $\bigoplus_{i<0} \mathcal{D}_i$, i.e. the ideal generated by $\bigoplus_{i<0}\mathcal{D}_i$ and $d(\mathcal{D}_{-1})$. This is in fact a functor as a map $f\colon \mathcal{D}\to \mathcal{D}'$ is degree-preserving and thus sends the differential ideal spanned by $\bigoplus_{i<0} \mathcal{D}_i$ to the differential ideal spanned by $\bigoplus_{i<0} \mathcal{D}_{i<0}'$. 

\begin{ex}
As an illustration we give the following example from \cite[Section 6.2]{KK18}. We have shifted their numbering by one in order to not have a vertex called $0$. The quasi-hereditary algebra $\Lambda$ is given by the path algebra of the quiver
\[\begin{tikzcd}
\mathtt{3}\arrow[bend right, yshift=-0.5ex]{r}{a_2}&\mathtt{2}\arrow[bend right]{l}{c_2}\arrow[bend right, yshift=-0.5ex]{r}{a_1}\arrow[bend right, yshift=-1ex]{r}[swap]{k_2}&\mathtt{1}\arrow[bend right, yshift=0.5ex]{l}{c_1}
\end{tikzcd}\]
with relations $c_2a_2=0$, $a_2c_2=c_1k_2$, and $c_1a_1=0$. One can check that in this case $\Ext^{\bullet}_{\Lambda}(\Delta,\Delta)$ is given by the following quiver where the solid arrows correspond to basis vectors of $\Ext^1_{\Lambda}(\Delta,\Delta)$ while the dashed arrows represent basis vectors of $\Hom_{\Lambda}(\Delta,\Delta)$. 
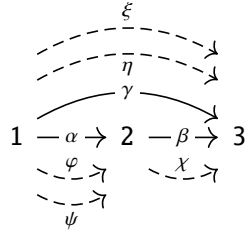
\[\begin{tikzcd}
\mathtt{1}\arrow{r}[description]{\alpha}\arrow[yshift=-1.5ex, dashed, bend right]{r}{\varphi}\arrow[yshift=-3.5ex, dashed, bend right]{r}[swap]{\psi}\arrow[bend left, dashed, yshift=3.5ex]{rr}[swap]{\eta}\arrow[bend left, dashed, yshift=5.5ex]{rr}{\xi}\arrow[yshift=0.5ex, bend left]{rr}[description]{\gamma}&\mathtt{2}\arrow{r}[description]{\beta}\arrow[dashed, bend right, yshift=-1.5ex]{r}{\chi}&\mathtt{3}
\end{tikzcd}
\]
with $m_2(\chi,\varphi)=\eta$, $m_2(\chi,\psi)=\xi$, $m_2(\chi,\alpha)=\gamma$, $m_2(\beta,\varphi)=\gamma$, $m_2(\beta,\psi)=0$, $m_2(\beta,\alpha)=0$. Since there are no paths of length greater than or equal to two, it follows that all the $m_n$ for $n\geq 2$ vanish.   
\end{ex}

Furthermore, we recall a theorem by Ro\u\i ter, \cite[Proposition 2]{Roi79}, which was written in a functorial way by Brzezi\'nski \cite{Brz13}. This associates to a semifree differential graded algebra a bocs. Recall that a differential graded algebra $\mathcal{D}$ is called \emphbf{semifree} if $\mathcal{D}\cong T_A(\overline{V})$ as graded algebras for an algebra $A$ and an $A$-bimodule $\overline{V}$.  

\begin{thm}\label{Brzmain}
There is an equivalence of categories between the category of semifree differential graded algebras over $\mathbb{L}$ and normal bocses. 

In one direction, it associates to a normal bocs $(A,V)$ with group-like $\omega$ the tensor algebra $T_A(\shift^{-1}\overline{V})$ with differential $\partial(a)=\shift^{-1}(\omega a-a\omega)$ in degree $0$ and $\partial(\shift^{-1}v)=(\shift\otimes\shift)^{-1}(\mu(v)-\omega\otimes v-v\otimes \omega)$ in degree $1$. For a morphism of bocses $(f_0,f_1)$, there is a morphism of semifree differential graded categories given in degree $0$ by $f_0$ and in degree $1$ by $\shift^{-1}f_1\shift$. 

In the reverse direction, it associates to a semifree differential graded algebra $(T_A(\shift^{-1}\overline{V}),\partial)$ the bocs with $V=A\omega\oplus \overline{V}$ as a left module and with right module structure given by $(a\omega+v)a'=aa'\omega+a\shift\partial(a')+va'$, comultiplication given by $\mu(a\omega+v)=a\omega\otimes \omega+\omega\otimes v+v\otimes \omega+(\shift\otimes\shift)\partial(\shift^{-1}v)$ and counit given by $\varepsilon(a\omega+v)=a$. For a morphism of differential graded algebras $f\colon (T_A(\shift^{-1}\overline{V}),\partial)\to (T_B(\shift^{-1}\overline{W}),\partial')$ a morphism of bocses $(f_0,f_1)$ is defined by $f_0=f|_A$ and $f_1(a\omega+v)=f_0(a)\omega+\shift f|_{\shift^{-1}\overline{V}}(\shift^{-1}v)$.
\end{thm}

We remark the following corollary, which is immediate from functoriality of the above constructions.

\begin{cor}\label{extisobocsiso}
Let $f\colon \mathcal{E}\to \mathcal{E}'$ be an isomorphism of minimal $A_\infty$-algebras. Then the corresponding bocses are isomorphic.
\end{cor}

The exact Borel subalgebras arising from directed bocses are quite special as shown more generally by Burt and Butler \cite{BB91}.

\begin{thm}[{\cite[Extension closure theorem (3.8)]{BB91}}]\label{extnisomorphic}
Let $(A,V)$ be a bocs with projective kernel with right algebra $R$. Then the functor $R\otimes_A-\colon \modu A\to \modu R$ induces $\Bbbk$-linear maps $\Ext^n_A(Z,X)\to \Ext^n(R\otimes_A Z,R\otimes_A X)$ for all $A$-modules $Z,X$. For $n=1$, this is an epimorphism and for $n\geq 2$, these maps are isomorphisms. Furthermore, the functor $R\otimes_A -\colon \modu A\to \modu R$ is faithful and dense.
\end{thm}

There is in general no isomorphism between $\Ext^1_A(Z,X)$ and $\Ext^1_R(R\otimes_A Z,R\otimes_A X)$. However, in the case where $Z$ and $X$ are simple modules we show in Lemma \ref{regularbocsesequivalence} that it can always be assumed by a process called regularisation. As noted in \cite[Appendix A.2]{KKO14}, this is an obstruction to uniqueness of bocses giving rise to the same quasi-hereditary algebra.

\begin{prop}\label{regularisation}
  Let $\mathfrak{A}=(A,V)$ be a bocs with $\partial(a)=\lambda\psi+\sum_{i}c_i\psi_ib_i$ where $a,b_i,c_i\in A$, $\psi\neq \psi_i$ are generators of $\overline{V}$ and $0\neq \lambda\in \Bbbk$. Then there is a bocs $\tilde{\mathfrak{B}}=(\tilde{B}, \tilde{W})$ with $\tilde{A}=A/(a)$ and $\tilde{V}=\tilde{A}\otimes_A V/(\psi)\otimes_A \tilde{A}$ such that the following statements hold:
    \begin{enumerate}[(i)]
    \item There is an equivalence of categories $\modu \mathfrak{A}\cong \modu \tilde{\mathfrak{A}}$.
    \item If $\mathfrak{A}=(A,V)$ is directed (resp. normal), then $\tilde{\mathfrak{A}}$ is directed (resp. normal).
    \item If $\mathfrak{A}$ is directed, then the right algebra of $\mathfrak{A}$ is Morita equivalent to the right algebra of $\tilde{\mathfrak{A}}$.
    \end{enumerate}
\end{prop}

\begin{proof}
The first claim is already contained in \cite[Proposition 2]{KR77}. It is immediate that if $\mathfrak{A}$ has a projective kernel, then so does $\tilde{\mathfrak{A}}$. Since the biquiver of $\tilde{\mathfrak{A}}$ is obtained from that of  $\mathfrak{A}$ by removing one solid and one dashed edge, directedness of $\tilde{\mathfrak{A}}$  follows from directedness of $\mathfrak{A}$. It is straightforward to check that the image of the group-like element $\omega$ is still a group-like in $\tilde{V}$. The last claim follows from the Dlab--Ringel standardisation theorem \cite[Theorem 2]{DR92} as both right algebras are quasi-hereditary and have equivalent categories of filtered modules $\modu \mathfrak{A}$ and $\modu \tilde{\mathfrak{A}}$.
  \end{proof}

Proceeding by induction on the total number of arrows in the biquiver of $\mathfrak{B}$, it is clear that for every directed bocs there is another directed bocs with equivalent module category which cannot be regularised further. This leads to the following notion introduced by Kleiner and Ro\u\i ter in \cite{KR77} (in the equivalent language of differential tensor algebras). 

\begin{defn}\label{regular}
A bocs $\mathfrak{A}=(A,V)$ is called \emphbf{regular} if $\partial(J)\subseteq \radoperator_{A\otimes A^{\op}} V$, where $J\subseteq A$ denotes the Jacobson radical of $A$. 
\end{defn}

From our point of view, the following definition seems reasonable since we will prove that (at least for directed bocses) it yields a unique, in some sense minimal, representative of each `Morita equivalence class' of bocses:

\begin{defn}\label{basic}
A bocs $\mathfrak{A}=(A,V)$ is called \emphbf{basic} if it is normal, regular and $A$ is basic.
\end{defn}

\subsection{Extensions between modules over a bocs}\label{sec:extensions}

In this subsection, we show that the possible definitions of $\Ext$ in the category of modules for a directed bocs agree. As was shown in \cite{KKO14}, the category of modules over a directed bocs is exact, hence one possible definition of $\Ext^j$ is of course the definition as equivalence classes of exact sequences of length $j+2$. Three others derive from the description of $\Hom$ in the category of modules over a bocs as 
\[\Hom_A(V\otimes_A M,N)\cong \Hom_A(M,\Hom_A(V,N))\cong \Hom_{A\otimes A^{\op}}(V,\Hom_{\Bbbk}(M,N))\]
where $M$ and $N$ are modules over the bocs $\mathfrak{A}=(A,V)$, see Remark \ref{equivalentdefinitions}. The last one makes it possible to describe $\Ext_{\mathfrak{A}}^1(M,N)$ as a quotient of $\Ext_{A}^1(M,N)$ by applying $\Hom_{A\otimes A^{\op}}(-,\Hom_{\Bbbk}(M,N))$ to the exact sequence $0\to \overline{V}\to V\overset{\varepsilon}{\to} A\to 0$. For $j\geq 2$, Burt and Butler \cite{BB91} showed that $\Ext^j_\mathfrak{A}(M,N)\cong \Ext^j_A(M,N)$, c.f. Theorem \ref{extnisomorphic}. 

\begin{lem}
Let $A$ be an algebra. Let $V$ be an $A$-$A$-bimodule which is projective as a left and as a right $A$-module. Let $M$, $N$ be $A$-modules. Then
\[\Ext^j_A(V\otimes_A M, N)\cong \Ext^j_A(M,\Hom_A(V,N)).\]
\end{lem}

\begin{proof}
Let $P^{\bullet}$ be a projective resolution of $M$ as an $A$-module. Since $V$ is projective as a left and as a right module, $V\otimes_A P^{\bullet}$ is a projective resolution of $V\otimes_A M$. Thus, $\Ext^j_A(V\otimes_A M,N)$ can be computed as the homology of the complex $\Hom_A(V\otimes_A P^{\bullet},N)$. But this complex is naturally isomorphic to the complex $\Hom_A(P^{\bullet},\Hom_A(V,N))$ whose homology is $\Ext^j_A(M,\Hom_A(V,N))$. 
\end{proof}

This proves the equivalence of two of the possible definitions of $\Ext$. The third one is captured in the following proposition.

\begin{prop}\label{extbimoduledescription}
Let $V$ be an $A$-$A$-bimodule which is projective as a left and as a right $A$-module. Let $M$ and $N$ be $A$-modules. Then
\[\Ext^j_{A\otimes A^{\op}}(V,\Hom_{\Bbbk}(M,N))\cong \Ext^j_A(V\otimes_A M,N).\] 
\end{prop}

\begin{proof}
Let $P^{\bullet}\to V$ be a projective resolution of $V$ as an $A$-$A$-bimodule. As $V$ is projective as a right module, the resolution is composed of split exact sequences as right $A$-modules. Thus, it is still exact when tensored with $M$ over $A$. Also, a module of the form $P\otimes_A M$ for $P$ a projective bimodule is projective as a left $A$-module. Thus, $P^{\bullet}\otimes M$ is a projective resolution of $V\otimes_A M$ as a left $A$-module. We conclude that $\Ext^j_{A\otimes A^{\op}}(V,\Hom_{\Bbbk}(M,N))$, which is by definition the homology of the complex $\Hom_{A\otimes A^{\op}}(P^{\bullet},\Hom_{\Bbbk}(M,N))\cong \Hom_A(P^{\bullet}\otimes_A M,N)$, is isomorphic to $\Ext^j_A(V\otimes M,N)$.    
\end{proof}

Applying $\Hom_{A\otimes A^{\op}}(-,\Hom_{\Bbbk}(M,N))$ to the exact sequence $0\to \overline{V}\to V\overset{\varepsilon}{\to} A\to 0$ yields the long exact sequence

\hspace*{-7cm}\begin{tikzcd}[row sep=0.5ex]
0\arrow{r}& \Hom_A(M,N)\arrow{r}&\Hom_{\mathfrak{A}}(M,N)\arrow{r} &\Hom_{A\otimes A^{\op}}(\overline{V},\Hom_{\Bbbk}(M,N))\arrow[out=0,in=180, looseness=2]{dll}\\
&\Ext^1_A(M,N)\arrow{r}&\Ext^1_{\mathfrak{A}}(M,N)\arrow{r}&\Ext^1_{A\otimes A^{\op}}(\overline{V},\Hom_{\Bbbk}(M,N))\to \dots
\end{tikzcd}

but as $\overline{V}$ is a projective bimodule, $\Ext^1_{A\otimes A^{\op}}(\overline{V},\Hom_{\Bbbk}(M,N))=0$. Thus, $\Ext^1_{\mathfrak{A}}(M,N)$ can also be defined as a quotient of $\Ext^1_A(M,N)$ by the image of $\Hom_{A\otimes A^{\op}}(\overline{V},\Hom_{\Bbbk}(M,N))$ under the map which can, via the standard adjunction, be described as sending a homomorphism $f\colon  \overline{V}\otimes_A M\to N$ to the pushout sequence of the  diagram

\begin{equation}
\begin{tikzcd}
0\arrow{r}&\overline{V}\otimes_A M\arrow{d}{f}\arrow{r}{i} &V\otimes_A M\arrow{d}\arrow{r}& A\otimes_A M\arrow[equals]{d}\arrow{r}&0\\
0\arrow{r}&N\arrow{r}&E\arrow{r}{p}&A\otimes_A M\arrow{r}&0
\end{tikzcd}\label{eqn:pushout}
\end{equation}

We now show that the derived functor definitions 
\[\Ext^j_{A\otimes A^{\op}}(V,\Hom_{\Bbbk}(M,N))\cong \Ext^j_A(V\otimes_A M,N)\cong \Ext^j_A(M,\Hom_A(V,N))\]
are equivalent to the definition in terms of Yoneda-$\Ext$. The proof of the following proposition is taken from unpublished notes of Ovsienko, cf. \cite{Ovs93}.

\begin{prop}
Let $\mathfrak{A}=(A,V)$ be a directed bocs, then $\Ext^j_\mathfrak{A}(M,N)\cong \Ext^j_A(V\otimes_A M,N)$ can equivalently be described as equivalence classes of exact sequences of length $j+2$ in $\modu \mathfrak{A}$ by the standard equivalence relation.
\end{prop}

\begin{proof}
As $\modu \mathfrak{A}$ has enough projectives, by dimension shifting it suffices to prove the statement for $j=1$. By additivity, it also suffices to prove that a sequence $\eta\in \Ext^1_A(M,N)$ belongs to the image of the map from $\Hom_{A\otimes A^{\op}}(\overline{V},\Hom_{\Bbbk}(M,N))$ if and only if it is equivalent to a split exact sequence. 

To prove that an element of the image is equivalent to a split exact sequence, consider an element $\eta$ in the image. Such an element is obtained as a pushout as in \eqref{eqn:pushout}. The middle vertical homomorphism $V\otimes_A M\to E$ corresponds to a morphism $s\in \Hom_{\mathfrak{A}}(M,E)$. The map $p \colon E\to A\otimes_A M$ is $A$-linear, hence can also be regarded as a morphism in $\modu \mathfrak{A}$. Their composition $p\circ s$ is the map $\varepsilon\otimes_A M$, which is the identity in $\modu \mathfrak{A}$, so the lower row in \eqref{eqn:pushout} is indeed split in $\modu \mathfrak{A}$. 

For the reverse direction, let $s\colon V\otimes_A M\to E$ correspond to a splitting of the map $p\colon E\to M$ in $\modu \mathfrak{A}$. This map can be fitted to obtain the right-most square in \eqref{eqn:pushout}. To finish, it suffices to give a map $f$ such that the diagram \eqref{eqn:pushout} commutes. Denote by $i\colon \overline{V}\otimes_AM\to V\otimes_AM$ the inclusion morphism. Then $p\circ s\circ i=0$ in $\modu A$. Thus, there exists $f\colon \overline{V}\otimes_A M\to N$ making the diagram commutative.
\end{proof}

\subsection{$A_\infty$-coalgebra structures on projective resolutions}
\label{sec:projectiveresolution}

Recall that Kadeishvili's theorem, Theorem \ref{kadeishvilitheorem}, tells us that an $A_\infty$-structure on a vector space induces an $A_\infty$-structure on homology. If one considers projective resolutions, one can go in the other direction and produce an $A_\infty$-coalgebra structure on a projective resolution from a coalgebra structure on its homology. In the analogous case of $A_\infty$-algebras over commutative rings (which are not necessarily fields) such a result was obtained by Jesse Burke in \cite[Theorem 3.1]{Bur18}. 
We will use this extensively in Section \ref{sec:uniquenessofbocses}. A special case can be found in \cite[Theorem 3.1.1]{NVW18}.

\begin{lem}\label{AinfinitycoalgebraobjectP}
Let $(A,V,\mu,\varepsilon)$ be a bocs. Let $(\mathcal{P}^i)_ {i\leq 0}$ be a projective resolution of $V$ as an $A$-$A$-bimodule. Then, on $\mathcal{P}=\bigoplus_{i\leq 0} \mathcal{P}^i$, there exists the structure of an $A_\infty$-coalgebra object in $\Mod A\otimes A^{\op}$ with the given differential $\mu_1$ such that the projection map $\mathcal{P}^0\to V$ provides a quasi-isomorphism of $A_\infty$-coalgebra objects $\mathcal{P}\to V$. 
\end{lem}

\begin{proof}
We define the maps $\mu_k^n\colon \mathcal{P}^n\to \bigoplus_{\sum i_j=n+2-k} \mathcal{P}^{i_1}\otimes_A \dots\otimes_A \mathcal{P}^{i_k}$ inductively via the universal property of projective modules. Let $\pi$ be the morphism $\mathcal{P}\to V$ defined by $\mathcal{P}^0\to V$ being the given projection and $\mathcal{P}^n\to V$ being zero for $n\geq 1$. Let $(\mu_2^n)_{n\leq 0}$ be a lift of $\mu\colon V\to V\otimes_A V$ to $\mathcal{P}^n$. 
Suppose $\mu_l$ has already been defined for $l<n$. We change the differential on $\mathcal{P}^{\otimes n}$ slightly to $(-1)^{n-1}\sum_{r+1+t=n} (\id^{\otimes r}\otimes \mu_1\otimes \id^{\otimes t})$. We will check that the diagram 
\[
\begin{tikzcd}
\dots\arrow{d}\arrow{r} &\mathcal{P}^{-2}\arrow{d}{\tilde{\mu}^2_n}\arrow{r}\arrow{d}&\mathcal{P}^{-1}\arrow{r}\arrow{d}{\tilde{\mu}^1_n} &\mathcal{P}^0\arrow{r}\arrow{d}{\tilde{\mu}^0_n}&V\arrow{d}{0}\arrow{r}&0\\
\dots\arrow{r} &(\mathcal{P}^{\otimes n})^{0-n}\arrow{r}&(\mathcal{P}^{\otimes n})^{1-n}\arrow{r}&(\mathcal{P}^{\otimes n})^{2-n}\arrow{r}&(\mathcal{P}^{\otimes n})^{3-n}\arrow{r}&\dots
\end{tikzcd}
\]
commutes, where $\tilde{\mu}_n=\sum\limits_{\substack{p+q+r=n\\q\neq n,1}}(-1)^{pq+r} (\id^{\otimes p}\otimes \mu_q\otimes \id^{\otimes r})\mu_{p+1+r}$. We check that $\tilde{\mu}_n$ is a morphism of complexes by computing
\begin{align*}
\tilde{\mu}_n\mu_1&=\sum_{\substack{p+q+r=n\\q\neq n,1}}(-1)^{pq+r} (\id^{\otimes p}\otimes \mu_q\otimes \id^{\otimes r})\mu_{p+1+r}\mu_1\\
&=\sum_{\substack{p+q+r=n\\q\neq n,1}}(-1)^{pq+r}(\id^{\otimes p}\otimes \mu_q\otimes \id^{\otimes r})\sum_{\substack{a+b+c=p+r+1\\b\neq p+r+1}}(-1)^{ab+c+1}(\id^{\otimes a}\otimes \mu_b\otimes \id^{\otimes c})\mu_{a+1+c}\\
&=\sum_{m\neq 1}\sum_{\substack{p+q+r=n\\q\neq n,1}}\sum_{\substack{a+b+c=p+r+1\\a+c+1=m}}(-1)^{pq+r+ab+c+1}(\id^{\otimes p}\otimes \mu_q\otimes \id^{\otimes r})(\id^{\otimes a}\otimes \mu_b\otimes \id^{\otimes c})\mu_m\\
&=\sum_{m\neq 1}\sum_{p+1+r=n}\sum_{\substack{a+b+c=p+r+1\\a+c+1=m}} (-1)^{p+r+ab+c}(\id^{\otimes p}\otimes \mu_1\otimes \id^{\otimes r})(\id^{\otimes a}\otimes \mu_b\otimes \id^{\otimes c})\mu_m\\
&=\sum_{p+1+r=n}(-1)^{p+r}(\id^{\otimes p}\otimes \mu_1\otimes \id^{\otimes r})\tilde{\mu}_n=\sum_{p+1+r=n}(-1)^{n-1}(\id^{\otimes p}\otimes \mu_1\otimes \id^{\otimes r})\tilde{\mu}_n
\end{align*}
where the second to last equality follows from the Lemma \ref{generalised}.
Thus, $\tilde{\mu}_n$ is a lift of the zero map and thus null-homotopic, the homotopy being given by $\mu_{n+1}$.
\end{proof}

We will apply this construction in the case of a basic directed bocs. The bar bimodule resolution $\mathcal{B}$ of the algebra $A$ induces a projective bimodule resolution $\mathcal{B}\oplus \overline{V}$ of $V$. We use this resolution to construct an $A_\infty$-coalgebra structure in $\Mod A\otimes A^{\op}$ on $\mathcal{B}\oplus \overline{V}$. This induces an $A_\infty$-algebra structure on $\Hom_{A\otimes A^{\op}}(\mathcal{B}\oplus \overline{V},\Hom_{\Bbbk}(\mathbb{L},\mathbb{L}))$ by Lemma \ref{dualtocoalgebraobject}. Note that the homology of this complex is $\Ext_{A\otimes A^{\op}}^{\bullet}(V,\Hom_{\Bbbk}(\mathbb{L},\mathbb{L}))\cong \Ext_{\mathfrak{A}}^{\bullet}(\mathbb{L},\mathbb{L})$. Thus, by a special case of Kadeishvili's theorem, Corollary \ref{kadeishvilitheoremhomology}, we obtain an induced $A_\infty$-structure on $\Ext^{\bullet}_{\mathfrak{A}}(\mathbb{L},\mathbb{L})$. The construction is explicit in low degrees which makes it possible to compute the restriction of the $A_\infty$-structure on $\Ext$ in low degrees. For $\overline{V}=0$, this recovers the result of Keller \cite[Proposition 2]{K02} (see also \cite[Theorem A]{LPWZ09}, \cite[Theorem 5.2]{HL05}), that an $A_\infty$-structure on the $\Ext$-algebra of simple modules can be defined from a presentation of the algebra as a quiver with relations.

Let $(A,V,\mu,\varepsilon)$ be a basic directed bocs. Denote by $J$ the Jacobson radical of $A$. Denote the chosen group-like in $V$ by $\omega$ and set  $\omega_{\mathtt{i}}=e_{\mathtt{i}}\omega e_{\mathtt{i}}$. Note that this implies $\varepsilon(\omega_{\mathtt{i}})=e_{\mathtt{i}}\in A$. 

Let $\eta\colon 0\to \overline{V}\to V\to A\to 0$ be the defining exact sequence of $\overline{V}:=\ker \varepsilon$. In the following, we will write $\otimes$ to mean $\otimes_\mathbb{L}$ and just write $xy$ for $x\otimes_A y$. Recall that the bimodule bar resolution $\mathcal{B}$ of $A$ as a bimodule over itself is of the form 

\[\dots\to A\otimes J\otimes J\otimes J\otimes A\stackrel{\delta^{-3}}{\to} A\otimes J\otimes J\otimes A\stackrel{\delta^{-2}}{\to} A\otimes J\otimes A\stackrel{\delta^{-1}}{\to} A\otimes \mathbb{L}\otimes A\stackrel{\delta^0}{\to} A\to 0,\]\label{Vbimoduleresolution}

where all tensors are over $\mathbb{L}$ and  $\delta^{-k}$ is given by

\begin{align*}
\delta^{-k}(1\otimes (a_1\otimes \dots\otimes a_k)\otimes 1)&=a_1\otimes (a_2\otimes \dots\otimes a_k)\otimes 1\\
&+\sum_{j=1}^{k-1}(-1)^j\cdot 1\otimes (a_1\otimes \dots\otimes a_{j-1}\otimes a_ja_{j+1}\otimes a_{j+2}\otimes \dots\otimes a_k)\otimes 1\\
&+(-1)^{k}1\otimes (a_2\otimes \dots\otimes a_{k-1})\otimes a_k.
\end{align*}

To make the shifts more explicit we write the bimodule bar resolution as a dg module 
\[\mathcal{B}=\bigoplus_{j\geq 0} A\otimes (\shift J)^{\otimes j}\otimes A.\]
Using this notation, the differential is given on direct summands by

\begin{align*}
\delta^{-k}(1\otimes (\shift a_1\otimes \dots\otimes \shift a_k)\otimes 1)&=a_1\otimes (\shift a_2\otimes \dots\otimes \shift a_k)\otimes 1\\
&+\sum_{j=1}^{k-1}(-1)^j\cdot 1\otimes (\shift a_1\otimes \dots\otimes \shift a_{j-1}\otimes \shift(a_ja_{j+1})\otimes \shift a_{j+2}\otimes \dots\otimes \shift a_k)\otimes 1\\
&+(-1)^{k}1\otimes (\shift a_2\otimes \dots\otimes \shift a_{k-1})\otimes a_k.
\end{align*}

By assumption, $\overline{V}$ is projective as a bimodule. Recall that $\partial\colon A\to \shift^{-1}\overline{V}$ was defined as $\partial(a)=\shift^{-1}(a\omega-\omega a)$, where $\omega=\sum \omega_{\mathtt{i}}$. Let $p$ denote the map $A\otimes \mathbb{L}\otimes A\to V, \hat{\omega}_{\mathtt{i}}\mapsto \omega_{\mathtt{i}}$ where $\hat{\omega}_{\mathtt{i}}$ denotes the generator $e_\mathtt{i}\otimes 1_\mathbb{L}\otimes e_\mathtt{i}$ of the projective bimodule $Ae_\mathtt{i}\otimes \mathbb{L}\otimes  e_{\mathtt{i}}A$. Write $\hat{\omega}=\sum \hat\omega_{\mathtt{i}}=1_A\otimes 1_\mathbb{L} \otimes 1_A$ Let $\tilde{\delta}^{-1}\colon A\otimes\shift J\otimes A\to (A\otimes \mathbb{L}\otimes A)\oplus \overline{V}$ be the map defined by $1\otimes \shift a\otimes 1\mapsto 1\otimes e_\mathtt{j}\otimes a - a\otimes e_\mathtt{i}\otimes 1- \shift\partial(a)$ for $a\in J(\mathtt{i},\mathtt{j})$, or, equivalently,  $1\otimes \shift a\otimes 1\mapsto \hat{\omega}a -a\hat{\omega} -\shift\partial(a)$.   By the horseshoe lemma, the following is a projective bimodule resolution of $V$, written as a dg module

\begin{equation}\label{projbimodres}\mathcal{P} = \left(\bigoplus_{j\geq 1} A\otimes (\shift J)^j\otimes A\right) \oplus \underbrace{(A\otimes \mathbb{L}\otimes A)\oplus \overline{V}}_{\text{in degree $0$}}
\end{equation}
with differential given by $\tilde{\delta}^{-1}$ starting  in degree $-1$ and $\delta^{-k}$ starting  in degree less than $-1$. There is a map of dg modules $\mathcal{P}\to V$ given by $(p,\iota)$ in degree $0$ where $\iota$ denotes the inclusion $\overline{V}\to V$.

The following lemma shows that regular bocses are in some sense the analogue of admissible ideals in the sense that $\Ext^1_R(\Delta(\mathtt{i}),\Delta(\mathtt{j}))$ corresponds precisely to the number of solid arrows from $\mathtt{i}$ to $\mathtt{j}$, and likewise $\Hom_R(\Delta(\mathtt{i}),\Delta(\mathtt{j}))$ corresponds to the number of dashed arrows from $\mathtt{i}$ to $\mathtt{j}$. After proving it, the authors found that it was already stated by Ovsienko in unpublished notes. It could be well-known in the Kiev school.

\begin{lem}\label{regularbocsesequivalence}
Let $\mathfrak{A}=(A,V)$ be a directed bocs and $R$ its right algebra. Then the following conditions are equivalent:
\begin{enumerate}[(1)]
\item $\mathfrak{A}$ is regular.
\item $\Ext^1_A(L,L)\to \Ext^1_R(\Delta,\Delta)$ is injective (and hence an isomorphism).
\item For $\mathtt{i}\neq \mathtt{j}$ we have $\dim \Hom_R(\Delta(\mathtt{i}),\Delta(\mathtt{j}))=\#\{\begin{tikzcd}\mathtt{i}\arrow[dashed]{r} &\mathtt{j}\end{tikzcd}\}$.
\end{enumerate}
\end{lem}

\begin{proof}
Consider the long exact sequence

\hspace{-5cm}\begin{tikzcd}
0\arrow{r}& \Hom_A(L(\mathtt{i}),L(\mathtt{j}))\arrow{r}&\Hom_{R}(\Delta(\mathtt{i}),\Delta(\mathtt{j}))\arrow{r}{g} &\Hom_{A\otimes A^{\op}}(\overline{V},\Hom_{\Bbbk}(L(\mathtt{i}),L(\mathtt{j})))\arrow[out=0,in=180, looseness=1.2]{dll}[description]{c}\\
&\Ext^1_A(L(\mathtt{i}),L(\mathtt{j}))\arrow{r}{h}&\Ext^1_R(\Delta(\mathtt{i}),\Delta(\mathtt{j}))\arrow{r}&\Ext^1_{A\otimes A^{\op}}(\overline{V},\Hom_{\Bbbk}(L(\mathtt{i}),L(\mathtt{j})))\to \dots
\end{tikzcd}

from \eqref{eqn:pushout} (using that $\Ext^k_\mathfrak{A}(L(\mathtt{i}),L(\mathtt{j}))\cong \Ext^k_R(\Delta(\mathtt{i}),\Delta(\mathtt{j}))$ for $k=0,1$. We see that $h$ is injective if and only if $c=0$ if and only if $g$ is surjective. Furthermore, $g$ is surjective if and only if $\dim \image g=\dim \Hom_{A\otimes A^{\op}}(\overline{V}, \Hom_{\Bbbk} (L(\mathtt{i}), L(\mathtt{j})))$. For $\mathtt{i}\neq \mathtt{j}$, since $\Hom_A(L(\mathtt{i}),L(\mathtt{j}))=0$, $\dim \image g$ is equal to $\dim \Hom_R(\Delta(\mathtt{i}), \Delta(\mathtt{j}))$. Notice furthermore that $\Hom_\Bbbk(L(\mathtt{i}),L(\mathtt{j}))$ is isomorphic to the simple top of the indecomposable projective module $Ae_\mathtt{j}\otimes e_\mathtt{i}A$. Thus, 
$\dim \Hom_{A\otimes A^{\op}}(\overline{V}, \Hom_{\Bbbk} (L(\mathtt{i}), L(\mathtt{j})))$ is equal to the number of indecomposable summands of $\overline{V}$ isomorphic to $Ae_\mathtt{j}\otimes e_\mathtt{i}A$. This number agrees with the number of dashed arrows in the biquiver associated to $\mathfrak{A}$ by definition. This shows that (2) is equivalent to (3). 

To show that (1) is equivalent to (2) note that using the explicit bimodule resolution $\mathcal{P}$, the map $c$ can, when identifying its codomain with a subspace of $\Hom_{A\otimes A^{\op}}(A\otimes J\otimes A, \Hom_\Bbbk(L(\mathtt{i}),L(\mathtt{j})))$, be expressed as the map 
\begin{align*}
\Hom_{A\otimes A^{\op}}(\overline{V},\Hom_\Bbbk(L(\mathtt{i}),L(\mathtt{j})))&\to \Hom_{A\otimes A^{\op}}(A\otimes J\otimes A, \Hom_\Bbbk(L(\mathtt{i}),L(\mathtt{j})))\\
f&\mapsto f\circ \tilde{\delta}^{-1}_{\overline{V}}
\end{align*} 
where $\tilde{\delta}^{-1}_{\overline{V}}$ denotes the composition of $\tilde{\delta}^{-1}$ with the projection to $\overline{V}$. This follows from the explicit definition of the connecting homomorphism using the horsehoe lemma and the snake lemma. As the domain is spanned by projections to the top of the projective bimodule $\overline{V}$, the map $c$ being zero is equivalent to the image of the map $\tilde{\delta}^{-1}_{\overline{V}}$ being contained in the radical of $\overline{V}$. But this is equivalent to $\mathfrak{A}$ being regular as this map is given by $1\otimes a\otimes 1\mapsto -\shift\partial(a)$. \qedhere
\end{proof}

\subsection{Uniqueness of Borel subalgebras}
\label{sec:uniquenessofexactborelsubalgebras}

In this subsection, we prove that for any two basic directed bocses $(A,V)$ and $(B,W)$ whose right algebras are Morita equivalent as quasi-hereditary algebras there is an isomorphism $A\cong B$. Furthermore, we also show that in this case the corresponding right algebras are isomorphic. This is a partial result of Theorem \ref{propermainthm} since it does not imply that the two isomorphisms are compatible. However, it has a much shorter proof which we give here. Generalisations of this theorem are considered in \cite{Thu20, CK21}.

\begin{thm}\label{sameainfty}
Let $\mathfrak{B}=(B,W)$ be a basic directed bocs. Then the functor $R\otimes_B-$ provides an isomorphism $\Ext^{\geq 1}_B(\mathbb{L},\mathbb{L})\cong \Ext^{\geq 1}_R(\Delta,\Delta)$ of $A_\infty$-categories. 
\end{thm}

\begin{proof}
According to Lemma \ref{regularbocsesequivalence}, the functor $R\otimes_B-\colon \modu B\to \modu R$ induces an isomorphism of (graded) vector spaces $\bigoplus_{n\geq 1} \Ext^n_B(\mathbb{L},\mathbb{L})\to \bigoplus_{n\geq 1}\Ext^n_R(R\otimes_B \mathbb{L},R\otimes_B \mathbb{L})=\bigoplus_{n\geq 1}\Ext^n_R(\Delta,\Delta)$. Since $R$ is projective over $B$, the functor $R\otimes_B-$ is also a functor on the derived category, thus the isomorphism is in fact one of (non-unital) $\Bbbk$-algebras.

The next step to prove is to prove some version of functoriality of the $A_\infty$-construction. 
We transfer Merkulov's construction (Theorem \ref{merkulovsconstruction}) of the $A_\infty$-algebra structure via $R\otimes_B -$. To achieve this, let $P^{\bullet}$ be a projective resolution of $\mathbb{L}$. Then $\mathcal{A}=\End_B^{\bullet}(P^{\bullet})$ possesses the structure of a differential graded algebra. To construct the $A_\infty$-structure on $H^{\bullet}(\mathcal{A})\cong \Ext_B^{\bullet}(\mathbb{L},\mathbb{L})$, we apply the strategy of the proof of \ref{merkulovsconstruction}. This amounts to first choosing a vector space complement $\mathcal{L}$ of the cycles $\mathcal{Z}$ of $\mathcal{A}$. Inside $\mathcal{Z}$, one then chooses a vector space complement $\mathcal{H}$ of the boundaries $\mathcal{B}$. The aim is now to construct such a decomposition for the differential graded algebra $\End_R^{\bullet}(R\otimes_B P^{\bullet})$ (with differential $\tilde{\partial}$ compatible with the induction functor). To this end, note that by assumption $R\otimes_B \mathcal{H}=\tilde{\mathcal{H}}$ is a complement of $\tilde{\mathcal{B}}$, the boundaries for $\tilde{\partial}$. The next step is to prove that $\{R\otimes_B f\,|\,f \in \mathcal{B}\}\subseteq \tilde{\mathcal{B}}$. Indeed, recall that $\tilde{\partial}(\tilde{f})=\tilde{f}\circ \tilde{d}-(-1)^{|\tilde{f}|} \tilde{d}\circ \tilde{f}$ where $\tilde{d}=R\otimes_B d$. Functoriality of $R\otimes_B -$ implies that the diagram 
\[
\begin{tikzcd}
\bigoplus_{i} \Hom_B(P_i,P_{i+n})\arrow{d}{R\otimes_B -}\arrow{r}{\partial} &\bigoplus_i \Hom_B(P_i,P_{i+n+1})\arrow{d}{R\otimes_B-}\\
\bigoplus_i \Hom_R(R\otimes P_i, R\otimes P_{i+n})\arrow{r}{\tilde{\partial}}&\bigoplus_i \Hom_R(R\otimes_B P_i, R\otimes P_{i+n+1})
\end{tikzcd}
\]
commutes.
This proves the claim. Faithfulness of $R\otimes_B -$, see Theorem \ref{extnisomorphic}, implies that a complement $\tilde{\mathcal{L}}$ of $\tilde{\mathcal{B}}\oplus \tilde{\mathcal{H}}$ can be chosen such that $\{R\otimes_B f\,|\,f\in \mathcal{L}\}\subseteq \tilde{\mathcal{L}}$. Choose complements $\hat{\mathcal{B}}$, $\hat{\mathcal{L}}$ such that $\{R\otimes_B f\,|\,f\in \mathcal{B}\}\oplus \hat{\mathcal{B}}=\tilde{\mathcal{B}}$ and $\{R\otimes_B f\,|\,f\in \mathcal{L}\}\oplus \hat{\mathcal{L}}=\tilde{\mathcal{L}}$. Define $\tilde{h}$ as a splitting of $\tilde{\partial}$ as follows: It is $0$ on $\hat{\mathcal{L}}$. For $R\otimes_B g\in \{R\otimes_B f\,|\,f\in \mathcal{B}\oplus \mathcal{H}\oplus \mathcal{L}\}$ define $\tilde{G}(R\otimes_B g):= R\otimes_B h(g)$ and define $\tilde{h}$ on $\hat{\mathcal{B}}$ such that $\tilde{h}$ is a splitting of $\tilde{\partial}$. Applying Merkulov's construction to the chosen splitting $\tilde{h}$ implies the theorem since, as noted above, $R\otimes_B \mathcal{H}=\tilde{\mathcal{H}}$. 
\end{proof}

As a first corollary, we obtain a special case of the following theorem due to Keller in the case where the $M_{\mathtt{i}}$ are standard modules. We would like to thank Dag Madsen for making us aware of this theorem. For convenience of the reader, we give a proof in the general case since we couldn't find it in the literature.

\begin{thm}
Let $\Lambda$ be a finite dimensional algebra and let $M_\mathtt{1},\dots,M_\mathtt{n}$ be modules such that $\mathcal{F}(M_\mathtt{1},\dots,M_\mathtt{n})$ is a resolving subcategory. Let $M=\bigoplus_{\mathtt{i}=\mathtt{1}}^{\mathtt{n}}M_\mathtt{i}$. Then $\Ext^{\geq 1}_\Lambda(M,M)$ is generated by $\Ext^1_\Lambda(M,M)$ as an $A_\infty$-category. 
\end{thm}

\begin{proof}
Let $\beta\in \Ext^s_\Lambda(N_0,N_{s+1})$ with $N_0,N_{s+1}\in \mathcal{F}(M_{\mathtt{1}},\dots,M_{\mathtt{n}})$. Since $\mathcal{F}(M_\mathtt{1},\dots,M_\mathtt{n})$ is resolving, $\beta$ is represented by a sequence
\[0\to N_{s+1}\to E_s\to \dots\to E_1\to N_0\to 0\]
with $E_1,\dots,E_s\in \mathcal{F}(M_\mathtt{1},\dots,M_\mathtt{n})$. As $\mathcal{F}(M_{\mathtt{1}},\dots,M_{\mathtt{n}})$ is resolving and thus, in particular, contains the projectives and is closed under kernels of epimorphisms there exist short exact sequences
\[\alpha_i\colon 0\to N_{i}\to E_i\to N_{i-1}\to 0,\quad i=1,\dots,{n+1}\]
in $\mathcal{F}(M_\mathtt{1},\dots,M_\mathtt{n})$ such that $\beta$ is equal to the Yoneda product $\beta=\alpha_{n+1}\cdots \alpha_1$ of these short exact sequences. 

As explained in Section \ref{sec:twistedcomplexes}, there is an equivalence of categories 
\[\tria(M_{\mathtt{1}},\dots,M_{\mathtt{n}})\cong H^0(\twcom \Ext_\Lambda^{\bullet}(M,M)).\] 
Furthermore, the extensions $\alpha_i$ correspond to morphisms $\alpha_i\colon N_{i-1}\to N_i[1]$ in the triangulated hull of $M_\mathtt{1},\dots,M_\mathtt{n}$. By definition, such a morphism is given by a matrix in the $\Ext^1_\Lambda(M_\mathtt{i},M_\mathtt{j})$. Recall that 
\[m_2^{\twcom \mathcal{A}}=\sum_{t=0}^\infty\sum \pm m_{2+t}^\mathcal{A}(\id^{\otimes i_1}\otimes \alpha^{\otimes j_1}\otimes \dots\otimes \id^{\otimes i_r}\otimes \alpha^{\otimes j_r}).\]
As $\beta$ is the Yoneda product of the $\alpha_i$, i.e. the product of the $\alpha_i$ in $H^0(\twcom \Ext^{\bullet}_\Lambda(M,M))$ using $m_2^{\twcom \mathcal{A}}$, the element $\beta$ can be generated by the $\alpha_i$ by applying $m_n^{\mathcal{A}}$ for $n\geq 2$. Since the $\alpha_i$ can be described using only $\Ext^1(M_\mathtt{i},M_\mathtt{j})$, the element $\beta$ is generated in degree $1$. 
\end{proof}

\begin{cor}\label{borel_unique}
Let $\mathfrak{A}=(A,V)$ and $\mathfrak{B}=(B,W)$ be two basic directed bocses with Morita equivalent  quasi-hereditary algebras $R$ and $R'$, then $A\cong B$ as algebras.
\end{cor}

\begin{proof}
By the Theorem \ref{sameainfty}, $\Ext_A^{\geq 1}(\mathbb{L},\mathbb{L})\cong \Ext_{B}^{\geq 1}(\mathbb{L},\mathbb{L})$ as $A_\infty$-algebras. However, by \cite[Proposition 2]{K02}, this $A_\infty$-structure determines the quivers with relations of $A$ and $B$ up to isomorphism. Hence, they are isomorphic.
\end{proof}

\begin{lem}\label{multiplicities}
Let $\mathfrak{A}=(A,V)$ be a regular directed bocs with right algebra $R$. Let $T$ be the decomposition matrix of $R$, i.e. the matrix with entries $t_{ij}=[P_R(\mathtt{j}):\Delta_R(\mathtt{i})]$, where $P_R(\mathtt{j})$ denotes the projective cover of $L_R(\mathtt{j})$. Then the vector of multiplicities of $P_R(\mathtt{j})$ as direct summands of $Q_R(\mathtt{i}):=R\otimes_A P_A(\mathtt{i})$ is given by $T^{-1}\underline{\dim} P_A(\mathtt{i})$, where $\underline{\dim}$ denotes the dimension vector. 
\end{lem}

\begin{proof}
It is well known that the $\Delta$-dimension vector $\dim_\Delta M$ of a module $M$ in $\mathcal{F}(\Delta)$, i.e. the vector of multiplicities of $\Delta(\mathtt{j})$ in a $\Delta$-filtration of $M$, is well-defined, see e.g. \cite[Lemma 2.3]{DR92}. As $R\otimes_A -$ is exact, and the induction of a simple $A$-module is a standard $R$-module, $\dim_\Delta Q(\mathtt{i})=\underline{\dim} P_A(\mathtt{i})$. But decomposing $Q(\mathtt{i})\cong \bigoplus_{\mathtt{j}=\mathtt{1}}^{\mathtt{n}} P_R(\mathtt{j})^{a_{\mathtt{i}\mathtt{j}}}$ we see that $\dim_\Delta Q(\mathtt{i})=T\cdot (a_{\mathtt{i}\mathtt{j}})_{\mathtt{j}=\mathtt{1},\dots,\mathtt{n}}$. The claim follows as $T$ is invertible.
\end{proof}

A recursive formula for the entries of the decomposition matrix is given in \cite{Con21}.

\begin{cor}\label{moritaiso}
Let $\mathfrak{A}=(A,V)$ and $\mathfrak{B}=(B,W)$ be basic directed bocses. Assume that the associated right algebras $R$ and $S$ are Morita equivalent as quasi-hereditary algebras. Then $R$ and $S$ are isomorphic.
\end{cor}

\section{$A_\infty$-structures on Ext-algebras and uniqueness of bocses}\label{sec:uniquenessofbocses}

Throughout this section, we consider a basic directed bocs $(A,V)$. We explicitly compute parts of the $A_\infty$-coalgebra structure on the projective resolution of $V$ given in  \eqref{projbimodres} (Section \ref{sec:explicitcalc}), relate the resulting comultiplications to the differential on the semifree dg algebra associate to $(A,V)$ (Section \ref{sec:muanddel}), and construct some well-behaved splittings of certain natural projections appearing in a presentation of $A$ by quiver and relations (Section \ref{splittings}). Finally, in Section \ref{sec:mainthm} we construct an explicit $A_\infty$-algebra structure on the Ext-algebra of standard modules of the associated quasi-hereditary algebra and use this to prove uniqueness of the bocs up to isomorphism.

\subsection{An explicit $A_\infty$-coalgebra stucture on $\mathcal{P}$}
\label{sec:explicitcalc}

In this subsection, using Lemma \ref{AinfinitycoalgebraobjectP}, we will inductively construct the first terms of an $A_\infty$-coalgebra object structure on $\mathcal{P}$ from \eqref{projbimodres}. Proofs are partially relegated to the appendix.

Consider the comultiplication $\mu\colon V\to V\otimes_A V$. Recall that the projective bimodule resolution $\mathcal{P}$ of $V$ from \eqref{projbimodres} has components $\mathcal{P}^0=(A\otimes \mathbb{L}\otimes A)\oplus \overline{V}$ and $\mathcal{P}^{-i}=A\otimes (\shift J)^{\otimes i}\otimes A$ for all $i\geq 1$. Recall that $\Phi$ was chosen as a complement to $J\overline{V}+\overline{V}J$ in $\overline{V}$.

Using the canonical identification of $A\otimes (\shift J)^{\otimes i}\otimes A\otimes_A A\otimes (\shift J)^{\otimes j}\otimes A\cong A\otimes (\shift J)^{\otimes i}\otimes A\otimes (\shift J)^{\otimes j}\otimes A$, thus always omitting tensors over $A$, $\mu_1$ is simply the differential on $\mathcal{P}$, explicitly given by 
\begin{align*}
\mu_1^0&=0\\
\mu_1^1(1\otimes \shift a\otimes 1)&= \hat{\omega} a - a\hat{\omega} -\shift\partial(a)\\
\mu_1^i(1\otimes \shift a_1\otimes \dots\otimes \shift a_i\otimes 1)&= a_1\otimes \shift a_2\otimes \dots\otimes \shift a_i\otimes 1\\
&+\sum_{j=1}^{i-1}(-1)^j1\otimes \shift a_1\otimes \dots\otimes \shift a_{j-1}\otimes \shift (a_ja_{j+1})\otimes \shift a_{j+2}\otimes \dots\otimes \shift a_i\otimes 1\\
&+(-1)^i 1\otimes \shift a_1\otimes \dots\otimes \shift a_{i-1}\otimes a_i.
\end{align*}

The next two lemmas provide explicit formulae for the component maps of $\mu_2$. 

\begin{lem}\label{mu201}

For $\varphi\in \Phi$ and $a\in J$ with, using Sweedler notation, $\partial(a)=a^1(\shift^{-1}a^2)a^3$ with $a^1,a^3\in J\cup \mathbb{L}$ and $a^2\in \overline{V}$ we obtain
\begin{align*}
\mu_2^0(\hat{\omega})&=\hat{\omega}\hat{\omega},\\
\mu_2^0(\varphi)&=\hat{\omega}\varphi+\varphi\hat{\omega}+(\shift\otimes\shift)\partial(\shift^{-1}\varphi),\\
\mu_2^1(1\otimes \shift a\otimes 1)&=\hat{\omega}\otimes \shift a\otimes 1+1\otimes \shift a\otimes \hat{\omega}-a^1a^2\otimes \shift a^3\otimes 1 +1\otimes \shift a^1\otimes a^2a^3,
\end{align*}
where the last equation is to be understood as summing only over those summands  where $a^3\in J$ for the third, resp. $a^1\in J$ for the fourth summand.   
\end{lem}

\begin{proof}
To prove the formula for $\mu_2^0$ we have to show that $\mu\circ (p,\iota)=((p,\iota)\otimes (p,\iota))\mu_2^0$. Indeed,
\begin{align*}
\mu\circ(p,\iota)(\hat{\omega})&=\mu(\omega)=\omega\omega\\
\mu\circ(p,\iota)(\varphi)&=\mu(\varphi)=\omega \varphi+\varphi \omega+(\shift \otimes \shift)\partial(\shift^{-1}\varphi)
\end{align*}
and
\begin{align*}
((p,\iota)\otimes (p,\iota))\mu_2^0(\hat{\omega})&=((p,\iota)\otimes (p,\iota))\hat{\omega} \hat{\omega}\\
&=\omega\otimes \omega\\
((p,\iota)\otimes (p,\iota))\mu_2^0(\varphi)&=((p,\iota)\otimes (p,\iota))(\hat{\omega} \varphi+\varphi\hat{\omega} +(\shift \otimes \shift)\partial(\shift^{-1}\varphi))\\
&=\omega\varphi+\varphi \omega+(\shift\otimes\shift)\partial(\shift^{-1}\varphi).
\end{align*}

Furthermore, in degree $1$ we obtain
\begin{align*}
\mu_2^0\mu_1^1(1\otimes \shift a\otimes 1)&=\mu_2^0\left(\hat{\omega}a-a\hat{\omega}-a^1a^2a^3\right)\\
&=\underbrace{\hat{\omega}\hat{\omega} a}_{1}-\underbrace{a\hat{\omega} \hat{\omega}}_{2}-\underbrace{a^1\hat{\omega} a^2a^3}_{3}- \underbrace{a^1a^2\hat{\omega} a^3}_{4}-\underbrace{a^1(\shift \otimes \shift)\partial(\shift^{-1}a^2)a^3}_{5} 
\end{align*}
and
\begin{align*}
&-(\mu_1\otimes 1+1\otimes \mu_1)\mu_2^1(1\otimes \shift a\otimes 1)\\
=\,&-(\mu_1\otimes 1+1\otimes \mu_1)\bigg(\hat{\omega}\otimes \shift a\otimes 1+1\otimes \shift a\otimes \hat{\omega}-a^1a^2\otimes \shift a^3\otimes 1+1\otimes \shift a^1\otimes a^2a^3\bigg)\\
=\,&-\underbrace{\hat{\omega}\hat{\omega}a}_{1}+\underbrace{\hat{\omega}a\hat{\omega}}_{6}+\underbrace{\hat{\omega} a^1a^2a^3}_{8}-\underbrace{\hat{\omega}a\hat{\omega}}_{6}+\underbrace{a\hat{\omega} \hat{\omega}}_{2}+\underbrace{a^1a^2a^3\hat{\omega}}_{7}\\
&+\underbrace{a^1a^2\hat{\omega} a^3}_{4}-\underbrace{a^1a^2a^3\hat{\omega}}_{7}-\underbrace{a^1a^2 \shift\partial(a^3)}_{5}-\underbrace{\hat{\omega} a^1a^2a^3}_{8}+\underbrace{a^1\hat{\omega} a^2a^3}_{3}+\underbrace{\shift\partial(a^1)a^2a^3}_{5}.
\end{align*}
To see that the two expressions sum to zero, we compare the parts with the same number, noting that the three terms marked $5$ add to zero as $\partial$ is a differential and therefore
\begin{align*}
0&=(\shift \otimes \shift)\partial^2(a)=(\shift\otimes\shift)\partial(a^1 \shift^{-1}(a^2) a^3)\\
&=(\shift\otimes\shift)(\partial(a^1)\shift^{-1}(a^2)a^3+a^1\partial(\shift^{-1}a^2)a^3-a^1\shift^{-1}(a^2)\partial(a^3)\\
&=-\shift\partial(a^1)a^2a^3+a^1(\shift\otimes\shift)\partial(\shift^{-1}a^2)a^3+a^1a^2\shift\partial(a^3).
\end{align*}.
\end{proof}

We continue by inductively constructing $\mu_2^i$. 

\begin{lem}\label{mu2i}
For $i\geq 2$ and $a_1,\dots,a_i\in J$, $\mu_2^i$ can be chosen as 
\begin{align*}
\mu_2^i(1\otimes \shift a_1\otimes \dots\otimes  \shift a_i\otimes 1)=\,&\hat{\omega}\otimes  \shift a_1\otimes \dots\otimes  \shift a_i\otimes 1+1\otimes  \shift a_1\otimes \dots\otimes  \shift a_i\otimes \hat{\omega}\\
&-\sum_{j=1}^{i-1}  1\otimes  \shift a_1\otimes \dots\otimes  \shift a_j\otimes 1\otimes  \shift a_{j+1}\otimes \dots\otimes  \shift a_i\otimes 1\\
&-a^1_{1}a^2_{1}\otimes  \shift a^3_{1}\otimes  \shift a_2\otimes \dots\otimes  \shift a_i\otimes 1\\
&+1\otimes  \shift a_1\otimes \dots\otimes  \shift a_{i-1}\otimes  \shift a_{i}^1\otimes a_{i}^2a_{i}^3
\end{align*}
\end{lem}

\begin{proof}
This proof is given in Appendix \ref{appendixA}.
\end{proof}

Note that 
\[(\mathcal{P}^{\otimes n})^k = \underset{k_1+\cdots +k_n=k}{\bigoplus}\mathcal{P}^{k_1}\otimes\cdots\otimes\mathcal{P}^{k_n}.\] Denoting $\mathcal{P}^{0,0}=A\otimes \mathbb{L}\otimes A$ and  $\mathcal{P}^{0,1}=\overline{V}$, we can refine this sum further. Let 
\[{}^{\mathbb{L}}\mathcal{P}^{\otimes n}=\underset{\substack{k_1+\cdots +k_n=k\\ k_j=0 }}{\sum}\mathcal{P}^{k_1}\otimes\cdots\otimes \mathcal{P}^{k_{j-1}}\otimes\mathcal{P}^{0,0}\otimes  \mathcal{P}^{k_{j+1}}\otimes\cdots\otimes\mathcal{P}^{k_n}\]

\begin{lem}\label{projLzero}
For $n\geq 3$ and any $i$, the composition of $\mu_n^i$ with the projection to ${}^{\mathbb{L}}\mathcal{P}^{\otimes n}$ is zero. 
\end{lem}

\begin{proof}
Direct inspection shows that $\mu_1$ and $\mu_2$ as computed in Lemmas \ref{mu201} and \ref{mu2i} satisfy the hypotheses of Proposition \ref{counitality}. This implies that the homotopies inductively defining the higher comultiplications in Lemma \ref{AinfinitycoalgebraobjectP} can be chosen in such a way that the resulting $A_\infty$-coalgebra structure on $\mathcal{P}$ becomes strictly counital. Summands in $\mu_n^i(\underline{a})$ potentially violating the statement of the lemma can be written as $\underline{x}b_1\hat{\omega}b_2\underline{y}$ where $b_1,b_2\in A$, and $\underline{x}$ (resp. $\underline{y}$) is either empty or ends (respectively starts) in the symbol $\otimes$ or a bimodule generator $\varphi$ of $\overline{V}$. Counitality implies that $\underline{x}b_1b_2\underline{y} =0$, which, given the form of  $\underline{x}$ and $\underline{y}$ described above, requires $b_1b_2=0$ in $A$.

Now we inductively assume that the statement of the lemma holds for $\mu_j$, whenever $j=3,\dots,n-1$, and  for $\mu_n^k$, whenever $k<i$. Using
\[
0=\sum (-1)^{rs+t} (\id^{\otimes r}\otimes \mu_s\otimes \id^{\otimes t})\mu_{r+1+t}\]
we write  
\begin{align}
\label{ainfinityABCDE}\underbrace{\sum (-1)^{r+t} (\id^{\otimes r}\otimes \mu_1\otimes \id^{\otimes t})\mu_{n}^{i}}_{(A)}=&\underbrace{\sum (-1)^{2r+t+1} (\id^{\otimes r}\otimes \mu_2\otimes \id^{\otimes t})\mu_{n-1}}_{(B)}\\
&+
\underbrace{\sum_{s\geq 3, r+1+t\geq 3} (-1)^{rs+t+1} (\id^{\otimes r}\otimes \mu_s\otimes \id^{\otimes t})\mu_{r+1+t}}_{(C)}\notag\\
&+
\underbrace{(\mu_{n-1}\otimes \id+(-1)^{n}\id\otimes \mu_{n-1}) \mu_2}_{(D)}-\underbrace{\mu_n^{i-1}\mu_1}_{(E)}\notag,
\end{align}
where we omit superscripts on all $\mu_j$ for $j<n$. 

By the induction hypothesis, terms $(C)$ and $(E)$ in \eqref{ainfinityABCDE} do not have summands of the form $\underline{x}b_1\hat{\omega}b_2\underline{y}$. We thus need to consider only the terms labelled $(B)$ and $(D)$. 

In $(B)$, by looking at the explicit formulae for $\mu_2^k$ we see that the only terms of the form $\underline{x}b_1\hat{\omega}b_2\underline{y}$ appearing satisfy $b_1=b_2=1$, contradicting the need for $b_1b_2=0$, and thus $(B)$ does not contribute any summands of the form $\underline{x}b_1\hat{\omega}b_2\underline{y}$.

For $(D)$, recall that $\mu_2(\underline{a})=\hat{\omega} \underline{a}+\underline{a}\hat{\omega}+y$ for some $y$ which is annihilated by the projection from $\mathcal{P}^{\otimes 2}$ to $\mathcal{P}^{0,0}\otimes \mathcal{P}+\mathcal{P}\otimes \mathcal{P}^{0,0}$. Since, by the induction hypothesis, composition of $\mu_{n-1}$ with projection onto ${}^{\mathbb{L}}\mathcal{P}^{\otimes n-1}$ is zero, any $\hat{\omega}$ appearing in 
$(D)$ must come from applying $\mu_{n-1}\otimes \id$ to $\underline{a}\hat{\omega}$ or from applying $(-1)^{n}\id\otimes \mu_{n-1}$ to $\hat{\omega} \underline{a}$. In these cases, they appear at the end or at the beginning of the summand, hence with $b_2=1$ or $b_1=1$, which again contradicts $b_1b_2=0$ in $A$. 

Therefore, the right hand side of \eqref{ainfinityABCDE} does not have any summands of the form $\underline{x}b_1\hat{\omega}b_2\underline{y}$ and the composition of  $\sum (-1)^{r+t} (\id^{\otimes r}\otimes \mu_1\otimes \id^{\otimes t})\mu_{n}^{i}$ with projection onto ${}^{\mathbb{L}}\mathcal{P}^{\otimes n}$ is zero. Since $\mathcal{P}$ is concentrated in nonpositive homological degrees, this implies that composition of  $\mu_{n}^i$ with projection onto ${}^{\mathbb{L}}\mathcal{P}^{\otimes n}$ is also zero, as required.
\end{proof}

To prove our main theorem, we will partially compute the $A_\infty$-algebra structure on $\Ext_R^\bullet(\Delta,\Delta)$ (where $R$ is the associated right algebra of the bocs $(A,V)$) by using the description of $\Ext_R^\bullet(\Delta,\Delta)$ as $\Hom_{A\otimes A^{\mathrm{op}}}(\mathcal{P}, \Hom_{\Bbbk}(\mathbb{L}, \mathbb{L}))$ from Section \ref{sec:projectiveresolution} and dualising the $A_\infty$-coalgebra structure on $\mathcal{P}$. Since we are interested in the quotient of the cobar resolution of $\Ext_R^\bullet(\Delta,\Delta)^\#$ by the differential ideal generated by negative degrees, we only need partial results to determine the differential on the associated semifree dg algebra. More precisely, the generators of the cobar resolution of $\Ext_R^\bullet(\Delta,\Delta)^\#$, as a tensor algebra over $ \mathbb{L}$, in degrees $-1,0$ and $1$ are given by $\shift^{-1}(\Ext_R^2(\Delta,\Delta)^\#)$, $\shift^{-1}(\Ext_R^1(\Delta,\Delta)^\#)$ and $\shift^{-1}(\radoperator_R(\Delta,\Delta)^\#)$, respectively, so we are only interested in the comultiplications $\mu_n^j$ in $\mathcal{P}$ for $j=0,1,2$. Similarly, since any term in the evaluation of a $\mu_n^j$ with a tensor factor corresponding to an $\Ext_R^2(\Delta,\Delta)^\#$ is sent to zero in the quotient semifree dg algebra, we only need to consider comultiplications $\mu_n^j$ with nonzero corestriction to direct summands of $\mathcal{P}^{\otimes n}$ all of whose tensor factors live in degree $0$ or $-1$.

Motivated, by this, in the following, we denote by ${}^{\small \textcircled{{\tiny 2}}}\mathcal{P}^{\otimes n}$ the subspace of $\mathcal{P}^{\otimes n}$ not involving tensor factors of the form $A\otimes  (\shift J)^{\otimes m}\otimes A$ for $m\geq 2$, that is
\[{}^{\small \textcircled{{\tiny 2}}}\mathcal{P}^{\otimes n} = \underset{k\geq 0}{\bigoplus}\;\;{\bigoplus_{\substack{k_1+\cdots +k_n=k\\k_i \geq -1 \;\;\forall i}}}\mathcal{P}^{k_1}\otimes\cdots\otimes\mathcal{P}^{k_n}.\]
We denote by $\pi$ the projection of $\mathcal{P}^{\otimes n}$ onto ${}^{\small \textcircled{{\tiny 2}}}\mathcal{P}^{\otimes n}$.

In the appendices, we compute parts of the $A_\infty$-coalgebra structure on $\mathcal{P}$. The computations and the resulting expressions for the relevant $\mu_i$ are long, however, in the following, only their projections onto ${}^{\small \textcircled{{\tiny 2}}}\mathcal{P}^{\otimes n}$ are relevant, which are given by the following proposition and consist of very few terms.

\begin{prop}\label{mu26}
Let $\mathcal{P}$ be the projective resolution of $V$ defined in \eqref{projbimodres}. Then the lifting $\mu_2$ of the comultiplication $\mu$ and the higher comultiplications $\mu_n$ (for $n\geq 3$) can be chosen in such a way that their projections $\pi\mu_j$ to ${}^{\small \textcircled{{\tiny 2}}}\mathcal{P}^{\otimes n}$ satisfy:
\begin{align*}
\pi\mu_1^0&=0\\
\pi\mu_1^1(1\otimes  \shift a\otimes 1)&=\hat{\omega} a-a\hat{\omega} -  \shift \partial(a)\\
\pi\mu_1^2(1\otimes  \shift a_1\otimes  \shift a_2\otimes 1)&=a_1\otimes  \shift a_2\otimes 1-1\otimes  \shift (a_1a_2)\otimes 1+1\otimes  \shift a_1\otimes a_2\\
\pi\mu_1^i(1\otimes  \shift a_1\otimes \dots\otimes  \shift a_i\otimes 1)&=0 \text{ for $i\geq 3$}\\
\pi\mu_2^0(\varphi)&=\hat{\omega}\varphi+\varphi\hat{\omega}+( \shift \otimes  \shift )\partial( \shift^{-1} \varphi)\\
\pi\mu_2^1(1\otimes  \shift a\otimes 1)&=\hat{\omega}\otimes  \shift a\otimes 1+1\otimes  \shift a\otimes \hat{\omega}-a^1a^2\otimes  \shift a^3\otimes 1+1\otimes  \shift a^1\otimes a^2a^3\\
\pi\mu_2^2(1\otimes  \shift a_1\otimes  \shift a_2\otimes 1)&=-1\otimes  \shift a_1\otimes 1\otimes  \shift a_2\otimes 1\\
\pi\mu_2^n(1\otimes  \shift a_1\otimes \dots\otimes  \shift a_i\otimes 1)&=0\text{ for $i\geq 3$}\\
\pi\mu_3^0(\varphi)&=1\otimes \shift \varphi^1\otimes \varphi^2\dots \varphi^5-\varphi^1\varphi^2\otimes \shift \varphi^3\otimes \varphi^4\varphi^5+\varphi^1\dots \varphi^4\otimes \shift \varphi^5\otimes 1 \\
\pi\mu_3^1(1\otimes  \shift a\otimes 1)&=1\otimes  \shift a^1\otimes a^2\otimes  \shift a^3\otimes 1 \\
\pi\mu_3^i(1\otimes \shift a_1\otimes \dots\otimes \shift a_i\otimes 1)&=0 \text{ for $i\geq 2$}\\
\pi\mu_4^0(\varphi)&=\varphi^1\varphi^2\otimes  \shift \varphi^3\otimes \varphi^4\otimes  \shift \varphi^5\otimes 1-1\otimes  \shift \varphi^1\otimes \varphi^2\varphi^3\varphi^4\otimes \shift  \varphi^5\otimes 1 \\
&\;\;\;+1\otimes  \shift \varphi^1\otimes \varphi^2\otimes  \shift \varphi^3\otimes \varphi^4\varphi^5 \\
\pi\mu_4^1(1\otimes \shift a\otimes 1)&=0\\
\pi\mu_4^2(1\otimes \shift a_1\otimes \shift a_2\otimes 1)&=0\\
\pi\mu_5^0(\varphi)&=- 1\otimes  \shift \varphi^1\otimes \varphi^2\otimes  \shift \varphi^3\otimes \varphi^4\otimes  \shift \varphi^5\otimes 1 \\
\pi\mu_5^1(1\otimes \shift a\otimes 1)&=0\\
\pi\mu_6^0(\varphi)&=0
\end{align*}
\end{prop}

\begin{proof}
This follows from the discussion at the beginning of this section, Lemmas \ref{mu201},  \ref{mu2i}, and Appendices \ref{appendixB} to \ref{appendixE}.
\end{proof}

\begin{lem}\label{property2}
The $A_\infty$-coalgebra structure on $\mathcal{P}$ can be chosen such that the only comultiplications $\mu_n^i$ with non-zero projections to ${}^{\small \textcircled{{\tiny 2}}}\mathcal{P}^{\otimes n}$ are $\mu_1^0,\dots,\mu_5^0, \mu_1^1,\dots, \mu_3^1, \mu_1^2, \mu_2^2$.  
\end{lem}

\begin{proof}
For notational simplicity we denote the following property by \textcircled{\raisebox{-0.9pt}{\small 2}}:
\begin{enumerate}
\item[\textcircled{\raisebox{-0.9pt}{\small 2}}] The projection of $\mu_n^i$ to ${}^{\small \textcircled{{\tiny 2}}}\mathcal{P}^{\otimes n}$
is zero.
\end{enumerate}
We thus have to prove that any $\mu_n^i$ not in the list satisfies \textcircled{\raisebox{-0.9pt}{\small 2}}. By Proposition \ref{mu26}, $\mu_6^0$, $\mu_5^1$, $\mu_4^1$, $\mu_3^2$, $\mu_2^3$, and $\mu_1^3$ satisfy \textcircled{\raisebox{-0.9pt}{\small 2}}.

Using the notation from \eqref{ainfinityABCDE}, by definition of $\mu_1$, the comultiplication $\mu_n$ will satisfy \textcircled{\raisebox{-0.9pt}{\small 2}} unless $(A)$ does not satisfy \textcircled{\raisebox{-0.9pt}{\small 2}} and preimages are taken under $(\id^{\otimes r}\otimes \mu_1^0\otimes \id^{\otimes t})$ for some $r,t$. In order to take a preimage under $(\id^{\otimes r}\otimes \mu_1^0\otimes \id^{\otimes t})$ an $\hat{\omega}$ in position $r+1$ has to exist. By Lemma \ref{projLzero}, $\hat{\omega}$ can only appear in $(B)$ and $(D)$. Let $n_0=6, n_1=4, n_2=3, n_j=0$ for $j\geq 3$. We consider $\mu_n^i$ for $n>n_i$. 

We start our analysis with summand (B). By induction, $\mu_{n-1}^i$ has property \textcircled{\raisebox{-0.9pt}{\small 2}}. Recalling 
\begin{align*}
\mu_2^0(\omega)=\,&\omega\omega\\
\mu_2^0(\varphi)=\,&\omega\varphi+\varphi\omega+\partial(\varphi)\\
\mu_2^i(1\otimes a_1\otimes \dots\otimes a_i\otimes 1)=\,&\omega\otimes a_1\otimes \dots\otimes a_i\otimes 1+1\otimes a_1\otimes \dots\otimes a_i\otimes \omega\\
&-\sum_j 1\otimes a_1\otimes \dots\otimes a_j\otimes 1\otimes a_{j+1}\otimes \dots\otimes a_i\otimes 1\\
&+a_1^1a_1^2\otimes a_1^3\otimes a_2\otimes\dots\otimes a_i\otimes 1\\
&-1\otimes a_1\otimes \dots\otimes a_{i-1}\otimes a_i^1\otimes a_i^2a_i^3,
\end{align*}
we see that the only way of creating summands without property \textcircled{\raisebox{-0.9pt}{\small 2}} is via the third term of $(\id^{\otimes r}\otimes \mu_2^2\otimes \id^{\otimes t})$ applied to the only tensor factor of the form $A\otimes J^{\otimes 2}\otimes A$. However, this term does not contain any $\hat{\omega}$. Thus, $(B)$ has terms with property \textcircled{\raisebox{-0.9pt}{\small 2}} and terms without property \textcircled{\raisebox{-0.9pt}{\small 2}} but also without $\hat{\omega}$. 

We next analyse summand (D). The comultiplication $\mu_2^i$ has terms in $A\otimes J^{\otimes m}\otimes A\otimes J^{\otimes i-m}\otimes A$, and $A\otimes (\mathbb{L}\oplus\Phi)\otimes A\otimes J^{\otimes i}\otimes A$ and $A\otimes J^{\otimes i}\otimes A\otimes (\mathbb{L}\oplus\Phi)\otimes A$. By abuse of notation, we will identify the case $m=0$ with the term $A\otimes (\mathbb{L}\oplus\Phi)\otimes A\otimes J^{\otimes i}\otimes A$. By symmetry of the terms in $\mu_2^i$, it suffices to apply $\mu_{n-1}^m\otimes 1$ to terms in $A\otimes J^{\otimes m}\otimes A\otimes J^{\otimes i-m}\otimes A$. If $i-m\geq 2$, the result will satisfy \textcircled{\raisebox{-0.9pt}{\small 2}}. Else, first consider $n=4$ and $i\geq 3$ (the case $i=2$ satisfies \textcircled{\raisebox{-0.9pt}{\small 2}} by Proposition \ref{mu26}). In this case, we apply $\mu_3^{m}$ with $m\geq i-1\geq 2$ which satisfies \textcircled{\raisebox{-0.9pt}{\small 2}} by Proposition \ref{mu26}. Next, for $n=5$ and $m\geq i-1\geq 1$, we apply $\mu_4^{m}$ which again satisfies \textcircled{\raisebox{-0.9pt}{\small 2}} by induction. For $n= 6$ and $m\geq i-1\geq 1$, we apply $\mu_5^m$ which has property \textcircled{\raisebox{-0.9pt}{\small 2}} by induction (the case $m=0$ has been treated in Proposition \ref{mu26}). For $n\geq 7$ and $m\geq i-1\geq 0$, we apply $\mu_{n-1}^m$ which satisfies property \textcircled{\raisebox{-0.9pt}{\small 2}} by induction. Thus, $(D)$ has property \textcircled{\raisebox{-0.9pt}{\small 2}}. 

Hence, all terms involving $\hat{\omega}$ have property \textcircled{\raisebox{-0.9pt}{\small 2}} and therefore every preimage under $(\id^{\otimes r}\otimes \mu_1\otimes \id^{\otimes t})$ has property \textcircled{\raisebox{-0.9pt}{\small 2}}. 
\end{proof}

\begin{lem}\label{classicalcase}
The projection of $\mu_n^i$ onto $\mathcal{P}^{-1}\otimes\cdots\otimes\mathcal{P}^{-1}$ is given by $1\otimes  \shift a_1\otimes 1\otimes  \shift a_2\otimes 1$ for $\mu_2^2$ and is zero for all other $n$ and $i$.
\end{lem}

\begin{proof}
The equation $k_1+\cdots +k_n=-i+2-n$ in this case yields $-n=-i+2-n$ and hence $i=2$. Moreover, if the projection onto $\mathcal{P}^{-1}\otimes\cdots\otimes\mathcal{P}^{-1}$ is nonzero, in particular the projection onto ${}^{\small \textcircled{{\tiny 2}}}\mathcal{P}^{\otimes n}$ is nonzero, and by Proposition \ref{property2} these cases are completely described in Proposition \ref{mu26}. Direct inspection of the formulae in Proposition \ref{mu26} shows that the projection onto $\mathcal{P}^{-1}\otimes\cdots\otimes\mathcal{P}^{-1}$ is given by is given by $1\otimes  \shift a_1\otimes 1\otimes  \shift a_2\otimes 1$ for $\mu_2^2$ and is zero for all other $n$.
\end{proof}

\subsection{Truncated multiplications and the differential on the semifree dg algebra}\label{sec:muanddel}
We will later be interested in $\mu_n^{\mathbb{L}}:=\mathbb{L}\otimes_A \mu_n \otimes_A \mathbb{L}$. We write $\mu_n^{i,\mathbb{L}}:=\mathbb{L}\otimes_A \mu_n^i \otimes_A \mathbb{L}$ and, in this section, establish a connection between $\mu_n^{\mathbb{L}}$ and $\partial$.

Consider the isomorphism $\rho\colon A\otimes \shift^{-1}\Phi\otimes A\to \shift^{-1}\overline{V}$ defined by $a\otimes \shift^{-1}\varphi\otimes b\mapsto a\shift^{-1}(\varphi)b$ using the embedding of $\Phi$ into $\overline{V}$. 

Notice that $\rho^{-1}\circ \partial_0(J)\subseteq  A\otimes \shift^{-1}\Phi \otimes J +  J\otimes\shift^{-1}\Phi \otimes A$ and observe that $A\otimes\shift^{-1} \Phi \otimes J +  J\otimes\shift^{-1}\Phi \otimes A = J\otimes \shift^{-1}\Phi \otimes J \oplus  J\otimes\shift^{-1}\Phi \oplus \shift^{-1}\Phi \otimes J$, where we omit tensor factors of the form $\mathbb{L}$ as usual.  We denote by $\mathtt{p}^1_2,\mathtt{p}^1_3$ the projections onto $ J\otimes\shift^{-1}\Phi \oplus \shift^{-1}\Phi \otimes J$ and $J\otimes \shift^{-1}\Phi \otimes J $ respectively.

\begin{lem}\label{claims910}
\begin{enumerate}[(i)]
\item\label{claim9} The projection to $\Phi \otimes \shift J \oplus  \shift J\otimes\Phi$ of $\mu_2^{1,\mathbb{L}}$  and $(\shift\otimes\shift)\mathtt{p}^1_2\rho^{-1}\partial_0 \shift^{-1}\vert_{\shift J}$ are equal.
\item\label{claim10}The projections to $\shift J\otimes\Phi \otimes \shift J$ of $\mu_3^{1,\mathbb{L}}$  and $(\shift\otimes\shift\otimes\shift)\mathtt{p}^1_3\rho^{-1}\partial_0 \shift^{-1}\vert_{\shift J} $ differ by a factor of $-1$.\end{enumerate}
\end{lem}

\begin{proof}
\begin{enumerate}[(i)]
\item  Consider the diagram
\begin{equation*}
\begin{tikzcd}
\shift J\ar{r}{\mu_2^{1,\mathbb{L}}}\ar{d}{\shift^{-1}} & \Phi\otimes \shift J\oplus \shift J\otimes \Phi\\
J\ar{r}{\mathtt{p}^1_2\rho^{-1}\partial_0}& \shift^{-1} \Phi\otimes  J\oplus J \otimes \shift^{-1} \Phi\ar{u}{\shift\otimes\shift}
\end{tikzcd}
\end{equation*}
which commutes by the explicit description of $\mu_2^1$ in Lemma \ref{mu26}, since 
\[(\shift\otimes\shift)\mathtt{p}^1_2\rho^{-1}\partial_0 \shift^{-1} (\shift a) = 
(\shift\otimes\shift)(a^1\shift^{-1} a^2\otimes a^3 + a^1\otimes \shift^{-1} a^2a^3) = -a^1a^2\otimes \shift a^3 +\shift a^1\otimes a^2a^3.\]
\item  Consider the diagram
\begin{equation*}
\begin{tikzcd}
\shift J\ar{r}{\mu_3^{1,\mathbb{L}}}\ar{d}{\shift^{-1}} & \shift J\otimes \Phi \otimes \shift J\\
J\ar{r}{-\mathtt{p}^1_3\rho^{-1}\partial_0}& J \otimes \shift^{-1} \Phi\otimes J.\ar{u}{\shift\otimes\shift\otimes\shift}
\end{tikzcd}
\end{equation*}
This commutes by the explicit description of $\mu_3^1$ in Lemma \ref{mu26} since 
\[(\shift\otimes\shift\otimes\shift)\mathtt{p}^1_2\rho^{-1}\partial_0 \shift^{-1} (\shift a) = 
(\shift\otimes\shift\otimes\shift)(a^1\otimes \shift^{-1} a^2\otimes a^3) = -\shift a^1\otimes a^2\otimes\shift a^3.\hfill\qedhere\]
\end{enumerate}
\end{proof}

Next note that 
\begin{equation*}
\begin{split}
\shift^{-1}\overline{V}\otimes_A\shift^{-1}\overline{V}\cong& J\otimes\shift^{-1}\Phi \otimes J\otimes\shift^{-1} \Phi \otimes J \\
&\oplus \shift^{-1}\Phi \otimes J\otimes \shift^{-1}\Phi \otimes J
\oplus J\otimes\shift^{-1}\Phi \otimes \shift^{-1}\Phi \otimes J
\oplus J\otimes\shift^{-1}\Phi \otimes J\otimes \shift^{-1}\Phi \\
&\oplus J\otimes\shift^{-1}\Phi \otimes \shift^{-1}\Phi
\oplus \shift^{-1}\Phi \otimes J\otimes \shift^{-1}\Phi
\oplus \shift^{-1}\Phi \otimes \shift^{-1}\Phi \otimes J\\
&\oplus \shift^{-1}\Phi \otimes \shift^{-1}\Phi
\end{split}
\end{equation*}
We denote by $\mathtt{p}^0_i$ the projection onto the direct sum of subspaces with $i$ tensor factors.

\begin{lem}\label{claims131415}
\begin{enumerate}[(i)]
\item\label{claim15} The projection to 
\[\shift J\otimes\Phi \otimes \Phi 
\oplus \Phi \otimes \shift J\otimes \Phi \oplus    \Phi \otimes \Phi \otimes\shift J\]
of $\mu_3^{0, \mathbb{L}}$  and $(\shift\otimes\shift\otimes\shift)\mathtt{p}^0_3(\rho^{-1}\otimes_A\rho^{-1})\partial_1\shift^{-1}\vert_{\Phi}$ differ by a factor of $-1$.
\item\label{claim14} The projection to 
\[\Phi \otimes \shift J\otimes \Phi \otimes \shift J   
\oplus \shift J\otimes\Phi \otimes \Phi \otimes \shift J\oplus  \shift J\otimes\Phi \otimes \shift J\otimes \Phi \]
 of $\mu_4^{0, \mathbb{L}}$  and $(\shift\otimes\shift\otimes\shift\otimes\shift)\mathtt{p}^0_4(\rho^{-1}\otimes_A\rho^{-1})\partial_1\shift^{-1}\vert_{\Phi}$ are equal.
\item\label{claim13} The projection to $\shift J \otimes \Phi \otimes \shift J \otimes \Phi \otimes \shift J$ of $\mu_5^{0, \mathbb{L}}$  and $(\shift\otimes\shift\otimes\shift\otimes\shift\otimes\shift)\mathtt{p}^0_5(\rho^{-1}\otimes_A\rho^{-1})\partial_1 \shift^{-1}\vert_{\Phi} $ differ by a factor $-1$.\\
\end{enumerate}
\end{lem}
\begin{proof}
\begin{enumerate}[(i)]
\item Consider the diagram
\begin{equation*}
\begin{tikzcd}
\Phi\ar{rr}{\mu_3^{0,\mathbb{L}}}\ar{d}{\shift^{-1}} &&\shift J\otimes \Phi \otimes  \Phi \oplus \Phi \otimes \shift J \otimes \Phi \oplus \Phi  \otimes  \Phi\otimes\shift J\\
\shift^{-1}\Phi\ar{rr}{-\mathtt{p}^0_3(\rho^{-1}\otimes_A\rho^{-1})\partial_1}&&  J \otimes \shift^{-1} \Phi\otimes \shift^{-1} \Phi \oplus  \shift^{-1} \Phi\otimes J\otimes \shift^{-1} \Phi\oplus  \shift^{-1} \Phi\otimes J\otimes \shift^{-1} \Phi\ar{u}{\shift\otimes\shift\otimes\shift}
\end{tikzcd}
\end{equation*}
This commutes by the explicit description of $\mu_3^0$ in Lemma \ref{mu26}, since 
\begin{equation*}
\begin{split}
(\shift&\otimes\shift\otimes\shift)\mathtt{p}^0_3(\rho^{-1}\otimes_A\rho^{-1})\partial_1\shift^{-1} (\varphi) \\
&= (\shift\otimes\shift\otimes\shift)(\varphi^1\otimes \shift^{-1} \varphi^2 \varphi^3 \shift^{-1} \varphi^4 \varphi^5 +\varphi^1\shift^{-1} \varphi^2 \otimes \varphi^3 \otimes \shift^{-1} \varphi^4 \varphi^5+ \varphi^1\shift^{-1} \varphi^2 \varphi^3 \shift^{-1} \varphi^4\otimes \varphi^5) \\
&= -\shift \varphi^1\otimes \varphi^2\varphi^3 \varphi^4\varphi^5 +\varphi^1 \varphi^2\otimes \shift\varphi^3 \otimes \varphi^4\varphi^5 -\varphi^1 \varphi^2\varphi^3 \otimes\varphi^4\otimes\shift \varphi^5.
\end{split}
\end{equation*}
\item Consider the diagram
\begin{equation*}
\begin{tikzcd}
\Phi\ar{rr}{\mu_4^{0,\mathbb{L}}}\ar{d}{\shift^{-1}} &&\Phi \otimes \shift J\otimes \Phi \otimes \shift J   
\oplus \shift J\otimes\Phi \otimes \Phi \otimes \shift J\oplus     \shift J\otimes\Phi \otimes \shift J\otimes \Phi \\
\shift^{-1}\Phi\ar{rr}{(\rho^{-1}\otimes_A\rho^{-1})\partial_1}&& \shift^{-1} \Phi\otimes J \otimes\shift^{-1} \Phi\otimes J   \oplus J \otimes \shift^{-1} \Phi\otimes \shift^{-1} \Phi\otimes J \oplus    J \otimes \shift^{-1} \Phi\otimes J\otimes \shift^{-1} \Phi .\ar{u}{\shift\otimes\shift\otimes\shift\otimes\shift}
\end{tikzcd}
\end{equation*}
This commutes by the explicit description of $\mu_4^0$ in Lemma \ref{mu26}, since 
\begin{equation*}
\begin{split}
(\shift&\otimes\shift\otimes\shift\otimes\shift)(\rho^{-1}\otimes_A\rho^{-1})\partial_1\shift^{-1} (\varphi) \\
&= (\shift\otimes\shift\otimes\shift\otimes\shift)(\varphi^1\shift^{-1} \varphi^2\otimes \varphi^3 \otimes \shift^{-1} \varphi^4\otimes \varphi^5  +  \varphi^1\otimes \shift^{-1} \varphi^2 \varphi^3 \shift^{-1} \varphi^4\otimes \varphi^5+    \varphi^1\otimes \shift^{-1} \varphi^2 \otimes \varphi^3 \otimes \shift^{-1} \varphi^4 \varphi^5) \\
&=\varphi^1 \varphi^2\otimes\shift \varphi^3 \otimes\varphi^4\otimes\shift \varphi^5   -\shift \varphi^1\otimes \varphi^2\varphi^3 \varphi^4\otimes\shift \varphi^5+    \shift \varphi^1\otimes \varphi^2\otimes \shift\varphi^3 \otimes \varphi^4\varphi^5.
\end{split}
\end{equation*}
\item Consider the diagram
\begin{equation*}
\begin{tikzcd}
\Phi\ar{rr}{\mu_5^{0,\mathbb{L}}}\ar{d}{\shift^{-1}} && \shift J\otimes \Phi \otimes \shift J\otimes \Phi \otimes \shift J\\
\shift^{-1}\Phi\ar{rr}{-\mathtt{p}^0_5(\rho^{-1}\otimes_A\rho^{-1})\partial_1}&& J \otimes \shift^{-1} \Phi\otimes J\otimes \shift^{-1} \Phi\otimes J.\ar{u}{\shift\otimes\shift\otimes\shift\otimes\shift\otimes\shift}
\end{tikzcd}
\end{equation*}
This commutes by the explicit description of $\mu_5^0$ in Lemma \ref{mu26} since 
\begin{equation*}
\begin{split}
(\shift\otimes\shift\otimes\shift\otimes\shift\otimes\shift)\mathtt{p}^0_5(\rho^{-1}\otimes_A\rho^{-1})\partial_1 \shift^{-1} (\varphi) &= 
(\shift\otimes\shift\otimes\shift\otimes\shift\otimes\shift)(\varphi^1\otimes \shift^{-1} \varphi^2\otimes \varphi^3\otimes \shift^{-1} \varphi^4\otimes \varphi^5) \\
&= \shift \varphi^1\otimes \varphi^2\otimes\shift \varphi^3\otimes \varphi^4\otimes\shift \varphi^5.\qedhere
\end{split}
\end{equation*}
\end{enumerate}
\end{proof}

\subsection{Well-behaved splittings}\label{splittings}
Consider a presentation of $A$ as $A=\Bbbk Q/I$ and denote by $Q_+$ the augmentation ideal of $\Bbbk Q$. As usual, $Q_1$ denotes the vector space of arrows. In this subsection, we will construct well behaved splittings of the natural projections $\Bbbk Q \twoheadrightarrow A$ and $I\twoheadrightarrow I/(IQ_++Q_+I)$ which are needed for the main theorem.

Observe that $A=\Bbbk Q/I$ is a filtered algebra, where the degree of each element is induced by path length, i.e. $a\in A_{\leq i}$ if $a=\alpha_1\cdots \alpha_i$ with $\alpha_j\in Q_1$ is an arrow for all $j=1,\dots,i$. Let $\mult\colon Q_1\otimes_\Bbbk A\twoheadrightarrow J$ be the multiplication map. Note that $\mult$ is compatible with the filtration, i.e. $\mult\colon Q_1\otimes_\Bbbk A_{\leq i}\twoheadrightarrow J_{\leq i+1}$. We construct a splitting of $\mult$ compatible with the filtration. We proceed by induction. Let $C_1=Q_1$ and choose $C_{i+1}$ to be a complement to $\mult(Q_1\otimes C_i)\cap A_{\leq i}$ in $\mult(Q_1\otimes C_i)$ spanned by paths. 

By construction, $J=\bigoplus_{i\geq 1} C_i$ and $J^2=\bigoplus_{i\geq 2} C_i$ as vector spaces. 

Choose $\xi_{1,i}\colon C_{i+1}\to Q_1\otimes C_i$ to be a splitting of the multiplication map sending a basis vector (which is a path) to $\alpha\otimes q$ for a path $q\in C_i$ and $\alpha\in Q_1$, and set $\xi_1=\bigoplus_{i\geq 1} \xi_{1,i}\colon J^2\to Q_1\otimes J$ using the decomposition $J^2=\bigoplus_{i\geq 2}C_i$.

We furthermore define $\xi_{n+1}$ inductively as the composition
\[J^2\stackrel{\xi_1}{\to} Q_1\otimes J\twoheadrightarrow Q_1\otimes J^2\stackrel{\id\otimes \xi_n}{\to} Q_1^{\otimes (n+1)}\otimes J,\]
where the middle map is the projection onto the second summand for the direct sum decomposition 
$Q_1\otimes J= Q_1\otimes (Q_1\oplus J^2)=Q_1^{\otimes 2}\oplus (Q_1\otimes J^2)$.

Set $C_0:=\mathbb{L}$. We now define $\theta\colon A\to \Bbbk Q$ on the components $C_i$ using the decomposition $A=\bigoplus_{i\geq 0} C_i$. We let $\theta_0\colon C_0=\mathbb{L}\hookrightarrow \Bbbk Q$ and $\theta_1\colon C_1=Q_1\hookrightarrow \Bbbk Q$ be the canonical inclusions 
and define $\theta_{i+1}\colon C_{i+1}\to Q_1^{\otimes (i+1)}$ 
such that the diagram
\[\begin{tikzcd}
C_{i+1}\arrow{d}[swap]{\xi_{1,i}}\arrow{r}{\theta_{i+1}}&Q_1^{\otimes (i+1)}\\
Q_1\otimes C_i\arrow{ru}[swap]{1\otimes \theta_i}
\end{tikzcd}
\]
commutes.
We furthermore denote by $r_i$ the projection from $A$ to $C_i$ and define $\pi_i\colon \Bbbk Q\to C_i$ as the composition $r_i\circ \pi$.

\begin{lem}\label{pithetasplitting}
The morphism $\theta=(\theta_i)_{i\geq 0}$ is a splitting of the canonical projection $\pi\colon \Bbbk Q \twoheadrightarrow A$. In other words, $\pi_i\theta_i=\id_{C_i}$ for $i\geq 0$, and $\pi_j\theta_i=0$ for $i\neq j$.
\end{lem}

\begin{proof}
We proceed by induction, the claim being obvious for $i=0,1$. Now assume the statement holds for $k\leq i$.
To show that $\pi_{i+1}\theta_{i+1}=\id_{C_{i+1}}$, we compute
\begin{equation*}
\begin{split}
\pi_{i+1}\theta_{i+1}& =r_{i+1} \pi\theta_{i+1}=r_{i+1}\mathrm{mult}(\id\otimes \pi)(\id\otimes\theta_{i})\xi_{1,i}\\
\end{split}
\end{equation*}
where the second equality uses the definition of $\theta_i$ and the fact that $\pi$ factors as $\mathrm{mult}(\id\otimes \pi)$.
By the inductive assumption that $\theta_i$ is annihilated by all $\pi_j$ for $j\neq i$, hence we compute
\begin{equation*}
\begin{split}
\pi_{i+1}\theta_{i+1}&=r_{i+1}\mathrm{mult}(\id\otimes \pi)(\id\otimes\theta_{i})\xi_{1,i}= r_{i+1}\mathrm{mult}(\id\otimes \pi_i)(\id\otimes\theta_{i})\xi_{1,i}
= r_{i+1}\mathrm{mult}\xi_{1,i}.
\end{split}
\end{equation*}
Since $\mathrm{mult}\xi_{1,i} = \id_A$, we see that $r_{i+1}\mathrm{mult}\xi_{1,i}$ restricted to $C_{i+1}$ is the identity on $C_{i+1}$ as claimed.
\end{proof}

The proof of the main theorem will require two technical lemmas, which we will now state.

\begin{lem}\label{xithetacompat}
For $n\geq 2$, there is a commutative diagram
\begin{equation}
\label{xitheta}
\begin{tikzcd}
J^2\arrow{r}{\xi_{n-1}} &Q_1^{\otimes (n-1)}\otimes J\\
C_n\arrow[hookrightarrow]{u}\arrow{r}{\theta_n}&Q_1^{\otimes n}\arrow[hookrightarrow]{u}
\end{tikzcd}
\end{equation}
where the upwards arrows are split monomorphisms given by the decomposition $J=\bigoplus_{i\geq 1} C_i$.
\end{lem}

\begin{proof}
We proceed by induction on $n$. For $n=2$, 
\begin{equation*}
\begin{tikzcd}
J^2\arrow{r}{\xi_1}&Q_1\otimes J\\
C_2\arrow[hookrightarrow]{u}\arrow{r}{\theta_2}&Q_1\otimes Q_1\arrow[hookrightarrow]{u}
\end{tikzcd}
\end{equation*}
commutes since $\theta_2 = (1\otimes \theta_1)\xi_{1,i}$, the image of $\xi_{1,i}$ lies in the component $Q_1\otimes Q_1$ and $\theta_1$ corestricted to the subspace $Q_1$ of $\Bbbk Q$ is just the identity.

Expanding the diagram \eqref{xitheta} by using the definitions of $\xi_{n-1}$ and $\theta_n$, we obtain, for $n\geq 3$,
\[
\begin{tikzcd}
J^2\arrow{r}{\xi_1}&Q_1\otimes J\arrow[twoheadrightarrow]{r}&Q_1\otimes J^2\arrow{r}{\id\otimes \xi_{n-2}}&Q_1^{\otimes (n-1)}\otimes J\\
C_n\arrow[hookrightarrow]{u}\arrow{r}{\xi_{1,n-1}}&Q_1\otimes C_{n-1}\arrow[equals]{r}\arrow[hookrightarrow]{u} &Q_1\otimes C_{n-1}\arrow{r}{\id\otimes \theta_{n-1}}\arrow[hookrightarrow]{u}&Q_1^{\otimes n}\arrow[hookrightarrow]{u}
\end{tikzcd}
\]
in which the first and second square commute by definition and the rightmost square by the induction hypothesis. This completes the proof.
\end{proof}

\begin{lem}\label{rxipi}
Let $r\colon J\to Q_1$ denote the canonical projection. 
Then  for $n\geq 2$, 
\[(\id_{Q_1^{\otimes(n-1)}}\otimes r)\circ\xi_{n-1}\circ \pi = (\id_{Q_1^{\otimes(n-1)}}\otimes r)\circ\xi_{n-1}\circ \pi_n.\]
\end{lem}

\begin{proof}
Note that by definition $\pi=\sum_{j\geq 0}\pi_j$.
We show by induction on $n$ that the restriction of $\xi_{n-1}$ to $\bigoplus_{i=1}^{n-1} C_i \subset J$ is zero, extending $\xi_1$ to $C_1$ by zero for convenience. Indeed, for $n=2$, this is true by definition. Assume inductively that the statement is true for $k<n$. Then the image of $C_i$ under $\xi_{1,i}$ is contained in $Q_1\otimes C_{i-1}$, and hence by the induction hypothesis annihilated by $\id\otimes\xi_j$ for all $j\geq i-1$. Thus $C_i$ is annihilated by $\xi_j$ for $j\geq i$.

We further claim that the image of the restriction of $\xi_{n-1}$ to $\bigoplus_{i\geq n+1} C_i$ is contained in $Q_1^{\otimes(n-1)}\otimes J^2$ and thus annihilated by $\id_{Q_1^{\otimes(n-1)}}\otimes r$. We again proceed by induction on $n$. The image of $C_j$, for $j\geq 3$, under $\xi_{1,j}$ is contained in $Q_1\otimes C_{j-1}$, and, by definition, $C_{j-1}\subset J^2$ for $j-1\geq 2$. Assume inductively that  the claim holds for $k<n$. For $i\geq n+1$, the image of $C_i$ under 
$\xi_{1,i}$ is contained in $Q_1\otimes C_{i-1}$, and, by the inductive hypothesis the image of $C_{i-1}$ under $\xi_{n-2}$ is contained in $Q_1^{\otimes(n-2)}\otimes J^2$, thus the claim follows.

This implies that the only summand of $\xi_{n-1}\circ \pi$ not annihilated by $\id_{Q_1^{\otimes(n-1)}}\otimes r$ is $\xi_{n-1}\circ \pi_n$ and completes the proof.
\end{proof}

The following lemma, which states that a splitting of the natural projection $I\to I/(IQ_++Q_+I)$ can be chosen to be compatible with $\theta$, will be important in the proof of Theorem \ref{propermainthm}\eqref{mainthm1}.

\begin{lem}\label{goodzeta}
There exists a splitting  of the projection $I\to I/(IQ_++Q_+I)$ whose image is contained in $Q_1\otimes \image(\theta)$.
\end{lem}

\begin{proof}
By construction, we have vector space decompositions $Q_+=Q_1\oplus (Q_1\otimes Q_+)$ and $Q_+= \theta(J)\oplus I$.
Combining these, we obtain a vector space decomposition
\begin{equation*}\begin{split}
Q_+& =Q_1\oplus (Q_1\otimes Q_+)= Q_1\oplus \left(Q_1\otimes \theta(J)\right) \oplus \left(Q_1\otimes I\right)
=Q_1\oplus \left( Q_1\otimes \bigoplus_{i\geq 1}\theta(C_i) \right)\oplus \left(Q_1\otimes I\right).
\end{split}\end{equation*}

Let $\zeta'\colon I/(IQ_++Q_+I)\to I$ be a splitting of the projection $I\to I/(IQ_++Q_+I)$. Let $r\in I/(IQ_++Q_+I)$ and consider $\zeta'(r)$. Since $I$ is admissible, the image of $\zeta'$ is contained in $\left( Q_1\otimes \bigoplus_{i\geq 1}\theta(C_i) \right)\oplus \left(Q_1\otimes I\right)$ in the above decomposition. Since $Q_1\otimes I$ is in the kernel of the projection $I\to I/(IQ_++Q_+I)$, composing $\zeta'$ with the projection onto $\left( Q_1\otimes \bigoplus_{i\geq 1}\theta(C_i) \right)$ (and followed again by the inclusion) is also a splitting, whose image is indeed in the image of $ \id_{Q_1}\otimes \theta$.
\end{proof}

\subsection{The main result}\label{sec:mainthm}
We define 
\begin{itemize}
\item $\zeta \colon I/(IQ_++Q_+I) \to Q_1\otimes \image(\theta) \hookrightarrow Q_1\otimes Q_+$ as the corestriction of the splitting constructed in Lemma \ref{goodzeta};
\item the image $\Gamma$ of $\zeta$, which is a complement to $IQ_++Q_+I$ in $I$;
\item  the projection $q_n\colon Q_+  \twoheadrightarrow  Q_1^{\otimes n}$.
\end{itemize}

We further recall
\begin{itemize}
\item the isomorphism $\rho\colon A\otimes_{\mathbb{L}}\shift^{-1} \Phi\otimes_{\mathbb{L}} A\to \shift^{-1}\overline{V},$ given by $ a\otimes \shift^{-1}\varphi\otimes b\mapsto a\shift^{-1}(\varphi) b$;
\item the section $\theta$ of the projection map $\pi\colon \Bbbk Q\twoheadrightarrow A$ constructed before Lemma \ref{xithetacompat};
\item the projection $r_i\colon A\twoheadrightarrow C_i$.
\end{itemize}

We now restate Theorem \ref{mainthm} with this precise notation and the relevant gradings.

\begin{thm}\label{propermainthm}
Let $\mathfrak{A}=(A=\Bbbk Q/I,V)$ be a basic directed bocs. Let $R$ be its right algebra with standard modules $\Delta$. 
Then there is an $A_\infty$-structure on $\Ext^{\bullet}_R(\Delta,\Delta)$ such that 
\begin{enumerate}[(i)]
\item\label{mainthm1} for the morphism 
\[m_n\colon \Ext^1_R(\Delta,\Delta)^{\otimes n}\to \Ext^2_R(\Delta,\Delta),\] 
$\iota_n^{-1}m_n^\#$ can be identified with the  map 
\[ \shift^2(I/(IQ_++Q_+I)) \stackrel{\zeta \shift^{-2}}{\to} Q_1\otimes Q_+ \hookrightarrow Q_+  \stackrel{\shift^{\otimes n}q_n}{\twoheadrightarrow}  (\shift Q_1)^{\otimes n};\] 
\item\label{mainthm2} for the morphism 
\[m_n\colon \Ext^1_R(\Delta,\Delta)^{\otimes i}\otimes \radoperator_R(\Delta,\Delta)\otimes \Ext^1_R(\Delta,\Delta)^{\otimes j}\to \Ext^1_R(\Delta,\Delta)\]
$\iota_{(\shift Q_1)^\#,\dots, (\shift Q_1)^\#, \Phi^\#,(\shift Q_1)^\#,\dots,(\shift Q_1)^\#}^{-1}m_n^\#$
 can be identified with the map 
 \[\shift Q_1 \xrightarrow{\shift^{-1}}A\stackrel{\partial_0}{\to} \shift^{-1}\overline{V}\xrightarrow{\rho^{-1}} A\otimes \shift^{-1}\Phi \otimes A \xrightarrow{r_i\otimes 1\otimes r_j} C_i\otimes \shift^{-1}\Phi\otimes C_j 
\xrightarrow{\shift^{\otimes n}(\theta_i\otimes 1\otimes \theta_j)} (\shift Q_1)^{\otimes i}\otimes \Phi\otimes (\shift Q_1)^{\otimes j};\] 
\item\label{mainthm3} for the morphism 
\[m_n\colon \Ext^1_R(\Delta,\Delta)^{\otimes i}\otimes \radoperator_{R}(\Delta,\Delta)\otimes \Ext^1_R(\Delta,\Delta)^{\otimes j}\otimes \radoperator_{R}(\Delta,\Delta)\otimes \Ext^1_R(\Delta,\Delta)^{\otimes k}\to \radoperator_R(\Delta,\Delta)\] 
$\iota_{(\shift Q_1)^\#,\dots, (\shift Q_1)^\#, \Phi^\#,(\shift Q_1)^\#,\dots,(\shift Q_1)^\#, \Phi^\#,(\shift Q_1)^\#,\dots,(\shift Q_1)^\#}^{-1}m_n^\#$ can be identifed with the map\\
\hspace*{-7cm}
\begin{tikzcd}
\Phi\arrow[hookrightarrow]{r}{\shift^{-1}} &\shift^{-1}\overline{V}\arrow{r}{\partial_1} &\shift^{-1}\overline{V}\otimes_A \shift^{-1}\overline{V}\arrow[out=0, in=180]{dll}[description]{\rho^{-1}\otimes_A\rho^{-1}} \\
A\otimes \shift^{-1}\Phi\otimes A\otimes \shift^{-1}\Phi\otimes A
\arrow{r}{r_i\otimes 1\otimes r_j\otimes r_k} &C_i\otimes \shift^{-1}\Phi\otimes C_j \otimes \shift^{-1}\Phi\otimes C_k 
\arrow{r}{\shift^{\otimes n}(\theta_i\otimes 1\otimes \theta_j\otimes 1\otimes \theta_k)} &(\shift Q_1)^{\otimes i}\otimes \Phi\otimes (\shift Q_1)^{\otimes j}\otimes \Phi\otimes (\shift Q_1)^{\otimes k}.
\end{tikzcd}

\end{enumerate} 
\end{thm}
\begin{proof}[Proof of Theorem \ref{propermainthm}]

Recall the projective bimodule resolution \eqref{projbimodres} of $V$ and its structure of $A_\infty$-coalgebra object in $\mathrm{mod}\,A\otimes A^{\op}$. As explained in Section \ref{sec:projectiveresolution}, we are interested in the $A_\infty$-category $\Hom_{A\otimes A^{\op}}(\mathcal{P},\Hom_{\Bbbk}(\mathbb{L},\mathbb{L}))$ obtained from Lemma \ref{dualtocoalgebraobject}, whose multiplications we will denote by $\tilde m_n^i$, since its homology is isomorphic to $\Ext_{\mathfrak{A}}^{\bullet}(\mathbb{L},\mathbb{L})$.
Consider the isomorphisms
\begin{equation}\begin{split}\label{dualisos}
\Hom_{A\otimes A^{\op}}(A\otimes \Psi\otimes A, \Hom_\Bbbk(\mathbb{L},\mathbb{L}))&\cong\Hom_{\mathbb{L}\otimes \mathbb{L}^{\op}}(\Psi, \Hom_{\Bbbk}(\mathbb{L},\mathbb{L}))\\
&\cong \Hom_\mathbb{L}(\Psi\otimes \mathbb{L}, \mathbb{L})
\cong \Hom_\mathbb{L}(\Psi,\mathbb{L}).
\end{split}\end{equation}
of vector spaces. Notice that the $\mathbb{L}$-$\mathbb{L}$-bimodule structure on $\Hom_{A\otimes A^{\op}}(A\otimes \Psi\otimes A, \Hom_\Bbbk(\mathbb{L},\mathbb{L}))$ is induced by composition in $\Hom_\Bbbk(\mathbb{L},\mathbb{L})$ (cf. Definition \ref{directedbocs}\eqref{morofbocsmod}). Thus, there is an equality of vector spaces $e_\mathtt{j}\Hom_{A\otimes A^{\op}}(A\otimes \Psi\otimes A, \Hom_\Bbbk(\mathbb{L},\mathbb{L})) e_\mathtt{i}= \Hom_{A\otimes A^{\op}}(A\otimes \Psi\otimes A, \Hom_\Bbbk(\mathbb{L}e_\mathtt{i},\mathbb{L}e_\mathtt{j}))$. Tracing this through the adjunctions yields $e_\mathtt{j}\Hom_\mathbb{L}(\Psi,\mathbb{L})e_\mathtt{i} = \Hom_\mathbb{L}(\Psi e_\mathtt{i},\mathbb{L}e_\mathtt{j})$, thus we obtain 
\[\Hom_{A\otimes A^{\op}}(A\otimes \Psi\otimes A, \Hom_\Bbbk(\mathbb{L},\mathbb{L}))\cong \Psi^\flat\cong \Psi^\#\] as an $\mathbb{L}$-$\mathbb{L}$-bimodule, using Lemma \ref{flathash}. Similarly to Remark  \ref{compflat}, we see that this isomorphism translates the multiplications $\tilde m_\ell = (-1)^\ell \Hom_{A\otimes A^{\op}}(\mu_\ell,\Hom_{\Bbbk}(\mathbb{L},\mathbb{L}))\iota_\ell'$ to $ (-1)^\ell (\mu_\ell^{\mathbb{L}})^\#\iota_\ell.$

Therefore, the resolution $\Hom_{A\otimes A^{\op}}(\mathcal{P},\Hom_\Bbbk(\mathbb{L},\mathbb{L}))$ is isomorphic to

\begin{equation}\label{dualbar}
 0\to (\mathbb{L}\oplus \Phi)^\#\stackrel{\partial^1}{\to} J^\#\stackrel{\partial^2}{\to} (J\otimes J)^\#\stackrel{\partial^3}{\to} \dots
\end{equation}

whose differential is given by $\partial^i=-(d^{-i})^\#$, where $d^{-i}$ is the differential of the complex 

\begin{equation}\label{reducedbar}\dots \stackrel{d^{-3}}{\to} J\otimes J\stackrel{d^{-2}}{\to} J\stackrel{d^{-1}}{\to} \mathbb{L}\oplus \Phi\to 0\end{equation}

defined via $d^{-i}(a_1\otimes \dots\otimes a_n)=\sum_{i=1}^{n-1} (-1)^{i}a_1\otimes \dots\otimes a_{i}a_{i+1}\otimes \dots\otimes a_n$ for $a_1,\dots a_n\in J$. 

Note that $\partial^1\colon (\mathbb{L}\oplus \Phi)^\#\to J^\#$ is the zero map because $(A,V)$ is a regular bocs. Indeed, for a regular bocs, $\partial(J)\subseteq JV+VJ$ and thus, 
\[(d^{-1})^\#(f)(a\otimes x\otimes b)=(-1)^{|f|}f(d^{-1}(a\otimes x\otimes b))=(-1)^{|f|}af(\omega)xb-axf(\omega)b+af(\partial(x))b=0\]
for $a,b\in A, x\in J$.
Thus, in position $0$, the homology is obviously $(L\oplus \Phi)^\#$. Furthermore, in position $1$, this map being zero also implies that the homology of the complex can be identified with $(J/J^2)^\#$, which is isomorphic to $Q_1^\#$, as in the classical case of simple modules (see e.g. \cite{LPWZ09, HL05}). From position $2$ on, everything is as in the classical case, so in particular the homology in position $2$ can be identified with $(I/(IQ_++Q_+I))^\#$.

\begin{description}
\item[Claim 1] There is a vector space decomposition
\[J\otimes J=\lefteqn{\overbrace{\phantom{\image \delta^{-3}\oplus \Gamma}}^{\ker d^{-2}}}\image d^{-3}\oplus\underbrace{\Gamma\oplus \xi_1(J^2)}_{\subseteq Q_1\otimes J}.\]
Proof:
Of course, $\ker d^{-2}=\image d^{-3}\oplus \Gamma$. It suffices to prove that $\xi_1(J^2)$ is a complement to $\ker d^{-2}$. This follows from the fact that the map $d^{-2}\colon J\otimes J\to J^2\subset J$ is a split epimorphism with right inverse $-\xi_1$. Thus, $J\otimes J=\xi_1(J^2)\oplus \ker d^{-2}$. 

\item[Claim 2] Consider \eqref{dualbar}, the dual of the reduced bar complex. Then
\begin{enumerate}[(a)]
\item $\image\partial^2\cong(\xi_1(J^2))^\#$ and $\ker \partial^3\cong(\xi_1(J^2))^\#\oplus \Gamma^\#$. 
\item $(J\otimes J)^\#\cong(\xi_1(J^2))^\#\oplus \Gamma^\#\oplus (\image d^{-3})^\#$ 
\end{enumerate}
Proof: It is an elementary fact of linear algebra that for a sequence $U\stackrel{f}{\to} V\stackrel{g}{\to} W$ with $gf=0$ decomposing a component $V=(\image f)\oplus (\ker g/\image f) \oplus L$, the corresponding decomposition of $V^\#$ in the dual sequence is $V^\#\cong\underbrace{\image g^\#}_{\cong L^\#}\oplus \underbrace{(\ker g/\image f)^\#}_{\cong(\ker f^\#)/(\image g^\#)}\oplus (\image f)^\#$.  Applying this to the reduced bar complex the claim follows. 
\item[Claim 3] The projection  $(J\otimes J)^\# \twoheadrightarrow \Gamma^\#$ annihilates $(J^2\otimes J)^\#$.\\
Proof: By construction $\Gamma$ embeds into the direct summand $Q_1\otimes J$ of the decomposition $J\otimes J=(Q_1\otimes J)\oplus (J^2\otimes J)$. 
\end{description}
For $n\geq 1$, define $\bar\xi_{n}=\shift^{\otimes n+1} \xi_n \shift^{-1}\colon \shift J^2 \to (\shift Q_1)^{\otimes n}\otimes \shift J$, which is a morphism of degree $n$.
\begin{description}
\item[Claim 4] Let $h$ be a homotopy between $\id$ and $\mathbf{i}\mathbf{p}$ where $\mathbf{p}$ and $\mathbf{i}$ are the projection and inclusion of the homology into \eqref{dualbar}, the dual of the reduced bar complex, i.e. $\mathbf{i}\mathbf{p} - \id = \partial h + h\partial$. Then $h^1$ can be chosen to be zero and $h^2$ can be chosen to be $-\bar\xi_1^\#$.\\
Proof: That $h^1$ can be chosen to be zero follows from the fact that the complex equals its homology in position $0$. In position $1$, the differential $\partial^2$ is the dual of multiplication in $A$, and hence the homotopy can be chosen as $-\bar\xi_1^\#$. Indeed,
\[-\bar\xi_1^\#\partial_2=\bar\xi_1^\#(d^{-2})^\#=-(d^{-2}\bar\xi_1)^\#=\id_{J^\#}.\] 
\end{description}

Now we apply the construction by Kontsevich and So\u\i belman given in the proof of  Theorem \ref{kadeishvilitheorem} to compute the projection to degree $0$, $1$, and $2$ of the restriction of the $A_\infty$-structure on $\Ext_{\mathfrak{A}}(\mathbb{L},\mathbb{L})$ to degrees $0$ and $1$.
As $\tilde m_n$ has degree $2-n$, at most two of the input tensor factors can have degree $0$. Since it is the essential step for all the cases, we first discuss the case where all the input tensor factors have degree $1$. This is essentially the same as in \cite{LPWZ09, HL05}.

Recall that  $m_n = \mathbf{p}\lambda_n \mathbf{i}^{\otimes n}$ and that by Definition \ref{lambdakerneldef} that $\lambda_2$ is defined as the multiplication on the dual of the reduced bar complex and $\lambda_n$ is constructed inductively as 
\[\lambda_n=\sum_{\substack{\ell\neq 1\\j_1+\dots+j_\ell=n}}(-1)^{\sum_{i=1}^{\ell-1}(\ell-i)(j_i-1)}m_\ell(h\lambda_{j_1}\otimes \dots\otimes h\lambda_{j_\ell}).\]
We thus have
$\iota_n^{-1}m_n^\# = \iota_n^{-1}(i^{\otimes n})^\# \lambda_n^\# p^\#$ and 
\begin{equation}\label{lambdahash}
\lambda_n^\# = \sum_{\substack{\ell\neq 1\\j_1+\dots+j_\ell=n}} (-1)^{\sum_{i=1}^\ell(\ell-i)(j_i-1)}(-1)^{(n-\ell)\ell} (h\lambda_{j_1}\otimes \cdots \otimes h\lambda_{j_\ell})^\#\tilde m_\ell^\#.
\end{equation}

We need the following claim:
\begin{description}
\item[Claim 5]
\begin{enumerate}[(a)]
\item $h\lambda_n|_{((\shift Q_1)^\#)^{\otimes n}} \subseteq (\shift J^2)^\#$ for $n\geq 2$. 
\item $h\lambda_n|_{((\shift Q_1)^\#)^{\otimes n}}=h\lambda_2(h\lambda_1\otimes h\lambda_{n-1})|_{((\shift Q_1)^\#)^{\otimes n}}$.
\item $h\lambda_n|_{((\shift Q_1)^\#)^{\otimes (n-1)}\otimes (\shift J)^\#}=(-1)^{n-1}\bar\xi_{n-1}^\#\iota_n$.
\item $\lambda_n|_{((\shift Q_1)^\#)^{\otimes (n-1)}\otimes (\shift J)^\#} = (-1)^{n-2}(\id\otimes \bar\xi_{n-2})^\#\iota_n,$
\end{enumerate}
Proof: 
\begin{enumerate}[(a)]
\item Note that $\lambda_n(((\shift Q_1)^\#)^{\otimes n})\subseteq (\shift J\otimes \shift J)^\#$. Since $h^2$ is chosen to be $-\bar\xi_1^\#$, whose image is contained in $(\shift J^2)^\#$, the claim follows.

\item 
Since the codomain of $\lambda_n|_{((\shift Q_1)^\#)^{\otimes n}}$ has degree $2$, the degrees of the codomains of the $h\lambda_j$ in a nonzero term $\tilde m_\ell(h\lambda_{j_1}\otimes \cdots \otimes h\lambda_{j_\ell})$ have to sum to $\ell$. Since each such degree is non-negative, and $h^1=0$, the codomain of each $h\lambda_j$ has strictly positive degree, and this is only possible if each ends in degree $1$. The only 
$\tilde m_\ell^\#$ appearing in \eqref{lambdahash} therefore correspond to $\mu_\ell^{2,\mathbb{L}}$ such that $\mu_\ell$ has nonzero projection onto $\mathcal{P}^{-1}\otimes\cdots\otimes\mathcal{P}^{-1}$, which, by Lemma \ref{classicalcase} necessitates $\ell=2$.
Therefore, 
\[\lambda_n|_{((\shift Q_1)^\#)^{\otimes n}}=\sum_{ j_1+j_2 =n} (-1)^{j_1-1} \tilde m_2(h\lambda_{j_1}\otimes h\lambda_{j_2})|_{((\shift Q_1)^\#)^{\otimes n}}.\]
If $j_1>1$, $h\lambda_{j_1}\subseteq (\shift J^2)^\#$ by part (a) of this claim. In this case 
\[h\tilde m_2(h\lambda_{j_1}\otimes h\lambda_{j_2})({((\shift Q_1)^\#)^{\otimes n}})\subseteq h\lambda_2((\shift J^2)^\#\otimes \shift J^\#)= h((\shift J^2\otimes \shift J)^\#)=-\bar\xi_1^\#((\shift J^2\otimes \shift J)^\#)=0,\]
where the first equality uses the fact that $\lambda_2$, restricted to the degree $1$ part, is equal to $\iota_2\colon \shift J^\#\otimes \shift J^\#\to (\shift J\otimes \shift J)^\#$, and last equality uses the splitting of $J=Q_1\oplus J^2$ and the fact that the image of $\xi_1$ is contained in $Q_1\otimes J$. This proves the claim.

\item For $n=2$, this follows from the fact that the multiplication $\lambda_2$, restricted to the degree $1$ part, is equal to $\iota_2\colon \shift J^\#\otimes \shift J^\#\to (\shift J\otimes \shift J)^\#$, and $h$ in this case is
 $h^2=-\bar\xi_1^\#$, leading to $h\lambda_2= -\bar\xi_1^\#\iota_2$.
 
Notice that 
\begin{equation}\label{xicomp}
\begin{split}
\bar\xi_{n-1}^\# &= (\shift^{\otimes n}\xi_{n-1}\shift^{-1})^\# =  (\shift^{\otimes n}(\id\otimes \xi_{n-2})\xi_1\shift^{-1})^\#
= ((\shift\otimes (\shift^{\otimes n-1}\xi_{n-2}))(\shift\otimes\shift)^{-1}\bar \xi_1)^\#\\
&=-((\shift\otimes (\shift^{\otimes n-1}\xi_{n-2}))(\shift^{-1}\otimes \shift^{-1})\bar\xi_1)^\#
=(-1)^n((\id\otimes \bar\xi_{n-2})\bar\xi_1)^\#
=\bar\xi_1^\#(\id\otimes \bar\xi_{n-2})^\#.
\end{split}
\end{equation}

For $n>2$, using part (b), we then obtain 
\begin{equation*}
\begin{split}
h\lambda_n &=h\lambda_2(h\lambda_1\otimes h\lambda_{n-1})
=-\bar{\xi}_1^\#\iota_2(\id\otimes((-1)^{n-2}\bar\xi_{n-2}^\#\iota_{n-1}))\\
&=(-1)^{n-1}\bar\xi_1^\#(\id\otimes \bar\xi_{n-2})^\#\iota_n
= (-1)^{n-1} \bar\xi_{n-1}^\#\iota_n 
\end{split}
\end{equation*}
where the third equality follows from commutativity of 
\[
\begin{tikzcd}
((\shift Q_1)^\#)^{\otimes (n-1)}\otimes \shift J^\#\arrow{r}{\id\otimes \iota_{n-1}}\arrow{rd}[swap]{\iota_n} &(\shift Q_1)^\#\otimes ((\shift Q_1)^{\otimes (n-2)}\otimes \shift J)^\#\arrow{d}{\iota_2}\arrow{r}{\id\otimes \bar\xi_{n-2}^\#} &(\shift Q_1)^\#\otimes (\shift J^2)^\#\arrow{d}{\iota_2} \arrow{d}\\
&{}((\shift Q_1){\otimes (n-1)}\otimes \shift J)^\#\arrow{r}{(\id\otimes \bar\xi_{n-2})^\#} &(\shift Q_1\otimes \shift J^2)^\#\arrow{r}{\xi_1^\#} &(\shift J)^\#.&{}
\end{tikzcd}
\]

\item By part (c), and using \eqref{xicomp} for the last equality,  we obtain
\begin{equation*}\begin{split}
h\lambda_n|_{((\shift Q_1)^\#)^{\otimes (n-1)}\otimes (\shift J)^\#}&=(-1)^{n-1}\bar\xi_{n-1}^\#\iota_n = (-1)^{n-1}\bar\xi_1^\#(\id\otimes \bar\xi_{n-2})^\#\iota_n.
\end{split}
\end{equation*}

Since $\bar\xi_1^\#$ restricted to the images of $\lambda_n|_{((\shift Q_1)^\#)^{\otimes (n-1)}\otimes (\shift J)^\#}$ and $(\id\otimes \bar\xi_{n-2})^\#$ is a monomorphism, we deduce that
$\lambda_n|_{((\shift Q_1)^\#)^{\otimes (n-1)}\otimes (\shift J)^\#} = (-1)^{n-2}(\id\otimes \bar\xi_{n-2})^\#\iota_n$. 
\end{enumerate}
\end{description}

Recall from Lemma \ref{lem:lambdakernel} that the multiplication $m_n$ on homology is given by $\mathbf{p}\lambda_n \mathbf{i}^{\otimes n}$ where $\mathbf{p}$ is projection onto homology and $\mathbf{i}$ is the inclusion of the homology into the reduced dual bar complex in \eqref{reducedbar}. 

The homology of \eqref{reducedbar}  in position $-2$ is $\shift^2 (I/(IQ_++Q_+I))$. This can be embedded into $\shift J\otimes \shift J$ by first using $\zeta \colon I/(IQ_++Q_+I) \to Q_1\otimes \im\theta \subset Q_1\otimes Q_+$ as defined following Lemma \ref{goodzeta}, then defining $\bar \zeta = (\shift\otimes\shift) \zeta \shift^{-2}$ and finally taking the embedding to be the composition
\begin{equation}\label{phash}
\shift^2 (I/(IQ_++Q_+I))\stackrel{\bar\zeta}{\to} \shift Q_1\otimes\shift Q_+\stackrel{\id\otimes \bar \pi}\twoheadrightarrow \shift Q_1\otimes \shift J\hookrightarrow \shift J\otimes \shift J,
\end{equation}
where $\bar \pi = \shift \pi \shift^{-1}$.
This is indeed an embedding thanks to the kernel of the projection having zero intersection with the image of $\zeta$ by Lemma \ref{goodzeta}. Thus, in degree two, we choose $\mathbf{p}$ dual to the morphism given in \eqref{phash}. In degree one, $\mathbf{i}$ is dual to the projection $\bar r\colon \shift r \shift^{-1}\colon \shift J\twoheadrightarrow \shift Q_1$, that is $\mathbf{i}^\# = \bar r$.

We wish to compute the dual of $m_n$ restricted to $((\shift Q_1)^\#)^{\otimes n}$, for which
\begin{equation*}
\begin{split}
\iota_n^{-1}m_n^\# &= \iota_n^{-1}(\mathbf{i}^{\otimes n})^\#\lambda_n^\# \mathbf{p}^\#=(-1)^{n-2}(\mathbf{i}^\#)^{\otimes n}(\id\otimes \bar\xi_{n-2})\mathbf{p}^\#
\end{split}
\end{equation*}
by Claim 5(d).
Thus, $\iota_n^{-1}m_n^\# = (-1)^{n-2}\phi$ where $\phi$ is given by the composition
\[\phi\colon \shift^2 (I/(IQ_++Q_+I)) \stackrel{\bar\zeta}{\to}  \shift Q_1\otimes \shift Q_+ \stackrel{\id\otimes \bar \pi}{\twoheadrightarrow} \shift Q_1\otimes \shift J \stackrel{\id\otimes \bar\xi_{n-2}}{\to} (\shift Q_1)^{\otimes (n-1)}\otimes \shift J \stackrel{\id\otimes \bar r}{\twoheadrightarrow}  (\shift Q_1)^{\otimes n}.\]

By Lemma \ref{rxipi}, this is equal to the composition
\[\phi \colon \shift^2 (I/(IQ_++Q_+I)) \stackrel{\bar\zeta}{\to}  \shift Q_1\otimes \shift Q_+\stackrel{\id\otimes \bar \pi_{n-1}}{\twoheadrightarrow} \shift Q_1\otimes \shift C_{n-1} \stackrel{\id\otimes \bar\xi_{n-2}}{\to} (\shift Q_1)^{\otimes (n-1)}\otimes \shift J \stackrel{\id\otimes \bar r}{\twoheadrightarrow}  (\shift Q_1)^{\otimes n}\]
where again $\bar \pi_{n-1} = \shift \pi_{n-1} \shift^{-1}$.

Consider the diagram
\begin{equation}\label{xitotheta}
\begin{tikzcd}
\shift^2 (I/(IQ_++Q_+I))\arrow{r}{\bar\zeta}&     \shift Q_1\otimes \shift Q_+\arrow[twoheadrightarrow]{dr}{\id\otimes \bar \pi_{n-1}}&\shift Q_1\otimes \shift J\arrow{r}{\id\otimes \bar\xi_{n-2}}&(\shift Q_1)^{\otimes (n-1)}\otimes \shift J\arrow[twoheadrightarrow]{r}{\id\otimes \bar r}   &(\shift Q_1)^{\otimes n}\arrow[equals]{d}\\
&&\shift Q_1\otimes \shift C_{n-1}\arrow{r}{\id\otimes \bar\theta_{n-1}}\arrow[hookrightarrow]{u}&(\shift Q_1)^{\otimes n}\arrow[equals]{r}\arrow[hookrightarrow]{u}&(\shift Q_1)^{\otimes n},
\end{tikzcd}
\end{equation}
where the middle square commutes by  \eqref{xithetacompat} and $\bar\theta_{n-1} = \shift^{\otimes n-1}\theta_{n-1} \shift^{-1}$.

Thus $\phi$ is equal to the composition
\[
\shift^2 (I/(IQ_++Q_+I))\stackrel{\bar\zeta}{\to}\shift Q_1\otimes \shift Q_+\stackrel{\id\otimes \bar \pi_{n-1}}{\twoheadrightarrow} \shift Q_1\otimes \shift C_{n-1} \stackrel{\id\otimes \bar\theta_{n-1}}{\to} (\shift Q_1)^{\otimes n}.
\]
Now, $\bar\theta_{n-1} \bar\pi_{n-1}= \shift^{\otimes n-1}\theta_{n-1} \pi_{n-1} \shift^{-1}$. By construction, the image of $\bar\zeta$ is containend $\shift Q_1\otimes \shift \image(\theta)$. On the image of $\theta$, $\theta_{n-1} \pi_{n-1}$ equals projection onto $Q_1^{\otimes n-1}$. 

Recalling that the projection $Q_+  \twoheadrightarrow  Q_1^{\otimes n}$ is denoted by $q_n$ and setting $\bar q_n = \shift^{\otimes n}q_n\shift^{-1}$, this proves that $\phi$ is given by the composition 
\[\shift^2 (I/(IQ_++Q_+I)) \stackrel{\bar\zeta}{\to}  \shift Q_1\otimes \shift Q_+  \stackrel{\id\otimes \bar q_{n-1}}\twoheadrightarrow  (\shift Q_1)^{\otimes n}.\]
Finally 
\[(\id\otimes \bar q_{n-1}) \bar\zeta =( \id \otimes  \shift^{\otimes n-1}q_{n-1}\shift^{-1} ) (\shift \otimes \shift)\zeta \shift^{-2} = (-1)^n \shift^{\otimes n}(\id\otimes q_{n-1})\zeta \shift^{-2} =   (-1)^{n-2} \shift^{\otimes n}q_n \zeta\shift^{-2}\]
so $\iota_n^{-1}m_n^\# =  (-1)^{n-2}\phi = \shift^{\otimes n}q_n \zeta\shift^{-2}$
as claimed. This concludes the proof of Theorem \ref{propermainthm}\eqref{mainthm1}.

We now prove Theorem \ref{propermainthm}\eqref{mainthm2}.

Define the map $\psi_{ij}$ given by the composition
\[\shift Q_1 \stackrel{\shift^{-1}}{\hookrightarrow} J\stackrel{\partial_0}{\to} \shift^{-1}\overline{V}\stackrel{\rho^{-1}}{\to} A\otimes
\shift^{-1} \Phi \otimes A \xrightarrow{\pi_i\otimes 1\otimes \pi_j}  C_i\otimes \shift^{-1}\Phi\otimes C_j 
\xrightarrow{\shift^{\otimes n}(\theta_i\otimes 1\otimes \theta_j)} (\shift Q_1)^{\otimes i}\otimes \Phi\otimes (\shift Q_1)^{\otimes j}.\]

Our aim is to show that for the morphism $m_n$ restricted to $\Ext^1_R(\Delta,\Delta)^{\otimes i}\otimes \radoperator_R(\Delta,\Delta)\otimes \Ext^1_R(\Delta,\Delta)^{\otimes j}$, using the identifications $\Ext^1_R(\Delta,\Delta)\cong (\shift Q_1)^\#$ and $ \radoperator_R(\Delta,\Delta)\cong \Phi^\#$, as well as those of spaces with their double duals, the morphism $\iota_{(\shift Q_1)^\#, \dots, (\shift Q_1)^\#,\Phi^\#,(\shift Q_1)^\#,\dots,(\shift Q_1)^\#}^{-1}m_n^\#$
can be identified with $\psi_{ij}$. To improve readability, we will write $\iota_{\blacksquare}$ for  $\iota_{(\shift Q_1)^\#, \dots, (\shift Q_1)^\#,\Phi^\#,(\shift Q_1)^\#,\dots,(\shift Q_1)^\#}$. We thus compute 
$\iota_{\blacksquare}^{-1}m_n^\#.$

As before $m_n=\mathbf{p}\lambda_n \mathbf{i}^{\otimes n}$.  The relevant tensor factors of $\mathbf{i}^{\otimes n}$ in this case are dual to the projections $\bar r\colon \shift J\twoheadrightarrow\shift Q_1$ in all but the $(i+1)$st factor as before and the  projection of $\mathbb{L}\oplus \Phi$ onto $\Phi$ in the middle. Since the direct summand $\mathbb{L}$ does not feature in any computations, we omit it and interpret $\mathbf{i}$ simply as the identity on $\Phi^\#$. More precisely, 
\[\mathbf{i}^{\otimes n} = (\bar r^\#)^{\otimes i}\otimes \id\otimes (\bar r^\#)^{\otimes j}\colon (\shift Q_1^\#)^{\otimes i}\otimes \Phi^\#\otimes (\shift Q_1^\#)^{\otimes j}  {\longrightarrow}  (\shift J^\#)^{\otimes i}\otimes \Phi^\#\otimes (\shift J^\#)^{\otimes j}.\]
The projection $\mathbf{p}$ is dual to the inclusion of $\shift Q_1$ into $\shift J$. 

Moreover, as before, $\lambda_n = \sum_{\substack{\ell\neq 1\\j_1+\dots+j_\ell=n}} (-1)^{\sum_{i=1}^\ell(\ell-i)(j_i-1)}  \tilde m_\ell(h\lambda_{j_1}\otimes \cdots \otimes h\lambda_{j_\ell}),$ and we now prove several claims to iteratively describe these.

\begin{description} 
\item[Claim 6] $h\lambda_n|_{((\shift Q_1)^\#)^{\otimes i}\otimes \Phi^\#\otimes ((\shift Q_1)^\#)^{\otimes j}}=0$ for $n\geq 2$. \\
Proof: Observe that $\lambda_n$ is of degree $2-n$ while the domain is in degree $n-1$, hence the codomain of $\lambda_n|_{((\shift Q_1)^\#)^{\otimes i}\otimes \Phi^\#\otimes ((\shift Q_1)^\#)^{\otimes j}}$ is in degree $1$. Since $h^1=0$, the claim follows.

\item[Claim 7] In the expansion 
\[\lambda_n|_{((\shift Q_1)^\#)^{\otimes i}\otimes \Phi^\#\otimes ((\shift Q_1)^\#)^{\otimes j}}=\sum_{\substack{\ell\neq 1\\j_1+\dots+j_\ell=n}}(-1)^{\sum_{i=1}^\ell(\ell-i)(j_i-1)}\tilde m_\ell(h\lambda_{j_1}\otimes \dots\otimes h\lambda_{j_\ell}),\] the only $\tilde m_\ell$ that contribute are $\tilde m_2$ (in case $i=0$ or $j=0$) and $\tilde m_3$ (in case $i>0$ and $j>0$). \\
Proof: By Lemma \ref{property2}, the only $\mu_n^1$ having a non-zero projection  to ${}^{\small \textcircled{{\tiny 2}}}\mathcal{P}^{\otimes n}$ are $\mu_1^1, \mu_2^1, \mu_3^1$. As $\tilde m_\ell= (-1)^\ell\Hom_{A\otimes A^{\op}}(\mu_\ell,\Hom_{\Bbbk}(\mathbb{L},\mathbb{L}))\iota_\ell$, only $\tilde m_2$ and $\tilde m_3$ contribute in 
\[\lambda_n|_{(\shift Q_1^\#)^{\otimes i}\otimes \Phi^\#\otimes (\shift Q_1^\#)^{\otimes j}}=\sum_{\substack{\ell\neq 1\\j_1+\dots+j_\ell=n}}(-1)^{\sum_{i=1}^\ell(\ell-i)(j_i-1)}\tilde m_\ell(h\lambda_{j_1}\otimes \dots\otimes h\lambda_{j_\ell}).\] 
Inspecting the explicit presentation of $\mu_2^1$ and $\mu_3^1$ in Lemma \ref{mu26} using the fact that the isomorphism $\Hom_{A\otimes A^{\op}}(A\otimes \Psi\otimes A, \Hom_\Bbbk(\mathbb{L},\mathbb{L}))\cong \Psi^\#$ kills elements whose left- or rightmost term is in $J$ we see that $\mu_2^1$ only contributes for $i=0$ or $j=0$ and $\mu_3^1$ only contributes for $i>0$ and $j>0$.  
\item[Claim 8] 
\[\lambda_n|_{((\shift Q_1)^\#)^{\otimes i}\otimes \Phi^\#\otimes ((\shift Q_1)^\#)^{\otimes j}}=
\begin{cases} 
\tilde m_3(h\lambda_{i}\otimes \id \otimes h\lambda_{j}) & \text{ if $i\neq 0$ and $j\neq 0$}\\
(-1)^{i-1} \tilde m_2(h\lambda_{i}\otimes \id)& \text{ if $i\neq 0$ and $j= 0$}\\
\tilde m_2(\id \otimes h\lambda_{j})& \text{ if $i=0$ and $j\neq 0$}
\end{cases}.\]
Proof: This follows from Claims 6 and 7.
\end{description}

Observe that the $h\lambda_i$ appearing in Claim 8 have domain $((\shift Q_1)^\#)^{\otimes i}$, and were hence described in Claim 5(c). 

We wish to use Claim 8 to compute 
\begin{equation*}
\begin{split}
\iota_{\blacksquare}^{-1}m_n^\#
&= \iota_{\blacksquare}^{-1}(\mathbf{p}\lambda_n\mathbf{i}^{\otimes n})^\#  
= \iota_{\blacksquare}^{-1}(\mathbf{i}^{\otimes n})^\#\lambda_n^\#\mathbf{p}^\#.
\end{split}
\end{equation*}

 Since the $\tilde m_\ell$ appearing in Claim 8 have domain $(\shift J)^{\#}\otimes \Phi^\#\otimes (\shift J)^{\#}$, $(\shift J)^{\#}\otimes \Phi^\#$ and $ \Phi^\#\otimes (\shift J)^{\#}$ in the different cases, respectively, the relevant parts of $\iota_\ell$ for the computation of $\tilde m_\ell^\#$ as $(-1)^\ell\iota_\ell^\#\mu_\ell^{\mathbb{L}}$ are $\iota_{(\shift J)^{\#},\Phi^\#, (\shift J)^{\#}}$, $\iota_{(\shift J)^{\#}, \Phi^\#}$ and $ \iota_{\Phi^\#, (\shift J)^{\#}}$, respectively. Identifying the domains of their duals with their double duals, e.g. identifying the domain  $(\shift J \otimes \Phi\otimes \shift J)^{\#\#}$ of $\iota^\#_{(\shift J)^{\#},\Phi^\#, (\shift J)^{\#}}$ with $\shift J \otimes \Phi\otimes \shift J$, we obtain that the relevant components of the 
$\tilde m_\ell^\# $ are 
\[(-1)^3\iota_{(\shift J)^{\#}\otimes \Phi^\#\otimes (\shift J)^{\#}}^\#\mu_3^{1,\mathbb{L}}, 
\qquad (-1)^2\iota_{(\shift J)^{\#}\otimes \Phi^\#}^\#\mu_2^{1,\mathbb{L}} \quad \text{ and } \quad
(-1)^2\iota_{\Phi^\#\otimes (\shift J)^{\#}}^\#\mu_2^{1,\mathbb{L}},\] 
respectively, where we keep the projections of $\mu_\ell^{1,\mathbb{L}}$ to corresponding direct summands implicit.

We first consider the case $i\neq 0\neq j$, in which case 
\begin{equation*}
\begin{split}
\iota_{\blacksquare}^{-1}m_n^\#
&= \iota_{\blacksquare}^{-1}(\mathbf{i}^{\otimes n})^\#\lambda_n^\#\mathbf{p}^\# 
=\iota_{\blacksquare}^{-1}(\mathbf{i}^{\otimes n})^\#(\tilde m_3(h\lambda_{i}\otimes \id \otimes h\lambda_{j}))^\#\mathbf{p}^\#\\
&=(-1)^{i+j}\iota_{\blacksquare}^{-1}(\mathbf{i}^{\otimes n})^\#(h\lambda_{i}\otimes \id \otimes h\lambda_{j})^\#\tilde m_3^\#\mathbf{p}^\#\\
&=(-1)^{i+j+1}\iota_{\blacksquare}^{-1}(\mathbf{i}^{\otimes n})^\#(h\lambda_{i}\otimes \id \otimes h\lambda_{j})^\#\iota_{(\shift J)^{\#}\otimes \Phi^\#\otimes (\shift J)^{\#}}^\#\mu_3^{1,\mathbb{L}}\mathbf{p}^\#\\
&=(-1)^{n}\iota_{\blacksquare}^{-1}(\mathbf{i}^{\otimes n})^\#(h\lambda_{i}\otimes \id \otimes h\lambda_{j})^\#\iota_{(\shift J)^{\#}\otimes \Phi^\#\otimes (\shift J)^{\#}}^\#\mu_3^{1,\mathbb{L}}\mathbf{p}^\#.
\end{split}
\end{equation*}
Recall that by Claim 5(c) 
\[ (h\lambda_{i}\otimes \id \otimes h\lambda_{j})^\# = (-1)^{i+j}(\bar\xi_{i-1}^\#\iota_{i} \otimes \id\otimes \bar\xi_{j-1}^\#\iota_{j})^\#. \]

Moreover, Lemma \ref{rxipi} and commutativity of \eqref{xitotheta} imply $r^{\otimes i}\xi_{i-1}\pi= r^{\otimes i}\xi_{i-1}\pi_i=\theta_i \pi_i = \theta_i r_i \pi$, which implies $r^{\otimes i}\xi_{i-1} =\theta_i r_i $ since $\pi$ is an epimorphism. 

Setting $\bar r_i = \shift r_i \shift^{-1}$, consider the commutative diagram
\begin{equation*}
\begin{tikzcd}
((\shift J)^{\#}\otimes \Phi^\#\otimes (\shift J)^{\#})^\#\ar{d}{\iota_{(\shift J)^{\#}, \Phi^\#, (\shift J)^{\#}}^{-1}} \ar{rr}{(h\lambda_{i}\otimes \id \otimes h\lambda_{j})^\#}
&&(((\shift J)^\#)^{\otimes i}\otimes \Phi^\#\otimes ((\shift J)^\#)^{\otimes j})^\# \ar{r}{(\mathbf{i}^{\otimes n})^\#}\ar{d}{\iota_{((\shift J)^{\#})^{\otimes i}, \Phi^\#, ((\shift J)^{\#})^{\otimes j}}^{-1}}
& (((\shift Q_1)^\#)^{\otimes i}\otimes \Phi^\#\otimes ((\shift Q_1)^\#)^{\otimes j})^\#\ar{ddd}{\iota_{\blacksquare}^{-1}}
\\
(\shift J)^{\#\#}\otimes \Phi^{\#\#}\otimes (\shift J)^{\#\#}\ar{drr}[description]{(-1)^{i+j}\bar \xi_{i-1}^{\#\#}\otimes \id \otimes \xi_{j-1}^{\#\#}}\ar{rr}{(h\lambda_{i})^\#\otimes \id \otimes (h\lambda_{j})^\#}\ar{dd}{\bar r_i^{\#\#}\otimes \id \otimes\bar r_j^{\#\#}}
&& (((\shift J)^\#)^{\otimes i})^\#\otimes \Phi^\#\otimes (((\shift J)^\#)^{\otimes j})^\#\ar{d}{\iota_i^{-1}\otimes \id\otimes\iota_j^{-1}}& \\
&&((\shift J)^{\#\#})^{\otimes i}\otimes \Phi^{\#\#}\otimes ((\shift J)^{\#\#})^{\otimes j} \ar{dr}{((\bar r)^{\#\#})^{\otimes i}\otimes \id \otimes ((\bar r)^{\#\#})^{\otimes j}}&\\
(\shift C_i)^{\#\#}\otimes \Phi^{\#\#}\otimes (\shift C_i)^{\#\#}\ar{rrr}{(-1)^{i+j}\bar \theta_i^{\#\#}\otimes \id \otimes \bar\theta_j^{\#\#}}&&& ((\shift Q_1)^{\#\#})^{\otimes i}\otimes \Phi^{\#\#}\otimes ((\shift Q_1)^{\#\#})^{\otimes j}
\end{tikzcd}
\end{equation*}
where the middle triangle commutes by Claim 5(c), and the bottom triangle commutes by the considerations in the paragraph preceding it.

Using this,  we obtain
\begin{equation}\label{iotaiota}
\begin{split}
\iota_{\blacksquare}^{-1}m_n^\#&= (-1)^{n}\iota_{\blacksquare}^{-1}(\mathbf{i}^{\otimes n})^\#(h\lambda_{i}\otimes \id \otimes h\lambda_{j})^\#\iota_{(\shift J)^{\#}\otimes \Phi^\#\otimes (\shift J)^{\#}}^\#\mu_3^{1,\mathbb{L}}\mathbf{p}^\#\\
&=-(\bar \theta_i\otimes \id \otimes \bar\theta_j)(\bar r_i\otimes \id \otimes\bar r_j)\iota_{(\shift J)^{\#}, \Phi^\#, (\shift J)^{\#}}^{-1}\iota_{(\shift J)^{\#}\otimes \Phi^\#\otimes (\shift J)^{\#}}^\#\mu_3^{1,\mathbb{L}} \mathbf{p}^\#.\\
\end{split}
\end{equation}
Noticing that $\iota_{(\shift J)^{\#}, \Phi^\#, (\shift J)^{\#}}^{-1}\iota_{(\shift J)^{\#}\otimes \Phi^\#\otimes (\shift J)^{\#}}^\#$ is just the canonical isomorphism
$(\shift J \otimes \Phi\otimes \shift J)^{\#\#} \to (\shift J)^{\#\#}\otimes \Phi^{\#\#}\otimes (\shift J)^{\#\#}$, we can rewrite this to obtain
\[\iota_{\blacksquare}^{-1}m_n^\# =  -(\bar \theta_i\otimes \id \otimes \bar\theta_j)(\bar r_i\otimes \id \otimes\bar r_j)\mu_3^{1,\mathbb{L}} \mathbf{p}^\#.\]
By Lemma \ref{claims910}\eqref{claim10}, we further compute
\begin{equation*}
\begin{split}
\iota_{\blacksquare}^{-1}m_n^\# &=  -(\bar \theta_i\otimes \id \otimes \bar\theta_j)(\bar r_i\otimes \id \otimes\bar r_j)\mu_3^{1,\mathbb{L}} \mathbf{p}^\#\\
&=(\bar \theta_i\otimes \id \otimes \bar\theta_j)(\bar r_i\otimes \id \otimes\bar r_j)(\shift\otimes\shift\otimes\shift)\mathtt{p}^1_3\rho^{-1}\partial_0 \shift^{-1}\mathbf{p}^\#\\
&=(\bar \theta_i\otimes \id \otimes \bar\theta_j)(\bar r_i\otimes \id \otimes\bar r_j)(\shift\otimes\shift\otimes\shift)\mathtt{p}^1_3\rho^{-1}\partial_0 \shift^{-1}\mathbf{p}^\#.
\end{split}
\end{equation*}

Simplifying shifts finally yields
\begin{equation*}
\begin{split}
\iota_{\blacksquare}^{-1}m_n^\#&=(\bar \theta_i\otimes \id \otimes \bar\theta_j)(\bar r_i\otimes \id \otimes\bar r_j)(\shift\otimes\shift\otimes\shift)\mathtt{p}^1_3\rho^{-1}\partial_0 \shift^{-1}\mathbf{p}^\#\\
&=((\shift^{\otimes i} \theta_i \shift^{-1})\otimes \id \otimes (\shift^{\otimes j} \theta_j\shift^{-1}))(\shift r_i\shift^{-1}\otimes \id \otimes\shift r_j\shift^{-1})(\shift\otimes\shift\otimes\shift)\mathtt{p}^1_3\rho^{-1}\partial_0 \shift^{-1}\mathbf{p}^\#\\
&=(\shift^{\otimes i} \theta_i  r_i\shift^{-1}  )\otimes \id \otimes (\shift^{\otimes j} \theta_j r_j\shift^{-1})  (\shift\otimes\shift\otimes\shift)\mathtt{p}^1_3\rho^{-1}\partial_0 \shift^{-1}\mathbf{p}^\# \\
&=\shift^{\otimes n}( \theta_i \otimes \id \otimes\theta_j)( r_i\otimes \id \otimes r_j)\mathtt{p}^1_3\rho^{-1}\partial_0 \shift^{-1}\mathbf{p}^\#\\
&=\psi_{ij}.
\end{split}
\end{equation*}

We now consider the case $i\neq 0, j=0$. In both cases where one of $i$ or $j$ is zero, we will be using Lemma \ref{claims910}\eqref{claim9}, and it will be clear from context, which of the relevant direct summands of the projection $\mathtt{p}^1_2$ will be used.

In the same way as before, we see that
\begin{equation*}
\begin{split}
\iota_{\blacksquare}^{-1}m_n^\# 
&= \iota_{\blacksquare}^{-1}(\mathbf{i}^{\otimes n})^\#\lambda_n^\#\mathbf{p}^\# 
=(-1)^{i-1}\iota_{\blacksquare}^{-1}(\mathbf{i}^{\otimes n})^\#(\tilde m_2(h\lambda_{i}\otimes \id ))^\#\mathbf{p}^\#\\
&=(-1)^{i-1}\iota_{\blacksquare}^{-1}(\mathbf{i}^{\otimes n})^\#(h\lambda_{i}\otimes \id )^\#\tilde m_2^\#\mathbf{p}^\#
=(-1)^{i-1}\iota_{\blacksquare}^{-1}(\mathbf{i}^{\otimes n})^\#(h\lambda_{i}\otimes \id )^\#\iota_{(\shift J)^{\#}\otimes \Phi^\#}^\#\mu_2^{1,\mathbb{L}}\mathbf{p}^\#\\
&=(\bar \theta_i\otimes \id )(\bar r_i\otimes \id)\mu_2^{1,\mathbb{L}}\mathbf{p}^\#.\\
\end{split}
\end{equation*}
By Lemma \ref{claims910}\eqref{claim9},
\begin{equation*}
\begin{split}
\iota_{\blacksquare}^{-1}m_n^\# &=(\bar \theta_i\otimes \id )(\bar r_i\otimes \id)\mu_2^{1,\mathbb{L}}\mathbf{p}^\#
=(\bar \theta_i\otimes \id )(\bar r_i\otimes \id)(\shift\otimes\shift)\mathtt{p}^1_2\rho^{-1}\partial_0\shift^{-1}\\
&=( \shift^{\otimes i}\theta_i\otimes \id )( r_i \shift^{-1}\otimes \id)(\shift\otimes\shift)\mathtt{p}^1_2\rho^{-1}\partial_0\shift^{-1}\\
&=\shift^{\otimes n}( \theta_i\otimes \id )( r_i\otimes \id)\mathtt{p}^1_2\rho^{-1}\partial_0\shift^{-1}=\psi_{ij}.
\end{split}
\end{equation*}

Finally, we consider the case $i=0, j\neq 0$. Again, we compute
\begin{equation*}
\begin{split}
\iota_{\blacksquare}^{-1}m_n^\#
&= \iota_{\blacksquare}^{-1}(\mathbf{i}^{\otimes n})^\#\lambda_n^\#\mathbf{p}^\# 
=\iota_{\blacksquare}^{-1}(\mathbf{i}^{\otimes n})^\#(\tilde m_2(\id \otimes h\lambda_{j}))^\#\mathbf{p}^\#\\
&=\iota_{\blacksquare}^{-1}(\mathbf{i}^{\otimes n})^\#(\id \otimes h\lambda_{j})^\#\tilde m_2^\#\mathbf{p}^\#
=\iota_{\blacksquare}^{-1}(\mathbf{i}^{\otimes n})^\#(\id \otimes h\lambda_{j})^\#\tilde m_2^\#\mathbf{p}^\#\\
&=\iota_{\blacksquare}^{-1}(\mathbf{i}^{\otimes n})^\#(\id \otimes h\lambda_{j})^\#\iota_{\Phi^\#\otimes (\shift J)^{\#}}^\#\mu_2^{1,\mathbb{L}}\mathbf{p}^\#
=(-1)^{j-1}(\id \otimes \bar \theta_j)(\id \otimes \bar r_j)\mu_2^{1,\mathbb{L}}\mathbf{p}^\#\\
\end{split}
\end{equation*}
by similar arguments as above. Using Lemma \ref{claims910}\eqref{claim9} again and simplifying shifts, we deduce that 
\begin{equation*}
\begin{split}
\iota_{\blacksquare}^{-1}m_n^\#&=(-1)^{j-1}(\id \otimes \bar \theta_j)(\id \otimes \bar r_j)\mu_2^{1,\mathbb{L}}\mathbf{p}^\#
=(-1)^{j-1}(\id \otimes \shift^{\otimes j} \theta_j)(\id\otimes  r_j\shift^{-1})(\shift\otimes\shift)\mathtt{p}^1_2\rho^{-1}\partial_0\shift^{-1}\\
&= \shift^{\otimes n}(\id \otimes  \theta_j)(\id\otimes  r_j)\mathtt{p}^1_2\rho^{-1}\partial_0\shift^{-1}=\psi_{ij}.
\end{split}
\end{equation*}

We finally prove Theorem \ref{mainthm}\eqref{mainthm3}.
We define the map $\chi_{ijk}$ as the composition\\
\hspace*{-7cm}
\begin{tikzcd}
\Phi\arrow[hookrightarrow]{r}{\shift^{-1}} &\shift^{-1}\overline{V}\arrow{r}{\partial_1} &\shift^{-1}\overline{V}\otimes_A \shift^{-1}\overline{V}
\arrow[out=0, in=180]{dll}[description]{\rho^{-1}\otimes_A\rho^{-1}}\\
A\otimes \shift^{-1}\Phi\otimes A\otimes \shift^{-1}\Phi\otimes A
\arrow{r}{ r_i\otimes 1\otimes  r_j\otimes  r_k} &C_i\otimes \shift^{-1}\Phi\otimes C_j \otimes \shift^{-1}\Phi\otimes C_k 
\arrow{r}{\shift^{\otimes n}(\theta_i\otimes 1\otimes \theta_j\otimes 1\otimes \theta_k)} &(\shift Q_1)^{\otimes i}\otimes \Phi\otimes (\shift Q_1)^{\otimes j}\otimes \Phi\otimes (\shift Q_1)^{\otimes k}.
\end{tikzcd}

\begin{description}
\item[Claim 9] In the expansion 
\[\lambda_n|_{((\shift Q_1)^\#)^{\otimes i}\otimes \Phi^\#\otimes ((\shift Q_1)^\#)^{\otimes j}\otimes \Phi^\#\otimes ((\shift Q_1)^\#)^{\otimes k}}=\sum_{\substack{\ell\neq 1\\j_1+\dots+j_\ell=n}}(-1)^{\sum_{i=1}^\ell(\ell-i)(j_i-1)}\tilde m_\ell(h\lambda_{j_1}\otimes \dots\otimes h\lambda_{j_\ell}),\]
the only $\tilde m_\ell$ that contribute are those for $\ell=2,\dots, 5$. \\
Proof: By Lemma \ref{property2}, the only $\mu_n^0$ having a non-zero projection to ${}^{\small \textcircled{{\tiny 2}}}\mathcal{P}^{\otimes n}$ are $\mu_2^0, \dots, \mu_5^0$. Again, since $\tilde m_\ell= (-1)^\ell\Hom_{A\otimes A^{\op}}(\mu_\ell,\Hom_{\Bbbk}(\mathbb{L},\mathbb{L}))\iota_\ell$, only $\tilde m_2, \dots, \tilde m_5$ contribute in 
\[\lambda_n|_{(\shift Q_1^\#)^{\otimes i}\otimes \Phi^\#\otimes (\shift Q_1^\#)^{\otimes j}\otimes \Phi^\#\otimes ((\shift Q_1)^\#)^{\otimes k}}=\sum_{\substack{\ell\neq 1\\j_1+\dots+j_\ell=n}}(-1)^{\sum_{i=1}^\ell(\ell-i)(j_i-1)}\tilde m_\ell(h\lambda_{j_1}\otimes \dots\otimes h\lambda_{j_\ell}).\]

\item[Claim 10] 
\[\lambda_n|_{((\shift Q_1)^\#)^{\otimes i}\otimes \Phi^\#\otimes ((\shift Q_1)^\#)^{\otimes j}}=
\begin{cases} 
\tilde m_5(h\lambda_{i}\otimes \id \otimes h\lambda_{j}  \otimes \id \otimes h\lambda_{k}) & \text{ if $i\neq 0, j\neq 0$ and $k\neq 0$}\\
(-1)^{i+j}\tilde m_4(h\lambda_{i}\otimes \id \otimes h\lambda_{j}  \otimes \id ) & \text{ if $i\neq 0, j\neq 0$ and $k= 0$}\\
(-1)^{i-1} \tilde m_4(h\lambda_{i}\otimes \id \otimes \id \otimes h\lambda_{k}) & \text{ if $i\neq 0, j= 0$ and $k\neq 0$}\\
\tilde m_4(\id \otimes h\lambda_{j}  \otimes \id \otimes h\lambda_{k}) & \text{ if $i= 0, j\neq 0$ and $k\neq 0$}\\
\tilde m_3(h\lambda_{i}\otimes \id \otimes \id ) & \text{ if $i\neq 0$ and $ j= k= 0$}\\
(-1)^{j-1}\tilde  m_3(\id \otimes h\lambda_{j}  \otimes \id) & \text{ if $i= 0, j\neq 0$ and $k= 0$}\\
\tilde m_3(\id \otimes \id \otimes h\lambda_{k}) & \text{ if $i= j= 0$ and $k\neq 0$}\\
\tilde m_2 &\text{ if $i=j=k=0$.}
\end{cases}.\]
Proof: This follows from Claims 6 and 9. 
\end{description}

As before, $m_n = \mathbf{p}\lambda_n \mathbf{i}^{\otimes n}$, and similarly to part (ii), 
\[\mathbf{i}^{\otimes n} = (\bar r^\#)^{\otimes i}\otimes \id\otimes (\bar r^\#)^{\otimes j}\otimes \id\otimes (\bar r^\#)^{\otimes k}\colon (\shift Q_1^\#)^{\otimes i}\otimes \Phi^\#\otimes (\shift Q_1^\#)^{\otimes j} \otimes \Phi^\#\otimes (\shift Q_1^\#)^{\otimes k}  {\longrightarrow}  (\shift J^\#)^{\otimes i}\otimes \Phi^\#\otimes (\shift J^\#)^{\otimes j}\otimes \Phi^\#\otimes (\shift J^\#)^{\otimes k}.\]

We set $\iota_{\blacksquare}^{-1}:=\iota_{(\shift Q_1)^\#,\dots, (\shift Q_1)^\#, \Phi^\#,(\shift Q_1)^\#,\dots,(\shift Q_1)^\#, \Phi^\#,(\shift Q_1)^\#,\dots,(\shift Q_1)^\#}^{-1}$ and again treat the cases separately. In each case, one part of Lemma \ref{claims131415} will be used, and we will now keep the relevant projections $\mathtt{p}^0_i$ entirely implicit.

We first consider the generic case of $i\neq 0,j\neq 0,k\neq 0$. In this case, we compute 

\begin{equation*}
\begin{split}
\iota_{\blacksquare}^{-1}m_n^\# 
&= \iota_{\blacksquare}^{-1}(\mathbf{i}^{\otimes n})^\#\lambda_n^\#\mathbf{p}^\# \\
&=\iota_{\blacksquare}^{-1}(\mathbf{i}^{\otimes n})^\#(\tilde m_5(h\lambda_{i}\otimes \id \otimes h\lambda_{j}\otimes \id \otimes h\lambda_{k}))^\#\mathbf{p}^\#\\
&=(-1)^{i+j+k-1}\iota_{\blacksquare}^{-1}(\mathbf{i}^{\otimes n})^\#(h\lambda_{i}\otimes \id \otimes h\lambda_{j}\otimes \id \otimes h\lambda_{k}))^\#\tilde m_5^\#\mathbf{p}^\#\\
&=(-1)^{n-1}\iota_{\blacksquare}^{-1}(\mathbf{i}^{\otimes n})^\#(h\lambda_{i}\otimes \id \otimes h\lambda_{j}\otimes \id \otimes h\lambda_{k})^\#\tilde m_5^\#\mathbf{p}^\#\\
&=(-1)^{n}\iota_{\blacksquare}^{-1}(\mathbf{i}^{\otimes n})^\#(h\lambda_{i}\otimes \id \otimes h\lambda_{j}\otimes \id \otimes h\lambda_{k})^\#\iota_5^\#\mu_5^{\mathbb{L}}\mathbf{p}^\#.
\end{split}
\end{equation*}
Similarly to above, we obtain 
\[\iota_{\blacksquare}^{-1}(\mathbf{i}^{\otimes n})^\#(h\lambda_{i}\otimes \id \otimes h\lambda_{j}\otimes \id \otimes h\lambda_{k})^\#= (-1)^{i+j+k-1}(\bar \theta_i\otimes \id \otimes \bar\theta_j\otimes \id \otimes \bar\theta_k)(\bar r_i\otimes \id \otimes\bar r_j\otimes \id \otimes \bar r_k)(\iota_{(\shift J)^{\#}, \Phi^\#, (\shift J)^{\#},\Phi^\#, (\shift J)^{\#} }^{-1})^{-1}\]
and putting this together, while using Lemma \ref{claims131415}\eqref{claim13} and cancelling occurrences of $\iota$ as in \eqref{iotaiota}, we obtain
\begin{equation*}
\begin{split}
\iota_{\blacksquare}^{-1}m_n^\# &=-(\bar \theta_i\otimes \id \otimes \bar\theta_j\otimes \id \otimes \bar\theta_k)(\bar r_i\otimes \id \otimes\bar r_j\otimes \id \otimes \bar r_k)\mu_5^{\mathbb{L}}\mathbf{p}^\#\\
&=(\bar \theta_i\otimes \id \otimes \bar\theta_j\otimes \id \otimes \bar\theta_k)(\bar r_i\otimes \id \otimes\bar r_j\otimes \id \otimes \bar r_k)(\shift\otimes\shift\otimes\shift\otimes\shift\otimes\shift)(\rho^{-1}\otimes_A\rho^{-1})\partial_1 \shift^{-1}\mathbf{p}^\#\\
&=(\shift^{\otimes n})(\theta_i\otimes \id \otimes \theta_j\otimes \id \otimes \theta_k)( r_i\otimes \id \otimes r_j\otimes \id \otimes r_k)(\rho^{-1}\otimes_A\rho^{-1})\partial_1 \shift^{-1}\mathbf{p}^\#\\
&=\chi_{ijk}.
\end{split}
\end{equation*}

Secondly, we treat the case of $i\neq 0,j\neq 0,k= 0$, which, using Lemma \ref{claims131415}\eqref{claim14}, yields
\begin{equation*}
\begin{split}
 \iota_{\blacksquare}^{-1}m_n^\# 
&=(-1)^{i+j}\iota_{\blacksquare}^{-1}(\mathbf{i}^{\otimes n})^\#(\tilde m_4(h\lambda_{i}\otimes \id \otimes h\lambda_{j} \otimes  \id) )^\#\mathbf{p}^\#\\
&=(-1)^{i+j}\iota_{\blacksquare}^{-1}(\mathbf{i}^{\otimes n})^\#(h\lambda_{i}\otimes \id \otimes h\lambda_{j} \otimes  \id)^\#\tilde m_4^\#\mathbf{p}^\#\\
&=(-1)^{i+j}\iota_{\blacksquare}^{-1}(\mathbf{i}^{\otimes n})^\#(h\lambda_{i}\otimes \id \otimes h\lambda_{j} \otimes  \id)^\#\iota_4^\#\mu_4^{\mathbb{L}}\mathbf{p}^\#\\
&=(\bar \theta_i\otimes \id \otimes \bar\theta_j\otimes  \id )(\bar r_i\otimes \id\otimes \bar r_j \otimes \id )\mu_4^{\mathbb{L}}\mathbf{p}^\#\\
&= (\shift^{\otimes i}\theta_i r_i\shift^{-1}\otimes \id \otimes \shift^{\otimes j}\theta_j r_j\shift^{-1}\otimes  \id )\mu_4^{\mathbb{L}}\mathbf{p}^\#\\
&= (\shift^{\otimes i}\theta_i r_i\shift^{-1}\otimes \id \otimes \shift^{\otimes j}\theta_j r_j\shift^{-1}\otimes  \id )(\shift\otimes\shift\otimes\shift\otimes\shift)(\rho^{-1}\otimes_A\rho^{-1})\partial_1 \shift^{-1}\mathbf{p}^\#\\
&= (\shift^{\otimes n})(\theta_i\otimes \id \otimes \theta_j\otimes \id )( r_i\otimes \id \otimes r_j\otimes \id )(\rho^{-1}\otimes_A\rho^{-1})\partial_1 \shift^{-1}\mathbf{p}^\#\\
&=\chi_{ijk}
.\end{split}
\end{equation*}

Next, we treat the case of  $i\neq 0,j= 0,k\neq 0$, where, using Lemma \ref{claims131415}\eqref{claim14}, we obtain 
\begin{equation*}
\begin{split}
 \iota_{\blacksquare}^{-1}m_n^\# 
&= \iota_{\blacksquare}^{-1}(\mathbf{i}^{\otimes n})^\#\lambda_n^\#\mathbf{p}^\# \\
&=(-1)^{i-1}\iota_{\blacksquare}^{-1}(\mathbf{i}^{\otimes n})^\#(\tilde m_4(h\lambda_{i}\otimes \id \otimes  \id \otimes h\lambda_{k}))^\#\mathbf{p}^\#\\
&=(-1)^{i-1}\iota_{\blacksquare}^{-1}(\mathbf{i}^{\otimes n})^\#(h\lambda_{i}\otimes \id \otimes \id \otimes h\lambda_{k})^\#\tilde m_4^\#\mathbf{p}^\#\\
&=(-1)^{i-1}\iota_{\blacksquare}^{-1}(\mathbf{i}^{\otimes n})^\#(h\lambda_{i}\otimes \id \otimes \id \otimes h\lambda_{k})^\#\iota_4^\#\mu_4^{\mathbb{L}}\mathbf{p}^\#\\
&= (-1)^{k-1}(\bar \theta_i\otimes \id \otimes  \id \otimes \bar\theta_k)(\bar r_i\otimes \id \otimes \id \otimes \bar r_k)\mu_4^{\mathbb{L}}\mathbf{p}^\#\\
&= (-1)^{k-1}(\shift^{\otimes i}\theta_i r_i\shift^{-1}\otimes \id \otimes  \id \otimes \shift^{\otimes k}\theta_k r_k\shift^{-1})\mu_4^{\mathbb{L}}\mathbf{p}^\#\\
&= (-1)^{k-1}(\shift^{\otimes i}\theta_i r_i\shift^{-1}\otimes \id \otimes  \id \otimes \shift^{\otimes k}\theta_k r_k\shift^{-1})(\shift\otimes\shift\otimes\shift\otimes\shift)(\rho^{-1}\otimes_A\rho^{-1})\partial_1 \shift^{-1}\mathbf{p}^\#\\
&=(\shift^{\otimes n})(\theta_i\otimes \id \otimes \id \otimes \theta_k)( r_i\otimes \id \otimes \id \otimes r_k)(\rho^{-1}\otimes_A\rho^{-1})\partial_1 \shift^{-1}\mathbf{p}^\#\\
&=\chi_{ijk}
\end{split}
\end{equation*}

For the case $i= 0,j\neq 0,k\neq 0$, we, again using Lemma \ref{claims131415}\eqref{claim14}, compute
\begin{equation*}
\begin{split}
\iota_{\blacksquare}^{-1}m_n^\# 
&=\iota_{\blacksquare}^{-1}(\mathbf{i}^{\otimes n})^\#(\tilde m_4( \id \otimes h\lambda_{j} \otimes  \id \otimes h\lambda_{k}))^\#\mathbf{p}^\#\\
&=\iota_{\blacksquare}^{-1}(\mathbf{i}^{\otimes n})^\#( \id \otimes h\lambda_{j} \otimes  \id \otimes h\lambda_{k})^\#\tilde m_4^\#\mathbf{p}^\#\\
&=\iota_{\blacksquare}^{-1}(\mathbf{i}^{\otimes n})^\#( \id \otimes h\lambda_{j} \otimes  \id \otimes h\lambda_{k})^\#\iota_4^\#\mu_4^{\mathbb{L}}\mathbf{p}^\#\\
&= (-1)^{j+k}( \id \otimes \bar\theta_j\otimes  \id \otimes\bar \theta_k)( \id\otimes \bar r_j \otimes \id \otimes \bar r_k)\mu_4^{\mathbb{L}}\mathbf{p}^\#\\
&= (-1)^{j+k}( \id \otimes \shift^{\otimes j}\theta_j r_j\shift^{-1}\otimes  \id \otimes\shift^{\otimes k}\theta_k r_k\shift^{-1})(\shift\otimes\shift\otimes\shift\otimes\shift)(\rho^{-1}\otimes_A\rho^{-1})\partial_1 \shift^{-1}\mathbf{p}^\#\\
&=(\shift^{\otimes n})( \id \otimes \theta_j\otimes \id \otimes\theta_k)( \id \otimes r_j\otimes \id \otimes r_k)(\rho^{-1}\otimes_A\rho^{-1})\partial_1 \shift^{-1}\mathbf{p}^\#.\\
&=\chi_{ijk}
\end{split}
\end{equation*}

In case  $i\neq 0,j=0,k= 0$, the calculation, now using Lemma \ref{claims131415}\eqref{claim15}, yields
\begin{equation*}
\begin{split}
\iota_{\blacksquare}^{-1}m_n^\# 
&= \iota_{\blacksquare}^{-1}(\mathbf{i}^{\otimes n})^\#\lambda_n^\#\mathbf{p}^\# \\
&=\iota_{\blacksquare}^{-1}(\mathbf{i}^{\otimes n})^\#(\tilde m_3(h\lambda_{i}\otimes \id \otimes  \id ))^\#\mathbf{p}^\#\\
&=(-1)^{i-1}\iota_{\blacksquare}^{-1}(\mathbf{i}^{\otimes n})^\#(h\lambda_{i}\otimes \id \otimes \id )^\#\tilde m_3^\#\mathbf{p}^\#\\
&=(-1)^{i}\iota_{\blacksquare}^{-1}(\mathbf{i}^{\otimes n})^\#(h\lambda_{i}\otimes \id \otimes \id )^\#\iota_3^\#\mu_3^{\mathbb{L}}\mathbf{p}^\#\\
&= -(\bar \theta_i\otimes \id \otimes  \id )(\bar r_i\otimes \id \otimes \id )\mu_3^{\mathbb{L}}\mathbf{p}^\#\\
&= -(\shift^{\otimes i}\theta_i r_i\shift^{-1}\otimes \id \otimes  \id )\mu_3^{\mathbb{L}}\mathbf{p}^\#\\
&=(\shift^{\otimes i}\theta_i r_i\shift^{-1}\otimes \id \otimes  \id )(\shift\otimes\shift\otimes\shift)(\rho^{-1}\otimes_A\rho^{-1})\partial_1 \shift^{-1}\mathbf{p}^\#\\
&=(\shift^{\otimes n})(\theta_i\otimes \id \otimes \id )( r_i\otimes \id \otimes \id)(\rho^{-1}\otimes_A\rho^{-1})\partial_1 \shift^{-1}\mathbf{p}^\#\\
&=\chi_{ijk}
.\end{split}
\end{equation*}

For the case $i= 0,j\neq 0,k= 0$, we, again using Lemma \ref{claims131415}\eqref{claim15}, compute
\begin{equation*}
\begin{split}
 \iota_{\blacksquare}^{-1}m_n^\# 
&= \iota_{\blacksquare}^{-1}(\mathbf{i}^{\otimes n})^\#\lambda_n^\#\mathbf{p}^\# \\
&=(-1)^{j-1}\iota_{\blacksquare}^{-1}(\mathbf{i}^{\otimes n})^\#(\tilde m_3( \id \otimes h\lambda_{j} \otimes  \id ))^\#\mathbf{p}^\#\\
&=\iota_{\blacksquare}^{-1}(\mathbf{i}^{\otimes n})^\#( \id \otimes h\lambda_{j} \otimes  \id )^\#\tilde m_3^\#\mathbf{p}^\#\\
&=-\iota_{\blacksquare}^{-1}(\mathbf{i}^{\otimes n})^\#( \id \otimes h\lambda_{j} \otimes  \id )^\#\iota_3^\#\mu_3^{\mathbb{L}}\mathbf{p}^\#\\
&= (-1)^{j}( \id \otimes \bar\theta_j\otimes  \id )( \id\otimes \bar r_j \otimes \id )\mu_3^{\mathbb{L}}\mathbf{p}^\#\\
&=(-1)^{j}( \id \otimes \shift^{\otimes j}\theta_j r_j\shift^{-1}\otimes  \id )\mu_3^{\mathbb{L}}\mathbf{p}^\#\\
&=(-1)^{j+1}( \id \otimes \shift^{\otimes j}\theta_j r_j\shift^{-1}\otimes  \id )(\shift\otimes\shift\otimes\shift)(\rho^{-1}\otimes_A\rho^{-1})\partial_1 \shift^{-1}\mathbf{p}^\#\\
&=(\shift^{\otimes n})( \id \otimes \theta_j\otimes \id)( \id \otimes r_j\otimes \id)(\rho^{-1}\otimes_A\rho^{-1})\partial_1 \shift^{-1}\mathbf{p}^\#\\
&=\chi_{ijk}
.\end{split}
\end{equation*}

For the case $i= 0,j= 0,k\neq 0$, we again use Lemma \ref{claims131415}\eqref{claim15} to compute
\begin{equation*}
\begin{split}
\iota_{\blacksquare}^{-1}m_n^\# 
&=\iota_{\blacksquare}^{-1}(\mathbf{i}^{\otimes n})^\#(\tilde m_3( \id \otimes   \id \otimes h\lambda_{k}))^\#\mathbf{p}^\#\\
&=(-1)^{k-1}\iota_{\blacksquare}^{-1}(\mathbf{i}^{\otimes n})^\#( \id \otimes  \id \otimes h\lambda_{k})^\#\tilde m_3^\#\mathbf{p}^\#\\
&=(-1)^{k}\iota_{\blacksquare}^{-1}(\mathbf{i}^{\otimes n})^\#( \id \otimes  \id \otimes h\lambda_{k})^\#\iota_3^\#\mu_3^{\mathbb{L}}\mathbf{p}^\#\\
&= -( \id \otimes  \id \otimes\bar \theta_k)( \id\otimes  \id \otimes \bar r_k)\mu_3^{\mathbb{L}}\mathbf{p}^\#\\
&=( \id \otimes  \id \otimes\bar \theta_k)( \id\otimes  \id \otimes \bar r_k)(\shift\otimes\shift\otimes\shift)(\rho^{-1}\otimes_A\rho^{-1})\partial_1 \shift^{-1}\mathbf{p}^\#\\
&=( \id \otimes  \id \otimes\shift^{\otimes k}\theta_k r_k\shift^{-1})(\shift\otimes\shift\otimes\shift)(\rho^{-1}\otimes_A\rho^{-1})\partial_1 \shift^{-1}\mathbf{p}^\#\\
&=(\shift^{\otimes n})( \id \otimes \id \otimes\theta_k)( \id \otimes \id \otimes r_k)(\rho^{-1}\otimes_A\rho^{-1})\partial_1 \shift^{-1}\mathbf{p}^\#\\
&=\chi_{ijk}.\qedhere
\end{split}
\end{equation*}
\end{proof}

\begin{cor}\label{inverse}
Let $\mathfrak{A}=(A,V)$ be a basic directed bocs, let $R$ be its right algebra. Equip $\Ext^{\bullet}_R(\Delta,\Delta)$ with the $A_\infty$-structure  $(m_n)_{n\in \mathbb{N}}$ obtained in Theorem \ref{mainthm}. Applying the construction in Theorem \ref{KKOmain} to this $A_\infty$-algebra yields a bocs isomorphic to $\mathfrak{A}$.
\end{cor}

\begin{proof}
By Theorem \ref{Brzmain}, it suffices to show that the quotient $\mathcal{D}$ of $\mathbf{B}(\Ext^{\bullet}_R(\Delta,\Delta))^\#$ by the differential ideal generated by elements of negative degree, is isomorphic to $T_A(\overline{V})$. Equivalently, using the isomorphism in Lemma \ref{barcobariso}\eqref{duality1}, we can work with $\mathbf{\Omega}((\Ext^{\bullet}_R(\Delta,\Delta))^{\#})$. Using \eqref{bocscoreq}, 
we see that this has differential with components
\[(-1)^n(\shift^{\otimes n})^{-1} \iota_n^{-1}m_n^{\#}\shift.\]
Since by Lemma \ref{extisobocsiso} isomorphic $A_\infty$-structures on $\Ext^{\bullet}_R(\Delta,\Delta)$  yield isomorphic associated bocses and, by Lemma \ref{changesigns}, replacing the multiplications $m_n$ by $(-1)^nm_n$ produces an isomorphic $A_\infty$-structure on $\Ext^{\bullet}_R(\Delta,\Delta)$, we can equivalently consider $\mathbf{\Omega}((\Ext^{\bullet}_R(\Delta,\Delta))^{\#})$ to be equipped with differential $d$ given by
\[d_n=(\shift^{\otimes n})^{-1} \iota_n^{-1}m_n^{\#}\shift.\]

We first check that $\mathcal{D}_0$ is isomorphic to $A$. For this, note that 
\[\mathcal{D}_0\cong T(\shift^{-1}(\Ext^1(\Delta,\Delta)^\#))/(d(\shift^{-2}(\Ext^2(\Delta,\Delta)^\#))).\] By Theorem \ref{propermainthm} \eqref{mainthm1}, $\iota_n^{-1}m_n^{\#}$ in this case is given by
\[ \shift^2(I/(IQ_++Q_+I)) \stackrel{\zeta \shift^{-2}}{\to} Q_1\otimes Q_+ \hookrightarrow Q_+  \stackrel{\shift^{\otimes n}q_n}{\twoheadrightarrow}  (\shift Q_1)^{\otimes n};\]
hence $d_n$ in this degree is given by $q_n\zeta \shift^{-1}$, where again $\zeta \colon I/(IQ_++Q_+I) \to Q_1\otimes \image(\theta) \hookrightarrow Q_1\otimes Q_+$ is the corestriction of the splitting constructed in Lemma \ref{goodzeta} and $q_n\colon Q_+  \twoheadrightarrow  Q_1^{\otimes n}$ is the natural projection.
Thus, the ideal spanned by the image of $d$ in $\Bbbk Q$ is equal to $I$. 

On the other hand, it follows from the proof of the main theorem of \cite{KKO14} that $\mathcal{D}$ is a semifree dg algebra. Its degree $1$ part is the projective bimodule spanned by $\shift^{-1} \radoperator_R(\Delta,\Delta)\cong \shift^{-1} \Phi$. Thus, it suffices to check that the maps $\partial_0'$ and $\partial_1'$ on $\mathcal{D}$ coincide with the maps $\partial_0$ and $\partial_1$ we started with. 

The differential on $\mathcal{D}$ is induced by the following diagram
\[
\begin{tikzcd}
T(\shift^{-1}(\Ext^1(\Delta,\Delta)^\#))\arrow{d}\arrow{r}{d} &T(\shift^{-1}(\Ext^1(\Delta,\Delta)^\#))\arrow{d}\\
\mathcal{D}\arrow{r}{\partial'} &\mathcal{D}.
\end{tikzcd}
\]
where the vertical arrows are the quotient maps modulo the differential ideal spanned by the negative degree part. 

Note that $\partial_0'$ is determined by the restriction $\partial'$ to the image of $\shift^{-1}(\Ext^1(\Delta,\Delta)^\#)$, which we identify with $Q_1$. By Theorem \ref{propermainthm}(ii), we can identify $d_n  =(\shift^{\otimes n})^{-1} \iota_n^{-1}m_n^{\#}\shift$ restricted to $Q_1$  
with 
\[Q_1 \hookrightarrow A\stackrel{\partial_0}{\to} \shift^{-1}\overline{V}\xrightarrow{\rho^{-1}} A\otimes \shift^{-1}\Phi \otimes A \xrightarrow{\bigoplus_{i+j+1=n}\theta_ir_i\otimes 1\otimes \theta_jr_j} \bigoplus_{i+j+1=n}(Q_1)^{\otimes i}\otimes \shift^{-1}\Phi\otimes (Q_1)^{\otimes j}.\]
Hence the differential on $\mathcal{D}$ restricted to (the image of) $Q_1$ is given by 
\[Q_1 \hookrightarrow A\stackrel{\partial_0}{\to} \shift^{-1}\overline{V}\xrightarrow{\rho^{-1}} A\otimes \shift^{-1}\Phi \otimes A \xrightarrow{\bigoplus_{n \geq 1}\bigoplus_{i+j+1=n}\pi_i\theta_ir_i\otimes 1\otimes \pi_j\theta_jr_j} \bigoplus_{n \geq 1}\bigoplus_{i+j+1=n}C_i\otimes \shift^{-1}\Phi\otimes C_j.\]
Using $\pi_i\theta_i = \id_{C_i}$ (see Lemma \ref{pithetasplitting}), $\bigoplus_{n \geq 1}\bigoplus_{i+j+1=n}C_i\otimes \shift^{-1}\Phi\otimes C_j = A\otimes \shift^{-1}\Phi \otimes A $ and 
\[\bigoplus_{n \geq 1}\bigoplus_{i+j+1=n}r_i\otimes 1\otimes r_j = \id_{A\otimes \shift^{-1}\Phi \otimes A},\]
we can identify this with
\[Q_1 \hookrightarrow A\stackrel{\partial_0}{\to} \shift^{-1}\overline{V}\xrightarrow{\rho^{-1}} A\otimes \shift^{-1}\Phi \otimes A \]
as claimed.

Analogously, $\partial_1'$ is determined by the restriction of $d$ to $\shift^{-1} \Phi$. By Theorem \ref{propermainthm}(iii), we can identify $d_n$ 
restricted to $\Phi$ with
\begin{equation*}
\begin{split}
\shift^{-1} \Phi\hookrightarrow\shift^{-1}\overline{V}&\xrightarrow{\partial_1} \shift^{-1}\overline{V}\otimes_A \shift^{-1}\overline{V}\xrightarrow{\rho^{-1}\otimes_A\rho^{-1}} A\otimes \shift^{-1}\Phi\otimes A\otimes \shift^{-1}\Phi\otimes A\\
&\xrightarrow{\bigoplus_{i+j+k+2=n}(\theta_ir_i\otimes 1\otimes \theta_jr_j\otimes 1\otimes \theta_kr_k)} Q_1^{\otimes i}\otimes \Phi\otimes Q_1^{\otimes j}\otimes \shift^{-1}\Phi\otimes Q_1^{\otimes k}.
\end{split}
\end{equation*}
Summing over $n$, projecting down to $\mathcal{D}$ and using the same arguments as for $\partial_0'$, we see that $\partial_1'$ is indeed given by 
\[\shift^{-1} \Phi\hookrightarrow\shift^{-1}\overline{V}\xrightarrow{\partial_1} \shift^{-1}\overline{V}\otimes_A \shift^{-1}\overline{V}\xrightarrow{\rho^{-1}\otimes_A\rho^{-1}} A\otimes \shift^{-1}\Phi\otimes A\otimes \shift^{-1}\Phi\otimes A,
\]  
and hence coincides with $\partial_1$, as required.
\end{proof}

\begin{thm}
Let $\mathfrak{A}=(A,V)$ and $\mathfrak{B}=(B,W)$ be basic directed bocses  such that $R_\mathfrak{A}$ and $R_\mathfrak{B}$ are Morita equivalent as quasi-hereditary algebras. Then $\mathfrak{A}$ and $\mathfrak{B}$ are isomorphic. Moreover, there exists an isomorphism between $R_\mathfrak{A}$ and $R_\mathfrak{B}$, which restricts to an isomorphism between $A$ and $B$.
\end{thm}

\begin{proof}
Set $R=R_\mathfrak{A}$ and $S=R_\mathfrak{B}$. If $R$ and $S$ are Morita equivalent as quasi-hereditary algebras, then endowing $\Ext^{\bullet}_R(\Delta_R, \Delta_R)$ and $\Ext^{\bullet}_S(\Delta_S, \Delta_S)$ with the $A_\infty$-structure as in Theorem \ref{propermainthm} yields
isomorphic $A_\infty$-algebras.

From Corollary \ref{inverse}, it follows that $\mathfrak{A}=(A,V)$ and $\mathfrak{B}=(B,W)$ are isomorphic. Hence, $R$ and $S$ are also isomorphic via an isomorphism restricting to the given isomorphism between $A$ and $B$.
\end{proof}

\bibliographystyle{alpha}
\bibliography{publication}

\begin{thebibliography}{LPWZ09}

\bibitem[BB91]{BB91}
William~L. Burt and Michael Charles~Richard Butler.
\newblock {Almost split sequences for bocses}.
\newblock In {\em {Representations of finite-dimensional algebras ({T}sukuba
  1990)}}, volume~11 of {\em {CMS Conference Proceedings}}, pages 89--121.
  American Mathematical Society, Providence, RI, 1991.

\bibitem[Be{\u{\i}}78]{Bei78}
A.~A. Be{\u{\i}}linson.
\newblock Coherent sheaves on {${\bf P}^{n}$} and problems in linear algebra.
\newblock {\em Akademiya Nauk SSSR. Funktsional'ny{\u{\i}} Analiz i ego
  Prilozheniya}, 12(3):68--69, 1978.

\bibitem[Be{\u\i}84]{Bei84}
Alexander~A. Be{\u\i}linson.
\newblock {The derived category of coherent sheaves on {${\bf P}^n$}}.
\newblock {\em Selecta Mathematica Sovietica}, 3(3):233--237, 1984.
\newblock Selected translations.

\bibitem[BKK20]{BKK20}
Tomasz Brzezi{\'{n}}ski, Steffen Koenig, and Julian K{\"{u}}lshammer.
\newblock From quasi-hereditary algebras with exact {B}orel subalgebras to
  directed bocses.
\newblock {\em Bulletin of the London Mathematical Society}, 52(2):367--378,
  2020.

\bibitem[BLM08]{BLM08}
Yu. Bespalov, V.~Lyubashenko, and O.~Manzyuk.
\newblock {\em Pretriangulated {$A_\infty$}-categories}, volume~76 of {\em
  Proceedings of Institute of Matheamtics of NAS of Ukraine. Mathematics and
  Applications}.
\newblock National{\'n}a Akadami ya Nauk Ukraini, Kiev, 2008.

\bibitem[Bon83]{Bon83}
Klaus Bongartz.
\newblock Algebras and quadratic forms.
\newblock {\em Journal of the London Mathematical Society. Second Series},
  28(3):461--469, 1983.

\bibitem[Bon89]{Bon89}
Alexey Bondal.
\newblock {Representations of associative algebras and coherent sheaves}.
\newblock {\em Izvestiia Akademii Nauk SSSR. Seriya Matematicheskaya},
  53(1):25--44, 1989.

\bibitem[Boo21]{Boo21}
Matt Booth.
\newblock The derived deformation theory of a point.
\newblock arXiv: 2009.09150, to appear in: Mathematische Zeitschrift, 2021.

\bibitem[Brz13]{Brz13}
Tomasz Brzezi{\'n}ski.
\newblock {Curved differential graded algebras and corings}.
\newblock {\em Bulletin of the Belgian Mathematical Society. Simon Stevin},
  20(5):909--936, 2013.

\bibitem[BSZ09]{BSZ09}
Raymundo Bautista, Leonardo Salmer{\'o}n, and Rita Zuazua.
\newblock {\em {Differential tensor categories and their module categories}},
  volume 362 of {\em {London Mathematical Society Lecture Note Series}}.
\newblock Cambridge University Press, Cambridge, 2009.
\newblock x+452 pp.

\bibitem[Bur05]{Bur05}
William~L. Burt.
\newblock {Almost Split Sequences and {BOCS}es}.
\newblock unpublished manuscript, August 2005.

\bibitem[Bur18]{Bur18}
Jesse Burke.
\newblock {Transfer of A-infinity structures to projective resolutions}.
\newblock preprint, arXiv: 1801.08933, 2018.

\bibitem[BW03]{BW03}
Tomasz Brzezinski and Robert Wisbauer.
\newblock {\em {Corings and comodules}}, volume 309 of {\em {London
  Mathematical Society Lecture Note Series}}.
\newblock Cambridge University Press, Cambridge, 2003.
\newblock xii+476 pp.

\bibitem[CK21]{CK21}
Teresa Conde and Julian K{\"u}lshammer.
\newblock Homological embedding, exact borel subalgebras, and
  {$A_\infty$}-algebras.
\newblock in preparation, 2021.

\bibitem[Con11]{Con11}
Andrew~Brondos Conner.
\newblock {\em {{$A_\infty$}-structures, generalized {K}oszul properties and
  combinatorial topology}}.
\newblock PhD thesis, University of Oregon, 2011.

\bibitem[Con21]{Con21}
Teresa Conde.
\newblock All quasihereditary algebras with a regular exact {B}orel subalgebra.
\newblock {\em Advances in Mathematics}, 384:Paper No. 107751, 45, 2021.

\bibitem[CPS88]{CPS88}
Edward Cline, Brian~J. Parshall, and L{\'e}onard~Lewy Scott.
\newblock {Finite-dimensional algebras and highest weight categories}.
\newblock {\em Journal f{\"u}r die Reine und Angewandte Mathematik. [Crelle's
  Journal]}, 391:85--99, 1988.

\bibitem[DR92]{DR92}
Vastimil Dlab and Claus~Michael Ringel.
\newblock {The module theoretical approach to quasi-hereditary algebras}.
\newblock In {\em {Representations of algebras and related topics ({K}yoto,
  1990)}}, volume 168 of {\em {London Mathematical Society Lecture Note
  Series}}, pages 200--224. Cambridge University Press, Cambridge, 1992.

\bibitem[Dro80]{D80}
Yurii~Anatolievich Drozd.
\newblock {Tame and wild matrix problems}.
\newblock In {\em {Representation theory, {II} ({P}roceedings of the {S}econd
  {I}nternational {C}onference, {C}arleton {U}niversity, {O}ttawa, {ON},
  1979)}}, volume 832 of {\em {Lecture Notes in Mathematics}}, pages 242--258.
  Springer, Berlin, 1980.

\bibitem[EL17]{EL17}
Tobias Ekholm and Yanki Lekili.
\newblock {Duality between {L}agrangian and {L}egendrian invariants}.
\newblock Preprint, arXiv:1701.01284, 2017.

\bibitem[Gov73]{Gov73}
V.~E. Govorov.
\newblock The global dimension of algebras.
\newblock {\em Akademiya Nauk SSSR. Matematicheskie Zametki}, 14:399--406,
  1973.

\bibitem[Her18]{Her18}
Estanislao Herscovich.
\newblock {Using torsion theory to compute the algebraic structure of
  {H}ochschild (co)homology}.
\newblock {\em Homology, Homotopy and Applications}, 20(1):117--139, 2018.

\bibitem[Iya03]{Iya03}
Osamu Iyama.
\newblock {Finiteness of representation dimension}.
\newblock {\em Proc. Amer. Math. Soc.}, 131(4):1011--1014, 2003.

\bibitem[Kad80]{Kad80}
T.~V. Kadeishvili.
\newblock {On the theory of homology of fiber spaces}.
\newblock {\em Akademiya Nauk SSSR i Moskovskoe matematicheskoe Obshchestvo.
  Uspekhi Matematicheskikh Nauk}, 35(3(213)):183--188, 1980.
\newblock International Topology Conference (Moscow State Univ., Moscow, 1979).

\bibitem[Kad82]{Kad82}
T.~V. Kadeishvili.
\newblock {The algebraic structure in the homology of an
  {$A(\infty)$}-algebra}.
\newblock {\em Soobshcheniya Akademii Nauk Gruzinskoui SSR}, 108(2):249--252
  (1983), 1982.

\bibitem[Kel01]{K01}
Bernhard Keller.
\newblock {Introduction to {$A$}-infinity algebras and modules}.
\newblock {\em Homology, Homotopy and Applications}, 3(1):1--35, 2001.

\bibitem[Kel02]{K02}
Bernhard Keller.
\newblock {{$A$}-infinity algebras in representation theory}.
\newblock In {\em {Representations of algebra. {V}ol. {I}, {II}}}, pages
  74--86. Beijing Normal University Press, 2002.

\bibitem[Kel06]{K06}
Bernhard Keller.
\newblock {{$A$}-infinity algebras, modules and functor categories}.
\newblock In {\em {Trends in representation theory of algebras and related
  topics}}, volume 406 of {\em {Contemporary Mathematics}}, pages 67--93.
  American Mathematical Society, Providence, RI, 2006.

\bibitem[KK99]{KK99}
Michael Klucznik and Steffen Koenig.
\newblock {Characteristic tilting modules over quasi hereditary algebras}.
\newblock Sonderforschungsbereich Diskrete Strukturen in der Mathematik
  (Bielefeld): Erg{\"a}nzungsreihe ; 99,004, 1999.

\bibitem[KK18]{KK18}
Martin Kalck and Joseph Karmazyn.
\newblock Ringel duality for certain strongly quasi-hereditary algebras.
\newblock {\em European Journal of Mathematics}, 4(3):1100--1140, 2018.

\bibitem[KKO14]{KKO14}
Steffen Koenig, Julian K{\"u}lshammer, and Sergiy Ovsienko.
\newblock {Quasi-hereditary algebras, exact {B}orel subalgebras,
  {$A_\infty$}-categories and boxes}.
\newblock {\em Advances in Mathematics}, 262:546--592, 2014.

\bibitem[KM23]{KM23}
Julian K{\"u}lshammer and Vanessa Miemietz.
\newblock Uniqueness of exact {B}orel subalgebras and bocses - {S}age{M}ath
  code, 2023.
\newblock https://zenodo.org/badge/latestdoi/678389337.

\bibitem[K{\"o}n95]{Koe95}
Steffen K{\"o}nig.
\newblock {Exact {B}orel subalgebras of quasi-hereditary algebras. {I}}.
\newblock {\em Mathematische Zeitschrift}, 220(3):399--426, 1995.
\newblock With an appendix by Leonard Scott.

\bibitem[Kop17]{Kop17}
Jakub Kop{\v{r}}iva.
\newblock On the homotopy transfer of {$A_\infty$} structures.
\newblock {\em Archivum Mathematicum}, 53(5):267--312, 2017.

\bibitem[KR77]{KR77}
Mark~M. Kl{\u\i}ner and A.~V. Ro{\u\i}ter.
\newblock {Representations of differential graded categories}.
\newblock In {\em {Matrix problems ({R}ussian)}}, pages 5--70. Akad. Nauk
  Ukrain. SSR, Inst. Mat., Kiev, 1977.

\bibitem[KS01]{KS01}
Maxim Kontsevich and Yan~S. So{\u\i}bel'man.
\newblock {Homological mirror symmetry and torus fibrations}.
\newblock In {\em {Symplectic geometry and mirror symmetry ({S}eoul, 2000)}},
  pages 203--263. World Sci. Publ., River Edge, NJ, 2001.

\bibitem[K{\"u}l17]{Kul17}
Julian K{\"u}lshammer.
\newblock {In the bocs seat: {Q}uasi-hereditary algebras and representation
  type.}
\newblock In {\em {Representation theory - current trends and perspectives}},
  {EMS Ser. Congr. Rep.}, pages 375--426. Eur. Math. Soc., Z{\"u}rich, 2017.

\bibitem[LDM05]{HL05}
Ji-Wei~He Lu and Di-Ming.
\newblock {Higher {K}oszul algebras and {$A$}-infinity algebras}.
\newblock {\em Journal of Algebra}, 293(2):335--362, 2005.

\bibitem[LH03]{L-H03}
Kenji Lef{\`e}vre-Hasegawa.
\newblock {\em {Sur les {$A_\infty$}-cat{\'e}gories}}.
\newblock PhD thesis, Universit{\'e} Paris 7, November 2003.
\newblock arXiv: math/0310337.

\bibitem[LPWZ04]{LPWZ04}
Di~Ming Lu, John~H. Palmieri, Quan~Shui Wu, and James~J. Zhang.
\newblock {{$A_\infty$}-algebras for ring theorists}.
\newblock In {\em {Proceedings of the {I}nternational {C}onference on
  {A}lgebra}}, volume~11, pages 91--128, 2004.

\bibitem[LPWZ08]{LPWZ08}
Di~Ming Lu, John~H. Palmieri, Quan~Shui Wu, and James~J. Zhang.
\newblock {Koszul equivalences in {$A_\infty$}-algebras}.
\newblock {\em New York Journal of Mathematics}, 14:325--378, 2008.

\bibitem[LPWZ09]{LPWZ09}
Di~Ming Lu, John~H. Palmieri, Quan~Shui Wu, and James~J. Zhang.
\newblock {{$A$}-infinity structure on {E}xt-algebras}.
\newblock {\em Journal of Pure and Applied Algebra}, 213(11):2017--2037, 2009.

\bibitem[LV12]{LV12}
Jean-Louis Loday and Bruno Vallette.
\newblock {\em Algebraic operads}, volume 346 of {\em Grundlehren der
  Mathematischen Wissenschaften [Fundamental Principles of Mathematical
  Sciences]}.
\newblock Springer, Heidelberg, 2012.

\bibitem[Mad02]{Mad02}
Dag Madsen.
\newblock {\em {Homological aspects in representation theory}}.
\newblock PhD thesis, Norwegian University of Science and Technology, February
  2002.

\bibitem[Mar06]{Mar06}
Martin Markl.
\newblock {Transferring {$A_\infty$} (strongly homotopy associative)
  structures}.
\newblock {\em Rendiconti del Circolo Matematico di Palermo. Serie II.
  Supplemento}, (79):139--151, 2006.

\bibitem[Mer99]{Mer99}
S.~A. Merkulov.
\newblock {Strong homotopy algebras of a {K}{\"a}hler manifold}.
\newblock {\em International Mathematics Research Notices}, (3):153--164, 1999.

\bibitem[Mor58]{Mor58}
Kiiti Morita.
\newblock Duality for modules and its applications to the theory of rings with
  minimum condition.
\newblock {\em Science Reports of the Tokyo Kyoiku Daigaku. Section A},
  6:83--142, 1958.

\bibitem[NVW20]{NVW18}
Cris Negron, Yury Volkov, and Sarah Witherspoon.
\newblock {$A_\infty$}-coderivations and the {G}erstenhaber bracket on
  {H}ochschild cohomology.
\newblock {\em Journal of Noncommutative Geometry}, 14(2):531--565, 2020.

\bibitem[Orl18]{Orl18}
Dmitri~O. Orlov.
\newblock Derived non-commutative schemes, geometric realizations, and
  finite-dimensional algebras.
\newblock {\em Uspekhi Matematicheskikh Nauk}, 73(5(443)):123--182, 2018.

\bibitem[Ovs93]{Ovs93}
Sergiy Ovsienko.
\newblock {Generic representations of free bocses}.
\newblock preprint 93-010, Universit{\"a}t Bielefeld, 1993.

\bibitem[Ro{\u\i}79]{Roi79}
A.~V. Ro{\u\i}ter.
\newblock {Matrix problems and representions of {BOCS}es}.
\newblock In {\em {Representations and quadratic forms ({R}ussian)}}, pages
  3--38, 154. Akad. Nauk Ukrain. SSR, Inst. Mat., Kiev, 1979.

\bibitem[Sco87]{Sco87}
Leonard~L. Scott.
\newblock {Simulating algebraic geometry with algebra. {I}. {T}he algebraic
  theory of derived categories}.
\newblock In {\em {The {A}rcate {C}onference on {R}epresentations of {F}inite
  {G}roups ({A}rcata, {C}alif., 1986}}, volume~47 of {\em {Proc. Sympos. Pure
  Math.}}, pages 271--281. American Mathematical Society, Providence, RI, 1987.

\bibitem[Seg08]{Seg08}
Ed~Segal.
\newblock {The {$A_\infty$} deformation theory of a point and the derived
  categories of local {C}alabi-{Y}aus}.
\newblock {\em Journal of Algebra}, 320(8):3232--3268, 2008.

\bibitem[Sei08]{Sei08}
Paul Seidel.
\newblock {\em {Fukaya categories and {P}icard-{L}efschetz Theory}}.
\newblock {Zurich Lectures in Advanced Mathematics}. European Mathematical
  Society (EMS), Z{\"u}rich, 2008.

\bibitem[Sta63]{Sta63}
James~Dillon Stasheff.
\newblock Homotopy associativity of {$H$}-spaces. {I}, {II}.
\newblock {\em Transactions of the American Mathematical Society},
  108:275--292, 293--312, 1963.

\bibitem[Thu20]{Thu20}
Markus Thuresson.
\newblock The {$\operatorname{Ext}$}-algebra of standard modules over dual
  extension algebras.
\newblock arXiv: 2011.11107, to appear in: Journal of Algebra, 2020.

\bibitem[Val14]{Val14}
Bruno Vallette.
\newblock Algebra + homotopy = operad.
\newblock In {\em Symplectic, {P}oisson, and noncommutative geometry},
  volume~62 of {\em Math. Sci. Res. Inst. Publ.}, pages 229--290. Cambridge
  University Press, New York, 2014.

\bibitem[Wal21]{Wal21}
Carl~Felix Waller.
\newblock Formality of {$A_\infty$}-algebras.
\newblock Master's thesis, 2021.
\newblock arXiv: 2107.03195.

\bibitem[Zha00]{Zha00}
Yuehui Zhang.
\newblock {Strong exact {B}orel subalgebras of quasi-hereditary algebras}.
\newblock {\em Journal of Algebra}, 231(2):463--472, 2000.

\end{thebibliography}

\begin{appendices}

In this appendix, we provide proofs for Lemma \ref{mu2i} and Proposition \ref{mu26}. All computations were done by hand by computing the already known terms in the $A_\infty$-equation and then finding an appropriate expression for the missing term, so that applying the terms $1^{\otimes r} \otimes \mu_1 \otimes 1^{\otimes t}$ to it achieves cancellation of all terms. Parts of these computations were checked using SageMath. More precisely, assuming that we have correctly applied the relevant comultiplications $\mu_n$ and identified the terms which vanish by one of Lemmas \ref{yellow}, \ref{green}, or \ref{gg}, we used SageMath to verify that the remaining terms add up to zero. In doing so, we have omitted summations in Appendix \ref{appendixA} and Lemma \ref{mu3i} (c) since these are relatively few and easily checked by hand. The code for this can be found in \cite{KM23}.

\section{The lift of the comultiplication $\mu_2$ on $\mathcal{P}$}
\label{appendixA}
In this appendix, we provide the proof of Lemma \ref{mu2i} to inductively compute the comultiplications $\mu_2^i$. 

\begin{proof}
We make use of the $A_\infty$-equation (specialised to degree $i$):
\[\mu_2^{i-1}\mu_1^i-(\mu_1^i\otimes 1+1\otimes \mu_1^i)\mu_2^i=0\]
Applying this equation to $(1\otimes \shift a_1\otimes \dots \shift a_i\otimes 1)$ we verify that it is satisfied by the stated formula. 
Separately computing the individual summands of $\mu_2^{i-1}\mu_1^i(1\otimes \shift a_1\otimes \dots\otimes  \shift a_i\otimes 1)$ we obtain 
\begin{align*}
\mu_2^{i-1}(a_1\otimes (\shift a_2\otimes \dots\otimes \shift a_i)\otimes 1)&=\underbrace{a_1\hat{\omega}\otimes  \shift a_2\otimes \dots\otimes \shift a_i\otimes 1}_{8}+\underbrace{a_1\otimes \shift a_2\otimes \dots\otimes \shift a_i\otimes \hat{\omega}}_{5}\\
& -\underbrace{\sum_{j=2}^{i-1}a_1\otimes (\shift a_2\otimes \dots\otimes \shift a_j)\otimes 1\otimes (\shift a_{j+1}\otimes \dots\otimes \shift a_i)\otimes 1}_{10}\\
& -\underbrace{a_1 a_{2}^1a_{2}^2\otimes \shift a_{2}^3\otimes \shift a_3\otimes \dots\otimes \shift a_i\otimes 1}_{13}\\
&+\underbrace{a_1\otimes \shift a_2\otimes \dots\otimes \shift a_{i-1}\otimes \shift a_{i}^1\otimes a_{i}^2a_{i}^3}_{15},
\end{align*}
and
\begin{align*}
&\mu_2^{i-1}((-1)^i 1\otimes (\shift a_1\otimes \dots\otimes \shift a_{i-1})\otimes a_i)\\
=\,&\underbrace{(-1)^i \hat{\omega}\otimes \shift a_1\otimes \dots\otimes \shift a_{i-1}\otimes a_i}_{4}+\underbrace{(-1)^i 1\otimes \shift a_1\otimes \dots\otimes \shift a_{i-1}\otimes \hat{\omega} a_i}_{9}\\
&+\underbrace{\sum_{j=1}^{i-2} (-1)^{i+1} 1\otimes \shift a_1\otimes \dots\otimes \shift a_j\otimes 1\otimes \shift a_{j+1}\otimes \dots\otimes \shift a_{i-1}\otimes a_i}_{17}\\
&+\underbrace{ (-1)^{i+1} a_{1}^1a_{1}^2\otimes \shift a_{1}^3\otimes \shift a_2\otimes \dots\otimes \shift a_{i-1}\otimes a_i}_{16}+\underbrace{(-1)^{i} 1\otimes \shift a_1\otimes \dots\otimes \shift a_{i-2}\otimes \shift a_{i-1}^1\otimes a_{i-1}^2 a_{i-1}^3 a_i}_{14},
\end{align*}
and
\begin{align*}
&\mu_2^{i-1}\left(\sum_{j=1}^{i-1}(-1)^j1\otimes \shift a_1\otimes \dots\otimes \shift a_{j-1}\otimes \shift (a_ja_{j+1})\otimes \shift a_{j+2}\otimes \dots\otimes \shift a_i\otimes 1\right)\\
=\,&\underbrace{\sum_{j=1}^{i-1}(-1)^j\hat{\omega} \otimes \shift a_1\otimes \dots\otimes \shift a_{j-1}\otimes \shift (a_ja_{j+1})\otimes \shift a_{j+2}\otimes \dots\otimes \shift a_i\otimes 1}_{6}\\
&+\underbrace{\sum_{j=1}^{i-1}(-1)^j 1\otimes \shift a_1\otimes \dots\otimes \shift a_{j-1}\otimes \shift (a_ja_{j+1})\otimes \shift a_{j+2}\otimes \dots\otimes \shift a_i\otimes \hat{\omega}}_{7}\\
&+\underbrace{\sum_{j=2}^{i-1}\sum_{k=1}^{j-1} (-1)^{j+1}1\otimes \shift a_1\otimes \dots\otimes \shift a_k\otimes 1\otimes \shift a_{k+1} \otimes \dots\otimes \shift a_{j-1}\otimes \shift (a_{j}a_{j+1})\otimes \shift a_{j+2}\otimes  \dots\otimes \shift a_i\otimes 1}_{23}\\
&+\underbrace{\sum_{j=1}^{i-1}\sum_{k=j+1}^{i-1}  (-1)^{j+1} 1\otimes \shift a_1\otimes \dots\otimes \shift a_{j-1}\otimes \shift (a_ja_{j+1})\otimes \shift a_{j+2}\otimes  \dots\otimes \shift a_k\otimes 1\otimes \shift a_{k+1}\otimes \dots\otimes \shift a_i\otimes 1}_{22}\\
&+\underbrace{\sum_{j=2}^{i-1}(-1)^{j+1} a_{1}^1 a^2_{1}\otimes \shift a_{1}^3\otimes \shift a_2\otimes \dots\otimes \shift a_{j-1}\otimes \shift (a_ja_{j+1})\otimes \shift a_{j+2}\otimes \dots\otimes \shift a_i\otimes 1}_{21}\\
&\underbrace{+a_{1}^1a_{1}^2\otimes \shift (a_{1}^3a_2)\otimes \shift a_3\otimes \dots\otimes \shift a_i\otimes 1}_{18}+\underbrace{a_1a_{2}^1a_{2}^2\otimes \shift a_{2}^3\otimes \shift a_3\otimes \dots\otimes\shift  a_i\otimes 1}_{13}\\
&+\underbrace{\sum_{j=1}^{i-2}(-1)^j 1\otimes \shift a_1\otimes \dots\otimes \shift a_{j-1}\otimes \shift (a_ja_{j+1})\otimes \shift a_{j+2}\otimes \dots\otimes \shift a_{i-1}\otimes \shift a_{i}^1\otimes a_{i}^2a_{i}^3}_{20}\\
&+\underbrace{(-1)^{i+1} 1\otimes \shift a_1\otimes \dots\otimes \shift a_{i-2}\otimes \shift a_{i-1}^1\otimes a_{i-1}^2a_{i-1}^3a_i}_{14}+\underbrace{(-1)^{i+1}1\otimes \shift a_1\otimes \dots\otimes \shift a_{i-2}\otimes \shift (a_{i-1} a_{i}^1)\otimes a_{i}^2 a_{i}^3.}_{19}
\end{align*}
On the other hand,
\begin{align*}
-(1\otimes \mu_1^i)(\hat{\omega}\otimes \shift a_1\otimes \dots\otimes \shift a_i\otimes 1)
=\,&-\underbrace{\hat{\omega} a_1\otimes \shift a_2\otimes \dots\otimes \shift a_i\otimes 1}_{1}\\
&-\underbrace{\sum_j (-1)^j\hat{\omega}\otimes \shift a_1\otimes \dots\otimes\shift  a_{j-1}\otimes \shift (a_ja_{j+1})\otimes \shift a_{j+2}\otimes \dots\otimes \shift a_i\otimes 1}_{6}\\
&-\underbrace{(-1)^i \hat{\omega}\otimes \shift a_1\otimes \dots\otimes \shift a_{i-1}\otimes a_i}_{4}
\end{align*}
and
\begin{align*}
-(\mu_1^i\otimes 1)(1\otimes \shift a_1\otimes \dots\otimes \shift a_i\otimes \hat{\omega})
=\,&-\underbrace{a_1\otimes \shift a_2\otimes \dots\otimes \shift a_i\otimes \hat{\omega}}_{5}\\
&-\underbrace{\sum (-1)^j 1\otimes \shift a_1\otimes \dots \otimes \shift a_{j-1}\otimes \shift (a_ja_{j+1})\otimes \shift a_{j+2}\otimes \dots\otimes \shift a_i\otimes \hat{\omega}}_{7}\\
&-\underbrace{(-1)^i1\otimes \shift a_1\otimes \dots\otimes \shift a_{i-1}\otimes a_i\hat{\omega},}_{2}
\end{align*}
and
\begin{align*}
-(\mu_1^1\otimes 1)( -1\otimes \shift a_1 \otimes 1\otimes \shift a_2 \otimes \dots \otimes \shift a_i\otimes 1)
=\,& +\underbrace{\hat{\omega}a_1 \otimes \shift a_2 \otimes \dots \otimes \shift a_i\otimes 1}_{1} - \underbrace{a_1\hat{\omega} \otimes \shift a_2 \otimes \dots \otimes \shift a_i\otimes 1}_{8} \\
& - \underbrace{\shift \partial(a_1) \otimes \shift a_2 \otimes \dots \otimes \shift a_i\otimes 1}_{11} 
\end{align*}
and
\begin{align*}
&\sum_{j=2}^{i-1}(-\mu_1^j\otimes 1) ( -1\otimes \shift a_1\otimes \dots\otimes \shift a_j\otimes 1\otimes \shift a_{j+1}\otimes \dots\otimes \shift a_i\otimes 1)\\
&=\underbrace{\sum_{j=2}^{i-1}  a_1\otimes \shift a_2\otimes \dots\otimes \shift a_j\otimes 1\otimes \shift a_{j+1}\otimes \dots\otimes \shift a_i\otimes 1}_{10}\\
&+\underbrace{\sum_{j=2}^{i-1}\sum_{k=1}^{j-1}  (-1)^{k} 1\otimes \shift a_1\otimes \dots\otimes \shift a_{k-1}\otimes \shift (a_ka_{k+1})\otimes \shift a_{k+2}\otimes \dots\otimes \shift a_j\otimes 1\otimes \shift a_{j+1}\otimes \dots\otimes \shift a_i\otimes 1}_{22}\\
&+\underbrace{\sum_{j=2}^{i-1}(-1)^{j} 1\otimes (\shift a_1\otimes \dots\otimes \shift a_{j-1})\otimes a_j\otimes (\shift a_{j+1}\otimes \dots\otimes \shift a_i)\otimes 1,}_{3}
\end{align*}
and
\begin{align*}
&\sum_{j=2}^{i-1}(-1\otimes \mu_1^{i-j})(-1\otimes \shift a_1\otimes \dots\otimes\shift  a_j\otimes 1\otimes \shift a_{j+1}\otimes \dots\otimes \shift a_i\otimes 1)\\
=\,&\underbrace{\sum_{j=1}^{i-2}  (-1)^{j} 1\otimes (\shift a_1\otimes \dots\otimes \shift a_j)\otimes a_{j+1}\otimes (\shift a_{j+2}\otimes \dots\otimes \shift a_i)\otimes 1}_{3}\\
&+\underbrace{\sum_{j=1}^{i-2}\sum_{k=j+1}^{i-1} (-1)^{j}(-1)^{k-j} 1\otimes \shift a_1\otimes \dots\otimes \shift a_{j}\otimes 1\otimes \shift a_{j+1}\otimes \dots\otimes \shift (a_ka_{k+1})\otimes \dots\otimes \shift a_i\otimes 1}_{23}\\
&+\underbrace{\sum_{j=1}^{i-2}  (-1)^{j}(-1)^{i-j} 1\otimes (\shift a_1\otimes \dots\otimes \shift a_j)\otimes 1\otimes (\shift a_{j+1}\otimes \dots\otimes \shift a_{i-1})\otimes a_i}_{17}\\
&+\underbrace{ (-1)^{i-1} 1\otimes \shift a_1\otimes \dots\otimes \shift a_{i-1}\otimes \hat{\omega} a_i}_{9}+\underbrace{(-1)^{i}1\otimes \shift a_1\otimes \dots\otimes \shift a_{i-1}\otimes a_i\hat{\omega}}_{2}\\
&+\underbrace{ (-1)^{i} 1\otimes \shift a_1\otimes \dots\otimes \shift a_{i-1}\otimes \shift \partial(a_i)}_{12}
\end{align*}
and
\begin{align*}
&-(1\otimes \mu_1^{i})(-a_1^1a_1^2\otimes \shift a_1^3\otimes \shift a_2\otimes \dots \otimes \shift a_i\otimes 1)\\
=\,&\underbrace{\shift \partial(a_1)\otimes \shift a_2\otimes \dots\otimes\shift  a_i\otimes 1}_{11}-\underbrace{a_1^1a_1^2\otimes \shift (a_1^3a_2)\otimes \shift a_3\otimes \dots \otimes \shift a_i\otimes 1}_{18}\\
&+\underbrace{\sum_{j=2}^{i-1} (-1)^ja_1^1a_1^2\otimes \shift a_1^3\otimes \shift a_2\otimes \dots \otimes \shift (a_ja_{j+1})\otimes \dots \otimes \shift a_i\otimes 1}_{21}\\
&+\underbrace{(-1)^i a_1^1a_1^2\otimes \shift a_1^3\otimes\shift  a_2\otimes \dots \otimes \shift a_{i-1}\otimes a_i}_{16}
\end{align*}
and
\begin{align*}
&-(\mu_1^i\otimes 1)(1\otimes\shift  a_1\otimes \dots \otimes \shift a_{i-1}\otimes \shift a_i^1\otimes a_i^2a_i^3)\\
=\,&\underbrace{-a_1\otimes \shift a_2\otimes \dots \otimes \shift a_{i-1}\otimes \shift a_i^1\otimes a_i^2a_i^3}_{15}\\
&+\underbrace{\sum_{j=1}^{i-2}(-1)^{j+1}1\otimes \shift a_1\otimes \dots \otimes \shift (a_ja_{j+1})\otimes \dots \otimes \shift a_{i-1}\otimes \shift a_i^1\otimes a_i^2a_i^3}_{20}\\
&+\underbrace{(-1)^{i}1\otimes \shift a_1\otimes \dots \otimes \shift a_{i-2}\otimes \shift (a_{i-1}a_i^1)\otimes a_i^2a_i^3}_{19}+\underbrace{(-1)^{i+1}1\otimes \shift a_1\otimes \dots \otimes \shift a_{i-1}\otimes \shift \partial(a_i)}_{12}.
\end{align*}
Again, the summands with the same number add to zero, proving the claim.
\end{proof}

\section{The homotopy $\mu_3$ on $\mathcal{P}$ up to which $\mu_2$ is coassociative}
\label{appendixB}

The lifting $\mu_2$ of $\mu$ to $\mathcal{P}$ is coassociative up to homotopy. In this section, we provide explicit formulae for this homotopy $\mu_3$.  In all of our computations, whenever we compute expressions of the form $(f_1+\dots+f_r)(a_1+\dots+a_s)$ we order the result and mark the beginning of the application of an $f_i$ to a new $a_j$ by a boldface plus (or minus) symbol and the application of a new $f_i$ as follows:
\[f_1(a_1)\+f_1(a_2)\+\dots\+f_1(a_s)\uwave{\+}f_2(a_1)\+f_2(a_2)\+\dots\+f_2(a_s)\uwave{\+}\dots\uwave{\+}f_r(a_1)\+f_r(a_2)\+\dots\+f_r(a_s)\]

We need the following lemma to simplify calculations by identifying vanishing summands.

\begin{lem}\label{yellow}
Let $a\in A$, then, using Sweedler notation,
\begin{align*}
a^{11}\otimes \shift^{-1}(a^{12})\otimes a^{13}\otimes \shift^{-1}(a^{2})\otimes a^3&=0 ,\\
a^1a^{21}\otimes \shift^{-1}(a^{22})\otimes a^{23}\otimes \shift^{-1}(a^{24})\otimes a^{25}a^3&=0,\\
a^1\otimes \shift^{-1}(a^2)\otimes a^{31}\otimes \shift^{-1}(a^{32})\otimes a^{33}&=0.
\end{align*}

For $\varphi\in \Phi$, we have

\begin{align*}
\varphi^{11}\otimes(\shift^{-1}\varphi^{12})\otimes\varphi^{13}\otimes(\shift^{-1}\varphi^2)\otimes\varphi^3\otimes(\shift^{-1}\varphi^4)\otimes\varphi^5 &=0,\\
\varphi^1\varphi^{21}\otimes(\shift^{-1}\varphi^{22})\otimes\varphi^{23}\otimes(\shift^{-1}\varphi^{24})\otimes\varphi^{25}\varphi^3\otimes(\shift^{-1}\varphi^4)\otimes\varphi^5&=0,\\
\varphi^1\otimes(\shift^{-1}\varphi^2)\otimes\varphi^{31}\otimes(\shift^{-1}\varphi^{32})\otimes\varphi^{33}\otimes(\shift^{-1}\varphi^4)\otimes\varphi^5 &=0,\\
\varphi^1\otimes(\shift^{-1}\varphi^2)\otimes\varphi^3\varphi^{41}\otimes(\shift^{-1}\varphi^{42})\otimes\varphi^{43}\otimes(\shift^{-1}\varphi^{44})\otimes\varphi^{45}\varphi^5 &=0,\\
\varphi^1\otimes(\shift^{-1}\varphi^2)\otimes\varphi^3\otimes(\shift^{-1}\varphi^4)\otimes\varphi^{51}\otimes(\shift^{-1}\varphi^{52})\otimes\varphi^{53} &=0.
\end{align*}
These vanishings also hold true for an arbitrary number of tensor symbols removed.
\end{lem}

\begin{proof}

As $\partial$ is a differential, we obtain
\[0=\partial^2(a)=a^{11} \shift^{-1} (a^{12}) a^{13}\shift^{-1}  (a^{2}) a^3+
a^1a^{21}\shift^{-1}(a^{22}) a^{23} \shift^{-1}(a^{24}) a^{25}a^3-a^1\shift^{-1}(a^2)a^{31}\shift^{-1}(a^{32})a^{33}
\]

and 
\begin{equation*}
\begin{split}
0= &-\partial^2(\varphi)\\
= &\partial(\varphi^1(\shift^{-1}\varphi^2)\varphi^3(\shift^{-1}\varphi^4)\varphi^5)\\
= &\varphi^{11}(\shift^{-1}\varphi^{12})\varphi^{13}(\shift^{-1}\varphi^2)\varphi^3(\shift^{-1}\varphi^4)\varphi^5 
+ \varphi^1\varphi^{21}(\shift^{-1}\varphi^{22})\varphi^{23}(\shift^{-1}\varphi^{24})\varphi^{25}\varphi^3(\shift^{-1}\varphi^4)\varphi^5\\
&- \varphi^1(\shift^{-1}\varphi^2)\varphi^{31}(\shift^{-1}\varphi^{32})\varphi^{33}(\shift^{-1}\varphi^4)\varphi^5 
- \varphi^1(\shift^{-1}\varphi^2)\varphi^3\varphi^{41}(\shift^{-1}\varphi^{42})\varphi^{43}(\shift^{-1}\varphi^{44})\varphi^{45}\varphi^5\\
&+\varphi^1(\shift^{-1}\varphi^2)\varphi^3(\shift^{-1}\varphi^4)\varphi^{51}(\shift^{-1}\varphi^{52})\varphi^{53}
\end{split}
\end{equation*}
As the biquiver of $\overline{V}$ is directed, it again follows that each of the summands is zero. By projectivity of $\overline{V}$, we have $A\otimes \Phi\otimes A\cong \overline{V}$ and 
$\overline{V}\otimes_A\overline{V}\otimes_A\overline{V}\cong A\otimes \Phi\otimes A\otimes \Phi\otimes A\otimes \Phi\otimes A$.  Hence introducing tensors does not change the vanishing of the terms.
\end{proof}

Whenever we identify a vanishing summand using the preceding lemma in the sequel, we mark it with a $y$. 

\begin{lem}\label{mu3i}
Let $\varphi\in \Phi$, $i\in \mathbb{N}$ and $a,a_1,\dots,a_i\in J$. Then, using Sweedler notation (writing $\partial(a)=a^1\shift^{-1}(a^2)a^3$ and furthermore, $\partial(\shift^{-1}a^2)=a^{21}\shift^{-1}(a^{22}) a^{23} \shift^{-1}(a^{24}) a^{25}$, we obtain, 
\begin{enumerate}[(a)]
\item\label{mu3:0} \begin{align*}\mu_3^0(\varphi)=\,&1\otimes \shift\varphi^1\otimes\varphi^2\varphi^3\varphi^4\varphi^5-\varphi^1\varphi^2\otimes \shift\varphi^3\otimes \varphi^4\varphi^5\\&+\varphi^1\varphi^2\varphi^3\varphi^4\otimes \shift\varphi^5\otimes 1
\end{align*}
\item\label{mu3:1} 
\begin{align*}
\mu_3^1(1\otimes \shift a\otimes 1)=\,&1\otimes \shift a^1\otimes a^2\otimes \shift a^3\otimes 1+a^1a^{21}\dots a^{24}\otimes \shift a^{25}\otimes \shift a^3\otimes 1\\+&1\otimes \shift a^1\otimes \shift a^{21}\otimes a^{22}\dots a^{25} a^3
\end{align*}
\item\label{mu3:i} 
\begin{align*}
\mu_3^i(1\otimes\shift  a_1\otimes \dots\otimes \shift a_i\otimes 1)=\,&\sum_{j=1}^{i} (-1)^{j+1} 1\otimes \shift a_1\otimes \dots\otimes \shift a_{j-1}\otimes \shift a_j^1\otimes a_j^2\otimes \shift a_j^3\otimes \shift a_{j+1}\otimes\dots\otimes \shift a_i\otimes 1\\
&+a_1^1a_1^{21}\dots a_1^{24}\otimes \shift a_1^{25}\otimes \shift a_1^3\otimes \shift a_2\otimes \dots\otimes\shift  a_i\otimes 1\\
&+(-1)^{i+1} 1\otimes \shift a_1\otimes \dots\otimes \shift a_{i-1}\otimes \shift a_i^1\otimes \shift a_i^{21}\otimes a_i^{22}\dots a_i^{25}a_i^3.
\end{align*}
\end{enumerate}
Again, we use the convention of only summing over those $a_1^3, a_i^1$ which are in $J$. 
\end{lem}

\begin{proof}
For the $A_\infty$-relations to hold, we need
\[(\mu_1\otimes 1^{\otimes 2}+1\otimes \mu_1\otimes 1+1^{\otimes 2}\otimes \mu_1)\mu_3+\mu_3\mu_1+(1\otimes \mu_2-\mu_2\otimes 1)\mu_2=0.\]
Using the preceding lemma, we construct the homotopy inductively by verifying that this is satisfied by the stated formulae.

\eqref{mu3:0}: We first compute $(1\otimes \mu_2-\mu_2\otimes 1)\mu_2$ in degree $0$. 
\begin{align*}
(1\otimes \mu_2-\mu_2\otimes 1)\mu_2(\varphi)=\,&(1\otimes \mu_2-\mu_2\otimes 1)\left(\hat{\omega}\varphi+\varphi\hat{\omega}+(\shift\otimes\shift)\partial(\shift^{-1}\varphi)\right)\\
=\,&(1\otimes \mu_2-\mu_2\otimes 1)\left(\hat{\omega}\varphi+\varphi\hat{\omega}-\varphi^1\varphi^2\varphi^3\varphi^4\varphi^5\right)\\
=\,&\underbrace{\hat{\omega}\hat{\omega}\varphi}_{1}+\underbrace{\hat{\omega}\varphi\hat{\omega}}_2+\underbrace{\hat{\omega}(\shift\otimes\shift)\partial(\shift^{-1}\varphi)}_3\+\underbrace{\varphi\hat{\omega}\hat{\omega}}_4\\
&\-\underbrace{\varphi^1\varphi^2\varphi^3\hat{\omega}\varphi^4\varphi^5}_6-\underbrace{\varphi^1\varphi^2\varphi^3\varphi^4\hat{\omega}\varphi^5}_7-\underbrace{
\varphi^1\varphi^2\varphi^3(\shift\otimes\shift)\partial(\shift^{-1}\varphi^4)\varphi^5}_y\\
&\boldsymbol{\uwave{-}}\underbrace{\hat{\omega}\hat{\omega}\varphi}_1\-\underbrace{\hat{\omega}\varphi\hat{\omega}}_2-\underbrace{\varphi\hat{\omega}\hat{\omega}}_4-\underbrace{(\shift\otimes\shift)\partial(\shift^{-1}\varphi)\hat{\omega}}_5\\
&\+\underbrace{\varphi^1\hat{\omega}\varphi^2\varphi^3\varphi^4\varphi^5}_8+\underbrace{\varphi^1\varphi^2\hat{\omega}\varphi^3\varphi^4\varphi^5}_9+\underbrace{\varphi^1(\shift\otimes\shift)\partial(\shift^{-1}\varphi^2)\varphi^3\varphi^4\varphi^5}_y
\end{align*}

Furthermore,
\begin{align*}
&(\mu_1\otimes 1^{\otimes 2}+1\otimes \mu_1\otimes 1+1^{\otimes 2}\otimes \mu_1)\mu_3(\varphi)\\
=\,&(\mu_1\otimes 1^{\otimes 2}+1\otimes \mu_1\otimes 1+1^{\otimes 2}\otimes \mu_1)(1\otimes \shift \varphi^1\otimes \varphi^2\varphi^3\varphi^4\varphi^5-\varphi^1\varphi^2\otimes \shift\varphi^3\otimes \varphi^4\varphi^5+\varphi^1\dots\varphi^4\otimes \shift\varphi^5\otimes 1)\\
=\,&-\underbrace{\hat{\omega}(\shift\otimes\shift)\partial(\shift^{-1}\varphi)}_3-\underbrace{\varphi^1\hat{\omega}\varphi^2\varphi^3\varphi^4 \varphi^5}_8-\underbrace{\shift\partial(\varphi^1)\varphi^2\varphi^3\varphi^4\varphi^5}_y\+0\+0\boldsymbol{\uwave{-}}0\-\underbrace{\varphi^1\varphi^2\hat{\omega}\varphi^3\varphi^4\varphi^5}_9+\underbrace{\varphi^1\varphi^2\varphi^3\hat{\omega}\varphi^4\varphi^5}_6\\
&+\underbrace{\varphi^1\varphi^2\shift\partial(\varphi^3)\varphi^4\varphi^5}_y\+0\uwave{\+}0\+0+\underbrace{\varphi^1\varphi^2\varphi^3 \varphi^4\hat{\omega}\varphi^5}_7+\underbrace{(\shift\otimes\shift)\partial(\shift^{-1}\varphi)\hat{\omega}}_5-\underbrace{\varphi^1\varphi^2\varphi^3\varphi^4\shift\partial(\varphi^5)}_y
\end{align*}

The terms marked with the same number cancel while the terms marked with a $y$ vanish because of the preceding lemma. As $\mu_1^0=0$, this $A_\infty$-relation does hold.

\eqref{mu3:1}: We check that $\mu_3^1$ satisfies the required property:
 
\begin{align*}
&(1\otimes \mu_2-\mu_2\otimes 1)\mu_2(1\otimes \shift a\otimes 1)\\
=\,&(1\otimes \mu_2-\mu_2\otimes 1)(\hat{\omega}\otimes \shift a\otimes 1+1\otimes \shift a\otimes \hat{\omega}-a^1a^2\otimes \shift a^3\otimes 1+1\otimes\shift  a^1\otimes a^2a^3)\\
=\,&\underbrace{\hat{\omega}\hat{\omega}\otimes \shift a\otimes 1}_1+\underbrace{\hat{\omega}\otimes \shift a\otimes \hat{\omega}}_2-\underbrace{\hat{\omega}a^1a^2\otimes \shift a^3\otimes 1}_5+\underbrace{\hat{\omega}\otimes \shift a^1\otimes a^2a^3}_8\+\underbrace{1\otimes \shift a\otimes \hat{\omega}\hat{\omega}}_3\\
&\-\underbrace{a^1a^2\hat{\omega}\otimes \shift a^3\otimes 1}_{12}-\underbrace{a^1a^2\otimes \shift a^3\otimes \hat{\omega}}_9+\underbrace{a^1a^2a^{31}a^{32}\otimes \shift a^{33}\otimes 1}_y-\underbrace{a^1a^2\otimes \shift a^{31}\otimes a^{32}a^{33}}_y\\
&\+\underbrace{1\otimes \shift a^1\otimes \hat{\omega}a^2a^3}_{11}+\underbrace{1\otimes \shift a^1\otimes a^2\hat{\omega}a^3}_7+\underbrace{1\otimes \shift a^1\otimes (\shift\otimes\shift)\partial(\shift^{-1} a^2)a^3}_{13}\\
&\boldsymbol{\uwave{-}}\underbrace{\hat{\omega}\hat{\omega}\otimes\shift  a\otimes 1}_1\-\underbrace{\hat{\omega}\otimes \shift a\otimes \hat{\omega}}_2-\underbrace{1\otimes \shift a\otimes \hat{\omega}\hat{\omega}}_3+\underbrace{a^1a^2\otimes\shift  a^3\otimes \hat{\omega}}_9-\underbrace{1\otimes \shift a^1\otimes a^2a^3\hat{\omega}}_4\\
&\+\underbrace{a^1\hat{\omega}a^2\otimes \shift a^3\otimes 1}_6+\underbrace{a^1a^2\hat{\omega}\otimes \shift a^3\otimes 1}_{12}+\underbrace{a^1 (\shift\otimes\shift)\partial(\shift^{-1} a^2)\otimes \shift a^3\otimes 1}_{14}\\
&\-\underbrace{\hat{\omega}\otimes \shift a^1\otimes a^2a^3}_8-\underbrace{1\otimes \shift a^1\otimes \hat{\omega}a^2a^3}_{11}+\underbrace{a^{11}a^{12}\otimes\shift  a^{13}\otimes a^2a^3}_y-\underbrace{1\otimes \shift a^{11}\otimes a^{12}a^{13}a^2a^3}_y
\end{align*} 
and
\begin{align*}
\mu_3\mu_1(1\otimes \shift a\otimes 1)=\,&\mu_3(\hat{\omega}a-a\hat{\omega}-a^1a^2a^3)\\
=\,0\+0&\-\underbrace{a^1\otimes \shift a^{21}\otimes a^{22}\dots a^{25}a^3}_{15}+\underbrace{a^1a^{21}a^{22}\otimes \shift a^{23}\otimes a^{24}a^{25}a^3}_y-\underbrace{a^1a^{21}\dots a^{24}\otimes \shift a^{25}\otimes a^3}_{10}
\end{align*}
and, for our asserted formula for $\mu_3^1$,
\begin{align*}
&(\mu_1\otimes 1^{\otimes 2}+1\otimes \mu_1\otimes 1+1^{\otimes 2}\otimes \mu_1)\mu_3(1\otimes \shift a\otimes 1)\\
=\,&(\mu_1\otimes 1^{\otimes 2}+1\otimes \mu_1\otimes 1+1^{\otimes 2}\otimes \mu_1)
(1\otimes \shift a^1\otimes a^2\otimes \shift a^3\otimes 1+a^1a^{21}\dots a^{24}\otimes \shift a^{25}\otimes \shift a^3\otimes 1\\&+1\otimes \shift a^1\otimes \shift a^{21}\otimes a^{22}\dots a^{25} a^3)\\
=\,&+\underbrace{\hat{\omega} a^1a^2\otimes \shift a^3\otimes 1}_5-\underbrace{a^1\hat{\omega}a^2\otimes\shift  a^3\otimes 1}_6-\underbrace{\shift \partial(a^1)a^2\otimes \shift a^3\otimes 1}_y\+0\\
&\+\underbrace{a^1\otimes \shift a^{21}\otimes a^{22}\dots a^{25}a^3}_{15}-\underbrace{1\otimes \shift (a^1a^{21})\otimes a^{22}\dots a^{25}a^3}_y-\underbrace{1\otimes \shift a^1\otimes (\shift \otimes \shift )\partial(\shift^{-1} a^2)a^3}_{13}\uwave{\+}0\+0\+0\\
&\boldsymbol{\uwave{-}}\underbrace{1\otimes \shift a^1\otimes a^2\hat{\omega}a^3}_7+\underbrace{1\otimes \shift a^1\otimes a^2a^3\hat{\omega}}_4+\underbrace{1\otimes \shift a^1\otimes a^2\shift \partial(a^3)}_y\\
&\-\underbrace{a^1(\shift \otimes \shift )\partial(\shift^{-1} a^2)\otimes \shift a^3\otimes 1}_{14}-\underbrace{a^1a^{21}\dots a^{24}\otimes \shift (a^{25}a^3)\otimes 1}_y+\underbrace{a^1a^{21}\dots a^{24}\otimes \shift a^{25}\otimes a^3}_{10}\+0.
\end{align*}
Again, the terms with the same number add up to zero while the terms marked with $y$ vanish by the Lemma \ref{yellow}.

\eqref{mu3:i} We proceed by induction on $i$. 
\begin{align*}
&(-\mu_2\otimes_A 1+1\otimes_A \mu_2)\mu_2(1\otimes \shift a_1\otimes \dots\otimes \shift a_i\otimes 1)\\
=\,&(-\mu_2\otimes_A 1+1\otimes_A \mu_2)(\hat{\omega}\otimes \shift a_1\otimes \dots\otimes \shift a_i\otimes 1+1\otimes \shift a_1\otimes \dots\otimes \shift a_i\otimes \hat{\omega}\\
&-\sum_{j=1}^{i-1} 1\otimes \shift a_1\otimes \dots\otimes \shift a_j\otimes 1\otimes \shift a_{j+1}\otimes \dots\otimes \shift a_i\otimes 1\\
&-a_1^1a_1^2\otimes \shift a_1^3\otimes \shift a_2\otimes \dots\otimes \shift a_i\otimes 1 +1\otimes \shift a_1\otimes \dots\otimes \shift a_{i-1}\otimes \shift a_i^1\otimes a_i^2a_i^3)
\end{align*}
We compute the terms separately:
\begin{align*}
&(-\mu_2\otimes_A 1+1\otimes_A \mu_2)(\hat{\omega}\otimes \shift a_1\otimes \dots\otimes \shift a_i\otimes 1+1\otimes \shift a_1\otimes \dots\otimes \shift a_i\otimes \hat{\omega})\\
=\,&\underbrace{-\hat{\omega}\hat{\omega}\otimes \shift a_1\otimes \dots\otimes \shift a_i\otimes 1}_{13}\underbrace{\-\hat{\omega}\otimes \shift a_1\otimes \dots\otimes \shift a_i\otimes \hat{\omega}}_{14}-\underbrace{1\otimes \shift a_1\otimes \dots\otimes \shift a_i\otimes \hat{\omega}\hat{\omega}}_{15}\\
&+\underbrace{\sum_{j=1}^{i-1} 1\otimes \shift a_1\otimes \dots\otimes \shift a_j\otimes 1\otimes \shift a_{j+1}\otimes \dots\otimes \shift a_i\otimes \hat{\omega}}_{17}\\
&+\underbrace{a_1^1a_1^2\otimes \shift a_1^3\otimes \shift a_2\otimes \dots\otimes \shift a_i\otimes \hat{\omega}}_{25}-\underbrace{1\otimes \shift a_1\otimes \dots\otimes \shift a_{i-1}\otimes \shift a_i^1\otimes a_i^2a_i^3\hat{\omega}}_{29}\\
&\uwave{\+}\underbrace{\hat{\omega}\hat{\omega}\otimes \shift a_1\otimes \dots\otimes \shift a_i\otimes 1}_{13}
+\underbrace{\hat{\omega}\otimes \shift a_1\otimes \dots\otimes \shift a_i\otimes \hat{\omega}}_{14}
\\
&\underbrace{-\sum_{j=1}^{i-1}\hat{\omega}\otimes \shift a_1\otimes \dots\otimes \shift a_j\otimes 1\otimes \shift a_{j+1}\otimes \dots\otimes \shift a_i\otimes 1}_{18}\\
&-\underbrace{\hat{\omega} a_1^1a_1^2\otimes \shift a_1^3\otimes \shift a_2\otimes \dots\otimes \shift a_i\otimes 1}_{28}+\underbrace{\hat{\omega}\otimes \shift a_1\otimes \dots\otimes \shift a_{i-1}\otimes \shift a_i^1\otimes a_i^2a_i^3}_{24}\\
&\+\underbrace{1\otimes \shift a_1\otimes \dots\otimes \shift a_i\otimes \hat{\omega}\hat{\omega}}_{15}
\end{align*}
and
\begin{align*}
&(-\mu_2\otimes_A 1+1\otimes_A \mu_2)\left(-\sum_{j=1}^{i-1} 1\otimes \shift a_1\otimes \dots\otimes \shift a_j\otimes 1\otimes \shift a_{j+1}\otimes \dots\otimes \shift a_i\otimes 1\right)\\
=\,&\underbrace{\sum_{j=1}^{i-1}\hat{\omega}\otimes \shift a_1\otimes \dots\otimes \shift a_j\otimes 1\otimes \shift a_{j+1}\otimes \dots\otimes \shift a_i\otimes 1}_{18}+\underbrace{\sum_{j=1}^{i-1}1\otimes \shift a_1\otimes \dots\otimes \shift a_j\otimes \hat{\omega}\otimes \shift a_{j+1}\otimes \dots\otimes \shift a_i\otimes 1}_{16}\\
&-\underbrace{\sum_{j=1}^{i-1}\sum_{k=1}^{j-1} 1\otimes \shift a_1\otimes \dots\otimes \shift a_k\otimes 1\otimes \shift a_{k+1}\otimes\dots\otimes \shift a_j\otimes 1\otimes \shift a_{j+1}\otimes \dots\otimes \shift a_i\otimes 1}_{19}\\
&-\underbrace{a_1^1a_1^2\otimes \shift a_1^3\otimes 1\otimes \shift a_2\otimes \dots\otimes \shift a_i\otimes 1}_{34}-\underbrace{\sum_{j=2}^{i-1}a_1^1a_1^2\otimes \shift a_1^3\otimes \shift a_2\otimes \dots\otimes \shift a_j\otimes 1\otimes \shift a_{j+1}\otimes \dots\otimes \shift a_i\otimes 1}_{26}\\
&+\underbrace{\sum_{j=1}^{i-1} 1\otimes \shift a_1\otimes \dots\otimes \shift a_{j-1}\otimes \shift a_j^1\otimes a_j^2a_j^3\otimes \shift a_{j+1}\otimes \dots\otimes \shift a_i\otimes 1}_{35}\\
&\boldsymbol{\uwave{-}}\underbrace{\sum_{j=1}^{i-1}1\otimes \shift a_1\otimes \dots\otimes \shift a_{j}\otimes \hat{\omega}\otimes \shift a_{j+1}\otimes \dots\otimes \shift a_i\otimes 1}_{16}-\underbrace{\sum_{j=1}^{i-1}1\otimes \shift a_1\otimes \dots\otimes \shift a_j\otimes 1\otimes \shift a_{j+1}\otimes \dots\otimes \shift a_i\otimes \hat{\omega}}_{17}\\
&+ \underbrace{\sum_{j=1}^{i-1}\sum_{k=j+1}^{i-1} 1\otimes \shift a_1\otimes \dots\otimes \shift a_j\otimes 1\otimes \shift a_{j+1}\otimes \dots\otimes \shift a_k\otimes 1\otimes \shift a_{k+1}\otimes \dots\otimes \shift a_i\otimes 1}_{19}\\
&+\underbrace{\sum_{j=1}^{i-1}1\otimes \shift a_1\otimes \dots\otimes \shift a_{j}\otimes a_{j+1}^1a_{j+1}^2\otimes \shift a_{j+1}^3\otimes \shift a_{j+2}\otimes \dots\otimes \shift a_i\otimes 1}_{36}\\
&-\underbrace{\left(\sum_{j=1}^{i-2} 1\otimes \shift a_1\otimes \dots\otimes \shift a_{j}\otimes 1\otimes \shift a_{j+1}\otimes \dots\otimes \shift a_{i-1}\otimes \shift a_i^1\otimes a_i^2a_i^3\right)}_{27}-\underbrace{1\otimes \shift a_1\otimes \dots\otimes \shift a_{i-1}\otimes 1\otimes \shift a_i^1\otimes a_i^2a_i^3}_{33}
\end{align*}
and
\begin{align*}
&(-\mu_2\otimes_A 1+1\otimes_A \mu_2)(-a_1^1a_1^2\otimes \shift a_1^3\otimes \shift a_2\otimes \dots\otimes \shift a_i\otimes 1+1\otimes \shift a_1\otimes \dots\otimes \shift a_{i-1}\otimes \shift a_i^1\otimes a_i^2a_i^3)\\
=\,&\underbrace{a_1^1\hat{\omega}a_1^2\otimes \shift a_1^3\otimes \shift a_2\otimes \dots\otimes \shift a_i\otimes 1}_{20}+\underbrace{a_1^1a_1^2\hat{\omega}\otimes \shift a_1^3\otimes \shift a_2\otimes \dots\otimes \shift a_i\otimes 1}_{22}\\
&+\underbrace{a_1^1(\shift \otimes \shift) \partial(\shift^{-1}a_1^2)\otimes \shift a_1^3\otimes \shift a_2\otimes \dots\otimes \shift a_i\otimes 1}_{21}\-\underbrace{\hat{\omega}\otimes \shift a_1\otimes \dots\otimes \shift a_{i-1}\otimes \shift a_i^1\otimes a_i^2a_i^3}_{24}\\
&-\underbrace{1\otimes \shift a_1\otimes \dots\otimes \shift a_{i-1}\otimes \shift a_i^1\otimes \hat{\omega}a_i^2a_i^3}_{30}\\
&+\underbrace{\sum_{j=1}^{i-2}1\otimes \shift a_1\otimes \dots\otimes \shift a_j\otimes 1\otimes \shift a_{j+1}\otimes \dots\otimes \shift a_{i-1}\otimes \shift a_i^1\otimes a_i^2a_i^3}_{27}+\underbrace{1\otimes \shift a_1\otimes \dots\otimes \shift a_{i-1}\otimes 1\otimes \shift a_i^1\otimes a_i^2a_i^3}_{33}\\
&+\underbrace{a_1^1a_1^2\otimes \shift a_1^3\otimes \shift a_2\otimes \dots\otimes \shift a_{i-1}\otimes \shift a_i^1\otimes a_i^2a_i^3}_{23}-\underbrace{1\otimes \shift a_1\otimes \dots\otimes \shift a_{i-1}\otimes \shift a_i^{11}\otimes a_i^{12}a_i^{13}a_i^2a_i^3}_y\\
&\boldsymbol{\uwave{-}}\underbrace{a_1^1a_1^2\hat{\omega}\otimes \shift a_1^3\otimes \shift a_2\otimes \dots\otimes \shift a_i\otimes 1}_{22}-\underbrace{a_1^1a_1^2\otimes \shift a_1^3\otimes \shift a_2\otimes \dots\otimes \shift a_i\otimes \hat{\omega}}_{25}\\
&+\underbrace{a_1^1a_1^2\otimes \shift a_1^3\otimes 1\otimes \shift a_2\otimes \dots\otimes \shift a_i\otimes 1}_{34}+\underbrace{\sum_{j=2}^{i-1}a_1^1a_1^2\otimes \shift a_1^3\otimes \shift a_2\otimes \dots\otimes \shift a_j\otimes 1\otimes \shift a_{j+1}\otimes \dots\otimes \shift a_i\otimes 1}_{26}\\
&+\underbrace{a_1^1a_1^2a_1^{31}a_1^{32}\otimes \shift a_1^{33}\otimes \shift a_2\otimes \dots\otimes \shift a_n\otimes 1}_y-\underbrace{a_1^1a_1^2\otimes \shift a_1^3\otimes \shift a_2\otimes \dots\otimes \shift a_{i-1}\otimes \shift a_i^1\otimes a_i^2a_i^3}_{23}\\
&\+\underbrace{1\otimes \shift a_1\otimes \dots\otimes \shift a_{i-1}\otimes \shift a_i^1\otimes \hat{\omega}a_i^2a_i^3}_{30}+\underbrace{1\otimes \shift a_1\otimes\dots\otimes \shift a_{i-1}\otimes \shift a_i^1\otimes a_i^2\hat{\omega}a_i^3}_{31}\\
&+\underbrace{1\otimes \shift a_1\otimes \dots\otimes \shift a_{i-1}\otimes \shift a_i^1\otimes (\shift \otimes \shift)\partial(\shift^{-1} a_i^2)a_i^3}_{32}
\end{align*}

We next compute $\mu_3\mu_1(1\otimes \shift a_1\otimes \dots\otimes \shift a_i\otimes 1)$ by applying $\mu_3^{i-1}$ to each summand in the bar resolution of $1\otimes \shift a_1\otimes \dots\otimes \shift a_i\otimes 1$ separately:
\begin{align*}
&\mu_3(a_1\otimes \shift a_2\otimes \dots\otimes \shift a_i\otimes 1)\\
=\,&\underbrace{\sum_{j=2}^{i-1} (-1)^j a_1\otimes \shift a_2\otimes \dots\otimes \shift a_{j-1}\otimes \shift a_j^1\otimes a_j^2\otimes \shift a_j^3\otimes \shift a_{j+1}\otimes \dots\otimes \shift a_i\otimes 1}_{39}\\
&+\underbrace{(-1)^{i} a_1\otimes \shift a_2\otimes \dots\otimes \shift a_{i-1}\otimes \shift a_i^1\otimes a_i^2\otimes \shift a_i^3\otimes 1}_{10}\\
&+\underbrace{a_1a_2^1a_2^{21}\dots a_2^{24}\otimes \shift a_2^{25}\otimes \shift a_2^3\otimes \shift a_3\otimes \dots\otimes \shift a_i\otimes 1}_{43}\\
&+\underbrace{(-1)^{i} a_1\otimes \shift a_2\otimes \dots\otimes \shift a_{i-1}\otimes \shift a_i^1\otimes \shift a_i^{21}\otimes a_i^{22}\dots a_i^{25}a_i^3}_{8}
\end{align*}

\begin{align*}
&\mu_3(-1\otimes \shift (a_1a_2)\otimes \shift a_3\otimes \dots\otimes \shift a_i\otimes 1)\\
=\,&-\underbrace{1\otimes \shift a_1^1\otimes a_1^2\otimes \shift (a_1^3a_2)\otimes \shift a_3\otimes \dots\otimes \shift a_i\otimes 1}_6-\underbrace{1\otimes \shift (a_1a_2^1)\otimes a_2^2\otimes \shift a_2^3\otimes \shift a_3\otimes \dots\otimes \shift a_i\otimes 1}_{40}\\
&+\underbrace{\sum_{k=3}^{i-1} (-1)^{k+1} 1\otimes \shift (a_1a_2)\otimes \dots\otimes \shift a_{k-1}\otimes \shift a_k^1\otimes a_k^2\otimes \shift a_k^3\otimes \shift a_{k+1}\otimes \dots \otimes \shift a_i\otimes 1}_{37}\\
&+\underbrace{(-1)^{i+1} 1\otimes \shift (a_1a_2)\otimes \dots\otimes \shift a_{i-1}\otimes \shift a_i^1\otimes a_i^2\otimes \shift a_i^3\otimes 1}_1\\
&-\underbrace{a_1^1a_1^{21}\dots a_1^{24}\otimes \shift a_1^{25}\otimes \shift (a_1^3a_2)\otimes \shift a_3\otimes \dots\otimes \shift a_i\otimes 1}_{11}\\
&-\underbrace{a_1a_2^1a_2^{21}\dots a_{2}^{24}\otimes \shift a_2^{25}\otimes \shift a_2^3\otimes \shift a_3\otimes \dots\otimes \shift a_i\otimes 1}_{43}\\
&+\underbrace{(-1)^{i+1} 1\otimes \shift (a_1a_2)\otimes \shift a_3\otimes \dots\otimes \shift a_{i-1}\otimes \shift a_i^1\otimes \shift a_i^{21}\otimes a_i^{22}\dots a_i^{25}a_i^3}_3
\end{align*}

\begin{align*}
&\mu_3((-1)^{i-1}1\otimes \shift a_1\otimes \dots\otimes \shift a_{i-2}\otimes \shift (a_{i-1}a_i)\otimes 1)\\
=\,&\underbrace{(-1)^{i-1}1\otimes \shift a_1^1\otimes a_1^2\otimes \shift a_1^3\otimes \shift a_2\otimes \dots\otimes \shift a_{i-2}\otimes \shift (a_{i-1}a_i)\otimes 1}_2\\
&-\underbrace{\sum_{k=2}^{i-2} (-1)^{k+i-1} 1\otimes \shift a_1\otimes \dots\otimes \shift a_{k-1}\otimes \shift a_k^1\otimes a_k^2\otimes \shift a_k^3\otimes \shift a_{k+1}\otimes\dots\otimes \shift a_{i-2}\otimes \shift (a_{i-1}a_i)\otimes 1}_{38}\\
&-\underbrace{1\otimes \shift a_1\otimes \dots\otimes \shift a_{i-2}\otimes \shift a_{i-1}^1\otimes a_{i-1}^2\otimes \shift (a_{i-1}^3 a_i)\otimes 1}_{42}\\
&-\underbrace{1\otimes \shift a_1\otimes \dots\otimes \shift a_{i-2}\otimes \shift (a_{i-1}a_i^1)\otimes a_i^2\otimes \shift a_i^3\otimes 1}_5\\
&-\underbrace{(-1)^i a_1^1a_1^{21}\dots a_1^{24} \otimes \shift a_1^{25}\otimes \shift a_1^3\otimes \shift a_2\otimes \dots\otimes \shift a_{i-2}\otimes \shift (a_{i-1}a_i)\otimes 1}_{4}\\
&-\underbrace{1\otimes \shift a_1\otimes \dots\otimes \shift a_{i-2}\otimes \shift a_{i-1}^1\otimes \shift a_{i-1}^{21}\otimes a_{i-1}^{22}\dots a_{i-1}^{25}a_{i-1}^3a_i}_{44}\\
&-\underbrace{1\otimes \shift a_1\otimes \dots\otimes \shift a_{i-2}\otimes \shift (a_{i-1}a_i^1)\otimes \shift a_i^{21}\otimes a_i^{22}\dots a_i^{25}a_i^3}_{12}
\end{align*}

\begin{align*}
&\mu_3\left(\sum_{j=2}^{i-2}(-1)^j1\otimes \shift a_1\otimes \dots\otimes \shift a_{j-1}\otimes \shift (a_ja_{j+1})\otimes \shift a_{j+2}\otimes \dots\otimes \shift a_i\otimes 1\right)\\
=&\underbrace{\sum_{j=2}^{i-2}\sum_{k=2}^{j-1} (-1)^{j+k+1} 1\otimes \shift a_1\otimes \dots\otimes \shift a_{k-1}\otimes \shift a_k^1\otimes a_k^2\otimes \shift a_k^3\otimes \shift a_{k+1}\otimes \dots\otimes \shift (a_ja_{j+1})\otimes \dots\otimes \shift a_i\otimes 1}_{38}\\
&+\underbrace{\sum_{j=2}^{i-2}\sum_{k=j+2}^{i-1} (-1)^{j+k} 1\otimes \shift a_1\otimes \dots\otimes \shift (a_ja_{j+1})\otimes \dots\otimes \shift a_{k-1}\otimes \shift a_k^1\otimes a_k^2\otimes \shift a_k^3\otimes \shift a_{k+1}\otimes \dots\otimes \shift a_i\otimes 1}_{37}\\
&-\underbrace{\sum_{j=2}^{i-2} 1\otimes \shift a_1\otimes \dots\otimes \shift a_{j-1}\otimes \shift a_j^1\otimes a_j^2\otimes \shift (a_j^3a_{j+1})\otimes \dots\otimes \shift a_i\otimes 1}_{42}\\
&-\underbrace{\sum_{j=2}^{i-2} 1\otimes \shift a_1\otimes \dots\otimes \shift a_{j-1}\otimes \shift (a_ja_{j+1}^1)\otimes a_{j+1}^2\otimes \shift a_{j+1}^3\otimes \dots\otimes \shift a_i\otimes 1}_{40}\\
&+\underbrace{\sum_{j=2}^{i-2} (-1)^j 1\otimes \shift a_1^1\otimes a_1^2\otimes \shift a_1^3\otimes \shift a_2\otimes \dots\otimes \shift (a_ja_{j+1})\otimes \dots \otimes \shift a_i\otimes 1}_2\\
&+\underbrace{\sum_{j=2}^{i-2} (-1)^{j+i} 1\otimes \shift a_1\otimes \dots\otimes \shift (a_ja_{j+1})\otimes \dots\otimes \shift a_{i-1}\otimes \shift a_i^1\otimes a_i^2\otimes \shift a_i^3\otimes 1}_1\\
&+\underbrace{\sum_{j=2}^{i-2} (-1)^{j} a_1^1 a_1^{21}\dots a_1^{24}\otimes \shift a_1^{25}\otimes \shift a_1^3\otimes \shift a_2\otimes \dots\otimes \shift (a_ja_{j+1})\otimes \dots\otimes \shift a_i\otimes 1}_{4}\\
&+\underbrace{\sum_{j=2}^{i-2}(-1)^{j+i} 1\otimes \shift a_1\otimes \dots\otimes \shift (a_ja_{j+1})\otimes \dots\otimes \shift a_{i-1}\otimes \shift a_i^1\otimes \shift a_i^{21}\otimes a_i^{22}\dots a_i^{25}a_i^3}_3
\end{align*}
and,
\begin{align*}
&\mu_3((-1)^i 1\otimes \shift a_1\otimes \dots\otimes \shift a_{i-1}\otimes a_i)\\
=\,&+\underbrace{(-1)^i 1\otimes \shift a_1^1\otimes a_1^2\otimes \shift a_1^3\otimes \shift a_2\otimes \dots\otimes \shift a_{i-1}\otimes a_i}_9\\
&+\underbrace{1\otimes \shift a_1\otimes \dots\otimes \shift a_{i-2}\otimes \shift a_{i-1}^1\otimes a_{i-1}^2\otimes \shift a_{i-1}^3\otimes a_i}_{41}\\
&+\underbrace{\sum_{j=2}^{i-2}(-1)^{i+j+1} 1\otimes \shift a_1\otimes \dots\otimes \shift a_{j-1}\otimes \shift a_j^1\otimes a_j^2\otimes \shift a_j^3\otimes \shift a_{j+1}\otimes \dots \otimes \shift a_{i-1}\otimes a_i}_{41}\\
&+\underbrace{(-1)^i a_1^1a_1^{21}\dots a_1^{24}\otimes \shift a_1^{25}\otimes \shift a_1^3\otimes a_2\otimes \dots\otimes \shift a_{i-1}\otimes a_i}_7\\
&+\underbrace{1\otimes \shift a_1\otimes \dots\otimes \shift a_{i-2}\otimes \shift a_{i-1}^1\otimes \shift a_{i-1}^{21}\otimes a_{i-1}^{22}\dots a_{i-1}^{25}a_{i-1}^3a_i}_{44}
\end{align*}

Finally, to compute $(\mu_1\otimes 1^{\otimes 2}+1\otimes \mu_1\otimes 1+1^{\otimes 2}\otimes \mu_1)\mu_3(1\otimes a_1\otimes \dots\otimes a_i\otimes 1)$ for our asserted formula for $\mu_3^i$, we consider the action of $(\mu_1\otimes 1^{\otimes 2}+1\otimes \mu_1\otimes 1+1^{\otimes 2}\otimes \mu_1)$ on each summand in turn and deal with the cases $j=1$ and $j=i$ separately. For $j=1$ we obtain
\begin{align*}
&(\mu_1\otimes 1^{\otimes 2}+1\otimes \mu_1\otimes 1+1^{\otimes 2}\otimes \mu_1)(1\otimes \shift a_1^1\otimes a_1^2\otimes \shift a_1^3\otimes \shift a_2\otimes \dots\otimes \shift a_i\otimes 1)\\
=\,&\underbrace{\hat{\omega}a_1^1a_1^2\otimes \shift a_1^3\otimes \shift a_2\otimes \dots\otimes \shift a_i\otimes 1}_{28}-\underbrace{a_1^1\hat{\omega} a_1^2\otimes \shift a_1^3\otimes \shift a_2\otimes \dots\otimes \shift a_i\otimes 1}_{20}\\
&-\underbrace{\shift \partial(a_1^1)a_1^2\otimes \shift a_1^3\otimes \shift a_2\otimes \dots\otimes \shift a_i\otimes 1}_y\uwave{\+}0\\
&\boldsymbol{\uwave{-}}\underbrace{1\otimes \shift a_1^1\otimes a_1^2a_1^3\otimes \shift a_2\otimes \dots\otimes \shift a_i\otimes 1}_{35}+\underbrace{1\otimes \shift a_1^1\otimes a_1^2\otimes \shift (a_1^3a_2)\otimes \shift a_3\otimes \dots\otimes \shift a_i\otimes 1}_6\\
&+\underbrace{\sum_{j=2}^{i-1}(-1)^{j+1}1\otimes \shift a_1^1\otimes a_1^2\otimes a_1^3\otimes \shift a_2\otimes \dots \otimes \shift a_{j-1} \otimes \shift (a_ja_{j+1})\otimes \shift a_{j+2}\otimes \dots\otimes \shift a_i\otimes 1}_2\\
&+\underbrace{(-1)^{i+1}1\otimes \shift a_1^1\otimes a_1^2\otimes \shift a_1^3\otimes \shift a_2\otimes \dots \otimes \shift a_{i-1}\otimes a_i}_9
\end{align*}
For $j=i$ we obtain
\begin{align*}
&(\mu_1\otimes 1^{\otimes 2}+1\otimes \mu_1\otimes 1+1^{\otimes 2}\otimes \mu_1)((-1)^{i+1}1\otimes \shift a_1\otimes \dots\otimes \shift a_{i-1}\otimes \shift a_i^1\otimes a_i^2\otimes \shift a_i^3\otimes 1)\\
=\,&\underbrace{(-1)^{i+1}a_1\otimes \shift a_2\otimes \dots\otimes \shift a_{i-1}\otimes \shift a_i^1\otimes a_i^2\otimes \shift a_i^3\otimes 1}_{10}\\
&+\underbrace{\sum_{j=1}^{i-2}(-1)^{j+i+1}1\otimes \shift a_1\otimes \dots\otimes \shift a_{j-1}\otimes \shift (a_ja_{j+1})\otimes \shift a_{j+2}\otimes \dots\otimes \shift a_{i-1}\otimes \shift a_i^1\otimes a_i^2\otimes a_i^3\otimes 1}_1\\
&+\underbrace{1\otimes \shift a_1\otimes \dots\otimes \shift a_{i-2}\otimes \shift (a_{i-1}a_i^1)\otimes a_i^2\otimes \shift a_i^3\otimes 1}_5-\underbrace{1\otimes \shift a_1\otimes \dots\otimes \shift a_{i-1}\otimes a_i^1a_i^2\otimes \shift a_i^3\otimes 1}_{36}\uwave{\+}0\\
&\boldsymbol{\uwave{-}}\underbrace{1\otimes \shift a_1\otimes \dots\otimes \shift a_{i-1}\otimes \shift a_i ^1\otimes a_i^2\hat{\omega}a_i^3}_{31}+\underbrace{1\otimes \shift a_1\otimes \dots\otimes \shift a_{i-1}\otimes \shift a_i^1\otimes a_i^2a_i^3\hat{\omega}}_{29}\\
&+\underbrace{1\otimes \shift a_1\otimes \dots\otimes \shift a_{i-1}\otimes \shift a_i^1\otimes a_i^2\shift \partial(a_i^3).}_y
\end{align*}
Now, for $j\neq 1,i$ we obtain
\begin{align*}
&(\mu_1\otimes 1^{\otimes 2}+1\otimes \mu_1\otimes 1+1^{\otimes 2}\otimes \mu_1)\left(\sum_{j=2}^{i-1} (-1)^{j+1} 1\otimes \shift a_1\otimes \dots\otimes \shift a_{j-1}\otimes \shift a_j^1\otimes a_j^2\otimes \shift a_j^3\otimes \shift a_{j+1}\otimes\dots\otimes \shift a_i\otimes 1\right)\\
=\,&\underbrace{\sum_{j=2}^{i-1}(-1)^{j+1}a_1\otimes \shift a_2\otimes \dots\otimes \shift a_{j-1}\otimes \shift a_j^1\otimes a_j^2\otimes \shift a_j^3\otimes \shift a_{j+1}\otimes \dots\otimes \shift a_i\otimes 1}_{39}\\
&+\underbrace{\sum_{j=2}^{i-1}\sum_{k=1}^{j-2}(-1)^{j+k+1} 1\otimes \shift a_1\otimes \dots\otimes \shift (a_ka_{k+1})\otimes \dots\otimes \shift a_{j-1}\otimes \shift a_j^1\otimes a_j^2\otimes \shift a_j^3\otimes a_{j+1}\otimes \dots\otimes \shift a_i\otimes 1}_{37}\\
&+\underbrace{\sum_{j=2}^{i-1} 1\otimes \shift a_1\otimes \dots\otimes \shift (a_{j-1}a_j^1)\otimes a_j^2\otimes a_j^3\otimes \shift a_{j+1}\otimes \dots\otimes \shift a_i\otimes 1}_{40}\\
&-\underbrace{\sum_{j=2}^{i-1} 1\otimes \shift a_1\otimes \dots\otimes \shift a_{j-1}\otimes a_j^1a_j^2\otimes \shift a_j^3\otimes \shift a_{j+1}\otimes \dots\otimes \shift a_i\otimes 1}_{36}\uwave{\+}0\\
&\boldsymbol{\uwave{-}}\underbrace{\sum_{j=2}^{i-1}1\otimes \shift a_1\otimes \dots\otimes \shift a_{j-1}\otimes \shift a_j^1\otimes a_j^2a_j^3\otimes \shift a_{j+1}\otimes \dots\otimes \shift a_i\otimes 1}_{35}\\
&+\underbrace{\sum_{j=2}^{i-1} 1\otimes \shift a_1\otimes \dots\otimes \shift a_{j-1}\otimes \shift a_j^1\otimes a_j^2\otimes \shift(a_j^3a_{j+1})\otimes \dots\otimes \shift a_i\otimes 1}_{42}\\
&+\underbrace{\sum_{j=2}^{i-1}\sum_{k=j+1}^{i-1} (-1)^{j+k}1\otimes \shift a_1\otimes \dots\otimes \shift a_{j-1}\otimes \shift a_j^1\otimes a_j^2\otimes \shift a_j^3\otimes \shift a_{j+1}\otimes \dots\otimes \shift a_{k-1}\otimes \shift (a_ka_{k+1})\otimes \shift a_{k+2}\otimes \dots\otimes \shift a_i\otimes 1}_{38}\\
&+\underbrace{\sum_{j=2}^{i-1}(-1)^{i+j} 1\otimes \shift a_1\otimes \dots\otimes \shift a_{j-1}\otimes \shift a_j^1\otimes a_j^2\otimes \shift a_j^3\otimes \shift a_{j+1}\otimes \dots\otimes \shift a_{i-1}\otimes a_i.}_{41}
\end{align*}
Furthermore,
\begin{align*}
&(\mu_1\otimes 1^{\otimes 2}+1\otimes \mu_1\otimes 1+1^{\otimes 2}\otimes \mu_1)(a_1^1a_1^{21}\dots a_1^{24}\otimes \shift a_1^{25}\otimes \shift a_1^3\otimes \shift a_2\otimes \dots\otimes \shift a_i\otimes 1)\\
=&0\uwave{\+}0\boldsymbol{\uwave{-}}\underbrace{a_1^1(\shift \otimes \shift)\partial(\shift^{-1}a_1^2)\otimes \shift a_1^3\otimes \shift a_2\otimes \dots\otimes \shift a_i\otimes 1}_{21}-\underbrace{a_1^1a_1^{21}\dots a_1^{24}\otimes \shift (a_1^{25}a_1^3)\otimes a_2\otimes \dots\otimes \shift a_i\otimes 1}_y\\
&+\underbrace{a_1^1a_1^{21}\dots a_1^{24}\otimes \shift a_1^{25}\otimes \shift(a_1^3a_2)\otimes \shift a_3\otimes \dots\otimes \shift a_i\otimes 1}_{11}\\
&-\underbrace{\sum_{j=2}^{i-1} (-1)^j a_1^1a_1^{21}\dots a_1^{24}\otimes \shift a_1^{25} \otimes \shift a_1^3\otimes a_2\otimes \dots\otimes \shift a_{j-1}\otimes \shift (a_ja_{j+1})\otimes \shift a_{j+2}\otimes \dots\otimes \shift a_i\otimes 1}_4\\
&-\underbrace{(-1)^i a_1^1a_1^{21}\dots a_1^{24}\otimes \shift a_1^{25}\otimes \shift a_1^3\otimes a_2\otimes \dots\otimes \shift a_{i-1}\otimes a_i}_7
\end{align*}
and, finally,
\begin{align*}
&(\mu_1\otimes 1^{\otimes 2}+1\otimes \mu_1\otimes 1+1^{\otimes 2}\otimes \mu_1)((-1)^{i+1} 1\otimes \shift a_1\otimes \dots\otimes \shift a_{i-1}\otimes \shift a_i^1\otimes a_i^{21}\otimes a_i^{22}\dots a_i^{25}a_i^3)\\
=\,&\underbrace{(-1)^{i+1} a_1\otimes \shift a_2\otimes \dots\otimes \shift a_{i-1}\otimes \shift a_i^1\otimes \shift a_i^{21}\otimes a_i^{22}\dots a_i^{25}a_i^3}_8\\
&+\underbrace{\sum_{j=1}^{i-2} (-1)^{j+i+1} 1\otimes \shift a_1\otimes \dots\otimes \shift a_{j-1}\otimes \shift (a_ja_{j+1})\otimes \shift a_{j+2}\otimes \dots\otimes \shift a_{i-1}\otimes \shift a_i^1\otimes \shift a_i^{21}\otimes a_i^{22}\dots a_i^{25}a_i^3}_3\\
&+\underbrace{1\otimes \shift a_1\otimes \dots\otimes \shift a_{i-2}\otimes \shift (a_{i-1}a_i^1)\otimes \shift a_i^{21}\otimes a_i^{22}\dots a_i^{25} a_i^3}_{12}\\
&-\underbrace{1\otimes \shift a_1\otimes \dots\otimes \shift a_{i-1}\otimes \shift (a_i^1a_i^{21})\otimes a_i^{22}\dots a_i^{25} a_i^3}_y\\
&-\underbrace{1\otimes \shift a_1\otimes \dots\otimes \shift a_{i-1}\otimes \shift a_i^1\otimes (\shift \otimes \shift)\partial(\shift^{-1}a_i^2) a_i^3}_{32}\uwave{\+}0\uwave{\+}0.
\end{align*}
Again, terms with the same label add up to zero, completing the proof.
\end{proof}

This concludes the computation for $\mu_3$.

\section{The maps $\mu_4^0, \mu_4^1, \mu_4^2$ on $\mathcal{P}$}
\label{appendixC}

Similarly to Lemma \ref{yellow}, we start by providing another useful tool to identify vanishing terms in the remainder.

\begin{lem}\label{green}
For $a \in A$, the following equalities hold:

\begin{align*}
a^{111}\otimes \shift^{-1}(a^{112})\otimes a^{113}\otimes \shift^{-1}(a^{12})\otimes a^{13}\otimes \shift^{-1}(a^{2})\otimes a^3&=0,\\
a^{11}\otimes a^{121}\otimes \shift^{-1}(a^{122})\otimes a^{123}\otimes \shift^{-1}(a^{124})\otimes a^{125}\otimes a^{13}\otimes \shift^{-1}(a^{2})\otimes a^3&=0,\\
a^{11}\otimes \shift^{-1}(a^{12})\otimes a^{131}\otimes \shift^{-1}(a^{132})\otimes a^{133}\otimes \shift^{-1}(a^{2})\otimes a^3&=0,\\
a^{11}\otimes \shift^{-1}(a^{12})\otimes a^{13}\otimes a^{21}\otimes \shift^{-1}(a^{22})\otimes a^{23}\otimes \shift^{-1}(a^{24})\otimes a^{25}\otimes a^3&=0,\\
a^{11}\otimes \shift^{-1}(a^{12})\otimes a^{13}\otimes \shift^{-1}(a^{2})\otimes a^{31}\otimes \shift^{-1}(a^{32})\otimes a^{33}&=0;
\end{align*}

\begin{align*}
a^{11}\otimes  \shift^{-1}(a^{12})\otimes a^{13}a^{21}\otimes \shift^{-1}(a^{22})\otimes a^{23}\otimes \shift^{-1}(a^{24})\otimes a^{25}a^3&=0,\\
a^1a^{211}\otimes \shift^{-1}(a^{212})\otimes a^{213}\otimes \shift^{-1}(a^{22})\otimes a^{23}\otimes \shift^{-1}(a^{24})\otimes a^{25}a^3&=0,\\
a^1a^{21}\otimes a^{221}\otimes \shift^{-1}(a^{222})\otimes a^{223}\shift^{-1}(a^{224})\otimes a^{225}\otimes a^{23}\otimes \shift^{-1}(a^{24})\otimes a^{25}a^3&=0,\\
a^1a^{21}\otimes \shift^{-1}(a^{22})\otimes a^{231}\otimes \shift^{-1}(a^{232})\otimes a^{233}\otimes \shift^{-1}(a^{24})\otimes a^{25}a^3&=0,\\
a^1a^{21}\otimes \shift^{-1}(a^{22})\otimes a^{23}\otimes a^{241}\otimes \shift^{-1}(a^{242})\otimes a^{243}\otimes \shift^{-1}(a^{244})\otimes a^{245}\otimes a^{25}a^3&=0,\\
a^1a^{21}\otimes \shift^{-1}(a^{22})\otimes a^{23}\otimes \shift^{-1}(a^{24})\otimes a^{251}\otimes \shift^{-1}(a^{252})\otimes a^{253}a^3&=0,\\
a^1a^{21}\otimes \shift^{-1}(a^{22})\otimes a^{23}\otimes \shift^{-1}(a^{24})\otimes a^{25}a^{31}\otimes \shift^{-1}(a^{32})\otimes a^{33}&=0;
\end{align*}

\begin{align*}
a^{11} \otimes \shift^{-1}(a^{12})\otimes a^{13}\otimes \shift^{-1}(a^2)\otimes a^{31}\otimes \shift^{-1}(a^{32})\otimes a^{33}&=0,\\
a^1\otimes a^{21}\otimes \shift^{-1}(a^{22})\otimes a^{23}\otimes \shift^{-1}(a^{24})\otimes a^{25}\otimes a^{31}\otimes \shift^{-1}(a^{32})\otimes a^{33}&=0,\\
a^1\otimes \shift^{-1}(a^2)\otimes a^{311}\otimes \shift^{-1}(a^{312})\otimes a^{313}\otimes \shift^{-1}(a^{32})\otimes a^{33}&=0,\\
a^1\otimes \shift^{-1}(a^2)\otimes a^{31}\otimes a^{321}\otimes \shift^{-1}(a^{322})\otimes a^{323}\otimes \shift^{-1}(a^{324})\otimes a^{325}\otimes a^{33}&=0,\\
a^1\otimes \shift^{-1}(a^2)\otimes a^{31}\otimes \shift^{-1}(a^{32})\otimes a^{331}\otimes \shift^{-1}(a^{332})\otimes a^{333}&=0.
\end{align*}

For $\varphi\in\Phi$, the following equalities hold:

\begin{align*}
\varphi^{111}\otimes (\shift^{-1}\varphi^{112})\otimes \varphi^{113}\otimes(\shift^{-1}\varphi^{12})\otimes\varphi^{13}\otimes(\shift^{-1}\varphi^2)\otimes\varphi^3\otimes(\shift^{-1}\varphi^4)\otimes\varphi^5&=0,\\
\varphi^{11}\otimes\varphi^{121}\otimes (\shift^{-1}\varphi^{122})\otimes \varphi^{123}\otimes (\shift^{-1}\varphi^{124})\otimes \varphi^{125}\otimes\varphi^{13}\otimes(\shift^{-1}\varphi^2)\otimes\varphi^3\otimes(\shift^{-1}\varphi^4)\otimes\varphi^5 &=0,\\
\varphi^{11}\otimes(\shift^{-1}\varphi^{12})\otimes\varphi^{131}\otimes (\shift^{-1}\varphi^{132})\varphi^{133}\otimes \otimes(\shift^{-1}\varphi^2)\otimes\varphi^3\otimes(\shift^{-1}\varphi^4)\otimes\varphi^5 &=0,\\
\varphi^{11}\otimes(\shift^{-1}\varphi^{12})\otimes\varphi^{13}\otimes\varphi^{21}\otimes (\shift^{-1}\varphi^{22})\otimes \varphi^{23}\otimes (\shift^{-1}\varphi^{24})\otimes \varphi^{25}\otimes\varphi^3\otimes(\shift^{-1}\varphi^4)\otimes\varphi^5 &=0,\\
\varphi^{11}\otimes(\shift^{-1}\varphi^{12})\otimes \varphi^{13} \otimes(\shift^{-1}\varphi^2)\otimes\varphi^{31}\otimes (\shift^{-1}\otimes \varphi^{32})\varphi^{33}\otimes(\shift^{-1}\varphi^4)\otimes\varphi^5 &=0,\\
\varphi^{11}\otimes(\shift^{-1}\varphi^{12})\otimes\varphi^{13}\otimes(\shift^{-1}\varphi^2)\otimes\varphi^3\otimes\varphi^{41}\otimes (\shift^{-1}\varphi^{42})\varphi^{43}\otimes (\shift^{-1}\varphi^{44})\otimes \varphi^{45}\otimes\varphi^5 &=0,\\
\varphi^{11}\otimes(\shift^{-1}\varphi^{12})\otimes\varphi^{13}\otimes(\shift^{-1}\varphi^2)\otimes\varphi^3\otimes(\shift^{-1}\varphi^4)\otimes\varphi^{51}\otimes (\shift^{-1}\varphi^{52})\otimes \varphi^{53} &=0;
\end{align*}

\begin{align*}
\varphi^{11}\otimes (\shift^{-1}\varphi^{12})\otimes \varphi^{13}\varphi^{21}\otimes(\shift^{-1}\varphi^{22})\otimes\varphi^{23}\otimes(\shift^{-1}\varphi^{24})\otimes\varphi^{25}\varphi^3\otimes(\shift^{-1}\varphi^4)\otimes\varphi^5 &=0,\\
\varphi^1\varphi^{211}\otimes (\shift^{-1}\varphi^{212})\otimes \varphi^{213}\otimes(\shift^{-1}\varphi^{22})\otimes\varphi^{23}\otimes(\shift^{-1}\varphi^{24})\otimes\varphi^{25}\varphi^3\otimes(\shift^{-1}\varphi^4)\otimes\varphi^5 &=0,\\
\varphi^1\varphi^{21}\otimes\varphi^{221}\otimes (\shift^{-1}\varphi^{222})\otimes \varphi^{223}\otimes (\shift^{-1}\varphi^{224})\otimes \varphi^{225}\otimes\varphi^{23}\otimes(\shift^{-1}\varphi^{24})\otimes\varphi^{25}\varphi^3\otimes(\shift^{-1}\varphi^4)\otimes\varphi^5 &=0,\\
\varphi^1\varphi^{21}\otimes(\shift^{-1}\varphi^{22})\otimes\varphi^{231}\otimes (\shift^{-1}\varphi^{232})\otimes \varphi^{233}\otimes(\shift^{-1}\varphi^{24})\otimes\varphi^{25}\varphi^3\otimes(\shift^{-1}\varphi^4)\otimes\varphi^5 &=0,\\
\varphi^1\varphi^{21}\otimes(\shift^{-1}\varphi^{22})\otimes\varphi^{23}\otimes\varphi^{241}\otimes (\shift^{-1}\varphi^{242})\otimes \varphi^{243}\otimes (\shift^{-1}\varphi^{244})\otimes \varphi^{245}\otimes\varphi^{25}\varphi^3\otimes(\shift^{-1}\varphi^4)\otimes\varphi^5 &=0,\\
\varphi^1\varphi^{21}\otimes(\shift^{-1}\varphi^{22})\otimes\varphi^{23}\otimes(\shift^{-1}\varphi^{24})\otimes\varphi^{251}\otimes (\shift^{-1}\varphi^{252})\otimes \varphi^{253}\varphi^3\otimes(\shift^{-1}\varphi^4)\otimes\varphi^5 &=0,\\
\varphi^1\varphi^{21}\otimes(\shift^{-1}\varphi^{22})\otimes\varphi^{23}\otimes(\shift^{-1}\varphi^{24})\otimes\varphi^{25}\varphi^{31}\otimes (\shift^{-1}\varphi^{32})\otimes \varphi^{33}\otimes(\shift^{-1}\varphi^4)\otimes\varphi^5 &=0,\\
\varphi^1\varphi^{21}\otimes(\shift^{-1}\varphi^{22})\otimes\varphi^{23}\otimes(\shift^{-1}\varphi^{24})\otimes\varphi^{25}\varphi^3\otimes\varphi^{41}\otimes (\shift^{-1}\varphi^{42})\otimes \varphi^{43}\otimes (\shift^{-1}\varphi^{44})\otimes \varphi^{45}\otimes\varphi^5 &=0,\\
\varphi^1\varphi^{21}\otimes(\shift^{-1}\varphi^{22})\otimes\varphi^{23}\otimes(\shift^{-1}\varphi^{24})\otimes\varphi^{25}\varphi^3\otimes(\shift^{-1}\varphi^4)\otimes\varphi^{51}\otimes (\shift^{-1}\varphi^{52})\otimes \varphi^{53} &=0;
\end{align*}

\begin{align*}
\varphi^{11}\otimes (\shift^{-1}\varphi^{12})\otimes \varphi^{13}\otimes(\shift^{-1}\varphi^2)\otimes\varphi^{31}\otimes(\shift^{-1}\varphi^{32})\otimes\varphi^{33}\otimes(\shift^{-1}\varphi^4)\otimes\varphi^5  &=0,\\
\varphi^1\otimes\varphi^{21}\otimes (\shift^{-1}\varphi^{22})\otimes \varphi^{23}\otimes (\shift^{-1}\varphi^{24})\otimes \varphi^{25}\otimes\varphi^{31}\otimes(\shift^{-1}\varphi^{32})\otimes\varphi^{33}\otimes(\shift^{-1}\varphi^4)\otimes\varphi^5  &=0,\\
\varphi^1\otimes(\shift^{-1}\varphi^2)\otimes\varphi^{311}\otimes (\shift^{-1}\varphi^{312})\otimes \varphi^{313}\otimes(\shift^{-1}\varphi^{32})\otimes\varphi^{33}\otimes(\shift^{-1}\varphi^4)\otimes\varphi^5  &=0,\\
\varphi^1\otimes(\shift^{-1}\varphi^2)\otimes\varphi^{31}\otimes\varphi^{321}\otimes (\shift^{-1}\varphi^{322})\otimes \varphi^{323}\otimes (\shift^{-1}\varphi^{324})\otimes \varphi^{325}\otimes\varphi^{33}\otimes(\shift^{-1}\varphi^4)\otimes\varphi^5  &=0,\\
\varphi^1\otimes(\shift^{-1}\varphi^2)\otimes\varphi^{31}\otimes(\shift^{-1}\varphi^{32})\otimes\varphi^{331}\otimes (\shift^{-1}\varphi^{332})\otimes \varphi^{333}\otimes(\shift^{-1}\varphi^4)\otimes\varphi^5  &=0,\\
\varphi^1\otimes(\shift^{-1}\varphi^2)\otimes\varphi^{31}\otimes(\shift^{-1}\varphi^{32})\otimes\varphi^{33}\otimes\varphi^{41}\otimes (\shift^{-1}\varphi^{42})\otimes \varphi^{43}\otimes (\shift^{-1}\varphi^{44})\otimes \varphi^{45}\otimes\varphi^5  &=0,\\
\varphi^1\otimes(\shift^{-1}\varphi^2)\otimes\varphi^{31}\otimes(\shift^{-1}\varphi^{32})\otimes\varphi^{33}\otimes(\shift^{-1}\varphi^4)\otimes\varphi^{51}\otimes (\shift^{-1}\varphi^{52})\otimes \varphi^{53}  &=0;
\end{align*}

\begin{align*}
\varphi^{11}\otimes (\shift^{-1}\varphi^{12})\otimes \varphi^{13}\otimes(\shift^{-1}\varphi^2)\otimes\varphi^3\varphi^{41}\otimes(\shift^{-1}\varphi^{42})\otimes\varphi^{43}\otimes(\shift^{-1}\varphi^{44})\otimes\varphi^{45}\varphi^5 &=0,\\
\varphi^1\otimes\varphi^{21}\otimes (\shift^{-1}\varphi^{22})\otimes \varphi^{23}\otimes (\shift^{-1}\varphi^{24})\otimes \varphi^{25}\otimes\varphi^3\varphi^{41}\otimes(\shift^{-1}\varphi^{42})\otimes\varphi^{43}\otimes(\shift^{-1}\varphi^{44})\otimes\varphi^{45}\varphi^5 &=0,\\
\varphi^1\otimes(\shift^{-1}\varphi^2)\otimes\varphi^{31}\otimes (\shift^{-1}\varphi^{32})\otimes \varphi^{33}\varphi^{41}\otimes(\shift^{-1}\varphi^{42})\otimes\varphi^{43}\otimes(\shift^{-1}\varphi^{44})\otimes\varphi^{45}\varphi^5 &=0,\\
\varphi^1\otimes(\shift^{-1}\varphi^2)\otimes\varphi^3\varphi^{411}\otimes (\shift^{-1}\varphi^{412})\otimes \varphi^{413}\otimes(\shift^{-1}\varphi^{42})\otimes\varphi^{43}\otimes(\shift^{-1}\varphi^{44})\otimes\varphi^{45}\varphi^5 &=0,\\
\varphi^1\otimes(\shift^{-1}\varphi^2)\otimes\varphi^3\varphi^{41}\otimes\varphi^{421}\otimes (\shift^{-1}\varphi^{422})\otimes \varphi^{423}\otimes (\shift^{-1}\varphi^{424})\otimes \varphi^{425}\otimes\varphi^{43}\otimes(\shift^{-1}\varphi^{44})\otimes\varphi^{45}\varphi^5 &=0,\\
\varphi^1\otimes(\shift^{-1}\varphi^2)\otimes\varphi^3\varphi^{41}\otimes(\shift^{-1}\varphi^{42})\otimes\varphi^{431}\otimes (\shift^{-1}\varphi^{432})\otimes \varphi^{433}\otimes(\shift^{-1}\varphi^{44})\otimes\varphi^{45}\varphi^5 &=0,\\
\varphi^1\otimes(\shift^{-1}\varphi^2)\otimes\varphi^3\varphi^{41}\otimes(\shift^{-1}\varphi^{42})\otimes\varphi^{43}\otimes\varphi^{441}\otimes (\shift^{-1}\varphi^{442})\otimes \varphi^{443}\otimes (\shift^{-1}\varphi^{444})\otimes \varphi^{445}\otimes\varphi^{45}\varphi^5 &=0,\\
\varphi^1\otimes(\shift^{-1}\varphi^2)\otimes\varphi^3\varphi^{41}\otimes(\shift^{-1}\varphi^{42})\otimes\varphi^{43}\otimes(\shift^{-1}\varphi^{44})\otimes\varphi^{451}\otimes (\shift^{-1}\varphi^{452})\otimes \varphi^{453} \varphi^5 &=0,\\
\varphi^1\otimes(\shift^{-1}\varphi^2)\otimes\varphi^3\varphi^{41}\otimes(\shift^{-1}\varphi^{42})\otimes\varphi^{43}\otimes(\shift^{-1}\varphi^{44})\otimes\varphi^{45}\varphi^{51}\otimes (\shift^{-1}\varphi^{52})\otimes \varphi^{53}  &=0;
\end{align*}

\begin{align*}
\varphi^{11}\otimes (\shift^{-1}\varphi^{12})\otimes \varphi^{13}\otimes(\shift^{-1}\varphi^2)\otimes\varphi^3\otimes(\shift^{-1}\varphi^4)\otimes\varphi^{51}\otimes(\shift^{-1}\varphi^{52})\otimes\varphi^{53} &=0,\\
\varphi^1\otimes\varphi^{21}\otimes (\shift^{-1}\varphi^{22})\otimes \varphi^{23}\otimes (\shift^{-1}\varphi^{24})\otimes \varphi^{25}\otimes\varphi^3\otimes(\shift^{-1}\varphi^4)\otimes\varphi^{51}\otimes(\shift^{-1}\varphi^{52})\otimes\varphi^{53} &=0,\\
\varphi^1\otimes(\shift^{-1}\varphi^2)\otimes\varphi^{31}\otimes (\shift^{-1}\varphi^{32})\otimes \varphi^{33}\otimes(\shift^{-1}\varphi^4)\otimes\varphi^{51}\otimes(\shift^{-1}\varphi^{52})\otimes\varphi^{53} &=0,\\
\varphi^1\otimes(\shift^{-1}\varphi^2)\otimes\varphi^3\otimes\varphi^{41}\otimes (\shift^{-1}\varphi^{42})\otimes \varphi^{43}\otimes (\shift^{-1}\varphi^{44})\otimes \varphi^{45}\otimes\varphi^{51}\otimes(\shift^{-1}\varphi^{52})\otimes\varphi^{53} &=0,\\
\varphi^1\otimes(\shift^{-1}\varphi^2)\otimes\varphi^3\otimes(\shift^{-1}\varphi^4)\otimes\varphi^{511}\otimes (\shift^{-1}\varphi^{512})\otimes \varphi^{513} \otimes(\shift^{-1}\varphi^{52})\otimes\varphi^{53} &=0,\\
\varphi^1\otimes(\shift^{-1}\varphi^2)\otimes\varphi^3\otimes(\shift^{-1}\varphi^4)\otimes\varphi^{51}\otimes\varphi^{521}\otimes (\shift^{-1}\varphi^{522})\otimes \varphi^{523}\otimes (\shift^{-1}\varphi^{524})\otimes \varphi^{525}\otimes\varphi^{53} &=0,\\
\varphi^1\otimes(\shift^{-1}\varphi^2)\otimes\varphi^3\otimes(\shift^{-1}\varphi^4)\otimes\varphi^{51}\otimes(\shift^{-1}\varphi^{52})\otimes\varphi^{531}\otimes (\shift^{-1}\varphi^{532})\otimes \varphi^{533} &=0.
\end{align*}

All statements again hold with any number of tensor symbols removed.
\end{lem}
\begin{proof}

Differentiating the terms proven to be zero in Lemma \ref{yellow}, and using the same argument of directedness of the biquiver as there, we obtain that every summand in these differentials is zero. Listing them, we obtain 

\begin{align*}
a^{111}\shift^{-1}(a^{112})a^{113}\otimes \shift^{-1}(a^{12})\otimes a^{13}\otimes \shift^{-1}(a^{2})\otimes a^3&=0,\\
a^{11}\otimes a^{121}\shift^{-1}(a^{122})a^{123}\shift^{-1}(a^{124})a^{125}\otimes a^{13}\otimes \shift^{-1}(a^{2})\otimes a^3&=0,\\
a^{11}\otimes \shift^{-1}(a^{12})\otimes a^{131}\shift^{-1}(a^{132})a^{133}\otimes \shift^{-1}(a^{2})\otimes a^3&=0,\\
a^{11}\otimes \shift^{-1}(a^{12})\otimes a^{13}\otimes a^{21}\shift^{-1}(a^{22})a^{23}\shift^{-1}(a^{24})a^{25}\otimes a^3&=0,\\
a^{11}\otimes \shift^{-1}(a^{12})\otimes a^{13}\otimes \shift^{-1}(a^{2})\otimes a^{31}\shift^{-1}(a^{32})a^{33}&=0;
\end{align*}

\begin{align*}
a^{11} \shift^{-1}(a^{12})a^{13}a^{21}\otimes \shift^{-1}(a^{22})\otimes a^{23}\otimes \shift^{-1}(a^{24})\otimes a^{25}a^3&=0,\\
a^1a^{211}\shift^{-1}(a^{212})a^{213}\otimes \shift^{-1}(a^{22})\otimes a^{23}\otimes \shift^{-1}(a^{24})\otimes a^{25}a^3&=0,\\
a^1a^{21}\otimes a^{221}\shift^{-1}(a^{222})a^{223}\shift^{-1}(a^{224})a^{225}\otimes a^{23}\otimes \shift^{-1}(a^{24})\otimes a^{25}a^3&=0,\\
a^1a^{21}\otimes \shift^{-1}(a^{22})\otimes a^{231}\shift^{-1}(a^{232})a^{233}\otimes \shift^{-1}(a^{24})\otimes a^{25}a^3&=0,\\
a^1a^{21}\otimes \shift^{-1}(a^{22})\otimes a^{23}\otimes a^{241}\shift^{-1}(a^{242})a^{243}\shift^{-1}(a^{244})a^{245}\otimes a^{25}a^3&=0,\\
a^1a^{21}\otimes \shift^{-1}(a^{22})\otimes a^{23}\otimes \shift^{-1}(a^{24})\otimes a^{251}\shift^{-1}(a^{252})a^{253}a^3&=0,\\
a^1a^{21}\otimes \shift^{-1}(a^{22})\otimes a^{23}\otimes \shift^{-1}(a^{24})\otimes a^{25}a^{31}\shift^{-1}(a^{32})a^{33}&=0;
\end{align*}

\begin{align*}
a^{11} \shift^{-1}(a^{12})a^{13}\otimes \shift^{-1}(a^2)\otimes a^{31}\otimes \shift^{-1}(a^{32})\otimes a^{33}&=0,\\
a^1\otimes a^{21}\shift^{-1}(a^{22})a^{23}\shift^{-1}(a^{24})a^{25}\otimes a^{31}\otimes \shift^{-1}(a^{32})\otimes a^{33}&=0,\\
a^1\otimes \shift^{-1}(a^2)\otimes a^{311}\shift^{-1}(a^{312})a^{313}\otimes \shift^{-1}(a^{32})\otimes a^{33}&=0,\\
a^1\otimes \shift^{-1}(a^2)\otimes a^{31}\otimes a^{321}\shift^{-1}(a^{322})a^{323}\shift^{-1}(a^{324})a^{325}\otimes a^{33}&=0,\\
a^1\otimes \shift^{-1}(a^2)\otimes a^{31}\otimes \shift^{-1}(a^{32})\otimes a^{331}\shift^{-1}(a^{332})a^{333}&=0.
\end{align*}

For $\varphi\in \Phi$, using the same arguments, we obtain, from

\[\varphi^{11}\otimes(\shift^{-1}\varphi^{12})\otimes\varphi^{13}\otimes(\shift^{-1}\varphi^2)\otimes\varphi^3\otimes(\shift^{-1}\varphi^4)\otimes\varphi^5 =0.\]

\begin{align*}
\varphi^{111}(\shift^{-1}\varphi^{112})\varphi^{113}\otimes(\shift^{-1}\varphi^{12})\otimes\varphi^{13}\otimes(\shift^{-1}\varphi^2)\otimes\varphi^3\otimes(\shift^{-1}\varphi^4)\otimes\varphi^5 &=0,\\
\varphi^{11}\otimes\varphi^{121}(\shift^{-1}\varphi^{122})\varphi^{123}(\shift^{-1}\varphi^{124})\varphi^{125}\otimes\varphi^{13}\otimes(\shift^{-1}\varphi^2)\otimes\varphi^3\otimes(\shift^{-1}\varphi^4)\otimes\varphi^5 &=0,\\
\varphi^{11}\otimes(\shift^{-1}\varphi^{12})\otimes\varphi^{131}(\shift^{-1}\varphi^{132})\varphi^{133}\otimes(\shift^{-1}\varphi^2)\otimes\varphi^3\otimes(\shift^{-1}\varphi^4)\otimes\varphi^5 &=0,\\
\varphi^{11}\otimes(\shift^{-1}\varphi^{12})\otimes\varphi^{13}\otimes\varphi^{21}(\shift^{-1}\varphi^{22})\varphi^{23}(\shift^{-1}\varphi^{24})\varphi^{25}\otimes\varphi^3\otimes(\shift^{-1}\varphi^4)\otimes\varphi^5 &=0,\\
\varphi^{11}\otimes(\shift^{-1}\varphi^{12})\otimes \varphi^{13} \otimes(\shift^{-1}\varphi^2)\otimes\varphi^{31}(\shift^{-1}\varphi^{32})\varphi^{33}\otimes(\shift^{-1}\varphi^4)\otimes\varphi^5 &=0,\\
\varphi^{11}\otimes(\shift^{-1}\varphi^{12})\otimes\varphi^{13}\otimes(\shift^{-1}\varphi^2)\otimes\varphi^3\otimes\varphi^{41}(\shift^{-1}\varphi^{42})\varphi^{43}(\shift^{-1}\varphi^{44})\varphi^{45}\otimes\varphi^5 &=0,\\
\varphi^{11}\otimes(\shift^{-1}\varphi^{12})\otimes\varphi^{13}\otimes(\shift^{-1}\varphi^2)\otimes\varphi^3\otimes(\shift^{-1}\varphi^4)\otimes\varphi^{51}(\shift^{-1}\varphi^{52})\varphi^{53} &=0;
\end{align*}

from \[\varphi^1\varphi^{21}\otimes(\shift^{-1}\varphi^{22})\otimes\varphi^{23}\otimes(\shift^{-1}\varphi^{24})\otimes\varphi^{25}\varphi^3\otimes(\shift^{-1}\varphi^4)\otimes\varphi^5 =0.\]

\begin{align*}
\varphi^{11}(\shift^{-1}\varphi^{12})\varphi^{13}\varphi^{21}\otimes(\shift^{-1}\varphi^{22})\otimes\varphi^{23}\otimes(\shift^{-1}\varphi^{24})\otimes\varphi^{25}\varphi^3\otimes(\shift^{-1}\varphi^4)\otimes\varphi^5 &=0,\\
\varphi^1\varphi^{211}(\shift^{-1}\varphi^{212})\varphi^{213}\otimes(\shift^{-1}\varphi^{22})\otimes\varphi^{23}\otimes(\shift^{-1}\varphi^{24})\otimes\varphi^{25}\varphi^3\otimes(\shift^{-1}\varphi^4)\otimes\varphi^5 &=0,\\
\varphi^1\varphi^{21}\otimes\varphi^{221}(\shift^{-1}\varphi^{222})\varphi^{223}(\shift^{-1}\varphi^{224})\varphi^{225}\otimes\varphi^{23}\otimes(\shift^{-1}\varphi^{24})\otimes\varphi^{25}\varphi^3\otimes(\shift^{-1}\varphi^4)\otimes\varphi^5 &=0,\\
\varphi^1\varphi^{21}\otimes(\shift^{-1}\varphi^{22})\otimes\varphi^{231}(\shift^{-1}\varphi^{232})\varphi^{233}\otimes(\shift^{-1}\varphi^{24})\otimes\varphi^{25}\varphi^3\otimes(\shift^{-1}\varphi^4)\otimes\varphi^5 &=0,\\
\varphi^1\varphi^{21}\otimes(\shift^{-1}\varphi^{22})\otimes\varphi^{23}\otimes\varphi^{241}(\shift^{-1}\varphi^{242})\varphi^{243}(\shift^{-1}\varphi^{244})\varphi^{245}\otimes\varphi^{25}\varphi^3\otimes(\shift^{-1}\varphi^4)\otimes\varphi^5 &=0,\\
\varphi^1\varphi^{21}\otimes(\shift^{-1}\varphi^{22})\otimes\varphi^{23}\otimes(\shift^{-1}\varphi^{24})\otimes\varphi^{251}(\shift^{-1}\varphi^{252})\varphi^{253}\varphi^3\otimes(\shift^{-1}\varphi^4)\otimes\varphi^5 &=0,\\
\varphi^1\varphi^{21}\otimes(\shift^{-1}\varphi^{22})\otimes\varphi^{23}\otimes(\shift^{-1}\varphi^{24})\otimes\varphi^{25}\varphi^{31}(\shift^{-1}\varphi^{32})\varphi^{33}\otimes(\shift^{-1}\varphi^4)\otimes\varphi^5 &=0,\\
\varphi^1\varphi^{21}\otimes(\shift^{-1}\varphi^{22})\otimes\varphi^{23}\otimes(\shift^{-1}\varphi^{24})\otimes\varphi^{25}\varphi^3\otimes\varphi^{41}(\shift^{-1}\varphi^{42})\varphi^{43}(\shift^{-1}\varphi^{44})\varphi^{45}\otimes\varphi^5 &=0,\\
\varphi^1\varphi^{21}\otimes(\shift^{-1}\varphi^{22})\otimes\varphi^{23}\otimes(\shift^{-1}\varphi^{24})\otimes\varphi^{25}\varphi^3\otimes(\shift^{-1}\varphi^4)\otimes\varphi^{51}(\shift^{-1}\varphi^{52})\varphi^{53} &=0;
\end{align*}

from \[\varphi^1\otimes(\shift^{-1}\varphi^2)\otimes\varphi^{31}\otimes(\shift^{-1}\varphi^{32})\otimes\varphi^{33}\otimes(\shift^{-1}\varphi^4)\otimes\varphi^5  =0.\]

\begin{align*}
\varphi^{11}(\shift^{-1}\varphi^{12})\varphi^{13}\otimes(\shift^{-1}\varphi^2)\otimes\varphi^{31}\otimes(\shift^{-1}\varphi^{32})\otimes\varphi^{33}\otimes(\shift^{-1}\varphi^4)\otimes\varphi^5  &=0,\\
\varphi^1\otimes\varphi^{21}(\shift^{-1}\varphi^{22})\varphi^{23}(\shift^{-1}\varphi^{24})\varphi^{25}\otimes\varphi^{31}\otimes(\shift^{-1}\varphi^{32})\otimes\varphi^{33}\otimes(\shift^{-1}\varphi^4)\otimes\varphi^5  &=0,\\
\varphi^1\otimes(\shift^{-1}\varphi^2)\otimes\varphi^{311}(\shift^{-1}\varphi^{312})\varphi^{313}\otimes(\shift^{-1}\varphi^{32})\otimes\varphi^{33}\otimes(\shift^{-1}\varphi^4)\otimes\varphi^5  &=0,\\
\varphi^1\otimes(\shift^{-1}\varphi^2)\otimes\varphi^{31}\otimes\varphi^{321}(\shift^{-1}\varphi^{322})\varphi^{323}(\shift^{-1}\varphi^{324})\varphi^{325}\otimes\varphi^{33}\otimes(\shift^{-1}\varphi^4)\otimes\varphi^5  &=0,\\
\varphi^1\otimes(\shift^{-1}\varphi^2)\otimes\varphi^{31}\otimes(\shift^{-1}\varphi^{32})\otimes\varphi^{331}(\shift^{-1}\varphi^{332})\varphi^{333}\otimes(\shift^{-1}\varphi^4)\otimes\varphi^5  &=0,\\
\varphi^1\otimes(\shift^{-1}\varphi^2)\otimes\varphi^{31}\otimes(\shift^{-1}\varphi^{32})\otimes\varphi^{33}\otimes\varphi^{41}(\shift^{-1}\varphi^{42})\varphi^{43}(\shift^{-1}\varphi^{44})\varphi^{45}\otimes\varphi^5  &=0,\\
\varphi^1\otimes(\shift^{-1}\varphi^2)\otimes\varphi^{31}\otimes(\shift^{-1}\varphi^{32})\otimes\varphi^{33}\otimes(\shift^{-1}\varphi^4)\otimes\varphi^{51}(\shift^{-1}\varphi^{52})\varphi^{53}  &=0;
\end{align*}

from \[\varphi^1\otimes(\shift^{-1}\varphi^2)\otimes\varphi^3\varphi^{41}\otimes(\shift^{-1}\varphi^{42})\otimes\varphi^{43}\otimes(\shift^{-1}\varphi^{44})\otimes\varphi^{45}\varphi^5 =0.\]

\begin{align*}
\varphi^{11}(\shift^{-1}\varphi^{12})\varphi^{13}\otimes(\shift^{-1}\varphi^2)\otimes\varphi^3\varphi^{41}\otimes(\shift^{-1}\varphi^{42})\otimes\varphi^{43}\otimes(\shift^{-1}\varphi^{44})\otimes\varphi^{45}\varphi^5 &=0,\\
\varphi^1\otimes\varphi^{21}(\shift^{-1}\varphi^{22})\varphi^{23}(\shift^{-1}\varphi^{24})\varphi^{25}\otimes\varphi^3\varphi^{41}\otimes(\shift^{-1}\varphi^{42})\otimes\varphi^{43}\otimes(\shift^{-1}\varphi^{44})\otimes\varphi^{45}\varphi^5 &=0,\\
\varphi^1\otimes(\shift^{-1}\varphi^2)\otimes\varphi^{31}(\shift^{-1}\varphi^{32})\varphi^{33}\varphi^{41}\otimes(\shift^{-1}\varphi^{42})\otimes\varphi^{43}\otimes(\shift^{-1}\varphi^{44})\otimes\varphi^{45}\varphi^5 &=0,\\
\varphi^1\otimes(\shift^{-1}\varphi^2)\otimes\varphi^3\varphi^{411}(\shift^{-1}\varphi^{412})\varphi^{413}\otimes(\shift^{-1}\varphi^{42})\otimes\varphi^{43}\otimes(\shift^{-1}\varphi^{44})\otimes\varphi^{45}\varphi^5 &=0,\\
\varphi^1\otimes(\shift^{-1}\varphi^2)\otimes\varphi^3\varphi^{41}\otimes\varphi^{421}(\shift^{-1}\varphi^{422})\varphi^{423}(\shift^{-1}\varphi^{424})\varphi^{425}\otimes\varphi^{43}\otimes(\shift^{-1}\varphi^{44})\otimes\varphi^{45}\varphi^5 &=0,\\
\varphi^1\otimes(\shift^{-1}\varphi^2)\otimes\varphi^3\varphi^{41}\otimes(\shift^{-1}\varphi^{42})\otimes\varphi^{431}(\shift^{-1}\varphi^{432})\varphi^{433}\otimes(\shift^{-1}\varphi^{44})\otimes\varphi^{45}\varphi^5 &=0,\\
\varphi^1\otimes(\shift^{-1}\varphi^2)\otimes\varphi^3\varphi^{41}\otimes(\shift^{-1}\varphi^{42})\otimes\varphi^{43}\otimes\varphi^{441}(\shift^{-1}\varphi^{442})\varphi^{443}(\shift^{-1}\varphi^{444})\varphi^{445}\otimes\varphi^{45}\varphi^5 &=0,\\
\varphi^1\otimes(\shift^{-1}\varphi^2)\otimes\varphi^3\varphi^{41}\otimes(\shift^{-1}\varphi^{42})\otimes\varphi^{43}\otimes(\shift^{-1}\varphi^{44})\otimes\varphi^{451}(\shift^{-1}\varphi^{452})\varphi^{453} \varphi^5 &=0,\\
\varphi^1\otimes(\shift^{-1}\varphi^2)\otimes\varphi^3\varphi^{41}\otimes(\shift^{-1}\varphi^{42})\otimes\varphi^{43}\otimes(\shift^{-1}\varphi^{44})\otimes\varphi^{45}\varphi^{51}(\shift^{-1}\varphi^{52})\varphi^{53}  &=0;
\end{align*}

and from \[\varphi^1\otimes(\shift^{-1}\varphi^2)\otimes\varphi^3\otimes(\shift^{-1}\varphi^4)\otimes\varphi^{51}\otimes(\shift^{-1}\varphi^{52})\otimes\varphi^{53} =0.\]

\begin{align*}
\varphi^{11}(\shift^{-1}\varphi^{12})\varphi^{13}\otimes(\shift^{-1}\varphi^2)\otimes\varphi^3\otimes(\shift^{-1}\varphi^4)\otimes\varphi^{51}\otimes(\shift^{-1}\varphi^{52})\otimes\varphi^{53} &=0,\\
\varphi^1\otimes\varphi^{21}(\shift^{-1}\varphi^{22})\varphi^{23}(\shift^{-1}\varphi^{24})\varphi^{25}\otimes\varphi^3\otimes(\shift^{-1}\varphi^4)\otimes\varphi^{51}\otimes(\shift^{-1}\varphi^{52})\otimes\varphi^{53} &=0,\\
\varphi^1\otimes(\shift^{-1}\varphi^2)\otimes\varphi^{31}(\shift^{-1}\varphi^{32})\varphi^{33}\otimes(\shift^{-1}\varphi^4)\otimes\varphi^{51}\otimes(\shift^{-1}\varphi^{52})\otimes\varphi^{53} &=0,\\
\varphi^1\otimes(\shift^{-1}\varphi^2)\otimes\varphi^3\otimes\varphi^{41}(\shift^{-1}\varphi^{42})\varphi^{43}(\shift^{-1}\varphi^{44})\varphi^{45}\otimes\varphi^{51}\otimes(\shift^{-1}\varphi^{52})\otimes\varphi^{53} &=0,\\
\varphi^1\otimes(\shift^{-1}\varphi^2)\otimes\varphi^3\otimes(\shift^{-1}\varphi^4)\otimes\varphi^{511}(\shift^{-1}\varphi^{512})\varphi^{513} \otimes(\shift^{-1}\varphi^{52})\otimes\varphi^{53} &=0,\\
\varphi^1\otimes(\shift^{-1}\varphi^2)\otimes\varphi^3\otimes(\shift^{-1}\varphi^4)\otimes\varphi^{51}\otimes\varphi^{521}(\shift^{-1}\varphi^{522})\varphi^{523}(\shift^{-1}\varphi^{524})\varphi^{525}\otimes\varphi^{53} &=0,\\
\varphi^1\otimes(\shift^{-1}\varphi^2)\otimes\varphi^3\otimes(\shift^{-1}\varphi^4)\otimes\varphi^{51}\otimes(\shift^{-1}\varphi^{52})\otimes\varphi^{531}(\shift^{-1}\varphi^{532})\varphi^{533} &=0.
\end{align*}

We again use the isomorphism $\overline{V}\cong A\otimes \Phi\otimes A$ to enable us to insert tensor in the appropriate places in the statement of the lemma.
\end{proof}

Whenever we use this lemma from hereon to identify a vanishing term, we will label this term with a $g$. 

\begin{lem}
Let $\varphi\in \Phi$, and $a,a_1,a_2\in J$. Then, using Sweedler notation as above, we obtain
\begin{enumerate}[(a)]
\item\label{mu4:0}
\begin{align*}
\mu_4^0(\varphi)=\,&\varphi^1\varphi^2\otimes \shift \varphi^3\otimes \varphi^4\otimes \shift \varphi^5\otimes 1-1\otimes \shift \varphi^1\otimes \varphi^2\varphi^3\varphi^4\otimes \shift \varphi^5\otimes 1+1\otimes \shift \varphi^1\otimes \varphi^2\otimes \shift \varphi^3\otimes \varphi^4\varphi^5\\
&+\varphi^1\varphi^{21}\dots \varphi^{24}\otimes \shift \varphi^{25}\otimes \shift \varphi^3\otimes \varphi^{4}\varphi^5+1\otimes \shift \varphi^1\otimes \shift \varphi^{21}\otimes \varphi^{22}\dots \varphi^{25}\varphi^3\varphi^4\varphi^5\\
&+\varphi^1\varphi^2\otimes \shift \varphi^3\otimes \shift \varphi^{41}\otimes \varphi^{42}\dots\varphi^{45}\varphi^{5}+\varphi^1\varphi^2\varphi^3\varphi^{41}\dots\varphi^{44}\otimes \shift \varphi^{45}\otimes \shift \varphi^5\otimes 1,
\end{align*}
\item\label{mu4:1}
\begin{align*}
\mu_4^1(1\otimes \shift a\otimes 1)=&-1\otimes \shift a^1\otimes \shift a^{21}\otimes \shift a^{221}\otimes a^{222}\dots a^{225} a^{23}\dots a^{25}a^3\\
&+a^1a^{21}\dots a^{23}a^{241}\dots a^{244}\otimes \shift a^{245}\otimes \shift a^{25}\otimes \shift a^3\otimes 1\\
&+a^1 a^{21}a^{22}\otimes \shift a^{23}\otimes a^{24}\otimes \shift a^{25}\otimes \shift a^3\otimes 1\\
&-1\otimes \shift a^1\otimes \shift a^{21}\otimes a^{22}\otimes \shift a^{23}\otimes a^{24}a^{25}a^3\\
&+1\otimes \shift a^1\otimes \shift a^{21}\otimes a^{22}\dots a^{25}\otimes \shift a^3\otimes 1\\
&-1\otimes \shift (a^1 a^{21})\otimes a^{22}\dots a^{24}\otimes \shift a^{25}\otimes \shift a^3\otimes 1\\
&+1\otimes \shift a^1\otimes \shift a^{21}\otimes a^{22}\dots a^{24}\otimes \shift a^{25}\otimes a^3,
\end{align*}
\item\label{mu4:2}
\begin{align*}
\mu_4^2(1\otimes \shift a_1\otimes \shift a_2\otimes 1)=&-1\otimes \shift a_1\otimes \shift (a_2^1 a_2^{21})\otimes a_2^{22}\dots a_2^{24}\otimes \shift a_2^{25}\otimes \shift a_2^3\otimes 1\\
&-1\otimes \shift a_1^1\otimes \shift a_1^{21}\otimes a_1^{22}\dots a_1^{25}\otimes \shift a_1^3\otimes \shift a_2\otimes 1\\
&-a_1\otimes \shift a_2^1\otimes \shift a_2^{21}\otimes a_2^{22}\dots a_2^{24}\otimes \shift a_2^{25}\otimes \shift a_2^3\otimes 1\\
&+1\otimes \shift a_1\otimes \shift (a_2^1 a_2^{21})\otimes a_2^{22}\dots a_2^{24}\otimes \shift a_2^{25}\otimes \shift a_2^3\otimes 1\\
&+1\otimes \shift a_1\otimes \shift a_2^1\otimes a_2^{21}\dots a_2^{24}\otimes \shift a_2^{25}\otimes \shift a_2^3\otimes 1\\
&+1\otimes \shift (a_1 a_2^1)\otimes \shift a_2^{21} \otimes a_2^{22} \dots a_2^{24}\otimes \shift a_2^{25} \otimes \shift a_2^3\otimes 1\\
&+1\otimes \shift (a_1^1 a_1^{21})\otimes a_1^{22}\dots a_1^{24}\otimes \shift a_1^{25}\otimes \shift a_1^3\otimes \shift a_2\otimes 1\\
&-a_1a_1^{21}a_1^{22}\otimes \shift a_1^{23}\otimes a_1^{24}\otimes \shift a_1^{25}\otimes \shift a_1^3\otimes \shift a_2\otimes 1\\
&-1\otimes \shift a_1\otimes \shift a_2^1\otimes \shift a_2^{21}\otimes a_2^{22}\dots a_2^{24}\otimes \shift (a_2^{25} a_2^3)\otimes 1\\
&+ 1\otimes \shift  a_1\otimes \shift a_2^1\otimes \shift a_2^{21}\otimes a_2^{22}\otimes \shift a_2^{23}\otimes a_2^{24} a_2^{25} a_2^3\\
&-1\otimes \shift a_1\otimes \shift a_2^1\otimes \shift a_2^{21}\otimes \shift a_2^{221}\otimes a_2^{222}\dots a_2^{225}a_2^{23}\dots a_2^{25}a_2^3\\
&+a_1^1a_1^{21}\dots a_1^{23}a_1^{241}\dots a_1^{244}\otimes \shift a_1^{245}\otimes \shift a_1^{25}\otimes \shift a_1^3\otimes \shift a_2\otimes 1.
\end{align*}
\end{enumerate}
\end{lem}

\begin{proof}
To compute $\mu_4$, we use the $A_\infty$-equation
\[\mu_4\mu_1-(\mu_1\otimes 1^{\otimes 3}+1\otimes \mu_1\otimes 1^{\otimes 2}+1^{\otimes 2}\otimes \mu_1\otimes 1+1^{\otimes 3}\otimes \mu_1)\mu_4-(\mu_3\otimes 1+1\otimes \mu_3)\mu_2+(\mu_2\otimes 1^{\otimes 2}-1\otimes \mu_2\otimes 1+1^{\otimes 2}\otimes \mu_2)\mu_3=0\]
and verify that it is satisfied by the stated formulae.

\eqref{mu4:0}: Note that
\begin{align*}
-(\mu_3\otimes 1)\mu_2(\varphi)=\,&-(\mu_3\otimes 1)(\hat{\omega}\varphi+\varphi\hat{\omega}-\varphi^1\varphi^2\varphi^3\varphi^4\varphi^5)\\
=\,&0\-\underbrace{1\otimes \shift \varphi^1\otimes \varphi^2\dots \varphi^5\hat{\omega}}_{25}+\underbrace{\varphi^1\varphi^2\otimes \shift \varphi^3\otimes \varphi^4\varphi^5\hat{\omega}}_{21}-\underbrace{\varphi^1\dots\varphi^4\otimes \shift \varphi^5\otimes \hat{\omega}}_{19}\\
&\+\underbrace{\varphi^1\otimes \shift \varphi^{21}\otimes \varphi^{22}\dots \varphi^{25}\varphi^3\varphi^4\varphi^5}_6-\underbrace{\varphi^1\varphi^{21}\varphi^{22}\otimes \shift \varphi ^{23}\otimes \varphi^{24}\varphi^{25}\varphi^3\varphi^4\varphi^5}_y\\
&+\underbrace{\varphi^1\varphi^{21}\dots\varphi^{24}\otimes \shift \varphi^{25}\otimes \varphi^3\varphi^4\varphi^5}_{15}
\end{align*}
and
\begin{align*}
-(1\otimes \mu_3)\mu_2(\varphi)=\,&-(1\otimes \mu_3)(\hat{\omega}\varphi+\varphi\hat{\omega}-\varphi^1\varphi^2\varphi^3\varphi^4\varphi^5)\\
=\,&-\underbrace{\hat{\omega}\otimes \shift \varphi^1\otimes \varphi^2\dots\varphi^5}_1+\underbrace{\hat{\omega}\varphi^1\varphi^2\otimes \shift \varphi^3\otimes \varphi^4\varphi^5}_2-\underbrace{\hat{\omega}\varphi^1\dots\varphi^4\otimes \shift \varphi^5\otimes 1}_{4}\+0\\
&\+\underbrace{\varphi^1\varphi^2\varphi^3\otimes \shift \varphi^{41}\otimes \varphi^{42}\dots \varphi^{45}\varphi^5}_{10}-\underbrace{\varphi^1\varphi^2\varphi^3\varphi^{41}\varphi^{42}\otimes \shift \varphi^{43}\otimes \varphi^{44}\varphi^{45}\varphi^5}_y\\
&+\underbrace{\varphi^1\varphi^2\varphi^3\varphi^{41}\dots \varphi^{44}\otimes \shift \varphi^{45}\otimes \varphi^5}_{23}.
\end{align*}
Furthermore,
\begin{align*}
(\mu_2\otimes 1^{\otimes 2})\mu_3(\varphi)=\,&(\mu_2\otimes 1^{\otimes 2})(1\otimes \shift \varphi^1\otimes \varphi^2\dots \varphi^5-\varphi^1\varphi^2\otimes \shift \varphi^3\otimes \varphi^4\varphi^5+\varphi^1\dots \varphi^4\otimes \shift \varphi^5\otimes 1)\\
=\,&\underbrace{\hat{\omega}\otimes \shift \varphi^1\otimes \varphi^2\dots \varphi^5}_1+\underbrace{1\otimes \shift \varphi^1\otimes \hat{\omega}\varphi^2\dots \varphi^5}_{26}\\
&-\underbrace{\varphi^{11}\varphi^{12}\otimes \shift \varphi^{13}\otimes \varphi^{2}\dots \varphi^5}_y+\underbrace{1\otimes \shift \varphi^{11}\otimes \varphi^{12}\varphi^{13}\varphi^2\dots \varphi^5}_y\\
&\-\underbrace{\varphi^1\hat{\omega}\varphi^2\otimes \shift \varphi^3\otimes \varphi^4\varphi^5}_3-\underbrace{\varphi^1\varphi^2\hat{\omega}\otimes \shift \varphi^3\otimes \varphi^4\varphi^5}_{16}-\underbrace{\varphi^1(\shift \otimes \shift)\partial(\shift^{-1}\varphi^2)\otimes \shift \varphi^3\otimes \varphi^4\varphi^5}_{14}\\
&\+\underbrace{\varphi^1\hat{\omega}\varphi^2\varphi^3\varphi^4\otimes \shift \varphi^5\otimes 1}_5+\underbrace{\varphi^1\varphi^2\hat{\omega}\varphi^3\varphi^4\otimes \shift \varphi^5\otimes 1}_8+\underbrace{\varphi^1(\shift \otimes \shift)\partial(\shift^{-1}\varphi^2)\varphi^3\varphi^4\otimes \shift \varphi^5\otimes 1}_y
\end{align*}
and
\begin{align*}
-(1\otimes \mu_2\otimes 1)\mu_3(\varphi)=\,&-\underbrace{1\otimes \shift \varphi^1\otimes \hat{\omega}\varphi^2\dots\varphi^5}_{26}-\underbrace{1\otimes \shift \varphi^1\otimes \varphi^2\hat{\omega}\varphi^3\varphi^4\varphi^5}_{12}-\underbrace{1\otimes \shift \varphi^1\otimes (\shift \otimes \shift)\partial(\shift^{-1}\varphi^2)\varphi^3\varphi^4\varphi^5}_7\\
&\+\underbrace{\varphi^1\varphi^2\otimes \shift \varphi^3\otimes \hat{\omega}\varphi^4\varphi^5}_{17}+\underbrace{\varphi^1\varphi^2\hat{\omega}\otimes \shift \varphi^3\otimes \varphi^4\varphi^5}_{16}\\
&-\underbrace{\varphi^1\varphi^2\varphi^{31}\varphi^{32}\otimes \shift \varphi^{33}\otimes \varphi^{4}\varphi^5}_y+\underbrace{\varphi^1\varphi^2\otimes \shift \varphi^{31}\otimes \varphi^{32}\varphi^{33}\varphi^4\varphi^5}_y\\
&\-\underbrace{\varphi^1\dots \varphi^4\hat{\omega}\otimes \shift \varphi^5\otimes 1}_{18}-\underbrace{\varphi^1\varphi^2\varphi^3\hat{\omega}\varphi^4\otimes \shift \varphi^5\otimes 1}_9-\underbrace{\varphi^1\varphi^2\varphi^3(\shift \otimes \shift)\partial(\shift^{-1}\varphi^4)\otimes \shift \varphi^5\otimes 1}_{22}
\end{align*}
and
\begin{align*}
(1^{\otimes 2}\otimes \mu_2)\mu_3(\varphi)=\,&\underbrace{1\otimes \shift \varphi^1\otimes \varphi^2\varphi^3\hat{\omega}\varphi^4\varphi^5}_{13}+\underbrace{1\otimes \shift \varphi^1\otimes \varphi^2\varphi^3\varphi^4\hat{\omega}\varphi^5}_{24}+\underbrace{1\otimes \shift \varphi^1\otimes \varphi^2\varphi^3(\shift\otimes\shift)\partial(\shift^{-1}\varphi^4)\varphi^5}_y\\
&\-\underbrace{\varphi^1\varphi^2\otimes \shift \varphi^3\otimes \hat{\omega}\varphi^4\varphi^5}_{17}-\underbrace{\varphi^1\varphi^2\otimes \shift \varphi^3\otimes \varphi^4\hat{\omega}\varphi^5}_{20}-\underbrace{\varphi^1\varphi^2\otimes \shift \varphi^3\otimes (\shift\otimes\shift)\partial(\shift^{-1}\varphi^4)\varphi^5}_{11}\\
&\+\underbrace{\varphi^1\dots \varphi^4\hat{\omega}\otimes \shift \varphi^5\otimes 1}_{18}+\underbrace{\varphi^1\dots \varphi^4\otimes \shift \varphi^5\otimes \hat{\omega}}_{19}\\
&-\underbrace{\varphi^1\dots \varphi^4\varphi^{51}\varphi^{52}\otimes \shift \varphi^{53}\otimes 1}_y+\underbrace{\varphi^1\dots \varphi^4\otimes \shift \varphi^{51}\otimes \varphi^{52}\varphi^{53}}_y.
\end{align*}
Finally, we write $\mu_1^{\otimes}=\mu_1\otimes 1^{\otimes 3}+1\otimes \mu_1\otimes 1^{\otimes 2}+1^{\otimes 2}\otimes \mu_1\otimes 1+1^{\otimes 3}\otimes \mu_1$. Applying this to each individual summand of $\mu_4^0(\varphi)$ we obtain

\begin{align*}
-\mu_1^{\otimes}(\varphi^1\varphi^2\otimes \shift \varphi^3\otimes \varphi^4\otimes \shift \varphi^5\otimes 1)=\,0&\boldsymbol{\uwave{-}}\underbrace{\varphi^1\varphi^2\hat{\omega}\varphi^3\varphi^4\otimes \shift \varphi^5\otimes 1}_8+\underbrace{\varphi^1\varphi^2\varphi^3\hat{\omega}\varphi^4\otimes \shift \varphi^5\otimes 1}_9\\
&+\underbrace{\varphi^1\varphi^2\shift\partial(\varphi^3)\varphi^4\otimes \shift \varphi^5\otimes 1}_y\boldsymbol{\uwave{+}} 0 \boldsymbol{\uwave{+}}
\underbrace{\varphi^1\varphi^2\otimes \shift \varphi^3\otimes \varphi^4\hat{\omega}\varphi^5}_{20}\\
&-\underbrace{\varphi^1\varphi^2\otimes \shift \varphi^3\otimes \varphi^4\varphi^5\hat{\omega}}_{21}-\underbrace{\varphi^1\varphi^2\otimes \shift \varphi^3\otimes \varphi^4\shift \partial(\varphi^5)}_y
\end{align*}
and
\begin{align*}
-\mu_1^{\otimes}(-1\otimes \shift \varphi^1\otimes \varphi^2\varphi^3\varphi^4\otimes \shift \varphi^5\otimes 1)=\,&
\underbrace{\hat{\omega}\varphi^1\dots \varphi^4\otimes \shift \varphi^5\otimes 1}_4-\underbrace{\varphi^1\hat{\omega}\varphi^2\varphi^3\varphi^4\otimes \shift \varphi^5\otimes 1}_5\\ &-\underbrace{\shift \partial(\varphi^1)\varphi^2\varphi^3\varphi^4\otimes \shift \varphi^5\otimes 1}_y\boldsymbol{\uwave{+}} 0\boldsymbol{\uwave{+}} 0 
\boldsymbol{\uwave{-}}\underbrace{1\otimes \shift \varphi^1\otimes \varphi^2\varphi^3\varphi^4\hat{\omega}\varphi^5}_{24}\\
&+\underbrace{1\otimes \shift \varphi^1\otimes \varphi^2\dots \varphi^5\hat{\omega}}_{25}+\underbrace{1\otimes \shift \varphi^1\otimes \varphi^2\varphi^3\varphi^4\shift\partial(\varphi^5)}_y
\end{align*}
and 
\begin{align*}
-\mu_1^{\otimes}(1\otimes \shift \varphi^1\otimes \varphi^2\otimes \shift \varphi^3\otimes \varphi^4\varphi^5))=\,&-\underbrace{\hat{\omega}\varphi^1\varphi^2\otimes \shift \varphi^3\otimes \varphi^4\varphi^5}_2+\underbrace{\varphi^1\hat{\omega}\varphi^2\otimes \shift \varphi^3\otimes \varphi^4\varphi^5}_3 \\&
+\underbrace{\shift\partial(\varphi^1)\varphi^2\otimes \shift \varphi^3\otimes \varphi^4\varphi^5}_y\boldsymbol{\uwave{+}} 0 \boldsymbol{\uwave{+}}\underbrace{1\otimes \shift \varphi^1\otimes \varphi^2\hat{\omega}\varphi^3\varphi^4\varphi^5}_{12}\\ &-\underbrace{1\otimes \shift \varphi^1\otimes \varphi^2\varphi^3\hat{\omega}\varphi^4\varphi^5}_{13}-\underbrace{1\otimes \shift \varphi^1\otimes \varphi^2\partial(\varphi^3)\varphi^4\varphi^5}_y\boldsymbol{\uwave{+}} 0
\end{align*}
and
\begin{align*}
-\mu_1^{\otimes}(\varphi^1\varphi^{21}\dots \varphi^{24}\otimes \shift \varphi^{25}\otimes \shift \varphi^3\otimes \varphi^{4}\varphi^5) =&\, 0\boldsymbol{\uwave{+}} 0\boldsymbol{\uwave{+}}\underbrace{\varphi^1(\shift\otimes\shift)\partial(\shift^{-1}\varphi^2)\otimes \shift \varphi^3\otimes \varphi^4\varphi^5}_{14}\\
&+\underbrace{\varphi^1\varphi^{21}\dots \varphi^{24}\otimes \shift (\varphi^{25}\varphi^3)\otimes \varphi^4\varphi^5}_y\\
&-\underbrace{\varphi^1\varphi^{21}\dots \varphi^{24}\otimes \shift \varphi^{25}\otimes \varphi^3\varphi^4\varphi^5}_{15}\boldsymbol{\uwave{+}} 0
\end{align*}
and
\begin{align*}
-\mu_1^{\otimes}(1\otimes \shift \varphi^1\otimes \shift \varphi^{21}\otimes \varphi^{22}\dots \varphi^{25}\varphi^3\varphi^4\varphi^5)=
&\, -\underbrace{\varphi^1\otimes \shift \varphi^{21}\otimes \varphi^{22}\dots \varphi^{25}\varphi^3\varphi^4\varphi^5}_6\\&+\underbrace{1\otimes \shift (\varphi^1\varphi^{21})\otimes \varphi^{22}\dots \varphi^{25}\varphi^3\varphi^4\varphi^5}_y\\
&+\underbrace{1\otimes \shift \varphi^1\otimes (\shift \otimes \shift)\partial(\shift^{-1}\varphi^2)\varphi^3\varphi^4\varphi^5}_7\boldsymbol{\uwave{+}}0\boldsymbol{\uwave{+}}0\boldsymbol{\uwave{+}}0
\end{align*}
and
\begin{align*}
-\mu_1^{\otimes}(\varphi^1\varphi^2\otimes \shift \varphi^3\otimes \shift \varphi^{41}\otimes \varphi^{42}\dots\varphi^{45}\varphi^{5})=&\,0\boldsymbol{\uwave{-}}\underbrace{\varphi^1\varphi^2\varphi^3\otimes \shift \varphi^{41}\otimes \varphi^{42}\dots \varphi^{45}\varphi^5}_{10}\\&+\underbrace{\varphi^1\varphi^2\otimes \shift (\varphi^3\varphi^{41})\otimes \varphi^{42}\dots \varphi^{45}\varphi^5}_y\\
&+\underbrace{\varphi^1\varphi^2\otimes \shift \varphi^3\otimes (\shift\otimes\shift) \partial(\shift^{-1}\varphi^4)\varphi^5}_{11}\boldsymbol{\uwave{+}}0\boldsymbol{\uwave{+}}0
\end{align*}
and 
\begin{align*}
-\mu_1^{\otimes}(\varphi^1\varphi^2\varphi^3\varphi^{41}\dots\varphi^{44}\otimes \shift \varphi^{45}\otimes \shift \varphi^5\otimes 1)=&\,0\boldsymbol{\uwave{+}}0\boldsymbol{\uwave{+}}0\boldsymbol{\uwave{+}}\underbrace{\varphi^1\varphi^2\varphi^3(\shift \otimes \shift)\partial(\shift^{-1}\varphi^4)\otimes \shift \varphi^5\otimes 1}_{22}\\&+\underbrace{\varphi^1\varphi^2\varphi^3\varphi^{41}\dots \varphi^{44}\otimes \shift (\varphi^{45}\varphi^5)\otimes 1}_y\\
&-\underbrace{\varphi^1\varphi^2\varphi^3\varphi^{41}\dots \varphi^{44}\otimes \shift \varphi^{45}\otimes \varphi^5}_{23}.
\end{align*}
This finishes the computation for $\mu_4^0$. 

\eqref{mu4:1}: For $\mu_4^1$ we compute:
\begin{align*}
\mu_4^0\mu_1^1(1\otimes \shift a\otimes 1)=\,&\mu_4^0(\hat{\omega}a-a\hat{\omega}-a^1a^2a^3)\\
=\,&0\+0\-\underbrace{a^1a^{21}a^{22}\otimes \shift a^{23}\otimes a^{24}\otimes \shift a^{25}\otimes a^3}_{18}+\underbrace{a^1\otimes \shift a^{21}\otimes a^{22}a^{23}a^{24}\otimes \shift a^{25}\otimes a^3}_{27}\\
&-\underbrace{a^1\otimes \shift a^{21}\otimes a^{22}\otimes \shift a^{23}\otimes a^{24}a^{25}a^3}_{20}-\underbrace{a^1a^{21}a^{221}\dots a^{224}\otimes \shift a^{225}\otimes \shift a^{23}\otimes a^{24}a^{25}a^3}_g\\
&-\underbrace{a^1\otimes \shift a^{21}\otimes \shift a^{221}\otimes a^{222}\dots a^{225}a^{23}a^{24}a^{25}a^3}_{26}\\
&-\underbrace{a^1a^{21}a^{22}\otimes \shift a^{23}\otimes \shift a^{241}\otimes a^{242}\dots a^{245}a^{25}a^3}_g\\
&-\underbrace{a^1a^{21}a^{22}a^{23}a^{241}\dots a^{244}\otimes \shift a^{245}\otimes \shift a^{25}\otimes a^3}_{33}
\end{align*}
and
\begin{align*}
-(\mu_3\otimes 1)\mu_2(1\otimes \shift a\otimes 1)=\,&-(\mu_3\otimes 1)(\hat{\omega}\otimes \shift a\otimes 1+1\otimes \shift a\otimes \hat{\omega}-a^1a^2\otimes \shift a^3\otimes 1+1\otimes \shift a^1\otimes a^2a^3)\\
=\,&0\-\underbrace{1\otimes \shift a^1\otimes a^2\otimes \shift a^3\otimes \hat{\omega}}_1-\underbrace{1\otimes \shift a^1\otimes \shift a^{21}\otimes a^{22}\dots a^{25}a^3\hat{\omega}}_{12}\\
&-\underbrace{a^1a^{21}\dots a^{24}\otimes \shift a^{25}\otimes \shift a^3\otimes \hat{\omega}}_3 \+\underbrace{a^1\otimes \shift a^{21}\otimes a^{22}\dots a^{25}\otimes \shift a^3\otimes 1}_{13}\\
&-\underbrace{a^1a^{21}a^{22}\otimes \shift a^{23}\otimes a^{24}a^{25}\otimes \shift a^3\otimes 1}_{19}+\underbrace{a^1a^{21}\dots a^{24}\otimes \shift a^{25}\otimes 1\otimes \shift a^3\otimes 1}_{10}\\
&\-\underbrace{1\otimes \shift a^{11}\otimes a^{12}\otimes \shift a^{13}\otimes a^2a^3}_{y}-\underbrace{1\otimes \shift a^{11}\otimes \shift a^{121}\otimes a^{122}\dots a^{125}a^{13}a^2a^3}_g\\
&-\underbrace{a^{11}a^{121}\dots a^{124}\otimes \shift a^{125}\otimes \shift a^{13}\otimes a^2a^3}_g
\end{align*}
and
\begin{align*}
-(1\otimes \mu_3)\mu_2(1\otimes \shift a\otimes 1)=\,&-\underbrace{\hat{\omega}\otimes \shift a^1\otimes a^2\otimes \shift a^3\otimes 1}_2-\underbrace{\hat{\omega}\otimes \shift a^1\otimes \shift a^{21}\otimes a^{22}\dots a^{25}a^3}_4\\
&-\underbrace{\hat{\omega}a^1a^{21}\dots a^{24}\otimes \shift a^{25}\otimes \shift a^3\otimes 1}_{14}\+0\+\underbrace{a^1a^2\otimes \shift a^{31}\otimes a^{32}\otimes \shift a^{33}\otimes 1}_y\\
&+\underbrace{a^1a^2\otimes \shift a^{31}\otimes \shift a^{321}\otimes a^{322}\dots a^{325}a^{33}}_g\\
&+\underbrace{a^1a^2a^{31}a^{321}\dots a^{324}\otimes \shift a^{325}\otimes \shift a^{33}\otimes 1}_g
\+\underbrace{1\otimes \shift a^1\otimes 1\otimes \shift a^{21}\otimes a^{22}\dots a^{25}a^3}_9\\
&-\underbrace{1\otimes \shift a^1\otimes a^{21}a^{22}\otimes \shift a^{23}\otimes a^{24}a^{25}a^3}_{21}
+\underbrace{1\otimes \shift a^1\otimes a^{21}\dots a^{24}\otimes \shift a^{25}\otimes a^3}_{15}
\end{align*}
and, since 
\[\mu_3^1(1\otimes \shift a\otimes 1)=1\otimes \shift a^1\otimes a^2\otimes \shift a^3\otimes 1+1\otimes \shift a^1\otimes \shift a^{21}\otimes a^{22}\dots a^{25}a^3+a^1a^{21}\dots a^{24}\otimes \shift a^{25}\otimes \shift a^3\otimes 1\]
we compute
\begin{align*}
&(\mu_2\otimes 1^{\otimes 2}-1\otimes \mu_2\otimes 1+1^{\otimes 2}\otimes \mu_2)(1\otimes \shift a^1\otimes a^2\otimes \shift a^3\otimes 1)\\
=\,&\underbrace{\hat{\omega}\otimes \shift a^1\otimes a^2\otimes \shift a^3\otimes 1}_2+\underbrace{1\otimes \shift a^1\otimes \hat{\omega} a^2\otimes \shift a^3\otimes 1}_5\\
&-\underbrace{a^{11}a^{12}\otimes \shift a^{13}\otimes a^2\otimes \shift a^3\otimes 1}_y+\underbrace{1\otimes \shift a^{11}\otimes a^{12}a^{13}a^2\otimes \shift a^3\otimes 1}_y\\
&\boldsymbol{\uwave{-}}\underbrace{1\otimes \shift a^1\otimes \hat{\omega}a^2\otimes \shift a^3\otimes 1}_5-\underbrace{1\otimes \shift a^1\otimes a^2\hat{\omega}\otimes \shift a^3\otimes 1}_6-\underbrace{1\otimes \shift a^1\otimes (\shift \otimes \shift)\partial(\shift^{-1}a^2)\otimes \shift a^3\otimes 1}_{11}\\
&\boldsymbol{\uwave{+}}\underbrace{1\otimes \shift a^1\otimes a^2\hat{\omega}\otimes \shift a^3\otimes 1}_{6}+\underbrace{1\otimes \shift a^1\otimes a^2\otimes \shift a^3\otimes \hat{\omega}}_1\\
&-\underbrace{1\otimes \shift a^1\otimes a^2a^{31}a^{32}\otimes \shift a^{33}\otimes 1}_y+\underbrace{1\otimes \shift a^1\otimes a^2\otimes \shift a^{31}\otimes a^{32}a^{33}}_y
\end{align*}
and
\begin{align*}
&(\mu_2\otimes 1^{\otimes 2}-1\otimes \mu_2\otimes 1+1^{\otimes 2}\otimes \mu_2)(1\otimes \shift a^1\otimes \shift a^{21}\otimes a^{22}\dots a^{25}a^3)\\
=\,&\underbrace{\hat{\omega}\otimes \shift a^1\otimes \shift a^{21}\otimes a^{22}\dots a^{25}a^3}_4
+\underbrace{1\otimes \shift a^1\otimes \shift a^{21}\otimes \hat{\omega}a^{22}\dots a^{25}a^3}_7
-\underbrace{1\otimes \shift a^1\otimes 1\otimes \shift a^{21}\otimes a^{22}\dots a^{25}a^3}_9\\
&-\underbrace{a^{11}a^{12}\otimes \shift a^{13}\otimes \shift a^{21}\otimes a^{22}\dots a^{25}a^3}_g
+\underbrace{1\otimes \shift a^1\otimes \shift a^{211}\otimes a^{212}a^{213}a^{22}\dots a^{25}a^3}_y\\
&\boldsymbol{\uwave{-}}\underbrace{1\otimes \shift a^1\otimes \shift a^{21}\otimes \hat{\omega}a^{22}\dots a^{25}a^3}_7
-\underbrace{1\otimes \shift a^1\otimes \shift a^{21}\otimes a^{22}\hat{\omega} a^{23}a^{24} a^{25}a^3}_{22}\\
&-\underbrace{1\otimes \shift a^1\otimes \shift a^{21}\otimes (\shift \otimes \shift)\partial(\shift^{-1}a^{22})a^{23}\dots a^{25}a^3}_{25}
\boldsymbol{\uwave{+}}\underbrace{1\otimes \shift a^1\otimes \shift a^{21}\otimes a^{22}a^{23}\hat{\omega}a^{24}a^{25}a^3}_{23}\\
&+\underbrace{1\otimes \shift a^1\otimes \shift a^{21}\otimes a^{22}\dots a^{24}\hat{\omega}a^{25}a^3}_{28}
+\underbrace{1\otimes \shift a^1\otimes \shift a^{21}\otimes a^{22}a^{23}(\shift \otimes \shift)\partial(\shift^{-1}a^{24})a^{25}a^3}_y
\end{align*}
and
\begin{align*}
&(\mu_2\otimes 1^{\otimes 2}-1\otimes \mu_2\otimes 1+1^{\otimes 2}\otimes \mu_2)(a^1a^{21}\dots a^{24}\otimes \shift a^{25}\otimes \shift a^3\otimes 1)\\
=\,&\underbrace{a^1a^{21}\hat{\omega}a^{22}\dots a^{24}\otimes \shift a^{25}\otimes \shift a^3\otimes 1}_{24}
+\underbrace{a^1a^{21}a^{22}\hat{\omega}a^{23}a^{24}\otimes \shift a^{25}\otimes \shift a^3\otimes 1}_{16}\\
&+\underbrace{a^1a^{21}(\shift \otimes \shift)\partial(\shift^{-1}a^{22})a^{23}a^{24}\otimes \shift a^{25}\otimes \shift a^3\otimes 1}_y
\boldsymbol{\uwave{-}}\underbrace{a^1a^{21}a^{22} a^{23}\hat{\omega}a^{24}\otimes \shift a^{25}\otimes \shift a^3\otimes 1}_{17}\\
&-\underbrace{a^1a^{21}\dots a^{24}\hat{\omega}\otimes \shift a^{25}\otimes \shift a^3\otimes 1}_8
-\underbrace{a^1a^{21}\dots a^{23}(\shift \otimes \shift)\partial(\shift^{-1}a^{24})\otimes \shift a^{25}\otimes \shift a^3\otimes 1}_{32}\\
&\boldsymbol{\uwave{+}}\underbrace{a^1a^{21}\dots a^{24}\hat{\omega}\otimes \shift a^{25}\otimes \shift a^3\otimes 1}_8
+\underbrace{a^1a^{21}\dots a^{24}\otimes \shift a^{25}\otimes \shift a^{3}\otimes \hat{\omega}}_3\\
&-\underbrace{a^1a^{21}\dots a^{24}\otimes \shift a^{25}\otimes 1\otimes \shift a^3\otimes 1}_{10}
-\underbrace{a^1a^{21}\dots a^{24}a^{251}a^{252}\otimes \shift a^{253}\otimes \shift a^3\otimes 1}_{y}\\
&+\underbrace{a^1a^{21}\dots a^{24}\otimes \shift a^{25}\otimes \shift a^{31}\otimes a^{32}a^{33}}_g
\end{align*}
To simplify notation, in the following we write $\mu_1^{\otimes}$ for $(\mu_1\otimes 1^{\otimes 3}+1\otimes \mu_1\otimes 1^{\otimes 2}+1^{\otimes 2}\otimes \mu_1\otimes 1+1^{\otimes 3}\otimes \mu_1)$. Applying this to each summand of $\mu_4^1(1\otimes \shift a \otimes 1)$ individually, we obtain:
\begin{align*}
&-\mu_1^{\otimes}(-1\otimes \shift a^1\otimes \shift a^{21}\otimes \shift a^{221}\otimes a^{222}\dots a^{225}a^{23}\dots a^{25}a^3)\\
=\,&\underbrace{a^1\otimes \shift a^{21}\otimes \shift a^{221}\otimes a^{222}\dots a^{225}a^{23}a^{24}a^{25}a^3}_{26}
-\underbrace{1\otimes \shift (a^1a^{21})\otimes \shift a^{221}\otimes a^{222}\dots a^{225}a^{23}\dots a^{25}a^3}_g\\
&+\underbrace{1\otimes \shift a^1\otimes \shift (a^{21}a^{221})\otimes a^{222}\dots a^{225}a^{23}\dots a^{25}a^3}_y
+\underbrace{1\otimes \shift a^1\otimes \shift a^{21}\otimes (\shift \otimes \shift)\partial(\shift^{-1}a^{22})a^{23}\dots a^{25}a^3}_{25}\boldsymbol{\uwave{+}}0\boldsymbol{\uwave{+}}0\boldsymbol{\uwave{+}}0
\end{align*}
and
\begin{align*}
&-\mu_1^{\otimes}(1\otimes \shift a^1\otimes \shift a^{21}\otimes a^{22}\dots a^{25}\otimes \shift a^3\otimes 1)\\
=\,&-\underbrace{a^1\otimes \shift a^{21}\otimes a^{22}\dots a^{25}\otimes \shift a^3\otimes 1}_{13}
+\underbrace{1\otimes \shift (a^1a^{21})\otimes a^{22}\dots a^{25}\otimes \shift a^3\otimes 1}_{31}\\
&+\underbrace{1\otimes \shift a^1\otimes (\shift \otimes \shift)\partial(\shift^{-1}a^2)\otimes \shift a^3\otimes 1}_{11}\boldsymbol{\uwave{+}}0\boldsymbol{\uwave{+}}0
\boldsymbol{\uwave{-}}\underbrace{1\otimes \shift a^1\otimes \shift a^{21}\otimes a^{22}\dots a^{25}\hat{\omega}a^3}_{29}\\
&+\underbrace{1\otimes \shift a^1\otimes \shift a^{21}\otimes a^{22}\dots a^{25}a^3\hat{\omega}}_{12}
+\underbrace{1\otimes \shift a^1\otimes \shift a^{21}\otimes a^{22}\dots a^{25}\shift\partial(a^3)}_g
\end{align*}
and
\begin{align*}
&-\mu_1^{\otimes}(-1\otimes \shift (a^1a^{21})\otimes a^{22}\dots a^{24}\otimes \shift a^{25}\otimes a^3\otimes 1)\\
=\,&\underbrace{\hat{\omega}a^1a^{21}\dots a^{24}\otimes \shift a^{25}\otimes \shift a^3\otimes 1}_{14}
-\underbrace{a^1a^{21}\hat{\omega}a^{22}\dots a^{24}\otimes \shift a^{25}\otimes \shift a^3\otimes 1}_{24}\\
&-\underbrace{\shift \partial(a^1)a^{21}\dots a^{24}\otimes \shift a^{25}\otimes \shift a^3\otimes 1}_g
-\underbrace{a^1\shift \partial(a^{21})a^{22}\dots a^{24}\otimes \shift a^{25}\otimes \shift a^3\otimes 1}_y\boldsymbol{\uwave{+}}0\boldsymbol{\uwave{+}}0\\
&\boldsymbol{\uwave{-}}\underbrace{1\otimes \shift (a^1a^{21})\otimes a^{22}\dots a^{25}\otimes \shift a^3\otimes 1}_{31}
+\underbrace{1\otimes \shift (a^1a^{21})\otimes a^{22}\dots a^{24}\otimes \shift (a^{25}a^3)\otimes 1}_y\\
&-\underbrace{1\otimes \shift (a^1a^{21})\otimes a^{22}\dots a^{24}\otimes \shift a^{25}\otimes a^3}_{30}
\end{align*}
and
\begin{align*}
&-\mu_1^{\otimes}(a^1a^{21}a^{22}\otimes \shift a^{23}\otimes a^{24}\otimes \shift a^{25}\otimes \shift a^3\otimes 1)\\
=\,0&\boldsymbol{\uwave{-}}\underbrace{a^1a^{21}a^{22}\hat{\omega}a^{23}a^{24}\otimes \shift a^{25}\otimes \shift a^3\otimes 1}_{16}
+\underbrace{a^1a^{21}a^{22}a^{23}\hat{\omega}a^{24}\otimes \shift a^{25}\otimes \shift a^3\otimes 1}_{17}\\
&+\underbrace{a^1a^{21}a^{22}\shift \partial(a^{23})a^{24}\otimes \shift a^{25}\otimes \shift a^3\otimes 1}_y\boldsymbol{\uwave{+}}0
\boldsymbol{\uwave{+}}\underbrace{a^1a^{21}a^{22}\otimes \shift a^{23}\otimes a^{24}a^{25}\otimes \shift a^3\otimes 1}_{19}\\
&-\underbrace{a^1a^{21}a^{22}\otimes \shift a^{23}\otimes a^{24}\otimes \shift (a^{25}a^3)\otimes 1}_y
+\underbrace{a^1a^{21}a^{22}\otimes \shift a^{23}\otimes a^{24}\otimes \shift a^{25}\otimes a^3}_{18}
\end{align*}
and
\begin{align*}
&-\mu_1^{\otimes}(-1\otimes \shift a^1\otimes \shift a^{21}\otimes a^{22}\otimes \shift a^{23}\otimes a^{24}a^{25}a^3)\\
=\,&\underbrace{a^1\otimes \shift a^{21}\otimes a^{22}\otimes \shift a^{23}\otimes a^{24}a^{25}a^3}_{20}
-\underbrace{1\otimes \shift (a^1a^{21})\otimes a^{22}\otimes \shift a^{23}\otimes a^{24}a^{25}a^3}_y\\
&+\underbrace{1\otimes \shift a^1\otimes a^{21}a^{22}\otimes \shift a^{23}\otimes a^{24}a^{25}a^3}_{21}\boldsymbol{\uwave{+}}0
\boldsymbol{\uwave{+}}\underbrace{1\otimes \shift a^1\otimes \shift a^{21}\otimes a^{22}\hat{\omega}a^{23}a^{24}a^{25}a^3}_{22}\\
&-\underbrace{1\otimes \shift a^1\otimes \shift a^{21}\otimes a^{22}a^{23}\hat{\omega}a^{24}a^{25}a^3}_{23}
-\underbrace{1\otimes \shift a^1\otimes a^{21}\otimes a^{22}\shift \partial(a^{23})a^{24}a^{25}a^3}_y\boldsymbol{\uwave{+}}0
\end{align*}
and
\begin{align*}
&-\mu_1^{\otimes}(1\otimes \shift a^1\otimes \shift a^{21}\otimes a^{22}\dots a^{24}\otimes \shift a^{25}\otimes a^3)\\
=\,&-\underbrace{a^1\otimes \shift a^{21}\otimes a^{22}\dots a^{24}\otimes \shift a^{25}\otimes a^3}_{27}
+\underbrace{1\otimes \shift (a^1a^{21})\otimes a^{22}a^{23} a^{24}\otimes \shift a^{25}\otimes a^3}_{30}\\
&-\underbrace{1\otimes \shift a^1\otimes a^{21}\dots a^{24}\otimes \shift a^{25}\otimes a^3}_{15}\boldsymbol{\uwave{+}}0\boldsymbol{\uwave{+}}0
\boldsymbol{\uwave{-}}\underbrace{1\otimes \shift a^1\otimes \shift a^{21}\otimes a^{22}\dots a^{24}\hat{\omega}a^{25}a^3}_{28}\\
&+\underbrace{1\otimes \shift a^1\otimes \shift a^{21}\otimes a^{22}\dots a^{25}\hat{\omega}a^3}_{29}
+\underbrace{1\otimes \shift a^1\otimes \shift a^{21}\otimes a^{22}\dots a^{24}\shift \partial(a^{25})a^3}_y
\end{align*}
and
\begin{align*}
&-\mu_1^{\otimes}(a^1a^{21}\dots a^{23}a^{241}\dots a^{244}\otimes \shift a^{245}\otimes \shift a^{25}\otimes \shift a^3\otimes 1)\\
=\,&0\boldsymbol{\uwave{+}}0\boldsymbol{\uwave{+}}0\boldsymbol{\uwave{+}}\underbrace{a^1a^{21}\dots a^{23}(\shift \otimes \shift)\partial(\shift^{-1}a^{24})\otimes \shift a^{25}\otimes \shift a^3\otimes 1}_{32}
+\underbrace{a^1a^{21}\dots a^{23}a^{241}\dots a^{244}\otimes \shift (a^{245}a^{25})\otimes \shift a^3\otimes 1}_y\\
&-\underbrace{a^1a^{21}\dots a^{23}a^{241}\dots a^{244}\otimes \shift a^{245}\otimes \shift (a^{25}a^3)\otimes 1}_g
+\underbrace{a^1a^{21}a^{22}a^{23}a^{241}\dots a^{244}\otimes \shift a^{245}\otimes \shift a^{25}\otimes a^3}_{33}
\end{align*}

\eqref{mu4:2}: For $\mu_4^2$ we compute:
\begin{align*}
&-(\mu_3\otimes 1)\mu_2^2(1\otimes \shift a_1\otimes \shift a_2\otimes 1)\\
=\,&-(\mu_3\otimes 1)(\hat{\omega}\otimes \shift a_1\otimes\shift  a_2\otimes 1+1\otimes \shift a_1\otimes \shift a_2\otimes \hat{\omega}-1\otimes \shift a_1\otimes 1\otimes \shift a_2\otimes 1\\
&-a_1^1a_1^2\otimes \shift a_1^3\otimes \shift a_2\otimes 1+1\otimes \shift a_1\otimes \shift a_2^1\otimes a_2^2a_2^3)\\
=\,&0\-\underbrace{1\otimes \shift a_1^1\otimes a_1^2\otimes \shift a_1^3\otimes \shift a_2\otimes \hat{\omega}}_5\\
&+\underbrace{1\otimes \shift a_1\otimes \shift a_2^1\otimes a_2^2\otimes \shift a_2^3\otimes \hat{\omega}}_6
-\underbrace{a_1^1a_1^{21}\dots a_1^{24}\otimes \shift a_1^{25}\otimes \shift a_1^3\otimes \shift a_2\otimes \hat{\omega}}_7\\
&+\underbrace{1\otimes \shift a_1\otimes \shift a_2^1\otimes \shift a_2^{21}\otimes a_2^{22}\dots a_2^{25}a_2^3\hat{\omega}}_{29}
\+\underbrace{1\otimes \shift a_1^1\otimes a_1^2\otimes \shift a_1^3\otimes 1\otimes \shift a_2\otimes 1}_{10}\\
&+\underbrace{a_1^1a_1^{21}\dots a_1^{24}\otimes \shift a_1^{25}\otimes \shift a_1^3\otimes 1\otimes \shift a_2\otimes 1}_{11}
+\underbrace{1\otimes \shift a_1^1\otimes \shift a_1^{21}\otimes a_1^{22}\dots a_1^{25}a_1^3\otimes \shift a_2\otimes 1}_{39}\\
&\+\underbrace{a_1^1\otimes \shift a_1^{21}\otimes a_1^{22}\dots a_1^{25}\otimes \shift a_1^3\otimes \shift a_2\otimes 1}_{40}
-\underbrace{a_1^1a_1^{21}a_1^{22}\otimes \shift a_1^{23}\otimes a_1^{24}a_1^{25}\otimes \shift a_1^3\otimes \shift a_2\otimes 1}_{27}\\
&+\underbrace{a_1^1a_1^{21}\dots a_1^{24}\otimes \shift a_1^{25}\otimes 1\otimes \shift a_1^3\otimes \shift a_2\otimes 1}_{17}
\-\underbrace{1\otimes \shift a_1^1\otimes a_1^2\otimes \shift a_1^3\otimes \shift a_2^1\otimes a_2^2a_2^3}_{18}\\
&+\underbrace{1\otimes \shift a_1\otimes \shift a_2^{11}\otimes a_2^{12}\otimes \shift a_2^{13}\otimes a_1^2a_1^3}_y
-\underbrace{a_1^1a_1^{21}\dots a_1^{24}\otimes \shift a_1^{25}\otimes \shift a_1^3\otimes \shift a_2^1\otimes a_2^2a_2^3}_{20}\\
&+\underbrace{1\otimes \shift a_1\otimes \shift a_2^{11}\otimes \shift a_2^{121}\otimes a_2^{122}\dots a_2^{125}a_2^{13}a_2^2a_2^3}_g
\end{align*}
and
\begin{align*}
&-(1\otimes \mu_3)\mu_2^2(1\otimes \shift a_1\otimes \shift a_2\otimes 1)\\
=\,&-(1\otimes \mu_3)(\hat{\omega}\otimes \shift a_1\otimes\shift  a_2\otimes 1+1\otimes \shift a_1\otimes \shift a_2\otimes \hat{\omega}-1\otimes \shift a_1\otimes 1\otimes \shift a_2\otimes 1\\
&-a_1^1a_1^2\otimes \shift a_1^3\otimes \shift a_2\otimes 1+1\otimes \shift a_1\otimes \shift a_2^1\otimes a_2^2a_2^3)\\
=\,&-\underbrace{\hat{\omega}\otimes \shift a_1^1\otimes a_1^2\otimes \shift a_1^3\otimes \shift a_2\otimes 1}_{1}
+\underbrace{\hat{\omega}\otimes \shift a_1\otimes \shift a_2^1\otimes a_2^2\otimes \shift a_2^3\otimes 1}_{3}\\
&-\underbrace{\hat{\omega}a_1^1a_1^{21}\dots a_1^{24}\otimes \shift a_1^{25}\otimes \shift a_1^3\otimes \shift a_2\otimes 1}_{23}
+\underbrace{\hat{\omega}\otimes \shift a_1\otimes \shift a_2^1\otimes \shift a_2^{21}\otimes a_2^{22}\dots a_2^{25}a_2^3}_4\\
&\+0\-\underbrace{1\otimes \shift a_1\otimes 1\otimes \shift a_2^1\otimes a_2^2\otimes \shift a_2^3\otimes 1}_8
-\underbrace{1\otimes \shift a_1\otimes a_2^1a_2^{21}\dots a_2^{24}\otimes \shift a_2^{25}\otimes \shift a_2^3\otimes 1}_{45}\\
&-\underbrace{1\otimes \shift a_1\otimes 1\otimes \shift a_2^1\otimes \shift a_2^{21}\otimes a_2^{22}\dots a_2^{25}a_2^3}_9
\+\underbrace{a_1^1a_1^2\otimes \shift a_1^{31}\otimes a_1^{32}\otimes \shift a_1^{33}\otimes \shift a_2\otimes 1}_y\\
&-\underbrace{a_1^1a_1^2\otimes \shift a_1^3\otimes \shift a_2^1\otimes a_2^2\otimes \shift a_2^3\otimes 1}_{14}
+\underbrace{a_1^1a_1^2a_1^{31}a_1^{321}\dots a_1^{324}\otimes \shift a_1^{325}\otimes \shift a_1^{33}\otimes \shift a_2\otimes 1}_g\\
&-\underbrace{a_1^1a_1^2\otimes \shift a_1^3\otimes \shift a_2^1\otimes \shift a_2^{21}\otimes a_2^{22}\dots a_2^{25}a_2^3}_{19}
\-\underbrace{1\otimes \shift a_1\otimes \shift a_2^1\otimes 1\otimes \shift a_2^{21}\otimes a_2^{22}\dots a_2^{25}a_2^3}_{15}\\
&+\underbrace{1\otimes \shift a_1\otimes \shift a_2^1\otimes a_2^{21}a_2^{22}\otimes \shift a_2^{23}\otimes a_2^{24}a_2^{25}a_2^3}_{32}
-\underbrace{1\otimes \shift a_1\otimes \shift a_2^1\otimes a_2^{21}\dots a_2^{24}\otimes \shift a_2^{25}\otimes a_2^3}_{49}
\end{align*}
and
\begin{align*}
&(\mu_2\otimes 1^{\otimes 2})\mu_3(1\otimes \shift a_1\otimes \shift a_2\otimes 1)\\
=\,&(\mu_2\otimes 1^{\otimes 2})(1\otimes \shift a_1^1\otimes a_1^2\otimes \shift a_1^3\otimes \shift a_2\otimes 1-1\otimes \shift a_1\otimes \shift a_2^1\otimes a_2^2\otimes \shift a_2^3\otimes 1\\
&+a_1^1a_1^{21}\dots a_1^{24}\otimes \shift a_1^{25}\otimes \shift a_1^3\otimes \shift a_2\otimes 1-1\otimes \shift a_1\otimes \shift a_2^1\otimes \shift a_2^{21}\otimes a_2^{22}\dots a_2^{25}a_2^3)\\
=\,&\underbrace{\hat{\omega}\otimes \shift a_1^1\otimes a_1^2\otimes \shift a_1^3\otimes \shift a_2\otimes 1}_1
+\underbrace{1\otimes \shift a_1^1\otimes \hat{\omega}a_1^2\otimes \shift a_1^3\otimes \shift a_2\otimes 1}_{13}\\
&-\underbrace{a_1^{11}a_1^{12}\otimes \shift a_1^{13}\otimes a_1^2\otimes \shift a_1^3\otimes \shift a_2\otimes 1}_y
+\underbrace{1\otimes \shift a_1^{11}\otimes a_1^{12}a_1^{13}a_1^2\otimes \shift a_1^3\otimes \shift a_2\otimes 1}_y&\\
&\-\underbrace{\hat{\omega}\otimes \shift a_1\otimes \shift a_2^1\otimes a_2^2\otimes \shift a_2^3\otimes 1}_3
-\underbrace{1\otimes \shift a_1\otimes \shift a_2^1\otimes \hat{\omega}a_2^2\otimes \shift a_2^3\otimes 1}_{2}\\
&+\underbrace{1\otimes \shift a_1\otimes 1\otimes \shift a_2^1\otimes a_2^2\otimes \shift a_2^3\otimes 1}_8
+\underbrace{a_1^1a_1^2\otimes \shift a_1^3\otimes \shift a_2^1\otimes a_2^2\otimes \shift a_2^3\otimes 1}_{14}\\
&-\underbrace{1\otimes \shift a_1\otimes \shift a_2^{11}\otimes a_2^{12}a_2^{13}a_2^2\otimes \shift a_2^3\otimes 1}_{y}
\+\underbrace{a_1^1a_1^{21}\hat{\omega}a_1^{22}\dots a_1^{24}\otimes \shift a_1^{25}\otimes \shift a_1^3\otimes \shift a_2\otimes 1}_{24}\\
&+\underbrace{a_1^1a_1^{21}a_1^{22}\hat{\omega}a_1^{23}a_1^{24}\otimes \shift a_1^{25}\otimes \shift a_1^3\otimes \shift a_2\otimes 1}_{25}
+\underbrace{a_1^1a_1^{21}(\shift \otimes \shift)\partial(\shift^{-1}a_1^{22})a_1^{23}a_1^{24}\otimes \shift a_1^{25}\otimes \shift a_1^3\otimes \shift a_2\otimes 1}_y\\
&\-\underbrace{\hat{\omega}\otimes \shift a_1\otimes \shift a_2^1\otimes \shift a_2^{21}\otimes a_2^{22}\dots a_2^{25}a_2^3}_4
-\underbrace{1\otimes \shift a_1\otimes \shift a_2^1\otimes \shift a_2^{21}\otimes \hat{\omega}a_2^{22}\dots a_2^{25}a_2^3}_{21}\\
&+\underbrace{1\otimes \shift a_1\otimes 1\otimes \shift a_2^1\otimes \shift a_2^{21}\otimes a_2^{22}\dots a_2^{25}a_2^3}_9
+\underbrace{1\otimes \shift a_1\otimes \shift a_2^1\otimes 1\otimes \shift a_2^{21}\otimes a_2^{22}\dots a_2^{25}a_2^3}_{15}\\
&+\underbrace{a_1^1a_1^2\otimes \shift a_1^3\otimes \shift a_2^1\otimes \shift a_2^{21}\otimes a_2^{22}\dots a_2^{25}a_2^3}_{19}
-\underbrace{1\otimes \shift a_1\otimes \shift a_2^1\otimes \shift a_2^{211}\otimes a_2^{212}a_2^{213}a_2^{22}\dots a_2^{25}a_2^3}_{y} 
\end{align*}
and
\begin{align*}
&-(1\otimes \mu_2\otimes 1)\mu_3(1\otimes \shift a_1\otimes \shift a_2\otimes 1)\\
=\,&-(1\otimes \mu_2\otimes 1)(1\otimes \shift a_1^1\otimes a_1^2\otimes \shift a_1^3\otimes \shift a_2\otimes 1-1\otimes \shift a_1\otimes \shift a_2^1\otimes a_2^2\otimes \shift a_2^3\otimes 1\\
&+a_1^1a_1^{21}\dots a_1^{24}\otimes \shift a_1^{25}\otimes \shift a_1^3\otimes \shift a_2\otimes 1-1\otimes \shift a_1\otimes \shift a_2^1\otimes \shift a_2^{21}\otimes a_2^{22}\dots a_2^{25}a_2^3)\\
=\,&-\underbrace{1\otimes \shift a_1^1\otimes \hat{\omega}a_1^2\otimes \shift a_1^3\otimes \shift a_2\otimes 1}_{13}
-\underbrace{1\otimes \shift a_1^1\otimes a_1^2\hat{\omega}\otimes \shift a_1^3\otimes \shift a_2\otimes 1}_{12}\\
&-\underbrace{1\otimes \shift a_1^1\otimes (\shift \otimes \shift)\partial(\shift^{-1}a_1^2)\otimes \shift a_1^3\otimes \shift a_2\otimes 1}_{36}
\+\underbrace{1\otimes \shift a_1\otimes \shift a_2^1\otimes \hat{\omega}a_2^2\otimes \shift a_2^3\otimes 1}_{2}\\
&+\underbrace{1\otimes \shift a_1\otimes \shift a_2^1\otimes a_2^2\hat{\omega}\otimes \shift a_2^3\otimes 1}_{16}
+\underbrace{1\otimes \shift a_1\otimes \shift a_2^1\otimes (\shift \otimes \shift)\partial(\shift^{-1}a_2^2)\otimes \shift a_2^3\otimes 1}_{50}\\
&\-\underbrace{a_1^1a_1^{21}\dots a_1^{23}\hat{\omega}a_1^{24}\otimes \shift a_1^{25}\otimes \shift a_1^3\otimes \shift a_2\otimes 1}_{26}
-\underbrace{a_1^1a_1^{21}\dots a_1^{24}\hat{\omega}\otimes \shift a_1^{25}\otimes \shift a_1^3\otimes \shift a_2\otimes 1}_{22}\\
&-\underbrace{a_1^1a_1^{21}a_1^{22} a_1^{23}(\shift \otimes \shift)\partial(\shift^{-1}a_1^{24})\otimes \shift a_1^{25}\otimes \shift a_1^3\otimes \shift a_2\otimes 1}_{75}
\+\underbrace{1\otimes \shift a_1\otimes \shift a_2^1\otimes \shift a_2^{21}\otimes \hat{\omega}a_2^{22}\dots a_2^{25}a_2^3}_{21}\\
&+\underbrace{1\otimes \shift a_1\otimes \shift a_2^1\otimes \shift a_2^{21}\otimes a_2^{22}\hat{\omega}a_2^{23}a_2^{24}a_2^{25}a_2^3}_{30}
+\underbrace{1\otimes \shift a_1\otimes \shift a_2^1 \otimes \shift a_2^{21}\otimes (\shift \otimes \shift)\partial(\shift^{-1}a_2^{22})a_2^{23}a_2^{24} a_2^{25}a_2^3}_{74}
\end{align*}
and
\begin{align*}
&(1^{\otimes 2}\otimes \mu_2)\mu_3(1\otimes \shift a_1\otimes \shift a_2\otimes 1)\\
=\,&(1^{\otimes 2}\otimes \mu_2)(1\otimes \shift a_1^1\otimes a_1^2\otimes \shift a_1^3\otimes \shift a_2\otimes 1-1\otimes \shift a_1\otimes \shift a_2^1\otimes a_2^2\otimes \shift a_2^3\otimes 1\\
&+a_1^1a_1^{21}\dots a_1^{24}\otimes \shift a_1^{25}\otimes \shift a_1^3\otimes \shift a_2\otimes 1
-1\otimes \shift a_1\otimes \shift a_2^1\otimes \shift a_2^{21}\otimes a_2^{22}\dots a_2^{25}a_2^3)\\
=\,&\underbrace{1\otimes \shift a_1^1\otimes a_1^2\hat{\omega}\otimes \shift a_1^3\otimes \shift a_2\otimes 1}_{12}
+\underbrace{1\otimes \shift a_1^1\otimes a_1^2\otimes \shift a_1^3\otimes \shift a_2\otimes \hat{\omega}}_5\\
&-\underbrace{1\otimes \shift a_1^1\otimes a_1^2\otimes \shift a_1^3\otimes 1\otimes \shift a_2\otimes 1}_{10}
-\underbrace{1\otimes \shift a_1^1\otimes a_1^2a_1^{31}a_1^{32}\otimes \shift a_1^{33}\otimes \shift a_2\otimes 1}_y\\
&+\underbrace{1\otimes \shift a_1^1\otimes a_1^2\otimes \shift a_1^3\otimes \shift a_2^1\otimes a_2^2a_2^3}_{18}
\-\underbrace{1\otimes \shift a_1\otimes \shift a_2^1\otimes a_2^2\hat{\omega}\otimes \shift a_2^3\otimes 1}_{16}\\
&-\underbrace{1\otimes \shift a_1\otimes \shift a_2^1\otimes a_2^2\otimes \shift a_2^3\otimes \hat{\omega}}_6
+\underbrace{1\otimes \shift a_1\otimes \shift a_2^1\otimes a_2^2a_2^{31}a_2^{32}\otimes \shift a_2^{33}\otimes 1}_y\\
&-\underbrace{1\otimes \shift a_1\otimes \shift a_2^1\otimes a_2^2\otimes \shift a_2^{31}\otimes a_2^{32}a_2^{33}}_y
\+\underbrace{a_1^1a_1^{21}\dots a_1^{24}\hat{\omega}\otimes \shift a_1^{25}\otimes \shift a_1^3\otimes \shift a_2\otimes 1}_{22}\\
&+\underbrace{a_1^1a_1^{21}\dots a_1^{24}\otimes \shift a_1^{25}\otimes \shift a_1^3\otimes \shift a_2\otimes \hat{\omega}}_7
-\underbrace{a_1^1a_1^{21}\dots a_1^{24}\otimes \shift a_1^{25}\otimes 1\otimes \shift a_1^3\otimes \shift a_2\otimes 1}_{17}\\
&-\underbrace{a_1^1a_1^{21}\dots a_1^{24}\otimes \shift a_1^{25}\otimes \shift a_1^3\otimes 1\otimes \shift a_2\otimes 1}_{11}
-\underbrace{a_1^1a_1^{21}\dots a_1^{24}a_1^{251}a_1^{252}\otimes \shift a_1^{253}\otimes \shift a_1^3\otimes \shift a_2\otimes 1}_{y}\\
&+\underbrace{a_1^1a_1^{21}\dots a_1^{24}\otimes \shift a_1^{25}\otimes \shift a_1^3\otimes \shift a_2^1\otimes a_2^2a_2^3}_{20}
\-\underbrace{1\otimes \shift a_1\otimes \shift a_2^1\otimes \shift a_2^{21}\otimes a_2^{22}a_2^{23}\hat{\omega}a_2^{24}a_2^{25}a_2^3}_{31}\\
&-\underbrace{1\otimes \shift a_1\otimes \shift a_2^1\otimes \shift a_2^{21}\otimes a_2^{22}\dots a_2^{24}\hat{\omega}a_2^{25}a_2^3}_{28}
-\underbrace{1\otimes \shift a_1\otimes \shift a_2^1\otimes \shift a_2^{21}\otimes a_2^{22}a_2^{23}(\shift \otimes \shift)\partial(\shift^{-1}a_2^{24})a_2^{25}a_2^3}_y
\end{align*}
and
\begin{align*}
&\mu_4^1\mu_1^2(1\otimes \shift a_1\otimes \shift a_2\otimes 1)\\
=\,&\mu_4^1(a_1\otimes \shift a_2\otimes 1-1\otimes \shift (a_1a_2)\otimes 1+1\otimes \shift a_1\otimes a_2)\\
=\,&-\underbrace{a_1\otimes \shift a_2^1\otimes \shift a_2^{21}\otimes \shift a_2^{221}\otimes a_2^{222}\dots a_2^{225}a_2^{23}\dots a_2^{25}a_2^3}_{63}\\
&+\underbrace{a_1a_2^1a_2^{21}\dots a_2^{23}a_2^{241}\dots a_2^{244}\otimes \shift a_2^{245}\otimes \shift a_2^{25}\otimes \shift a_2^3\otimes 1}_{69}\\
&+\underbrace{a_1a_2^1a_2^{21}a_2^{22}\otimes \shift a_2^{23}\otimes a_2^{24}\otimes \shift a_2^{25}\otimes \shift a_2^3\otimes 1}_{53}
-\underbrace{a_1\otimes \shift a_2^1\otimes \shift a_2^{21}\otimes a_2^{22}\otimes \shift a_2^{23}\otimes a_2^{24}a_2^{25}a_2^3}_{33}\\
&+\underbrace{a_1\otimes \shift a_2^1\otimes \shift a_2^{21}\otimes a_2^{22}\dots a_2^{25}\otimes \shift a_2^3\otimes 1}_{59}
-\underbrace{a_1\otimes \shift (a_2^1a_2^{21})\otimes a_2^{22}\dots a_2^{24}\otimes \shift a_2^{25}\otimes \shift a_2^3\otimes 1}_{46}\\
&+\underbrace{a_1\otimes \shift a_2^1\otimes \shift a_2^{21}\otimes a_2^{22}\dots a_2^{24}\otimes \shift a_2^{25}\otimes a_2^3}_{58}
\-\underbrace{a_1a_2^1a_2^{21}a_2^{22}\otimes \shift a_2^{23}\otimes a_2^{24}\otimes \shift a_2^{25}\otimes \shift a_2^3\otimes 1}_{53}\\
&-\underbrace{a_1^1a_1^{21}a_1^{22}\otimes \shift a_1^{23}\otimes a_1^{24}\otimes \shift a_1^{25}\otimes \shift (a_1^3a_2)\otimes 1}_{48}+\underbrace{1\otimes \shift (a_1a_2^1)\otimes \shift a_2^{21}\otimes a_2^{22}\otimes \shift a_2^{23}\otimes a_2^{24}a_2^{25}a_2^3}_{47}\\
&+\underbrace{1\otimes \shift a_1^1\otimes \shift a_1^{21}\otimes a_1^{22}\otimes \shift a_1^{23}\otimes a_1^{24}a_1^{25}a_1^3a_2}_{52}-\underbrace{1\otimes \shift (a_1a_2^1)\otimes \shift a_2^{21}\otimes a_2^{22}\dots a_2^{25}\otimes \shift a_2^3\otimes 1}_{55}\\
&-\underbrace{1\otimes \shift a_1^1\otimes \shift a_1^{21}\otimes a_1^{22}\dots a_1^{25}\otimes \shift (a_1^3a_2)\otimes 1}_{38}+\underbrace{1\otimes \shift (a_1a_2^1a_2^{21})\otimes a_2^{22}\dots a_2^{24}\otimes \shift a_2^{25}\otimes \shift a_2^3\otimes 1}_{44}\\
&+\underbrace{1\otimes \shift (a_1^1a_1^{21})\otimes a_1^{22}\dots a_1^{24}\otimes \shift a_1^{25}\otimes \shift (a_1^3a_2)\otimes 1}_{41}-\underbrace{1\otimes \shift (a_1a_2^1)\otimes \shift a_2^{21}\otimes a_2^{22}\dots a_2^{24}\otimes \shift a_2^{25}\otimes a_2^3}_{54}\\
&-\underbrace{1\otimes \shift a_1^1\otimes \shift a_1^{21}\otimes a_1^{22}\dots a_1^{24}\otimes \shift a_1^{25}\otimes a_1^3a_2}_{42}\\
&+\underbrace{1\otimes \shift a_1^1\otimes \shift a_1^{21}\otimes \shift a_1^{221}\otimes a_1^{222}\dots a_1^{225}a_1^{23}\dots a_1^{25}a_1^3a_2}_{70}\\
&+\underbrace{1\otimes \shift (a_1a_2^1)\otimes \shift a_2^{21}\otimes \shift a_2^{221}\otimes a_2^{222}\dots a_2^{225}a_2^{23}\dots a_2^{25}a_2^3}_{71}\\
&-\underbrace{a_1^1a_1^{21}\dots a_1^{23}a_1^{241}\dots a_1^{244}\otimes \shift a_1^{245}\otimes \shift a_1^{25}\otimes \shift (a_1^3a_2)\otimes 1}_{72}\\
&-\underbrace{a_1a_2^1a_2^{21}\dots a_2^{23}a_2^{241}\dots a_2^{244}\otimes \shift a_2^{245}\otimes \shift a_2^{25}\otimes \shift a_2^3\otimes 1}_{69}\\
&\-\underbrace{1\otimes \shift a_1^1\otimes \shift a_1^{21}\otimes \shift a_1^{221}\otimes a_1^{222}\dots a_1^{225}a_1^{23}\dots a_1^{25}a_1^3a_2}_{70}\\
&+\underbrace{a_1^1a_1^{21}\dots a_1^{23}a_1^{241}\dots a_1^{244}\otimes \shift a_1^{245}\otimes \shift a_1^{25}\otimes \shift a_1^3\otimes a_2}_{73}
+\underbrace{a_1^1a_1^{21}a_1^{22}\otimes \shift a_1^{23}\otimes a_1^{24}\otimes \shift a_1^{25}\otimes \shift a_1^3\otimes a_2}_{34}\\
&-\underbrace{1\otimes \shift a_1^1\otimes \shift a_1^{21}\otimes a_1^{22}\otimes \shift a_1^{23}\otimes a_1^{24}a_1^{25}a_1^3a_2}_{52}
+\underbrace{1\otimes \shift a_1^1\otimes \shift a_1^{21}\otimes a_1^{22}\dots a_1^{25}\otimes \shift a_1^3\otimes a_2}_{37}\\
&-\underbrace{1\otimes \shift (a_1^1a_1^{21})\otimes a_1^{22}\dots a_1^{24}\otimes \shift a_1^{25}\otimes \shift a_1^3\otimes a_2}_{43}
+\underbrace{1\otimes \shift a_1^1\otimes \shift a_1^{21}\otimes a_1^{22}\dots a_1^{24}\otimes \shift a_1^{25}\otimes a_1^3a_2}_{42}
\end{align*}

Again, we write $\mu_1^{\otimes}:=(\mu_1\otimes 1^{\otimes 3}+1\otimes \mu_1\otimes 1^{\otimes 2}+1^{\otimes 2}\otimes \mu_1\otimes 1+1^{\otimes 3}\otimes \mu_1)$. Applying this to each individual summand of $\mu_4^2(1\otimes \shift a_1\otimes \shift a_2\otimes 1)$ we obtain:
\begin{align*}
&-\mu_1^{\otimes}(1\otimes \shift a_1\otimes \shift (a_2^1a_2^{21})\otimes a_2^{22}\dots a_2^{24}\otimes \shift a_2^{25}\otimes \shift a_2^3\otimes 1)\\
=\,&-\underbrace{a_1\otimes \shift (a_2^1a_2^{21})\otimes a_2^{22}\dots a_2^{24}\otimes \shift a_2^{25}\otimes \shift a_2^3\otimes 1}_{67}+\underbrace{1\otimes \shift (a_1a_2^1a_2^{21})\otimes a_2^{22}\dots a_2^{24}\otimes \shift a_2^{25}\otimes \shift a_2^3\otimes 1}_{66}\\
&-\underbrace{1\otimes \shift a_1\otimes a_2^1a_2^{21}\dots a_2^{24}\otimes \shift a_2^{25}\otimes \shift a_2^3\otimes 1}_{65}\boldsymbol{\uwave{+}}0\boldsymbol{\uwave{+}}0\boldsymbol{\uwave{-}}
\underbrace{1\otimes \shift a_1\otimes \shift (a_2^1a_2^{21})\otimes a_2^{22}\dots a_2^{25}\otimes \shift a_2^3\otimes 1}_{68}\\
&+\underbrace{1\otimes \shift a_1\otimes \shift (a_2^1a_2^{21})\otimes a_2^{22}\dots a_2^{24}\otimes \shift (a_2^{25}a_2^3)\otimes 1}_y
-\underbrace{1\otimes \shift a_1\otimes \shift (a_2^1a_2^{21})\otimes a_2^{22}\dots a_2^{24}\otimes \shift a_2^{25}\otimes a_2^3}_{64} \end{align*}
and
\begin{align*}
&-\mu_1^{\otimes}(1\otimes \shift a_1^1\otimes \shift a_1^{21}\otimes a_1^{22}\dots a_1^{25}\otimes \shift a_1^3\otimes \shift a_2\otimes 1)\\
=\,&-\underbrace{a_1^1\otimes \shift a_1^{21}\otimes a_1^{22}\dots a_1^{25}\otimes \shift a_1^3\otimes \shift a_2\otimes 1}_{40}
+\underbrace{1\otimes \shift (a_1^1a_1^{21})\otimes a_1^{22}\dots a_1^{25}\otimes \shift a_1^3\otimes \shift a_2\otimes 1}_{35}\\
&+\underbrace{1\otimes \shift a_1^1\otimes (\shift \otimes \shift)\partial(\shift^{-1}a_1^2)\otimes \shift a_1^3\otimes \shift a_2\otimes 1}_{36}
\boldsymbol{\uwave{+}}0\boldsymbol{\uwave{+}}0\boldsymbol{\uwave{-}}\underbrace{1\otimes \shift a_1^1\otimes \shift a_1^{21}\otimes a_1^{22}\dots a_1^{25}a_1^3\otimes \shift a_2\otimes 1}_{39}\\
&+\underbrace{1\otimes \shift a_1^1\otimes \shift a_1^{21}\otimes a_1^{22}\dots a_1^{25}\otimes \shift (a_1^3a_2)\otimes 1}_{38}
-\underbrace{1\otimes \shift a_1^1\otimes \shift a_1^{21}\otimes a_1^{22}\dots a_1^{25}\otimes \shift a_1^3\otimes a_2}_{37} 
\end{align*}
and
\begin{align*}
&-\mu_1^{\otimes}(a_1\otimes \shift a_2^1\otimes \shift a_2^{21}\otimes a_2^{22}\dots a_2^{24}\otimes \shift a_2^{25}\otimes \shift a_2^3\otimes 1)\\
=\,&-\underbrace{a_1a_2^1\otimes \shift a_2^{21}\otimes a_2^{22}\dots a_2^{24}\otimes \shift a_2^{25}\otimes \shift a_2^3\otimes 1}_{61}
+\underbrace{a_1\otimes \shift (a_2^1a_2^{21})\otimes a_2^{22}\dots a_2^{24}\otimes \shift a_2^{25}\otimes \shift a_2^3\otimes 1}_{67}\\
&-\underbrace{a_1\otimes \shift a_2^1\otimes a_2^{21}\dots a_2^{24}\otimes \shift a_2^{25}\otimes \shift a_2^3\otimes 1}_{62}\boldsymbol{\uwave{+}}0\boldsymbol{\uwave{+}}0
\boldsymbol{\uwave{-}}\underbrace{a_1\otimes \shift a_2^1\otimes \shift a_2^{21}\otimes a_2^{22}\dots a_2^{25}\otimes \shift a_2^3\otimes 1}_{59}\\
&+\underbrace{a_1\otimes \shift a_2^1\otimes \shift a_2^{21}\otimes a_2^{22}\dots a_2^{24}\otimes \shift (a_2^{25}a_2^3)\otimes 1}_{60}
-\underbrace{a_1\otimes \shift a_2^1\otimes \shift a_2^{21}\otimes a_2^{22}\dots a_2^{24}\otimes \shift a_2^{25}\otimes a_2^3}_{58}
\end{align*}
and
\begin{align*}
&-\mu_1^{\otimes}(-1\otimes \shift a_1\otimes \shift (a_2^1a_2^{21})\otimes a_2^{22}\dots a_2^{24}\otimes \shift a_2^{25}\otimes \shift a_2^3\otimes 1)\\
=\,&\underbrace{a_1\otimes \shift (a_2^1a_2^{21})\otimes a_2^{22}\dots a_2^{24}\otimes \shift a_2^{25}\otimes \shift a_2^3\otimes 1}_{46}
-\underbrace{1\otimes \shift (a_1a_2^1a_2^{21})\otimes a_2^{22}\dots a_2^{24}\otimes \shift a_2^{25}\otimes \shift a_2^3\otimes 1}_{44}\\
&+\underbrace{1\otimes \shift a_1\otimes a_2^1a_2^{21}\dots a_2^{24}\otimes \shift a_2^{25}\otimes \shift a_2^3\otimes 1}_{45}\boldsymbol{\uwave{+}}0\boldsymbol{\uwave{+}}0
\boldsymbol{\uwave{+}}\underbrace{1\otimes \shift a_1\otimes \shift (a_2^1a_2^{21})\otimes a_2^{22}\dots a_2^{25}\otimes \shift a_2^3\otimes 1}_{68}\\
&-\underbrace{1\otimes \shift a_1\otimes \shift (a_2^1a_2^{21})\otimes a_2^{22}\dots a_2^{24}\otimes \shift (a_2^{25}a_2^3)\otimes 1}_y
+\underbrace{1\otimes \shift a_1\otimes \shift (a_2^1a_2^{21})\otimes a_2^{22}\dots a_2^{24}\otimes \shift a_2^{25}\otimes a_2^3}_{64}
\end{align*}
and
\begin{align*}
&-\mu_1^{\otimes}(-1\otimes \shift a_1\otimes \shift a_2^1\otimes a_2^{21}\dots a_2^{24}\otimes \shift a_2^{25}\otimes \shift a_2^3\otimes 1)\\
=\,&\underbrace{a_1\otimes \shift a_2^1\otimes a_2^{21}\dots a_2^{24}\otimes \shift a_2^{25}\otimes \shift a_2^3\otimes 1}_{62}
-\underbrace{1\otimes \shift (a_1a_2^1)\otimes a_2^{21}\dots a_2^{24}\otimes \shift a_2^{25}\otimes \shift a_2^3\otimes 1}_{57}\\
&+\underbrace{1\otimes \shift a_1\otimes a_2^1a_2^{21}\dots a_2^{24}\otimes \shift a_2^{25}\otimes \shift a_2^3\otimes 1}_{65}\boldsymbol{\uwave{+}}0\boldsymbol{\uwave{+}}0
\boldsymbol{\uwave{-}}\underbrace{1\otimes \shift a_1\otimes \shift a_2^1\otimes (\shift \otimes \shift)\partial(\shift^{-1}a_2^2)\otimes \shift a_2^3\otimes 1}_{50}\\
&-\underbrace{1\otimes \shift a_1\otimes \shift a_2^1\otimes a_2^{21}\dots a_2^{24}\otimes \shift (a_2^{25}a_2^3)\otimes 1}_{51}
+\underbrace{1\otimes \shift a_1\otimes \shift a_2^1\otimes a_2^{21}\dots a_2^{24}\otimes \shift a_2^{25}\otimes a_2^3}_{49}
\end{align*}
and
\begin{align*}
&-\mu_1^{\otimes}(-1\otimes \shift (a_1a_2^1)\otimes \shift a_2^{21}\otimes a_2^{22}\dots a_2^{24}\otimes \shift a_2^{25}\otimes \shift a_2^3\otimes 1)\\
=\,&\underbrace{a_1a_2^1\otimes \shift a_2^{21}\otimes a_2^{22}\dots a_2^{24}\otimes \shift a_2^{25}\otimes \shift a_2^3\otimes 1}_{61}
-\underbrace{1\otimes \shift (a_1a_2^1a_2^{21})\otimes a_2^{22}\dots a_2^{24}\otimes \shift a_2^{25}\otimes \shift a_2^3\otimes 1}_{66}\\
&+\underbrace{1\otimes \shift (a_1a_2^1)\otimes a_2^{21}\dots a_2^{24}\otimes \shift a_2^{25}\otimes \shift a_2^3\otimes 1}_{57}\boldsymbol{\uwave{+}}0\boldsymbol{\uwave{+}}0
\boldsymbol{\uwave{+}}\underbrace{1\otimes \shift (a_1a_2^1)\otimes \shift a_2^{21}\otimes a_2^{22}\dots a_2^{25}\otimes \shift a_2^3\otimes 1}_{55}\\
&-\underbrace{1\otimes \shift (a_1a_2^1)\otimes  \shift a_2^{21}\otimes a_2^{22}\dots a_2^{24}\otimes \shift (a_2^{25}a_2^3)\otimes 1}_{56}
+\underbrace{1\otimes \shift (a_1a_2^1)\otimes \shift a_2^{21}\otimes a_2^{22}\dots a_2^{24}\otimes \shift a_2^{25}\otimes a_2^3}_{54}
\end{align*}
and
\begin{align*}
&-\mu_1^{\otimes}(-1\otimes \shift (a_1^1a_1^{21})\otimes a_1^{22}\dots a_1^{24}\otimes \shift a_1^{25}\otimes \shift a_1^3\otimes \shift a_2\otimes 1)\\
=\,&\underbrace{\hat{\omega}a_1^1a_1^{21}\dots a_1^{24}\otimes \shift a_1^{25}\otimes \shift a_1^3\otimes \shift a_2\otimes 1}_{23}
-\underbrace{a_1^1a_1^{21}\hat{\omega}a_1^{22}\dots a_1^{24}\otimes \shift a_1^{25}\otimes \shift a_1^3\otimes \shift a_2\otimes 1}_{24}\\
&-\underbrace{\shift \partial(a_1^1)a_1^{21}\dots a_1^{24}\otimes \shift a_1^{25}\otimes \shift a_1^3\otimes \shift a_2\otimes 1}_g
-\underbrace{a_1^1\shift \partial(a_1^{21})a_1^{22}\dots a_1^{24}\otimes \shift a_1^{25}\otimes \shift a_1^3\otimes \shift a_2\otimes 1}_y\\
&\boldsymbol{\uwave{+}}0\boldsymbol{\uwave{+}}0\boldsymbol{\uwave{-}}\underbrace{1\otimes \shift (a_1^1a_1^{21})\otimes a_1^{22}\dots a_1^{25}\otimes \shift a_1^3\otimes \shift a_2\otimes 1}_{35}
+\underbrace{1\otimes \shift (a_1^1a_1^{21})\otimes a_1^{22}\dots a_1^{24}\otimes \shift (a_1^{25}a_1^3)\otimes \shift a_2\otimes 1}_y\\
&-\underbrace{1\otimes \shift (a_1^1a_1^{21})\otimes a_1^{22}\dots a_1^{24}\otimes \shift a_1^{25}\otimes \shift (a_1^3a_2)\otimes 1}_{41}
+\underbrace{1\otimes \shift (a_1^1a_1^{21})\otimes a_1^{22}\dots a_1^{24}\otimes \shift a_1^{25}\otimes \shift a_1^3\otimes a_2}_{43}\end{align*}
and
\begin{align*}
&-\mu_1^{\otimes}(a_1^1a_1^{21}a_1^{22}\otimes \shift a_1^{23}\otimes a_1^{24}\otimes \shift a_1^{25}\otimes \shift a_1^3\otimes \shift a_2\otimes 1)\\
=\,&0\boldsymbol{\uwave{-}}\underbrace{a_1^1a_1^{21}a_1^{22}\hat{\omega}a_1^{23}a_1^{24}\otimes \shift a_1^{25}\otimes \shift a_1^3\otimes \shift a_2\otimes 1}_{25}
+\underbrace{a_1^1a_1^{21}\dots a_1^{23}\hat{\omega}a_1^{24}\otimes \shift a_1^{25}\otimes \shift a_1^3\otimes \shift a_2\otimes 1}_{26}\\
&+\underbrace{a_1^1a_1^{21}a_1^{22}\shift \partial(a_1^{23})a_1^{24}\otimes \shift a_1^{25}\otimes \shift a_1^3\otimes \shift a_2\otimes 1}_y\boldsymbol{\uwave{+}}0\boldsymbol{\uwave{+}}
\underbrace{a_1^1a_1^{21}a_1^{22}\otimes \shift a_1^{23}\otimes a_1^{24}a_1^{25}\otimes \shift a_1^3\otimes \shift a_2\otimes 1}_{27}\\
&-\underbrace{a_1^1a_1^{21}a_1^{22}\otimes \shift a_1^{23}\otimes a_1^{24}\otimes \shift (a_1^{25}a_1^3)\otimes \shift a_2\otimes 1}_y
+\underbrace{a_1^1a_1^{21}a_1^{22}\otimes \shift a_1^{23}\otimes a_1^{24}\otimes \shift a_1^{25}\otimes \shift (a_1^3a_2)\otimes 1}_{48}\\
&-\underbrace{a_1^1a_1^{21}a_1^{22}\otimes \shift a_1^{23}\otimes a_1^{24}\otimes \shift a_1^{25}\otimes \shift a_1^3\otimes a_2}_{34}
\end{align*}
and
\begin{align*}
&-\mu_1^{\otimes}(1\otimes \shift a_1\otimes \shift a_2^1\otimes \shift a_2^{21}\otimes a_2^{22}\dots a_2^{24}\otimes \shift (a_2^{25}a_2^3)\otimes 1)\\
=\,&-\underbrace{a_1\otimes \shift a_2^1\otimes \shift a_2^{21}\otimes a_2^{22}\dots a_2^{24}\otimes \shift (a_2^{25}a_2^3)\otimes 1}_{60}
+\underbrace{1\otimes \shift (a_1a_2^1)\otimes \shift a_2^{21}\otimes a_2^{22}\dots a_2^{24}\otimes \shift (a_2^{25}a_2^3)\otimes 1}_{56}\\
&-\underbrace{1\otimes \shift a_1\otimes \shift (a_2^1a_2^{21})\otimes a_2^{22}\dots a_2^{24}\otimes \shift (a_2^{25}a_2^3)\otimes 1}_y
+\underbrace{1\otimes \shift a_1\otimes \shift a_2^1\otimes a_2^{21}\dots a_2^{24}\otimes \shift (a_2^{25}a_2^3)\otimes 1}_{51}\boldsymbol{\uwave{+}}0\boldsymbol{\uwave{+}}0\\
&\boldsymbol{\uwave{+}}\underbrace{1\otimes \shift a_1\otimes \shift a_2^1\otimes \shift a_2^{21}\otimes a_2^{22}\dots a_2^{24}\hat{\omega}a_2^{25}a_2^3}_{28}
-\underbrace{1\otimes \shift a_1\otimes \shift a_2^1\otimes \shift a_2^{21}\otimes a_2^{22}\dots a_2^{25}a_2^3\hat{\omega}}_{29}\\
&-\underbrace{1\otimes \shift a_1\otimes \shift a_2^1\otimes \shift a_2^{21}\otimes a_2^{22}\dots a_2^{24}\shift \partial(a_2^{25})a_2^3}_y
-\underbrace{1\otimes \shift a_1\otimes \shift a_2^1\otimes \shift a_2^{21}\otimes a_2^{22}\dots a_2^{25}\shift \partial(a_2^3)}_g
\end{align*}
and
\begin{align*}
&-\mu_1^{\otimes}(-1\otimes \shift a_1\otimes \shift a_2^1\otimes \shift a_2^{21}\otimes a_2^{22}\otimes \shift a_2^{23}\otimes a_2^{24}a_2^{25}a_2^3)\\
=\,&\underbrace{a_1\otimes \shift a_2^1\otimes \shift a_2^{21}\otimes a_2^{22}\otimes \shift a_2^{23}\otimes a_2^{24}a_2^{25}a_2^3}_{33}
-\underbrace{1\otimes \shift (a_1a_2^1)\otimes \shift a_2^{21}\otimes a_2^{22}\otimes \shift a_2^{23}\otimes a_2^{24}a_2^{25}a_2^3}_{47}\\
&+\underbrace{1\otimes \shift a_1\otimes \shift (a_2^1a_2^{21})\otimes a_2^{22}\otimes \shift a_2^{23}\otimes a_2^{24}a_2^{25}a_2^3}_y
-\underbrace{1\otimes \shift a_1\otimes \shift a_2^1\otimes a_2^{21}a_2^{22}\otimes \shift a_2^{23}\otimes a_2^{24}a_2^{25}a_2^3}_{32}\boldsymbol{\uwave{+}}0\\
&\boldsymbol{\uwave{-}}\underbrace{1\otimes \shift a_1\otimes \shift a_2^1\otimes \shift a_2^{21}\otimes a_2^{22}\hat{\omega}a_2^{23}a_2^{24}a_2^{25}a_2^3}_{30}
+\underbrace{1\otimes \shift a_1\otimes \shift a_2^1\otimes \shift a_2^{21}\otimes a_2^{22}a_2^{23}\hat{\omega}a_2^{24}a_2^{25}a_2^3}_{31}\\
&+\underbrace{1\otimes \shift a_1\otimes \shift a_2^1\otimes \shift a_2^{21}\otimes a_2^{22}\shift \partial(a_2^{23})a_2^{24}a_2^{25}a_2^3}_y\boldsymbol{\uwave{+}}0
\end{align*}
and
\begin{align*}
&-\mu_1^{\otimes}(-1\otimes \shift a_1\otimes \shift a_2^1\otimes \shift a_2^{21}\otimes \shift a_2^{221}\otimes a_2^{222}\dots a_2^{225}a_2^{23}\dots a_2^{25}a_2^3)\\
=\,&+\underbrace{a_1\otimes \shift a_2^1\otimes \shift a_2^{21}\otimes \shift a_2^{221}\otimes a_2^{222}\dots a_2^{225}a_2^{23}\dots a_2^{25}a_2^3}_{63}\\
&
-\underbrace{1\otimes \shift (a_1a_2^1)\otimes \shift a_2^{21}\otimes \shift a_2^{221}\otimes a_2^{222}\dots a_2^{225}a_2^{23}\dots a_2^{25}a_2^3}_{71}\\
&+\underbrace{1\otimes \shift a_1\otimes \shift (a_2^1a_2^{21})\otimes \shift a_2^{221}\otimes a_2^{222}\dots a_2^{225}a_2^{23}\dots a_2^{25}a_2^3}_g\\
&
-\underbrace{1\otimes \shift a_1\otimes \shift a_2^1\otimes \shift (a_2^{21}a_2^{221})\otimes a_2^{222}\dots a_2^{225}a_2^{23}\dots a_2^{25}a_2^3}_y\\
&-\underbrace{1\otimes \shift a_1\otimes \shift a_2^1\otimes \shift a_2^{21}\otimes (\shift \otimes \shift)\partial(\shift^{-1}a_2^{22})a_2^{23}\dots a_2^{25}a_2^3}_{74}\boldsymbol{\uwave{+}}0\boldsymbol{\uwave{+}}0\boldsymbol{\uwave{+}}0
\end{align*}
and
\begin{align*}
&-\mu_1^{\otimes}(a_1a_1^{21}\dots a_1^{23}a_1^{241}\dots a_1^{244}\otimes \shift a_1^{245}\otimes \shift a_1^{25}\otimes \shift a_1^3\otimes \shift a_2\otimes 1)\\
=\,&0\boldsymbol{\uwave{+}}0\boldsymbol{\uwave{+}}0\boldsymbol{\uwave{+}}\underbrace{a_1^1a_1^{21}\dots a_1^{23}(\shift \otimes \shift)\partial(\shift^{-1}a_1^{24})\otimes \shift a_1^{25}\otimes \shift a_1^3\otimes \shift a_2\otimes 1}_{75}\\
&
+\underbrace{a_1a_1^{21}\dots a_1^{23}a_1^{241}\dots a_1^{244}\otimes \shift (a_1^{245}a_1^{25})\otimes \shift a_1^3\otimes \shift a_2\otimes 1}_y\\
&-\underbrace{a_1a_1^{21}\dots a_1^{23}a_1^{241}\dots a_1^{244}\otimes \shift a_1^{245}\otimes \shift (a_1^{25}a_1^3)\otimes \shift a_2\otimes 1}_g\\
&
+\underbrace{a_1^1a_1^{21}a_1{22}a_1^{23}a_1^{241}\dots a_1^{244}\otimes \shift a_1^{245}\otimes \shift a_1^{25}\otimes \shift (a_1^3a_2)\otimes 1}_{72}\\
&-\underbrace{a_1^1a_1^{21}a_1^{22} a_1^{23}a_1^{241}\dots a_1^{244}\otimes \shift a_1^{245}\otimes \shift a_1^{25}\otimes \shift a_1^3\otimes a_2}_{73}
\end{align*}
Noting that terms with the same label add up to zero, this finishes the proof for $\mu_4^2$.
\end{proof}

\section{The maps $\mu_5^0$ and $\mu_5^1$ on $\mathcal{P}$}
\label{appendixD}

We start by stating the following trivial lemma, analogous to Lemmas \ref{yellow} and \ref{green}.

\begin{lem}\label{gg}
Any term obtained from one of the zero terms in Lemmas \ref{yellow} and \ref{green} by applying $\partial$ to one of the tensor factors is zero. The same holds when inserting tensors in the factor that has newly been differentiated.
\end{lem}

Whenever we use this lemma from now on to identify a vanishing term, we will label the latter by $gg$.

We illustrate this by the following examples:
By Lemma \ref{yellow}, $\varphi^1\otimes\varphi^{2}\otimes \varphi^3\otimes \varphi^4\partial(\varphi^5)= 0$, and hence applying $\partial$ to $\varphi^2$, we obtain that $\varphi^1\otimes\varphi^{21}\dots \varphi^{24}\otimes \varphi^{25}\otimes \varphi^3\otimes \varphi^4\partial(\varphi^5)=0$ and hence so is
$\varphi^1\varphi^{21}\dots \varphi^{24}\otimes \varphi^{25}\otimes \varphi^3\otimes \varphi^4\partial(\varphi^5)$, which explains the first occurrence of $gg$ on page \pageref{ggocc1}. 

Similarly,
$1\otimes a^{11}\otimes a^{121}\otimes a^{122}\otimes a^{123}\dots a^{125}a^{13}a^2a^3=0$ by Lemma \ref{green}, and applying $\partial$ to $ a^{122}$, we obtain that $1\otimes a^{11}\otimes a^{121}\otimes a^{1221}\otimes a^{1222}\dots a^{1225}\otimes a^{123}\dots a^{125}a^{13}a^2a^3=0$ and thus also 
$1\otimes a^{11}\otimes a^{121}\otimes a^{1221}\otimes a^{1222}\dots a^{1225}a^{123}\dots a^{125}a^{13}a^2a^3=0$, which explains the first occurrence of $gg$ in $-(\mu_4\otimes 1)\mu_2^1(1\otimes \shift a\otimes 1)$.

\begin{lem}
For $\varphi\in \Phi$ and $a\in J$ we obtain:
\begin{enumerate}[(a)]
\item\label{mu5:0} \begin{align*}
\mu_5^0(\varphi)&=-1\otimes \shift \varphi^1\otimes \varphi^2\otimes \shift \varphi^3\otimes \varphi^4\otimes \shift \varphi^5\otimes 1\\
&-1\otimes \shift \varphi^1\otimes \varphi^2\varphi^3\varphi^{41}\dots \varphi^{44}\otimes \shift \varphi^{45}\otimes \shift \varphi^5\otimes 1\\
&-1\otimes \shift \varphi^1\otimes \shift \varphi^{21}\otimes \varphi^{22}\dots \varphi^{25}\varphi^3\varphi^4\otimes \shift \varphi^5\otimes 1\\
&-1\otimes \shift \varphi^1\otimes \varphi^{21}\dots \varphi^{24}\otimes \shift \varphi^{25}\otimes \shift \varphi^3\otimes \varphi^4\varphi^5\\
&-\varphi^1\varphi^2\otimes \shift \varphi^3\otimes \shift \varphi^{41}\otimes \varphi^{42}\dots \varphi^{45}\otimes \shift \varphi^5\otimes 1\\
&-\varphi^1\varphi^{21}\dots \varphi^{24}\otimes \shift \varphi^{25}\otimes \shift \varphi^3\otimes \varphi^4\otimes \shift \varphi^5\otimes 1\\
&-1\otimes \shift \varphi^1\otimes \varphi^2\otimes \shift \varphi^3\otimes \shift \varphi^{41}\otimes \varphi^{42}\dots \varphi^{45}\varphi^5\\
&-\varphi^1\otimes \shift \varphi^{21}\otimes \varphi^{22}\dots \varphi^{24}\otimes \shift \varphi^{25}\otimes \shift \varphi^3\otimes \varphi^4\varphi^5\\
&-\varphi^1\varphi^2\otimes \shift \varphi^3\otimes \shift \varphi^{41}\otimes \varphi^{42}\dots \varphi^{44}\otimes \shift \varphi^{45}\otimes \varphi^5\\
&+\varphi^1\varphi^2\otimes \shift (\varphi^3\varphi^{41})\otimes \varphi^{42}\dots \varphi^{44}\otimes \shift \varphi^{45}\otimes \shift \varphi^5\otimes 1\\
&+1\otimes \shift \varphi^1\otimes \shift \varphi^{21}\otimes \varphi^{22}\dots \varphi^{24}\otimes \shift (\varphi^{25}\varphi^3)\otimes \varphi^4\varphi^5\\
&+\varphi^1\varphi^{21}\varphi^{22}\otimes \shift \varphi^{23}\otimes \varphi^{24}\otimes \shift \varphi^{25}\otimes \shift \varphi^3\otimes \varphi^4\varphi^5\\
&+\varphi^1\varphi^2\otimes \shift \varphi^3\otimes \shift \varphi^{41}\otimes \varphi^{42}\otimes \shift \varphi^{43}\otimes \varphi^{44}\varphi^{45}\varphi^{5}\\
&-1\otimes \shift \varphi^1\otimes \shift \varphi^{21}\otimes \varphi^{22}\otimes \shift \varphi^{23}\otimes \varphi^{24}\varphi^{25}\varphi^3\varphi^4\varphi^5\\
&-\varphi^1\varphi^2\varphi^3\varphi^{41}\varphi^{42}\otimes \shift \varphi^{43}\otimes \varphi^{44}\otimes \shift \varphi^{45}\otimes \shift \varphi^5\otimes 1\\
&-\varphi^1\varphi^{21}\dots \varphi^{24}\otimes \shift \varphi^{25}\otimes \shift \varphi^3\otimes \shift \varphi^{41}\otimes \varphi^{42}\dots \varphi^{45} \varphi^5\\
&-\varphi^1\varphi^2\varphi^3\varphi^{41}\dots\varphi^{43}\varphi^{441}\dots \varphi^{444}\otimes \shift \varphi^{445}\otimes \shift \varphi^{45}\otimes \shift \varphi^5\otimes 1\\
&-1\otimes \shift \varphi^1\otimes \shift \varphi^{21}\otimes \shift \varphi^{221}\otimes \varphi^{222}\dots \varphi^{225}\varphi^{23}\dots \varphi^{25}\varphi^3\varphi^4\varphi^5\\
&+\varphi^1\varphi^{21}\dots \varphi^{23}\varphi^{241}\dots \varphi^{244}\otimes \shift \varphi^{245}\otimes \shift \varphi^{25}\otimes \shift \varphi^3\otimes \varphi^4\varphi^5\\
&+\varphi^1\varphi^2\otimes \shift\varphi^3\otimes \shift \varphi^{41}\otimes \shift \varphi^{421}\otimes \varphi^{422}\dots \varphi^{425} \varphi^{43}\dots \varphi^{45}\varphi^5
\end{align*}
\item\label{mu5:1} \begin{align*}
\mu_5^1(1\otimes \shift a\otimes 1)&=1\otimes \shift a^1\otimes \shift a^{21}\otimes a^{22}\dots a^{24}\otimes \shift a^{25}\otimes 1\otimes \shift a^3\otimes 1\\
&-1\otimes \shift (a^1a^{21})\otimes a^{22}\otimes \shift a^{23}\otimes a^{24}\otimes \shift a^{25}\otimes \shift a^3\otimes 1\\
&+1\otimes \shift (a^1a^{21})\otimes a^{22}a^{23}a^{241}\dots a^{244}\otimes \shift a^{245}\otimes \shift a^{25}\otimes \shift a^3\otimes 1\\
&-a^1a^{21}a^{22}\otimes \shift (a^{23}a^{241})\otimes a^{242}\dots a^{244}\otimes \shift a^{245}\otimes \shift a^{25}\otimes \shift a^3\otimes 1\\
&-1\otimes \shift a^1\otimes \shift a^{21}\otimes a^{22}\otimes \shift a^{23}\otimes a^{24}a^{25}\otimes \shift a^3\otimes 1\\
&-1\otimes \shift a^1\otimes \shift a^{21}\otimes a^{22}\otimes \shift a^{23}\otimes a^{24}\otimes \shift a^{25}\otimes a^3\\
&+a^1a^{21}\dots a^{23}a^{241}a^{242}\otimes \shift a^{243}\otimes a^{244}\otimes \shift a^{245}\otimes \shift a^{25}\otimes \shift a^3\otimes 1\\
&-1\otimes \shift a^1\otimes \shift a^{21}\otimes \shift a^{221}\otimes a^{222}\dots a^{224}\otimes \shift (a^{225}a^{23})\otimes a^{24}a^{25}a^3\\
&+1\otimes \shift a^1\otimes \shift a^{21}\otimes \shift a^{221}\otimes a^{222}\dots a^{225}a^{23}a^{24}\otimes \shift (a^{25}a^3)\otimes 1\\
&+1\otimes \shift a^1\otimes \shift a^{21}\otimes \shift a^{221}\otimes a^{222}\otimes \shift a^{223}\otimes a^{224}a^{225}a^{23}\dots a^{25}a^3\\
&+a^1\otimes \shift a^{21}\otimes \shift a^{221}\otimes a^{222}\dots a^{225}a^{23}a^{24}\otimes \shift a^{25}\otimes \shift a^3\otimes 1\\
&-1\otimes \shift a^1\otimes \shift a^{21}\otimes (\shift \otimes \shift)\partial(\shift^{-1} a^{22})a^{23}a^{24}\otimes \shift a^{25}\otimes \shift a^3\otimes 1\\
&+a^1a^{21}a^{221}\dots a^{224}\otimes \shift a^{225}\otimes \shift a^{23}\otimes a^{24}\otimes \shift a^{25}\otimes \shift a^3\otimes 1\\
&+1\otimes \shift (a^1a^{21})\otimes \shift a^{221}\otimes a^{222}\dots a^{224}\otimes \shift a^{225}\otimes \shift a^{23}\otimes a^{24}a^{25}a^3\\
&+1\otimes \shift a^1\otimes \shift a^{21}\otimes \shift a^{221}\dots a^{224}\otimes \shift a^{225}\otimes \shift a^{23}\otimes a^{24}a^{25}a^3\\
&+a^1a^{21}a^{22}\otimes \shift a^{23}\otimes \shift a^{241}\otimes a^{242}\dots a^{245}\otimes \shift a^{25}\otimes \shift a^3\otimes 1\\
&+1\otimes \shift a^1\otimes \shift a^{21}\otimes a^{22}a^{23}a^{241}\dots a^{244}\otimes \shift a^{245}\otimes \shift a^{25}\otimes a^3\\
&+1\otimes \shift a^1\otimes \shift a^{21}\otimes a^{22}\otimes \shift a^{23}\otimes \shift a^{241}\otimes a^{242}\dots a^{245}a^{25}a^3\\
&-1\otimes \shift a^1\otimes \shift a^{21}\otimes \shift a^{221}\otimes \shift a^{2221}\otimes a^{2222}\dots a^{2225}a^{223}\dots a^{225}a^{23}\dots a^{25}a^3\\
&-a^1a^{21}\dots a^{23}a^{241}\dots a^{243}a^{2441}\dots a^{2444}\otimes \shift a^{2445}\otimes \shift a^{245}\otimes \shift a^{25}\otimes \shift a^3\otimes 1
\end{align*}
\end{enumerate}
\end{lem}

\begin{proof}
We verify that the stated formulae verify the $A_\infty$-relation 
\begin{align*}
&\mu_5\mu_1+(\mu_1\otimes 1^{\otimes 4}+1\otimes \mu_1\otimes 1^{\otimes 3}+1^{\otimes 2}\otimes \mu_1\otimes 1^{\otimes 2}+1^{\otimes 3}\otimes \mu_1\otimes 1+1^{\otimes 4}\otimes \mu_1)\mu_5\\
&+(-\mu_2\otimes 1^{\otimes 3}+1\otimes \mu_2\otimes 1^{\otimes 2}-1^{\otimes 2}\otimes \mu_2\otimes 1+1^{\otimes 3}\otimes \mu_2)\mu_4+(-\mu_4\otimes 1+1\otimes \mu_4)\mu_2\\
&+(\mu_3\otimes 1^{\otimes 2}+1\otimes \mu_3\otimes 1+1^{\otimes 2}\otimes \mu_3)=0.
\end{align*}
\eqref{mu5:0}: For $\mu_5^0$ we compute:
\begin{align*}
&(-\mu_2\otimes 1^{\otimes 3})\mu_4^0(\varphi)\\
=\,&-\underbrace{\varphi^1\hat{\omega}\varphi^2\otimes \shift \varphi^3\otimes \varphi^4\otimes \shift \varphi^5\otimes 1}_{68}
-\underbrace{\varphi^1\varphi^2\hat{\omega}\otimes \shift \varphi^3\otimes \varphi^4\otimes \shift \varphi^5\otimes 1}_{10}\\
&-\underbrace{\varphi^1(\shift \otimes \shift)\partial(\shift^{-1}\varphi^2)\otimes \shift \varphi^3\otimes \varphi^4\otimes \shift \varphi^5\otimes 1}_{73}
\+\underbrace{\hat{\omega}\otimes \shift \varphi^1\otimes \varphi^2\varphi^3\varphi^4\otimes \shift \varphi^5\otimes 1}_{16}\\
&+\underbrace{1\otimes \shift \varphi^1\otimes \hat{\omega}\varphi^2\varphi^3\varphi^4\otimes \shift \varphi^5\otimes 1}_6
-\underbrace{\varphi^{11}\varphi^{12}\otimes \shift \varphi^{13}\otimes \varphi^2\varphi^3\varphi^4\otimes \shift \varphi^5\otimes 1}_y\\
&+\underbrace{1\otimes \shift \varphi^{11}\otimes \varphi^{12}\varphi^{13}\varphi^2\varphi^3\varphi^4\otimes \shift \varphi^5\otimes 1}_y
\-\underbrace{\hat{\omega}\otimes \shift \varphi^1\otimes \varphi^2\otimes \shift \varphi^3\otimes \varphi^4\varphi^5}_{14}\\
&-\underbrace{1\otimes \shift \varphi^1\otimes \hat{\omega}\varphi^2\otimes \shift \varphi^3\otimes \varphi^4\varphi^5}_1
+\underbrace{\varphi^{11}\varphi^{12}\otimes \shift \varphi^{13}\otimes \varphi^2\otimes \shift \varphi^3\otimes \varphi^4\varphi^5}_y\\
&-\underbrace{1\otimes \shift \varphi^{11}\otimes \varphi^{12}\varphi^{13}\varphi^2\otimes \shift \varphi^3\otimes \varphi^4\varphi^5}_y
\-\underbrace{\varphi^1\varphi^{21}\hat{\omega}\varphi^{22}\varphi^{23}\varphi^{24}\otimes \shift \varphi^{25}\otimes \shift \varphi^3\otimes \varphi^4\varphi^5}_{84}\\
&-\underbrace{\varphi^1\varphi^{21}\varphi^{22}\hat{\omega}\varphi^{23}\varphi^{24}\otimes \shift \varphi^{25}\otimes \shift \varphi^3\otimes \varphi^4\varphi^5}_{90}
-\underbrace{\varphi^1\varphi^{21}(\shift \otimes \shift)\partial(\shift^{-1}\varphi^{22})\varphi^{23}\varphi^{24}\otimes \shift \varphi^{25}\otimes \shift \varphi^3\otimes \varphi^4\varphi^5}_y\\
&\-\underbrace{\hat{\omega}\otimes \shift \varphi^1\otimes \shift \varphi^{21}\otimes \varphi^{22}\dots \varphi^{25}\varphi^3\varphi^4\varphi^5}_{21}
-\underbrace{1\otimes \shift \varphi^1\otimes \shift \varphi^{21}\otimes \hat{\omega}\varphi^{22}\dots \varphi^{25}\varphi^3\varphi^4\varphi^5}_{24}\\
&+\underbrace{1\otimes \shift \varphi^1\otimes 1\otimes \shift \varphi^{21}\otimes \varphi^{22}\dots \varphi^{25}\varphi^3\varphi^4\varphi^5}_{18}
+\underbrace{\varphi^{11}\varphi^{12}\otimes \shift \varphi^{13}\otimes \shift \varphi^{21}\otimes \varphi^{22}\dots \varphi^{25}\varphi^3\varphi^4\varphi^5}_g\\
&-\underbrace{1\otimes \shift \varphi^1\otimes \shift \varphi^{211}\otimes \varphi^{212}\varphi^{213}\varphi^{22}\dots \varphi^{25}\varphi^3\varphi^4\varphi^5}_y
\-\underbrace{\varphi^1\hat{\omega}\varphi^2\otimes \shift \varphi^3\otimes \shift \varphi^{41}\otimes \varphi^{42}\dots \varphi^{45}\varphi^5}_{92}\\
&-\underbrace{\varphi^1\varphi^2\hat{\omega}\otimes \shift \varphi^3\otimes \shift \varphi^{41}\otimes \varphi^{42}\dots \varphi^{45}\varphi^5}_7
-\underbrace{\varphi^1(\shift \otimes \shift)\partial(\shift^{-1}\varphi^2)\otimes \shift \varphi^3\otimes \shift \varphi^{41}\otimes \varphi^{42}\dots \varphi^{45}\varphi^5}_{45}\\
&\-\underbrace{\varphi^1\hat{\omega}\varphi^2\varphi^3\varphi^{41}\dots \varphi^{44}\otimes \shift \varphi^{45}\otimes \shift \varphi^5\otimes 1}_{76}
-\underbrace{\varphi^1\varphi^2\hat{\omega}\varphi^3\varphi^{41}\dots \varphi^{44}\otimes \shift \varphi^{45}\otimes \shift \varphi^5\otimes 1}_{87}\\
&-\underbrace{\varphi^1(\shift \otimes \shift)\partial(\shift^{-1}\varphi^2)\varphi^3\varphi^{41}\dots \varphi^{44}\otimes \varphi^{45}\otimes \varphi^5\otimes 1}_g
\end{align*}
and
\begin{align*}
(1\otimes& \mu_2\otimes 1^{\otimes 2})\mu_4^0(\varphi)\\
=\,&+\underbrace{\varphi^1\varphi^2\hat{\omega}\otimes \shift \varphi^3\otimes \varphi^4\otimes \shift \varphi^5\otimes 1}_{10}
+\underbrace{\varphi^1\varphi^2\otimes \shift \varphi^3\otimes \hat{\omega}\varphi^4\otimes \shift \varphi^5\otimes 1}_3\\
&-\underbrace{\varphi^1\varphi^2\varphi^{31}\varphi^{32}\otimes \shift \varphi^{33}\otimes \varphi^4\otimes \shift \varphi^5\otimes 1}_y
+\underbrace{\varphi^1\varphi^2\otimes \shift \varphi^{31}\otimes \varphi^{32}\varphi^{33}\varphi^4\otimes \shift \varphi^5\otimes 1}_y\\
&\-\underbrace{1\otimes \shift \varphi^1\otimes \hat{\omega}\varphi^2\varphi^3\varphi^4\otimes \shift \varphi^5\otimes 1}_6
-\underbrace{1\otimes \shift \varphi^1\otimes \varphi^2\hat{\omega}\varphi^3\varphi^4\otimes \shift \varphi^5\otimes 1}_{69}\\
&-\underbrace{1\otimes \shift \varphi^1\otimes (\shift \otimes \shift)\partial(\shift^{-1}\varphi^2)\varphi^3\varphi^4\otimes \shift \varphi^5\otimes 1}_{54}
\+\underbrace{1\otimes \shift \varphi^1\otimes \hat{\omega}\varphi^2\otimes \shift \varphi^3\otimes \varphi^4\varphi^5}_1\\
&+\underbrace{1\otimes \shift \varphi^1\otimes \varphi^2\hat{\omega}\otimes \shift \varphi^3\otimes \varphi^4\varphi^5}_5
+\underbrace{1\otimes \shift \varphi^1\otimes (\shift \otimes \shift)\partial(\shift^{-1}\varphi^2)\otimes \shift\varphi^3\otimes \varphi^4\varphi^5}_{81}\\
&\+\underbrace{\varphi^1\varphi^{21}\varphi^{22}\varphi^{23}\hat{\omega}\varphi^{24}\otimes \shift \varphi^{25}\otimes \shift \varphi^3\otimes \varphi^4\varphi^5}_{38}
+\underbrace{\varphi^1\varphi^{21}\dots \varphi^{24}\hat{\omega}\otimes \shift \varphi^{25}\otimes \shift \varphi^3\otimes \varphi^4\varphi^5}_{12}\\
&+\underbrace{\varphi^1\varphi^{21}\varphi^{22}\varphi^{23}(\shift \otimes \shift)\partial(\shift^{-1}\varphi^{24})\otimes \shift \varphi^{25}\otimes \shift \varphi^3\otimes \varphi^4\varphi^5}_{48}
\+\underbrace{1\otimes \shift \varphi^1\otimes \shift \varphi^{21}\otimes \hat{\omega}\varphi^{22}\dots \varphi^{25}\varphi^3\varphi^4\varphi^5}_{24}\\
&+\underbrace{1\otimes \shift \varphi^1\otimes \shift \varphi^{21}\otimes \varphi^{22}\hat{\omega}\varphi^{23}\varphi^{24}\varphi^{25}\varphi^3\varphi^4\varphi^5}_{43}
+\underbrace{1\otimes \shift \varphi^1\otimes \shift \varphi^{21}\otimes (\shift \otimes \shift )\partial(\shift^{-1}\varphi^{22})\varphi^{23}\dots \varphi^{25}\varphi^3\varphi^4\varphi^5}_{33}\\
&\+\underbrace{\varphi^1\varphi^2\hat{\omega}\otimes \shift \varphi^3\otimes \shift \varphi^{41}\otimes \varphi^{42}\dots \varphi^{45}\varphi^5}_7
+\underbrace{\varphi^1\varphi^2\otimes \shift \varphi^3\otimes \shift \varphi^{41}\otimes \hat{\omega}\varphi^{42}\dots \varphi^{45}\varphi^5}_{11}\\
&-\underbrace{\varphi^1\varphi^2\otimes \shift \varphi^3\otimes 1\otimes \shift \varphi^{41}\otimes \varphi^{42}\dots \varphi^{45}\varphi^5}_{19}
-\underbrace{\varphi^1\varphi^2\varphi^{31}\varphi^{32}\otimes \shift \varphi^{33}\otimes \shift \varphi^{41}\otimes \varphi^{42}\dots \varphi^{45}\varphi^5}_g\\
&+\underbrace{\varphi^1\varphi^2\otimes \shift \varphi^3\otimes \shift \varphi^{411}\otimes \varphi^{412}\varphi^{413}\varphi^{42}\dots \varphi^{45}\varphi^5}_{y}
\+\underbrace{\varphi^1\varphi^2\varphi^3\varphi^{41}\hat{\omega}\varphi^{42}\dots \varphi^{44}\otimes \shift \varphi^{45}\otimes \shift \varphi^5\otimes 1}_{88}\\
&+\underbrace{\varphi^1\varphi^2\varphi^3\varphi^{41}\varphi^{42}\hat{\omega}\varphi^{43}\varphi^{44}\otimes \shift \varphi^{45}\otimes \shift \varphi^{5}\otimes 1}_{37}
+\underbrace{\varphi^1\varphi^2\varphi^3\varphi^{41}(\shift \otimes \shift)\partial(\shift^{-1} \varphi^{42})\varphi^{43}\varphi^{44}\otimes \varphi^{45}\otimes \shift \varphi^5\otimes 1}_y
\end{align*}
and
\begin{align*}
&-(1^{\otimes 2}\otimes \mu_2\otimes 1)\mu_4^0(\varphi)\\
=\,&-\underbrace{\varphi^1\varphi^2\otimes \shift \varphi^3\otimes \hat{\omega}\varphi^4\otimes \shift \varphi^5\otimes 1}_3
-\underbrace{\varphi^1\varphi^2\otimes \shift \varphi^3\otimes \varphi^4\hat{\omega}\otimes \shift \varphi^5\otimes 1}_2\\
&-\underbrace{\varphi^1\varphi^2\otimes \shift \varphi^3\otimes (\shift \otimes \shift)\partial(\shift^{-1}\varphi^4)\otimes \shift \varphi^5\otimes 1}_{58}
\+\underbrace{1\otimes \shift \varphi^1\otimes \varphi^2\varphi^3\hat{\omega}\varphi^4\otimes \shift \varphi^5\otimes 1}_{71}\\
&+\underbrace{1\otimes \shift \varphi^1\otimes \varphi^2\varphi^3\varphi^4\hat{\omega}\otimes \shift \varphi^5\otimes 1}_4
+\underbrace{1\otimes \shift \varphi^1\otimes \varphi^2\varphi^3(\shift \otimes \shift)\partial(\shift^{-1}\varphi^4)\otimes \shift \varphi^5\otimes 1}_{77}\\
&\-\underbrace{1\otimes \shift \varphi^1\otimes \varphi^2\hat{\omega}\otimes \shift \varphi^3\otimes \varphi^4\varphi^5}_5
-\underbrace{1\otimes \shift \varphi^1\otimes \varphi^2\otimes \shift \varphi^3\otimes \hat{\omega}\varphi^4\varphi^5}_9\\
&+\underbrace{1\otimes \shift \varphi^1\otimes \varphi^2\varphi^{31}\varphi^{32}\otimes \shift \varphi^{33}\otimes \varphi^4\varphi^5}_y
-\underbrace{1\otimes \shift \varphi^1\otimes \varphi^2\otimes \shift \varphi^{31}\otimes \varphi^{32}\varphi^{33}\varphi^4\varphi^5}_y\\
&\-\underbrace{\varphi^1\varphi^{21}\dots \varphi^{24}\hat{\omega}\otimes \shift \varphi^{25}\otimes \shift \varphi^3\otimes \varphi^4\varphi^5}_{12}
-\underbrace{\varphi^1\varphi^{21}\dots \varphi^{24}\otimes \shift \varphi^{25}\otimes \shift \varphi^3\otimes \hat{\omega}\varphi^4\varphi^5}_8\\
&+\underbrace{\varphi^1\varphi^{21}\dots \varphi^{24}\otimes \shift \varphi^{25}\otimes 1\otimes \shift \varphi^3\otimes \varphi^4\varphi^5}_{20}
+\underbrace{\varphi^1\varphi^{21}\dots \varphi^{24}\varphi^{251}\varphi^{252}\otimes \shift \varphi^{253}\otimes \shift \varphi^3\otimes \varphi^4\varphi^5}_{y}\\
&-\underbrace{\varphi^1\varphi^{21}\dots \varphi^{24}\otimes \shift \varphi^{25}\otimes \shift \varphi^{31}\otimes \varphi^{32}\varphi^{33}\varphi^4\varphi^5}_g
\-\underbrace{1\otimes \shift \varphi^1\otimes \shift \varphi^{21}\otimes \varphi^{22}\varphi^{23}\hat{\omega}\varphi^{24}\varphi^{25}\varphi^3\varphi^4\varphi^5}_{50}\\
&-\underbrace{1\otimes \shift \varphi^1\otimes \shift \varphi^{21}\otimes \varphi^{22}\dots \varphi^{24}\hat{\omega}\varphi^{25}\varphi^3\varphi^4\varphi^5}_{65}
-\underbrace{1\otimes \shift \varphi^1\otimes \shift \varphi^{21}\otimes \varphi^{22}\varphi^{23}(\shift \otimes \shift)\partial(\shift^{-1}\varphi^{24})\varphi^{25}\varphi^3\varphi^4\varphi^5}_y\\
&\-\underbrace{\varphi^1\varphi^2\otimes \shift \varphi^3\otimes \shift \varphi^{41}\otimes \hat{\omega}\varphi^{42}\dots \varphi^{45}\varphi^5}_{11}
-\underbrace{\varphi^1\varphi^2\otimes \shift \varphi^3\otimes \shift \varphi^{41}\otimes \varphi^{42}\hat{\omega}\varphi^{43}\varphi^{44} \varphi^{45}\varphi^5}_{25}\\
&-\underbrace{\varphi^1\varphi^2\otimes \shift \varphi^3\otimes \shift \varphi^{41}\otimes (\shift \otimes \shift)\partial(\shift^{-1}\varphi^{42})\varphi^{43}\dots \varphi^{45}\varphi^5}_{35}
\-\underbrace{\varphi^1\varphi^2\varphi^3\varphi^{41}\varphi^{42} \varphi^{43}\hat{\omega}\varphi^{44}\otimes \shift \varphi^{45}\otimes \shift \varphi^5\otimes 1}_{30}\\
&-\underbrace{\varphi^1\varphi^2\varphi^3\varphi^{41}\dots \varphi^{44}\hat{\omega}\otimes \shift \varphi^{45}\otimes \shift \varphi^5\otimes 1}_{23}
-\underbrace{\varphi^1\varphi^2\varphi^3\varphi^{41}\dots \varphi^{43}(\shift \otimes \shift)\partial(\shift^{-1}\varphi^{44})\otimes \shift \varphi^{45}\otimes \shift \varphi^5\otimes 1}_{46}
\end{align*}
and
\begin{align*}
&(1^{\otimes 3}\otimes \mu_2)\mu_4^0(\varphi)\\
=\,&+\underbrace{\varphi^1\varphi^2\otimes \shift \varphi^3\otimes \varphi^4\hat{\omega}\otimes \shift \varphi^5\otimes 1}_2
+\underbrace{\varphi^1\varphi^2\otimes \shift \varphi^3\otimes \varphi^4\otimes \shift \varphi^5\otimes \hat{\omega}}_{15}\\
&-\underbrace{\varphi^1\varphi^2\otimes \shift \varphi^3\otimes \varphi^4\varphi^{51}\varphi^{52}\otimes \shift \varphi^{53}\otimes 1}_y
+\underbrace{\varphi^1\varphi^2\otimes \shift \varphi^3\otimes \varphi^4\otimes \shift \varphi^{51}\otimes \varphi^{52}\varphi^{53}}_y\\
&\-\underbrace{1\otimes \shift \varphi^1\otimes \varphi^2\varphi^3\varphi^4\hat{\omega}\otimes \shift \varphi^5\otimes 1}_4
-\underbrace{1\otimes \shift \varphi^1\otimes \varphi^2\varphi^3\varphi^4\otimes \shift \varphi^5\otimes \hat{\omega}}_{13}\\
&+\underbrace{1\otimes \shift \varphi^1\otimes \varphi^2\varphi^3\varphi^4\varphi^{51}\varphi^{52}\otimes \shift \varphi^{53}\otimes 1}_y
-\underbrace{1\otimes \shift \varphi^1\otimes \varphi^2\varphi^3\varphi^4\otimes \shift \varphi^{51}\otimes \varphi^{52}\varphi^{53}}_y\\
&\+\underbrace{1\otimes \shift \varphi^1\otimes \varphi^2\otimes \shift \varphi^3\otimes \hat{\omega}\varphi^4\varphi^5}_9
+\underbrace{1\otimes \shift \varphi^1\otimes \varphi^2\otimes \shift \varphi^3\otimes \varphi^4\hat{\omega}\varphi^5}_{70}\\
&+\underbrace{1\otimes \shift \varphi^1\otimes \varphi^2\otimes \shift \varphi^3\otimes (\shift \otimes \shift)\partial(\shift^{-1}\varphi^4)\varphi^5}_{93}
\+\underbrace{\varphi^1\varphi^{21}\dots \varphi^{24}\otimes \shift \varphi^{25}\otimes \shift \varphi^3\otimes \hat{\omega}\varphi^4\varphi^5}_8\\
&+\underbrace{\varphi^1\varphi^{21}\dots \varphi^{24}\otimes \shift \varphi^{25}\otimes \shift \varphi^3\otimes \varphi^4\hat{\omega}\varphi^5}_{51}
+\underbrace{\varphi^1\varphi^{21}\dots \varphi^{24}\otimes \shift \varphi^{25}\otimes \shift \varphi^3\otimes (\shift\otimes\shift)\partial(\shift^{-1}\varphi^4)\varphi^5}_{32}\\
&\+\underbrace{1\otimes \shift \varphi^1\otimes \shift \varphi^{21}\otimes \varphi^{22}\dots \varphi^{25}\varphi^3\hat{\omega}\varphi^4\varphi^5}_{64}
+\underbrace{1\otimes \shift \varphi^1\otimes \shift \varphi^{21}\otimes \varphi^{22}\dots \varphi^{25}\varphi^3\varphi^4\hat{\omega}\varphi^5}_{53}\\
&+\underbrace{1\otimes \shift \varphi^1\otimes \shift \varphi^{21}\otimes \varphi^{22}\dots \varphi^{25}\varphi^3(\shift \otimes \shift)\partial(\shift^{-1}\varphi^4)\varphi^5}_{g}
\+\underbrace{\varphi^1\varphi^2\otimes \shift \varphi^3\otimes \shift \varphi^{41}\otimes \varphi^{42}\varphi^{43}\hat{\omega}\varphi^{44}\varphi^{45}\varphi^5}_{67}\\
&+\underbrace{\varphi^1\varphi^2\otimes \shift \varphi^3\otimes \shift \varphi^{41}\otimes \varphi^{42}\dots \varphi^{44}\hat{\omega}\varphi^{45}\varphi^5}_{61}
+\underbrace{\varphi^1\varphi^2\otimes \shift \varphi^3\otimes \shift \varphi^{41}\otimes \varphi^{42}\varphi^{43}(\shift \otimes \shift)\partial(\shift^{-1}\varphi^{44})\varphi^{45}\varphi^5}_y\\
&\+\underbrace{\varphi^1\varphi^2\varphi^3\varphi^{41}\dots \varphi^{44}\hat{\omega}\otimes \shift \varphi^{45}\otimes \shift \varphi^5\otimes 1}_{23}
+\underbrace{\varphi^1\varphi^2\varphi^3\varphi^{41}\dots \varphi^{44}\otimes \shift \varphi^{45}\otimes \shift \varphi^5\otimes \hat{\omega}}_{22}\\
&-\underbrace{\varphi^1\varphi^2\varphi^3\varphi^{41}\dots \varphi^{44}\otimes \shift \varphi^{45}\otimes 1\otimes \shift \varphi^5\otimes 1}_{17}
-\underbrace{\varphi^1\varphi^2\varphi^3\varphi^{41}\dots \varphi^{44}\varphi^{451}\varphi^{452}\otimes \shift \varphi^{453}\otimes \shift \varphi^5\otimes 1}_y\\
&+\underbrace{\varphi^1\varphi^2\varphi^3\varphi^{41}\dots \varphi^{44}\otimes \shift \varphi^{45}\otimes \shift \varphi^{51}\otimes \varphi^{52}\varphi^{53}}_g
\end{align*}
and
\begin{align*}
-(\mu_4\otimes 1)\mu_2(\varphi)=\,&-(\mu_4\otimes 1)(\hat{\omega}\varphi+\varphi\hat{\omega}+(\shift \otimes \shift)\partial(\shift^{-1}\varphi))\\
=\,&0\-\underbrace{1\otimes \shift \varphi^1\otimes \varphi^2\otimes \shift \varphi^3\otimes \varphi^4\varphi^5\hat{\omega}}_{72}
-\underbrace{\varphi^1\varphi^2\otimes \shift \varphi^3\otimes \varphi^4\otimes \shift \varphi^5\otimes \hat{\omega}}_{15}\\
&+\underbrace{1\otimes \shift \varphi^1\otimes \varphi^2\varphi^3\varphi^4\otimes \shift \varphi^5\otimes \hat{\omega}}_{13}
-\underbrace{\varphi^1\varphi^{21}\dots \varphi^{24}\otimes \shift \varphi^{25}\otimes \shift \varphi^3\otimes \varphi^4\varphi^5\hat{\omega}}_{91}\\
&-\underbrace{\varphi^1\varphi^2\otimes \shift \varphi^3\otimes \shift \varphi^{41}\otimes \varphi^{42}\dots \varphi^{45}\varphi^5\hat{\omega}}_{57}
-\underbrace{1\otimes \shift \varphi^1\otimes \shift \varphi^{21}\otimes \varphi^{22}\dots \varphi^{25}\varphi^3\varphi^4\varphi^5\hat{\omega}}_{55}\\
&-\underbrace{\varphi^1\varphi^2\varphi^3\varphi^{41}\dots \varphi^{44}\otimes \shift \varphi^{45}\otimes \shift \varphi^5\otimes \hat{\omega}}_{22}
\+\underbrace{\varphi^1\otimes \shift \varphi^{21}\otimes \varphi^{22}\otimes \shift \varphi^{23}\otimes \varphi^{24}\varphi^{25}\varphi^3\varphi^4\varphi^5}_{42}\\
&+\underbrace{\varphi^1\varphi^{21}\varphi^{22}\otimes \shift \varphi^{23}\otimes \varphi^{24}\otimes \shift \varphi^{25}\otimes \varphi^3\varphi^4\varphi^5}_{39}
-\underbrace{\varphi^1\otimes \shift \varphi^{21}\otimes \varphi^{22}\dots \varphi^{24}\otimes \shift \varphi^{25}\otimes \varphi^3\varphi^4\varphi^5}_{86}\\
&+\underbrace{\varphi^1\varphi^{21}\varphi^{221}\dots \varphi^{224}\otimes \shift \varphi^{225}\otimes \shift \varphi^{23}\otimes\varphi^{24}\varphi^{25}\varphi^3\varphi^4\varphi^5}_{g}\\
&+\underbrace{\varphi^1\varphi^{21}\varphi^{22}\otimes \shift \varphi^{23}\otimes \shift \varphi^{241}\otimes \varphi^{242}\dots \varphi^{245}\varphi^{25}\varphi^3\varphi^4\varphi^5}_g\\
&+\underbrace{\varphi^1\otimes \shift \varphi^{21}\otimes \shift \varphi^{221}\otimes \varphi^{222}\dots \varphi^{225}\varphi^{23}\dots \varphi^{25}\varphi^3\varphi^4\varphi^5}_{34}\\
&+\underbrace{\varphi^1\varphi^{21}\varphi^{22} \varphi^{23}\varphi^{241}\dots \varphi^{244}\otimes \shift \varphi^{245}\otimes \shift \varphi^{25}\otimes \varphi^3\varphi^4\varphi^5}_{49}
\end{align*}
and
\begin{align*}
(1\otimes \mu_4)\mu_2(\varphi)=\,&(1\otimes \mu_4)(\hat{\omega}\varphi+\varphi\hat{\omega}+(\shift \otimes \shift)\partial(\shift^{-1}\varphi))\\
=\,&\underbrace{\hat{\omega}\otimes \shift \varphi^1\otimes \varphi^2\otimes \shift \varphi^3\otimes \varphi^4\varphi^5}_{14}
+\underbrace{\hat{\omega}\varphi^1\varphi^2\otimes \shift \varphi^3\otimes \varphi^4\otimes \shift \varphi^5\otimes 1}_{94}\\
&-\underbrace{\hat{\omega}\otimes \shift \varphi^1\otimes \varphi^2\varphi^3\varphi^4\otimes \shift \varphi^5\otimes 1}_{16}
+\underbrace{\hat{\omega}\varphi^1\varphi^{21}\dots \varphi^{24}\otimes \shift \varphi^{25}\otimes \shift \varphi^3\otimes \varphi^4\varphi^5}_{80}\\
&+\underbrace{\hat{\omega}\varphi^1\varphi^2\otimes \shift \varphi^3\otimes \shift \varphi^{41}\otimes \varphi^{42}\dots \varphi^{45}\varphi^5}_{74}
+\underbrace{\hat{\omega}\otimes \shift \varphi^1\otimes \shift \varphi^{21}\otimes \varphi^{22}\dots \varphi^{25}\varphi^3\varphi^4\varphi^5}_{21}\\
&+\underbrace{\hat{\omega}\varphi^1\varphi^2\varphi^3\varphi^{41}\dots \varphi^{44}\otimes \shift \varphi^{45}\otimes \shift \varphi^5\otimes 1}_{78}\+0\\
&\-\underbrace{\varphi^1\varphi^2\varphi^3\otimes \shift \varphi^{41}\otimes \varphi^{42}\otimes \shift \varphi^{43}\otimes \varphi^{44}\varphi^{45}\varphi^5}_{26}
-\underbrace{\varphi^1\varphi^2\varphi^3\varphi^{41}\varphi^{42}\otimes \shift \varphi^{43}\otimes \varphi^{44}\otimes \shift \varphi^{45}\otimes \varphi^5}_{29}\\
&+\underbrace{\varphi^1\varphi^2\varphi^3\otimes \shift \varphi^{41}\otimes \varphi^{42}\dots \varphi^{44}\otimes \shift \varphi^{45}\otimes \varphi^5}_{63}\\
&-\underbrace{\varphi^1\varphi^2\varphi^3\varphi^{41}\varphi^{421}\dots \varphi^{424}\otimes \shift \varphi^{425}\otimes \shift \varphi^{43}\otimes \varphi^{44}\otimes \varphi^{45}\varphi^5}_g\\
&-\underbrace{\varphi^1\varphi^2\varphi^3\varphi^{41}\varphi^{42}\otimes \shift \varphi^{43}\otimes \shift \varphi^{441}\otimes \varphi^{442}\dots \varphi^{445}\varphi^{45}\varphi^5}_g\\
&-\underbrace{\varphi^1\varphi^2\varphi^3\otimes \shift \varphi^{41}\otimes \shift \varphi^{421}\otimes \varphi^{422}\dots \varphi^{425}\varphi^{43}\varphi^{44}\varphi^{45}\varphi^5}_{36}\\
&-\underbrace{\varphi^1\varphi^2\varphi^3\varphi^{41}\varphi^{42}\varphi^{43}\varphi^{441}\dots \varphi^{444}\otimes \shift \varphi^{445}\otimes \shift \varphi^{45}\otimes \varphi^5}_{47}
\end{align*}
and
\begin{align*}
(\mu_3\otimes 1^{\otimes 2})\mu_3(\varphi)=\,&(\mu_3\otimes 1^{\otimes 2})(1\otimes \shift\varphi^1\otimes\varphi^2\varphi^3\varphi^4\varphi^5-\varphi^1\varphi^2\otimes \shift\varphi^3\otimes \varphi^4\varphi^5\\&+\varphi^1\varphi^2\varphi^3\varphi^4\otimes \shift\varphi^5\otimes 1)\\
=\,&\underbrace{1\otimes \shift \varphi^{11}\otimes \varphi^{12}\otimes \shift \varphi^{13}\otimes \varphi^2\dots \varphi^5}_y
+\underbrace{1\otimes \shift \varphi^{11}\otimes \shift \varphi^{121}\otimes \varphi^{122}\dots \varphi^{125}\varphi^{13}\varphi^2\dots \varphi^5}_g\\
&+\underbrace{\varphi^{11}\varphi^{121}\dots \varphi^{124}\otimes \shift \varphi^{125}\otimes \shift \varphi^{13}\otimes \varphi^2\dots \varphi^5}_g
\-\underbrace{\varphi^1\otimes \shift \varphi^{21}\otimes \varphi^{22}\dots \varphi^{25}\otimes \shift \varphi^3\otimes \varphi^4\varphi^5}_{85}\\
&+\underbrace{\varphi^1\varphi^{21}\varphi^{22}\otimes \shift \varphi^{23}\otimes \varphi^{24}\varphi^{25}\otimes \shift \varphi^3\otimes \varphi^4\varphi^5}_{40}
-\underbrace{\varphi^1\varphi^{21}\dots \varphi^{24}\otimes \shift \varphi^{25}\otimes 1\otimes \shift \varphi^3\otimes \varphi^4\varphi^5}_{20}\\
&\+\underbrace{\varphi^1\otimes \shift \varphi^{21}\otimes \varphi^{22}\dots \varphi^{25}\varphi^3\varphi^4\otimes \shift \varphi^5\otimes 1}_{56}
-\underbrace{\varphi^1\varphi^{21}\varphi^{22}\otimes \shift \varphi^{23}\otimes \varphi^{24}\varphi^{25}\varphi^3\varphi^4\otimes \shift \varphi^5\otimes 1}_y\\
&+\underbrace{\varphi^1\varphi^{21}\dots \varphi^{24}\otimes \shift \varphi^{25}\otimes \varphi^3\varphi^4\otimes \shift \varphi^5\otimes 1}_{52}
\end{align*}
and
\begin{align*}
(1\otimes \mu_3\otimes 1)\mu_3(\varphi)=\,&(1\otimes \mu_3\otimes 1)(1\otimes \shift\varphi^1\otimes\varphi^2\varphi^3\varphi^4\varphi^5-\varphi^1\varphi^2\otimes \shift\varphi^3\otimes \varphi^4\varphi^5\\&+\varphi^1\varphi^2\varphi^3\varphi^4\otimes \shift\varphi^5\otimes 1)\\
=\,&-\underbrace{1\otimes \shift \varphi^1\otimes 1\otimes \shift \varphi^{21}\otimes \varphi^{22}\dots \varphi^{25}\varphi^3\varphi^4\varphi^5}_{18}
+\underbrace{1\otimes \shift \varphi^1\otimes \varphi^{21}\varphi^{22}\otimes \shift \varphi^{23}\otimes \varphi^{24}\varphi^{25}\varphi^3\varphi^4\varphi^5}_{44}\\
&-\underbrace{1\otimes \shift \varphi^1\otimes \varphi^{21}\dots \varphi^{24}\otimes \shift \varphi^{25}\otimes \varphi^3\varphi^4\varphi^5}_{82}
\-\underbrace{\varphi^1\varphi^2\otimes \shift \varphi^{31}\otimes \varphi^{32}\otimes \shift \varphi^{33}\otimes \varphi^4\varphi^5}_y\\
&-\underbrace{\varphi^1\varphi^2\otimes \shift \varphi^{31}\otimes \shift \varphi^{321}\otimes \varphi^{322}\dots \varphi^{325}\varphi^{33}\varphi^4\varphi^5}_g\\
&
-\underbrace{\varphi^1\varphi^2\varphi^{31}\varphi^{321}\dots \varphi^{324}\otimes \shift \varphi^{325}\otimes \shift \varphi^{33}\otimes \varphi^4\varphi^5}_g\\
&\+\underbrace{\varphi^1\varphi^2\varphi^3\otimes \shift \varphi^{41}\otimes \varphi^{42}\dots \varphi^{45}\otimes \shift \varphi^5\otimes 1}_{59}
-\underbrace{\varphi^1\varphi^2\varphi^3\varphi^{41}\varphi^{42}\otimes \shift \varphi^{43}\otimes \varphi^{44}\varphi^{45}\otimes \shift \varphi^{5}\otimes 1}_{31}\\
&+\underbrace{\varphi^1\varphi^2\varphi^3\varphi^{41}\dots \varphi^{44}\otimes \shift \varphi^{45}\otimes 1\otimes \shift \varphi^5\otimes 1}_{17} 
\end{align*}
and
\begin{align*}
(1^{\otimes 2}\otimes \mu_3)\mu_3(\varphi)=\,&(1^{\otimes 2}\otimes \mu_3)(1\otimes \shift\varphi^1\otimes\varphi^2\varphi^3\varphi^4\varphi^5-\varphi^1\varphi^2\otimes \shift\varphi^3\otimes \varphi^4\varphi^5\\&+\varphi^1\varphi^2\varphi^3\varphi^4\otimes \shift\varphi^5\otimes 1)\\
=\,&-\underbrace{1\otimes \shift \varphi^1\otimes \varphi^2\varphi^3\otimes \shift \varphi^{41}\otimes \varphi^{42}\dots \varphi^{45}\varphi^5}_{75}
+\underbrace{1\otimes \shift \varphi^1\otimes \varphi^2\varphi^3\varphi^{41}\varphi^{42}\otimes \shift \varphi^{43}\otimes \varphi^{44}\varphi^{45}\varphi^5}_y\\
&-\underbrace{1\otimes \shift \varphi^1\otimes \varphi^2\varphi^3\varphi^{41}\dots \varphi^{44}\otimes \shift \varphi^{45}\otimes \varphi^5}_{79}
\+\underbrace{\varphi^1\varphi^2\otimes \shift \varphi^3\otimes 1\otimes \shift \varphi^{41}\otimes \varphi^{42}\dots \varphi^{45}\varphi^5}_{19}\\
&-\underbrace{\varphi^1\varphi^2\otimes \shift \varphi^3\otimes \varphi^{41}\varphi^{42}\otimes \shift \varphi^{43}\otimes \varphi^{44}\varphi^{45}\varphi^5}_{27}
+\underbrace{\varphi^1\varphi^2\otimes \shift \varphi^3\otimes \varphi^{41}\dots \varphi^{44}\otimes \shift \varphi^{45}\otimes \varphi^5}_{62}\\
&\+\underbrace{\varphi^1\dots \varphi^4\otimes \shift \varphi^{51}\otimes \varphi^{52}\otimes \shift \varphi^{53}\otimes 1}_y
+\underbrace{\varphi^1\dots \varphi^4\otimes \shift \varphi^{51}\otimes \shift \varphi^{521}\otimes \varphi^{522}\dots \varphi^{525}\varphi^{53}}_g\\
&+\underbrace{\varphi^1\varphi^2\varphi^3\varphi^4\varphi^{51}\varphi^{521}\dots \varphi^{524}\otimes \shift \varphi^{525}\otimes \shift \varphi^{53}\otimes 1}_g
\end{align*}
For simplicity, we use $\mu_1^{\otimes}:=(\mu_1\otimes 1^{\otimes 4}+1\otimes \mu_1\otimes 1^{\otimes 3}+1^{\otimes 2}\otimes \mu_1\otimes 1^{\otimes 2}+1^{\otimes 3}\otimes \mu_1\otimes 1+1^{\otimes 4}\otimes \mu_1)$. Applying this to each summand of our asserted $\mu_5^0(\varphi)$ individually, we obtain

\begin{align*}
\mu_1^{\otimes}(-1\otimes \shift \varphi^1\otimes \varphi^2\otimes \shift \varphi^3\otimes \varphi^4\otimes \shift \varphi^5\otimes 1)
&=\,-\underbrace{\hat{\omega}\varphi^1\varphi^2\otimes \shift \varphi^3\otimes \varphi^4\otimes \shift \varphi^5\otimes 1}_{94}
+\underbrace{\varphi^1\hat{\omega}\varphi^2\otimes \shift \varphi^3\otimes \varphi^4\otimes \shift \varphi^5\otimes 1}_{68}\\
&+\underbrace{\shift \partial(\varphi^1)\varphi^2\otimes \shift \varphi^3\otimes \varphi^4\otimes \shift \varphi^5\otimes 1}_y\boldsymbol{\uwave{+}}0
\boldsymbol{\uwave{+}}\underbrace{1\otimes \shift \varphi^1\otimes \varphi^2\hat{\omega}\varphi^3\varphi^4\otimes \shift \varphi^5\otimes 1}_{69}\\
&-\underbrace{1\otimes \shift \varphi^1\otimes \varphi^2\varphi^3\hat{\omega}\varphi^4\otimes \shift \varphi^5\otimes 1}_{71}
-\underbrace{1\otimes \shift \varphi^1\otimes \varphi^2\shift \partial(\varphi^3)\varphi^4\otimes \shift \varphi^5\otimes 1}_y\boldsymbol{\uwave{+}}0\\
&\boldsymbol{\uwave{-}}\underbrace{1\otimes \shift \varphi^1\otimes \varphi^2\otimes \shift \varphi^3\otimes \varphi^4\hat{\omega}\varphi^5}_{70}
+\underbrace{1\otimes \shift \varphi^1\otimes \varphi^2\otimes \shift \varphi^3\otimes \varphi^4\varphi^5\hat{\omega}}_{72}\\
&-\underbrace{1\otimes \shift \varphi^1\otimes \varphi^2\otimes \shift \varphi^3\otimes \varphi^4\shift \partial(\varphi^5)}_y
\end{align*}
and
\begin{align*}
&\mu_1^{\otimes}(-1\otimes \shift \varphi^1\otimes \varphi^2\varphi^3\varphi^{41}\dots \varphi^{44}\otimes \shift \varphi^{45}\otimes \shift \varphi^5\otimes 1)\\
=\,&-\underbrace{\hat{\omega}\varphi^1\varphi^2\varphi^3\varphi^{41}\dots \varphi^{44}\otimes \shift \varphi^{45}\otimes \shift \varphi^5\otimes 1}_{78}
+\underbrace{\varphi^1\hat{\omega}\varphi^2\varphi^3\varphi^{41}\dots \varphi^{44}\otimes \shift \varphi^{45}\otimes \shift \varphi^5\otimes 1}_{76}\\
&+\underbrace{\shift \partial(\varphi^1)\varphi^2\varphi^3\varphi^{41}\dots \varphi^{44}\otimes \shift \varphi^{45}\otimes \shift \varphi^5\otimes 1}_g\boldsymbol{\uwave{+}}0\boldsymbol{\uwave{+}}0\boldsymbol{\uwave{+}}0
\boldsymbol{\uwave{-}}\underbrace{1\otimes \shift \varphi^1\otimes \varphi^2\varphi^3(\shift \otimes \shift)\partial(\shift^{-1} \varphi^4)\otimes \shift \varphi^5\otimes 1}_{77}\\
&-\underbrace{1\otimes \shift \varphi^1\otimes \varphi^2\varphi^3\varphi^{41}\dots \varphi^{44}\otimes \shift(\varphi^{45}\varphi^5)\otimes 1}_y
+\underbrace{1\otimes \shift \varphi^1\otimes \varphi^2\varphi^3\varphi^{41}\dots \varphi^{44}\otimes \shift \varphi^{45}\otimes \varphi^5}_{79}
\end{align*}
and
\begin{align*}
&\mu_1^{\otimes}(-1\otimes \shift \varphi^1\otimes \shift \varphi^{21}\otimes \varphi^{22}\dots \varphi^{25}\varphi^3\varphi^4\otimes \shift \varphi^5\otimes 1)\\
=\,&-\underbrace{\varphi^1\otimes \shift \varphi^{21}\otimes \varphi^{22}\dots \varphi^{25}\varphi^3\varphi^4\otimes \shift \varphi^5\otimes 1}_{56}
+\underbrace{1\otimes \shift(\varphi^1\varphi^{21})\otimes \varphi^{22}\dots \varphi^{25}\varphi^3\varphi^4\otimes \shift \varphi^5\otimes 1}_y\\
&+\underbrace{1\otimes \shift \varphi^1\otimes (\shift \otimes \shift)\partial(\shift^{-1} \varphi^2)\varphi^3\varphi^4\otimes \shift \varphi^5\otimes 1}_{54}\boldsymbol{\uwave{+}}0\boldsymbol{\uwave{+}}0\boldsymbol{\uwave{+}}0
\boldsymbol{\uwave{-}}\underbrace{1\otimes \shift \varphi^1\otimes \shift \varphi^{21}\otimes \varphi^{22}\dots \varphi^{25}\varphi^3\varphi^4\hat{\omega}\varphi^5}_{53}\\
&+\underbrace{1\otimes \shift \varphi^1\otimes \shift \varphi^{21}\otimes \varphi^{22}\dots \varphi^{25}\varphi^3\varphi^4\varphi^5\hat{\omega}}_{55}
+\underbrace{1\otimes \shift \varphi^1\otimes \shift \varphi^{21}\otimes \varphi^{22}\dots \varphi^{25}\varphi^3\varphi^4\shift \partial(\varphi^5)}_g
\end{align*}
and
\begin{align*}
&\mu_1^{\otimes}(-1\otimes \shift \varphi^1\otimes \varphi^{21}\dots \varphi^{24}\otimes \shift \varphi^{25}\otimes \shift \varphi^3\otimes \varphi^4\varphi^5)\\
=\,&- \underbrace{\hat{\omega}\varphi^1\varphi^{21}\dots \varphi^{24}\otimes \shift \varphi^{25}\otimes \shift \varphi^3\otimes \varphi^4\varphi^5}_{80}
+\underbrace{\varphi^1\hat{\omega}\varphi^{21}\dots \varphi^{24}\otimes \shift \varphi^{25}\otimes \shift \varphi^3\otimes \varphi^4\varphi^5}_{83}\\
&+\underbrace{\shift \partial(\varphi^1)\varphi^{21}\dots \varphi^{24}\otimes \shift \varphi^{25}\otimes \shift \varphi^3\otimes \varphi^4\varphi^5}_g\boldsymbol{\uwave{+}}0\boldsymbol{\uwave{+}}0
\boldsymbol{\uwave{-}}\underbrace{1\otimes \shift \varphi^1\otimes (\shift \otimes \shift)\partial(\shift^{-1} \varphi^2)\otimes \shift \varphi^3\otimes \varphi^4\varphi^5}_{81}\\
&-\underbrace{1\otimes \shift \varphi^1\otimes \varphi^{21}\dots \varphi^{24}\otimes \shift(\varphi^{25}\varphi^3)\otimes \varphi^4\varphi^5}_{66}
+\underbrace{1\otimes \shift \varphi^1\otimes \varphi^{21}\dots \varphi^{24}\otimes \shift \varphi^{25}\otimes \varphi^3\varphi^4\varphi^5}_{82}\boldsymbol{\uwave{+}}0
\end{align*}
and
\begin{align*}
&\mu_1^{\otimes}(-\varphi^1\varphi^2\otimes \shift \varphi^3\otimes \shift \varphi^{41}\otimes \varphi^{42}\dots \varphi^{45}\otimes \shift \varphi^5\otimes 1)\\
=\,&0\boldsymbol{\uwave{-}}\underbrace{\varphi^1\varphi^2\varphi^3\otimes \shift \varphi^{41}\otimes \varphi^{42}\dots \varphi^{45}\otimes \shift \varphi^5\otimes 1}_{59}
+\underbrace{\varphi^1\varphi^2\otimes \shift(\varphi^3\varphi^{41})\otimes \varphi^{42}\dots \varphi^{45}\otimes \shift \varphi^5\otimes 1}_{89}\\
&+\underbrace{\varphi^1\varphi^2\otimes \shift \varphi^3\otimes (\shift \otimes \shift)\partial(\shift^{-1} \varphi^4)\otimes \shift \varphi^5\otimes 1}_{58}\boldsymbol{\uwave{+}}0\boldsymbol{\uwave{+}}0
\boldsymbol{\uwave{-}}\underbrace{\varphi^1\varphi^2\otimes \shift \varphi^3\otimes \shift \varphi^{41}\otimes \varphi^{42}\dots \varphi^{45}\hat{\omega}\varphi^5}_{60}\\
&+\underbrace{\varphi^1\varphi^2\otimes \shift \varphi^3\otimes \shift \varphi^{41}\otimes \varphi^{42}\dots \varphi^{45}\varphi^5\hat{\omega}}_{57}
+\underbrace{\varphi^1\varphi^2\otimes \shift \varphi^3\otimes \shift \varphi^{41}\otimes \varphi^{42}\dots \varphi^{45}\shift \partial(\varphi^5)}_g
\end{align*}
and
\label{ggocc1}
\begin{align*}
&\mu_1^{\otimes}(-\varphi^1\varphi^{21}\dots \varphi^{24}\otimes \shift \varphi^{25}\otimes \shift \varphi^3\otimes \varphi^4\otimes \shift \varphi^5\otimes 1)\\
=\,&0\boldsymbol{\uwave{+}}0\boldsymbol{\uwave{+}}\underbrace{\varphi^1(\shift \otimes \shift)\partial(\shift ^{-1}\varphi^2)\otimes \shift \varphi^3\otimes \varphi^4\otimes \shift \varphi^5\otimes 1}_{73}
+\underbrace{\varphi^1\varphi^{21}\dots \varphi^{24}\otimes \shift(\varphi^{25}\varphi^3)\otimes \varphi^4\otimes \shift \varphi^5\otimes 1}_g\\
&-\underbrace{\varphi^1\varphi^{21}\dots \varphi^{24}\otimes \shift \varphi^{25}\otimes \varphi^3\varphi^4\otimes \shift \varphi^5\otimes 1}_{52}\boldsymbol{\uwave{+}}0
\boldsymbol{\uwave{-}}\underbrace{\varphi^1\varphi^{21}\dots \varphi^{24}\otimes \shift \varphi^{25}\otimes \shift \varphi^3\otimes \varphi^4\hat{\omega}\varphi^5}_{51}\\
&+\underbrace{\varphi^1\varphi^{21}\dots \varphi^{24}\otimes \shift \varphi^{25}\otimes \shift \varphi^3\otimes \varphi^4\varphi^5\hat{\omega}}_{91}
+\underbrace{\varphi^1\varphi^{21}\dots \varphi^{24}\otimes \shift \varphi^{25}\otimes \shift \varphi^3\otimes \varphi^4\shift \partial(\varphi^5)}_{gg}
\end{align*}
and
\begin{align*}
&\mu_1^{\otimes}(-1\otimes \shift \varphi^1\otimes \varphi^2\otimes \shift \varphi^3\otimes \shift \varphi^{41}\otimes \varphi^{42}\dots \varphi^{45}\varphi^5)\\
=\,&-\underbrace{\hat{\omega}\varphi^1\varphi^2\otimes \shift \varphi^3\otimes \shift \varphi^{41}\otimes \varphi^{42}\dots \varphi^{45}\varphi^5}_{74}
+\underbrace{\varphi^1\hat{\omega}\varphi^2\otimes \shift \varphi^3\otimes \shift \varphi^{41}\otimes \varphi^{42}\dots \varphi^{45}\varphi^5}_{92}\\
&+\underbrace{\shift \partial(\varphi^1)\varphi^2\otimes \shift \varphi^3\otimes \shift \varphi^{41}\otimes \varphi^{42}\dots \varphi^{45}\varphi^5}_{g}\boldsymbol{\uwave{+}}0
\boldsymbol{\uwave{+}}\underbrace{1\otimes \shift \varphi^1\otimes \varphi^2\varphi^3\otimes \shift \varphi^{41}\otimes \varphi^{42}\dots \varphi^{45}\varphi^5}_{75}\\
&-\underbrace{1\otimes \shift \varphi^1\otimes \varphi^2\otimes \shift(\varphi^3\varphi^{41})\otimes \varphi^{42}\dots \varphi^{45}\varphi^5}_{g}
-\underbrace{1\otimes \shift \varphi^1\otimes \varphi^2\otimes \shift \varphi^3\otimes (\shift \otimes \shift)\partial(\shift^{-1}\varphi^4)\varphi^5}_{93}\boldsymbol{\uwave{+}}0\boldsymbol{\uwave{+}}0
\end{align*}
and
\begin{align*}
&\mu_1^{\otimes}(-\varphi^1\otimes \shift \varphi^{21}\otimes \varphi^{22}\dots \varphi^{24}\otimes \shift \varphi^{25}\otimes \shift \varphi^3\otimes \varphi^4\varphi^5)\\
=\,&-\underbrace{\varphi^1\hat{\omega}\varphi^{21}\dots \varphi^{24}\otimes \shift \varphi^{25}\otimes \shift \varphi^3\otimes \varphi^4\varphi^5}_{83}
+\underbrace{\varphi^1\varphi^{21}\hat{\omega}\varphi^{22}\varphi^{23} \varphi^{24}\otimes \shift \varphi^{25}\otimes \shift \varphi^3\otimes \varphi^4\varphi^5}_{84}\\
&+\underbrace{\varphi^1\shift \partial(\varphi^{21})\varphi^{22}\dots \varphi^{24}\otimes \shift \varphi^{25}\otimes \shift \varphi^3\otimes \varphi^4\varphi^5}_y\boldsymbol{\uwave{+}}0\boldsymbol{\uwave{+}}0
\boldsymbol{\uwave{+}}\underbrace{\varphi^1\otimes \shift \varphi^{21}\otimes \varphi^{22}\dots \varphi^{25}\otimes \shift \varphi^3\otimes \varphi^4\varphi^5}_{85}\\
&-\underbrace{\varphi^1\otimes \shift \varphi^{21}\otimes \varphi^{22}\dots \varphi^{24}\otimes \shift(\varphi^{25}\varphi^3)\otimes \varphi^4\varphi^5}_{41}
+\underbrace{\varphi^1\otimes \shift \varphi^{21}\otimes \varphi^{22}\dots \varphi^{24}\otimes \shift \varphi^{25}\otimes \varphi^{3}\varphi^4\varphi^5}_{86}\boldsymbol{\uwave{+}}0
\end{align*}
and
\begin{align*}
&\mu_1^{\otimes}(-\varphi^1\varphi^2\otimes \shift \varphi^3\otimes \shift \varphi^{41}\otimes \varphi^{42}\dots \varphi^{44}\otimes \shift \varphi^{45}\otimes \varphi^5)\\
=\,&0\boldsymbol{\uwave{-}}\underbrace{\varphi^1\varphi^2\varphi^3\otimes \shift \varphi^{41}\otimes \varphi^{42}\dots \varphi^{44}\otimes \shift \varphi^{45}\otimes \varphi^5}_{63}
+\underbrace{\varphi^1\varphi^2\otimes \shift(\varphi^3\varphi^{41})\otimes \varphi^{42}\dots \varphi^{44}\otimes \shift \varphi^{45}\otimes \varphi^5}_{28}\\
&-\underbrace{\varphi^1\varphi^2\otimes \shift \varphi^3\otimes \varphi^{41}\dots \varphi^{44}\otimes \shift \varphi^{45}\otimes \varphi^5}_{62}\boldsymbol{\uwave{+}}0\boldsymbol{\uwave{+}}0
\boldsymbol{\uwave{-}}\underbrace{\varphi^1\varphi^2\otimes \shift \varphi^3\otimes \shift \varphi^{41}\otimes \varphi^{42}\dots \varphi^{44}\hat{\omega}\varphi^{45}\varphi^5}_{61}\\
&+\underbrace{\varphi^1\varphi^2\otimes \shift \varphi^3\otimes \shift \varphi^{41}\otimes \varphi^{42}\dots \varphi^{45}\hat{\omega}\varphi^5}_{60}
+\underbrace{\varphi^1\varphi^2\otimes \shift \varphi^3\otimes \shift \varphi^{41}\otimes \varphi^{42}\dots \varphi^{44}\shift \partial(\varphi^{45})\varphi^5}_y
\end{align*}
and
\begin{align*}
&\mu_1^{\otimes}(\varphi^1\varphi^2\otimes \shift(\varphi^3\varphi^{41})\otimes \varphi^{42}\dots \varphi^{44}\otimes \shift \varphi^{45}\otimes \shift \varphi^5\otimes 1)\\
=\,&0\boldsymbol{\uwave{+}}\underbrace{\varphi^1\varphi^2\hat{\omega}\varphi^3\varphi^{41}\dots \varphi^{44}\otimes \shift \varphi^{45}\otimes \shift \varphi^5\otimes 1}_{87}
-\underbrace{\varphi^1\varphi^2\varphi^3\varphi^{41}\hat{\omega}\varphi^{42}\dots \varphi^{44}\otimes \shift \varphi^{45}\otimes \shift \varphi^5\otimes 1}_{88}\\
&-\underbrace{\varphi^1\varphi^2\shift \partial(\varphi^3)\varphi^{41}\dots \varphi^{44}\otimes \shift \varphi^{45}\otimes \shift \varphi^5\otimes 1}_g
-\underbrace{\varphi^1\varphi^2\varphi^3\shift \partial(\varphi^{41})\varphi^{42}\dots \varphi^{44}\otimes \shift \varphi^{45}\otimes \shift \varphi^5\otimes 1}_y\boldsymbol{\uwave{+}}0\boldsymbol{\uwave{+}}0\\
&\boldsymbol{\uwave{-}}\underbrace{\varphi^1\varphi^2\otimes \shift(\varphi^3\varphi^{41})\otimes \varphi^{42}\dots \varphi^{45}\otimes \shift \varphi^5\otimes 1}_{89}+\underbrace{\varphi^1\varphi^2\otimes \shift(\varphi^3\varphi^{41})\otimes \varphi^{42}\dots \varphi^{44}\otimes \shift(\varphi^{45}\varphi^5)\otimes 1}_y\\
&-\underbrace{\varphi^1\varphi^2\otimes \shift(\varphi^3\varphi^{41})\otimes \varphi^{42}\dots \varphi^{44}\otimes \shift \varphi^{45}\otimes \varphi^5}_{28}
\end{align*}
and
\begin{align*}
&\mu_1^{\otimes}(1\otimes \shift \varphi^1\otimes \shift \varphi^{21}\otimes \varphi^{22}\dots \varphi^{24}\otimes \shift(\varphi^{25}\varphi^3)\otimes \varphi^4\varphi^5)\\
=\,&\underbrace{\varphi^1\otimes \shift \varphi^{21}\otimes \varphi^{22}\dots \varphi^{24}\otimes \shift(\varphi^{25}\varphi^3)\otimes \varphi^4\varphi^5}_{41}
-\underbrace{1\otimes \shift(\varphi^1\varphi^{21})\otimes \varphi^{22}\dots \varphi^{24}\otimes \shift(\varphi^{25}\varphi^3)\otimes \varphi^4\varphi^5}_y\\
&+\underbrace{1\otimes\shift \varphi^1\otimes \varphi^{21}\dots \varphi^{24}\otimes \shift(\varphi^{25}\varphi^3)\otimes \varphi^4\varphi^5}_{66}\boldsymbol{\uwave{+}}0\boldsymbol{\uwave{+}}0
\boldsymbol{\uwave{+}}\underbrace{1\otimes \shift \varphi^1\otimes\shift  \varphi^{21}\otimes \varphi^{22}\dots \varphi^{24}\hat{\omega}\varphi^{25}\varphi^3\varphi^4\varphi^5}_{65}\\
&-\underbrace{1\otimes \shift \varphi^1\otimes \shift \varphi^{21}\otimes \varphi^{22}\dots \varphi^{25}\varphi^3\hat{\omega}\varphi^4\varphi^5}_{64}
-\underbrace{1\otimes \shift \varphi^1\otimes \shift \varphi^{21}\otimes \varphi^{22}\dots \varphi^{24}\shift \partial(\varphi^{25})\varphi^3\varphi^4\varphi^5}_y\\
&-\underbrace{1\otimes \shift \varphi^1\otimes \shift \varphi^{21}\otimes \varphi^{22}\dots \varphi^{25}\shift \partial(\varphi^3)\varphi^4\varphi^5}_g\boldsymbol{\uwave{+}}0
\end{align*}
and

\begin{align*}
&\mu_1^{\otimes}(\varphi^1\varphi^{21}\varphi^{22}\otimes \shift \varphi^{23}\otimes \varphi^{24}\otimes \shift \varphi^{25}\otimes \shift \varphi^3\otimes \varphi^4\varphi^5)\\
=\,&0\boldsymbol{\uwave{+}}\underbrace{\varphi^1\varphi^{21}\varphi^{22}\hat{\omega}\varphi^{23}\varphi^{24}\otimes \shift \varphi^{25}\otimes \shift \varphi^3\otimes \varphi^4\varphi^5}_{90}
-\underbrace{\varphi^1\varphi^{21}\varphi^{22}\varphi^{23}\hat{\omega}\varphi^{24}\otimes \shift \varphi^{25}\otimes \shift \varphi^3\otimes \varphi^4\varphi^5}_{38}\\
&-\underbrace{\varphi^1\varphi^{21}\varphi^{22}\shift \partial(\varphi^{23})\varphi^{24}\otimes\shift  \varphi^{25}\otimes \shift \varphi^3\otimes \varphi^4\varphi^5}_y\boldsymbol{\uwave{+}}0\boldsymbol{\uwave{-}}
\underbrace{\varphi^1\varphi^{21}\varphi^{22}\otimes \shift \varphi^{23}\otimes \varphi^{24}\varphi^{25}\otimes \shift \varphi^3\otimes \varphi^4\varphi^5}_{40}\\
&+\underbrace{\varphi^1\varphi^{21}\varphi^{22}\otimes \shift \varphi^{23}\otimes \varphi^{24}\otimes \shift(\varphi^{25}\varphi^3)\otimes \varphi^4\varphi^5}_y
-\underbrace{\varphi^1\varphi^{21}\varphi^{22}\otimes \shift \varphi^{23}\otimes \varphi^{24}\otimes \shift \varphi^{25}\otimes \varphi^3\varphi^4\varphi^5}_{39}\boldsymbol{\uwave{+}}0
\end{align*}
and
\begin{align*}
&\mu_1^{\otimes}(\varphi^1\varphi^2\otimes \shift \varphi^3\otimes \shift \varphi^{41}\otimes \varphi^{42}\otimes \shift \varphi^{43}\otimes \varphi^{44}\varphi^{45}\varphi^5)\\
=\,&0\boldsymbol{\uwave{+}}\underbrace{\varphi^1\varphi^2\varphi^3\otimes \shift \varphi^{41}\otimes \varphi^{42}\otimes \shift \varphi^{43}\otimes \varphi^{44}\varphi^{45}\varphi^5}_{26}
-\underbrace{\varphi^1\varphi^2\otimes \shift(\varphi^3\varphi^{41})\otimes \varphi^{42}\otimes \shift \varphi^{43}\otimes \varphi^{44}\varphi^{45}\varphi^5}_y\\
&+\underbrace{\varphi^1\varphi^2\otimes \shift \varphi^3\otimes \varphi^{41}\varphi^{42}\otimes\shift  \varphi^{43}\otimes \varphi^{44}\varphi^{45}\varphi^5}_{27}\boldsymbol{\uwave{+}}0
\boldsymbol{\uwave{+}}\underbrace{\varphi^1\varphi^2\otimes \shift \varphi^3\otimes \shift \varphi^{41}\otimes \varphi^{42}\hat{\omega}\varphi^{43}\varphi^{44}\varphi^{45}\varphi^5}_{25}\\
&-\underbrace{\varphi^1\varphi^2\otimes \shift \varphi^3\otimes \shift \varphi^{41}\otimes \varphi^{42}\varphi^{43}\hat{\omega}\varphi^{44}\varphi^{45}\varphi^5}_{67}
-\underbrace{\varphi^1\varphi^2\otimes \shift \varphi^3\otimes \shift \varphi^{41}\otimes \varphi^{42}\shift \partial(\varphi^{43})\varphi^{44}\varphi^{45}\varphi^5}_y\boldsymbol{\uwave{+}}0
\end{align*}
and
\begin{align*}
&\mu_1^{\otimes}(-1\otimes \shift \varphi^1\otimes \shift \varphi^{21}\otimes \varphi^{22}\otimes \shift \varphi^{23}\otimes \varphi^{24}\varphi^{25}\varphi^3\varphi^4\varphi^5)\\
=\,&-\underbrace{\varphi^1\otimes \shift \varphi^{21}\otimes \varphi^{22}\otimes \shift \varphi^{23}\otimes \varphi^{24}\varphi^{25}\varphi^3\varphi^4\varphi^5}_{42}
+\underbrace{1\otimes \shift(\varphi^1\varphi^{21})\otimes \varphi^{22}\otimes \shift \varphi^{23}\otimes \varphi^{24}\varphi^{25}\varphi^3\varphi^4\varphi^5}_y\\
&-\underbrace{1\otimes\shift  \varphi^1\otimes \varphi^{21}\varphi^{22}\otimes \shift \varphi^{23}\otimes \varphi^{24}\varphi^{25}\varphi^3\varphi^4\varphi^5}_{44}\boldsymbol{\uwave{+}}0
\boldsymbol{\uwave{-}}\underbrace{1\otimes \shift \varphi^1\otimes \shift \varphi^{21}\otimes \varphi^{22}\hat{\omega}\varphi^{23}\varphi^{24}\varphi^{25}\varphi^3\varphi^4\varphi^5}_{43}\\
&+\underbrace{1\otimes \shift \varphi^1\otimes \shift \varphi^{21}\otimes \varphi^{22}\varphi^{23}\hat{\omega}\varphi^{24}\varphi^{25}\varphi^3\varphi^4\varphi^5}_{50}
+\underbrace{1\otimes \shift \varphi^1\otimes \shift \varphi^{21}\otimes \varphi^{22}\shift \partial(\varphi^{23})\varphi^{24}\varphi^{25}\varphi^3\varphi^4\varphi^5}_y\boldsymbol{\uwave{+}}0\boldsymbol{\uwave{+}}0
\end{align*}
and
\begin{align*}
&\mu_1^{\otimes}(-\varphi^1\varphi^2\varphi^3\varphi^{41}\varphi^{42}\otimes \shift \varphi^{43}\otimes \varphi^{44}\otimes \shift \varphi^{45}\otimes \shift \varphi^5\otimes 1)\\
=\,&0\boldsymbol{\uwave{+}}0\boldsymbol{\uwave{-}}\underbrace{\varphi^1\varphi^2\varphi^3\varphi^{41}\varphi^{42}\hat{\omega}\varphi^{43}\varphi^{44}\otimes \shift \varphi^{45}\otimes \shift \varphi^5\otimes 1}_{37}
+\underbrace{\varphi^1\varphi^2\varphi^3\varphi^{41}\varphi^{42}\varphi^{43}\hat{\omega}\varphi^{44}\otimes \shift \varphi^{45}\otimes \shift \varphi^5\otimes 1}_{30}\\
&+\underbrace{\varphi^1\varphi^2\varphi^3\varphi^{41}\varphi^{42}\shift \partial(\varphi^{43})\varphi^{44}\otimes \shift \varphi^{45}\otimes \shift \varphi^5\otimes 1}_y\boldsymbol{\uwave{+}}0
\boldsymbol{\uwave{+}}\underbrace{\varphi^1\varphi^2\varphi^3\varphi^{41}\varphi^{42}\otimes \shift \varphi^{43}\otimes \varphi^{44}\varphi^{45}\otimes \shift \varphi^5\otimes 1}_{31}\\
&-\underbrace{\varphi^1\varphi^2\varphi^3\varphi^{41}\varphi^{42}\otimes \shift \varphi^{43}\otimes\varphi^{44}\otimes \shift(\varphi^{45}\varphi^5)\otimes 1}_y
+\underbrace{\varphi^1\varphi^2\varphi^3\varphi^{41}\varphi^{42}\otimes \shift \varphi^{43}\otimes \varphi^{44}\otimes \shift \varphi^{45}\otimes \varphi^5}_{29}
\end{align*}
and
\begin{align*}
&\mu_1^{\otimes}(-\varphi^1\varphi^{21}\dots \varphi^{24}\otimes \shift \varphi^{25}\otimes \shift \varphi^3\otimes \shift \varphi^{41}\otimes \varphi^{42}\dots \varphi^{45}\varphi^5)\\
=\,&0\boldsymbol{\uwave{+}}0\boldsymbol{\uwave{+}}\underbrace{\varphi^1(\shift \otimes \shift)\partial(\shift^{-1} \varphi^2)\otimes \shift \varphi^3\otimes \shift \varphi^{41}\otimes \varphi^{42}\dots \varphi^{45}\varphi^5}_{45}
+\underbrace{\varphi^1\varphi^{21}\dots \varphi^{24}\otimes \shift(\varphi^{25}\varphi^3)\otimes \shift \varphi^{41}\otimes \varphi^{42}\dots \varphi^{45}\varphi^5}_g\\
&-\underbrace{\varphi^1\varphi^{21}\dots \varphi^{24}\otimes \shift \varphi^{25}\otimes \shift(\varphi^3\varphi^{41})\otimes \varphi^{42}\dots \varphi^{45}\varphi^5}_g
-\underbrace{\varphi^1\varphi^{21}\dots \varphi^{24}\otimes \shift \varphi^{25}\otimes \shift \varphi^3\otimes (\shift \otimes \shift)\partial(\shift^{-1} \varphi^4)\varphi^5}_{32}\boldsymbol{\uwave{+}}0\boldsymbol{\uwave{+}}0
\end{align*}
and
\begin{align*}
&\mu_1^{\otimes}(-\varphi^1\varphi^2\varphi^3\varphi^{41}\dots \varphi^{43}\varphi^{441}\dots \varphi^{444}\otimes \shift \varphi^{445}\otimes \shift \varphi^{45}\otimes \shift \varphi^5\otimes 1)\\
=\,&0\boldsymbol{\uwave{+}}0\boldsymbol{\uwave{+}}0\boldsymbol{\uwave{+}}0\boldsymbol{\uwave{+}}\underbrace{\varphi^1\varphi^2\varphi^3\varphi^{41}\varphi^{42} \varphi^{43}(\shift \otimes \shift)\partial(\shift^{-1}\varphi^{44})\otimes \shift \varphi^{45}\otimes \shift \varphi^5\otimes 1}_{46}\\
&
+\underbrace{\varphi^1\varphi^2\varphi^3\varphi^{41}\dots \varphi^{43}\varphi^{441}\dots \varphi^{444}\otimes \shift(\varphi^{445}\varphi^{45})\otimes \shift \varphi^5\otimes 1}_y\\
&-\underbrace{\varphi^1\varphi^2\varphi^3\varphi^{41}\dots \varphi^{43}\varphi^{441}\dots \varphi^{444}\otimes \shift \varphi^{445}\otimes \shift (\varphi^{45}\varphi^5)\otimes 1}_g\\
&
+\underbrace{\varphi^1\varphi^2\varphi^3\varphi^{41}\varphi^{42} \varphi^{43}\varphi^{441}\dots \varphi^{444}\otimes \shift \varphi^{445}\otimes \shift \varphi^{45}\otimes \varphi^5}_{47}
\end{align*}
and
\begin{align*}
&\mu_1^{\otimes}(-1\otimes \shift \varphi^1\otimes \shift \varphi^{21}\otimes \shift \varphi^{221}\otimes \varphi^{222}\dots \varphi^{225}\varphi^{23}\dots \varphi^{25}\varphi^3\varphi^4\varphi^5)\\
=\,&-\underbrace{\varphi^1\otimes \shift \varphi^{21}\otimes \shift \varphi^{221}\otimes \varphi^{222}\dots \varphi^{225}\varphi^{23}\varphi^{24} \varphi^{25}\varphi^3\varphi^4\varphi^5}_{34}\\
&
+\underbrace{1\otimes \shift(\varphi^1\varphi^{21})\otimes \shift \varphi^{221}\otimes \varphi^{222}\dots \varphi^{225}\varphi^{23}\dots \varphi^{25}\varphi^3\varphi^4\varphi^5}_g\\
&-\underbrace{1\otimes \shift \varphi^1\otimes \shift(\varphi^{21}\varphi^{221})\otimes \varphi^{222}\dots \varphi^{225}\varphi^{23}\dots \varphi^{25}\varphi^3\varphi^4\varphi^5}_y\\
&
-\underbrace{1\otimes \shift \varphi^1\otimes \shift \varphi^{21}\otimes (\shift \otimes \shift)\partial(\shift^{-1} \varphi^{22})\varphi^{23}\varphi^{24} \varphi^{25}\varphi^3\varphi^4\varphi^5}_{33}\boldsymbol{\uwave{+}}0\boldsymbol{\uwave{+}}0\boldsymbol{\uwave{+}}0\boldsymbol{\uwave{+}}0
\end{align*}
and
\begin{align*}
&\mu_1^{\otimes}(\varphi^1\varphi^{21}\dots \varphi^{23}\varphi^{241}\dots \varphi^{244}\otimes \shift \varphi^{245}\otimes \shift \varphi^{25}\otimes \shift \varphi^3\otimes \varphi^4\varphi^5)\\
=\,&0\boldsymbol{\uwave{+}}0\boldsymbol{\uwave{+}}0\boldsymbol{\uwave{-}}\underbrace{\varphi^1\varphi^{21}\varphi^{22}\varphi^{23}(\shift \otimes \shift)\partial(\shift^{-1}\varphi^{24})\otimes \shift \varphi^{25}\otimes \shift \varphi^3\otimes \varphi^4\varphi^5}_{48}\\
&
-\underbrace{\varphi^1\varphi^{21}\varphi^{22}\varphi^{23}\varphi^{241}\dots \varphi^{244}\otimes \shift(\varphi^{245}\varphi^{25})\otimes \shift \varphi^3\otimes \varphi^4\varphi^5}_y\\
&+\underbrace{\varphi^1\varphi^{21}\dots \varphi^{23}\varphi^{241}\dots \varphi^{244}\otimes \shift \varphi^{245}\otimes \shift(\varphi^{25}\varphi^3)\otimes \varphi^4\varphi^5}_g\\
&
-\underbrace{\varphi^1\varphi^{21}\varphi^{22}\varphi^{23}\varphi^{241}\dots \varphi^{244}\otimes \shift \varphi^{245}\otimes \shift \varphi^{25}\otimes \varphi^3\varphi^4\varphi^5}_{49}\boldsymbol{\uwave{+}}0
\end{align*}
and
\begin{align*}
&\mu_1^{\otimes}(\varphi^1\varphi^2\otimes \shift \varphi^3\otimes \shift \varphi^{41}\otimes \shift \varphi^{421}\otimes \varphi^{422}\dots \varphi^{425}\varphi^{43}\dots \varphi^{45}\varphi^5)\\
=\,&0\boldsymbol{\uwave{+}}\underbrace{\varphi^1\varphi^2\varphi^3\otimes \shift \varphi^{41}\otimes \shift \varphi^{421}\otimes \varphi^{422}\dots \varphi^{425}\varphi^{43}\varphi^{44} \varphi^{45}\varphi^5}_{36}\\
&
-\underbrace{\varphi^1\varphi^2\otimes \shift(\varphi^3\varphi^{41})\otimes \shift \varphi^{421}\otimes \varphi^{422}\dots \varphi^{425}\varphi^{43}\dots \varphi^{45}\varphi^5}_g\\
&+\underbrace{\varphi^1\varphi^2\otimes \shift \varphi^3\otimes \shift(\varphi^{41}\varphi^{421})\otimes \varphi^{422}\dots \varphi^{425}\varphi^{43}\dots \varphi^{45}\varphi^5}_y\\
&
+\underbrace{\varphi^1\varphi^2\otimes \shift \varphi^3\otimes \shift \varphi^{41}\otimes (\shift \otimes \shift)\partial(\shift^{-1}\varphi^{42})\varphi^{43}\dots \varphi^{45}\varphi^5}_{35}\boldsymbol{\uwave{+}}0\boldsymbol{\uwave{+}}0\boldsymbol{\uwave{+}}0.
\end{align*}

As usual, terms with the same label sum to zero, finishing the calculation for $\mu_5^0$. 

\eqref{mu5:1}: For $\mu_5^1$ we compute:

\begin{align*}
&-(\mu_4\otimes 1)\mu_2^1(1\otimes \shift a\otimes 1)\\
=\,&-(\mu_4\otimes 1)(\hat{\omega}\otimes \shift a\otimes 1+1\otimes \shift a\otimes \hat{\omega}-a^1a^2\otimes \shift a^3\otimes 1+1\otimes \shift a^1\otimes a^2a^3)\\
=\,&0\+\underbrace{1\otimes \shift a^1\otimes \shift a^{21}\otimes \shift a^{221}\otimes a^{222}\dots a^{225}a^{23}\dots a^{25}a^3\hat{\omega}}_{69}\\
&
-\underbrace{a^1a^{21}\dots a^{23}a^{241}\dots a^{244}\otimes \shift a^{245}\otimes \shift a^{25}\otimes \shift a^3\otimes \hat{\omega}}_{20}\\
&-\underbrace{a^1a^{21}a^{22}\otimes \shift a^{23}\otimes a^{24}\otimes \shift a^{25}\otimes \shift a^3\otimes \hat{\omega}}_{22}
+\underbrace{1\otimes \shift a^1\otimes \shift a^{21}\otimes a^{22}\otimes \shift a^{23}\otimes a^{24}a^{25}a^3\hat{\omega}}_{52}\\
&-\underbrace{1\otimes \shift a^1\otimes \shift a^{21}\otimes a^{22}\dots a^{25}\otimes \shift a^3\otimes \hat{\omega}}_{21}
+\underbrace{1\otimes \shift (a^1a^{21})\otimes a^{22}\dots a^{24}\otimes \shift a^{25}\otimes \shift a^3\otimes \hat{\omega}}_{23}\\
&-\underbrace{1\otimes \shift a^1\otimes \shift a^{21}\otimes a^{22}\dots a^{24}\otimes \shift a^{25}\otimes a^3\hat{\omega}}_{31}
\+\underbrace{a^1a^{21}a^{22}\otimes \shift a^{23}\otimes a^{24}\otimes \shift a^{25}\otimes 1\otimes \shift a^3\otimes 1}_{26}\\
&
-\underbrace{a^1\otimes \shift a^{21}\otimes a^{22}\dots a^{24}\otimes \shift a^{25}\otimes 1\otimes \shift a^3\otimes 1}_{28}
+\underbrace{a^1\otimes \shift a^{21}\otimes a^{22}\otimes \shift a^{23}\otimes a^{24}a^{25}\otimes \shift a^3\otimes 1}_{48}\\
&+\underbrace{a^1a^{21}a^{221}\dots a^{224}\otimes \shift a^{225}\otimes \shift a^{23}\otimes a^{24}a^{25}\otimes \shift a^3\otimes 1}_{83}\\
&
+\underbrace{a^1\otimes \shift a^{21}\otimes \shift a^{221}\otimes a^{222}\dots a^{225}a^{23}\dots a^{25}\otimes \shift a^3\otimes 1}_{75}\\
&+\underbrace{a^1a^{21}a^{22}\otimes \shift a^{23}\otimes \shift a^{241}\otimes a^{242}\dots a^{245}a^{25}\otimes \shift a^3\otimes 1}_{98}\\
&
+\underbrace{a^1a^{21}\dots a^{23}a^{241}\dots a^{244}\otimes \shift a^{245}\otimes \shift a^{25}\otimes 1\otimes \shift a^3\otimes 1}_{27}\\
&\+\underbrace{1\otimes \shift a^{11}\otimes \shift a^{121}\otimes \shift a^{1221}\otimes a^{1222}\dots a^{1225}a^{123}\dots a^{125}a^{13}a^2a^3}_{gg}\\
&-\underbrace{a^{11}a^{121}\dots a^{123}a^{1241}\dots a^{1244}\otimes \shift a^{1245}\otimes \shift a^{125}\otimes \shift a^{13}\otimes a^2a^3}_{gg}\\
&-\underbrace{a^{11} a^{121}a^{122}\otimes \shift a^{123}\otimes a^{124}\otimes \shift a^{125}\otimes \shift a^{13}\otimes a^2a^3}_{g}
+\underbrace{1\otimes \shift a^{11}\otimes \shift a^{121}\otimes a^{122}\otimes \shift a^{123}\otimes a^{124}a^{125}a^{13}a^2a^3}_{ g}\\
&-\underbrace{1\otimes \shift a^{11}\otimes \shift a^{121}\otimes a^{122}\dots a^{125}\otimes \shift a^{13}\otimes a^2a^3}_{g}
+\underbrace{1\otimes \shift (a^{11} a^{121})\otimes a^{122}\dots a^{124}\otimes \shift a^{125}\otimes \shift a^{13}\otimes a^2a^3}_{g}\\
&-\underbrace{1\otimes \shift a^{11}\otimes \shift a^{121}\otimes a^{122}\dots a^{124}\otimes \shift a^{125}\otimes a^{13}a^2a^3}_{g}
\end{align*}
and
\begin{align*}
&(1\otimes \mu_4)\mu_2^1(1\otimes \shift a\otimes 1)\\
=\,&-\underbrace{\hat{\omega}\otimes \shift a^1\otimes \shift a^{21}\otimes \shift a^{221}\otimes a^{222}\dots a^{225}a^{23}\dots a^{25}a^3}_{2}\\
&
+\underbrace{\hat{\omega}a^1a^{21}\dots a^{23}a^{241}\dots a^{244}\otimes \shift a^{245}\otimes \shift a^{25}\otimes \shift a^3\otimes 1}_{42}\\
&+\underbrace{\hat{\omega}a^1a^{21}a^{22}\otimes \shift a^{23}\otimes a^{24}\otimes \shift a^{25}\otimes \shift a^3\otimes 1}_{38}
-\underbrace{\hat{\omega}\otimes \shift a^1\otimes \shift a^{21}\otimes a^{22}\otimes \shift a^{23}\otimes a^{24}a^{25}a^3}_5\\
&+\underbrace{\hat{\omega}\otimes \shift a^1\otimes \shift a^{21}\otimes a^{22}\dots a^{25}\otimes \shift a^3\otimes 1}_{10}
-\underbrace{\hat{\omega}\otimes \shift (a^1a^{21})\otimes a^{22}\dots a^{24}\otimes \shift a^{25}\otimes \shift a^3\otimes 1}_{8}\\
&+\underbrace{\hat{\omega}\otimes \shift a^1\otimes \shift a^{21}\otimes a^{22}\dots a^{24}\otimes \shift a^{25}\otimes a^3}_{14}\-0\\
&\+\underbrace{a^1a^2\otimes \shift a^{31}\otimes \shift a^{321}\otimes \shift a^{3221}\otimes a^{3222}\dots a^{3225}a^{323}\dots a^{325}a^{33}}_{gg}\\
&-\underbrace{a^1a^2a^{31}a^{321}\dots a^{323}a^{3241}\dots a^{3244}\otimes \shift a^{3245}\otimes \shift a^{325}\otimes \shift a^{33}\otimes 1}_{gg}\\
&-\underbrace{a^1a^2a^{31}a^{321}a^{322}\otimes \shift a^{323}\otimes \shift a^{324}\otimes \shift a^{325}\otimes a^{33}\otimes 1}_g
+\underbrace{a^1a^2\otimes \shift a^{31}\otimes \shift a^{321}\otimes a^{322}\otimes \shift a^{323}\otimes a^{324}a^{325}a^{33}}_g\\
&-\underbrace{a^1a^2\otimes \shift a^{31}\otimes \shift a^{321}\otimes a^{322}\dots a^{325}\otimes \shift a^{33}\otimes 1}_g
+\underbrace{a^1a^2\otimes \shift (a^{31}a^{321})\otimes a^{322}\dots a^{324}\otimes \shift a^{325}\otimes \shift a^{33}\otimes 1}_g\\
&-\underbrace{a^1a^2\otimes \shift a^{31}\otimes \shift a^{321}\otimes a^{322}\dots a^{324}\otimes \shift a^{325}\otimes a^{33}}_g
\+\underbrace{1\otimes \shift a^1\otimes a^{21}a^{22}\otimes \shift a^{23}\otimes a^{24}\otimes \shift a^{25}\otimes a^3}_{58}\\
&
-\underbrace{1\otimes \shift a^1\otimes 1\otimes \shift a^{21}\otimes a^{22}\dots a^{24}\otimes \shift a^{25}\otimes a^3}_{18}
+\underbrace{1\otimes \shift a^1\otimes 1\otimes \shift a^{21}\otimes a^{22}\otimes \shift a^{23}\otimes a^{24}a^{25}a^3}_{25}\\
&+\underbrace{1\otimes \shift a^1\otimes a^{21}a^{221}\dots a^{224}\otimes \shift a^{225}\otimes \shift a^{23}\otimes a^{24}a^{25}a^3}_{90}\\
&
+\underbrace{1\otimes \shift a^1\otimes 1\otimes \shift a^{21}\otimes \shift a^{221}\otimes a^{222}\dots a^{225}a^{23}\dots a^{25}a^3}_{24}\\
&+\underbrace{1\otimes \shift a^1\otimes a^{21}a^{22}\otimes \shift a^{23}\otimes \shift a^{241}\otimes a^{242}\dots a^{245}a^{25}a^3}_{105}\\
&
+\underbrace{1\otimes \shift a^1\otimes a^{21}\dots a^{23}a^{241}\dots a^{244}\otimes \shift a^{245}\otimes \shift a^{25}\otimes a^3}_{101}
\end{align*}
and
\begin{align*}
&-(\mu_2\otimes 1^{\otimes 3})\mu_4(1\otimes \shift a\otimes 1)\\
=\,&+\underbrace{\hat{\omega}\otimes \shift a^1\otimes \shift a^{21}\otimes \shift a^{221}\otimes a^{222}\dots a^{225}a^{23}\dots a^{25}a^3}_2\\
&+\underbrace{1\otimes \shift a^1\otimes \shift a^{21}\otimes \shift a^{221}\otimes \hat{\omega}a^{222}\dots a^{225}a^{23}\dots a^{25}a^3}_1\\
&-\underbrace{1\otimes \shift a^1\otimes 1\otimes \shift a^{21}\otimes \shift a^{221}\otimes a^{222}\dots a^{225}a^{23}\dots a^{25}a^3}_{24}\\
&-\underbrace{1\otimes \shift a^1\otimes \shift a^{21}\otimes 1\otimes \shift a^{221}\otimes a^{222}\dots a^{225}a^{23}\dots a^{25}a^3}_{37}\\
&-\underbrace{a^{11}a^{12}\otimes \shift a^{13}\otimes \shift a^{21}\otimes \shift a^{221}\otimes a^{222}\dots a^{225}a^{23}\dots a^{25}a^3}_{gg}\\
&+\underbrace{1\otimes \shift a^1\otimes \shift a^{21}\otimes \shift a^{2211}\otimes a^{2212}a^{2213}a^{222}\dots a^{225}a^{23}\dots a^{25}a^3}_{y}\\
&\-\underbrace{a^1a^{21}\hat{\omega}a^{22}a^{23}a^{241}\dots a^{244}\otimes \shift a^{245}\otimes \shift a^{25}\otimes \shift a^3\otimes 1}_{43}\\
&-\underbrace{a^1a^{21}a^{22}\hat{\omega}a^{23}a^{241}\dots a^{244}\otimes \shift a^{245}\otimes \shift a^{25}\otimes \shift a^3\otimes 1}_{45}\\
&-\underbrace{a^1 a^{21}(\shift \otimes \shift)\partial(\shift^{-1}a^{22})a^{23}a^{241}\dots a^{244}\otimes \shift a^{245}\otimes \shift a^{25}\otimes \shift a^3\otimes 1}_g\\
&\-\underbrace{a^1a^{21}\hat{\omega}a^{22}\otimes \shift a^{23}\otimes a^{24}\otimes \shift a^{25}\otimes \shift a^3\otimes 1}_{39}
-\underbrace{a^1a^{21}a^{22}\hat{\omega}\otimes \shift a^{23}\otimes a^{24}\otimes \shift a^{25}\otimes \shift a^3\otimes 1}_3\\
&-\underbrace{a^1a^{21}(\shift \otimes \shift)\partial(\shift^{-1}a^{22})\otimes \shift a^{23}\otimes a^{24}\otimes \shift a^{25}\otimes \shift a^3\otimes 1}_{82}
\+\underbrace{\hat{\omega}\otimes \shift a^1\otimes \shift a^{21}\otimes a^{22}\otimes \shift a^{23}\otimes a^{24}a^{25}a^3}_5\\
&+\underbrace{1\otimes \shift a^1\otimes \shift a^{21}\otimes \hat{\omega}a^{22}\otimes \shift a^{23}\otimes a^{24}a^{25}a^3}_4
-\underbrace{1\otimes \shift a^1\otimes 1\otimes \shift a^{21}\otimes a^{22}\otimes \shift a^{23}\otimes a^{24}a^{25}a^3}_{25}\\
&-\underbrace{a^{11}a^{12}\otimes \shift a^{13}\otimes \shift a^{21}\otimes a^{22}\otimes \shift a^{23}\otimes a^{24}a^{25}a^3}_g
+\underbrace{1\otimes \shift a^1\otimes \shift a^{211}\otimes a^{212}a^{213}a^{22}\otimes \shift a^{23}\otimes a^{24}a^{25}a^3}_y\\
&\-\underbrace{\hat{\omega}\otimes \shift a^1\otimes \shift a^{21}\otimes a^{22}\dots a^{25}\otimes \shift a^3\otimes 1}_{10}
-\underbrace{1\otimes \shift a^1\otimes \shift a^{21}\otimes \hat{\omega}a^{22}\dots a^{25}\otimes \shift a^3\otimes 1}_6\\
&+\underbrace{1\otimes \shift a^1\otimes 1\otimes \shift a^{21}\otimes a^{22}\dots a^{25}\otimes \shift a^3\otimes 1}_{36}
+\underbrace{a^{11}a^{12}\otimes \shift a^{13}\otimes \shift a^{21}\otimes a^{22}\dots a^{25}\otimes \shift a^3\otimes 1}_{g}\\
&-\underbrace{1\otimes \shift a^1\otimes \shift a^{211}\otimes a^{212}a^{213}a^{22}\dots a^{25}\otimes \shift a^3\otimes 1}_y
\+\underbrace{\hat{\omega}\otimes \shift (a^1a^{21})\otimes a^{22}\dots a^{24}\otimes \shift a^{25}\otimes \shift a^3\otimes 1}_8\\
&+\underbrace{1\otimes \shift (a^1a^{21})\otimes \hat{\omega}a^{22}\dots a^{24}\otimes \shift a^{25}\otimes \shift a^3\otimes 1}_7
-\underbrace{a^{11}a^{12}\otimes \shift (a^{13}a^{21})\otimes a^{22}\dots a^{24}\otimes \shift a^{25}\otimes \shift a^3\otimes 1}_g\\
&-\underbrace{a^1a^{211}a^{212}\otimes \shift a^{213}\otimes a^{22}\dots a^{24}\otimes \shift a^{25}\otimes \shift a^3\otimes 1}_y
+\underbrace{1\otimes \shift a^{11}\otimes a^{12}a^{13}a^{21}\dots a^{24}\otimes \shift a^{25}\otimes \shift a^3\otimes 1}_g\\
&+\underbrace{1\otimes \shift (a^1a^{211})\otimes a^{212}a^{213}a^{22}\dots a^{24}\otimes \shift a^{25}\otimes \shift a^3\otimes 1}_y
\-\underbrace{\hat{\omega}\otimes \shift a^1\otimes \shift a^{21}\otimes a^{22}\dots a^{24}\otimes \shift a^{25}\otimes a^3}_{14}\\
&-\underbrace{1\otimes \shift a^1\otimes \shift a^{21}\otimes \hat{\omega}a^{22}\dots a^{24}\otimes \shift a^{25}\otimes a^3}_9
+\underbrace{1\otimes \shift a^1\otimes 1\otimes \shift a^{21}\otimes a^{22}\dots a^{24}\otimes \shift a^{25}\otimes a^3}_{18}\\
&+\underbrace{a^{11}a^{12}\otimes \shift a^{13}\otimes \shift a^{21}\otimes a^{22}\dots a^{24}\otimes \shift a^{25}\otimes a^3}_g
-\underbrace{1\otimes \shift a^1\otimes \shift a^{211}\otimes a^{212}a^{213}a^{22}\dots a^{24}\otimes \shift a^{25}\otimes a^3}_y
\end{align*}
and
\begin{align*}
&(1\otimes \mu_2\otimes 1^{\otimes 2})\mu_4(1\otimes \shift a\otimes 1)\\
=\,&-\underbrace{1\otimes \shift a^1\otimes \shift a^{21}\otimes \shift a^{221}\otimes \hat{\omega}a^{222}\dots a^{225}a^{23}\dots a^{25}a^3}_1\\
&-\underbrace{1\otimes \shift a^1\otimes \shift a^{21}\otimes \shift a^{221}\otimes a^{222}\hat{\omega}a^{223}\dots a^{225}a^{23}\dots a^{25}a^3}_{71}\\
&-\underbrace{1\otimes \shift a^1\otimes \shift a^{21}\otimes \shift a^{221}\otimes (\shift \otimes \shift)\partial(\shift^{-1}a^{222})a^{223}\dots a^{225}a^{23}\dots a^{25}a^3}_{109}\\
&\+\underbrace{a^1a^{21}a^{22}a^{23}a^{241}\hat{\omega}a^{242}\dots a^{244}\otimes \shift a^{245}\otimes \shift a^{25}\otimes \shift a^3\otimes 1}_{46}\\
&+\underbrace{a^1a^{21}\dots a^{23}a^{241}a^{242}\hat{\omega}a^{243}a^{244}\otimes \shift a^{245}\otimes \shift a^{25}\otimes \shift a^3\otimes 1}_{61}\\
&+\underbrace{a^1a^{21}\dots a^{23}a^{241}(\shift \otimes \shift)\partial(\shift^{-1}a^{242})a^{243}a^{244}\otimes \shift a^{245}\otimes \shift a^{25}\otimes \shift a^3\otimes 1}_y\\
&\+\underbrace{a^1a^{21}a^{22}\hat{\omega}\otimes \shift a^{23}\otimes a^{24}\otimes \shift a^{25}\otimes \shift a^3\otimes 1}_3
+\underbrace{a^1a^{21}a^{22}\otimes \shift a^{23}\otimes \hat{\omega}a^{24}\otimes \shift a^{25}\otimes \shift a^3\otimes 1}_{17}\\
&-\underbrace{a^1a^{21}a^{22}a^{231}a^{232}\otimes \shift a^{233}\otimes a^{24}\otimes \shift a^{25}\otimes \shift a^3\otimes 1}_y
+\underbrace{a^1a^{21}a^{22}\otimes \shift a^{231}\otimes a^{232}a^{233}a^{24}\otimes \shift a^{25}\otimes \shift a^3\otimes 1}_y\\
&\-\underbrace{1\otimes \shift a^1\otimes \shift a^{21}\otimes \hat{\omega}a^{22}\otimes \shift a^{23}\otimes a^{24}a^{25}a^3}_4
-\underbrace{1\otimes \shift a^1\otimes \shift a^{21}\otimes a^{22}\hat{\omega}\otimes \shift a^{23}\otimes a^{24}a^{25}a^3}_{11}\\
&-\underbrace{1\otimes \shift a^1\otimes \shift a^{21}\otimes (\shift \otimes \shift)\partial(\shift^{-1}a^{22})\otimes \shift a^{23}\otimes a^{24}a^{25}a^3}_{91}
\+\underbrace{1\otimes \shift a^1\otimes \shift a^{21}\otimes \hat{\omega}a^{22}\dots a^{25}\otimes \shift a^3\otimes 1}_6\\
&+\underbrace{1\otimes \shift a^1\otimes \shift a^{21}\otimes a^{22}\hat{\omega}a^{23}a^{24} a^{25}\otimes \shift a^3\otimes 1}_{50}
+\underbrace{1\otimes \shift a^1\otimes \shift a^{21}\otimes (\shift \otimes \shift)\partial(\shift^{-1}a^{22})a^{23}\dots a^{25}\otimes \shift a^3\otimes 1}_{80}\\
&\-\underbrace{1\otimes \shift(a^1a^{21})\otimes \hat{\omega}a^{22}\dots a^{24}\otimes \shift a^{25}\otimes \shift a^3\otimes 1}_7\\
&-\underbrace{1\otimes \shift (a^1a^{21})\otimes a^{22}\hat{\omega}a^{23}a^{24}\otimes \shift a^{25}\otimes \shift a^3\otimes 1}_{40}
-\underbrace{1\otimes \shift (a^1a^{21})\otimes (\shift \otimes \shift)\partial(\shift^{-1}a^{22})a^{23}a^{24}\otimes \shift a^{25}\otimes \shift a^3\otimes 1}_{78}\\
&\+\underbrace{1\otimes \shift a^1\otimes \shift a^{21}\otimes \hat{\omega}a^{22}\dots a^{24}\otimes \shift a^{25}\otimes a^3}_9
+\underbrace{1\otimes \shift a^1\otimes \shift a^{21}\otimes a^{22}\hat{\omega}a^{23}a^{24}\otimes \shift a^{25}\otimes a^3}_{59}\\
&+\underbrace{1\otimes \shift a^1\otimes \shift a^{21}\otimes (\shift \otimes \shift)\partial(\shift^{-1}a^{22})a^{23}a^{24}\otimes \shift a^{25}\otimes a^3}_{81}
\end{align*}
and
\begin{align*}
&-(1^{\otimes 2}\otimes \mu_2\otimes 1)\mu_4(1\otimes \shift a\otimes 1)\\
=\,&\underbrace{1\otimes \shift a^1\otimes \shift a^{21}\otimes \shift a^{221}\otimes a^{222}a^{223}\hat{\omega}a^{224}a^{225}a^{23}\dots a^{25}a^3}_{72}\\
&+\underbrace{1\otimes \shift a^1\otimes \shift a^{21}\otimes \shift a^{221}\otimes a^{222}\dots a^{224}\hat{\omega}a^{225}a^{23}\dots a^{25}a^3}_{66}\\
&+\underbrace{1\otimes \shift a^1\otimes \shift a^{21}\otimes \shift a^{221}\otimes a^{222}a^{223}(\shift \otimes \shift)\partial(\shift^{-1}a^{224})a^{225}a^{23}\dots a^{25}a^3}_y\\
&\-\underbrace{a^1a^{21}\dots a^{23}a^{241}\dots a^{243}\hat{\omega}a^{244}\otimes \shift a^{245}\otimes \shift a^{25}\otimes \shift a^3\otimes 1}_{62}\\
&-\underbrace{a^1a^{21}\dots a^{23}a^{241}\dots a^{244}\hat{\omega}\otimes \shift a^{245}\otimes \shift a^{25}\otimes \shift a^3\otimes 1}_{12}\\
&-\underbrace{a^1a^{21}\dots a^{23}a^{241}\dots a^{243}(\shift \otimes \shift)\partial(\shift^{-1}a^{244})\otimes \shift a^{245}\otimes \shift a^{25}\otimes \shift a^3\otimes 1}_{110}\\
&\-\underbrace{a^1a^{21}a^{22}\otimes \shift a^{23}\otimes \hat{\omega}a^{24}\otimes \shift a^{25}\otimes \shift a^3\otimes 1}_{17}
-\underbrace{a^1a^{21}a^{22}\otimes \shift a^{23}\otimes a^{24}\hat{\omega}\otimes \shift a^{25}\otimes \shift a^3\otimes 1}_{16}\\
&-\underbrace{a^1a^{21}a^{22}\otimes \shift a^{23}\otimes (\shift \otimes \shift)\partial(\shift^{-1}a^{24})\otimes \shift a^{25}\otimes \shift a^3\otimes 1}_{96}
\+\underbrace{1\otimes \shift a^1\otimes \shift a^{21}\otimes a^{22}\hat{\omega}\otimes \shift a^{23}\otimes a^{24}a^{25}a^3}_{11}\\
&+\underbrace{1\otimes \shift a^1\otimes \shift a^{21}\otimes a^{22}\otimes \shift a^{23}\otimes \hat{\omega}a^{24}a^{25}a^3}_{19}
-\underbrace{1\otimes \shift a^1\otimes \shift a^{21}\otimes a^{22}a^{231}a^{232}\otimes \shift a^{233}\otimes a^{24}a^{25}a^3}_y\\
&+\underbrace{1\otimes \shift a^1\otimes \shift a^{21}\otimes a^{22}\otimes \shift a^{231}\otimes a^{232}a^{233}a^{24}a^{25}a^3}_y
\-\underbrace{1\otimes \shift a^1\otimes \shift a^{21}\otimes a^{22}a^{23}\hat{\omega}a^{24}a^{25}\otimes \shift a^3\otimes 1}_{51}\\
&-\underbrace{1\otimes \shift a^1\otimes \shift a^{21}\otimes a^{22}\dots a^{24}\hat{\omega}a^{25}\otimes \shift a^3\otimes 1}_{33}
-\underbrace{1\otimes \shift a^1\otimes \shift a^{21}\otimes a^{22}a^{23}(\shift \otimes \shift)\partial(\shift^{-1}a^{24})a^{25}\otimes a^3\otimes 1}_y\\
&\+\underbrace{1\otimes \shift (a^1a^{21})\otimes a^{22}a^{23}\hat{\omega}a^{24}\otimes \shift a^{25}\otimes \shift a^3\otimes 1}_{41}
+\underbrace{1\otimes \shift (a^1a^{21})\otimes a^{22}\dots a^{24}\hat{\omega}\otimes \shift a^{25}\otimes \shift a^3\otimes 1}_{15}\\
&+\underbrace{1\otimes \shift (a^1a^{21})\otimes a^{22}a^{23}(\shift \otimes \shift)\partial(\shift^{-1}a^{24})\otimes \shift a^{25}\otimes \shift a^3\otimes 1}_{44}
\-\underbrace{1\otimes \shift a^1\otimes \shift a^{21}\otimes a^{22}a^{23}\hat{\omega}a^{24}\otimes \shift a^{25}\otimes a^3}_{57}\\
&-\underbrace{1\otimes \shift a^1\otimes \shift a^{21}\otimes a^{22}\dots a^{24}\hat{\omega}\otimes \shift a^{25}\otimes a^3}_{13}
-\underbrace{1\otimes \shift a^1\otimes \shift a^{21}\otimes a^{22}a^{23}(\shift \otimes \shift)\partial(\shift^{-1}a^{24})\otimes \shift a^{25}\otimes a^3}_{102}
\end{align*}
and
\begin{align*}
&(1^{\otimes 3}\otimes \mu_2)\mu_4(1\otimes \shift a\otimes 1)\\
=\,&-\underbrace{1\otimes \shift a^1\otimes \shift a^{21}\otimes \shift a^{221}\otimes a^{222}\dots a^{225}a^{23}\hat{\omega}a^{24}a^{25}a^3}_{67}\\
&-\underbrace{1\otimes \shift a^1\otimes \shift a^{21}\otimes \shift a^{221}\otimes a^{222}\dots a^{225}a^{23}a^{24}\hat{\omega}a^{25}a^3}_{68}\\
&-\underbrace{1\otimes \shift a^1\otimes \shift a^{21}\otimes \shift a^{221}\otimes a^{222}\dots a^{225}a^{23}(\shift \otimes \shift)\partial(\shift^{-1}a^{24})a^{25}a^3}_g\\
&\+\underbrace{a^1a^{21}\dots a^{23}a^{241}\dots a^{244}\hat{\omega}\otimes \shift a^{245}\otimes \shift a^{25}\otimes \shift a^3\otimes 1}_{12}\\
&+\underbrace{a^1a^{21}\dots a^{23}a^{241}\dots a^{244}\otimes \shift a^{245}\otimes \shift a^{25}\otimes \shift a^3\otimes \hat{\omega}}_{20}\\
&-\underbrace{a^1a^{21}\dots a^{23}a^{241}\dots a^{244}\otimes \shift a^{245}\otimes 1\otimes \shift a^{25}\otimes \shift a^3\otimes 1}_{35}\\
&-\underbrace{a^1a^{21}\dots a^{23}a^{241}\dots a^{244}\otimes \shift a^{245}\otimes \shift a^{25}\otimes 1\otimes \shift a^3\otimes 1}_{27}\\
&-\underbrace{a^1a^{21}\dots a^{23}a^{241}\dots a^{244}a^{2451}a^{2452}\otimes \shift a^{2453}\otimes \shift a^{25}\otimes \shift a^3\otimes 1}_{y}\\
&+\underbrace{a^1a^{21}\dots a^{23}a^{241}\dots a^{244}\otimes \shift a^{245}\otimes \shift a^{25}\otimes \shift a^{31}\otimes a^{32}a^{33}}_{gg}\\
&\+\underbrace{a^1a^{21}a^{22}\otimes \shift a^{23}\otimes a^{24}\hat{\omega}\otimes \shift a^{25}\otimes \shift a^3\otimes 1}_{16}\\
&+\underbrace{a^1a^{21}a^{22}\otimes \shift a^{23}\otimes a^{24}\otimes \shift a^{25}\otimes \shift a^3\otimes \hat{\omega}}_{22}
-\underbrace{a^1a^{21}a^{22}\otimes \shift a^{23}\otimes a^{24}\otimes \shift a^{25}\otimes 1\otimes \shift a^3\otimes 1}_{26}\\
&-\underbrace{a^1a^{21}a^{22}\otimes \shift a^{23}\otimes a^{24}a^{251}a^{252}\otimes \shift a^{253}\otimes a^3\otimes 1}_y
+\underbrace{a^1a^{21}a^{22}\otimes \shift a^{23}\otimes a^{24}\otimes \shift a^{25}\otimes \shift a^{31}\otimes a^{32}a^{33}}_g\\
&\-\underbrace{1\otimes \shift a^1\otimes \shift a^{21}\otimes a^{22}\otimes \shift a^{23}\otimes \hat{\omega}a^{24}a^{25}a^3}_{19}
-\underbrace{1\otimes \shift a^1\otimes \shift a^{21}\otimes a^{22}\otimes \shift a^{23}\otimes a^{24}\hat{\omega}a^{25}a^3}_{54}\\
&-\underbrace{1\otimes \shift a^1\otimes \shift a^{21}\otimes a^{22}\otimes \shift a^{23}\otimes (\shift \otimes \shift)\partial(\shift^{-1}a^{24})a^{25}a^3}_{107}
\+\underbrace{1\otimes \shift a^1\otimes \shift a^{21}\otimes a^{22}\dots a^{25}\hat{\omega}\otimes \shift a^3\otimes 1}_{32}\\
&+\underbrace{1\otimes \shift a^1\otimes \shift a^{21}\otimes a^{22}\dots a^{25}\otimes \shift a^3\otimes \hat{\omega}}_{21}
-\underbrace{1\otimes \shift a^1\otimes \shift a^{21}\otimes a^{22}\dots a^{25}a^{31}a^{32}\otimes \shift a^{33}\otimes 1}_g\\
&+\underbrace{1\otimes \shift a^1\otimes \shift a^{21}\otimes a^{22}\dots a^{25}\otimes \shift a^{31}\otimes a^{32}a^{33}}_g
\-\underbrace{1\otimes \shift (a^1a^{21})\otimes a^{22}\dots a^{24}\hat{\omega}\otimes \shift a^{25}\otimes \shift a^3\otimes 1}_{15}\\
&-\underbrace{1\otimes \shift (a^1a^{21})\otimes a^{22}\dots a^{24}\otimes \shift a^{25}\otimes \shift a^3\otimes \hat{\omega}}_{23}
+\underbrace{1\otimes \shift (a^1a^{21})\otimes a^{22}\dots a^{24}\otimes \shift a^{25}\otimes 1\otimes \shift a^3\otimes 1}_{29}\\
&+\underbrace{1\otimes \shift (a^1a^{21})\otimes a^{22}\dots a^{24}a^{251}a^{252}\otimes \shift a^{253}\otimes \shift a^3\otimes 1}_y\\
&-\underbrace{1\otimes \shift (a^1a^{21})\otimes a^{22}\dots a^{24}\otimes \shift a^{25}\otimes \shift a^{31}\otimes a^{32}a^{33}}_g\\
&\+\underbrace{1\otimes \shift a^1\otimes \shift a^{21}\otimes a^{22}\dots a^{24}\hat{\omega}\otimes \shift a^{25}\otimes a^3}_{13}
+\underbrace{1\otimes \shift a^1\otimes \shift a^{21}\otimes a^{22}\dots a^{24}\otimes \shift a^{25}\otimes \hat{\omega}a^3}_{30}\\
&-\underbrace{1\otimes \shift a^1\otimes \shift a^{21}\otimes a^{22}\dots a^{24}a^{251}a^{252}\otimes \shift a^{253}\otimes a^3}_y\\
&+\underbrace{1\otimes \shift a^1\otimes \shift a^{21}\otimes a^{22}\dots a^{24}\otimes \shift a^{251}\otimes a^{252}a^{253}a^3}_y
\end{align*}

and
\begin{align*}
&(\mu_3\otimes 1^{\otimes 2})\mu_3(1\otimes \shift a\otimes 1)\\
=\,&\underbrace{1\otimes \shift a^{11}\otimes a^{12}\otimes \shift a^{13}\otimes a^2\otimes \shift a^3\otimes 1}_y
+\underbrace{a^{11}a^{121}\dots a^{124}\otimes \shift a^{125}\otimes \shift a^{13}\otimes a^2\otimes \shift a^3\otimes 1}_g\\
&+\underbrace{1\otimes \shift a^{11}\otimes \shift a^{121}\otimes a^{122}\dots a^{125}a^{13}a^2\otimes \shift a^3\otimes 1}_g\\
&\+\underbrace{a^1a^{21}\otimes \shift a^{221}\otimes a^{222}\dots a^{225}a^{23}a^{24}\otimes \shift a^{25}\otimes \shift a^3\otimes 1}_{76}\\
&-\underbrace{a^1a^{21}a^{221}a^{222}\otimes \shift a^{223}\otimes a^{224}a^{225}a^{23}a^{24}\otimes \shift a^{25}\otimes \shift a^3\otimes 1}_y\\
&+\underbrace{a^1a^{21}a^{221}\dots a^{224}\otimes \shift a^{225}\otimes a^{23}a^{24}\otimes \shift a^{25}\otimes \shift a^3\otimes 1}_{85}\\
&\+\underbrace{1\otimes \shift a^{11}\otimes a^{12}\otimes \shift a^{13}\otimes \shift a^{21}\otimes a^{22}\dots a^{25}a^3}_g
-\underbrace{1\otimes \shift a^1\otimes \shift a^{211}\otimes a^{212}\otimes \shift a^{213}\otimes a^{22}\dots a^{25}a^3}_y\\
&+\underbrace{a^{11}a^{121}\dots a^{124}\otimes \shift a^{125}\otimes \shift a^{13}\otimes \shift a^{21}\otimes a^{22}\dots a^{25}a^3}_{gg}\\
&+\underbrace{1\otimes \shift a^1\otimes \shift a^{211}\otimes \shift a^{2121}\otimes a^{2122}\dots a^{2125}a^{213}a^{22}\dots a^{25}a^3}_g 
\end{align*}
and
\begin{align*}
&(1\otimes \mu_3\otimes 1)\mu_3(1\otimes \shift a\otimes 1)\\
=\,&-\underbrace{1\otimes \shift a^1\otimes 1\otimes \shift a^{21}\otimes a^{22}\dots a^{25}\otimes \shift a^3\otimes 1}_{36}
+\underbrace{1\otimes \shift a^1\otimes a^{21}a^{22}\otimes \shift a^{23}\otimes a^{24}a^{25}\otimes \shift a^3\otimes 1}_{49}\\
&-\underbrace{1\otimes \shift a^1\otimes a^{21}\dots a^{24}\otimes \shift a^{25}\otimes 1\otimes \shift a^3\otimes 1}_{34}
\+\underbrace{a^1a^{21}\dots a^{23}\otimes \shift a^{241}\otimes a^{242}\dots a^{245}\otimes \shift a^{25}\otimes \shift a^3\otimes 1}_{94}\\
&-\underbrace{a^1a^{21}\dots a^{23}a^{241}a^{242}\otimes \shift a^{243}\otimes a^{244}a^{245}\otimes \shift a^{25}\otimes \shift a^3\otimes 1}_{63}\\
&+\underbrace{a^1a^{21}\dots a^{23}a^{241}\dots a^{244}\otimes \shift a^{245}\otimes 1\otimes \shift a^{25}\otimes \shift a^3\otimes 1}_{35}\\
&\+\underbrace{1\otimes \shift a^1\otimes \shift a^{21}\otimes 1\otimes \shift a^{221}\otimes a^{222}\dots a^{225}a^{23}\dots a^{25}a^3}_{37}\\
&-\underbrace{1\otimes \shift a^1\otimes \shift a^{21}\otimes a^{221}a^{222}\otimes \shift a^{223}\otimes a^{224}a^{225}a^{23}\dots a^{25}a^3}_{73}\\
&+\underbrace{1\otimes \shift a^1\otimes \shift a^{21}\otimes a^{221}\dots a^{224}\otimes \shift a^{225}\otimes a^{23}\dots a^{25}a^3}_{93}
\end{align*}
and
\begin{align*}
&(1^{\otimes 2}\otimes \mu_3)\mu_3(1\otimes \shift a\otimes 1)\\
=\,&-\underbrace{1\otimes \shift a^1\otimes a^2\otimes \shift a^{31}\otimes a^{32}\otimes \shift a^{33}\otimes 1}_y
-\underbrace{1\otimes \shift a^1\otimes a^2a^{31}a^{321}\dots a^{324}\otimes \shift a^{325}\otimes \shift a^{33}\otimes 1}_g\\
&-\underbrace{1\otimes \shift a^1\otimes a^2\otimes \shift a^{31}\otimes \shift a^{321}\otimes a^{322}\dots a^{325}a^{33}}_g
\+\underbrace{a^1a^{21}\dots a^{24}\otimes \shift a^{251}\otimes a^{252}\otimes \shift a^{253}\otimes \shift a^3\otimes 1}_y\\
&-\underbrace{a^1a^{21}\dots a^{24}\otimes \shift a^{25}\otimes \shift a^{31}\otimes a^{32}\otimes a^{33}\otimes 1}_g\\
&+\underbrace{a^1a^{21}\dots a^{24}a^{251}a^{2521}\dots a^{2524}\otimes \shift a^{2525}\otimes \shift a^{253}\otimes \shift a^3\otimes 1}_g\\
&+\underbrace{a^1a^{21}\dots a^{24}\otimes \shift a^{25}\otimes \shift a^{31}\otimes \shift a^{321}\otimes a^{322}\dots a^{325}a^{33}}_{gg}\\
&\+\underbrace{1\otimes \shift a^1\otimes \shift a^{21}\otimes a^{22}a^{23}\otimes \shift a^{241}\otimes a^{242}\dots a^{245}a^{25}a^3}_{106}\\
&-\underbrace{1\otimes \shift a^1\otimes \shift a^{21}\otimes a^{22}a^{23}a^{241}a^{242}\otimes \shift a^{243}\otimes a^{244}a^{245}a^{25}a^3}_y\\
&+\underbrace{1\otimes \shift a^1\otimes \shift a^{21}\otimes a^{22}a^{23}a^{241}\dots a^{244}\otimes \shift a^{245}\otimes a^{25}a^3}_{103}
\end{align*}
and
\begin{align*}
&\mu_5^0\mu_1^1(1\otimes \shift a\otimes 1)= \mu_5^0(\hat{\omega} a - a\hat{\omega} -a^1a^2a^3)\\
=\,&0\+0\+\underbrace{a^1\otimes \shift a^{21}\otimes a^{22}\otimes \shift a^{23}\otimes a^{24}\otimes \shift a^{25}\otimes a^3}_{53}
+\underbrace{a^1\otimes \shift a^{21}\otimes a^{22}a^{23}a^{241}\dots a^{244}\otimes \shift a^{245}\otimes \shift a^{25}\otimes a^3}_{99}\\
&+\underbrace{a^1\otimes \shift a^{21}\otimes \shift a^{221}\otimes a^{222}\dots a^{225}a^{23}a^{24}\otimes \shift a^{25}\otimes a^3}_{77}\\
&+\underbrace{a^1\otimes \shift a^{21}\otimes a^{221}\dots a^{224}\otimes \shift a^{225}\otimes \shift a^{23}\otimes a^{24}a^{25}a^3}_{89}\\
&+\underbrace{a^1a^{21}a^{22}\otimes \shift a^{23}\otimes \shift a^{241}\otimes a^{242}\dots a^{245}\otimes \shift a^{25}\otimes a^3}_{97}\\
&+\underbrace{a^1a^{21}a^{221}\dots a^{224}\otimes \shift a^{225}\otimes \shift a^{23}\otimes a^{24}\otimes \shift a^{25}\otimes a^3}_{84}\\
&+\underbrace{a^1\otimes \shift a^{21}\otimes a^{22}\otimes \shift a^{23}\otimes \shift a^{241}\otimes a^{242}\dots a^{245}a^{25}a^3}_{104}\\
&+\underbrace{a^1a^{21}\otimes \shift a^{221}\otimes a^{222}\dots a^{224}\otimes\shift  a^{225}\otimes \shift a^{23}\otimes a^{24}a^{25}a^3}_{87}\\
&+\underbrace{a^1a^{21}a^{22}\otimes \shift a^{23}\otimes \shift a^{241}\otimes a^{242}\dots a^{244}\otimes \shift a^{245}\otimes a^{25}a^3}_{g}\\
&-\underbrace{a^1a^{21}a^{22}\otimes \shift (a^{23}a^{241})\otimes a^{242}\dots a^{244}\otimes \shift a^{245}\otimes \shift a^{25}\otimes a^3}_{47}\\
&-\underbrace{a^1\otimes \shift a^{21}\otimes \shift a^{221}\otimes a^{222}\dots a^{224}\otimes \shift (a^{225}a^{23})\otimes a^{24}a^{25}a^3}_{65}\\
&-\underbrace{a^1a^{21}a^{221}a^{222}\otimes \shift a^{223}\otimes a^{224}\otimes \shift a^{225}\otimes \shift a^{23}\otimes a^{24}a^{25}a^3}_{g}\\
&-\underbrace{a^1a^{21}a^{22}\otimes \shift a^{23}\otimes \shift a^{241}\otimes a^{242}\otimes \shift a^{243}\otimes a^{244}a^{245}a^{25}a^3}_{g}\\
&+\underbrace{a^1\otimes \shift a^{21}\otimes \shift a^{221}\otimes a^{222}\otimes \shift a^{223}\otimes a^{224}a^{225}a^{23}a^{24}a^{25}a^3}_{70}\\
&+\underbrace{a^1a^{21}a^{22}a^{23}a^{241}a^{242}\otimes \shift a^{243}\otimes a^{244}\otimes \shift a^{245}\otimes \shift a^{25}\otimes a^3}_{64}\\
&+\underbrace{a^1a^{21}a^{221}\dots a^{224}\otimes \shift a^{225}\otimes \shift a^{23}\otimes \shift a^{241}\otimes a^{242}\dots a^{245}a^{25}a^3}_{gg}\\
&+\underbrace{a^1a^{21}a^{22}a^{23}a^{241}\dots a^{243}a^{2441}\dots a^{2444}\otimes \shift a^{2445}\otimes \shift a^{245}\otimes \shift a^{25}\otimes a^3}_{88}\\
&+\underbrace{a^1\otimes \shift a^{21}\otimes \shift a^{221}\otimes \shift a^{2221}\otimes a^{2222}\dots a^{2225}a^{223}\dots a^{225}a^{23}\dots a^{25}a^3}_{108}\\
&-\underbrace{a^1a^{21}a^{221}\dots a^{223}a^{2241}\dots a^{2244}\otimes \shift a^{2245}\otimes \shift a^{225}\otimes \shift a^{23}\otimes a^{24}a^{25}a^3}_{gg}\\
&-\underbrace{a^1a^{21}a^{22}\otimes \shift a^{23}\otimes \shift a^{241}\otimes \shift a^{2421}\otimes a^{2422}\dots a^{2425}a^{243}\dots a^{245}a^{25}a^3}_{gg}
\end{align*}

Writing $\mu_1^{\otimes}:=(\mu_1\otimes 1^{\otimes 4}+1\otimes \mu_1\otimes 1^{\otimes 3}+1^{\otimes 2}\otimes \mu_1\otimes 1^{\otimes 2}+1^{\otimes 3}\otimes \mu_1\otimes 1+1^{\otimes 4}\otimes \mu_1)$ and applying this individually to each summand of our asserted $\mu_5^1(1\otimes \shift a \otimes 1)$ we obtain
\begin{align*}
&\mu_1^{\otimes}(1\otimes \shift a^1\otimes \shift a^{21}\otimes a^{22}\dots a^{24}\otimes \shift a^{25}\otimes 1\otimes \shift a^3\otimes 1)\\
=\,&\underbrace{a^1\otimes \shift a^{21}\otimes a^{22}\dots a^{24}\otimes \shift a^{25}\otimes 1\otimes \shift a^3\otimes 1}_{28}
-\underbrace{1\otimes \shift (a^1a^{21})\otimes a^{22}\dots a^{24}\otimes \shift a^{25}\otimes 1\otimes \shift a^3\otimes 1}_{29}\\
&+\underbrace{1\otimes \shift a^1\otimes a^{21}\dots a^{24}\otimes \shift a^{25}\otimes 1\otimes \shift a^3\otimes 1}_{34}
\boldsymbol{\uwave{+}}0\boldsymbol{\uwave{+}}0\boldsymbol{\uwave{+}}\underbrace{1\otimes \shift a^1\otimes \shift a^{21}\otimes a^{22}\dots a^{24}\hat{\omega}a^{25}\otimes \shift a^3\otimes 1}_{33}\\
&-\underbrace{1\otimes \shift a^1\otimes \shift a^{21}\otimes a^{22}\dots a^{25}\hat{\omega}\otimes \shift a^3\otimes 1}_{32}
-\underbrace{1\otimes \shift a^1\otimes \shift a^{21}\otimes a^{22}\dots a^{24}\shift \partial(a^{25})\otimes \shift a^3\otimes 1}_y\\
&\boldsymbol{\uwave{-}}\underbrace{1\otimes \shift a^1\otimes \shift a^{21}\otimes a^{22}\dots a^{24}\otimes \shift a^{25}\otimes \hat{\omega}a^3}_{30}
+\underbrace{1\otimes \shift a^1\otimes \shift a^{21}\otimes a^{22}\dots a^{24}\otimes \shift a^{25}\otimes a^3\hat{\omega}}_{31}\\
&+\underbrace{1\otimes \shift a^1\otimes \shift a^{21}\otimes a^{22}\dots a^{24}\otimes \shift a^{25}\otimes \shift \partial(a^3)}_g
\end{align*}
and
\begin{align*}
&\mu_1^{\otimes}(-1\otimes \shift (a^1a^{21})\otimes a^{22}\otimes \shift a^{23}\otimes a^{24}\otimes \shift a^{25}\otimes \shift a^3\otimes 1)\\
=\,&-\underbrace{\hat{\omega}a^1a^{21}a^{22}\otimes \shift a^{23}\otimes a^{24}\otimes \shift a^{25}\otimes \shift a^3\otimes 1}_{38}
+\underbrace{a^1a^{21}\hat{\omega}a^{22}\otimes \shift a^{23}\otimes a^{24}\otimes \shift a^{25}\otimes \shift a^3\otimes 1}_{39}\\
&+\underbrace{\shift \partial(a^1)a^{21}a^{22}\otimes \shift a^{23}\otimes a^{24}\otimes \shift a^{25}\otimes \shift a^3\otimes 1}_g
+\underbrace{a^1\shift \partial(a^{21})a^{22}\otimes \shift a^{23}\otimes a^{24}\otimes \shift a^{25}\otimes \shift a^3\otimes 1}_y\\
&\boldsymbol{\uwave{+}}0\boldsymbol{\uwave{+}}\underbrace{1\otimes \shift (a^1a^{21})\otimes a^{22}\hat{\omega}a^{23}a^{24}\otimes \shift a^{25}\otimes \shift a^3\otimes 1}_{40}
-\underbrace{1\otimes \shift (a^1a^{21})\otimes a^{22}a^{23}\hat{\omega}a^{24}\otimes \shift a^{25}\otimes \shift a^3\otimes 1}_{41}\\
&-\underbrace{1\otimes \shift (a^1a^{21})\otimes a^{22}\shift \partial(a^{23})a^{24}\otimes \shift a^{25}\otimes \shift a^3\otimes 1}_y
\boldsymbol{\uwave{+}}0\boldsymbol{\uwave{-}}\underbrace{1\otimes \shift (a^1a^{21})\otimes a^{22}\otimes \shift a^{23}\otimes a^{24}a^{25}\otimes \shift a^3\otimes 1}_{60}\\
&+\underbrace{1\otimes \shift (a^1a^{21})\otimes a^{22}\otimes \shift a^{23}\otimes a^{24}\otimes \shift (a^{25}a^3)\otimes 1}_y
-\underbrace{1\otimes \shift (a^1a^{21})\otimes a^{22}\otimes \shift a^{23}\otimes a^{24}\otimes \shift a^{25}\otimes a^3}_{56}
\end{align*}
and
\begin{align*}
&\mu_1^{\otimes}(-1\otimes \shift (a^1a^{21})\otimes a^{22}a^{23}a^{241}\dots a^{244}\otimes \shift a^{245}\otimes \shift a^{25}\otimes \shift a^3\otimes 1)\\
=\,&-\underbrace{\hat{\omega}a^1a^{21}\dots a^{23}a^{241}\dots a^{244}\otimes \shift a^{245}\otimes \shift a^{25}\otimes \shift a^3\otimes 1}_{42}\\
&+\underbrace{a^1a^{21}\hat{\omega}a^{22}a^{23}a^{241}\dots a^{244}\otimes \shift a^{245}\otimes \shift a^{25}\otimes \shift a^3\otimes 1}_{43}\\
&+\underbrace{\shift \partial(a^1)a^{21}\dots a^{23}a^{241}\dots a^{244}\otimes \shift a^{245}\otimes \shift a^{25}\otimes \shift a^3\otimes 1}_{gg}\\
&+\underbrace{a^1\shift \partial(a^{21})a^{22}a^{23}a^{241}\dots a^{244}\otimes \shift a^{245}\otimes \shift a^{25}\otimes \shift a^3\otimes 1}_g\\
&\boldsymbol{\uwave{+}}0\boldsymbol{\uwave{+}}0\boldsymbol{\uwave{+}}0\boldsymbol{\uwave{-}}\underbrace{1\otimes \shift (a^1a^{21})\otimes a^{22}a^{23}(\shift \otimes \shift)\partial(\shift^{-1}a^{24})\otimes \shift a^{25}\otimes \shift a^3\otimes 1}_{44}\\
&-\underbrace{1\otimes \shift (a^1a^{21})\otimes a^{22}a^{23}a^{241}\dots a^{244}\otimes \shift (a^{245}a^{25})\otimes \shift a^3\otimes 1}_y\\
&+\underbrace{1\otimes \shift (a^1a^{21})\otimes a^{22}a^{23}a^{241}\dots a^{244}\otimes \shift a^{245}\otimes \shift (a^{25}a ^3)\otimes 1}_g\\
&-\underbrace{1\otimes \shift (a^1a^{21})\otimes a^{22}a^{23}a^{241}\dots a^{244}\otimes \shift a^{245}\otimes \shift a^{25}\otimes a^3}_{100}
\end{align*}
and
\begin{align*}
&\mu_1^{\otimes}(a^1a^{21}a^{22}\otimes \shift (a^{23}a^{241})\otimes a^{242}\dots a^{244}\otimes \shift a^{245}\otimes \shift a^{25}\otimes \shift a^3\otimes 1)\\
=\,&0\boldsymbol{\uwave{+}}\underbrace{a^1a^{21}a^{22}\hat{\omega}a^{23}a^{241}\dots a^{244}\otimes \shift a^{245}\otimes \shift a^{25}\otimes \shift a^3\otimes 1}_{45}\\
&-\underbrace{a^1a^{21}a^{22}a^{23}a^{241}\hat{\omega}a^{242}\dots a^{244}\otimes \shift a^{245}\otimes \shift a^{25}\otimes \shift a^3\otimes 1}_{46}\\
&-\underbrace{a^1a^{21}a^{22}\shift \partial(a^{23})a^{241}\dots a^{244}\otimes \shift a^{245}\otimes \shift a^{25}\otimes \shift a^3\otimes 1}_g\\
&-\underbrace{a^1a^{21}\dots a^{23}\shift \partial(a^{241})a^{242}\dots a^{244}\otimes \shift a^{245}\otimes \shift a^{25}\otimes \shift a^3\otimes 1}_y\\
&\boldsymbol{\uwave{+}}0\boldsymbol{\uwave{+}}0\boldsymbol{\uwave{-}}\underbrace{a^1a^{21}a^{22}\otimes \shift (a^{23}a^{241})\otimes a^{242}\dots a^{245}\otimes \shift a^{25}\otimes \shift a^3\otimes 1}_{95}\\
&+\underbrace{a^1a^{21}a^{22}\otimes \shift (a^{23}a^{241})\otimes a^{242}\dots a^{244}\otimes \shift (a^{245}a^{25})\otimes \shift a^3\otimes 1}_y\\
&-\underbrace{a^1a^{21}a^{22}\otimes \shift (a^{23}a^{241})\otimes a^{242}\dots a^{244}\otimes \shift a^{245}\otimes \shift (a^{25}a^3)\otimes 1}_g\\
&+\underbrace{a^1a^{21}a^{22}\otimes \shift (a^{23}a^{241})\otimes a^{242}\dots a^{244}\otimes \shift a^{245}\otimes \shift a^{25}\otimes a^3}_{47}
\end{align*}
and
\begin{align*}
&\mu_1^{\otimes}(-1\otimes \shift a^1\otimes \shift a^{21}\otimes a^{22}\otimes \shift a^{23}\otimes a^{24}a^{25}\otimes \shift a^3\otimes 1)\\
=\,&-\underbrace{a^1\otimes \shift a^{21}\otimes a^{22}\otimes \shift a^{23}\otimes a^{24}a^{25}\otimes \shift a^3\otimes 1}_{48}
+\underbrace{1\otimes \shift (a^1a^{21})\otimes a^{22}\otimes \shift a^{23}\otimes a^{24}a^{25}\otimes \shift a^3\otimes 1}_{60}\\
&-\underbrace{1\otimes \shift a^1\otimes a^{21}a^{22}\otimes \shift a^{23}\otimes a^{24}a^{25}\otimes \shift a^3\otimes 1}_{49}
\boldsymbol{\uwave{+}}0\boldsymbol{\uwave{-}}\underbrace{1\otimes \shift a^1\otimes \shift a^{21}\otimes a^{22}\hat{\omega}a^{23}a^{24}a^{25}\otimes \shift a^3\otimes 1}_{50}\\
&+\underbrace{1\otimes \shift a^1\otimes \shift a^{21}\otimes a^{22}a^{23}\hat{\omega}a^{24}a^{25}\otimes \shift a^3\otimes 1}_{51}
+\underbrace{1\otimes \shift a^1\otimes \shift a^{21}\otimes a^{22}\shift \partial(a^{23})a^{24}a^{25}\otimes \shift a^3\otimes 1}_y\\
&\boldsymbol{\uwave{+}}0\boldsymbol{\uwave{+}}\underbrace{1\otimes \shift a^1\otimes \shift a^{21}\otimes a^{22}\otimes \shift a^{23}\otimes a^{24}a^{25}\hat{\omega}a^3}_{55}
-\underbrace{1\otimes \shift a^1\otimes \shift a^{21}\otimes a^{22}\otimes \shift a^{23}\otimes a^{24}a^{25}a^3\hat{\omega}}_{52}\\
&-\underbrace{1\otimes \shift a^1\otimes \shift a^{21}\otimes a^{22}\otimes \shift a^{23}\otimes a^{24}a^{25}\shift \partial(a^3)}_g
\end{align*}
and
\begin{align*}
&\mu_1^{\otimes}(-1\otimes \shift a^1\otimes \shift a^{21}\otimes a^{22}\otimes \shift a^{23}\otimes a^{24}\otimes \shift a^{25}\otimes a^3)\\
=\,&-\underbrace{a^1\otimes \shift a^{21}\otimes a^{22}\otimes \shift a^{23}\otimes a^{24}\otimes \shift a^{25}\otimes a^3}_{53}
+\underbrace{1\otimes \shift (a^1a^{21})\otimes a^{22}\otimes \shift a^{23}\otimes a^{24}\otimes \shift a^{25}\otimes a^3}_{56}\\
&-\underbrace{1\otimes \shift a^1\otimes a^{21}a^{22}\otimes \shift a^{23}\otimes a^{24}\otimes \shift a^{25}\otimes a^3}_{58}
\boldsymbol{\uwave{+}}0\boldsymbol{\uwave{-}}\underbrace{1\otimes \shift a^1\otimes \shift a^{21}\otimes a^{22}\hat{\omega}a^{23}a^{24}\otimes \shift a^{25}\otimes a^3}_{59}\\
&+\underbrace{1\otimes \shift a^1\otimes \shift a^{21}\otimes a^{22}a^{23}\hat{\omega}a^{24}\otimes \shift a^{25}\otimes a^3}_{57}
+\underbrace{1\otimes \shift a^1\otimes \shift a^{21}\otimes a^{22}\shift \partial(a^{23})a^{24}\otimes \shift a^{25}\otimes a^3}_y\\
&\boldsymbol{\uwave{+}}0\boldsymbol{\uwave{+}}\underbrace{1\otimes \shift a^1\otimes \shift a^{21}\otimes a^{22}\otimes \shift a^{23}\otimes a^{24}\hat{\omega}a^{25}a^3}_{54}
-\underbrace{1\otimes \shift a^1\otimes \shift a^{21}\otimes a^{22}\otimes \shift a^{23}\otimes a^{24}a^{25}\hat{\omega}a^3}_{55}\\
&-\underbrace{1\otimes \shift a^1\otimes \shift a^{21}\otimes a^{22}\otimes \shift a^{23}\otimes a^{24}\shift \partial(a^{25})a^3}_y
\end{align*}
and
\begin{align*}
&\mu_1^{\otimes}(-a^1a^{21}\dots a^{23}a^{241}a^{242}\otimes \shift a^{243}\otimes a^{244}\otimes \shift a^{245}\otimes \shift a^{25}\otimes \shift a^3\otimes 1)\\
=\,&0\boldsymbol{\uwave{+}}0\boldsymbol{\uwave{-}}\underbrace{a^1a^{21}\dots a^{23}a^{241}a^{242}\hat{\omega}a^{243}a^{244}\otimes \shift a^{245}\otimes \shift a^{25}\otimes \shift a^3\otimes 1}_{61}\\
&+\underbrace{a^1a^{21}\dots a^{23}a^{241}\dots a^{243}\hat{\omega}a^{244}\otimes \shift a^{245}\otimes \shift a^{25}\otimes \shift a^3\otimes 1}_{62}\\
&+\underbrace{a^1a^{21}\dots a^{23}a^{241}a^{242}\shift \partial(a^{243})a^{244}\otimes \shift a^{245}\otimes \shift a^{25}\otimes \shift a^3\otimes 1}_y\\
&\boldsymbol{\uwave{+}}0\boldsymbol{\uwave{+}}\underbrace{a^1a^{21}\dots a^{23}a^{241}a^{242}\otimes \shift a^{243}\otimes a^{244}a^{245}\otimes \shift a^{25}\otimes \shift a^3\otimes 1}_{63}\\
&-\underbrace{a^1a^{21}\dots a^{23}a^{241}a^{242}\otimes \shift a^{243}\otimes a^{244}\otimes \shift (a^{245}a^{25})\otimes \shift a^3\otimes 1}_y\\
&+\underbrace{a^1a^{21}\dots a^{23}a^{241}a^{242}\otimes \shift a^{243}\otimes a^{244}\otimes \shift a^{245}\otimes \shift (a^{25}a^3)\otimes 1}_g\\
&-\underbrace{a^1a^{21}a^{22}a^{23}a^{241}a^{242}\otimes \shift a^{243}\otimes a^{244}\otimes \shift a^{245}\otimes \shift a^{25}\otimes a^3}_{64}
\end{align*}
and
\begin{align*}
&\mu_1^{\otimes}(1\otimes \shift a^1\otimes \shift a^{21}\otimes \shift a^{221}\otimes a^{222}\dots a^{224}\otimes \shift (a^{225}a^{23})\otimes a^{24}a^{25}a^3)\\
=\,&\underbrace{a^1\otimes \shift a^{21}\otimes \shift a^{221}\otimes a^{222}\dots a^{224}\otimes \shift (a^{225}a^{23})\otimes a^{24}a^{25}a^3}_{65}\\
&-\underbrace{1\otimes \shift (a^1a^{21})\otimes \shift a^{221}\otimes a^{222}\dots a^{224}\otimes \shift (a^{225}a^{23})\otimes a^{24}a^{25}a^3}_{86}\\
&+\underbrace{1\otimes \shift a^1\otimes \shift (a^{21}a^{221})\otimes a^{222}\dots a^{224}\otimes \shift (a^{225}a^{23})\otimes a^{24}a^{25}a^3}_y\\
&-\underbrace{1\otimes \shift a^1\otimes \shift a^{21}\otimes a^{221}\dots a^{224}\otimes \shift(a^{225}a^{23})\otimes a^{24}a^{25}a^3}_{92}\\
&\boldsymbol{\uwave{+}}0\boldsymbol{\uwave{+}}0\boldsymbol{\uwave{-}}\underbrace{1\otimes \shift a^1\otimes \shift a^{21}\otimes \shift a^{221}\otimes a^{222}\dots a^{224}\hat{\omega}a^{225}a^{23}\dots a^{25}a^3}_{66}\\
&+\underbrace{1\otimes \shift a^1\otimes \shift a^{21}\otimes \shift a^{221}\otimes a^{222}\dots a^{225}a^{23}\hat{\omega}a^{24}a^{25}a^3}_{67} \\
&+\underbrace{1\otimes \shift a^1\otimes \shift a^{21}\otimes \shift a^{221}\otimes a^{222}\dots a^{224}\shift \partial(a^{225})a^{23}\dots a^{25}a^3}_y\\
&+\underbrace{1\otimes \shift a^1\otimes \shift a^{21}\otimes \shift a^{221}\otimes a^{222}\dots a^{225}\shift \partial(a^{23})a^{24}a^{25}a^3}_g\boldsymbol{\uwave{+}}0
\end{align*}
and
\begin{align*}
&\mu_1^{\otimes}(-1\otimes \shift a^1\otimes \shift a^{21}\otimes \shift a^{221}\otimes a^{222}\dots a^{225}a^{23}a^{24}\otimes \shift (a^{25}a^3)\otimes 1)\\
=\,&-\underbrace{a^1\otimes \shift a^{21}\otimes \shift a^{221}\otimes a^{222}\dots a^{225}a^{23}a^{24}\otimes \shift (a^{25}a^3)\otimes 1}_{74}\\
&+\underbrace{1\otimes \shift(a^1a^{21})\otimes \shift a^{221}\otimes a^{222}\dots a^{225}a^{23}a^{24}\otimes \shift(a^{25}a^3)\otimes 1}_{g}\\
&-\underbrace{1\otimes \shift a^1\otimes \shift (a^{21}a^{221})\otimes a^{222}\dots a^{225}a^{23}a^{24}\otimes \shift (a^{25}a^3)\otimes 1}_y\\
&-\underbrace{1\otimes \shift a^1\otimes \shift a^{21}\otimes (\shift \otimes \shift)\partial\shift^{-1}(a^{22})a^{23}a^{24}\otimes \shift (a^{25}a^3)\otimes 1}_{111}\\
&\boldsymbol{\uwave{+}}0\boldsymbol{\uwave{+}}0\boldsymbol{\uwave{+}}0\boldsymbol{\uwave{+}}\underbrace{1\otimes \shift a^1\otimes \shift a^{21}\otimes \shift a^{221}\otimes a^{222}\dots a^{225}a^{23}a^{24}\hat{\omega}a^{25}a^3}_{68}\\
&-\underbrace{1\otimes \shift a^1\otimes \shift a^{21}\otimes \shift a^{221}\otimes a^{222}\dots a^{225}a^{23}\dots a^{25}a^3\hat{\omega}}_{69}\\
&-\underbrace{1\otimes \shift a^1\otimes \shift a^{21}\otimes \shift a^{221}\otimes a^{222}\dots a^{225}a^{23}a^{24}\shift \partial(a^{25})a^3}_g\\
&-\underbrace{1\otimes \shift a^1\otimes \shift a^{21}\otimes \shift a^{221}\otimes a^{222}\dots a^{225}a^{23}\dots a^{25}\shift \partial(a^3)}_{gg}
\end{align*}
and
\begin{align*}
&\mu_1^{\otimes}(-1\otimes \shift a^1\otimes \shift a^{21}\otimes \shift a^{221}\otimes a^{222}\otimes \shift a^{223}\otimes a^{224}a^{225}a^{23}\dots a^{25}a^3\\
=\,&-\underbrace{a^1\otimes \shift a^{21}\otimes \shift a^{221}\otimes a^{222}\otimes \shift a^{223}\otimes a^{224}a^{225}a^{23}a^{24} a^{25}a^3}_{70}\\
&+\underbrace{1\otimes \shift (a^1a^{21})\otimes \shift a^{221}\otimes a^{222}\otimes \shift a^{223}\otimes a^{224}a^{225}a^{23}\dots a^{25}a^3}_g\\
&-\underbrace{1\otimes \shift a^1\otimes \shift (a^{21}a^{221})\otimes a^{222}\otimes \shift a^{223}\otimes a^{224}a^{225}a^{23}\dots a^{25}a^3}_y\\
&+\underbrace{1\otimes \shift a^1\otimes \shift a^{21}\otimes a^{221}a^{222}\otimes \shift a^{223}\otimes a^{224}a^{225}a^{23}\dots a^{25}a^3}_{73}\\
&\boldsymbol{\uwave{+}}0\boldsymbol{\uwave{+}}\underbrace{1\otimes \shift a^1\otimes \shift a^{21}\otimes \shift a^{221}\otimes a^{222}\hat{\omega}a^{223}\dots a^{225}a^{23}\dots a^{25}a^3}_{71}\\
&-\underbrace{1\otimes \shift a^1\otimes \shift a^{21}\otimes \shift a^{221}\otimes a^{222}a^{223}\hat{\omega}a^{224}a^{225}a^{23}\dots a^{25}a^3}_{72}\\
&-\underbrace{1\otimes \shift a^1\otimes \shift a^{21}\otimes \shift a^{221}\otimes a^{222}\shift \partial(a^{223})a^{224}a^{225}a^{23}\dots a^{25}a^3}_y\boldsymbol{\uwave{+}}0\boldsymbol{\uwave{+}}0
\end{align*}
and
\begin{align*}
&\mu_1^{\otimes}(-a^1\otimes \shift a^{21}\otimes \shift a^{221}\otimes a^{222}\dots a^{225}a^{23}a^{24}\otimes \shift a^{25}\otimes \shift a^3\otimes 1)\\
=\,&-\underbrace{a^1a^{21}\otimes \shift a^{221}\otimes a^{222}\dots a^{225}a^{23}a^{24}\otimes \shift a^{25}\otimes \shift a^3\otimes 1}_{76}\\
&
+\underbrace{a^1\otimes \shift (a^{21}a^{221})\otimes a^{222}\dots a^{225}a^{23}a^{24}\otimes \shift a^{25}\otimes \shift a^3\otimes 1}_y\\
&+\underbrace{a^1\otimes \shift a^{21}\otimes (\shift \otimes \shift) \partial(\shift^{-1}a^{22})a^{23}a^{24}\otimes \shift a^{25}\otimes \shift a^3\otimes 1}_{79}\\
&
\boldsymbol{\uwave{+}}0\boldsymbol{\uwave{+}}0\boldsymbol{\uwave{+}}0\boldsymbol{\uwave{-}}\underbrace{a^1\otimes \shift a^{21}\otimes \shift a^{221}\otimes a^{222}\dots a^{225}a^{23}\dots a^{25}\otimes \shift a^3\otimes 1}_{75}\\
&+\underbrace{a^1\otimes \shift a^{21}\otimes \shift a^{221}\otimes a^{222}\dots a^{225}a^{23}a^{24}\otimes \shift(a^{25}a^3)\otimes 1}_{74}\\
&
-\underbrace{a^1\otimes \shift a^{21}\otimes \shift a^{221}\otimes a^{222}\dots a^{225}a^{23}a^{24}\otimes \shift a^{25}\otimes a^3}_{77}
\end{align*}
and
\begin{align*}
&\mu_1^{\otimes}(-a^1a^{21}a^{221}\dots a^{224}\otimes \shift a^{225}\otimes \shift a^{23}\otimes a^{24}\otimes \shift a^{25}\otimes \shift a^3\otimes 1)\\
=\,&0\boldsymbol{\uwave{+}}0\boldsymbol{\uwave{+}}\underbrace{a^1a^{21}(\shift \otimes \shift)\partial(\shift^{-1}a^{22})\otimes \shift a^{23}\otimes a^{24}\otimes \shift a^{25}\otimes \shift a^3\otimes 1}_{82}\\
&
+\underbrace{a^1a^{21}a^{221}\dots a^{224}\otimes \shift (a^{225}a^{23})\otimes a^{24}\otimes \shift a^{25}\otimes \shift a^3\otimes 1}_y\\
&-\underbrace{a^1a^{21}a^{221}\dots a^{224}\otimes \shift a^{225}\otimes a^{23}a^{24}\otimes \shift a^{25}\otimes \shift a^3\otimes 1}_{85}\\
&
\boldsymbol{\uwave{+}}0\boldsymbol{\uwave{-}}\underbrace{a^1a^{21}a^{221}\dots a^{224}\otimes \shift a^{225}\otimes \shift a^{23}\otimes a^{24}a^{25}\otimes \shift a^3\otimes 1}_{83}\\
&+\underbrace{a^1a^{21}a^{221}\dots a^{224}\otimes \shift a^{225}\otimes \shift a^{23}\otimes a^{24}\otimes \shift (a^{25}a^3)\otimes 1}_g\\
&
-\underbrace{a^1a^{21}a^{221}\dots a^{224}\otimes \shift a^{225}\otimes \shift a^{23}\otimes a^{24}\otimes \shift a^{25}\otimes a^3}_{84}
\end{align*}
and
\begin{align*}
&\mu_1^{\otimes}(-1\otimes \shift a^1\otimes \shift a^{21}\otimes (\shift \otimes \shift)\partial(\shift^{-1}a^{22})a^{23}a^{24}\otimes \shift a^{25}\otimes \shift a^3\otimes 1)\\
=\,&-\underbrace{a^1\otimes \shift a^{21}\otimes (\shift \otimes \shift)\partial(\shift^{-1}a^{22})a^{23}a^{24}\otimes \shift a^{25}\otimes \shift a^3\otimes 1}_{79}\\
&
+\underbrace{1\otimes \shift (a^1a^{21})\otimes (\shift \otimes \shift)\partial(\shift^{-1}a^{22})a^{23}a^{24}\otimes \shift a^{25}\otimes \shift a^3\otimes 1}_{78}\\
&-\underbrace{1\otimes \shift a^1\otimes a^{21}(\shift \otimes \shift)\partial(\shift^{-1}a^{22})a^{23}a^{24}\otimes \shift a^{25}\otimes \shift a^3\otimes 1}_y\\
&\boldsymbol{\uwave{+}}0\boldsymbol{\uwave{+}}0\boldsymbol{\uwave{+}}0
\boldsymbol{\uwave{-}}\underbrace{1\otimes \shift a^1\otimes \shift a^{21}\otimes (\shift \otimes \shift)\partial(\shift^{-1}a^{22})a^{23}\dots a^{25}\otimes \shift a^3\otimes 1}_{80}\\
&+\underbrace{1\otimes \shift a^1\otimes \shift a^{21}\otimes (\shift \otimes \shift)\partial(\shift^{-1}a^{22})a^{23}a^{24}\otimes \shift(a^{25}a^3)\otimes 1}_{111}\\
&
-\underbrace{1\otimes \shift a^1\otimes \shift a^{21}\otimes (\shift \otimes \shift)\partial(\shift^{-1}a^{22})a^{23}a^{24}\otimes \shift a^{25}\otimes a^3}_{81}
\end{align*}
and
\begin{align*}
&\mu_1^{\otimes}(-1\otimes \shift (a^1a^{21})\otimes \shift a^{221}\otimes a^{222}\dots a^{224}\otimes \shift a^{225}\otimes \shift a^{23}\otimes a^{24}a^{25}a^3)\\
=\,&-\underbrace{a^1a^{21}\otimes \shift a^{221}\otimes a^{222}\dots a^{224}\otimes \shift a^{225}\otimes \shift a^{23}\otimes a^{24}a^{25}a^3}_{87}\\
&+\underbrace{1\otimes \shift (a^1a^{21}a^{221})\otimes a^{222}\dots a^{224}\otimes \shift a^{225}\otimes \shift a^{23}\otimes a^{24}a^{25}a^3}_g\\
&-\underbrace{1\otimes \shift (a^1a^{21})\otimes a^{221}\dots a^{224}\otimes \shift a^{225}\otimes \shift a^{23}\otimes a^{24}a^{25}a^3}_g\\
&\boldsymbol{\uwave{+}}0\boldsymbol{\uwave{+}}0\boldsymbol{\uwave{-}}\underbrace{1\otimes \shift (a^1a^{21})\otimes \shift a^{221}\otimes a^{222}\dots a^{225}\otimes \shift a^{23}\otimes a^{24}a^{25}a^3}_{g}\\
&+\underbrace{1\otimes \shift (a^1a^{21})\otimes \shift a^{221}\otimes a^{222}\dots a^{224}\otimes \shift (a^{225}a^{23})\otimes a^{24}a^{25}a^3}_{86}\\
&-\underbrace{1\otimes \shift (a^1a^{21})\otimes \shift a^{221}\otimes a^{222}\dots a^{224}\otimes \shift a^{225}\otimes a^{23}\dots a^{25}a^3}_{g}\boldsymbol{\uwave{+}}0
\end{align*}
and
\begin{align*}
&\mu_1^{\otimes}(-1\otimes \shift a^1\otimes \shift a^{21}\otimes a^{221}\dots a^{224}\otimes \shift a^{225}\otimes \shift a^{23}\otimes a^{24}a^{25}a^3)\\
=\,&-\underbrace{a^1\otimes \shift a^{21}\otimes a^{221}\dots a^{224}\otimes \shift a^{225}\otimes \shift a^{23}\otimes a^{24}a^{25}a^3}_{89}\\
&
+\underbrace{1\otimes \shift (a^1a^{21})\otimes a^{221}\dots a^{224}\otimes \shift a^{225}\otimes \shift a^{23}\otimes a^{24}a^{25}a^3}_g\\
&-\underbrace{1\otimes \shift a^1\otimes a^{21}a^{221}\dots a^{224}\otimes \shift a^{225}\otimes \shift a^{23}\otimes a^{24}a^{25}a^3}_{90}\\
&\boldsymbol{\uwave{+}}0\boldsymbol{\uwave{+}}0
\boldsymbol{\uwave{+}}\underbrace{1\otimes \shift a^1\otimes \shift a^{21}\otimes (\shift \otimes \shift)\partial(\shift^{-1}a^{22})\otimes \shift a^{23}\otimes a^{24}a^{25}a^3}_{91}\\
&+\underbrace{1\otimes \shift a^1\otimes \shift a^{21}\otimes a^{221}\dots a^{224}\otimes \shift (a^{225}a^{23})\otimes a^{24}a^{25}a^3}_{92}\\
&
-\underbrace{1\otimes \shift a^1\otimes \shift a^{21}\otimes a^{221}\dots a^{224}\otimes \shift a^{225}\otimes a^{23}\dots a^{25}a^3}_{93}\boldsymbol{\uwave{+}}0
\end{align*}
and
\begin{align*}
&\mu_1^{\otimes}(-a^1a^{21}a^{22}\otimes \shift a^{23}\otimes \shift a^{241}\otimes a^{242}\dots a^{245}\otimes \shift a^{25}\otimes \shift a^3\otimes 1)\\
=\,&0\boldsymbol{\uwave{-}}\underbrace{a^1a^{21}\dots a^{23}\otimes \shift a^{241}\otimes a^{242}\dots a^{245}\otimes \shift a^{25}\otimes \shift a^3\otimes 1}_{94}\\
&
+\underbrace{a^1a^{21}a^{22}\otimes \shift (a^{23}a^{241})\otimes a^{242}\dots a^{245}\otimes \shift a^{25}\otimes \shift a^3\otimes 1}_{95}\\
&+\underbrace{a^1a^{21}a^{22}\otimes \shift a^{23}\otimes (\shift \otimes \shift)\partial(\shift^{-1}a^{24})\otimes \shift a^{25}\otimes \shift a^3\otimes 1}_{96}\\
&\boldsymbol{\uwave{+}}0\boldsymbol{\uwave{+}}0
\boldsymbol{\uwave{-}}\underbrace{a^1a^{21}a^{22}\otimes \shift a^{23}\otimes \shift a^{241}\otimes a^{242}\dots a^{245}a^{25}\otimes \shift a^3\otimes 1}_{98}\\
&+\underbrace{a^1a^{21}a^{22}\otimes \shift a^{23}\otimes \shift a^{241}\otimes a^{242}\dots a^{245}\otimes \shift (a^{25}a^3)\otimes 1}_g\\
&
-\underbrace{a^1a^{21}a^{22}\otimes \shift a^{23}\otimes \shift a^{241}\otimes a^{242}\dots a^{245}\otimes \shift a^{25}\otimes a^3}_{97}
\end{align*}
and
\begin{align*}
&\mu_1^{\otimes}(-1\otimes \shift a^1\otimes \shift a^{21}\otimes a^{22}a^{23}a^{241}\dots a^{244}\otimes \shift a^{245}\otimes \shift a^{25}\otimes a^3)\\
=\,&-\underbrace{a^1\otimes \shift a^{21}\otimes a^{22}a^{23}a^{241}\dots a^{244}\otimes \shift a^{245}\otimes \shift a^{25}\otimes a^3}_{99}\\
&
+\underbrace{1\otimes \shift (a^1a^{21})\otimes a^{22}a^{23}a^{241}\dots a^{244}\otimes \shift a^{245}\otimes \shift a^{25}\otimes a^3}_{100}\\
&-\underbrace{1\otimes \shift a^1\otimes a^{21}\dots a^{23}a^{241}\dots a^{244}\otimes \shift a^{245}\otimes \shift a^{25}\otimes a^3}_{101}\\
&\boldsymbol{\uwave{+}}0\boldsymbol{\uwave{+}}0\boldsymbol{\uwave{+}}0
\boldsymbol{\uwave{+}}\underbrace{1\otimes \shift a^1\otimes \shift a^{21}\otimes a^{22}a^{23}(\shift \otimes \shift)\partial(\shift^{-1}a^{24})\otimes \shift a^{25}\otimes a^3}_{102}\\
&+\underbrace{1\otimes \shift a^1\otimes \shift a^{21}\otimes a^{22}a^{23}a^{241}\dots a^{244}\otimes \shift (a^{245}a^{25})\otimes a^3}_y\\
&
-\underbrace{1\otimes \shift a^1\otimes \shift a^{21}\otimes a^{22}a^{23}a^{241}\dots a^{244}\otimes \shift a^{245}\otimes a^{25}a^3}_{103}
\end{align*}
and
\begin{align*}
&\mu_1^{\otimes}(-1\otimes \shift a^1\otimes \shift a^{21}\otimes a^{22}\otimes \shift a^{23}\otimes \shift a^{241}\otimes a^{242}\dots a^{245}a^{25}a^3)\\
=\,&-\underbrace{a^1\otimes \shift a^{21}\otimes a^{22}\otimes \shift a^{23}\otimes \shift a^{241}\otimes a^{242}\dots a^{245}a^{25}a^3}_{104}\\
&
+\underbrace{1\otimes \shift (a^1a^{21})\otimes a^{22}\otimes \shift a^{23}\otimes \shift a^{241}\otimes a^{242}\dots a^{245}a^{25}a^3}_g\\
&-\underbrace{1\otimes \shift a^1\otimes a^{21}a^{22}\otimes \shift a^{23}\otimes \shift a^{241}\otimes a^{242}\dots a^{245}a^{25}a^3}_{105}\\
&\boldsymbol{\uwave{+}}0
\boldsymbol{\uwave{-}}\underbrace{1\otimes \shift a^1\otimes \shift a^{21}\otimes a^{22}a^{23}\otimes \shift a^{241}\otimes a^{242}\dots a^{245}a^{25}a^3}_{106}\\
&+\underbrace{1\otimes \shift a^1\otimes \shift a^{21}\otimes a^{22}\otimes \shift (a^{23}a^{241})\otimes a^{242}\dots a^{245}a^{25}a^3}_y\\
&
+\underbrace{1\otimes \shift a^1\otimes \shift a^{21}\otimes a^{22}\otimes \shift a^{23}\otimes (\shift \otimes \shift)\partial(\shift^{-1}a^{24})a^{25}a^3}_{107}\boldsymbol{\uwave{+}}0\boldsymbol{\uwave{+}}0
\end{align*}
and
\begin{align*}
&\mu_1^{\otimes}(-1\otimes \shift a^1\otimes \shift a^{21}\otimes \shift a^{221}\otimes \shift a^{2221}\otimes a^{2222}\dots a^{2225}a^{223}\dots a^{225}a^{23}\dots a^{25}a^3)\\
=\,&-\underbrace{a^1\otimes \shift a^{21}\otimes \shift a^{221}\otimes \shift a^{2221}\otimes a^{2222}\dots a^{2225}a^{223}\dots a^{225}a^{23}\dots a^{25}a^3}_{108}\\
&+\underbrace{1\otimes \shift (a^1a^{21})\otimes \shift a^{221}\otimes \shift a^{2221}\otimes a^{2222}\dots a^{2225}a^{223}\dots a^{225}a^{23}\dots a^{25}a^3}_{gg}\\
&-\underbrace{1\otimes \shift a^1\otimes \shift (a^{21}a^{221})\otimes \shift a^{2221}\otimes a^{2222}\dots a^{2225}a^{223}\dots a^{225}a^{23}\dots a^{25}a^3}_g\\
&+\underbrace{1\otimes \shift a^1\otimes \shift a^{21}\otimes \shift (a^{221}a^{2221})\otimes a^{2222}\dots a^{2225}a^{223}\dots a^{225}a^{23}\dots a^{25}a^3}_y\\
&+\underbrace{1\otimes \shift a^1\otimes \shift a^{21}\otimes \shift a^{221}\otimes (\shift \otimes \shift)\partial(\shift^{-1}a^{222})a^{223}\dots a^{225}a^{23}\dots a^{25}a^3}_{109}\boldsymbol{\uwave{+}}0\boldsymbol{\uwave{+}}0\boldsymbol{\uwave{+}}0\boldsymbol{\uwave{+}}0
\end{align*}
and
\begin{align*}
&\mu_1^{\otimes}(-a^1a^{21}\dots a^{23}a^{241}\dots a^{243}a^{2441}\dots a^{2444}\otimes \shift a^{2445}\otimes \shift a^{245}\otimes \shift a^{25}\otimes \shift a^3\otimes 1)\\
=\,&0\boldsymbol{\uwave{+}}0\boldsymbol{\uwave{+}}0\boldsymbol{\uwave{+}}0\boldsymbol{\uwave{+}}\underbrace{a^1a^{21}\dots a^{23}a^{241}\dots a^{243}(\shift \otimes \shift)\partial(\shift^{-1}a^{244})\otimes \shift a^{245}\otimes \shift a^{25}\otimes \shift a^3\otimes 1}_{110}\\
&+\underbrace{a^1a^{21}\dots a^{23}a^{241}\dots a^{243}a^{2441}\dots a^{2444}\otimes \shift (a^{2445}a^{245})\otimes \shift a^{25}\otimes \shift a^3\otimes 1}_y\\
&-\underbrace{a^1a^{21}\dots a^{23}a^{241}\dots a^{243}a^{2441}\dots a^{2444}\otimes \shift a^{2445}\otimes \shift (a^{245}a^{25})\otimes \shift a^3\otimes 1}_g\\
&+\underbrace{a^1a^{21}\dots a^{23}a^{241}\dots a^{243}a^{2441}\dots a^{2444}\otimes \shift a^{2445}\otimes \shift a^{245}\otimes \shift (a^{25} a^3)\otimes 1}_{gg}\\
&-\underbrace{a^1a^{21}a^{22} a^{23}a^{241}\dots a^{243}a^{2441}\dots a^{2444}\otimes \shift a^{2445}\otimes \shift a^{245}\otimes \shift a^{25}\otimes a^3}_{88}
\end{align*}
Noting that terms with the same label add up to zero, this finishes the proof for $\mu_5^1$. 
\end{proof}

\section{The map $\mu_6^0$ on $\mathcal{P}$}\label{appendixE}

Finally, in this section we provide a formula for $\mu_6^0$, which is the final higher comultiplication we need for the proof of the main theorem. 

\begin{lem}
For $\varphi\in \overline{V}$ we have:
\begin{align*}
\mu_6^0(\varphi)
=\,&-1\otimes \shift \varphi^1\otimes \shift \varphi^{21}\otimes \varphi^{22}\varphi^{23}\varphi^{241}\dots \varphi^{244}\otimes \shift \varphi^{245}\otimes \shift \varphi^{25}\otimes \varphi^3\varphi^4\varphi^5\\
&-\varphi^1\varphi^{21}\dots \varphi^{23}\varphi^{241}\dots \varphi^{244}\otimes \shift \varphi^{245}\otimes \shift \varphi^{25}\otimes \shift \varphi^3\otimes \shift \varphi^{41}\otimes \varphi^{42}\dots \varphi^{45}\varphi^5\\
&-\varphi^1\varphi^2\varphi^3\varphi^{41}\varphi^{421}\dots \varphi^{424}\otimes \shift \varphi^{425}\otimes \shift \varphi^{43}\otimes \varphi^{44}\otimes \shift \varphi^{45}\otimes \shift \varphi^5\otimes 1\\
&-\varphi^1\varphi^{21}\varphi^{22}\otimes \shift \varphi^{23}\otimes \shift \varphi^{241}\otimes \varphi^{242}\dots \varphi^{245}\otimes \shift \varphi^{25}\otimes \shift \varphi^3\otimes \varphi^4\varphi^5\\
&-\varphi^1\varphi^2\otimes \shift \varphi^3\otimes \shift \varphi^{41}\otimes \varphi^{421}\dots \varphi^{424}\otimes \shift \varphi^{425}\otimes \shift \varphi^{43}\otimes \varphi^{44}\varphi^{45}\varphi^5\\
&-1\otimes \shift \varphi^1\otimes \shift \varphi^{21}\otimes \varphi^{22}\otimes \shift \varphi^{23}\otimes \shift \varphi^{241}\otimes \varphi^{242}\dots \varphi^{245}\varphi^{25}\varphi^3\varphi^4\varphi^5\\
&-\varphi^1\varphi^2\otimes \shift \varphi^3\otimes \shift \varphi^{41}\otimes \varphi^{42}\varphi^{43}\varphi^{441}\dots \varphi^{444}\otimes \shift \varphi^{445}\otimes \shift \varphi^{45}\otimes \varphi^5\\
&-\varphi^1\varphi^2\varphi^3\varphi^{41}\varphi^{42}\otimes \shift \varphi^{43}\otimes \shift \varphi^{441}\otimes \varphi^{442}\dots \varphi^{445}\otimes \shift \varphi^{45}\otimes \shift \varphi^5\otimes 1\\
&-\varphi^1\varphi^{21}\dots \varphi^{23}\varphi^{241}\dots \varphi^{243}\varphi^{2441}\dots \varphi^{2444}\otimes \shift \varphi^{2445}\otimes \shift \varphi^{245}\otimes \shift \varphi^{25}\otimes \shift \varphi^3\otimes \varphi^4\varphi^5\\
&-1\otimes \shift \varphi^1\otimes \varphi^2\otimes \shift \varphi^3\otimes \shift \varphi^{41}\otimes \varphi^{42}\otimes \shift \varphi^{43}\otimes \varphi^{44}\varphi^{45}\varphi^5\\
&-\varphi^1\otimes \shift \varphi^{21}\otimes \varphi^{22}\otimes \shift \varphi^{23}\otimes \varphi^{24}\otimes \shift \varphi^{25}\otimes \shift \varphi^3\otimes \varphi^4\varphi^5\\
&+\varphi^1\otimes \shift \varphi^{21}\otimes \varphi^{22}\dots \varphi^{24}\otimes \shift \varphi^{25}\otimes \shift \varphi^3\otimes \varphi^4\otimes \shift \varphi^5\otimes 1\\
&+1\otimes \shift \varphi^1\otimes \varphi^2\otimes \shift \varphi^3\otimes \shift \varphi^{41}\otimes \varphi^{42}\dots \varphi^{44}\otimes \shift \varphi^{45}\otimes \varphi^5\\
&-\varphi^1\varphi^2\otimes \shift \varphi^3\otimes \shift \varphi^{41}\otimes \varphi^{42}\otimes \shift \varphi^{43}\otimes \varphi^{44}\otimes \shift \varphi^{45}\otimes \varphi^5\\
&-\varphi^1\varphi^2\otimes \shift (\varphi^3\varphi^{41})\otimes \varphi^{42}\otimes \shift \varphi^{43}\otimes \varphi^{44}\otimes \shift \varphi^{45}\otimes \shift \varphi^5\otimes 1\\
&-1\otimes \shift \varphi^1\otimes \shift \varphi^{21}\otimes \varphi^{22}\otimes \shift \varphi^{23}\otimes \varphi^{24}\otimes \shift (\varphi^{25}\varphi^3)\otimes \varphi^4\varphi^5\\
&-\varphi^1\varphi^{21}\varphi^{22}\otimes \shift \varphi^{23}\otimes \varphi^{24}\otimes \shift \varphi^{25}\otimes \shift \varphi^3\otimes \varphi^4\otimes \shift \varphi^5\otimes 1\\
&-\varphi^1\varphi^{21}\dots \varphi^{23}\varphi^{241}\varphi^{242}\otimes \shift \varphi^{243}\otimes \varphi^{244}\otimes \shift \varphi^{245}\otimes \shift \varphi^{25}\otimes \shift \varphi^3\otimes \varphi^4\varphi^5\\
&-\varphi^1\varphi^2\varphi^3\varphi^{41}\dots \varphi^{43}\varphi^{441}\varphi^{442}\otimes \shift \varphi^{443}\otimes \varphi^{444}\otimes \shift \varphi^{445}\otimes \shift \varphi^{45}\otimes \shift \varphi^5\otimes 1\\
&-\varphi^1(\shift \otimes \shift)\partial(\shift^{-1}\varphi^2)\otimes \shift \varphi^3\otimes \shift \varphi^{41}\otimes \varphi^{42}\dots \varphi^{44}\otimes \shift \varphi^{45}\otimes \shift \varphi^5\otimes 1\\
&-\varphi^1\varphi^{21}\dots \varphi^{24}\otimes \shift \varphi^{25}\otimes \shift \varphi^3\otimes \varphi^{41}\dots \varphi^{44}\otimes \shift \varphi^{45}\otimes \shift \varphi^5\otimes 1\\
&-\varphi^1\varphi^{21}\varphi^{221}\dots \varphi^{224}\otimes \shift \varphi^{225}\otimes \shift \varphi^{23}\otimes \varphi^{24}\otimes \shift \varphi^{25}\otimes \shift \varphi^3\otimes \varphi^4\varphi^5\\
&-1\otimes \shift \varphi^1\otimes \shift \varphi^{21}\otimes \varphi^{22}\dots \varphi^{25}\varphi^3\varphi^{41}\dots \varphi^{44}\otimes \shift \varphi^{45}\otimes \shift \varphi^5\otimes 1\\
&-\varphi^1\varphi^{21}\dots \varphi^{24}\otimes \shift \varphi^{25}\otimes \shift \varphi^3\otimes \shift \varphi^{41}\otimes \shift \varphi^{421}\otimes \varphi^{422}\dots \varphi^{425}\varphi^{43}\dots \varphi^{45}\varphi^5\\
&-1\otimes \shift \varphi^1\otimes \shift \varphi^{21}\otimes \varphi^{22}\dots \varphi^{25}\otimes \shift \varphi^3\otimes \shift \varphi^{41}\otimes \varphi^{42}\dots \varphi^{45}\varphi^5\\
&-1\otimes \shift \varphi^1\otimes \shift \varphi^{21}\otimes \varphi^{22}\dots \varphi^{24}\otimes \shift \varphi^{25}\otimes \shift \varphi^3\otimes (\shift \otimes \shift)\partial(\shift^{-1}\varphi^4)\varphi^5\\
&-\varphi^1\varphi^2\otimes \shift \varphi^3\otimes \shift \varphi^{41}\otimes (\shift \otimes \shift)\partial(\shift^{-1}\varphi^{42})\varphi^{43}\varphi^{44}\otimes \shift \varphi^{45}\otimes \shift \varphi^5\otimes 1\\
&-\varphi^1\otimes \shift \varphi^{21}\otimes \shift \varphi^{221}\otimes \varphi^{222}\dots \varphi^{225}\varphi^{23}\varphi^{24}\otimes \shift \varphi^{25}\otimes \shift \varphi^3\otimes \varphi^4\varphi^5\\
&-\varphi^1\varphi^2\varphi^3\otimes \shift \varphi^{41}\otimes \shift \varphi^{421}\otimes \varphi^{422}\dots \varphi^{425}\varphi^{43}\varphi^{44}\otimes \shift \varphi^{45}\otimes \shift \varphi^5\otimes 1\\
&-1\otimes \shift \varphi^1\otimes \shift \varphi^{21}\otimes \varphi^{221}\dots \varphi^{224}\otimes \shift \varphi^{225}\otimes \shift \varphi^{23}\otimes \varphi^{24}\varphi^{25}\varphi^3\varphi^4\varphi^5\\
&-1\otimes \shift \varphi^1\otimes \shift \varphi^{21}\otimes \varphi^{22}\varphi^{23}(\shift \otimes \shift)\partial(\shift^{-1}\varphi^{24})\otimes \shift \varphi^{25}\otimes \shift \varphi^3\otimes \varphi^4\varphi^5\\
&-1\otimes \shift \varphi^1\otimes \varphi^{21}\varphi^{22}\otimes \shift \varphi^{23}\otimes \varphi^{24}\otimes \shift \varphi^{25}\otimes \shift \varphi^3\otimes \varphi^4\varphi^5\\
&+1\otimes \shift \varphi^1\otimes \varphi^2\varphi^3\varphi^{41}\varphi^{42}\otimes \shift \varphi^{43}\otimes \varphi^{44}\otimes \shift \varphi^{45}\otimes \shift \varphi^5\otimes 1\\
&+1\otimes \shift \varphi^1\otimes \varphi^2\otimes \shift \varphi^3\otimes \shift \varphi^{41}\otimes \varphi^{42}\dots \varphi^{45}\otimes \shift \varphi^5\otimes 1\\
&+1\otimes \shift (\varphi^1\varphi^{21})\otimes \varphi^{22}\dots \varphi^{24}\otimes \shift \varphi^{25}\otimes \shift \varphi^3\otimes \shift \varphi^{41}\otimes \varphi^{42}\dots \varphi^{45}\varphi^5\\
&-\varphi^1\varphi^2\otimes \shift \varphi^3\otimes \shift \varphi^{41}\otimes \varphi^{42}\otimes \shift \varphi^{43}\otimes \varphi^{44}\varphi^{45}\otimes \shift \varphi^5\otimes 1\\
&-\varphi^1\varphi^{21}\varphi^{22}\otimes \shift \varphi^{23}\otimes \varphi^{24}\otimes \shift \varphi^{25}\otimes \shift \varphi^3\otimes \shift \varphi^{41}\otimes \varphi^{42}\dots \varphi^{45}\varphi^5\\
&-\varphi^1\varphi^{21}\dots \varphi^{24}\otimes \shift \varphi^{25}\otimes \shift \varphi^3\otimes \shift \varphi^{41}\otimes \varphi^{42}\otimes \shift \varphi^{43}\otimes \varphi^{44}\varphi^{45}\varphi^5\\
&+\varphi^1\varphi^{21}\dots \varphi^{24}\otimes \shift \varphi^{25}\otimes \shift \varphi^3\otimes \shift \varphi^{41}\otimes \varphi^{42}\dots \varphi^{44}\otimes \shift (\varphi^{45}\varphi^5)\otimes 1\\
&-\varphi^1\varphi^2\otimes \shift (\varphi^3\varphi^{41})\otimes \varphi^{42}\varphi^{43}\varphi^{441}\dots \varphi^{444}\otimes \shift \varphi^{445}\otimes \shift \varphi^{45}\otimes \shift \varphi^5\otimes 1\\
&+\varphi^1\varphi^2\varphi^3\varphi^{41}\varphi^{42}\otimes \shift (\varphi^{43}\varphi^{441})\otimes \varphi^{442}\dots \varphi^{444}\otimes \shift \varphi^{445}\otimes \shift \varphi^{45}\otimes \shift \varphi^5\otimes 1\\
&-1\otimes \shift \varphi^1\otimes \shift \varphi^{21}\otimes \shift \varphi^{221}\otimes \varphi^{222}\otimes \shift \varphi^{223}\otimes \varphi^{224}\varphi^{225}\varphi^{23}\dots \varphi^{25}\varphi^3\varphi^4\varphi^5\\
&-1\otimes \shift \varphi^1\otimes \shift \varphi^{21}\otimes \shift \varphi^{221}\otimes \varphi^{222}\dots \varphi^{225}\varphi^{23}\varphi^{24}\otimes \shift (\varphi^{25}\varphi^3)\otimes \varphi^4\varphi^5\\
&+1\otimes \shift \varphi^1\otimes \shift \varphi^{21}\otimes \shift \varphi^{221}\otimes \varphi^{222}\dots \varphi^{224}\otimes \shift (\varphi^{225}\varphi^{23})\otimes \varphi^{24}\varphi^{25}\varphi^3\varphi^4\varphi^5\\
&+\varphi^1\varphi^{21}\varphi^{22}\otimes \shift (\varphi^{23}\varphi^{241})\otimes \varphi^{242}\dots \varphi^{244}\otimes \shift \varphi^{245}\otimes \shift \varphi^{25}\otimes \shift \varphi^3\otimes \varphi^4\varphi^5\\
&-\varphi^1\varphi^{21}\dots \varphi^{23}\varphi^{241}\dots \varphi^{244}\otimes \shift \varphi^{245}\otimes \shift \varphi^{25}\otimes \shift \varphi^3\otimes \varphi^4\otimes \shift \varphi^5\otimes 1\\
&-\varphi^1\varphi^2\otimes \shift \varphi^3\otimes \shift \varphi^{41}\otimes \shift \varphi^{421}\otimes \varphi^{422}\otimes \shift \varphi^{423}\otimes \varphi^{424}\varphi^{425}\varphi^{43}\dots \varphi^{45}\varphi^5\\
&+\varphi^1\varphi^2\otimes \shift \varphi^3\otimes \shift \varphi^{41}\otimes \shift \varphi^{421}\otimes \varphi^{422}\dots \varphi^{424}\otimes \shift (\varphi^{425}\varphi^{43})\otimes \varphi^{44}\varphi^{45}\varphi^5\\
&+1\otimes \shift \varphi^1\otimes \varphi^{21}\dots \varphi^{24}\otimes \shift \varphi^{25}\otimes \shift \varphi^3\otimes \varphi^4\otimes \shift \varphi^5\otimes 1\\
&-1\otimes \shift \varphi^1\otimes \varphi^2\otimes \shift (\varphi^3\varphi^{41})\otimes \varphi^{42}\dots \varphi^{44}\otimes \shift \varphi^{45}\otimes \shift \varphi^5\otimes 1\\
&+1\otimes \shift \varphi^1\otimes 1\otimes \shift \varphi^{21}\otimes \varphi^{22}\dots \varphi^{24}\otimes \shift \varphi^{25}\otimes \shift \varphi^3\otimes \varphi^4\varphi^5\\
&+1\otimes \shift \varphi^1\otimes \shift \varphi^{21}\otimes \varphi^{22}\otimes \shift \varphi^{23}\otimes \varphi^{24}\varphi^{25}\varphi^3\varphi^4\otimes \shift \varphi^5\otimes 1\\
&-1\otimes \shift \varphi^1\otimes \shift \varphi^{21}\otimes \varphi^{22}\varphi^{23}\varphi^{24}\otimes \shift (\varphi^{25}\varphi^3)\otimes \varphi^4\otimes \shift \varphi^5\otimes 1\\
&-\varphi^1\varphi^2\otimes \shift \varphi^3\otimes \shift \varphi^{41}\otimes \varphi^{42}\otimes \shift \varphi^{43}\otimes \shift \varphi^{441}\otimes \varphi^{442}\dots \varphi^{445}\varphi^{45}\varphi^5\\
&+\varphi^1\varphi^2\otimes \shift \varphi^3\otimes \shift \varphi^{41}\otimes \varphi^{42}\dots \varphi^{44}\otimes \shift \varphi^{45}\otimes 1\otimes \shift \varphi^5\otimes 1\\
&+1\otimes \shift \varphi^1\otimes \shift \varphi^{21}\otimes \shift \varphi^{221}\otimes \varphi^{222}\dots \varphi^{225}\varphi^{23}\dots \varphi^{25}\varphi^3\varphi^4\otimes \shift \varphi^5\otimes 1\\
&+1\otimes \shift \varphi^1\otimes \varphi^2\varphi^3\varphi^{41}\dots \varphi^{43}\varphi^{441}\dots \varphi^{444}\otimes \shift \varphi^{445}\otimes \shift \varphi^{45}\otimes \shift \varphi^5\otimes 1\\
&-1\otimes \shift \varphi^1\otimes \varphi^2\otimes \shift \varphi^3\otimes \shift \varphi^{41}\otimes \shift \varphi^{421}\otimes \varphi^{422}\dots \varphi^{425}\varphi^{43}\dots \varphi^{45}\varphi^5\\
&-1\otimes \shift (\varphi^1\varphi^{21})\otimes \varphi^{22}\varphi^{23}\varphi^{241}\dots \varphi^{244}\otimes \shift \varphi^{245}\otimes \shift \varphi^{25}\otimes \shift \varphi^3\otimes \varphi^4\varphi^5\\
&-\varphi^1\varphi^2\otimes \shift \varphi^3\otimes \shift \varphi^{41}\otimes \shift \varphi^{421}\otimes \varphi^{422}\dots \varphi^{425}\varphi^{43}\varphi^{44}\otimes \shift (\varphi^{45}\varphi^5)\otimes 1\\
&-\varphi^1\varphi^2\otimes \shift \varphi^3\otimes \shift \varphi^{41}\otimes \shift \varphi^{421}\otimes \shift \varphi^{4221}\otimes \varphi^{4222}\dots \varphi^{4225}\varphi^{423}\dots \varphi^{425}\varphi^{43}\dots \varphi^{45}\varphi^5\\
&-\varphi^1\varphi^2\varphi^3\varphi^{41}\dots \varphi^{43}\varphi^{441}\dots \varphi^{443}\varphi^{4441}\dots \varphi^{4444}\otimes \shift \varphi^{4445}\otimes \shift \varphi^{445}\otimes \shift \varphi^{45}\otimes \shift \varphi^5\otimes 1
\end{align*}
\end{lem}

\begin{proof}
Recall the $A_\infty$-equation given in the order in which we compute the terms below:
\begin{align*}
&-(\mu_5\otimes 1+1\otimes \mu_5)\mu_2+(\mu_2\otimes 1^{\otimes 4}-1\otimes \mu_2\otimes 1^{\otimes 3}+1^{\otimes 2}\otimes \mu_2\otimes 1^{\otimes 2}-1^{\otimes 3}\otimes \mu_2\otimes 1+1^{\otimes 4}\otimes \mu_2)\mu_5\\
&+(\mu_4\otimes 1^{\otimes 2}-1\otimes \mu_4\otimes 1+1^{\otimes 2}\otimes \mu_4)\mu_3-(\mu_3\otimes 1^{\otimes 3}+1\otimes \mu_3\otimes 1^{\otimes 2}+1^{\otimes 2}\otimes \mu_3\otimes 1+1^{\otimes 3}\otimes \mu_3)\mu_4\\
&-(\mu_1\otimes 1^{\otimes 5}+1\otimes \mu_1\otimes 1^{\otimes 4}+1^{\otimes 2}\otimes \mu_1\otimes 1^{\otimes 3}+1^{\otimes 3}\otimes \mu_1\otimes 1^{\otimes 2}+1^{\otimes 4}\otimes \mu_1\otimes 1+1^{\otimes 5}\otimes \mu_1)\mu_6+\mu_6\mu_1=0
\end{align*}

As $\mu_2(\varphi)=\hat{\omega}\varphi+\varphi\hat{\omega}+(\shift \otimes \shift)\partial(\shift^{-1}\varphi)$ we compute

\begin{align*}
&-(\mu_5\otimes 1)(\hat{\omega}\varphi+\varphi\hat{\omega})\\
=\,&0\+\underbrace{1\otimes \shift \varphi^1\otimes \varphi^2\otimes \shift \varphi^3\otimes \varphi^4\otimes \shift \varphi^5\otimes \hat{\omega}}_{44}
+\underbrace{1\otimes \shift \varphi^1\otimes \varphi^2\varphi^3\varphi^{41}\dots \varphi^{44}\otimes \shift \varphi^{45}\otimes \shift \varphi^5\otimes \hat{\omega}}_{45}\\
&+\underbrace{1\otimes \shift \varphi^1\otimes \shift \varphi^{21}\otimes \varphi^{22}\dots \varphi^{25}\varphi^3\varphi^4\otimes \shift \varphi^5\otimes \hat{\omega}}_{46}
+\underbrace{1\otimes \shift \varphi^1\otimes \varphi^{21}\dots \varphi^{24}\otimes \shift \varphi^{25}\otimes \shift \varphi^3\otimes \varphi^4\varphi^5\hat{\omega}}_{52}\\
&+\underbrace{\varphi^1\varphi^2\otimes \shift \varphi^3\otimes \shift \varphi^{41}\otimes \varphi^{42}\dots \varphi^{45}\otimes \shift \varphi^5\otimes \hat{\omega}}_{47}
+\underbrace{\varphi^1\varphi^{21}\dots \varphi^{24}\otimes \shift \varphi^{25}\otimes \shift \varphi^3\otimes \varphi^4\otimes \shift \varphi^5\otimes \hat{\omega}}_{48}\\
&+\underbrace{1\otimes \shift \varphi^1\otimes \varphi^2\otimes \shift \varphi^3\otimes \shift \varphi^{41}\otimes \varphi^{42}\dots \varphi^{45}\varphi^5\hat{\omega}}_{78}
+\underbrace{\varphi^1\otimes \shift \varphi^{21}\otimes \varphi^{22}\dots \varphi^{24}\otimes \shift \varphi^{25}\otimes \shift \varphi^3\otimes \varphi^4\varphi^5\hat{\omega}}_{71}\\
&+\underbrace{\varphi^1\varphi^2\otimes \shift \varphi^3\otimes \shift \varphi^{41}\otimes \varphi^{42}\varphi^{43}\varphi^{44}\otimes \shift \varphi^{45}\otimes \varphi^5\hat{\omega}}_{117}
-\underbrace{\varphi^1\varphi^2\otimes \shift (\varphi^3\varphi^{41})\otimes \varphi^{42}\dots \varphi^{44}\otimes \shift \varphi^{45}\otimes \shift \varphi^5\otimes \hat{\omega}}_{49}\\
&-\underbrace{1\otimes \shift \varphi^1\otimes \shift \varphi^{21}\otimes \varphi^{22}\varphi^{23} \varphi^{24}\otimes \shift (\varphi^{25}\varphi^3)\otimes \varphi^4\varphi^5\hat{\omega}}_{98}
-\underbrace{\varphi^1\varphi^{21}\varphi^{22}\otimes \shift \varphi^{23}\otimes \varphi^{24}\otimes \shift \varphi^{25}\otimes \shift \varphi^3\otimes \varphi^4\varphi^5\hat{\omega}}_{102}\\
&-\underbrace{\varphi^1\varphi^2\otimes \shift \varphi^3\otimes \shift \varphi^{41}\otimes \varphi^{42}\otimes \shift \varphi^{43}\otimes \varphi^{44}\varphi^{45}\varphi^5\hat{\omega}}_{110}
+\underbrace{1\otimes \shift \varphi^1\otimes \shift \varphi^{21}\otimes \varphi^{22}\otimes \shift \varphi^{23}\otimes \varphi^{24}\varphi^{25}\varphi^3\varphi^4\varphi^5\hat{\omega}}_{106}\\
&+\underbrace{\varphi^1\varphi^2\varphi^3\varphi^{41}\varphi^{42}\otimes \shift \varphi^{43}\otimes \varphi^{44}\otimes \shift \varphi^{45}\otimes \shift \varphi^5\otimes \hat{\omega}}_{50}
+\underbrace{\varphi^1\varphi^{21}\dots \varphi^{24}\otimes \shift \varphi^{25}\otimes \shift \varphi^3\otimes \shift \varphi^{41}\otimes \varphi^{42}\dots \varphi^{45}\varphi^5\hat{\omega}}_{121}\\
&+\underbrace{\varphi^1\varphi^2\varphi^3\varphi^{41}\dots \varphi^{43}\varphi^{441}\dots \varphi^{444}\otimes \shift \varphi^{445}\otimes \shift \varphi^{45}\otimes \shift \varphi^5\otimes \hat{\omega}}_{51}\\
&+\underbrace{1\otimes \shift \varphi^1\otimes \shift \varphi^{21}\otimes \shift \varphi^{221}\otimes \varphi^{222}\dots \varphi^{225}\varphi^{23}\varphi^{24} \varphi^{25}\varphi^3\varphi^4\varphi^5\hat{\omega}}_{300}\\
&-\underbrace{\varphi^1\varphi^{21}\varphi^{22}\varphi^{23}\varphi^{241}\dots \varphi^{244}\otimes \shift \varphi^{245}\otimes \shift \varphi^{25}\otimes \shift \varphi^3\otimes \varphi^4\varphi^5\hat{\omega}}_{147}\\
&-\underbrace{\varphi^1\varphi^2\otimes \shift \varphi^3\otimes \shift \varphi^{41}\otimes \shift \varphi^{421}\otimes \varphi^{422}\dots \varphi^{425}\varphi^{43}\varphi^{44}\varphi^{45}\varphi^5\hat{\omega}}_{308}
\end{align*}
and
\begin{align*}
&-(\mu_5\otimes 1)((\shift \otimes \shift)\partial(\shift^{-1}\varphi))=(\mu_5\otimes 1)(\varphi^1\varphi^2\varphi^3\varphi^4\varphi^5)\\
=\,&-\underbrace{\varphi^1\otimes \shift \varphi^{21}\otimes \varphi^{22}\otimes \shift \varphi^{23}\otimes \varphi^{24}\otimes \shift \varphi^{25}\otimes \varphi^3\varphi^4\varphi^5}_{201}\\
&
-\underbrace{\varphi^1\otimes \shift \varphi^{21}\otimes \varphi^{22}\varphi^{23}\varphi^{241}\dots \varphi^{244}\otimes \shift \varphi^{245}\otimes \shift \varphi^{25}\otimes \varphi^3\varphi^4\varphi^5}_{257}\\
&-\underbrace{\varphi^1\otimes \shift \varphi^{21}\otimes \shift \varphi^{221}\otimes \varphi^{222}\dots \varphi^{225}\varphi^{23}\varphi^{24}\otimes \shift \varphi^{25}\otimes \varphi^3\varphi^4\varphi^5}_{241}\\
&-\underbrace{\varphi^1\otimes \shift \varphi^{21}\otimes \varphi^{221}\dots \varphi^{224}\otimes \shift\varphi^{225}\otimes \shift \varphi^{23}\otimes \varphi^{24}\varphi^{25}\varphi^3\varphi^4\varphi^5}_{251}\\
&-\underbrace{\varphi^1\varphi^{21}\varphi^{22}\otimes \shift \varphi^{23}\otimes \shift \varphi^{241}\otimes \varphi^{242}\dots \varphi^{245}\otimes \shift \varphi^{25}\otimes \varphi^3\varphi^4\varphi^5}_{271}\\
&-\underbrace{\varphi^1\varphi^{21}\varphi^{221}\dots \varphi^{224}\otimes \shift \varphi^{225}\otimes \shift \varphi^{23}\otimes \varphi^{24}\otimes \shift \varphi^{25}\otimes \varphi^3\varphi^4\varphi^5}_{173}\\
&-\underbrace{\varphi^1\otimes \shift \varphi^{21}\otimes \varphi^{22}\otimes \shift \varphi^{23}\otimes \shift \varphi^{241}\otimes \varphi^{242}\dots \varphi^{245}\varphi^{25}\varphi^3\varphi^4\varphi^5}_{278}\\
&-\underbrace{\varphi^1\varphi^{21}\otimes \shift \varphi^{221}\otimes \varphi^{222}\dots \varphi^{224}\otimes \shift \varphi^{225}\otimes \shift \varphi^{23}\otimes \varphi^{24}\varphi^{25}\varphi^3\varphi^4\varphi^5}_g\\
&-\underbrace{\varphi^1\varphi^{21}\varphi^{22}\otimes \shift \varphi^{23}\otimes \shift \varphi^{241}\otimes \varphi^{242}\dots \varphi^{244}\otimes \shift \varphi^{245}\otimes \varphi^{25}\varphi^3\varphi^4\varphi^5}_g\\
&+\underbrace{\varphi^1\varphi^{21}\varphi^{22}\otimes \shift (\varphi^{23}\varphi^{241})\otimes \varphi^{242}\dots \varphi^{244}\otimes \shift \varphi^{245}\otimes \shift \varphi^{25}\otimes \varphi^3\varphi^4\varphi^5}_{290}\\
&+\underbrace{\varphi^1\otimes \shift \varphi^{21}\otimes \shift \varphi^{221}\otimes \varphi^{222}\dots \varphi^{224}\otimes \shift (\varphi^{225}\varphi^{23})\otimes \varphi^{24}\varphi^{25}\varphi^3\varphi^4\varphi^5}_{236}\\
&+\underbrace{\varphi^1\varphi^{21}\varphi^{221}\varphi^{222}\otimes \shift \varphi^{223}\otimes \varphi^{224}\otimes \shift \varphi^{225}\otimes \shift \varphi^{23}\otimes \varphi^{24}\varphi^{25}\varphi^3\varphi^4\varphi^5}_g\\
&+\underbrace{\varphi^1\varphi^{21}\varphi^{22}\otimes \shift \varphi^{23}\otimes \shift \varphi^{241}\otimes \varphi^{242}\otimes \shift \varphi^{243}\otimes \varphi^{244}\varphi^{245}\varphi^{25}\varphi^3\varphi^4\varphi^5}_g\\
&-\underbrace{\varphi^1\otimes \shift \varphi^{21}\otimes \shift \varphi^{221}\otimes \varphi^{222}\otimes \shift \varphi^{223}\otimes \varphi^{224}\varphi^{225}\varphi^{23}\varphi^{24}\varphi^{25}\varphi^3\varphi^4\varphi^5}_{218}\\
&-\underbrace{\varphi^1\varphi^{21}\varphi^{22}\varphi^{23}\varphi^{241}\varphi^{242}\otimes \shift \varphi^{243}\otimes \varphi^{244}\otimes \shift \varphi^{245}\otimes \shift \varphi^{25}\otimes \varphi^3\varphi^4\varphi^5}_{174}\\
&-\underbrace{\varphi^1\varphi^{21}\varphi^{221}\dots \varphi^{224}\otimes \shift \varphi^{225}\otimes \shift \varphi^{23}\otimes \shift \varphi^{241}\otimes \varphi^{242}\dots \varphi^{245}\varphi^{25}\varphi^3\varphi^4\varphi^5}_{gg}\\
&-\underbrace{\varphi^1\varphi^{21}\dots \varphi^{23}\varphi^{241}\dots \varphi^{243}\varphi^{2441}\dots \varphi^{2444}\otimes \shift \varphi^{2445}\otimes \shift \varphi^{245}\otimes \shift \varphi^{25}\otimes \varphi^3\varphi^4\varphi^5}_{296}\\
&-\underbrace{\varphi^1\otimes \shift \varphi^{21}\otimes \shift \varphi^{221}\otimes \shift \varphi^{2221}\otimes \varphi^{2222}\dots \varphi^{2225}\varphi^{223}\dots \varphi^{225}\varphi^{23}\dots \varphi^{25}\varphi^3\varphi^4\varphi^5}_{gg}\\
&+\underbrace{\varphi^1\varphi^{21}\varphi^{221}\dots \varphi^{223}\varphi^{2241}\dots \varphi^{2244}\otimes \shift \varphi^{2245}\otimes \shift \varphi^{225}\otimes \shift \varphi^{23}\otimes \varphi^{24}\varphi^{25}\varphi^3\varphi^4\varphi^5}_{gg}\\
&+\underbrace{\varphi^1\varphi^{21}\varphi^{22}\otimes \shift \varphi^{23}\otimes \shift \varphi^{241}\otimes \shift \varphi^{2421}\otimes \varphi^{2422}\dots \varphi^{2425}\varphi^{243}\dots \varphi^{245}\varphi^{25}\varphi^3\varphi^4\varphi^5}_{gg}
\end{align*}
and
\begin{align*}
&-(1\otimes \mu_5)(\hat{\omega}\varphi+\varphi\hat{\omega})\\
=\,&\underbrace{\hat{\omega}\otimes \shift \varphi^1\otimes \varphi^2\otimes \shift \varphi^3\otimes \varphi^4\otimes \shift \varphi^5\otimes 1}_1
+\underbrace{\hat{\omega}\otimes \shift \varphi^1\otimes \varphi^2\varphi^3\varphi^{41}\dots \varphi^{44}\otimes \shift \varphi^{45}\otimes \shift \varphi^5\otimes 1}_2\\
&+\underbrace{\hat{\omega}\otimes \shift \varphi^1\otimes \shift \varphi^{21}\otimes \varphi^{22}\dots \varphi^{25}\varphi^3\varphi^4\otimes \shift \varphi^5\otimes 1}_3
+\underbrace{\hat{\omega}\otimes \shift \varphi^1\otimes \varphi^{21}\dots \varphi^{24}\otimes \shift\varphi^{25}\otimes \shift \varphi^3\otimes \varphi^4\varphi^5}_4\\
&+\underbrace{\hat{\omega}\varphi^1\varphi^2\otimes \shift \varphi^3\otimes \shift \varphi^{41}\otimes \varphi^{42}\dots \varphi^{45}\otimes \shift \varphi^5\otimes 1}_{75}
+\underbrace{\hat{\omega}\varphi^1\varphi^{21}\dots \varphi^{24}\otimes \shift \varphi^{25}\otimes \shift \varphi^3\otimes \varphi^4\otimes \shift \varphi^5\otimes 1}_{54}\\
&+\underbrace{\hat{\omega}\otimes \shift \varphi^1\otimes \varphi^2\otimes \shift \varphi^3\otimes \shift \varphi^{41}\otimes \varphi^{42}\dots \varphi^{45}\varphi^5}_5
+\underbrace{\hat{\omega}\varphi^1\otimes \shift \varphi^{21}\otimes \varphi^22\varphi^{23} \varphi^{24}\otimes \shift \varphi^{25}\otimes \shift \varphi^3\otimes \varphi^4\varphi^5}_{152}\\
&+\underbrace{\hat{\omega}\varphi^1\varphi^2\otimes \shift \varphi^3\otimes \shift \varphi^{41}\otimes \varphi^{42}\dots \varphi^{44}\otimes \shift \varphi^{45}\otimes \varphi^5}_{79}
-\underbrace{\hat{\omega}\varphi^1\varphi^2\otimes \shift (\varphi^3\varphi^{41})\otimes \varphi^{42}\dots \varphi^{44}\otimes \shift \varphi^{45}\otimes \shift \varphi^5\otimes 1}_{83}\\
&-\underbrace{\hat{\omega}\otimes \shift \varphi^1\otimes \shift \varphi^{21}\otimes \varphi^{22}\dots \varphi^{24}\otimes \shift (\varphi^{25}\varphi^3)\otimes \varphi^4\varphi^5}_6
-\underbrace{\hat{\omega}\varphi^1\varphi^{21}\varphi^{22}\otimes \shift \varphi^{23}\otimes \varphi^{24}\otimes \shift \varphi^{25}\otimes \shift \varphi^3\otimes \varphi^4\varphi^5}_{128}\\
&-\underbrace{\hat{\omega}\varphi^1\varphi^2\otimes \shift \varphi^3\otimes \shift \varphi^{41}\otimes \varphi^{42}\otimes \shift \varphi^{43}\otimes \varphi^{44}\varphi^{45}\varphi^5}_{81}
+\underbrace{\hat{\omega}\otimes \shift \varphi^1\otimes \shift \varphi^{21}\otimes \varphi^{22}\otimes \shift \varphi^{23}\otimes \varphi^{24}\varphi^{25}\varphi^3\varphi^4\varphi^5}_7\\
&+\underbrace{\hat{\omega}\varphi^1\varphi^2\varphi^3\varphi^{41}\varphi^{42}\otimes \shift \varphi^{43}\otimes \varphi^{44}\otimes \shift \varphi^{45}\otimes \shift \varphi^5\otimes 1}_{82}
+\underbrace{\hat{\omega}\varphi^1\varphi^{21}\dots \varphi^{24}\otimes \shift \varphi^{25}\otimes \shift \varphi^3\otimes \shift \varphi^{41}\otimes \varphi^{42}\dots \varphi^{45}\varphi^5}_{123}\\
&+\underbrace{\hat{\omega}\varphi^1\varphi^2\varphi^3\varphi^{41}\varphi^{42}\varphi^{43}\varphi^{441}\dots \varphi^{444}\otimes \shift \varphi^{445}\otimes \shift \varphi^{45}\otimes \shift \varphi^5\otimes 1}_{301}\\
&+\underbrace{\hat{\omega}\otimes \shift \varphi^1\otimes \shift \varphi^{21}\otimes \shift \varphi^{221}\otimes \varphi^{222}\dots \varphi^{225}\varphi^{23}\dots \varphi^{25}\varphi^3\varphi^4\varphi^5}_8\\
&-\underbrace{\hat{\omega}\varphi^1\varphi^{21}\varphi^{22} \varphi^{23}\varphi^{241}\dots \varphi^{244}\otimes \shift \varphi^{245}\otimes \shift \varphi^{25}\otimes \shift \varphi^3\otimes \varphi^4\varphi^5}_{306}\\
&-\underbrace{\hat{\omega}\varphi^1\varphi^2\otimes \shift \varphi^3\otimes \shift \varphi^{41}\otimes \shift \varphi^{421}\otimes \varphi^{422}\dots \varphi^{425}\varphi^{43}\varphi^{44}\varphi^{45}\varphi^5}_{302}\+0
\end{align*}
and
\begin{align*}
&-(1\otimes \mu_5)((\shift \otimes \shift)\partial(\shift^{-1}\varphi))=(1\otimes \mu_5)(\varphi^1\varphi^2\varphi^3\varphi^4\varphi^5)\\
=\,&-\underbrace{\varphi^1\varphi^2\varphi^3\otimes \shift \varphi^{41}\otimes \varphi^{42}\otimes \shift \varphi^{43}\otimes \varphi^{44}\otimes \shift \varphi^{45}\otimes \varphi^5}_{207}\\
&
-\underbrace{\varphi^1\varphi^2\varphi^3\otimes \shift \varphi^{41}\otimes \varphi^{42}\varphi^{43}\varphi^{441}\dots \varphi^{444}\otimes \shift \varphi^{445}\otimes \shift \varphi^{45}\otimes \varphi^5}_{281}\\
&-\underbrace{\varphi^1\varphi^2\varphi^3\otimes \shift \varphi^{41}\otimes \shift \varphi^{421}\otimes \varphi^{422}\dots \varphi^{425}\varphi^{43}\varphi^{44}\otimes \shift \varphi^{45}\otimes \varphi^5}_{244}\\
&-\underbrace{\varphi^1\varphi^2\varphi^3\otimes \shift \varphi^{41}\otimes \varphi^{421}\dots \varphi^{424}\otimes \shift \varphi^{425}\otimes \shift \varphi^{43}\otimes \varphi^{44}\varphi^{45}\varphi^5}_{274}\\
&-\underbrace{\varphi^1\varphi^2\varphi^3\varphi^{41}\varphi^{42}\otimes \shift \varphi^{43}\otimes \shift \varphi^{441}\otimes \varphi^{442}\dots \varphi^{445}\otimes \shift \varphi^{45}\otimes \varphi^5}_{287}\\
&-\underbrace{\varphi^1\varphi^2\varphi^3\varphi^{41}\varphi^{421}\dots \varphi^{424}\otimes \shift \varphi^{425}\otimes \shift \varphi^{43}\otimes \varphi^{44}\otimes \shift \varphi^{45}\otimes \varphi^5}_{266}\\
&-\underbrace{\varphi^1\varphi^2\varphi^3\otimes \shift \varphi^{41}\otimes \varphi^{42}\otimes \shift \varphi^{43}\otimes \shift \varphi^{441}\otimes \varphi^{442}\dots \varphi^{445}\varphi^{45}\varphi^5}_{294}\\
&-\underbrace{\varphi^1\varphi^2\varphi^3\varphi^{41}\otimes \shift \varphi^{421}\otimes \varphi^{422}\dots \varphi^{424}\otimes \shift \varphi^{425}\otimes \shift \varphi^{43}\otimes \varphi^{44}\varphi^{45}\varphi^5}_g\\
&-\underbrace{\varphi^1\varphi^2\varphi^3\varphi^{41}\varphi^{42}\otimes \shift \varphi^{43}\otimes \shift \varphi^{441}\otimes \varphi^{442}\dots \varphi^{444}\otimes \shift \varphi^{445}\otimes \varphi^{45}\varphi^5}_g\\
&+\underbrace{\varphi^1\varphi^2\varphi^3\varphi^{41}\varphi^{42}\otimes \shift (\varphi^{43}\varphi^{441})\otimes \varphi^{442}\dots \varphi^{444}\otimes \shift \varphi^{445}\otimes \shift \varphi^{45}\otimes \varphi^5}_{217}\\
&+\underbrace{\varphi^1\varphi^2\varphi^3\otimes \shift \varphi^{41}\otimes \shift \varphi^{421}\otimes \varphi^{422}\varphi^{423}\varphi^{424}\otimes \shift (\varphi^{425}\varphi^{43})\otimes \varphi^{44}\varphi^{45}\varphi^5}_{289}\\
&+\underbrace{\varphi^1\varphi^2\varphi^3\varphi^{41}\varphi^{421}\varphi^{422}\otimes \shift \varphi^{423}\otimes \varphi^{424}\otimes \shift \varphi^{425}\otimes \shift \varphi^{43}\otimes \varphi^{44}\varphi^{45}\varphi^5}_g\\
&+\underbrace{\varphi^1\varphi^2\varphi^3\varphi^{41}\varphi^{42}\otimes \shift \varphi^{43}\otimes \shift \varphi^{441}\otimes \varphi^{442}\otimes \shift \varphi^{443}\otimes \varphi^{444}\varphi^{445}\varphi^{45}\varphi^5}_g\\
&-\underbrace{\varphi^1\varphi^2\varphi^3\otimes \shift \varphi^{41}\otimes \shift \varphi^{421}\otimes \varphi^{422}\otimes \shift \varphi^{423}\otimes \varphi^{424}\varphi^{425}\varphi^{43}\varphi^{44}\varphi^{45}\varphi^5}_{220}\\
&-\underbrace{\varphi^1\varphi^2\varphi^3\varphi^{41}\varphi^{42}\varphi^{43}\varphi^{441}\varphi^{442}\otimes \shift \varphi^{443}\otimes \varphi^{444}\otimes \shift \varphi^{445}\otimes \shift \varphi^{45}\otimes \varphi^5}_{176}\\
&-\underbrace{\varphi^1\varphi^2\varphi^3\varphi^{41}\varphi^{421}\dots \varphi^{424}\otimes \shift \varphi^{425}\otimes \shift \varphi^{43}\otimes \shift \varphi^{441}\otimes \varphi^{442}\dots \varphi^{445}\varphi^{45}\varphi^5}_{gg}\\
&-\underbrace{\varphi^1\varphi^2\varphi^3\varphi^{41}\varphi^{42}\varphi^{43}\varphi^{441}\varphi^{442} \varphi^{443}\varphi^{4441}\dots \varphi^{4444}\otimes \shift \varphi^{4445}\otimes \shift \varphi^{445}\otimes \shift \varphi^{45}\otimes \varphi^5}_{158}\\
&-\underbrace{\varphi^1\varphi^2\varphi^3\otimes \shift \varphi^{41}\otimes \shift \varphi^{421}\otimes \shift \varphi^{4221}\otimes \varphi^{4222}\dots \varphi^{4225}\varphi^{423}\varphi^{424}\varphi^{425}\varphi^{43}\varphi^{44} \varphi^{45}\varphi^5}_{298}\\
&+\underbrace{\varphi^1\varphi^2\varphi^3\varphi^{41}\varphi^{421}\dots \varphi^{423}\varphi^{4241}\dots \varphi^{4244}\otimes \shift \varphi^{4245}\otimes \shift \varphi^{425}\otimes \shift \varphi^{43}\otimes \varphi^{44}\varphi^{45}\varphi^5}_{gg}\\
&+\underbrace{\varphi^1\varphi^2\varphi^3\varphi^{41}\varphi^{42}\otimes \shift \varphi^{43}\otimes \shift \varphi^{441}\otimes \shift \varphi^{4421}\otimes \varphi^{4422}\dots \varphi^{4425}\varphi^{443}\dots \varphi^{445}\varphi^{45}\varphi^5}_{gg}
\end{align*}
and
\begin{align*}
&(\mu_2\otimes 1^{\otimes 4})\mu_5(\varphi)\\
=\,&-\underbrace{\hat{\omega}\otimes \shift \varphi^1\otimes \varphi^2\otimes \shift \varphi^3\otimes \varphi^4\otimes \shift \varphi^5\otimes 1}_1
-\underbrace{1\otimes \shift \varphi^1\otimes \hat{\omega}\varphi^2\otimes \shift \varphi^3\otimes \varphi^4\otimes \shift \varphi^5\otimes 1}_9\\
&+\underbrace{\varphi^{11}\varphi^{12}\otimes \shift \varphi^{13}\otimes \varphi^2\otimes \shift \varphi^3\otimes \varphi^4\otimes \shift \varphi^5\otimes 1}_y
-\underbrace{1\otimes \shift \varphi^{11}\otimes \varphi^{12}\varphi^{13}\varphi^2\otimes \shift \varphi^3\otimes \varphi^4\otimes \shift \varphi^5\otimes 1}_y\\
&\-\underbrace{\hat{\omega}\otimes \shift \varphi^1\otimes \varphi^2\varphi^3\varphi^{41}\dots \varphi^{44}\otimes \shift \varphi^{45}\otimes \shift \varphi^5\otimes 1}_2
-\underbrace{1\otimes \shift \varphi^1\otimes \hat{\omega}\varphi^2\varphi^3\varphi^{41}\dots \varphi^{44}\otimes \shift \varphi^{45}\otimes \shift \varphi^5\otimes 1}_{10}\\
&+\underbrace{\varphi^{11}\varphi^{12}\otimes \shift \varphi^{13}\otimes \varphi^2\varphi^3\varphi^{41}\dots \varphi^{44}\otimes \shift \varphi^{45}\otimes \shift \varphi^5\otimes 1}_g\\
&
-\underbrace{1\otimes \shift \varphi^{11}\otimes \varphi^{12}\varphi^{13}\varphi^2\varphi^3\varphi^{41}\dots \varphi^{44}\otimes \shift \varphi^{45}\otimes \shift \varphi^5\otimes 1}_g\\
&\-\underbrace{\hat{\omega}\otimes \shift \varphi^1\otimes \shift \varphi^{21}\otimes \varphi^{22}\dots \varphi^{25}\varphi^3\varphi^4\otimes \shift \varphi^5\otimes 1}_3
-\underbrace{1\otimes \shift \varphi^1\otimes \shift \varphi^{21}\otimes \hat{\omega}\varphi^{22}\dots \varphi^{25}\varphi^3\varphi^4\otimes \shift \varphi^5\otimes 1}_{11}\\
&+\underbrace{1\otimes \shift \varphi^1\otimes 1\otimes \shift \varphi^{21}\otimes \varphi^{22}\dots \varphi^{25}\varphi^3\varphi^4\otimes \shift \varphi^5\otimes 1}_{333}
+\underbrace{\varphi^{11}\varphi^{12}\otimes \shift \varphi^{13}\otimes \shift \varphi^{21}\otimes \varphi^{22}\dots \varphi^{25}\varphi^3\varphi^4\otimes \shift \varphi^5\otimes 1}_g\\
&-\underbrace{1\otimes \shift \varphi^1\otimes \shift \varphi^{211}\otimes \varphi^{212}\varphi^{213}\varphi^{22}\dots \varphi^{25}\varphi^3\varphi^4\otimes \shift \varphi^5\otimes 1}_y
\-\underbrace{\hat{\omega}\otimes \shift \varphi^1\otimes \varphi^{21}\dots \varphi^{24}\otimes \shift \varphi^{25}\otimes \shift \varphi^3\otimes \varphi^4\varphi^5}_4\\
&-\underbrace{1\otimes \shift \varphi^1\otimes \hat{\omega}\varphi^{21}\dots \varphi^{24}\otimes \shift \varphi^{25}\otimes \shift \varphi^3\otimes \varphi^4\varphi^5}_{96}
+\underbrace{\varphi^{11}\varphi^{12}\otimes \shift \varphi^{13}\otimes \varphi^{21}\dots \varphi^{24}\otimes \shift \varphi^{25}\otimes \shift \varphi^3\otimes \varphi^4\varphi^5}_g\\
&-\underbrace{1\otimes \shift \varphi^{11}\otimes \varphi^{12}\varphi^{13}\varphi^{21}\dots \varphi^{24}\otimes \shift \varphi^{25}\otimes \shift \varphi^3\otimes \varphi^4\varphi^5}_g
\-\underbrace{\varphi^1\hat{\omega}\varphi^2\otimes \shift \varphi^3\otimes \shift \varphi^{41}\otimes \varphi^{42}\dots \varphi^{45}\otimes \shift \varphi^5\otimes 1}_{76}\\
&-\underbrace{\varphi^1\varphi^2\hat{\omega}\otimes \shift \varphi^3\otimes \shift \varphi^{41}\otimes \varphi^{42}\dots \varphi^{45}\otimes \shift \varphi^5\otimes 1}_{12}\\
&
-\underbrace{\varphi^1(\shift \otimes \shift)\partial(\shift^{-1}\varphi^2)\otimes \shift \varphi^3\otimes \shift \varphi^{41}\otimes \varphi^{42}\dots \varphi^{45}\otimes \shift \varphi^5\otimes 1}_{161}\\
&\-\underbrace{\varphi^1\varphi^{21}\hat{\omega}\varphi^{22}\dots \varphi^{24}\otimes \shift \varphi^{25}\otimes \shift \varphi^3\otimes \varphi^4\otimes \shift \varphi^5\otimes 1}_{73}
-\underbrace{\varphi^1\varphi^{21}\varphi^{22}\hat{\omega}\varphi^{23}\varphi^{24}\otimes \shift \varphi^{25}\otimes \shift \varphi^3\otimes \varphi^4\otimes \shift \varphi^5\otimes 1}_{104}\\
&-\underbrace{\varphi^1\varphi^{21}(\shift \otimes \shift)\partial(\shift^{-1}\varphi^{22})\varphi^{23}\varphi^{24}\otimes \shift \varphi^{25}\otimes \shift \varphi^3\otimes \varphi^4\otimes \varphi^5\otimes 1}_y\\
&
\-\underbrace{\hat{\omega}\otimes \shift \varphi^1\otimes \varphi^2\otimes \shift \varphi^3\otimes \shift \varphi^{41}\otimes \varphi^{42}\dots \varphi^{45}\varphi^5}_5
-\underbrace{1\otimes \shift \varphi^1\otimes \hat{\omega}\varphi^2\otimes \shift \varphi^3\otimes \shift \varphi^{41}\otimes \varphi^{42}\dots \varphi^{45}\varphi^5}_{13}\\
&
+\underbrace{\varphi^{11}\varphi^{12}\otimes \shift \varphi^{13}\otimes \varphi^2\otimes \shift \varphi^3\otimes \shift \varphi^{41}\otimes \varphi^{42}\dots \varphi^{45}\varphi^5}_g\\
&-\underbrace{1\otimes \shift \varphi^{11}\otimes \varphi^{12}\varphi^{13}\varphi^2\otimes \shift \varphi^3\otimes \shift \varphi^{41}\otimes \varphi^{42}\dots \varphi^{45}\varphi^5}_g\\
&
\-\underbrace{\varphi^1\hat{\omega}\otimes \shift \varphi^{21}\otimes \varphi^{22}\varphi^{23}\varphi^{24}\otimes \shift \varphi^{25}\otimes \shift \varphi^3\otimes \varphi^4\varphi^5}_{94}\\
&-\underbrace{\varphi^1\otimes \shift \varphi^{21}\otimes \hat{\omega}\varphi^{22}\varphi^{23}\varphi^{24}\otimes \shift \varphi^{25}\otimes \shift \varphi^3\otimes \varphi^4\varphi^5}_{14}\\
&
+\underbrace{\varphi^1\varphi^{211}\varphi^{212}\otimes \shift \varphi^{213}\otimes \varphi^{22}\dots \varphi^{24}\otimes \shift \varphi^{25}\otimes \shift \varphi^3\otimes \varphi^4\varphi^5}_y\\
&-\underbrace{\varphi^1\otimes \shift \varphi^{211}\otimes \varphi^{212}\varphi^{213}\otimes \varphi^{22}\dots \varphi^{24}\otimes \shift \varphi^{25}\otimes \shift \varphi^3\otimes \varphi^4\varphi^5}_y\\
&
\-\underbrace{\varphi^1\hat{\omega}\varphi^2\otimes \shift \varphi^3\otimes \shift \varphi^{41}\otimes \varphi^{42}\dots \varphi^{44}\otimes \shift \varphi^{45}\otimes \varphi^5}_{80}\\
&-\underbrace{\varphi^1\varphi^2\hat{\omega}\otimes \shift \varphi^3\otimes \shift \varphi^{41}\otimes \varphi^{42}\dots \varphi^{44}\otimes \shift \varphi^{45}\otimes \varphi^5}_{15}\\
&
-\underbrace{\varphi^1(\shift \otimes \shift)\partial(\shift^{-1}\varphi^2)\otimes \shift \varphi^3\otimes \shift \varphi^{41}\otimes \varphi^{42}\dots \varphi^{44}\otimes \shift \varphi^{45}\otimes \varphi^5}_{163}\\
&\+\underbrace{\varphi^1\hat{\omega}\varphi^2\otimes \shift (\varphi^3\varphi^{41})\otimes \varphi^{42}\varphi^{43}\varphi^{44}\otimes \shift \varphi^{45}\otimes \shift \varphi^5\otimes 1}_{84}\\
&
+\underbrace{\varphi^1\varphi^2\hat{\omega}\otimes \shift (\varphi^3\varphi^{41})\otimes \varphi^{42}\dots \varphi^{44}\otimes \shift \varphi^{45}\otimes \shift \varphi^5\otimes 1}_{16}\\
&+\underbrace{\varphi^1(\shift \otimes \shift)\partial(\shift^{-1}\varphi^2)\otimes \shift (\varphi^3\varphi^{41})\otimes \varphi^{42}\dots \varphi^{44}\otimes \shift \varphi^{45}\otimes \shift \varphi^5\otimes 1}_{353}\\
&
\+\underbrace{\hat{\omega}\otimes \shift \varphi^1\otimes \shift \varphi^{21}\otimes \varphi^{22}\dots \varphi^{24}\otimes \shift (\varphi^{25}\varphi^3)\otimes \varphi^4\varphi^5}_6
+\underbrace{1\otimes \shift \varphi^1\otimes \shift \varphi^{21}\otimes \hat{\omega}\varphi^{22}\dots \varphi^{24}\otimes \shift (\varphi^{25}\varphi^3)\otimes \varphi^4\varphi^5}_{17}\\
&-\underbrace{1\otimes \shift \varphi^1\otimes 1\otimes \shift \varphi^{21}\otimes \varphi^{22}\varphi^{23} \varphi^{24}\otimes \shift (\varphi^{25}\varphi^3)\otimes \varphi^4\varphi^5}_{95}\\
&
-\underbrace{\varphi^{11}\varphi^{12}\otimes \shift \varphi^{13}\otimes \shift \varphi^{21}\otimes \varphi^{22}\dots \varphi^{24}\otimes \shift (\varphi^{25}\varphi^3)\otimes \varphi^4\varphi^5}_g\\
&+\underbrace{1\otimes \shift \varphi^1\otimes \shift \varphi^{211}\otimes \varphi^{212}\varphi^{213}\varphi^{22}\dots \varphi^{24}\otimes \shift (\varphi^{25}\varphi^3)\otimes \varphi^4\varphi^5}_y\\
&
\+\underbrace{\varphi^1\varphi^{21}\hat{\omega}\varphi^{22}\otimes \shift \varphi^{23}\otimes \varphi^{24}\otimes \shift \varphi^{25}\otimes \shift \varphi^3\otimes \varphi^4\varphi^5}_{131}\\
&+\underbrace{\varphi^1\varphi^{21}\varphi^{22}\hat{\omega}\otimes \shift \varphi^{23}\otimes \varphi^{24}\otimes \shift \varphi^{25}\otimes \shift \varphi^3\otimes \varphi^4\varphi^5}_{18}\\
&
+\underbrace{\varphi^1\varphi^{21}(\shift \otimes \shift)\partial(\shift^{-1}\varphi^{22})\otimes \shift \varphi^{23}\otimes \varphi^{24}\otimes \shift \varphi^{25}\otimes \shift \varphi^3\otimes \varphi^4\varphi^5}_{170}\\
&\+\underbrace{\varphi^1\hat{\omega}\varphi^2\otimes \shift \varphi^3\otimes \shift \varphi^{41}\otimes \varphi^{42}\otimes \shift \varphi^{43}\otimes \varphi^{44}\varphi^{45}\varphi^5}_{93}
+\underbrace{\varphi^1\varphi^2\hat{\omega}\otimes \shift \varphi^3\otimes \shift \varphi^{41}\otimes \varphi^{42}\otimes \shift \varphi^{43}\otimes \varphi^{44}\varphi^{45}\varphi^5}_{19}\\
&+\underbrace{\varphi^1(\shift \otimes \shift)\partial(\shift^{-1}\varphi^2)\otimes \shift \varphi^3\otimes \shift \varphi^{41}\otimes \varphi^{42}\otimes \shift \varphi^{43}\otimes \varphi^{44}\varphi^{45}\varphi^5}_{345}\\
&
\-\underbrace{\hat{\omega}\otimes \shift \varphi^1\otimes \shift \varphi^{21}\otimes \varphi^{22}\otimes \shift \varphi^{23}\otimes \varphi^{24}\varphi^{25}\varphi^3\varphi^4\varphi^5}_7\\
&-\underbrace{1\otimes \shift \varphi^1\otimes \shift \varphi^{21}\otimes \hat{\omega}\varphi^{22}\otimes \shift \varphi^{23}\otimes \varphi^{24}\varphi^{25}\varphi^3\varphi^4\varphi^5}_{20}
+\underbrace{1\otimes \shift \varphi^1\otimes 1\otimes \shift \varphi^{21}\otimes \varphi^{22}\otimes \shift \varphi^{23}\otimes \varphi^{24}\varphi^{25}\varphi^3\varphi^4\varphi^5}_{69}\\
&+\underbrace{\varphi^{11}\varphi^{12}\otimes \shift \varphi^{13}\otimes \shift \varphi^{21}\otimes \varphi^{22}\otimes \shift \varphi^{23}\otimes \varphi^{24}\varphi^{25}\varphi^3\varphi^4\varphi^5}_g\\
&
-\underbrace{1\otimes \shift \varphi^1\otimes \shift \varphi^{211}\otimes \varphi^{212}\varphi^{213}\varphi^{22}\otimes \shift \varphi^{23}\otimes \varphi^{24}\varphi^{25}\varphi^3\varphi^4\varphi^5}_y\\
&\-\underbrace{\varphi^1\hat{\omega}\varphi^2\varphi^3\varphi^{41}\varphi^{42}\otimes \shift \varphi^{43}\otimes \varphi^{44}\otimes \shift \varphi^{45}\otimes \shift \varphi^5\otimes 1}_{88}
-\underbrace{\varphi^1\varphi^2\hat{\omega}\varphi^3\varphi^{41}\varphi^{42}\otimes \shift \varphi^{43}\otimes \varphi^{44}\otimes \shift \varphi^{45}\otimes \shift \varphi^5\otimes 1}_{135}\\
&-\underbrace{\varphi^1(\shift \otimes \shift)\partial(\shift^{-1}\varphi^2)\varphi^3\varphi^{41}\varphi^{42}\otimes \shift \varphi^{43}\otimes \varphi^{44}\otimes \shift \varphi^{45}\otimes \shift \varphi^5\otimes 1}_g\\
&
\-\underbrace{\varphi^1\varphi^{21}\hat{\omega}\varphi^{22}\dots \varphi^{24}\otimes \shift \varphi^{25}\otimes \shift \varphi^3\otimes \shift \varphi^{41}\otimes \varphi^{42}\dots \varphi^{45}\varphi^5}_{124}\\
&-\underbrace{\varphi^1\varphi^{21}\varphi^{22}\hat{\omega}\varphi^{23}\varphi^{24}\otimes \shift \varphi^{25}\otimes \shift \varphi^3\otimes \shift \varphi^{41}\otimes \varphi^{42}\dots \varphi^{45}\varphi^5}_{133}\\
&
-\underbrace{\varphi^1\varphi^{21}(\shift \otimes \shift)\partial(\shift^{-1}\varphi^{22})\varphi^{23}\varphi^{24}\otimes \shift \varphi^{25}\otimes \shift \varphi^3\otimes \shift \varphi^{41}\otimes \varphi^{42}\dots \varphi^{45}\varphi^5}_y\\
&\-\underbrace{\varphi^1\hat{\omega}\varphi^2\varphi^3\varphi^{41}\varphi^{42} \varphi^{43}\varphi^{441}\dots \varphi^{444}\otimes \shift \varphi^{445}\otimes \shift \varphi^{45}\otimes \shift \varphi^5\otimes 1}_{304}\\
&
-\underbrace{\varphi^1\varphi^2\hat{\omega}\varphi^3\varphi^{41}\varphi^{42}\varphi^{43}\varphi^{441}\dots \varphi^{444}\otimes \shift \varphi^{445}\otimes \shift \varphi^{45}\otimes \shift \varphi^5\otimes 1}_{310}\\
&-\underbrace{\varphi^1(\shift \otimes \shift)\partial(\shift^{-1}\varphi^2)\varphi^3\varphi^{41}\dots \varphi^{43}\varphi^{441}\dots \varphi^{444}\otimes \shift \varphi^{445}\otimes \shift \varphi^{45}\otimes \shift \varphi^5\otimes 1}_{gg}\\
&
\-\underbrace{\hat{\omega}\otimes \shift \varphi^1\otimes \shift \varphi^{21}\otimes \shift \varphi^{221}\otimes \varphi^{222}\dots \varphi^{225}\varphi^{23}\dots \varphi^{25}\varphi^3\varphi^4\varphi^5}_8\\
&-\underbrace{1\otimes \shift \varphi^1\otimes \shift \varphi^{21}\otimes \shift \varphi^{221}\otimes \hat{\omega}\varphi^{222}\dots \varphi^{225}\varphi^{23}\dots \varphi^{25}\varphi^3\varphi^4\varphi^5}_{21}\\
&+\underbrace{1\otimes \shift \varphi^1\otimes 1\otimes \shift \varphi^{21}\otimes \shift \varphi^{221}\otimes \varphi^{222}\dots \varphi^{225}\varphi^{23}\varphi^{24} \varphi^{25}\varphi^3\varphi^4\varphi^5}_{70}\\
&+\underbrace{1\otimes \shift \varphi^1\otimes \shift \varphi^{21}\otimes 1\otimes \shift \varphi^{221}\otimes \varphi^{222}\dots \varphi^{225}\varphi^{23}\dots \varphi^{25}\varphi^3\varphi^4\varphi^5}_{334}\\
&+\underbrace{\varphi^{11}\varphi^{12}\otimes \shift \varphi^{13}\otimes \shift \varphi^{21}\otimes \shift \varphi^{221}\otimes \varphi^{222}\dots \varphi^{225}\varphi^{23}\dots \varphi^{25}\varphi^3\varphi^4\varphi^5}_{gg}\\
&-\underbrace{1\otimes \shift \varphi^1\otimes \shift \varphi^{21}\otimes \shift \varphi^{2211}\otimes \varphi^{2212}\varphi^{2213}\varphi^{222}\dots \varphi^{225}\varphi^{23}\dots \varphi^{25}\varphi^3\varphi^4\varphi^5}_{y}\\
&\+\underbrace{\varphi^1\varphi^{21}\hat{\omega}\varphi^{22}\varphi^{23}\varphi^{241}\dots \varphi^{244}\otimes \shift \varphi^{245}\otimes \shift \varphi^{25}\otimes \shift \varphi^3\otimes \varphi^4\varphi^5}_{307}\\
&+\underbrace{\varphi^1\varphi^{21}\varphi^{22}\hat{\omega}\varphi^{23}\varphi^{241}\dots \varphi^{244}\otimes \shift \varphi^{245}\otimes \shift \varphi^{25}\otimes \shift \varphi^3\otimes \varphi^4\varphi^5}_{312}\\
&
+\underbrace{\varphi^1\varphi^{21}(\shift \otimes \shift)\partial(\shift^{-1}\varphi^{22})\varphi^{23}\varphi^{241}\dots \varphi^{244}\otimes \shift \varphi^{245}\otimes \shift \varphi^{25}\otimes \shift \varphi^3\otimes \varphi^4\varphi^5}_{gg}\\
&\+\underbrace{\varphi^1\hat{\omega}\varphi^2\otimes \shift \varphi^3\otimes \shift \varphi^{41}\otimes \shift \varphi^{421}\otimes \varphi^{422}\dots \varphi^{425}\varphi^{43}\varphi^{44}\varphi^{45}\varphi^5}_{303}\\
&
+\underbrace{\varphi^1\varphi^2\hat{\omega}\otimes \shift \varphi^3\otimes \shift \varphi^{41}\otimes \shift \varphi^{421}\otimes \varphi^{422}\dots \varphi^{425}\varphi^{43}\dots \varphi^{45}\varphi^5}_{22}\\
&+\underbrace{\varphi^1(\shift \otimes \shift)\partial( \shift^{-1}\varphi^2)\otimes \shift \varphi^3\otimes \shift \varphi^{41}\otimes \shift \varphi^{421}\otimes \varphi^{422}\dots \varphi^{425}\varphi^{43}\dots \varphi^{45}\varphi^5}_{183}
\end{align*}
and
\begin{align*}
&-(1\otimes \mu_2\otimes 1^{\otimes 3})\mu_5(\varphi)\\
=\,&+\underbrace{1\otimes \shift \varphi^1\otimes \hat{\omega}\varphi^2\otimes \shift \varphi^3\otimes \varphi^4\otimes \shift \varphi^5\otimes 1}_9
+\underbrace{1\otimes \shift \varphi^1\otimes \varphi^2\hat{\omega}\otimes \shift \varphi^3\otimes \varphi^4\otimes \shift \varphi^5\otimes 1}_{23}\\
&+\underbrace{1\otimes \shift \varphi^1\otimes (\shift \otimes \shift)\partial(\shift^{-1}\varphi^2)\otimes \shift \varphi^3\otimes \varphi^4\otimes \shift \varphi^5\otimes 1}_{55}
\+\underbrace{1\otimes \shift \varphi^1\otimes \hat{\omega}\varphi^2\varphi^3\varphi^{41}\dots \varphi^{44}\otimes \shift \varphi^{45}\otimes \shift \varphi^5\otimes 1}_{10}\\
&+\underbrace{1\otimes \shift \varphi^1\otimes \varphi^2\hat{\omega}\varphi^3\varphi^{41}\dots \varphi^{44}\otimes \shift \varphi^{45}\otimes \shift \varphi^5\otimes 1}_{85}\\
&
+\underbrace{1\otimes \shift \varphi^1\otimes (\shift \otimes \shift)\partial(\shift^{-1}\varphi^2)\varphi^3\varphi^{41}\dots \varphi^{44}\otimes \shift \varphi^{45}\otimes \shift \varphi^5\otimes 1}_{179}\\
&\+\underbrace{1\otimes \shift \varphi^1\otimes \shift \varphi^{21}\otimes \hat{\omega}\varphi^{22}\dots \varphi^{25}\varphi^3\varphi^4\otimes \shift \varphi^5\otimes 1}_{11}
+\underbrace{1 \otimes \shift\varphi^1\otimes \shift \varphi^{21}\otimes \varphi^{22}\hat{\omega}\varphi^{23}\varphi^{24}\varphi^{25}\varphi^3\varphi^4\otimes \shift \varphi^5\otimes 1}_{108}\\
&+\underbrace{1\otimes \shift \varphi^1\otimes \shift \varphi^{21}\otimes (\shift \otimes \shift)\partial(\shift^{-1}\varphi^{22})\varphi^{23}\varphi^{24}\varphi^{25}\varphi^3\varphi^4\otimes \shift \varphi^5\otimes 1}_{178}\\
&
\+\underbrace{1\otimes \shift \varphi^1\otimes \varphi^{21}\hat{\omega}\varphi^{22}\varphi^{23} \varphi^{24}\otimes \shift \varphi^{25}\otimes \shift\varphi^3\otimes \varphi^4\varphi^5}_{97}\\
&+\underbrace{1\otimes \shift \varphi^1\otimes \varphi^{21}\varphi^{22}\hat{\omega}\varphi^{23}\varphi^{24}\otimes \shift \varphi^{25}\otimes \shift \varphi^3\otimes \varphi^4\varphi^5}_{126}\\
&
+\underbrace{1\otimes \shift \varphi^1\otimes \varphi^{21}(\shift \otimes \shift)\partial(\shift^{-1}\varphi^{22})\varphi^{23}\varphi^{24}\otimes \shift \varphi^{25}\otimes \shift \varphi^3\otimes \varphi^4\varphi^5}_y\\
&\+\underbrace{\varphi^1\varphi^2\hat{\omega}\otimes \shift \varphi^3\otimes \shift \varphi^{41}\otimes \varphi^{42}\dots \varphi^{45}\otimes \shift \varphi^5\otimes 1}_{12}
+\underbrace{\varphi^1\varphi^2\otimes \shift \varphi^3\otimes \shift \varphi^{41}\otimes \hat{\omega}\varphi^{42}\dots \varphi^{45}\otimes \shift \varphi^5\otimes 1}_{24}\\
&-\underbrace{\varphi^1\varphi^2\otimes \shift \varphi^3\otimes 1\otimes \shift \varphi^{41}\otimes \varphi^{42}\dots \varphi^{45}\otimes \shift \varphi^5\otimes 1}_{335}
-\underbrace{\varphi^1\varphi^2\varphi^{31}\varphi^{32}\otimes \shift \varphi^{33}\otimes \shift \varphi^{41}\otimes \varphi^{42}\dots \varphi^{45}\otimes \shift \varphi^5\otimes 1}_g\\
&+\underbrace{\varphi^1\varphi^2\otimes \shift \varphi^3\otimes \shift \varphi^{411}\otimes \varphi^{412}\varphi^{413}\varphi^{42}\dots \varphi^{45}\otimes \shift \varphi^5\otimes 1}_y\\
&
\+\underbrace{\varphi^1\varphi^{21}\dots \varphi^{23}\hat{\omega}\varphi^{24}\otimes \shift \varphi^{25}\otimes \shift \varphi^3\otimes \varphi^4\otimes \shift \varphi^5\otimes 1}_{105}\\
&+\underbrace{\varphi^1\varphi^{21}\dots \varphi^{24}\hat{\omega}\otimes \shift \varphi^{25}\otimes \shift \varphi^3\otimes \varphi^4\otimes \shift \varphi^5\otimes 1}_{25}\\
&
+\underbrace{\varphi^1\varphi^{21}\varphi^{22}\varphi^{23}(\shift \otimes \shift)\partial(\shift^{-1}\varphi^{24})\otimes \shift \varphi^{25}\otimes \shift \varphi^3\otimes \varphi^4\otimes \shift \varphi^5\otimes 1}_{154}\\
&\+\underbrace{1\otimes \shift \varphi^1\otimes \hat{\omega}\varphi^2\otimes \shift \varphi^3\otimes \shift \varphi^{41}\otimes \varphi^{42}\dots \varphi^{45}\varphi^5}_{13}
+\underbrace{1\otimes \shift \varphi^1\otimes \varphi^2\hat{\omega}\otimes \shift \varphi^3\otimes \shift \varphi^{41}\otimes \varphi^{42}\dots \varphi^{45}\varphi^5}_{26}\\
&+\underbrace{1\otimes \shift \varphi^1\otimes (\shift \otimes \shift)\partial(\shift^{-1}\varphi^2)\otimes \shift \varphi^3\otimes \shift \varphi^{41}\otimes \varphi^{42}\dots \varphi^{45}\varphi^5}_{185}\\
&
\+\underbrace{\varphi^1\otimes \shift \varphi^{21}\otimes \hat{\omega}\varphi^{22}\varphi^{23}\varphi^{24}\otimes \shift \varphi^{25}\otimes \shift \varphi^3\otimes \varphi^4\varphi^5}_{14}\\
&+\underbrace{\varphi^1\otimes \shift \varphi^{21}\otimes \varphi^{22}\hat{\omega}\varphi^{23}\varphi^{24}\otimes \shift \varphi^{25}\otimes \shift \varphi^3\otimes \varphi^4\varphi^5}_{129}\\
&
+\underbrace{\varphi^1\otimes \shift \varphi^{21}\otimes (\shift \otimes \shift)\partial(\shift^{-1}\varphi^{22})\varphi^{23}\varphi^{24}\otimes \shift \varphi^{25}\otimes \shift\varphi^3\otimes \varphi^4\varphi^5}_{237}\\
&\+\underbrace{\varphi^1\varphi^2\hat{\omega}\otimes \shift \varphi^3\otimes \shift \varphi^{41}\otimes \varphi^{42}\dots \varphi^{44}\otimes \shift \varphi^{45}\otimes \varphi^5}_{15}
+\underbrace{\varphi^1\varphi^2\otimes \shift \varphi^3\otimes \shift \varphi^{41}\otimes \hat{\omega}\varphi^{42}\dots \varphi^{44}\otimes \shift \varphi^{45}\otimes \varphi^5}_{27}\\
&-\underbrace{\varphi^1\varphi^2\otimes \shift \varphi^3\otimes 1\otimes \shift \varphi^{41}\otimes \varphi^{42}\dots \varphi^{44}\otimes \shift \varphi^{45}\otimes \varphi^5}_{349}\\
&
-\underbrace{\varphi^1\varphi^2\varphi^{31}\varphi^{32}\otimes \shift \varphi^{33}\otimes \shift \varphi^{41}\otimes \varphi^{42}\dots \varphi^{44}\otimes \shift \varphi^{45}\otimes \varphi^5}_g\\
&+\underbrace{\varphi^1\varphi^2\otimes \shift \varphi^3\otimes \shift \varphi^{411}\otimes \varphi^{412}\varphi^{413}\varphi^{42}\dots \varphi^{44}\otimes \shift \varphi^{45}\otimes \varphi^5}_y\\
&
\-\underbrace{\varphi^1\varphi^2\hat{\omega}\otimes \shift (\varphi^3\varphi^{41})\otimes \varphi^{42}\dots \varphi^{44}\otimes \shift \varphi^{45}\otimes \shift \varphi^5\otimes 1}_{16}\\
&-\underbrace{\varphi^1\varphi^2\otimes \shift (\varphi^3\varphi^{41})\otimes \hat{\omega}\varphi^{42}\dots \varphi^{44}\otimes \shift \varphi^{45}\otimes \shift \varphi^5\otimes 1}_{28}\\
&
+\underbrace{\varphi^1\varphi^2\varphi^{31}\varphi^{32}\otimes \shift (\varphi^{33}\varphi^{41})\otimes \varphi^{42}\dots \varphi^{44}\otimes \shift \varphi^{45}\otimes \shift \varphi^5\otimes 1}_g\\
&+\underbrace{\varphi^1\varphi^2\varphi^3\varphi^{411}\varphi^{412}\otimes \shift \varphi^{413}\otimes \varphi^{42}\dots \varphi^{44}\otimes \shift \varphi^{45}\otimes \shift \varphi^5\otimes 1}_y\\
&
-\underbrace{\varphi^1\varphi^2\otimes \shift \varphi^{31}\otimes \varphi^{32}\varphi^{33}\varphi^{41}\dots \varphi^{44}\otimes \shift \varphi^{45}\otimes \shift \varphi^5\otimes 1}_g\\
&-\underbrace{\varphi^1\varphi^2\otimes \shift (\varphi^3\varphi^{411})\otimes \varphi^{412}\varphi^{413}\varphi^{42}\dots \varphi^{44}\otimes \shift \varphi^{45}\otimes \shift \varphi^5\otimes 1}_y\\
&\-\underbrace{1\otimes \shift \varphi^1\otimes \shift \varphi^{21}\otimes \hat{\omega}\varphi^{22}\dots \varphi^{24}\otimes \shift (\varphi^{25}\varphi^3)\otimes \varphi^4\varphi^5}_{17}\\
&
-\underbrace{1\otimes \shift \varphi^1\otimes \shift \varphi^{21}\otimes \varphi^{22}\hat{\omega}\varphi^{23}\varphi^{24}\otimes \shift (\varphi^{25}\varphi^3)\otimes \varphi^4\varphi^5}_{143}\\
&-\underbrace{1\otimes \shift \varphi^1\otimes \shift \varphi^{21}\otimes (\shift \otimes \shift)\partial(\shift^{-1}\varphi^{22})\varphi^{23}\varphi^{24}\otimes \shift (\varphi^{25}\varphi^3)\otimes \varphi^4\varphi^5}_{225}\\
&
\-\underbrace{\varphi^1\varphi^{21}\varphi^{22}\hat{\omega}\otimes \shift \varphi^{23}\otimes \varphi^{24}\otimes \shift \varphi^{25}\otimes \shift \varphi^3\otimes \varphi^4\varphi^5}_{18}\\
&-\underbrace{\varphi^1\varphi^{21}\varphi^{22}\otimes \shift \varphi^{23}\otimes \hat{\omega}\varphi^{24}\otimes \shift \varphi^{25}\otimes \shift \varphi^3\otimes \varphi^4\varphi^5}_{29}
+\underbrace{\varphi^1\varphi^{21}\varphi^{22}\varphi^{231}\varphi^{232}\otimes \shift \varphi^{233}\otimes \varphi^{24}\otimes \shift \varphi^{25}\otimes \shift \varphi^3\otimes \varphi^4\varphi^5}_y\\
&-\underbrace{\varphi^1\varphi^{21}\varphi^{22}\otimes \shift \varphi^{231}\otimes \varphi^{232}\varphi^{233}\varphi^{24}\otimes \shift \varphi^{25}\otimes \shift \varphi^3\otimes \varphi^4\varphi^5}_y
\-\underbrace{\varphi^1\varphi^2\hat{\omega}\otimes \shift \varphi^3\otimes \shift \varphi^{41}\otimes \varphi^{42}\otimes \shift \varphi^{43}\otimes \varphi^{44}\varphi^{45}\varphi^5}_{19}\\
&-\underbrace{\varphi^1\varphi^2\otimes \shift \varphi^3\otimes \shift \varphi^{41}\otimes \hat{\omega}\varphi^{42}\otimes \shift \varphi^{43}\otimes \varphi^{44}\varphi^{45}\varphi^5}_{30}
+\underbrace{\varphi^1\varphi^2\otimes \shift \varphi^3\otimes 1\otimes \shift \varphi^{41}\otimes \varphi^{42}\otimes \shift \varphi^{43}\otimes \varphi^{44}\varphi^{45}\varphi^5}_{336}\\
&+\underbrace{\varphi^1\varphi^2\varphi^{31}\varphi^{32}\otimes \shift \varphi^{33}\otimes \shift \varphi^{41}\otimes \varphi^{42}\otimes \shift \varphi^{43}\otimes \varphi^{44}\varphi^{45}\varphi^5}_g\\
&
-\underbrace{\varphi^1\varphi^2\otimes \shift \varphi^3\otimes \shift \varphi^{411}\otimes \varphi^{412}\varphi^{413}\varphi^{42}\otimes \shift \varphi^{43}\otimes \varphi^{44}\varphi^{45}\varphi^5}_y\\
&\+\underbrace{1\otimes \shift \varphi^1\otimes \shift \varphi^{21}\otimes \hat{\omega}\varphi^{22}\otimes \shift \varphi^{23}\otimes \varphi^{24}\varphi^{25}\varphi^3\varphi^4\varphi^5}_{20}
+\underbrace{1\otimes \shift \varphi^1\otimes \shift \varphi^{21}\otimes \varphi^{22}\hat{\omega}\otimes \shift \varphi^{23}\otimes \varphi^{24}\varphi^{25}\varphi^3\varphi^4\varphi^5}_{31}\\
&+\underbrace{1\otimes \shift \varphi^1\otimes \shift \varphi^{21}\otimes (\shift \otimes \shift)\partial(\shift^{-1}\varphi^{22})\otimes \shift \varphi^{23}\otimes \varphi^{24}\varphi^{25}\varphi^3\varphi^4\varphi^5}_{247}\\
&
\+\underbrace{\varphi^1\varphi^2\varphi^3\varphi^{41}\hat{\omega}\varphi^{42}\otimes \shift \varphi^{43}\otimes \varphi^{44}\otimes \shift \varphi^{45}\otimes \shift \varphi^5\otimes 1}_{136}\\
&+\underbrace{\varphi^1\varphi^2\varphi^3\varphi^{41}\varphi^{42}\hat{\omega}\otimes \shift \varphi^{43}\otimes \varphi^{44}\otimes \shift \varphi^{45}\otimes \shift \varphi^5\otimes 1}_{32}\\
&
+\underbrace{\varphi^1\varphi^2\varphi^3\varphi^{41}(\shift \otimes \shift)\partial(\shift^{-1}\varphi^{42})\otimes \shift \varphi^{43}\otimes \varphi^{44}\otimes \shift \varphi^{45}\otimes \shift \varphi^5\otimes 1}_{263}\\
&\+\underbrace{\varphi^1\varphi^{21}\varphi^{22}\varphi^{23}\hat{\omega}\varphi^{24}\otimes \shift \varphi^{25}\otimes \shift \varphi^3\otimes \shift \varphi^{41}\otimes \varphi^{42}\dots \varphi^{45}\varphi^5}_{134}\\
&
+\underbrace{\varphi^1\varphi^{21}\dots \varphi^{24}\hat{\omega}\otimes \shift \varphi^{25}\otimes \shift \varphi^3\otimes \shift \varphi^{41}\otimes \varphi^{42}\dots \varphi^{45}\varphi^5}_{33}\\
&+\underbrace{\varphi^1\varphi^{21}\varphi^{22}\varphi^{23}(\shift \otimes \shift)\partial(\shift^{-1}\varphi^{24})\otimes \shift \varphi^{25}\otimes \shift \varphi^3\otimes \shift \varphi^{41}\otimes \varphi^{42}\dots \varphi^{45}\varphi^5}_{261}\\
&
\+\underbrace{\varphi^1\varphi^2\varphi^3\varphi^{41}\hat{\omega}\varphi^{42}\varphi^{43}\varphi^{441}\dots \varphi^{444}\otimes \shift \varphi^{445}\otimes \shift \varphi^{45}\otimes \shift \varphi^5\otimes 1}_{311}\\
&+\underbrace{\varphi^1\varphi^2\varphi^3\varphi^{41}\varphi^{42}\hat{\omega}\varphi^{43}\varphi^{441}\dots \varphi^{444}\otimes \shift \varphi^{445}\otimes \shift \varphi^{45}\otimes \shift \varphi^5\otimes 1}_{315}\\
&
+\underbrace{\varphi^1\varphi^2\varphi^3\varphi^{41}(\shift \otimes \shift)\partial(\shift^{-1}\varphi^{42})\varphi^{43}\varphi^{441}\dots \varphi^{444}\otimes \shift \varphi^{445}\otimes \shift \varphi^{45}\otimes \shift \varphi^5\otimes 1}_g\\
&\+\underbrace{1\otimes \shift \varphi^1\otimes \shift \varphi^{21}\otimes \shift \varphi^{221}\otimes \hat{\omega}\varphi^{222}\dots \varphi^{225}\varphi^{23}\dots \varphi^{25}\varphi^3\varphi^4\varphi^5}_{21}\\
&+\underbrace{1\otimes \shift \varphi^1\otimes \shift \varphi^{21}\otimes \shift \varphi^{221}\otimes \varphi^{222}\hat{\omega}\varphi^{223}\dots \varphi^{225}\varphi^{23}\dots \varphi^{25}\varphi^3\varphi^4\varphi^5}_{317}\\
&+\underbrace{1\otimes \shift \varphi^1\otimes \shift \varphi^{21}\otimes \shift \varphi^{221}\otimes (\shift \otimes \shift)\partial(\shift^{-1}\varphi^{222})\varphi^{223}\dots \varphi^{225}\varphi^{23}\dots \varphi^{25}\varphi^3\varphi^4\varphi^5}_g\\
&\-\underbrace{\varphi^1\varphi^{21}\dots \varphi^{23}\varphi^{241}\hat{\omega}\varphi^{242}\dots \varphi^{244}\otimes \shift \varphi^{245}\otimes \shift \varphi^{25}\otimes \shift \varphi^3\otimes \varphi^4\varphi^5}_{313}\\
&-\underbrace{\varphi^1\varphi^{21}\dots \varphi^{23}\varphi^{241}\varphi^{242}\hat{\omega}\varphi^{243}\varphi^{244}\otimes \shift \varphi^{245}\otimes \shift \varphi^{25}\otimes \shift \varphi^3\otimes \varphi^4\varphi^5}_{319}\\
&-\underbrace{\varphi^1\varphi^{21}\dots \varphi^{23}\varphi^{241}(\shift \otimes \shift)\partial(\shift^{-1}\varphi^{242})\varphi^{243}\varphi^{244}\otimes \shift \varphi^{245}\otimes \shift \varphi^{25}\otimes \shift \varphi^3\otimes \varphi^4\varphi^5}_y\\
&\-\underbrace{\varphi^1\varphi^2\hat{\omega}\otimes \shift \varphi^3\otimes \shift \varphi^{41}\otimes \shift \varphi^{421}\otimes \varphi^{422}\dots \varphi^{425}\varphi^{43}\dots \varphi^{45}\varphi^5}_{22}\\
&-\underbrace{\varphi^1\varphi^2\otimes \shift \varphi^3\otimes \shift \varphi^{41}\otimes \shift \varphi^{421}\otimes \hat{\omega}\varphi^{422}\dots \varphi^{425}\varphi^{43}\dots \varphi^{45}\varphi^5}_{314}\\
&+\underbrace{\varphi^1\varphi^2\otimes \shift \varphi^3\otimes 1\otimes \shift \varphi^{41}\otimes \shift \varphi^{421}\otimes \varphi^{422}\dots \varphi^{425}\varphi^{43}\dots \varphi^{45}\varphi^5}_{337}\\
&+\underbrace{\varphi^1\varphi^2\otimes \shift \varphi^3\otimes \shift \varphi^{41}\otimes 1\otimes \shift \varphi^{421}\otimes \varphi^{422}\dots \varphi^{425}\varphi^{43}\dots \varphi^{45}\varphi^5}_{338}\\
&+\underbrace{\varphi^1\varphi^2\varphi^{31}\varphi^{32}\otimes \shift \varphi^{33}\otimes \shift \varphi^{41}\otimes \shift \varphi^{421}\otimes \varphi^{422}\dots \varphi^{425}\varphi^{43}\dots \varphi^{45}\varphi^5}_g\\
&-\underbrace{\varphi^1\varphi^2\otimes \shift \varphi^3\otimes \shift \varphi^{41}\otimes \shift \varphi^{4211}\otimes \varphi^{4212}\varphi^{4213}\varphi^{422}\dots \varphi^{425}\varphi^{43}\dots \varphi^{45}\varphi^5}_y
\end{align*}
and
\begin{align*}
&(1^{\otimes 2}\otimes \mu_2\otimes 1^{\otimes 2})\mu_5(\varphi)\\
=\,&-\underbrace{1\otimes \shift \varphi^1\otimes \varphi^2\hat{\omega}\otimes \shift \varphi^3\otimes \varphi^4\otimes \shift \varphi^5\otimes 1}_{23}
-\underbrace{1\otimes \shift \varphi^1\otimes \varphi^2\otimes \shift \varphi^3\otimes \hat{\omega}\varphi^4\otimes \shift \varphi^5\otimes 1}_{34}\\
&+\underbrace{1\otimes \shift \varphi^1\otimes \varphi^2\varphi^{31}\varphi^{32}\otimes \shift \varphi^{33}\otimes \varphi^4\otimes \shift \varphi^5\otimes 1}_y
-\underbrace{1\otimes \shift \varphi^1\otimes \varphi^2\otimes \shift \varphi^{31}\otimes \varphi^{32}\varphi^{33}\varphi^4\otimes \shift \varphi^5\otimes 1}_y\\
&\-\underbrace{1\otimes \shift \varphi^1\otimes \varphi^2\varphi^3\varphi^{41}\hat{\omega}\varphi^{42}\varphi^{43}\varphi^{44}\otimes \shift \varphi^{45}\otimes \shift \varphi^5\otimes 1}_{86}
-\underbrace{1\otimes \shift \varphi^1\otimes \varphi^2\varphi^3\varphi^{41}\varphi^{42}\hat{\omega}\varphi^{43}\varphi^{44}\otimes \shift \varphi^{45}\otimes \shift \varphi^5\otimes 1}_{87}\\
&-\underbrace{1\otimes \shift \varphi^1\otimes \varphi^2\varphi^3\varphi^{41}(\shift \otimes \shift)\partial(\shift^{-1}\varphi^{42})\varphi^{43}\varphi^{44}\otimes \shift \varphi^{45}\otimes \shift \varphi^5\otimes 1}_y\\
&
\-\underbrace{1\otimes \shift \varphi^1\otimes \shift \varphi^{21}\otimes \varphi^{22}\varphi^{23}\hat{\omega}\varphi^{24}\varphi^{25}\varphi^3\varphi^4\otimes \shift \varphi^5\otimes 1}_{109}\\
&-\underbrace{1\otimes \shift \varphi^1\otimes \shift \varphi^{21}\otimes \varphi^{22}\varphi^{23}\varphi^{24}\hat{\omega}\varphi^{25}\varphi^3\varphi^4\otimes \shift \varphi^5\otimes 1}_{99}\\
&
-\underbrace{1\otimes \shift \varphi^1\otimes \shift \varphi^{21}\otimes \varphi^{22}\varphi^{23}(\shift \otimes \shift)\partial(\shift^{-1}\varphi^{24})\varphi^{25}\varphi^3\varphi^4\otimes \shift \varphi^5\otimes 1}_y\\
&\-\underbrace{1\otimes \shift \varphi^1\otimes \varphi^{21}\varphi^{22} \varphi^{23}\hat{\omega}\varphi^{24}\otimes \shift \varphi^{25}\otimes \shift\varphi^3\otimes \varphi^4\varphi^5}_{127}
-\underbrace{1\otimes \shift \varphi^1\otimes \varphi^{21}\dots \varphi^{24}\hat{\omega}\otimes \shift \varphi^{25}\otimes \shift \varphi^3\otimes \varphi^4\varphi^5}_{35}\\
&-\underbrace{1\otimes \shift \varphi^1\otimes \varphi^{21}\varphi^{22} \varphi^{23}(\shift \otimes \shift)\partial(\shift^{-1}\varphi^{24})\otimes \shift \varphi^{25}\otimes \shift\varphi^3\otimes \varphi^4\varphi^5}_{253}\\
&
\-\underbrace{\varphi^1\varphi^2\otimes \shift \varphi^3\otimes \shift \varphi^{41}\otimes \hat{\omega}\varphi^{42}\dots \varphi^{45}\otimes \shift \varphi^5\otimes 1}_{24}\\
&-\underbrace{\varphi^1\varphi^2\otimes \shift \varphi^3\otimes \shift \varphi^{41}\otimes \varphi^{42}\hat{\omega}\varphi^{43}\varphi^{44} \varphi^{45}\otimes \shift \varphi^5\otimes 1}_{113}\\
&
-\underbrace{\varphi^1\varphi^2\otimes \shift \varphi^3\otimes \shift \varphi^{41}\otimes (\shift \otimes \shift)\partial(\shift^{-1}\varphi^{42})\varphi^{43}\dots \varphi^{45}\otimes \shift \varphi^5\otimes 1}_{234}\\
&\-\underbrace{\varphi^1\varphi^{21}\dots \varphi^{24}\hat{\omega}\otimes \shift \varphi^{25}\otimes \shift \varphi^3\otimes \varphi^4\otimes \shift \varphi^5\otimes 1}_{25}
-\underbrace{\varphi^1\varphi^{21}\dots \varphi^{24}\otimes \shift \varphi^{25}\otimes \shift \varphi^3\otimes \hat{\omega}\varphi^4\otimes \shift \varphi^5\otimes 1}_{36}\\
&+\underbrace{\varphi^1\varphi^{21}\dots \varphi^{24}\otimes \shift \varphi^{25}\otimes 1\otimes \shift \varphi^3\otimes \varphi^4\otimes \shift \varphi^5\otimes 1}_{342}
+\underbrace{\varphi^1\varphi^{21}\dots \varphi^{24}\varphi^{251}\varphi^{252}\otimes \shift \varphi^{253}\otimes \shift \varphi^3\otimes \varphi^4\otimes \shift \varphi^5\otimes 1}_y\\
&-\underbrace{\varphi^1\varphi^{21}\dots \varphi^{24}\otimes \shift \varphi^{25}\otimes \shift \varphi^{31}\otimes \varphi^{32}\varphi^{33}\varphi^4\otimes \shift \varphi^5\otimes 1}_g
\-\underbrace{1\otimes \shift \varphi^1\otimes \varphi^2\hat{\omega}\otimes \shift \varphi^3\otimes \shift \varphi^{41}\otimes \varphi^{42}\dots \varphi^{45}\varphi^5}_{26}\\
&-\underbrace{1\otimes \shift \varphi^1\otimes \varphi^2\otimes \shift \varphi^3\otimes \shift \varphi^{41}\otimes \hat{\omega}\varphi^{42}\dots \varphi^{45}\varphi^5}_{37}
+\underbrace{1\otimes \shift \varphi^1\otimes \varphi^2\otimes \shift \varphi^3\otimes 1\otimes \shift \varphi^{41}\otimes \varphi^{42}\dots \varphi^{45}\varphi^5}_{339}\\
&+\underbrace{1\otimes \shift \varphi^1\otimes \varphi^2\varphi^{31}\varphi^{32}\otimes \shift \varphi^{33}\otimes \shift \varphi^{41}\otimes \varphi^{42}\dots \varphi^{45}\varphi^5}_g\\
&
-\underbrace{1\otimes \shift \varphi^1\otimes \varphi^2\otimes \shift \varphi^3\otimes \shift \varphi^{411}\otimes \varphi^{412}\varphi^{413}\varphi^{42}\dots \varphi^{45}\varphi^5}_y\\
&\-\underbrace{\varphi^1\otimes \shift \varphi^{21}\otimes \varphi^{22}\varphi^{23}\hat{\omega}\varphi^{24}\otimes \shift \varphi^{25}\otimes \shift \varphi^3\otimes \varphi^4\varphi^5}_{132}
-\underbrace{\varphi^1\otimes \shift \varphi^{21}\otimes \varphi^{22}\dots \varphi^{24}\hat{\omega}\otimes \shift \varphi^{25}\otimes \shift \varphi^3\otimes \varphi^4\varphi^5}_{38}\\
&-\underbrace{\varphi^1\otimes \shift \varphi^{21}\otimes \varphi^{22}\varphi^{23}(\shift \otimes \shift)\partial(\shift^{-1}\varphi^{24})\otimes \shift \varphi^{25}\otimes \shift \varphi^3\otimes \varphi^4\varphi^5}_{254}\\
&
\-\underbrace{\varphi^1\varphi^2\otimes \shift \varphi^3\otimes \shift \varphi^{41}\otimes \hat{\omega}\varphi^{42}\dots \varphi^{44}\otimes \shift \varphi^{45}\otimes \varphi^5}_{27}\\
&-\underbrace{\varphi^1\varphi^2\otimes \shift \varphi^3\otimes \shift \varphi^{41}\otimes \varphi^{42}\hat{\omega}\varphi^{43}\varphi^{44}\otimes \shift \varphi^{45}\otimes \varphi^5}_{115}\\
&
-\underbrace{\varphi^1\varphi^2\otimes \shift \varphi^3\otimes \shift \varphi^{41}\otimes (\shift \otimes \shift)\partial(\shift^{-1}\varphi^{42})\varphi^{43}\varphi^{44}\otimes \shift \varphi^{45}\otimes \varphi^5}_{232}\\
&\+\underbrace{\varphi^1\varphi^2\otimes \shift (\varphi^3\varphi^{41})\otimes \hat{\omega}\varphi^{42}\dots \varphi^{44}\otimes \shift \varphi^{45}\otimes \shift \varphi^5\otimes 1}_{28}
+\underbrace{\varphi^1\varphi^2\otimes \shift (\varphi^3\varphi^{41})\otimes \varphi^{42}\hat{\omega}\varphi^{43}\varphi^{44}\otimes \shift \varphi^{45}\otimes \shift \varphi^5\otimes 1}_{137}\\
&+\underbrace{\varphi^1\varphi^2\otimes \shift (\varphi^3\varphi^{41})\otimes (\shift \otimes \shift)\partial(\shift^{-1}\varphi^{42})\varphi^{43}\varphi^{44}\otimes \shift \varphi^{45}\otimes \shift \varphi^5\otimes 1}_{233}\\
&
\+\underbrace{1\otimes \shift \varphi^1\otimes \shift \varphi^{21}\otimes \varphi^{22}\varphi^{23}\hat{\omega}\varphi^{24}\otimes \shift (\varphi^{25}\varphi^3)\otimes \varphi^4\varphi^5}_{144}\\
&+\underbrace{1\otimes \shift \varphi^1\otimes \shift \varphi^{21}\otimes \varphi^{22}\varphi^{23}\varphi^{24}\hat{\omega}\otimes \shift (\varphi^{25}\varphi^3)\otimes \varphi^4\varphi^5}_{39}\\
&
+\underbrace{1\otimes \shift \varphi^1\otimes \shift \varphi^{21}\otimes \varphi^{22}\varphi^{23}(\shift \otimes \shift)\partial(\shift^{-1}\varphi^{24})\otimes \shift (\varphi^{25}\varphi^3)\otimes \varphi^4\varphi^5}_{255}\\
&
\+\underbrace{\varphi^1\varphi^{21}\varphi^{22}\otimes \shift \varphi^{23}\otimes \hat{\omega}\varphi^{24}\otimes \shift \varphi^{25}\otimes \shift \varphi^3\otimes \varphi^4\varphi^5}_{29}
+\underbrace{\varphi^1\varphi^{21}\varphi^{22}\otimes \shift \varphi^{23}\otimes \varphi^{24}\hat{\omega}\otimes \shift \varphi^{25}\otimes \shift \varphi^3\otimes \varphi^4\varphi^5}_{40}\\
&+\underbrace{\varphi^1\varphi^{21}\varphi^{22}\otimes \shift \varphi^{23}\otimes (\shift \otimes \shift)\partial(\shift^{-1}\varphi^{24})\otimes \shift \varphi^{25}\otimes \shift \varphi^3\otimes \varphi^4\varphi^5}_{267}\\
&
\+\underbrace{\varphi^1\varphi^2\otimes \shift \varphi^3\otimes \shift \varphi^{41}\otimes \hat{\omega}\varphi^{42}\otimes \shift \varphi^{43}\otimes \varphi^{44}\varphi^{45}\varphi^5}_{30}\\
&+\underbrace{\varphi^1\varphi^2\otimes \shift \varphi^3\otimes \shift \varphi^{41}\otimes \varphi^{42}\hat{\omega}\otimes \shift \varphi^{43}\otimes \varphi^{44}\varphi^{45}\varphi^5}_{151}\\
&
+\underbrace{\varphi^1\varphi^2\otimes \shift \varphi^3\otimes \shift \varphi^{41}\otimes (\shift \otimes \shift)\partial(\shift^{-1}\varphi^{42})\otimes \shift \varphi^{43}\otimes \varphi^{44}\varphi^{45}\varphi^5}_{273}\\
&\-\underbrace{1\otimes \shift \varphi^1\otimes \shift \varphi^{21}\otimes \varphi^{22}\hat{\omega}\otimes \shift \varphi^{23}\otimes \varphi^{24}\varphi^{25}\varphi^3\varphi^4\varphi^5}_{31}
-\underbrace{1\otimes \shift \varphi^1\otimes \shift \varphi^{21}\otimes \varphi^{22}\otimes \shift \varphi^{23}\otimes \hat{\omega}\varphi^{24}\varphi^{25}\varphi^3\varphi^4\varphi^5}_{41}\\
&+\underbrace{1\otimes \shift \varphi^1\otimes \shift \varphi^{21}\otimes \varphi^{22}\varphi^{231}\varphi^{232}\otimes \shift \varphi^{233}\otimes \varphi^{24}\varphi^{25}\varphi^3\varphi^4\varphi^5}_y\\
&
-\underbrace{1\otimes \shift \varphi^1\otimes \shift \varphi^{21}\otimes \varphi^{22}\otimes \shift \varphi^{231}\otimes \varphi^{232}\varphi^{233}\varphi^{24}\varphi^{25}\varphi^3\varphi^4\varphi^5}_y\\
&\-\underbrace{\varphi^1\varphi^2\varphi^3\varphi^{41}\varphi^{42}\hat{\omega}\otimes \shift \varphi^{43}\otimes \varphi^{44}\otimes \shift \varphi^{45}\otimes \shift \varphi^5\otimes 1}_{32}
-\underbrace{\varphi^1\varphi^2\varphi^3\varphi^{41}\varphi^{42}\otimes \shift \varphi^{43}\otimes \hat{\omega}\varphi^{44}\otimes \shift \varphi^{45}\otimes \shift \varphi^5\otimes 1}_{42}\\
&+\underbrace{\varphi^1\varphi^2\varphi^3\varphi^{41}\varphi^{42}\varphi^{431}\varphi^{432}\otimes \shift \varphi^{433}\otimes \varphi^{44}\otimes \shift \varphi^{45}\otimes \shift \varphi^5\otimes 1}_y\\
&
-\underbrace{\varphi^1\varphi^2\varphi^3\varphi^{41}\varphi^{42}\otimes \shift \varphi^{431}\otimes \varphi^{432}\varphi^{433}\varphi^{44}\otimes \shift \varphi^{45}\otimes \shift \varphi^5\otimes 1}_y\\
&\-\underbrace{\varphi^1\varphi^{21}\dots \varphi^{24}\hat{\omega}\otimes \shift \varphi^{25}\otimes \shift \varphi^3\otimes \shift \varphi^{41}\otimes \varphi^{42}\dots \varphi^{45}\varphi^5}_{33}\\
&
-\underbrace{\varphi^1\varphi^{21}\dots \varphi^{24}\otimes \shift \varphi^{25}\otimes \shift \varphi^3\otimes \shift \varphi^{41}\otimes \hat{\omega}\varphi^{42}\dots \varphi^{45}\varphi^5}_{43}\\
&+\underbrace{\varphi^1\varphi^{21}\dots \varphi^{24}\otimes \shift \varphi^{25}\otimes 1\otimes \shift \varphi^3\otimes \shift \varphi^{41}\otimes \varphi^{42}\dots \varphi^{45}\varphi^5}_{341}\\
&
+\underbrace{\varphi^1\varphi^{21}\dots \varphi^{24}\otimes \shift \varphi^{25}\otimes \shift \varphi^3\otimes 1\otimes \shift \varphi^{41}\otimes \varphi^{42}\dots \varphi^{45}\varphi^5}_{340}\\
&+\underbrace{\varphi^1\varphi^{21}\dots \varphi^{24}\varphi^{251}\varphi^{252}\otimes \shift \varphi^{253}\otimes \shift \varphi^3\otimes \shift \varphi^{41}\otimes \varphi^{42}\dots \varphi^{45}\varphi^5}_y\\
&-\underbrace{\varphi^1\varphi^{21}\dots \varphi^{24}\otimes \shift \varphi^{25}\otimes \shift \varphi^3\otimes \shift \varphi^{411}\otimes \varphi^{412}\varphi^{413}\varphi^{42}\dots \varphi^{45}\varphi^5}_y\\
&\-\underbrace{\varphi^1\varphi^2\varphi^3\varphi^{41}\dots \varphi^{43}\varphi^{441}\hat{\omega}\varphi^{442}\dots \varphi^{444}\otimes \shift \varphi^{445}\otimes \shift \varphi^{45}\otimes \shift \varphi^5\otimes 1}_{316}\\
&-\underbrace{\varphi^1\varphi^2\varphi^3\varphi^{41}\dots \varphi^{43}\varphi^{441}\varphi^{442}\hat{\omega}\varphi^{443}\varphi^{444}\otimes \shift \varphi^{445}\otimes \shift \varphi^{45}\otimes \shift \varphi^5\otimes 1}_{322}\\
&-\underbrace{\varphi^1\varphi^2\varphi^3\varphi^{41}\dots \varphi^{43}\varphi^{441}(\shift \otimes \shift)\partial(\shift^{-1}\varphi^{442})\varphi^{443}\varphi^{444}\otimes \shift \varphi^{445}\otimes \shift \varphi^{45}\otimes \shift \varphi^5\otimes 1}_y\\
&\-\underbrace{1\otimes \shift \varphi^1\otimes \shift \varphi^{21}\otimes \shift \varphi^{221}\otimes \varphi^{222}\varphi^{223}\hat{\omega}\varphi^{224}\varphi^{225}\varphi^{23}\dots \varphi^{25}\varphi^3\varphi^4\varphi^5}_{318}\\
&-\underbrace{1\otimes \shift \varphi^1\otimes \shift \varphi^{21}\otimes \shift \varphi^{221}\otimes \varphi^{222}\dots \varphi^{224}\hat{\omega}\varphi^{225}\varphi^{23}\dots \varphi^{25}\varphi^3\varphi^4\varphi^5}_{324}\\
&-\underbrace{1\otimes \shift \varphi^1\otimes \shift \varphi^{21}\otimes \shift \varphi^{221}\otimes \varphi^{222}\varphi^{223}(\shift \otimes \shift)\partial(\shift^{-1}\varphi^{224})\varphi^{225}\varphi^{23}\dots \varphi^{25}\varphi^3\varphi^4\varphi^5}_y\\
&\+\underbrace{\varphi^1\varphi^{21}\dots \varphi^{23}\varphi^{241}\dots \varphi^{243}\hat{\omega}\varphi^{244}\otimes \shift \varphi^{245}\otimes \shift \varphi^{25}\otimes \shift \varphi^3\otimes \varphi^4\varphi^5}_{320}\\
&+\underbrace{\varphi^1\varphi^{21}\dots \varphi^{23}\varphi^{241}\dots \varphi^{244}\hat{\omega}\otimes \shift \varphi^{245}\otimes \shift \varphi^{25}\otimes \shift \varphi^3\otimes \varphi^4\varphi^5}_{321}\\
&+\underbrace{\varphi^1\varphi^{21}\dots \varphi^{23}\varphi^{241}\dots \varphi^{243}(\shift \otimes \shift)\partial(\shift^{-1}\varphi^{244})\otimes \shift \varphi^{245}\otimes \shift \varphi^{25}\otimes \shift \varphi^3\otimes \varphi^4\varphi^5}_{291}\\
&\+\underbrace{\varphi^1\varphi^2\otimes \shift \varphi^3\otimes \shift \varphi^{41}\otimes \shift \varphi^{421}\otimes \hat{\omega}\varphi^{422}\dots \varphi^{425}\varphi^{43}\dots \varphi^{45}\varphi^5}_{314}\\
&+\underbrace{\varphi^1\varphi^2\otimes \shift \varphi^3\otimes \shift \varphi^{41}\otimes \shift \varphi^{421}\otimes \varphi^{422}\hat{\omega}\varphi^{423}\varphi^{424}\varphi^{425}\varphi^{43}\dots \varphi^{45}\varphi^5}_{326}\\
&+\underbrace{\varphi^1\varphi^2\otimes \shift \varphi^3\otimes \shift \varphi^{41}\otimes \shift \varphi^{421}\otimes (\shift \otimes \shift)\partial(\shift^{-1}\varphi^{422})\varphi^{423}\varphi^{424} \varphi^{425}\varphi^{43}\varphi^{44} \varphi^{45}\varphi^5}_{297}
\end{align*}
and
\begin{align*}
&-(1^{\otimes 3}\otimes \mu_2\otimes 1)\mu_5(\varphi)\\
=\,&+\underbrace{1\otimes \shift \varphi^1\otimes \varphi^2\otimes \shift \varphi^3\otimes \hat{\omega}\varphi^4\otimes \shift \varphi^5\otimes 1}_{34}
+\underbrace{1\otimes \shift \varphi^1\otimes \varphi^2\otimes \shift \varphi^3\otimes \varphi^4\hat{\omega}\otimes \shift \varphi^5\otimes 1}_{125}\\
&+\underbrace{1\otimes \shift \varphi^1\otimes \varphi^2\otimes \shift \varphi^3\otimes (\shift \otimes \shift)\partial(\shift^{-1}\varphi^4)\otimes \shift \varphi^5\otimes 1}_{148}
\+\underbrace{1\otimes \shift \varphi^1\otimes \varphi^2\varphi^3\varphi^{41}\varphi^{42}\varphi^{43}\hat{\omega}\varphi^{44}\otimes \shift \varphi^{45}\otimes \shift \varphi^5\otimes 1}_{89}\\
&+\underbrace{1\otimes \shift \varphi^1\otimes \varphi^2\varphi^3\varphi^{41}\dots \varphi^{44}\hat{\omega}\otimes \shift \varphi^{45}\otimes \shift \varphi^5\otimes 1}_{56}\\
&
+\underbrace{1\otimes \shift \varphi^1\otimes \varphi^2\varphi^3\varphi^{41}\varphi^{42}\varphi^{43}(\shift \otimes \shift)\partial(\shift^{-1}\varphi^{44})\otimes \shift \varphi^{45}\otimes \shift \varphi^5\otimes 1}_{227}\\
&\+\underbrace{1\otimes \shift \varphi^1\otimes \shift \varphi^{21}\otimes \varphi^{22}\dots \varphi^{25}\varphi^3\hat{\omega}\varphi^4\otimes \shift \varphi^5\otimes 1}_{100}
+\underbrace{1\otimes \shift \varphi^1\otimes \shift \varphi^{21}\otimes \varphi^{22}\dots \varphi^{25}\varphi^3\varphi^4\hat{\omega}\otimes \shift \varphi^5\otimes 1}_{57}\\
&+\underbrace{1\otimes \shift \varphi^1\otimes \shift \varphi^{21}\otimes \varphi^{22}\dots \varphi^{25}\varphi^3(\shift \otimes \shift)\partial(\shift^{-1}\varphi^4)\otimes \shift \varphi^5\otimes 1}_{181}\\
&
\+\underbrace{1\otimes \shift \varphi^1\otimes \varphi^{21}\dots \varphi^{24}\hat{\omega}\otimes \shift \varphi^{25}\otimes \shift \varphi^3\otimes \varphi^4\varphi^5}_{35}\\
&+\underbrace{1\otimes \shift \varphi^1\otimes \varphi^{21}\dots \varphi^{24}\otimes \shift \varphi^{25}\otimes \shift \varphi^3\otimes \hat{\omega}\varphi^4\varphi^5}_{58}
-\underbrace{1\otimes \shift \varphi^1\otimes \varphi^{21}\dots \varphi^{24}\otimes \shift \varphi^{25}\otimes 1\otimes \shift \varphi^3\otimes \varphi^4\varphi^5}_{347}\\
&-\underbrace{1\otimes \shift \varphi^1\otimes \varphi^{21}\dots \varphi^{24}\varphi^{251}\varphi^{252}\otimes \shift \varphi^{253}\otimes \shift \varphi^3\otimes \varphi^4\varphi^5}_y\\
&
+\underbrace{1\otimes \shift \varphi^1\otimes \varphi^{21}\dots \varphi^{24}\otimes \shift \varphi^{25}\otimes \shift \varphi^{31}\otimes \varphi^{32}\varphi^{33}\varphi^4\varphi^5}_g\\
&\+\underbrace{\varphi^1\varphi^2\otimes \shift \varphi^3\otimes \shift \varphi^{41}\otimes \varphi^{42}\varphi^{43}\hat{\omega}\varphi^{44}\varphi^{45}\otimes \shift \varphi^5\otimes 1}_{114}
+\underbrace{\varphi^1\varphi^2\otimes \shift \varphi^3\otimes \shift \varphi^{41}\otimes \varphi^{42}\varphi^{43} \varphi^{44}\hat{\omega}\varphi^{45}\otimes \shift \varphi^5\otimes 1}_{119}\\
&+\underbrace{\varphi^1\varphi^2\otimes \shift \varphi^3\otimes \shift \varphi^{41}\otimes \varphi^{42}\varphi^{43}(\shift \otimes \shift)\partial(\shift^{-1}\varphi^{44})\varphi^{45}\otimes \shift \varphi^5\otimes 1}_y\\
&
\+\underbrace{\varphi^1\varphi^{21}\dots \varphi^{24}\otimes \shift \varphi^{25}\otimes \shift \varphi^3\otimes \hat{\omega}\varphi^4\otimes \shift \varphi^5\otimes 1}_{36}\\
&+\underbrace{\varphi^1\varphi^{21}\dots \varphi^{24}\otimes \shift \varphi^{25}\otimes \shift \varphi^3\otimes \varphi^4\hat{\omega}\otimes \shift \varphi^5\otimes 1}_{59}\\
&
+\underbrace{\varphi^1\varphi^{21}\dots \varphi^{24}\otimes \shift \varphi^{25}\otimes \shift \varphi^3\otimes (\shift \otimes \shift)\partial(\shift^{-1}\varphi^4)\otimes \shift \varphi^5\otimes 1}_{166}\\
&\+\underbrace{1\otimes \shift \varphi^1\otimes \varphi^2\otimes \shift \varphi^3\otimes \shift \varphi^{41}\otimes \hat{\omega}\varphi^{42}\dots \varphi^{45}\varphi^5}_{37}
+\underbrace{1\otimes \shift \varphi^1\otimes \varphi^2\otimes \shift \varphi^3\otimes \shift \varphi^{41}\otimes \varphi^{42}\hat{\omega}\varphi^{43}\varphi^{44} \varphi^{45}\varphi^5}_{91}\\
&+\underbrace{1\otimes \shift \varphi^1\otimes \varphi^2\otimes \shift \varphi^3\otimes \shift \varphi^{41}\otimes (\shift \otimes \shift)\partial(\shift^{-1}\varphi^{42})\varphi^{43}\varphi^{44} \varphi^{45}\varphi^5}_{229}\\
&
\+\underbrace{\varphi^1\otimes \shift \varphi^{21}\otimes \varphi^{22}\dots \varphi^{24}\hat{\omega}\otimes \shift \varphi^{25}\otimes \shift \varphi^3\otimes \varphi^4\varphi^5}_{38}\\
&+\underbrace{\varphi^1\otimes \shift \varphi^{21}\otimes \varphi^{22}\dots \varphi^{24}\otimes \shift \varphi^{25}\otimes \shift \varphi^3\otimes \hat{\omega}\varphi^4\varphi^5}_{60}
-\underbrace{\varphi^1\otimes \shift \varphi^{21}\otimes \varphi^{22}\dots \varphi^{24}\otimes \shift \varphi^{25}\otimes 1\otimes \shift \varphi^3\otimes \varphi^4\varphi^5}_{346}\\
&-\underbrace{\varphi^1\otimes \shift \varphi^{21}\otimes \varphi^{22}\dots \varphi^{24}\varphi^{251}\varphi^{252}\otimes \shift \varphi^{253}\otimes \shift \varphi^3\otimes \varphi^4\varphi^5}_y\\
&
+\underbrace{\varphi^1\otimes \shift \varphi^{21}\otimes \varphi^{22}\dots \varphi^{24}\otimes \shift \varphi^{25}\otimes \shift \varphi^{31}\otimes \varphi^{32}\varphi^{33}\varphi^4\varphi^5}_g\\
&\+\underbrace{\varphi^1\varphi^2\otimes \shift \varphi^3\otimes \shift \varphi^{41}\otimes \varphi^{42}\varphi^{43}\hat{\omega}\varphi^{44}\otimes \shift \varphi^{45}\otimes \varphi^5}_{116}
+\underbrace{\varphi^1\varphi^2\otimes \shift \varphi^3\otimes \shift \varphi^{41}\otimes \varphi^{42}\varphi^{43}\varphi^{44}\hat{\omega}\otimes \shift \varphi^{45}\otimes \varphi^5}_{61}\\
&+\underbrace{\varphi^1\varphi^2\otimes \shift \varphi^3\otimes \shift \varphi^{41}\otimes \varphi^{42}\varphi^{43}(\shift \otimes \shift)\partial(\shift^{-1}\varphi^{44})\otimes \shift \varphi^{45}\otimes \varphi^5}_{280}\\
&
\-\underbrace{\varphi^1\varphi^2\otimes \shift (\varphi^3\varphi^{41})\otimes \varphi^{42}\varphi^{43}\hat{\omega}\varphi^{44}\otimes \shift \varphi^{45}\otimes \shift \varphi^5\otimes 1}_{138}\\
&-\underbrace{\varphi^1\varphi^2\otimes \shift (\varphi^3\varphi^{41})\otimes \varphi^{42}\dots \varphi^{44}\hat{\omega}\otimes \shift \varphi^{45}\otimes \shift \varphi^5\otimes 1}_{62}\\
&
-\underbrace{\varphi^1\varphi^2\otimes \shift (\varphi^3\varphi^{41})\otimes \varphi^{42}\varphi^{43}(\shift \otimes \shift)\partial(\shift^{-1}\varphi^{44})\otimes \shift \varphi^{45}\otimes \shift \varphi^5\otimes 1}_{235}\\
&\-\underbrace{1\otimes \shift \varphi^1\otimes \shift \varphi^{21}\otimes \varphi^{22}\varphi^{23}\varphi^{24}\hat{\omega}\otimes \shift (\varphi^{25}\varphi^3)\otimes \varphi^4\varphi^5}_{39}
-\underbrace{1\otimes \shift \varphi^1\otimes \shift \varphi^{21}\otimes \varphi^{22}\dots \varphi^{24}\otimes \shift (\varphi^{25}\varphi^3)\otimes \hat{\omega}\varphi^4\varphi^5}_{63}\\
&+\underbrace{1\otimes \shift \varphi^1\otimes \shift \varphi^{21}\otimes \varphi^{22}\dots \varphi^{24}\varphi^{251}\varphi^{252}\otimes \shift (\varphi^{253}\varphi^3)\otimes \varphi^4\varphi^5}_y\\
&+\underbrace{1\otimes \shift \varphi^1\otimes \shift \varphi^{21}\otimes \varphi^{22}\dots \varphi^{25}\varphi^{31}\varphi^{32}\otimes \shift \varphi^{33}\otimes \varphi^4\varphi^5}_g\\
&-\underbrace{1\otimes \shift \varphi^1\otimes \shift \varphi^{21}\otimes \varphi^{22}\dots \varphi^{24}\otimes \shift \varphi^{251}\otimes \varphi^{252}\varphi^{253}\varphi^3\varphi^4\varphi^5}_y\\
&
-\underbrace{1\otimes \shift \varphi^1\otimes \shift \varphi^{21}\otimes \varphi^{22}\dots \varphi^{24}\otimes \shift (\varphi^{25}\varphi^{31})\otimes \varphi^{32}\varphi^{33}\varphi^4\varphi^5}_g\\
&\-\underbrace{\varphi^1\varphi^{21}\varphi^{22}\otimes \shift \varphi^{23}\otimes \varphi^{24}\hat{\omega}\otimes \shift \varphi^{25}\otimes \shift \varphi^3\otimes \varphi^4\varphi^5}_{40}
-\underbrace{\varphi^1\varphi^{21}\varphi^{22}\otimes \shift \varphi^{23}\otimes \varphi^{24}\otimes \shift \varphi^{25}\otimes \shift \varphi^3\otimes \hat{\omega}\varphi^4\varphi^5}_{64}\\
&+\underbrace{\varphi^1\varphi^{21}\varphi^{22}\otimes \shift \varphi^{23}\otimes \varphi^{24}\otimes \shift \varphi^{25}\otimes 1\otimes \shift \varphi^3\otimes \varphi^4\varphi^5}_{332}
+\underbrace{\varphi^1\varphi^{21}\varphi^{22}\otimes \shift \varphi^{23}\otimes \varphi^{24}\varphi^{251}\varphi^{252}\otimes \shift \varphi^{253}\otimes \shift \varphi^3\otimes \varphi^4\varphi^5}_y\\
&-\underbrace{\varphi^1\varphi^{21}\varphi^{22}\otimes \shift \varphi^{23}\otimes \varphi^{24}\otimes \shift \varphi^{25}\otimes \shift \varphi^{31}\otimes \varphi^{32}\varphi^{33}\varphi^4\varphi^5}_g
\-\underbrace{\varphi^1\varphi^2\otimes \shift \varphi^3\otimes \shift \varphi^{41}\otimes \varphi^{42}\hat{\omega}\otimes \shift \varphi^{43}\otimes \varphi^{44}\varphi^{45}\varphi^5}_{151}\\
&-\underbrace{\varphi^1\varphi^2\otimes \shift \varphi^3\otimes \shift \varphi^{41}\otimes \varphi^{42}\otimes \shift \varphi^{43}\otimes \hat{\omega}\varphi^{44}\varphi^{45}\varphi^5}_{65}\\
&
+\underbrace{\varphi^1\varphi^2\otimes \shift \varphi^3\otimes \shift \varphi^{41}\otimes \varphi^{42}\varphi^{431}\varphi^{432}\otimes \shift \varphi^{433}\otimes \varphi^{44}\varphi^{45}\varphi^5}_y\\
&-\underbrace{\varphi^1\varphi^2\otimes \shift \varphi^3\otimes \shift \varphi^{41}\otimes \varphi^{42}\otimes \shift \varphi^{431}\otimes \varphi^{432}\varphi^{433}\varphi^{44}\varphi^{45}\varphi^5}_y\\
&
\+\underbrace{1\otimes \shift \varphi^1\otimes \shift \varphi^{21}\otimes \varphi^{22}\otimes \shift \varphi^{23}\otimes \hat{\omega}\varphi^{24}\varphi^{25}\varphi^3\varphi^4\varphi^5}_{41}\\
&+\underbrace{1\otimes \shift \varphi^1\otimes \shift \varphi^{21}\otimes \varphi^{22}\otimes \shift \varphi^{23}\otimes \varphi^{24}\hat{\omega}\varphi^{25}\varphi^3\varphi^4\varphi^5}_{141}
+\underbrace{1\otimes \shift \varphi^1\otimes \shift \varphi^{21}\otimes \varphi^{22}\otimes \shift \varphi^{23}\otimes (\shift\otimes \shift)\partial(\shift^{-1}\varphi^{24})\varphi^{25}\varphi^3\varphi^4\varphi^5}_{276}\\
&\+\underbrace{\varphi^1\varphi^2\varphi^3\varphi^{41}\varphi^{42}\otimes \shift \varphi^{43}\otimes \hat{\omega}\varphi^{44}\otimes \shift \varphi^{45}\otimes \shift \varphi^5\otimes 1}_{42}
+\underbrace{\varphi^1\varphi^2\varphi^3\varphi^{41}\varphi^{42}\otimes \shift \varphi^{43}\otimes \varphi^{44}\hat{\omega}\otimes \shift \varphi^{45}\otimes \shift \varphi^5\otimes 1}_{66}\\
&+\underbrace{\varphi^1\varphi^2\varphi^3\varphi^{41}\varphi^{42}\otimes \shift \varphi^{43}\otimes (\shift \otimes \shift)\partial(\shift^{-1}\varphi^{44})\otimes \shift \varphi^{45}\otimes \shift\varphi^5\otimes 1}_{285}\\
&
\+\underbrace{\varphi^1\varphi^{21}\dots \varphi^{24}\otimes \shift \varphi^{25}\otimes \shift \varphi^3\otimes \shift \varphi^{41}\otimes \hat{\omega}\varphi^{42}\dots \varphi^{45}\varphi^5}_{43}\\
&+\underbrace{\varphi^1\varphi^{21}\dots \varphi^{24}\otimes \shift \varphi^{25}\otimes \shift \varphi^3\otimes \shift \varphi^{41}\otimes \varphi^{42}\hat{\omega}\varphi^{43}\varphi^{44}\varphi^{45}\varphi^5}_{140}\\
&
+\underbrace{\varphi^1\varphi^{21}\dots \varphi^{24}\otimes \shift \varphi^{25}\otimes \shift \varphi^3\otimes \shift \varphi^{41}\otimes (\shift \otimes \shift)\partial(\shift^{-1}\varphi^{42})\varphi^{43}\dots \varphi^{45}\varphi^5}_{184}\\
&\+\underbrace{\varphi^1\varphi^2\varphi^3\varphi^{41}\dots\varphi^{43}\varphi^{441}\dots \varphi^{443}\hat{\omega}\varphi^{444}\otimes \shift \varphi^{445}\otimes \shift \varphi^{45}\otimes \shift \varphi^5\otimes 1}_{323}\\
&+\underbrace{\varphi^1\varphi^2\varphi^3\varphi^{41}\dots \varphi^{43}\varphi^{441}\dots \varphi^{444}\hat{\omega}\otimes \shift \varphi^{445}\otimes \shift \varphi^{45}\otimes \shift \varphi^5\otimes 1}_{67}\\
&+\underbrace{\varphi^1\varphi^2\varphi^3\varphi^{41}\varphi^{42} \varphi^{43}\varphi^{441}\varphi^{442} \varphi^{443}(\shift \otimes \shift)\partial(\shift^{-1}\varphi^{444})\otimes \shift \varphi^{445}\otimes \shift \varphi^{45}\otimes \shift \varphi^5\otimes 1}_{299}\\
&\+\underbrace{1\otimes \shift \varphi^1\otimes \shift \varphi^{21}\otimes \shift \varphi^{221}\otimes \varphi^{222}\dots \varphi^{225}\varphi^{23}\hat{\omega}\varphi^{24}\varphi^{25}\varphi^3\varphi^4\varphi^5}_{325}\\
&+\underbrace{1\otimes \shift \varphi^1\otimes \shift \varphi^{21}\otimes \shift \varphi^{221}\otimes \varphi^{222}\dots \varphi^{225}\varphi^{23}\varphi^{24}\hat{\omega}\varphi^{25}\varphi^3\varphi^4\varphi^5}_{328}\\
&+\underbrace{1\otimes \shift \varphi^1\otimes \shift \varphi^{21}\otimes \shift \varphi^{221}\otimes \varphi^{222}\dots \varphi^{225}\varphi^{23}(\shift \otimes \shift)\partial(\shift^{-1}\varphi^{24})\varphi^{25}\varphi^3\varphi^4\varphi^5}_g\\
&\-\underbrace{\varphi^1\varphi^{21}\dots \varphi^{23}\varphi^{241}\dots \varphi^{244}\hat{\omega}\otimes \shift \varphi^{245}\otimes \shift \varphi^{25}\otimes \shift \varphi^3\otimes \varphi^4\varphi^5}_{321}\\
&-\underbrace{\varphi^1\varphi^{21}\dots \varphi^{23}\varphi^{241}\dots \varphi^{244}\otimes \shift \varphi^{245}\otimes \shift \varphi^{25}\otimes \shift \varphi^3\otimes \hat{\omega}\varphi^4\varphi^5}_{68}\\
&+\underbrace{\varphi^1\varphi^{21}\dots \varphi^{23}\varphi^{241}\dots \varphi^{244}\otimes \shift \varphi^{245}\otimes 1\otimes \shift \varphi^{25}\otimes \shift \varphi^3\otimes \varphi^4\varphi^5}_{348}\\
&+\underbrace{\varphi^1\varphi^{21}\dots \varphi^{23}\varphi^{241}\dots \varphi^{244}\otimes \shift \varphi^{245}\otimes \shift \varphi^{25}\otimes 1\otimes \shift \varphi^3\otimes \varphi^4\varphi^5}_{352}\\
&+\underbrace{\varphi^1\varphi^{21}\dots \varphi^{23}\varphi^{241}\dots \varphi^{244}\varphi^{2451}\varphi^{2452}\otimes \shift \varphi^{2453}\otimes \shift \varphi^{25}\otimes \shift \varphi^3\otimes \varphi^4\varphi^5}_y\\
&-\underbrace{\varphi^1\varphi^{21}\dots \varphi^{23}\varphi^{241}\dots \varphi^{244}\otimes \shift \varphi^{245}\otimes \shift \varphi^{25}\otimes \shift \varphi^{31}\otimes \varphi^{32}\varphi^{33}\varphi^4\varphi^5}_{gg}\\
&\-\underbrace{\varphi^1\varphi^2\otimes \shift \varphi^3\otimes \shift \varphi^{41}\otimes \shift \varphi^{421}\otimes \varphi^{422}\varphi^{423}\hat{\omega}\varphi^{424}\varphi^{425}\varphi^{43}\varphi^{44}\varphi^{45}\varphi^5}_{327}\\
&-\underbrace{\varphi^1\varphi^2\otimes \shift \varphi^3\otimes \shift \varphi^{41}\otimes \shift \varphi^{421}\otimes \varphi^{422}\varphi^{423} \varphi^{424}\hat{\omega}\varphi^{425}\varphi^{43}\varphi^{44}\varphi^{45}\varphi^5}_{330}\\
&-\underbrace{\varphi^1\varphi^2\otimes \shift \varphi^3\otimes \shift \varphi^{41}\otimes \shift \varphi^{421}\otimes \varphi^{422}\varphi^{423}(\shift \otimes \shift)\partial(\shift^{-1}\varphi^{424})\varphi^{425}\varphi^{43}\dots \varphi^{45}\varphi^5}_y
\end{align*}
and
\begin{align*}
&(1^{\otimes 4}\otimes \mu_2)\mu_5(\varphi)\\
=\,&-\underbrace{1\otimes \shift \varphi^1\otimes \varphi^2\otimes \shift \varphi^3\otimes \varphi^4\hat{\omega}\otimes \shift \varphi^5\otimes 1}_{125}
-\underbrace{1\otimes \shift \varphi^1\otimes \varphi^2\otimes \shift \varphi^3\otimes \varphi^4\otimes \shift \varphi^5\otimes \hat{\omega}}_{44}\\
&+\underbrace{1\otimes \shift \varphi^1\otimes \varphi^2\otimes \shift \varphi^3\otimes \varphi^4\varphi^{51}\varphi^{52}\otimes \shift\varphi^{53}\otimes 1}_{y}
-\underbrace{1\otimes \shift \varphi^1\otimes \varphi^2\otimes \shift \varphi^3\otimes \varphi^4\otimes \shift \varphi^{51}\otimes \varphi^{52}\varphi^{53}}_y\\
&\-\underbrace{1\otimes \shift \varphi^1\otimes \varphi^2\varphi^3\varphi^{41}\dots \varphi^{44}\hat{\omega}\otimes \shift \varphi^{45}\otimes \shift \varphi^5\otimes 1}_{56}
-\underbrace{1\otimes \shift \varphi^1\otimes \varphi^2\varphi^3\varphi^{41}\dots \varphi^{44}\otimes \shift \varphi^{45}\otimes \shift \varphi^5\otimes \hat{\omega}}_{45}\\
&+\underbrace{1\otimes \shift \varphi^1\otimes \varphi^2\varphi^3\varphi^{41}\dots \varphi^{44}\otimes \shift \varphi^{45}\otimes 1\otimes \shift \varphi^5\otimes 1}_{350}\\
&
+\underbrace{1\otimes \shift \varphi^1\otimes \varphi^2\varphi^3\varphi^{41}\dots \varphi^{44}\varphi^{451}\varphi^{452}\otimes \shift \varphi^{453}\otimes \shift \varphi^5\otimes 1}_y\\
&-\underbrace{1\otimes \shift \varphi^1\otimes \varphi^2\varphi^3\varphi^{41}\dots \varphi^{44}\otimes \shift \varphi^{45}\otimes \shift \varphi^{51}\otimes \varphi^{52}\varphi^{53}}_g
\-\underbrace{1\otimes \shift \varphi^1\otimes \shift \varphi^{21}\otimes \varphi^{22}\dots \varphi^{25}\varphi^3\varphi^4\hat{\omega}\otimes \shift \varphi^5\otimes 1}_{57}\\
&-\underbrace{1\otimes \shift \varphi^1\otimes \shift \varphi^{21}\otimes \varphi^{22}\dots \varphi^{25}\varphi^3\varphi^4\otimes \shift \varphi^5\otimes \hat{\omega}}_{46}
+\underbrace{1\otimes \shift \varphi^1\otimes \shift \varphi^{21}\otimes \varphi^{22}\dots \varphi^{25}\varphi^3\varphi^4\varphi^{51}\varphi^{52}\otimes \shift \varphi^{53}\otimes 1}_g\\
&-\underbrace{1\otimes \shift \varphi^1\otimes \shift \varphi^{21}\otimes \varphi^{22}\dots \varphi^{25}\varphi^3\varphi^4\otimes \shift \varphi^{51}\otimes \varphi^{52}\varphi^{53}}_g
\-\underbrace{1\otimes \shift \varphi^1\otimes \varphi^{21}\dots \varphi^{24}\otimes \shift \varphi^{25}\otimes \shift \varphi^3\otimes \hat{\omega}\varphi^4\varphi^5}_{58}\\
&-\underbrace{1\otimes \shift \varphi^1\otimes \varphi^{21}\dots \varphi^{24}\otimes \shift \varphi^{25}\otimes \shift \varphi^3\otimes \varphi^4\hat{\omega}\varphi^5}_{53}\\
&
-\underbrace{1\otimes \shift \varphi^1\otimes \varphi^{21}\dots \varphi^{24}\otimes \shift \varphi^{25}\otimes \shift \varphi^3\otimes (\shift \otimes \shift)\partial(\shift^{-1}\varphi^4)\varphi^5}_{189}\\
&\-\underbrace{\varphi^1\varphi^2\otimes \shift \varphi^3\otimes \shift \varphi^{41}\otimes \varphi^{42}\dots \varphi^{45}\hat{\omega}\otimes \shift \varphi^5\otimes 1}_{120}
-\underbrace{\varphi^1\varphi^2\otimes \shift \varphi^3\otimes \shift \varphi^{41}\otimes \varphi^{42}\dots \varphi^{45}\otimes \shift \varphi^5\otimes \hat{\omega}}_{47}\\
&+\underbrace{\varphi^1\varphi^2\otimes \shift \varphi^3\otimes \shift \varphi^{41}\otimes \varphi^{42}\dots \varphi^{45}\varphi^{51}\varphi^{52}\otimes \shift\varphi^{53}\otimes 1}_{ g}\\
&
-\underbrace{\varphi^1\varphi^2\otimes \shift \varphi^3\otimes \shift \varphi^{41}\otimes \varphi^{42}\dots \varphi^{45}\otimes \shift \varphi^{51}\otimes \varphi^{52}\varphi^{53}}_g\\
&\-\underbrace{\varphi^1\varphi^{21}\dots \varphi^{24}\otimes \shift \varphi^{25}\otimes \shift \varphi^3\otimes \varphi^4\hat{\omega}\otimes \shift \varphi^5\otimes 1}_{59}
-\underbrace{\varphi^1\varphi^{21}\dots \varphi^{24}\otimes \shift \varphi^{25}\otimes \shift \varphi^3\otimes \varphi^4\otimes \shift \varphi^5\otimes \hat{\omega}}_{48}\\
&+\underbrace{\varphi^1\varphi^{21}\dots \varphi^{24}\otimes \shift \varphi^{25}\otimes \shift \varphi^3\otimes \varphi^4\varphi^{51}\varphi^{52}\otimes \shift\varphi^{53}\otimes 1}_{g}\\
&
-\underbrace{\varphi^1\varphi^{21}\dots \varphi^{24}\otimes \shift \varphi^{25}\otimes \shift \varphi^3\otimes \varphi^4\otimes \shift \varphi^{51}\otimes \varphi^{52}\varphi^{53}}_g\\
&\-\underbrace{1\otimes \shift \varphi^1\otimes \varphi^2\otimes \shift \varphi^3\otimes \shift \varphi^{41}\otimes \varphi^{42}\varphi^{43}\hat{\omega}\varphi^{44}\varphi^{45}\varphi^5}_{92}
-\underbrace{1\otimes \shift \varphi^1\otimes \varphi^2\otimes \shift \varphi^3\otimes \shift \varphi^{41}\otimes \varphi^{42}\dots \varphi^{44}\hat{\omega}\varphi^{45}\varphi^5}_{77}\\
&-\underbrace{1\otimes \shift \varphi^1\otimes \varphi^2\otimes \shift \varphi^3\otimes \shift \varphi^{41}\otimes \varphi^{42}\varphi^{43}(\shift \otimes \shift)\partial(\shift^{-1}\varphi^{44})\varphi^{45}\varphi^5}_y\\
&
\-\underbrace{\varphi^1\otimes \shift \varphi^{21}\otimes \varphi^{22}\dots \varphi^{24}\otimes \shift \varphi^{25}\otimes \shift \varphi^3\otimes \hat{\omega}\varphi^4\varphi^5}_{60}\\
&-\underbrace{\varphi^1\otimes \shift \varphi^{21}\otimes \varphi^{22}\dots \varphi^{24}\otimes \shift \varphi^{25}\otimes \shift \varphi^3\otimes \varphi^4\hat{\omega}\varphi^5}_{74}\\
&
-\underbrace{\varphi^1\otimes \shift \varphi^{21}\otimes \varphi^{22}\dots \varphi^{24}\otimes \shift \varphi^{25}\otimes \shift \varphi^3\otimes (\shift \otimes \shift)\partial(\shift^{-1}\varphi^4)\varphi^5}_{190}\\
&\-\underbrace{\varphi^1\varphi^2\otimes \shift \varphi^3\otimes \shift \varphi^{41}\otimes \varphi^{42}\varphi^{43}\varphi^{44}\hat{\omega}\otimes \shift \varphi^{45}\otimes \varphi^5}_{61}
-\underbrace{\varphi^1\varphi^2\otimes \shift \varphi^3\otimes \shift \varphi^{41}\otimes \varphi^{42}\varphi^{43} \varphi^{44}\otimes \shift \varphi^{45}\otimes \hat{\omega}\varphi^5}_{118}\\
&+\underbrace{\varphi^1\varphi^2\otimes \shift \varphi^3\otimes \shift \varphi^{41}\otimes \varphi^{42}\dots \varphi^{44}\varphi^{451}\varphi^{452}\otimes \shift \varphi^{453}\otimes \varphi^5}_y\\
&
-\underbrace{\varphi^1\varphi^2\otimes \shift \varphi^3\otimes \shift \varphi^{41}\otimes \varphi^{42}\dots \varphi^{44}\otimes \shift \varphi^{451}\otimes \varphi^{452}\varphi^{453}\varphi^5}_y\\
&\+\underbrace{\varphi^1\varphi^2\otimes \shift (\varphi^3\varphi^{41})\otimes \varphi^{42}\dots \varphi^{44}\hat{\omega}\otimes \shift \varphi^{45}\otimes \shift \varphi^5\otimes 1}_{62}
+\underbrace{\varphi^1\varphi^2\otimes \shift (\varphi^3\varphi^{41})\otimes \varphi^{42}\dots \varphi^{44}\otimes \shift \varphi^{45}\otimes \shift \varphi^5\otimes \hat{\omega}}_{49}\\
&-\underbrace{\varphi^1\varphi^2\otimes \shift (\varphi^3\varphi^{41})\otimes \varphi^{42}\varphi^{43} \varphi^{44}\otimes \shift \varphi^{45}\otimes 1\otimes \shift \varphi^5\otimes 1}_{354}\\
&
-\underbrace{\varphi^1\varphi^2\otimes \shift (\varphi^3\varphi^{41})\otimes \varphi^{42}\dots \varphi^{44}\varphi^{451}\varphi^{452}\otimes \shift \varphi^{453}\otimes \shift \varphi^5\otimes 1}_y\\
&+\underbrace{\varphi^1\varphi^2\otimes \shift (\varphi^3\varphi^{41})\otimes \varphi^{42}\dots \varphi^{44}\otimes \shift \varphi^{45}\otimes \shift \varphi^{51}\otimes \varphi^{52}\varphi^{53}}_g\\
&
\+\underbrace{1\otimes \shift \varphi^1\otimes \shift \varphi^{21}\otimes \varphi^{22}\dots \varphi^{24}\otimes \shift (\varphi^{25}\varphi^3)\otimes \hat{\omega}\varphi^4\varphi^5}_{63}\\
&+\underbrace{1\otimes \shift \varphi^1\otimes \shift \varphi^{21}\otimes \varphi^{22}\varphi^{23}\varphi^{24}\otimes \shift (\varphi^{25}\varphi^3)\otimes \varphi^4\hat{\omega}\varphi^5}_{101}\\
&
+\underbrace{1\otimes \shift \varphi^1\otimes \shift \varphi^{21}\otimes \varphi^{22}\dots \varphi^{24}\otimes \shift (\varphi^{25}\varphi^3)\otimes (\shift \otimes \shift)\partial(\shift^{-1}\varphi^4)\varphi^5}_{188}\\
&\+\underbrace{\varphi^1\varphi^{21}\varphi^{22}\otimes \shift \varphi^{23}\otimes \varphi^{24}\otimes \shift \varphi^{25}\otimes \shift \varphi^3\otimes \hat{\omega}\varphi^4\varphi^5}_{64}
+\underbrace{\varphi^1\varphi^{21}\varphi^{22}\otimes \shift \varphi^{23}\otimes \varphi^{24}\otimes \shift \varphi^{25}\otimes \shift \varphi^3\otimes \varphi^4\hat{\omega}\varphi^5}_{103}\\
&+\underbrace{\varphi^1\varphi^{21}\varphi^{22}\otimes \shift \varphi^{23}\otimes \varphi^{24}\otimes \shift \varphi^{25}\otimes \shift \varphi^3\otimes (\shift \otimes \shift)\partial(\shift^{-1}\varphi^4)\varphi^5}_{155}\\
&
\+\underbrace{\varphi^1\varphi^2\otimes \shift \varphi^3\otimes \shift \varphi^{41}\otimes \varphi^{42}\otimes \shift \varphi^{43}\otimes \hat{\omega}\varphi^{44}\varphi^{45}\varphi^5}_{65}\\
&+\underbrace{\varphi^1\varphi^2\otimes \shift \varphi^3\otimes \shift \varphi^{41}\otimes \varphi^{42}\otimes \shift \varphi^{43}\otimes \varphi^{44}\hat{\omega}\varphi^{45}\varphi^5}_{111}\\
&
+\underbrace{\varphi^1\varphi^2\otimes \shift \varphi^3\otimes \shift \varphi^{41}\otimes \varphi^{42}\otimes \shift \varphi^{43}\otimes (\shift \otimes \shift)\partial(\shift^{-1}\varphi^{44})\varphi^{45}\varphi^5}_{292}\\
&\-\underbrace{1\otimes \shift \varphi^1\otimes \shift \varphi^{21}\otimes \varphi^{22}\otimes \shift \varphi^{23}\otimes \varphi^{24}\varphi^{25}\varphi^3\hat{\omega}\varphi^4\varphi^5}_{142}
-\underbrace{1\otimes \shift \varphi^1\otimes \shift \varphi^{21}\otimes \varphi^{22}\otimes \shift \varphi^{23}\otimes \varphi^{24}\varphi^{25}\varphi^3\varphi^4\hat{\omega}\varphi^5}_{107}\\
&-\underbrace{1\otimes \shift \varphi^1\otimes \shift \varphi^{21}\otimes \varphi^{22}\otimes \shift \varphi^{23}\otimes \varphi^{24}\varphi^{25}\varphi^3(\shift \otimes \shift)\partial(\shift^{-1}\varphi^4)\varphi^5}_g\\
&
\-\underbrace{\varphi^1\varphi^2\varphi^3\varphi^{41}\varphi^{42}\otimes \shift \varphi^{43}\otimes \varphi^{44}\hat{\omega}\otimes \shift \varphi^{45}\otimes \shift \varphi^5\otimes 1}_{66}\\
&-\underbrace{\varphi^1\varphi^2\varphi^3\varphi^{41}\varphi^{42}\otimes \shift \varphi^{43}\otimes \varphi^{44}\otimes \shift \varphi^{45}\otimes \shift \varphi^5\otimes \hat{\omega}}_{50}
+\underbrace{\varphi^1\varphi^2\varphi^3\varphi^{41}\varphi^{42}\otimes \shift \varphi^{43}\otimes \varphi^{44}\otimes \shift \varphi^{45}\otimes 1\otimes \shift \varphi^5\otimes 1}_{351}\\
&+\underbrace{\varphi^1\varphi^2\varphi^3\varphi^{41}\varphi^{42}\otimes \shift \varphi^{43}\otimes \varphi^{44}\varphi^{451}\varphi^{452}\otimes \shift \varphi^{453}\otimes \shift \varphi^5\otimes 1}_y\\
&
-\underbrace{\varphi^1\varphi^2\varphi^3\varphi^{41}\varphi^{42}\otimes \shift \varphi^{43}\otimes \varphi^{44}\otimes \shift \varphi^{45}\otimes \shift \varphi^{51}\otimes \varphi^{52}\varphi^{53}}_g\\
&\-\underbrace{\varphi^1\varphi^{21}\dots \varphi^{24}\otimes \shift \varphi^{25}\otimes \shift \varphi^3\otimes \shift \varphi^{41}\otimes \varphi^{42}\varphi^{43}\hat{\omega}\varphi^{44}\varphi^{45}\varphi^5}_{139}\\
&
-\underbrace{\varphi^1\varphi^{21}\dots \varphi^{24}\otimes \shift \varphi^{25}\otimes \shift \varphi^3\otimes \shift \varphi^{41}\otimes \varphi^{42}\varphi^{43}\varphi^{44}\hat{\omega}\varphi^{45}\varphi^5}_{122}\\
&-\underbrace{\varphi^1\varphi^{21}\dots \varphi^{24}\otimes \shift \varphi^{25}\otimes \shift \varphi^3\otimes \shift \varphi^{41}\otimes \varphi^{42}\varphi^{43}(\shift \otimes \shift)\partial(\shift^{-1}\varphi^{44})\varphi^{45}}_y\\
&
\-\underbrace{\varphi^1\varphi^2\varphi^3\varphi^{41}\dots \varphi^{43}\varphi^{441}\dots \varphi^{444}\hat{\omega}\otimes \shift \varphi^{445}\otimes \shift \varphi^{45}\otimes \shift \varphi^5\otimes 1}_{67}\\
&-\underbrace{\varphi^1\varphi^2\varphi^3\varphi^{41}\dots \varphi^{43}\varphi^{441}\dots \varphi^{444}\otimes \shift \varphi^{445}\otimes \shift \varphi^{45}\otimes \shift \varphi^5\otimes \hat{\omega}}_{51}\\
&+\underbrace{\varphi^1\varphi^2\varphi^3\varphi^{41}\varphi^{42}\varphi^{43}\varphi^{441}\dots \varphi^{444}\otimes \shift \varphi^{445}\otimes 1\otimes \shift \varphi^{45}\otimes \shift \varphi^5\otimes 1}_{344}\\
&+\underbrace{\varphi^1\varphi^2\varphi^3\varphi^{41}\dots \varphi^{43}\varphi^{441}\dots \varphi^{444}\otimes \shift \varphi^{445}\otimes \shift \varphi^{45}\otimes 1\otimes \shift \varphi^5\otimes 1}_{343}\\
&+\underbrace{\varphi^1\varphi^2\varphi^3\varphi^{41}\dots \varphi^{43}\varphi^{441}\dots \varphi^{444}\varphi^{4451}\varphi^{4452}\otimes \shift \varphi^{4453}\otimes \shift \varphi^{45}\otimes \shift \varphi^5\otimes 1}_y\\
&-\underbrace{\varphi^1\varphi^2\varphi^3\varphi^{41}\dots \varphi^{43}\varphi^{441}\dots \varphi^{444}\otimes \shift \varphi^{445}\otimes \shift \varphi^{45}\otimes \shift \varphi^{51}\otimes \varphi^{52}\varphi^{53}}_{gg}\\
&\-\underbrace{1\otimes \shift \varphi^1\otimes \shift \varphi^{21}\otimes \shift \varphi^{221}\otimes \varphi^{222}\dots \varphi^{225}\varphi^{23}\varphi^{24} \varphi^{25}\varphi^3\hat{\omega}\varphi^4\varphi^5}_{329}\\
&-\underbrace{1\otimes \shift \varphi^1\otimes \shift \varphi^{21}\otimes \shift \varphi^{221}\otimes \varphi^{222}\dots \varphi^{225}\varphi^{23}\varphi^{24} \varphi^{25}\varphi^3\varphi^4\hat{\omega}\varphi^5}_{305}\\
&-\underbrace{1\otimes \shift \varphi^1\otimes \shift \varphi^{21}\otimes \shift \varphi^{221}\otimes \varphi^{222}\dots \varphi^{225}\varphi^{23}\dots \varphi^{25}\varphi^3(\shift \otimes \shift)\partial(\shift^{-1}\varphi^4)\varphi^5}_{gg}\\
&\+\underbrace{\varphi^1\varphi^{21}\dots \varphi^{23}\varphi^{241}\dots \varphi^{244}\otimes \shift \varphi^{245}\otimes \shift \varphi^{25}\otimes \shift \varphi^3\otimes \hat{\omega}\varphi^4\varphi^5}_{68}\\
&
+\underbrace{\varphi^1\varphi^{21}\varphi^{22}\varphi^{23}\varphi^{241}\dots \varphi^{244}\otimes \shift \varphi^{245}\otimes \shift \varphi^{25}\otimes \shift \varphi^3\otimes \varphi^4\hat{\omega}\varphi^5}_{146}\\
&+\underbrace{\varphi^1\varphi^{21}\dots \varphi^{23}\varphi^{241}\dots \varphi^{244}\otimes \shift \varphi^{245}\otimes \shift \varphi^{25}\otimes \shift \varphi^3\otimes (\shift \otimes \shift) \partial(\shift^{-1}\varphi^4)\varphi^5}_{262}\\
&
\+\underbrace{\varphi^1\varphi^2\otimes \shift \varphi^3\otimes \shift \varphi^{41}\otimes \shift \varphi^{421}\otimes \varphi^{422}\dots \varphi^{425}\varphi^{43}\hat{\omega}\varphi^{44}\varphi^{45}\varphi^5}_{331}\\
&+\underbrace{\varphi^1\varphi^2\otimes \shift \varphi^3\otimes \shift \varphi^{41}\otimes \shift \varphi^{421}\otimes \varphi^{422}\dots \varphi^{425}\varphi^{43}\varphi^{44}\hat{\omega}\varphi^{45}\varphi^5}_{309}\\
&
+\underbrace{\varphi^1\varphi^2\otimes \shift \varphi^3\otimes \shift \varphi^{41}\otimes \shift \varphi^{421}\otimes \varphi^{422}\dots \varphi^{425}\varphi^{43}(\shift \otimes \shift)\partial(\shift^{-1}\varphi^{44})\varphi^{45}\varphi^5}_g
\end{align*}
and
\begin{align*}
(\mu_4\otimes 1^{\otimes 2})\mu_3(\varphi) = \,&  (\mu_4\otimes 1^{\otimes 2})(1\otimes \shift\varphi^1\otimes\varphi^2\varphi^3\varphi^4\varphi^5-\varphi^1\varphi^2\otimes \shift\varphi^3\otimes \varphi^4\varphi^5+\varphi^1\varphi^2\varphi^3\varphi^4\otimes \shift\varphi^5\otimes 1)\\
=\,&-\underbrace{1\otimes \shift \varphi^{11}\otimes \shift \varphi^{121}\otimes \shift \varphi^{1221}\otimes \varphi^{1222}\dots \varphi^{1225}\varphi^{123}\dots \varphi^{125}\varphi^{13}\varphi^2\varphi^3\varphi^4\varphi^5}_{gg}\\
&+\underbrace{\varphi^{11}\varphi^{121}\dots \varphi^{123}\varphi^{1241}\dots \varphi^{1244}\otimes \shift \varphi^{1245}\otimes \shift \varphi^{125}\otimes \shift \varphi^{13}\otimes \varphi^2\varphi^3\varphi^4\varphi^5}_{gg}\\
&+\underbrace{\varphi^{11}\varphi^{121}\varphi^{122}\otimes \shift \varphi^{123}\otimes \varphi^{124}\otimes \shift \varphi^{125}\otimes \shift \varphi^{13}\otimes \varphi^2\varphi^3\varphi^4\varphi^5}_g\\
&-\underbrace{1\otimes \shift \varphi^{11}\otimes \shift \varphi^{121}\otimes \varphi^{122}\otimes \shift \varphi^{123}\otimes \varphi^{124}\varphi^{125}\varphi^{13}\varphi^2\varphi^3\varphi^4\varphi^5}_g\\
&+\underbrace{1\otimes \shift \varphi^{11}\otimes \shift \varphi^{121}\otimes \varphi^{122}\dots \varphi^{125}\otimes \shift \varphi^{13}\otimes \varphi^2\varphi^3\varphi^4\varphi^5}_g\\
&-\underbrace{1\otimes \shift (\varphi^{11}\varphi^{121})\otimes \varphi^{122}\dots \varphi^{124}\otimes \shift \varphi^{125}\otimes \shift \varphi^{13}\otimes \varphi^2\varphi^3\varphi^4\varphi^5}_g\\
&+\underbrace{1\otimes \shift \varphi^1\otimes \shift \varphi^{121}\otimes \varphi^{122}\dots \varphi^{124}\otimes \shift \varphi^{125}\otimes \varphi^{13}\varphi^2\varphi^3\varphi^4\varphi^5}_g\\
&\-\underbrace{\varphi^1\varphi^{21}\varphi^{22}\otimes \shift \varphi^{23}\otimes \varphi^{24}\otimes \shift \varphi^{25}\otimes 1\otimes \shift \varphi^3\otimes \varphi^4\varphi^5}_{332}\\
&+\underbrace{\varphi^1\otimes \shift \varphi^{21}\otimes \varphi^{22}\dots \varphi^{24}\otimes \shift \varphi^{25}\otimes 1\otimes \shift \varphi^3\otimes \varphi^4\varphi^5}_{346}\\
&-\underbrace{\varphi^1\otimes \shift \varphi^{21}\otimes \varphi^{22}\otimes \shift \varphi^{23}\otimes \varphi^{24}\varphi^{25}\otimes \shift \varphi^3\otimes \varphi^4\varphi^5}_{202}\\
&-\underbrace{\varphi^1\varphi^{21}\varphi^{221}\dots \varphi^{224}\otimes \shift \varphi^{225}\otimes \shift \varphi^{23}\otimes \varphi^{24}\varphi^{25}\otimes \shift \varphi^3\otimes \varphi^4\varphi^5}_{172}\\
&-\underbrace{\varphi^1\otimes \shift \varphi^{21}\otimes \shift \varphi^{221}\otimes \varphi^{222}\dots \varphi^{225}\varphi^{23}\dots \varphi^{25}\otimes \shift \varphi^3\otimes \varphi^4\varphi^5}_{239}\\
&-\underbrace{\varphi^1\varphi^{21}\varphi^{22}\otimes \shift \varphi^{23}\otimes \shift \varphi^{241}\otimes \varphi^{242}\dots \varphi^{245}\varphi^{25}\otimes \shift \varphi^3\otimes \varphi^4\varphi^5}_{269}\\
&-\underbrace{\varphi^1\varphi^{21}\dots \varphi^{23}\varphi^{241}\dots \varphi^{244}\otimes \shift \varphi^{245}\otimes \shift \varphi^{25}\otimes 1\otimes \shift \varphi^3\otimes \varphi^4\varphi^5}_{352}\\
&\+\underbrace{\varphi^1\varphi^{21}\varphi^{22}\otimes \shift \varphi^{23}\otimes \varphi^{24}\otimes \shift \varphi^{25}\otimes \varphi^3\varphi^4\otimes \shift \varphi^5\otimes 1}_{214}\\
&-\underbrace{\varphi^1\otimes \shift \varphi^{21}\otimes \varphi^{22}\dots \varphi^{24}\otimes \shift \varphi^{25}\otimes \varphi^3\varphi^4\otimes \shift \varphi^5\otimes 1}_{203}\\
&+\underbrace{\varphi^1\otimes \shift \varphi^{21}\otimes \varphi^{22}\otimes \shift \varphi^{23}\otimes \varphi^{24}\varphi^{25}\varphi^3\varphi^4\otimes \shift \varphi^5\otimes 1}_{196}\\
&+\underbrace{\varphi^1\varphi^{21}\varphi^{221}\dots \varphi^{224}\otimes \shift \varphi^{225}\otimes \shift \varphi^{23}\otimes \varphi^{24}\varphi^{25}\varphi^3\varphi^4\otimes \shift \varphi^5\otimes 1}_g\\
&+\underbrace{\varphi^1\otimes \shift \varphi^{21}\otimes \shift \varphi^{221}\otimes \varphi^{222}\dots \varphi^{225}\varphi^{23}\varphi^{24}\varphi^{25}\varphi^3\varphi^4\otimes \shift \varphi^5\otimes 1}_{226}\\
&+\underbrace{\varphi^1\varphi^{21}\varphi^{22}\otimes \shift \varphi^{23}\otimes \shift \varphi^{241}\otimes \varphi^{242}\dots \varphi^{245}\varphi^{25}\varphi^3\varphi^4\otimes \shift \varphi^5\otimes 1}_g\\
&+\underbrace{\varphi^1\varphi^{21}\varphi^{22}\varphi^{23}\varphi^{241}\dots \varphi^{244}\otimes \shift \varphi^{245}\otimes \shift \varphi^{25}\otimes \varphi^3\varphi^4\otimes \shift \varphi^5\otimes 1}_{270}
\end{align*}
and
\begin{align*}
-(1\otimes \mu_4\otimes 1)\mu_3(\varphi) 
=\,& -(1\otimes \mu_4\otimes 1)
(1\otimes \shift\varphi^1\otimes\varphi^2\varphi^3\varphi^4\varphi^5-\varphi^1\varphi^2\otimes \shift\varphi^3\otimes \varphi^4\varphi^5+\varphi^1\varphi^2\varphi^3\varphi^4\otimes \shift\varphi^5\otimes 1)
\\
=\,&
-\underbrace{1\otimes \shift \varphi^1\otimes \varphi^{21}\varphi^{22}\otimes \shift \varphi^{23}\otimes \varphi^{24}\otimes \shift \varphi^{25}\otimes \varphi^3\varphi^4\varphi^5}_{222}\\
&+\underbrace{1\otimes \shift \varphi^1\otimes 1\otimes \shift \varphi^{21}\otimes \varphi^{22}\varphi^{23}\varphi^{24}\otimes \shift \varphi^{25}\otimes \varphi^3\varphi^4\varphi^5}_{356}\\
&-\underbrace{1\otimes \shift \varphi^1\otimes 1\otimes \shift \varphi^{21}\otimes \varphi^{22}\otimes \shift \varphi^{23}\otimes \varphi^{24}\varphi^{25}\varphi^3\varphi^4\varphi^5}_{69}\\
&-\underbrace{1\otimes \shift \varphi^1\otimes \varphi^{21}\varphi^{221}\dots \varphi^{224}\otimes \shift \varphi^{225}\otimes \shift \varphi^{23}\otimes \varphi^{24}\varphi^{25}\varphi^3\varphi^4\varphi^5}_{250}\\
&-\underbrace{1\otimes \shift \varphi^1\otimes 1\otimes \shift \varphi^{21}\otimes \shift \varphi^{221}\otimes \varphi^{222}\dots \varphi^{225}\varphi^{23}\varphi^{24}\varphi^{25}\varphi^3\varphi^4\varphi^5}_{70}\\
&-\underbrace{1\otimes \shift \varphi^1\otimes \varphi^{21}\varphi^{22}\otimes \shift \varphi^{23}\otimes \shift \varphi^{241}\otimes \varphi^{242}\dots \varphi^{245}\varphi^{25}\varphi^3\varphi^4\varphi^5}_{279}\\
&-\underbrace{1\otimes \shift \varphi^1\otimes \varphi^{21}\varphi^{22}\varphi^{23}\varphi^{241}\dots \varphi^{244}\otimes \shift \varphi^{245}\otimes \shift \varphi^{25}\otimes \varphi^3\varphi^4\varphi^5}_{259}\\
&\-\underbrace{\varphi^1\varphi^2\otimes \shift \varphi^{31}\otimes \shift \varphi^{321}\otimes \shift \varphi^{3221}\otimes \varphi^{3222}\dots \varphi^{3225}\varphi^{323}\dots \varphi^{325}\varphi^{33}\varphi^4\varphi^5}_g\\
&+\underbrace{\varphi^1\varphi^2\varphi^{31}\varphi^{321}\dots \varphi^{323}\varphi^{3241}\dots \varphi^{3244}\otimes \shift \varphi^{3245}\otimes \shift \varphi^{325}\otimes \shift \varphi^{33}\otimes \varphi^4\varphi^5}_g\\
&+\underbrace{\varphi^1\varphi^2\varphi^{31} \varphi^{321}\varphi^{322}\otimes \shift \varphi^{323}\otimes \varphi^{324}\otimes \shift \varphi^{325}\otimes \shift \varphi^{33}\otimes \varphi^4\varphi^5}_{g}\\
&-\underbrace{\varphi^1\varphi^2\otimes \shift \varphi^{31}\otimes \shift \varphi^{321}\otimes \varphi^{322}\otimes \shift \varphi^{323}\otimes \varphi^{324}\varphi^{325}\varphi^{33}\varphi^4\varphi^5}_g\\
&+\underbrace{\varphi^1\varphi^2\otimes \shift \varphi^{31}\otimes \shift \varphi^{321}\otimes \varphi^{322}\dots \varphi^{325}\otimes \shift \varphi^{33}\otimes \varphi^4\varphi^5}_g\\
&-\underbrace{\varphi^1\varphi^2\otimes \shift (\varphi^{31}\varphi^{321})\otimes \varphi^{322}\dots \varphi^{324}\otimes \shift \varphi^{325}\otimes \shift \varphi^{33}\otimes \varphi^4\varphi^5}_g\\
&+\underbrace{\varphi^1\varphi^2\otimes \shift \varphi^{31}\otimes \shift \varphi^{321}\otimes \varphi^{322}\dots \varphi^{324}\otimes \shift \varphi^{325}\otimes \varphi^{33}\varphi^4\varphi^5}_g\\
&\-\underbrace{\varphi^1\varphi^2\varphi^3\varphi^{41}\varphi^{42}\otimes \shift \varphi^{43}\otimes \varphi^{44}\otimes \shift \varphi^{45}\otimes 1\otimes \shift \varphi^5\otimes 1}_{351}\\
&+\underbrace{\varphi^1\varphi^2\varphi^3\otimes \shift \varphi^{41}\otimes \varphi^{42}\varphi^{43}\varphi^{44}\otimes \shift \varphi^{45}\otimes 1\otimes \shift \varphi^5\otimes 1}_{150}\\
&-\underbrace{\varphi^1\varphi^2\varphi^3\otimes \shift \varphi^{41}\otimes \varphi^{42}\otimes \shift \varphi^{43}\otimes \varphi^{44}\varphi^{45}\otimes \shift \varphi^5\otimes 1}_{215}\\
&-\underbrace{\varphi^1\varphi^2\varphi^3\varphi^{41}\varphi^{421}\dots \varphi^{424}\otimes \shift \varphi^{425}\otimes \shift \varphi^{43}\otimes \varphi^{44}\varphi^{45}\otimes \shift \varphi^5\otimes 1}_{264}\\
&-\underbrace{\varphi^1\varphi^2\varphi^3\otimes \shift \varphi^{41}\otimes \shift \varphi^{421}\otimes \varphi^{422}\dots \varphi^{425}\varphi^{43}\dots \varphi^{45}\otimes \shift \varphi^5\otimes 1}_{246}\\
&-\underbrace{\varphi^1\varphi^2\varphi^3\varphi^{41}\varphi^{42}\otimes \shift \varphi^{43}\otimes \shift \varphi^{441}\otimes \varphi^{442}\dots \varphi^{445}\varphi^{45}\otimes \shift \varphi^5\otimes 1}_{286}\\
&-\underbrace{\varphi^1\varphi^2\varphi^3\varphi^{41}\dots \varphi^{43}\varphi^{441}\dots \varphi^{444}\otimes \shift \varphi^{445}\otimes \shift \varphi^{45}\otimes 1\otimes \shift \varphi^5\otimes 1}_{343}
\end{align*}
and
\begin{align*}
(1^{\otimes 2}\otimes \mu_4)\mu_3(\varphi)=\,&(1^{\otimes 2}\otimes \mu_4)
(1\otimes \shift\varphi^1\otimes\varphi^2\varphi^3\varphi^4\varphi^5-\varphi^1\varphi^2\otimes \shift\varphi^3\otimes \varphi^4\varphi^5+\varphi^1\varphi^2\varphi^3\varphi^4\otimes \shift\varphi^5\otimes 1)\\
=\,&
+\underbrace{1\otimes \shift \varphi^1\otimes \varphi^2\varphi^3\varphi^{41}\varphi^{42}\otimes \shift \varphi^{43}\otimes \varphi^{44}\otimes \shift \varphi^{45}\otimes \varphi^5}_{211}\\
&-\underbrace{1\otimes \shift \varphi^1\otimes \varphi^2\varphi^3\otimes \shift \varphi^{41}\otimes \varphi^{42}\dots \varphi^{44}\otimes \shift \varphi^{45}\otimes \varphi^5}_{206}\\
&+\underbrace{1\otimes \shift \varphi^1\otimes \varphi^2\varphi^3\otimes \shift \varphi^{41}\otimes \varphi^{42}\otimes \shift \varphi^{43}\otimes \varphi^{44}\varphi^{45}\varphi^5}_{198}\\
&+\underbrace{1\otimes \shift \varphi^1\otimes \varphi^2\varphi^3\varphi^{41}\varphi^{421}\dots \varphi^{424}\otimes \shift \varphi^{425}\otimes \shift \varphi^{43}\otimes \varphi^{44}\varphi^{45}\varphi^5}_g\\
&+\underbrace{1\otimes \shift \varphi^1\otimes \varphi^2\varphi^3\otimes \shift \varphi^{41}\otimes \shift \varphi^{421}\otimes \varphi^{422}\dots \varphi^{425}\varphi^{43}\varphi^{44}\varphi^{45}\varphi^5}_{157}\\
&+\underbrace{1\otimes \shift \varphi^1\otimes \varphi^2\varphi^3\varphi^{41}\varphi^{42}\otimes \shift \varphi^{43}\otimes \shift \varphi^{441}\otimes \varphi^{442}\dots \varphi^{445}\varphi^{45}\varphi^5}_g\\
&+\underbrace{1\otimes \shift \varphi^1\otimes \varphi^2\varphi^3\varphi^{41}\varphi^{42} \varphi^{43}\varphi^{441}\dots \varphi^{444}\otimes \shift \varphi^{445}\otimes \shift \varphi^{45}\otimes \varphi^5}_{228}\\
&\-\underbrace{\varphi^1\varphi^2\otimes \shift \varphi^3\otimes \varphi^{41}\varphi^{42}\otimes \shift \varphi^{43}\otimes \varphi^{44}\otimes \shift \varphi^{45}\otimes \varphi^5}_{221}\\
&+\underbrace{\varphi^1\varphi^2\otimes \shift \varphi^3\otimes 1\otimes \shift \varphi^{41}\otimes \varphi^{42}\dots \varphi^{44}\otimes \shift \varphi^{45}\otimes \varphi^5}_{349}\\
&-\underbrace{\varphi^1\varphi^2\otimes \shift \varphi^3\otimes 1\otimes \shift \varphi^{41}\otimes \varphi^{42}\otimes \shift \varphi^{43}\otimes \varphi^{44}\varphi^{45}\varphi^5}_{336}\\
&-\underbrace{\varphi^1\varphi^2\otimes \shift \varphi^3\otimes \varphi^{41}\varphi^{421}\dots \varphi^{424}\otimes \shift \varphi^{425}\otimes \shift \varphi^{43}\otimes \varphi^{44}\varphi^{45}\varphi^5}_{275}\\
&-\underbrace{\varphi^1\varphi^2\otimes \shift \varphi^3\otimes 1\otimes \shift \varphi^{41}\otimes \shift \varphi^{421}\otimes \varphi^{422}\dots \varphi^{425}\varphi^{43}\dots \varphi^{45}\varphi^5}_{337}\\
&-\underbrace{\varphi^1\varphi^2\otimes \shift \varphi^3\otimes \varphi^{41}\varphi^{42}\otimes \shift \varphi^{43}\otimes \shift \varphi^{441}\otimes \varphi^{442}\dots \varphi^{445}\varphi^{45}\varphi^5}_{295}\\
&-\underbrace{\varphi^1\varphi^2\otimes \shift \varphi^3\otimes \varphi^{41}\dots \varphi^{43}\varphi^{441}\dots \varphi^{444}\otimes \shift \varphi^{445}\otimes \shift \varphi^{45}\otimes \varphi^5}_{283}\\
&\-\underbrace{\varphi^1\dots \varphi^4\otimes \shift \varphi^{51}\otimes \shift \varphi^{521}\otimes \shift \varphi^{5221}\otimes \varphi^{5222}\dots \varphi^{5225}\varphi^{523}\dots \varphi^{525}\varphi^{53}}_{gg}\\
&+\underbrace{\varphi^1\dots \varphi^4\varphi^{51}\varphi^{521}\dots \varphi^{523}\varphi^{5241}\dots \varphi^{5244}\otimes \shift \varphi^{5245}\otimes \shift \varphi^{525}\otimes \shift \varphi^{53}\otimes 1}_{gg}\\
&+\underbrace{\varphi^1\dots \varphi^4\varphi^{51}\varphi^{521}\varphi^{522}\otimes \shift \varphi^{523}\otimes \varphi^{524}\otimes \shift \varphi^{525}\otimes \shift \varphi^{53}\otimes 1}_g\\
&-\underbrace{\varphi^1\dots \varphi^4\otimes \shift \varphi^{51}\otimes \shift \varphi^{521}\otimes \varphi^{522}\otimes \shift \varphi^{523}\otimes \varphi^{524}\varphi^{525}\varphi^{53}}_g\\
&+\underbrace{\varphi^1\dots \varphi^4\otimes \shift \varphi^{51}\otimes \shift \varphi^{521}\otimes \varphi^{522}\dots \varphi^{525}\otimes \shift \varphi^{53}\otimes 1}_g\\
&-\underbrace{\varphi^1\dots \varphi^4\otimes \shift (\varphi^{51}\varphi^{521})\otimes \varphi^{522}\dots \varphi^{524}\otimes \shift \varphi^{525}\otimes \shift \varphi^{53}\otimes 1}_g\\
&+\underbrace{\varphi^1\dots \varphi^4\otimes \shift \varphi^{51}\otimes \shift \varphi^{521}\otimes \varphi^{522}\dots \varphi^{524}\otimes \shift \varphi^{525}\otimes \varphi^{53}}_g
\end{align*}
and
\begin{align*}
&-(\mu_3\otimes 1^{\otimes 3})\mu_4(\varphi)\\
=\,&-\underbrace{\varphi^1\otimes \shift \varphi^{21}\otimes \varphi^{22}\dots \varphi^{25}\otimes \shift \varphi^3\otimes \varphi^4\otimes \shift \varphi^5\otimes 1}_{204}\\
&+\underbrace{\varphi^1\varphi^{21}\varphi^{22}\otimes \shift \varphi^{23}\otimes \varphi^{24}\varphi^{25}\otimes \shift \varphi^3\otimes \varphi^4\otimes \shift \varphi^5\otimes 1}_{213}\\
&-\underbrace{\varphi^1\varphi^{21}\dots \varphi^{24}\otimes \shift \varphi^{25}\otimes 1\otimes \shift \varphi^3\otimes \varphi^4\otimes \shift \varphi^5\otimes 1}_{342}\\
&\+\underbrace{1\otimes \shift \varphi^{11}\otimes \varphi^{12}\otimes \shift \varphi^{13}\otimes \varphi^2\varphi^3\varphi^4\otimes \shift \varphi^5\otimes 1}_g\\
&+\underbrace{\varphi^{11}\varphi^{121}\dots \varphi^{124}\otimes \shift \varphi^{125}\otimes \shift \varphi^{13}\otimes \varphi^2\varphi^3\varphi^4\otimes \shift \varphi^5\otimes 1}_g\\
&+\underbrace{1\otimes \shift \varphi^{11}\otimes \shift \varphi^{121}\otimes \varphi^{122}\dots \varphi^{125}\varphi^{13}\varphi^2\varphi^3\varphi^4\otimes \shift \varphi^5\otimes 1}_{g}\\
&\-\underbrace{1\otimes \shift \varphi^{11}\otimes \varphi^{12}\otimes \shift \varphi^{13}\otimes \varphi^2\otimes \shift \varphi^3\otimes \varphi^4\varphi^5}_y\\
&-\underbrace{\varphi^{11}\varphi^{121}\dots \varphi^{124}\otimes \shift \varphi^{125}\otimes \shift \varphi^{13}\otimes \varphi^2\otimes \shift \varphi^3\otimes \varphi^4\varphi^5}_g\\
&-\underbrace{1\otimes \shift \varphi^{11}\otimes \shift \varphi^{121}\otimes \varphi^{122}\dots \varphi^{125}\varphi^{13}\varphi^2\otimes \shift \varphi^3\otimes \varphi^4\varphi^5}_g\\
&\-\underbrace{\varphi^1\varphi^{21}\otimes \shift \varphi^{221}\otimes \varphi^{222}\dots \varphi^{225}\varphi^{23}\varphi^{24}\otimes \shift \varphi^{25}\otimes \shift \varphi^3\otimes \varphi^4\varphi^5}_{240}\\
&+\underbrace{\varphi^1\varphi^{21}\varphi^{221}\varphi^{222}\otimes \shift \varphi^{223}\otimes \varphi^{224}\varphi^{225}\varphi^{23}\varphi^{24}\otimes \shift \varphi^{25}\otimes \shift \varphi^3\otimes \varphi^4\varphi^5}_y\\
&-\underbrace{\varphi^1\varphi^{21}\varphi^{221}\dots \varphi^{224}\otimes \shift \varphi^{225}\otimes \varphi^{23}\varphi^{24}\otimes \shift \varphi^{25}\otimes \shift \varphi^3\otimes \varphi^4\varphi^5}_{171}\\
&\-\underbrace{1\otimes \shift \varphi^{11}\otimes \varphi^{12}\otimes \shift \varphi^{13}\otimes \shift \varphi^{21}\otimes \varphi^{22}\dots \varphi^{25}\varphi^3\varphi^4\varphi^5}_g\\
&+\underbrace{1\otimes \shift \varphi^1\otimes \shift \varphi^{211}\otimes \varphi^{212}\otimes \shift \varphi^{213}\otimes \varphi^{22}\dots \varphi^{25}\varphi^3\varphi^4\varphi^5}_y\\
&-\underbrace{\varphi^{11}\varphi^{121}\dots \varphi^{124}\otimes \shift \varphi^{125}\otimes \shift \varphi^{13}\otimes \shift \varphi^{21}\otimes \varphi^{22}\dots \varphi^{25}\varphi^3\varphi^4\varphi^5}_{gg}\\
&+\underbrace{1\otimes \shift \varphi^1\otimes \shift \varphi^{211}\otimes \shift \varphi^{2121}\otimes \varphi^{2122}\dots \varphi^{2125}\varphi^{213}\varphi^{22}\dots \varphi^{25}\varphi^3\varphi^4\varphi^5}_g\\
&\-\underbrace{\varphi^1\otimes \shift \varphi^{21}\otimes \varphi^{22}\dots \varphi^{25}\otimes \shift \varphi^3\otimes \shift \varphi^{41}\otimes \varphi^{42}\dots \varphi^{45}\varphi^5}_{187}\\
&
+\underbrace{\varphi^1\varphi^{21}\varphi^{22}\otimes \shift \varphi^{23}\otimes \varphi^{24}\varphi^{25}\otimes \shift \varphi^3\otimes \shift \varphi^{41}\otimes \varphi^{42}\dots \varphi^{45}\varphi^5}_{216}\\
&-\underbrace{\varphi^1\varphi^{21}\dots \varphi^{24}\otimes \shift \varphi^{25}\otimes 1\otimes \shift \varphi^3\otimes \shift \varphi^{41}\otimes \varphi^{42}\dots \varphi^{45}\varphi^5}_{341}\\
&
\-\underbrace{\varphi^1\otimes \shift \varphi^{21}\otimes \varphi^{22}\dots \varphi^{25}\varphi^3\varphi^{41}\dots \varphi^{44}\otimes \shift \varphi^{45}\otimes \shift \varphi^5\otimes 1}_{180}\\
&+\underbrace{\varphi^1\varphi^{21}\varphi^{22}\otimes \shift \varphi^{23}\otimes \varphi^{24}\varphi^{25}\varphi^3\varphi^{41}\dots \varphi^{44}\otimes \shift \varphi^{45}\otimes \shift \varphi^5\otimes 1}_g\\
&
-\underbrace{\varphi^1\varphi^{21}\dots \varphi^{24}\otimes \shift \varphi^{25}\otimes \varphi^3\varphi^{41}\dots \varphi^{44}\otimes \shift \varphi^{45}\otimes \shift \varphi^5\otimes 1}_{165}
\end{align*}
and
\begin{align*}
&-(1\otimes \mu_3\otimes 1^{\otimes 2})\mu_4(\varphi)\\
=\,&-\underbrace{\varphi^1\varphi^2\otimes \shift \varphi^{31}\otimes \varphi^{32}\otimes \shift \varphi^{33}\otimes \varphi^4\otimes \shift \varphi^5\otimes 1}_{y}\\
&-\underbrace{\varphi^1\varphi^2\varphi^{31}\varphi^{321}\dots \varphi^{324}\otimes \shift \varphi^{325}\otimes \shift \varphi^{33}\otimes \varphi^4\otimes \shift \varphi^5\otimes 1}_g\\
&
-\underbrace{\varphi^1\varphi^2\otimes \shift \varphi^{31}\otimes \shift \varphi^{321}\otimes \varphi^{322}\dots \varphi^{325}\varphi^{33}\varphi^4\otimes \shift \varphi^5\otimes 1}_g\\
&\-\underbrace{1\otimes \shift \varphi^1\otimes 1\otimes \shift \varphi^{21}\otimes \varphi^{22}\dots \varphi^{25}\varphi^3\varphi^4\otimes \shift \varphi^5\otimes 1}_{333}
+\underbrace{1\otimes \shift \varphi^1\otimes \varphi^{21}\varphi^{22}\otimes \shift \varphi^{23}\otimes \varphi^{24}\varphi^{25}\varphi^3\varphi^4\otimes \shift \varphi^5\otimes 1}_{197}\\
&-\underbrace{1\otimes \shift \varphi^1\otimes \varphi^{21}\dots \varphi^{24}\otimes \shift \varphi^{25}\otimes \varphi^3\varphi^4\otimes \shift \varphi^5\otimes 1}_{230}\\
&\+\underbrace{1\otimes \shift \varphi^1\otimes 1\otimes \shift \varphi^{21}\otimes \varphi^{22}\dots \varphi^{25}\otimes \shift \varphi^3\otimes \varphi^4\varphi^5}_{355}
-\underbrace{1\otimes \shift \varphi^1\otimes \varphi^{21}\varphi^{22}\otimes \shift \varphi^{23}\otimes \varphi^{24}\varphi^{25}\otimes \shift \varphi^3\otimes \varphi^4\varphi^5}_{145}\\
&+\underbrace{1\otimes \shift \varphi^1\otimes \varphi^{21}\dots \varphi^{24}\otimes \shift \varphi^{25}\otimes 1\otimes \shift \varphi^3\otimes \varphi^4\varphi^5}_{347}\\
&
\-\underbrace{\varphi^1\varphi^{21}\dots \varphi^{23}\otimes \shift \varphi^{241}\otimes \varphi^{242}\dots \varphi^{245}\otimes \shift \varphi^{25}\otimes \shift \varphi^3\otimes \varphi^4\varphi^5}_{272}\\
&+\underbrace{\varphi^1\varphi^{21}\dots \varphi^{23}\varphi^{241}\varphi^{242}\otimes \shift \varphi^{243}\otimes \varphi^{244}\varphi^{245}\otimes \shift \varphi^{25}\otimes \shift \varphi^3\otimes \varphi^4\varphi^5}_{175}\\
&-\underbrace{\varphi^1\varphi^{21}\dots \varphi^{23}\varphi^{241}\dots \varphi^{244}\otimes \shift \varphi^{245}\otimes 1\otimes \shift \varphi^{25}\otimes \shift \varphi^3\otimes \varphi^4\varphi^5}_{348}\\
&\-\underbrace{1\otimes \shift \varphi^1\otimes \shift \varphi^{21}\otimes 1\otimes \shift \varphi^{221}\otimes \varphi^{222}\dots \varphi^{225}\varphi^{23}\dots \varphi^{25}\varphi^3\varphi^4\varphi^5}_{334}\\
&+\underbrace{1\otimes \shift \varphi^1\otimes \shift \varphi^{21}\otimes \varphi^{221}\varphi^{222}\otimes \shift \varphi^{223}\otimes \varphi^{224}\varphi^{225}\varphi^{23}\dots \varphi^{25}\varphi^3\varphi^4\varphi^5}_{219}\\
&-\underbrace{1\otimes \shift \varphi^1\otimes \shift \varphi^{21}\otimes \varphi^{221}\dots \varphi^{224}\otimes \shift \varphi^{225}\otimes \varphi^{23}\dots \varphi^{25}\varphi^3\varphi^4\varphi^5}_{249}\\
&
\-\underbrace{\varphi^1\varphi^2\otimes \shift \varphi^{31}\otimes \varphi^{32}\otimes \shift \varphi^{33}\otimes \shift \varphi^{41}\otimes \varphi^{42}\dots \varphi^{45}\varphi^5}_g\\
&+\underbrace{\varphi^1\varphi^2\otimes \shift \varphi^3\otimes \shift \varphi^{411}\otimes \varphi^{412}\otimes \shift \varphi^{413}\otimes \varphi^{42}\dots \varphi^{45}\varphi^5}_y\\
&
-\underbrace{\varphi^1\varphi^2\varphi^{31}\varphi^{321}\dots \varphi^{324}\otimes \shift \varphi^{325}\otimes \shift \varphi^{33}\otimes \shift \varphi^{41}\otimes \varphi^{42}\dots \varphi^{45}\varphi^5}_{gg}\\
&+\underbrace{\varphi^1\varphi^2\otimes \shift \varphi^3\otimes \shift \varphi^{411}\otimes \shift \varphi^{4121}\otimes \varphi^{4122}\dots \varphi^{4125}\varphi^{413}\varphi^{42}\dots \varphi^{45}\varphi^5}_g\\
&\-\underbrace{\varphi^1\varphi^2\varphi^3\varphi^{41}\otimes \shift \varphi^{421}\otimes \varphi^{422}\dots \varphi^{425}\varphi^{43}\varphi^{44}\otimes \shift \varphi^{45}\otimes \shift \varphi^5\otimes 1}_{245}\\
&+\underbrace{\varphi^1\varphi^2\varphi^3\varphi^{41}\varphi^{421}\varphi^{422}\otimes \shift \varphi^{423}\otimes \varphi^{424}\varphi^{425}\varphi^{43}\varphi^{44}\otimes \shift \varphi^{45}\otimes \shift \varphi^5\otimes 1}_y\\
&-\underbrace{\varphi^1\varphi^2\varphi^3\varphi^{41}\varphi^{421}\dots \varphi^{424}\otimes \shift \varphi^{425}\otimes \varphi^{43}\varphi^{44}\otimes \shift \varphi^{45}\otimes \shift\varphi^5\otimes 1}_{265}
\end{align*}
and
\begin{align*}
&-(1^{\otimes 2}\otimes \mu_3\otimes 1)\mu_4(\varphi)\\
=\,&+\underbrace{\varphi^1\varphi^2\otimes \shift \varphi^3\otimes 1\otimes \shift \varphi^{41}\otimes \varphi^{42}\dots \varphi^{45}\otimes \shift \varphi^5\otimes 1}_{335}\\
&-\underbrace{\varphi^1\varphi^2\otimes \shift \varphi^3\otimes \varphi^{41}\varphi^{42}\otimes \shift \varphi^{43}\otimes \varphi^{44}\varphi^{45}\otimes \shift \varphi^5\otimes 1}_{223}
+\underbrace{\varphi^1\varphi^2\otimes \shift \varphi^3\otimes \varphi^{41}\dots \varphi^{44}\otimes \shift \varphi^{45}\otimes 1\otimes \shift \varphi^5\otimes 1}_{357}\\
&\-\underbrace{1\otimes \shift \varphi^1\otimes \varphi^2\varphi^3\otimes \shift \varphi^{41}\otimes \varphi^{42}\dots \varphi^{45}\otimes \shift \varphi^5\otimes 1}_{210}
+\underbrace{1\otimes \shift \varphi^1\otimes \varphi^2\varphi^3\varphi^{41}\varphi^{42}\otimes \shift \varphi^{43}\otimes \varphi^{44}\varphi^{45}\otimes \shift \varphi^5\otimes 1}_{212}\\
&-\underbrace{1\otimes \shift \varphi^1\otimes \varphi^2\varphi^3\varphi^{41}\dots \varphi^{44}\otimes \shift \varphi^{45}\otimes 1\otimes \shift \varphi^5\otimes 1}_{350}\\
&\+\underbrace{1\otimes \shift \varphi^1\otimes \varphi^2\otimes \shift \varphi^{31}\otimes \varphi^{32}\otimes \shift \varphi^{33}\otimes \varphi^4\varphi^5}_y
+\underbrace{1\otimes \shift \varphi^1\otimes \varphi^2\varphi^{31}\varphi^{321}\dots \varphi^{324}\otimes \shift \varphi^{325}\otimes \shift \varphi^{33}\otimes \varphi^4\varphi^5}_g\\
&+\underbrace{1\otimes \shift \varphi^1\otimes \varphi^2\otimes \shift \varphi^{31}\otimes \shift \varphi^{321}\otimes \varphi^{322}\dots \varphi^{325}\varphi^{33}\varphi^4\varphi^5}_g\\
&\-\underbrace{\varphi^1\varphi^{21}\dots \varphi^{24}\otimes \shift \varphi^{251}\otimes \varphi^{252}\otimes \shift \varphi^{253}\otimes \shift \varphi^3\otimes \varphi^4\varphi^5}_y\\
&+\underbrace{\varphi^1\varphi^{21}\dots \varphi^{24}\otimes \shift \varphi^{25}\otimes \shift \varphi^{31}\otimes \varphi^{32}\otimes \shift \varphi^{33}\otimes \varphi^4\varphi^5}_g\\
&-\underbrace{\varphi^1\varphi^{21}\dots \varphi^{24}\varphi^{251}\varphi^{2521}\dots \varphi^{2524}\otimes \shift \varphi^{2525}\otimes \shift \varphi^{253}\otimes \shift \varphi^3\otimes \varphi^4\varphi^5}_g\\
&+\underbrace{\varphi^1\varphi^{21}\dots \varphi^{24}\otimes \shift \varphi^{25}\otimes \shift \varphi^{31}\otimes \shift \varphi^{321}\otimes \varphi^{322}\dots \varphi^{325}\varphi^{33}\otimes \varphi^4\varphi^5}_{gg}\\
&\-\underbrace{1\otimes \shift \varphi^1\otimes \shift \varphi^{21}\otimes \varphi^{22}\varphi^{23}\otimes \shift \varphi^{241}\otimes \varphi^{242}\dots \varphi^{245}\varphi^{25}\varphi^3\varphi^4\varphi^5}_{277}\\
&+\underbrace{1\otimes \shift \varphi^1\otimes \shift \varphi^{21}\otimes \varphi^{22}\varphi^{23}\varphi^{241}\varphi^{242}\otimes \shift \varphi^{243}\otimes \varphi^{244}\varphi^{245}\varphi^3\varphi^4\varphi^5}_y\\
&-\underbrace{1\otimes \shift \varphi^1\otimes \shift \varphi^{21}\otimes \varphi^{22}\varphi^{23}\varphi^{241}\dots \varphi^{244}\otimes \shift \varphi^{245}\otimes \varphi^{25}\varphi^3\varphi^4\varphi^5}_{260}\\
&\-\underbrace{\varphi^1\varphi^2\otimes \shift \varphi^3\otimes \shift \varphi^{41}\otimes 1\otimes \shift \varphi^{421}\otimes \varphi^{422}\dots \varphi^{425}\varphi^{43}\dots \varphi^{45}\varphi^5}_{338}\\
&+\underbrace{\varphi^1\varphi^2\otimes \shift \varphi^3\otimes \shift \varphi^{41}\otimes \varphi^{421}\varphi^{422}\otimes \shift \varphi^{423}\otimes \varphi^{424}\varphi^{425}\varphi^{43}\varphi^{44}\varphi^{45}\varphi^5}_{156}\\
&-\underbrace{\varphi^1\varphi^2\otimes \shift \varphi^3\otimes \shift \varphi^{41}\otimes \varphi^{421}\dots \varphi^{424}\otimes \shift \varphi^{425}\otimes \varphi^{43}\dots \varphi^{45}\varphi^5}_{160}\\
&\-\underbrace{\varphi^1\varphi^2\varphi^3\varphi^{41}\dots \varphi^{43}\otimes \shift \varphi^{441}\otimes \varphi^{442}\dots \varphi^{445}\otimes \shift \varphi^{45}\otimes \shift \varphi^5\otimes 1}_g\\
&+\underbrace{\varphi^1\varphi^2\varphi^3\varphi^{41}\dots \varphi^{43}\varphi^{441}\varphi^{442}\otimes \shift \varphi^{443}\otimes \varphi^{444}\varphi^{445}\otimes \shift \varphi^{45}\otimes \shift \varphi^5\otimes 1}_{177} \\
&-\underbrace{\varphi^1\varphi^2\varphi^3\varphi^{41}\varphi^{42}\varphi^{43}\varphi^{441}\dots \varphi^{444}\otimes \shift \varphi^{445}\otimes 1\otimes \shift \varphi^{45}\otimes \shift \varphi^5\otimes 1}_{344}
\end{align*}
and
\begin{align*}
&-(1^{\otimes 3}\otimes \mu_3)\mu_4(\varphi)\\
=\,&+\underbrace{\varphi^1\varphi^2\otimes \shift \varphi^3\otimes \varphi^4\otimes \shift \varphi^{51}\otimes \varphi^{52}\otimes \shift \varphi^{53}\otimes 1}_y\\
&+\underbrace{\varphi^1\varphi^2\otimes \shift \varphi^3\otimes \varphi^4\varphi^{51}\varphi^{521}\dots \varphi^{524}\otimes \shift \varphi^{525}\otimes \shift \varphi^{53}\otimes 1}_g\\
&
+\underbrace{\varphi^1\varphi^2\otimes \shift \varphi^3\otimes \varphi^4\otimes \shift \varphi^{51}\otimes \shift \varphi^{521}\otimes \varphi^{522}\dots \varphi^{525}\varphi^{53}}_g\\
&\-\underbrace{1\otimes \shift \varphi^1\otimes \varphi^2\dots \varphi^4\otimes \shift \varphi^{51}\otimes \varphi^{52}\otimes \shift \varphi^{53}\otimes 1}_y
-\underbrace{1\otimes \shift \varphi^1\otimes \varphi^2\varphi^3\varphi^4\varphi^{51}\varphi^{521}\dots \varphi^{524}\otimes \shift \varphi^{525}\otimes \shift \varphi^{53}\otimes 1}_g\\
&-\underbrace{1\otimes \shift \varphi^1\otimes \varphi^2\varphi^3\varphi^4\otimes \shift \varphi^{51}\otimes \shift \varphi^{521}\otimes \varphi^{522}\dots \varphi^{525}\varphi^{53}}_g\\
&\-\underbrace{1\otimes \shift \varphi^1\otimes \varphi^2\otimes \shift \varphi^3\otimes 1\otimes \shift \varphi^{41}\otimes \varphi^{42}\dots \varphi^{45}\varphi^5}_{339}
+\underbrace{1\otimes \shift \varphi^1\otimes \varphi^2\otimes \shift \varphi^3\otimes \varphi^{41}\varphi^{42}\otimes \shift \varphi^{43}\otimes \varphi^{44}\varphi^{45}\varphi^5}_{199}\\
&-\underbrace{1\otimes \shift \varphi^1\otimes \varphi^2\otimes \shift \varphi^3\otimes \varphi^{41}\dots \varphi^{44}\otimes \shift \varphi^{45}\otimes \varphi^5}_{205}\\
&\-\underbrace{\varphi^1\varphi^{21}\dots \varphi^{24}\otimes \shift \varphi^{25}\otimes \shift \varphi^3\otimes 1\otimes \shift \varphi^{41}\otimes \varphi^{42}\dots \varphi^{45}\varphi^5}_{340}\\
&+\underbrace{\varphi^1\varphi^{21}\dots \varphi^{24}\otimes \shift \varphi^{25}\otimes \shift \varphi^3\otimes \varphi^{41}\varphi^{42}\otimes \shift \varphi^{43}\otimes \varphi^{44}\varphi^{45}\varphi^5}_{168}\\
&
-\underbrace{\varphi^1\varphi^{21}\dots \varphi^{24}\otimes \shift \varphi^{25}\otimes \shift \varphi^3\otimes \varphi^{41}\dots \varphi^{44}\otimes \shift \varphi^{45}\otimes \varphi^5}_{169}\\
&\-\underbrace{1\otimes \shift \varphi^1\otimes \shift \varphi^{21}\otimes \varphi^{22}\dots \varphi^{25}\varphi^3\otimes \shift \varphi^{41}\otimes \varphi^{42}\dots \varphi^{45}\varphi^5}_{193}\\
&
+\underbrace{1\otimes \shift \varphi^1\otimes \shift \varphi^{21}\otimes \varphi^{22}\dots \varphi^{25}\varphi^3\varphi^{41}\varphi^{42}\otimes \shift \varphi^{43}\otimes \varphi^{44}\varphi^{45}\varphi^5}_g\\
&-\underbrace{1\otimes \shift \varphi^1\otimes \shift \varphi^{21}\otimes \varphi^{22}\dots \varphi^{25}\varphi^3\varphi^{41}\dots \varphi^{44}\otimes \shift \varphi^{45}\otimes \varphi^5}_{182}\\
&\-\underbrace{\varphi^1\varphi^2\otimes \shift \varphi^3\otimes \shift \varphi^{41}\otimes \varphi^{42}\varphi^{43}\otimes \shift \varphi^{441}\otimes \varphi^{442}\dots \varphi^{445}\varphi^{45}\varphi^5}_{293}\\
&+\underbrace{\varphi^1\varphi^2\otimes \shift \varphi^3\otimes \shift \varphi^{41}\otimes \varphi^{42}\varphi^{43}\varphi^{441}\varphi^{442}\otimes \shift \varphi^{443}\otimes \varphi^{444}\varphi^{445}\varphi^{45}\varphi^5}_y\\
&-\underbrace{\varphi^1\varphi^2\otimes \shift \varphi^3\otimes \shift \varphi^{41}\otimes \varphi^{42}\varphi^{43}\varphi^{441}\dots \varphi^{444}\otimes \shift \varphi^{445}\otimes \varphi^{45}\varphi^5}_{282}\\
&\-\underbrace{\varphi^1\varphi^2\varphi^3\varphi^{41}\dots \varphi^{44}\otimes \shift \varphi^{451}\otimes \varphi^{452}\otimes \shift \varphi^{453}\otimes \shift \varphi^5\otimes 1}_y\\
&
+\underbrace{\varphi^1\varphi^2\varphi^3\varphi^{41}\dots \varphi^{44}\otimes \shift \varphi^{45}\otimes \shift \varphi^{51}\otimes \varphi^{52}\otimes \shift \varphi^{53}\otimes 1}_g\\
&-\underbrace{\varphi^1\varphi^2\varphi^3\varphi^{41}\dots \varphi^{44}\varphi^{451}\varphi^{4521}\dots \varphi^{4524}\otimes \shift \varphi^{4525}\otimes \shift \varphi^{453}\otimes \shift \varphi^5\otimes 1}_g\\
&+\underbrace{\varphi^1\varphi^2\varphi^3\varphi^{41}\dots \varphi^{44}\otimes \shift \varphi^{45}\otimes \shift \varphi^{51}\otimes \shift \varphi^{521}\otimes \varphi^{522}\dots \varphi^{525}\varphi^{53}\otimes 1}_{gg}
\end{align*}

Writing $\mu_1^{\otimes}=(\mu_1\otimes 1^{\otimes 5}+1\otimes \mu_1\otimes 1^{\otimes 4}+1^{\otimes 2}\otimes \mu_1\otimes 1^{\otimes 3}+1^{\otimes 3}\otimes \mu_1\otimes 1^{\otimes 2}+1^{\otimes 4}\otimes \mu_1\otimes 1+1^{\otimes 5}\otimes \mu_1)$ and applying it to the individual summands of $\mu_6^0(\varphi)$ we obtain

\begin{align*}
&-\mu_1^{\otimes}(-1\otimes \shift \varphi^1\otimes \shift \varphi^{21}\otimes \varphi^{22}\varphi^{23}\varphi^{241}\dots \varphi^{244}\otimes \shift \varphi^{245}\otimes \shift \varphi^{25}\otimes \varphi^3\varphi^4\varphi^5)\\
=\,&+\underbrace{\varphi^1\otimes \shift \varphi^{21}\otimes \varphi^{22}\varphi^{23}\varphi^{241}\dots \varphi^{244}\otimes \shift \varphi^{245}\otimes \shift \varphi^{25}\otimes \varphi^3\varphi^4\varphi^5}_{257}\\
&-\underbrace{1\otimes \shift (\varphi^1\varphi^{21})\otimes \varphi^{22}\varphi^{23}\varphi^{241}\dots \varphi^{244}\otimes \shift \varphi^{245}\otimes \shift \varphi^{25}\otimes \varphi^3\varphi^4\varphi^5}_{258}\\
&+\underbrace{1\otimes \shift \varphi^1\otimes \varphi^{21}\varphi^{22}\varphi^{23}\varphi^{241}\dots \varphi^{244}\otimes \shift \varphi^{245}\otimes \shift \varphi^{25}\otimes \varphi^3\varphi^4\varphi^5}_{259}\\
&\boldsymbol{\uwave{+}}0\boldsymbol{\uwave{+}}0\boldsymbol{\uwave{+}}0\boldsymbol{\uwave{-}}\underbrace{1\otimes \shift \varphi^1\otimes \shift \varphi^{21}\otimes \varphi^{22}\varphi^{23}(\shift \otimes \shift)\partial(\shift^{-1}\varphi^{24})\otimes \shift \varphi^{25}\otimes \varphi^3\varphi^4\varphi^5}_{256}\\
&-\underbrace{1\otimes \shift \varphi^1\otimes \shift \varphi^{21}\otimes \varphi^{22}\varphi^{23}\varphi^{241}\dots \varphi^{244}\otimes \shift (\varphi^{245}\varphi^{25})\otimes \varphi^3\varphi^4\varphi^5}_y\\
&+\underbrace{1\otimes \shift \varphi^1\otimes \shift \varphi^{21}\otimes \varphi^{22}\varphi^{23}\varphi^{241}\dots \varphi^{244}\otimes \shift \varphi^{245}\otimes \varphi^{25}\varphi^3\varphi^4\varphi^5}_{260}\boldsymbol{\uwave{+}}0
\end{align*}
and
\begin{align*}
&-\mu_1^{\otimes}(-\varphi^1\varphi^{21}\dots \varphi^{23}\varphi^{241}\dots \varphi^{244}\otimes \shift \varphi^{245}\otimes \shift \varphi^{25}\otimes \shift \varphi^3\otimes \shift \varphi^{41}\otimes \varphi^{42}\dots \varphi^{45}\varphi^5)\\
=\,&0\boldsymbol{\uwave{+}}0\boldsymbol{\uwave{+}}0\boldsymbol{\uwave{+}}0\boldsymbol{\uwave{-}}\underbrace{\varphi^1\varphi^{21}\varphi^{22}\varphi^{23}(\shift \otimes \shift)\partial(\shift^{-1}\varphi^{24})\otimes \shift \varphi^{25}\otimes \shift \varphi^3\otimes \shift \varphi^{41}\otimes \varphi^{42}\dots \varphi^{45}\varphi^5}_{261}\\
&-\underbrace{\varphi^1\varphi^{21}\dots \varphi^{23}\varphi^{241}\dots \varphi^{244}\otimes \shift (\varphi^{245}\varphi^{25})\otimes \shift \varphi^3\otimes \shift \varphi^{41}\otimes \varphi^{42}\dots \varphi^{45}\varphi^5}_y\\
&+\underbrace{\varphi^1\varphi^{21}\dots \varphi^{23}\varphi^{241}\dots \varphi^{244}\otimes \shift \varphi^{245}\otimes \shift (\varphi^{25}\varphi^3)\otimes \shift \varphi^{41}\otimes \varphi^{42}\dots \varphi^{45}\varphi^5}_{gg}\\
&-\underbrace{\varphi^1\varphi^{21}\dots \varphi^{23}\varphi^{241}\dots \varphi^{244}\otimes \shift \varphi^{245}\otimes \shift \varphi^{25}\otimes \shift (\varphi^3\varphi^{41})\otimes \varphi^{42}\dots \varphi^{45}\varphi^5}_{gg}\\
&-\underbrace{\varphi^1\varphi^{21}\dots \varphi^{23}\varphi^{241}\dots \varphi^{244}\otimes \shift \varphi^{245}\otimes \shift \varphi^{25} \otimes \shift \varphi^3\otimes (\shift \otimes \shift)\partial(\shift^{-1}\varphi^4)\varphi^5}_{262}\boldsymbol{\uwave{+}}0
\end{align*}
and
\begin{align*}
&-\mu_1^{\otimes}(-\varphi^1\varphi^2\varphi^3\varphi^{41}\varphi^{421}\dots \varphi^{424}\otimes \shift \varphi^{425}\otimes \shift \varphi^{43}\otimes \varphi^{44}\otimes \shift \varphi^{45}\otimes \shift \varphi^5\otimes 1)\\
=\,&0\boldsymbol{\uwave{+}}0\boldsymbol{\uwave{+}}0\boldsymbol{\uwave{-}}\underbrace{\varphi^1\varphi^2\varphi^3\varphi^{41}(\shift \otimes \shift)\partial(\shift^{-1}\varphi^{42})\otimes \shift \varphi^{43}\otimes \varphi^{44}\otimes \shift \varphi^{45}\otimes \shift \varphi^5\otimes 1}_{263}\\
&-\underbrace{\varphi^1\varphi^2\varphi^3\varphi^{41}\varphi^{421}\dots \varphi^{424}\otimes \shift (\varphi^{425}\varphi^{43})\otimes \varphi^{44}\otimes \shift \varphi^{45}\otimes \shift \varphi^5\otimes 1}_y\\
&+\underbrace{\varphi^1\varphi^2\varphi^3\varphi^{41}\varphi^{421}\dots \varphi^{424}\otimes \shift \varphi^{425}\otimes \varphi^{43}\varphi^{44}\otimes \shift \varphi^{45}\otimes \shift \varphi^5\otimes 1}_{265}\\
&\boldsymbol{\uwave{+}}0\boldsymbol{\uwave{+}}\underbrace{\varphi^1\varphi^2\varphi^3\varphi^{41}\varphi^{421}\dots \varphi^{424}\otimes \shift \varphi^{425}\otimes \shift \varphi^{43}\otimes \varphi^{44}\varphi^{45}\otimes \shift \varphi^5\otimes 1}_{264}\\
&-\underbrace{\varphi^1\varphi^2\varphi^3\varphi^{41}\varphi^{421}\dots \varphi^{424}\otimes \shift \varphi^{425}\otimes \shift \varphi^{43}\otimes \varphi^{44}\otimes \shift (\varphi^{45}\varphi^5) \otimes 1}_g\\
&+\underbrace{\varphi^1\varphi^2\varphi^3\varphi^{41}\varphi^{421}\dots \varphi^{424}\otimes \shift \varphi^{425}\otimes \shift \varphi^{43}\otimes \varphi^{44}\otimes \shift \varphi^{45}\otimes \varphi^5}_{266}
\end{align*}
and
\begin{align*}
&-\mu_1^{\otimes}(-\varphi^1\varphi^{21}\varphi^{22}\otimes \shift \varphi^{23}\otimes \shift \varphi^{241}\otimes \varphi^{242}\dots \varphi^{245}\otimes \shift \varphi^{25}\otimes \shift \varphi^3\otimes \varphi^4\varphi^5)\\
=\,&0\boldsymbol{\uwave{+}}\underbrace{\varphi^1\varphi^{21}\dots \varphi^{23}\otimes \shift \varphi^{241}\otimes \varphi^{242}\dots \varphi^{245}\otimes \shift \varphi^{25}\otimes \shift \varphi^3\otimes \varphi^4\varphi^5}_{272}\\
&-\underbrace{\varphi^1\varphi^{21}\varphi^{22}\otimes \shift (\varphi^{23}\varphi^{241})\otimes \varphi^{242}\dots \varphi^{245}\otimes \shift \varphi^{25}\otimes \shift \varphi^3\otimes \varphi^4\varphi^5}_{268}\\
&-\underbrace{\varphi^1\varphi^{21}\varphi^{22}\otimes \shift \varphi^{23}\otimes (\shift \otimes \shift)\partial(\shift^{-1}\varphi^{24})\otimes \shift \varphi^{25}\otimes \shift \varphi^3\otimes \varphi^4\varphi^5}_{267}\\
&\boldsymbol{\uwave{+}}0\boldsymbol{\uwave{+}}0\boldsymbol{\uwave{+}}\underbrace{\varphi^1\varphi^{21}\varphi^{22}\otimes \shift \varphi^{23}\otimes \shift \varphi^{241}\otimes \varphi^{242}\dots \varphi^{245}\varphi^{25}\otimes \shift \varphi^3\otimes \varphi^4\varphi^5}_{269}\\
&-\underbrace{\varphi^1\varphi^{21}\varphi^{22}\otimes \shift \varphi^{23}\otimes \shift \varphi^{241}\otimes \varphi^{242}\dots \varphi^{245}\otimes \shift (\varphi^{25}\varphi^3)\otimes \varphi^4\varphi^5}_g\\
&+\underbrace{\varphi^1\varphi^{21}\varphi^{22}\otimes \shift \varphi^{23}\otimes \shift \varphi^{241}\otimes \varphi^{242}\dots \varphi^{245}\otimes \shift \varphi^{25}\otimes \varphi^3\varphi^4\varphi^5}_{271}\boldsymbol{\uwave{+}}0
\end{align*}
and
\begin{align*}
&-\mu_1^{\otimes}(-\varphi^1\varphi^2\otimes \shift \varphi^3\otimes \shift \varphi^{41}\otimes \varphi^{421}\dots \varphi^{424}\otimes \shift \varphi^{425}\otimes \shift \varphi^{43}\otimes \varphi^{44}\varphi^{45}\varphi^5)\\
=\,&0\boldsymbol{\uwave{+}}\underbrace{\varphi^1\varphi^2\varphi^3\otimes \shift \varphi^{41}\otimes \varphi^{421}\dots \varphi^{424}\otimes \shift \varphi^{425}\otimes \shift \varphi^{43}\otimes \varphi^{44}\varphi^{45}\varphi^5}_{274}\\
&-\underbrace{\varphi^1\varphi^2\otimes \shift (\varphi^3\varphi^{41})\otimes \varphi^{421}\dots \varphi^{424}\otimes \shift \varphi^{425}\otimes \shift \varphi^{43}\otimes \varphi^{44}\varphi^{45}\varphi^5}_g\\
&+\underbrace{\varphi^1\varphi^2\otimes \shift \varphi^3\otimes \varphi^{41}\varphi^{421}\dots \varphi^{424}\otimes \shift \varphi^{425}\otimes \shift \varphi^{43}\otimes \varphi^{44}\varphi^{45}\varphi^5}_{275}\\
&\boldsymbol{\uwave{+}}0\boldsymbol{\uwave{+}}0\boldsymbol{\uwave{-}}\underbrace{\varphi^1\varphi^2\otimes \shift \varphi^3\otimes \shift \varphi^{41} \otimes (\shift \otimes \shift)\partial(\shift^{-1}\varphi^{42})\otimes \shift \varphi^{43}\otimes \varphi^{44}\varphi^{45}\varphi^5}_{273}\\
&-\underbrace{\varphi^1\varphi^2\otimes \shift \varphi^3\otimes \shift \varphi^{41}\otimes \varphi^{421}\dots \varphi^{424}\otimes \shift (\varphi^{425}\varphi^{43})\otimes \varphi^{44}\varphi^{45}\varphi^5}_{159}\\
&+\underbrace{\varphi^1\varphi^2\otimes \shift \varphi^3\otimes \shift \varphi^{41}\otimes \varphi^{421}\dots \varphi^{424}\otimes \shift \varphi^{425}\otimes \varphi^{43}\dots \varphi^{45}\varphi^5}_{160}\boldsymbol{\uwave{+}}0
\end{align*}
and
\begin{align*}
&-\mu_1^{\otimes}(-1\otimes \shift \varphi^1\otimes \shift \varphi^{21}\otimes \varphi^{22}\otimes \shift \varphi^{23}\otimes \shift \varphi^{241}\otimes \varphi^{242}\dots \varphi^{245}\varphi^{25}\varphi^3\varphi^4\varphi^5)\\
=\,&+\underbrace{\varphi^1\otimes \shift \varphi^{21}\otimes \varphi^{22}\otimes \shift \varphi^{23}\otimes \shift \varphi^{241}\otimes \varphi^{242}\dots \varphi^{245}\varphi^{25}\varphi^3\varphi^4\varphi^5}_{278}\\
&-\underbrace{1\otimes \shift (\varphi^1\varphi^{21})\otimes \varphi^{22}\otimes \shift \varphi^{23}\otimes \shift \varphi^{241}\otimes \varphi^{242}\dots \varphi^{245}\varphi^{25}\varphi^3\varphi^4\varphi^5}_g\\
&+\underbrace{1\otimes \shift \varphi^1\otimes \varphi^{21}\varphi^{22}\otimes \shift \varphi^{23}\otimes \shift \varphi^{241}\otimes \varphi^{242}\dots \varphi^{245}\varphi^{25}\varphi^3\varphi^4\varphi^5}_{279}\\
&\boldsymbol{\uwave{+}}0\boldsymbol{\uwave{+}}\underbrace{1\otimes \shift \varphi^1\otimes \shift \varphi^{21}\otimes \varphi^{22}\varphi^{23}\otimes \shift \varphi^{241}\otimes \varphi^{242}\dots \varphi^{245}\varphi^{25}\varphi^3\varphi^4\varphi^5}_{277}\\
&-\underbrace{1\otimes \shift \varphi^1\otimes \shift \varphi^{21}\otimes \varphi^{22}\otimes \shift (\varphi^{23}\varphi^{241})\otimes \varphi^{242}\dots \varphi^{245}\varphi^{25}\varphi^3\varphi^4\varphi^5}_y\\
&-\underbrace{1\otimes \shift \varphi^1\otimes \shift \varphi^{21}\otimes \varphi^{22}\otimes \shift \varphi^{23}\otimes (\shift \otimes \shift)\partial(\shift^{-1}\varphi^{24})\varphi^{25}\varphi^3\varphi^4\varphi^5}_{276}\boldsymbol{\uwave{+}}0\boldsymbol{\uwave{+}}0\boldsymbol{\uwave{+}}0
\end{align*}
and
\begin{align*}
&-\mu_1^{\otimes}(-\varphi^1\varphi^2\otimes \shift \varphi^3\otimes \shift \varphi^{41}\otimes \varphi^{42}\varphi^{43}\varphi^{441}\dots \varphi^{444}\otimes \shift \varphi^{445}\otimes \shift \varphi^{45}\otimes \varphi^5)\\
=\,&0\boldsymbol{\uwave{+}}\underbrace{\varphi^1\varphi^2\varphi^3\otimes \shift \varphi^{41}\otimes \varphi^{42}\varphi^{43}\varphi^{441}\dots \varphi^{444}\otimes \shift \varphi^{445}\otimes \shift \varphi^{45}\otimes \varphi^5}_{281}\\
&-\underbrace{\varphi^1\varphi^2\otimes \shift (\varphi^3\varphi^{41})\otimes \varphi^{42}\varphi^{43}\varphi^{441}\dots \varphi^{444}\otimes \shift \varphi^{445}\otimes \shift \varphi^{45}\otimes \varphi^5}_{284}\\
&+\underbrace{\varphi^1\varphi^2\otimes \shift \varphi^3\otimes \varphi^{41}\dots \varphi^{43}\varphi^{441}\dots \varphi^{444}\otimes \shift \varphi^{445}\otimes \shift \varphi^{45}\otimes \varphi^5}_{283}\\
&\boldsymbol{\uwave{+}}0\boldsymbol{\uwave{+}}0\boldsymbol{\uwave{+}}0\boldsymbol{\uwave{-}}\underbrace{\varphi^1\varphi^2\otimes \shift \varphi^3\otimes \shift \varphi^{41}\otimes \varphi^{42}\varphi^{43}(\shift \otimes \shift)\partial(\shift^{-1}\varphi^{44})\otimes \shift \varphi^{45}\otimes \varphi^5}_{280}\\
&-\underbrace{\varphi^1\varphi^2\otimes \shift \varphi^3\otimes \shift \varphi^{41}\otimes \varphi^{42}\varphi^{43}\varphi^{441}\dots \varphi^{444}\otimes \shift (\varphi^{445}\varphi^{45})\otimes \varphi^5}_y\\
&+\underbrace{\varphi^1\varphi^2\otimes \shift \varphi^3\otimes \shift \varphi^{41}\otimes \varphi^{42}\varphi^{43}\varphi^{441}\dots \varphi^{444}\otimes \shift \varphi^{445}\otimes \varphi^{45}\varphi^5}_{282}
\end{align*}
and
\begin{align*}
&-\mu_1^{\otimes}(-\varphi^1\varphi^2\varphi^3\varphi^{41}\varphi^{42}\otimes \shift \varphi^{43}\otimes \shift \varphi^{441}\otimes \varphi^{442}\dots \varphi^{445}\otimes \shift \varphi^{45}\otimes \shift \varphi^5\otimes 1)\\
=\,&0\boldsymbol{\uwave{+}}0\boldsymbol{\uwave{+}}\underbrace{\varphi^1\varphi^2\varphi^3\varphi^{41}\dots \varphi^{43}\otimes \shift \varphi^{441}\otimes \varphi^{442}\dots \varphi^{445}\otimes \shift \varphi^{45}\otimes \shift \varphi^5\otimes 1}_{g}\\
&-\underbrace{\varphi^1\varphi^2\varphi^3\varphi^{41}\varphi^{42}\otimes \shift (\varphi^{43}\varphi^{441})\otimes \varphi^{442}\dots \varphi^{445}\otimes \shift \varphi^{45}\otimes \shift \varphi^5\otimes 1}_{288}\\
&-\underbrace{\varphi^1\varphi^2\varphi^3\varphi^{41}\varphi^{42}\otimes \shift \varphi^{43}\otimes (\shift \otimes \shift)\partial(\shift^{-1}\varphi^{44})\otimes \shift \varphi^{45}\otimes \shift \varphi^5\otimes 1}_{285}\\
&\boldsymbol{\uwave{+}}0\boldsymbol{\uwave{+}}0\boldsymbol{\uwave{+}}\underbrace{\varphi^1\varphi^2\varphi^3\varphi^{41}\varphi^{42}\otimes \shift \varphi^{43}\otimes \shift \varphi^{441}\otimes \varphi^{442}\dots \varphi^{445}\varphi^{45}\otimes \shift \varphi^5\otimes 1}_{286}\\
&-\underbrace{\varphi^1\varphi^2\varphi^3\varphi^{41}\varphi^{42}\otimes \shift \varphi^{43}\otimes \shift \varphi^{441}\otimes \varphi^{442}\dots \varphi^{445}\otimes \shift (\varphi^{45}\varphi^5)\otimes 1}_g\\
&+\underbrace{\varphi^1\varphi^2\varphi^3\varphi^{41}\varphi^{42}\otimes \shift \varphi^{43}\otimes \shift \varphi^{441}\otimes \varphi^{442}\dots \varphi^{445}\otimes \shift \varphi^{45}\otimes \varphi^5}_{287}
\end{align*}
and
\begin{align*}
&-\mu_1^{\otimes}(-\varphi^1\varphi^{21}\dots \varphi^{23}\varphi^{241}\dots \varphi^{243}\varphi^{2441}\dots \varphi^{2444}\otimes \shift \varphi^{2445}\otimes \shift \varphi^{245}\otimes \shift \varphi^{25}\otimes \shift \varphi^3\otimes \varphi^4\varphi^5)\\
=\,&0\boldsymbol{\uwave{+}}0\boldsymbol{\uwave{+}}0\boldsymbol{\uwave{+}}0\boldsymbol{\uwave{-}}\underbrace{\varphi^1\varphi^{21}\dots \varphi^{23}\varphi^{241}\dots \varphi^{243}(\shift \otimes \shift)\partial(\shift^{-1}\varphi^{244})\otimes \shift \varphi^{245}\otimes \shift \varphi^{25}\otimes \shift \varphi^3\otimes \varphi^4\varphi^5}_{291}\\
&-\underbrace{\varphi^1\varphi^{21}\dots \varphi^{23}\varphi^{241}\dots \varphi^{243}\varphi^{2441}\dots \varphi^{2444}\otimes \shift (\varphi^{2445}\varphi^{245})\otimes \shift \varphi^{25}\otimes \shift \varphi^3\otimes \varphi^4\varphi^5}_y\\
&+\underbrace{\varphi^1\varphi^{21}\dots \varphi^{23}\varphi^{241}\dots \varphi^{243}\varphi^{2441}\dots \varphi^{2444}\otimes \shift \varphi^{2445}\otimes \shift (\varphi^{245}\varphi^{25})\otimes \shift \varphi^3\otimes \varphi^4\varphi^5}_g\\
&-\underbrace{\varphi^1\varphi^{21}\dots \varphi^{23}\varphi^{241}\dots \varphi^{243}\varphi^{2441}\dots \varphi^{2444}\otimes \shift \varphi^{2445}\otimes \shift \varphi^{245}\otimes \shift (\varphi^{25}\varphi^3)\otimes \varphi^4\varphi^5}_{gg}\\
&+\underbrace{\varphi^1\varphi^{21}\dots \varphi^{23}\varphi^{241}\dots \varphi^{243}\varphi^{2441}\dots \varphi^{2444}\otimes \shift \varphi^{2445}\otimes \shift \varphi^{245}\otimes \shift \varphi^{25}\otimes \varphi^3\varphi^4\varphi^5}_{296}\boldsymbol{\uwave{+}}0 
\end{align*}
and
\begin{align*}
&-\mu_1^{\otimes}(-1\otimes \shift \varphi^1\otimes \varphi^2\otimes \shift \varphi^3\otimes \shift \varphi^{41}\otimes \varphi^{42}\otimes \shift \varphi^{43}\otimes \varphi^{44}\varphi^{45}\varphi^5)\\
=\,&+\underbrace{\hat{\omega}\varphi^1\varphi^2\otimes \shift \varphi^3\otimes \shift \varphi^{41}\otimes \varphi^{42}\otimes \shift \varphi^{43}\otimes \varphi^{44}\varphi^{45}\varphi^5}_{81}\\
&-\underbrace{\varphi^1\hat{\omega}\varphi^2\otimes \shift \varphi^3\otimes \shift \varphi^{41}\otimes \varphi^{42}\otimes \shift \varphi^{43}\otimes \varphi^{44}\varphi^{45}\varphi^5}_{93}\\
&-\underbrace{\shift \partial(\varphi^1)\varphi^2\otimes \shift \varphi^3\otimes \shift \varphi^{41}\otimes \varphi^{42}\otimes \shift \varphi^{43}\otimes \varphi^{44}\varphi^{45}\varphi^5}_g\\
&\boldsymbol{\uwave{+}}0\boldsymbol{\uwave{-}}\underbrace{1\otimes \shift \varphi^1\otimes \varphi^2\varphi^3\otimes \shift \varphi^{41}\otimes \varphi^{42}\otimes \shift \varphi^{43}\otimes \varphi^{44}\varphi^{45}\varphi^5}_{198}\\
&+\underbrace{1\otimes \shift \varphi^1\otimes \varphi^2\otimes \shift (\varphi^3\varphi^{41})\otimes \varphi^{42}\otimes \shift \varphi^{43}\otimes \varphi^{44}\varphi^{45}\varphi^5}_y\\
&-\underbrace{1\otimes \shift \varphi^1\otimes \varphi^2\otimes \shift \varphi^3\otimes \varphi^{41}\varphi^{42}\otimes \shift \varphi^{43}\otimes \varphi^{44}\varphi^{45}\varphi^5}_{199}\\
&\boldsymbol{\uwave{+}}0\boldsymbol{\uwave{-}}\underbrace{1\otimes \shift \varphi^1\otimes \varphi^2\otimes \shift \varphi^3\otimes \shift \varphi^{41}\otimes \varphi^{42}\hat{\omega}\varphi^{43}\varphi^{44}\varphi^{45}\varphi^5}_{91}\\
&+\underbrace{1\otimes \shift \varphi^1\otimes \varphi^2\otimes \shift \varphi^3\otimes \shift \varphi^{41}\otimes \varphi^{42}\varphi^{43}\hat{\omega}\varphi^{44}\varphi^{45}\varphi^5}_{92}\\
&+\underbrace{1\otimes \shift \varphi^1\otimes \varphi^2\otimes \shift \varphi^3\otimes \shift \varphi^{41}\otimes \varphi^{42}\shift \partial(\varphi^{43})\varphi^{44}\varphi^{45}\varphi^5}_y\boldsymbol{\uwave{+}}0
\end{align*}
and
\begin{align*}
&-\mu_1^{\otimes}(-\varphi^1\otimes \shift \varphi^{21}\otimes \varphi^{22}\otimes \shift \varphi^{23}\otimes \varphi^{24}\otimes \shift \varphi^{25}\otimes \shift \varphi^3\otimes \varphi^4\varphi^5)\\
=\,&+\underbrace{\varphi^1\hat{\omega}\varphi^{21}\varphi^{22}\otimes \shift \varphi^{23}\otimes \varphi^{24}\otimes \shift \varphi^{25}\otimes \shift \varphi^3\otimes \varphi^4\varphi^5}_{130}\\
&-\underbrace{\varphi^1\varphi^{21}\hat{\omega}\varphi^{22}\otimes \shift \varphi^{23}\otimes \varphi^{24}\otimes \shift \varphi^{25}\otimes \shift \varphi^3\otimes \varphi^4\varphi^5}_{131}\\
&-\underbrace{\varphi^1\shift \partial(\varphi^{21})\varphi^{22}\otimes \shift \varphi^{23}\otimes \varphi^{24}\otimes \shift \varphi^{25}\otimes \shift \varphi^3\otimes \varphi^4\varphi^5}_y\\
&\boldsymbol{\uwave{+}}0\boldsymbol{\uwave{-}}\underbrace{\varphi^1\otimes \shift \varphi^{21}\otimes \varphi^{22}\hat{\omega}\varphi^{23}\varphi^{24}\otimes \shift \varphi^{25}\otimes \shift \varphi^3\otimes \varphi^4\varphi^5}_{129}\\
&+\underbrace{\varphi^1\otimes \shift \varphi^{21}\otimes \varphi^{22}\varphi^{23}\hat{\omega}\varphi^{24}\otimes \shift \varphi^{25}\otimes \shift \varphi^3\otimes \varphi^4\varphi^5}_{132}\\
&+\underbrace{\varphi^1\otimes \shift \varphi^{21}\otimes \varphi^{22}\shift \partial(\varphi^{23})\varphi^{24}\otimes \shift \varphi^{25}\otimes \shift \varphi^3\otimes \varphi^4\varphi^5}_y\\
&\boldsymbol{\uwave{+}}0\boldsymbol{\uwave{+}}\underbrace{\varphi^1\otimes \shift \varphi^{21}\otimes \varphi^{22}\otimes \shift \varphi^{23}\otimes \varphi^{24}\varphi^{25}\otimes \shift \varphi^3\otimes \varphi^4\varphi^5}_{202}\\
&-\underbrace{\varphi^1\otimes \shift \varphi^{21}\otimes \varphi^{22}\otimes \shift \varphi^{23}\otimes \varphi^{24}\otimes \shift (\varphi^{25}\varphi^3)\otimes \varphi^4\varphi^5}_{200}\\
&+\underbrace{\varphi^1\otimes \shift \varphi^{21}\otimes \varphi^{22} \otimes \shift \varphi^{23}\otimes \varphi^{24}\otimes \shift \varphi^{25}\otimes \varphi^3\varphi^4\varphi^5}_{201}\boldsymbol{\uwave{+}}0
\end{align*}
and
\begin{align*}
&-\mu_1^{\otimes}(\varphi^1\otimes \shift \varphi^{21}\otimes \varphi^{22}\dots \varphi^{24}\otimes \shift \varphi^{25}\otimes \shift \varphi^3\otimes \varphi^4\otimes \shift \varphi^5\otimes 1)\\
=\,&-\underbrace{\varphi^1\hat{\omega}\varphi^{21}\dots \varphi^{24}\otimes \shift \varphi^{25}\otimes \shift \varphi^3\otimes \varphi^4\otimes \shift \varphi^5\otimes 1}_{72}\\
&+\underbrace{\varphi^1\varphi^{21}\hat{\omega}\varphi^{22}\dots \varphi^{24}\otimes \shift \varphi^{25}\otimes \shift \varphi^3\otimes \varphi^4\otimes \shift \varphi^5\otimes 1}_{73}\\
&+\underbrace{\varphi^1\shift \partial(\varphi^{21})\varphi^{22}\dots \varphi^{24}\otimes \shift \varphi^{25}\otimes \shift \varphi^3\otimes \varphi^4\otimes \shift \varphi^5\otimes 1}_y\\
&\boldsymbol{\uwave{+}}0\boldsymbol{\uwave{+}}0\boldsymbol{\uwave{+}}\underbrace{\varphi^1\otimes \shift \varphi^{21}\otimes \varphi^{22}\dots \varphi^{25}\otimes \shift \varphi^3\otimes \varphi^4\otimes \shift \varphi^5\otimes 1}_{204}\\
&-\underbrace{\varphi^1\otimes \shift \varphi^{21}\otimes \varphi^{22}\varphi^{23}\varphi^{24}\otimes \shift (\varphi^{25}\varphi^3)\otimes \varphi^4\otimes \shift \varphi^5\otimes 1}_{153}\\
&+\underbrace{\varphi^1\otimes \shift \varphi^{21}\otimes \varphi^{22}\dots \varphi^{24}\otimes \shift \varphi^{25}\otimes \varphi^3\varphi^4\otimes \shift \varphi^5\otimes 1}_{203}\\
&\boldsymbol{\uwave{+}}0\boldsymbol{\uwave{+}}\underbrace{\varphi^1\otimes \shift \varphi^{21}\otimes \varphi^{22}\dots \varphi^{24}\otimes \shift \varphi^{25}\otimes \shift \varphi^3\otimes \varphi^4\hat{\omega}\varphi^5}_{74}\\
&-\underbrace{\varphi^1\otimes \shift \varphi^{21}\otimes \varphi^{22}\dots \varphi^{24}\otimes \shift \varphi^{25}\otimes \shift \varphi^3\otimes \varphi^4\varphi^5\hat{\omega}}_{71}\\
&-\underbrace{\varphi^1\otimes \shift \varphi^{21}\otimes \varphi^{22}\dots \varphi^{24}\otimes \shift \varphi^{25}\otimes \shift \varphi^3\otimes \varphi^4\shift \partial(\varphi^5)}_g
\end{align*}
and
\begin{align*}
&-\mu_1^{\otimes}(1\otimes \shift \varphi^1\otimes \varphi^2\otimes \shift \varphi^3\otimes \shift \varphi^{41}\otimes \varphi^{42}\dots \varphi^{44}\otimes \shift \varphi^{45}\otimes \varphi^5)\\
=\,&-\underbrace{\hat{\omega}\varphi^1\varphi^2\otimes \shift \varphi^3\otimes \shift \varphi^{41}\otimes \varphi^{42}\dots \varphi^{44}\otimes \shift \varphi^{45}\otimes \varphi^5}_{79}\\
&+\underbrace{\varphi^1\hat{\omega}\varphi^2\otimes \shift \varphi^3\otimes \shift \varphi^{41}\otimes \varphi^{42}\dots \varphi^{44}\otimes \shift \varphi^{45}\otimes \varphi^5}_{80}\\
&+\underbrace{\shift \partial(\varphi^1)\varphi^2\otimes \shift \varphi^3\otimes \shift \varphi^{41}\otimes \varphi^{42}\dots \varphi^{44}\otimes \shift \varphi^{45}\otimes \varphi^5}_g\\
&\boldsymbol{\uwave{+}}0\boldsymbol{\uwave{+}}\underbrace{1\otimes \shift \varphi^1\otimes \varphi^2\varphi^3\otimes \shift \varphi^{41}\otimes \varphi^{42}\dots \varphi^{44}\otimes \shift \varphi^{45}\otimes \varphi^5}_{206}\\
&-\underbrace{1\otimes \shift \varphi^1\otimes \varphi^2\otimes \shift (\varphi^3\varphi^{41})\otimes \varphi^{42}\varphi^{43}\varphi^{44}\otimes \shift \varphi^{45}\otimes \varphi^5}_{194}\\
&+\underbrace{1\otimes \shift \varphi^1\otimes \varphi^2\otimes \shift \varphi^3\otimes \varphi^{41}\dots \varphi^{44}\otimes \shift \varphi^{45}\otimes \varphi^5}_{205}\\
&\boldsymbol{\uwave{+}}0\boldsymbol{\uwave{+}}0\boldsymbol{\uwave{+}}\underbrace{1\otimes \shift \varphi^1\otimes \varphi^2\otimes \shift \varphi^3\otimes \shift \varphi^{41}\otimes \varphi^{42}\dots \varphi^{44}\hat{\omega}\varphi^{45}\varphi^5}_{77}\\
&-\underbrace{1\otimes \shift \varphi^1\otimes \varphi^2\otimes \shift \varphi^3\otimes \shift \varphi^{41}\otimes \varphi^{42}\dots \varphi^{45}\hat{\omega}\varphi^5}_{90}\\
&-\underbrace{1\otimes \shift \varphi^1\otimes \varphi^2\otimes \shift \varphi^3\otimes \shift \varphi^{41}\otimes \varphi^{42}\dots \varphi^{44}\otimes \varphi^{45}\shift \partial(\varphi^5)}_g
\end{align*}
and
\begin{align*}
&-\mu_1^\otimes(-\varphi^1\varphi^2\otimes \shift \varphi^3\otimes \shift \varphi^{41}\otimes \varphi^{42}\otimes \shift \varphi^{43}\otimes \varphi^{44}\otimes \shift \varphi^{45}\otimes \varphi^5)\\
=\,&0\boldsymbol{\uwave{+}}\underbrace{\varphi^1\varphi^2\varphi^3\otimes \shift \varphi^{41}\otimes \varphi^{42}\otimes \shift \varphi^{43}\otimes \varphi^{44}\otimes \shift \varphi^{45}\otimes \varphi^5}_{207}\\
&-\underbrace{\varphi^1\varphi^2\otimes \shift (\varphi^3\varphi^{41})\otimes \varphi^{42}\otimes \shift \varphi^{43}\otimes \varphi^{44}\otimes \shift \varphi^{45}\otimes \varphi^5}_{208}\\
&+\underbrace{\varphi^1\varphi^2\otimes \shift \varphi^3\otimes \varphi^{41}\varphi^{42}\otimes \shift \varphi^{43}\otimes \varphi^{44}\otimes \shift \varphi^{45}\otimes \varphi^5}_{221}\\
&0\boldsymbol{\uwave{+}}\underbrace{\varphi^1\varphi^2\otimes \shift \varphi^3\otimes \shift \varphi^{41}\otimes \varphi^{42}\hat{\omega}\varphi^{43}\varphi^{44}\otimes \shift \varphi^{45}\otimes \varphi^5}_{115}\\
&-\underbrace{\varphi^1\varphi^2\otimes \shift \varphi^3\otimes \shift \varphi^{41}\otimes \varphi^{42}\varphi^{43}\hat{\omega}\varphi^{44}\otimes \shift \varphi^{45}\otimes \varphi^5}_{116}\\
&-\underbrace{\varphi^1\varphi^2\otimes \shift \varphi^3\otimes \shift \varphi^{41}\otimes \varphi^{42}\shift \partial(\varphi^{43})\varphi^{44}\otimes \shift \varphi^{45}\otimes \varphi^5}_y\\
&\boldsymbol{\uwave{+}}0\boldsymbol{\uwave{-}}\underbrace{\varphi^1\varphi^2\otimes \shift \varphi^3\otimes \shift \varphi^{41}\otimes \varphi^{42}\otimes \shift \varphi^{43}\otimes \varphi^{44}\hat{\omega}\varphi^{45}\varphi^5}_{111}\\
&+\underbrace{\varphi^1\varphi^2\otimes \shift \varphi^3\otimes \shift \varphi^{41}\otimes \varphi^{42}\otimes \shift \varphi^{43}\otimes \varphi^{44}\varphi^{45}\hat{\omega}\varphi^5}_{112}\\
&+\underbrace{\varphi^1\varphi^2\otimes \shift \varphi^3\otimes \shift \varphi^{41}\otimes \varphi^{42}\otimes \shift \varphi^{43}\otimes \varphi^{44}\shift \partial(\varphi^{45})\varphi^5}_y
\end{align*}
and
\begin{align*}
&-\mu_1^\otimes(-\varphi^1\varphi^2\otimes \shift (\varphi^3\varphi^{41})\otimes \varphi^{42}\otimes \shift \varphi^{43}\otimes \varphi^{44}\otimes \shift \varphi^{45}\otimes \shift \varphi^5\otimes 1)\\
=\,&0\boldsymbol{\uwave{+}}\underbrace{\varphi^1\varphi^2\hat{\omega}\varphi^3\varphi^{41}\varphi^{42}\otimes \shift \varphi^{43}\otimes \varphi^{44}\otimes \shift \varphi^{45}\otimes \shift \varphi^5\otimes 1}_{135}\\
&-\underbrace{\varphi^1\varphi^2\varphi^3\varphi^{41}\hat{\omega}\varphi^{42}\otimes \shift \varphi^{43}\otimes \varphi^{44}\otimes \shift \varphi^{45}\otimes \shift \varphi^5\otimes 1}_{136}\\
&-\underbrace{\varphi^1\varphi^2\varphi^3\shift \partial(\varphi^{41})\varphi^{42}\otimes \shift \varphi^{43}\otimes \varphi^{44}\otimes \shift \varphi^{45}\otimes \shift \varphi^5\otimes 1}_y\\
&-\underbrace{\varphi^1\varphi^2\shift \partial(\varphi^3)\varphi^{41}\varphi^{42}\otimes \shift \varphi^{43}\otimes \varphi^{44}\otimes \shift \varphi^{45}\otimes \shift \varphi^5\otimes 1}_g\\
&\boldsymbol{\uwave{+}}0\boldsymbol{\uwave{-}}\underbrace{\varphi^1\varphi^2\otimes \shift (\varphi^3\varphi^{41})\otimes \varphi^{42}\hat{\omega}\varphi^{43}\varphi^{44}\otimes \shift \varphi^{45}\otimes \shift \varphi^5\otimes 1}_{137}\\
&+\underbrace{\varphi^1\varphi^2\otimes \shift (\varphi^3\varphi^{41})\otimes \varphi^{42}\varphi^{43}\hat{\omega}\varphi^{44}\otimes \shift \varphi^{45}\otimes \shift \varphi^5\otimes 1}_{138}\\
&+\underbrace{\varphi^1\varphi^2\otimes \shift (\varphi^3\varphi^{41})\otimes \varphi^{42}\shift \partial(\varphi^{43})\varphi^{44}\otimes \shift \varphi^{45}\otimes \shift \varphi^5\otimes 1}_y\\
&\boldsymbol{\uwave{+}}0\boldsymbol{\uwave{+}}\underbrace{\varphi^1\varphi^2\otimes \shift (\varphi^3\varphi^{41})\otimes \varphi^{42}\otimes \shift \varphi^{43}\otimes \varphi^{44}\varphi^{45}\otimes \shift \varphi^5\otimes 1}_{209}\\
&-\underbrace{\varphi^1\varphi^2\otimes \shift (\varphi^3\varphi^{41})\otimes \varphi^{42}\otimes \shift \varphi^{43}\otimes \varphi^{44}\otimes \shift (\varphi^{45}\varphi^5)\otimes 1}_y\\
&+\underbrace{\varphi^1\varphi^2\otimes \shift (\varphi^3\varphi^{41})\otimes \varphi^{42}\otimes \shift \varphi^{43}\otimes \varphi^{44}\otimes \shift \varphi^{45}\otimes \varphi^5}_{208}
\end{align*}
and
\begin{align*}
&-\mu_1^{\otimes}(-1\otimes \shift \varphi^1\otimes \shift \varphi^{21}\otimes \varphi^{22}\otimes \shift \varphi^{23}\otimes \varphi^{24}\otimes \shift (\varphi^{25}\varphi^3)\otimes \varphi^4\varphi^5)\\
=\,&\underbrace{\varphi^1\otimes \shift \varphi^{21}\otimes \varphi^{22}\otimes \shift \varphi^{23}\otimes \varphi^{24}\otimes \shift (\varphi^{25}\varphi^3)\otimes \varphi^4\varphi^5}_{200}\\
&-\underbrace{1\otimes \shift (\varphi^1\varphi^{21})\otimes \varphi^{22}\otimes \shift \varphi^{23}\otimes \varphi^{24}\otimes \shift (\varphi^{25}\varphi^3)\otimes \varphi^4\varphi^5}_y\\
&+\underbrace{1\otimes \shift \varphi^1\otimes \varphi^{21}\varphi^{22}\otimes \shift \varphi^{23}\otimes \varphi^{24}\otimes \shift (\varphi^{25}\varphi^3)\otimes \varphi^4\varphi^5}_{224}\\
&\boldsymbol{\uwave{+}}0\boldsymbol{\uwave{+}}\underbrace{1\otimes \shift \varphi^1\otimes \shift \varphi^{21}\otimes \varphi^{22}\hat{\omega}\varphi^{23}\varphi^{24}\otimes \shift (\varphi^{25}\varphi^3)\otimes \varphi^4\varphi^5}_{143}\\
&-\underbrace{1\otimes \shift \varphi^1\otimes \shift \varphi^{21}\otimes \varphi^{22}\varphi^{23}\hat{\omega}\varphi^{24}\otimes \shift (\varphi^{25}\varphi^3)\otimes \varphi^4\varphi^5}_{144}\\
&-\underbrace{1\otimes \shift \varphi^1\otimes \shift \varphi^{21}\otimes \varphi^{22}\shift \partial(\varphi^{23})\varphi^{24}\otimes \shift (\varphi^{25}\varphi^3)\otimes \varphi^4\varphi^5}_y\\
&\boldsymbol{\uwave{+}}0\boldsymbol{\uwave{-}}\underbrace{1\otimes \shift \varphi^1\otimes \shift \varphi^{21}\otimes \varphi^{22}\otimes \shift \varphi^{23}\otimes \varphi^{24}\hat{\omega}\varphi^{25}\varphi^3\varphi^4\varphi^5}_{141}\\
&+\underbrace{1\otimes \shift \varphi^1\otimes \shift \varphi^{21}\otimes \varphi^{22}\otimes \shift \varphi^{23}\otimes \varphi^{24}\varphi^{25}\varphi^3\hat{\omega}\varphi^4\varphi^5}_{142}\\
&+\underbrace{1\otimes \shift \varphi^1\otimes \shift \varphi^{21}\otimes \varphi^{22}\otimes  \shift \varphi^{23}\otimes \varphi^{24}\varphi^{25}\shift \partial(\varphi^3)\varphi^4\varphi^5}_y\\
&+\underbrace{1\otimes \shift \varphi^1\otimes \shift \varphi^{21}\otimes \varphi^{22}\otimes \shift \varphi^{23}\otimes \varphi^{24}\shift \partial(\varphi^{25})\varphi^3\varphi^4\varphi^5}_g\boldsymbol{\uwave{+}}0
\end{align*}
and
\begin{align*}
&-\mu_1^{\otimes}(-\varphi^1\varphi^{21}\varphi^{22}\otimes \shift \varphi^{23}\otimes \varphi^{24}\otimes \shift \varphi^{25}\otimes \shift \varphi^3\otimes \varphi^4\otimes \shift \varphi^5\otimes 1)\\
=\,&0\boldsymbol{\uwave{+}}\underbrace{\varphi^1\varphi^{21}\varphi^{22}\hat{\omega}\varphi^{23}\varphi^{24}\otimes \shift \varphi^{25}\otimes \shift \varphi^3\otimes \varphi^4\otimes \shift \varphi^5\otimes 1}_{104}
\\
&-\underbrace{\varphi^1\varphi^{21}\dots \varphi^{23}\hat{\omega}\varphi^{24}\otimes \shift \varphi^{25}\otimes \shift \varphi^3\otimes \varphi^4\otimes \shift \varphi^5\otimes 1}_{105}\\
&-\underbrace{\varphi^1\varphi^{21}\varphi^{22}\shift \partial(\varphi^{23})\varphi^{24}\otimes \shift \varphi^{25}\otimes \shift \varphi^3\otimes \varphi^4\otimes \shift \varphi^5\otimes 1}_y\\
&\boldsymbol{\uwave{+}}0\boldsymbol{\uwave{-}}\underbrace{\varphi^1\varphi^{21}\varphi^{22}\otimes \shift \varphi^{23}\otimes \varphi^{24}\varphi^{25}\otimes \shift \varphi^3\otimes \varphi^4\otimes \shift \varphi^5\otimes 1}_{213}\\
&+\underbrace{\varphi^1\varphi^{21}\varphi^{22}\otimes \shift \varphi^{23}\otimes \varphi^{24}\otimes \shift (\varphi^{25}\varphi^3)\otimes \varphi^4\otimes \shift \varphi^5\otimes 1}_y\\
&-\underbrace{\varphi^1\varphi^{21}\varphi^{22}\otimes \shift \varphi^{23}\otimes \varphi^{24}\otimes \shift \varphi^{25}\otimes \varphi^3\varphi^4\otimes \shift \varphi^5\otimes 1}_{214}\\
&\boldsymbol{\uwave{+}}0\boldsymbol{\uwave{-}}\underbrace{\varphi^1\varphi^{21}\varphi^{22}\otimes \shift \varphi^{23}\otimes \varphi^{24}\otimes \shift \varphi^{25}\otimes \shift \varphi^3\otimes \varphi^4\hat{\omega}\varphi^5}_{103}\\
&+\underbrace{\varphi^1\varphi^{21}\varphi^{22}\otimes \shift \varphi^{23}\otimes \varphi^{24}\otimes \shift \varphi^{25}\otimes \shift \varphi^3\otimes \varphi^4\varphi^5\hat{\omega}}_{102}\\
&+\underbrace{\varphi^1\varphi^{21}\varphi^{22}\otimes \shift \varphi^{23}\otimes \varphi^{24}\otimes \shift \varphi^{25}\otimes \shift \varphi^3\otimes \varphi^4\shift \partial(\varphi^5)}_y
\end{align*}
and
\begin{align*}
&-\mu_1^{\otimes}(-\varphi^1\varphi^{21}\dots \varphi^{23}\varphi^{241}\varphi^{242}\otimes \shift \varphi^{243}\otimes \varphi^{244}\otimes \shift \varphi^{245}\otimes \shift \varphi^{25}\otimes \shift \varphi^3\otimes \varphi^4\varphi^5)\\
=\,&0\boldsymbol{\uwave{+}}0\boldsymbol{\uwave{+}}\underbrace{\varphi^1\varphi^{21}\dots \varphi^{23}\varphi^{241}\varphi^{242}\hat{\omega}\varphi^{243}\varphi^{244}\otimes \shift \varphi^{245}\otimes \shift \varphi^{25}\otimes \shift \varphi^3\otimes \varphi^4\varphi^5}_{319}\\
&-\underbrace{\varphi^1\varphi^{21}\dots \varphi^{23}\varphi^{241}\dots \varphi^{243}\hat{\omega}\varphi^{244}\otimes \shift \varphi^{245}\otimes \shift \varphi^{25}\otimes \shift \varphi^3\otimes \varphi^4\varphi^5}_{320}\\
&-\underbrace{\varphi^1\varphi^{21}\dots \varphi^{23}\varphi^{241}\varphi^{242}\shift \partial(\varphi^{243})\varphi^{244}\otimes \shift \varphi^{245}\otimes \shift \varphi^{25}\otimes \shift \varphi^3\otimes \varphi^4\varphi^5}_y\\
&\boldsymbol{\uwave{+}}0\boldsymbol{\uwave{-}}\underbrace{\varphi^1\varphi^{21}\dots \varphi^{23}\varphi^{241}\varphi^{242}\otimes \shift \varphi^{243}\otimes \varphi^{244}\varphi^{245}\otimes \shift \varphi^{25}\otimes \shift \varphi^3\otimes \varphi^4\varphi^5}_{175}\\
&+\underbrace{\varphi^1\varphi^{21}\dots \varphi^{23}\varphi^{241}\varphi^{242}\otimes \shift \varphi^{243}\otimes \varphi^{244}\otimes \shift (\varphi^{245}\varphi^{25})\otimes \shift \varphi^3\otimes \varphi^4\varphi^5}_y\\
&-\underbrace{\varphi^1\varphi^{21}\dots \varphi^{23}\varphi^{241}\varphi^{242}\otimes \shift \varphi^{243}\otimes \varphi^{244}\otimes \shift \varphi^{245}\otimes \shift \varphi^{25}\varphi^3\otimes \varphi^4\varphi^5}_g\\
&+\underbrace{\varphi^1\varphi^{21}\varphi^{22} \varphi^{23}\varphi^{241}\varphi^{242}\otimes \shift \varphi^{243}\otimes \varphi^{244}\otimes \shift \varphi^{245}\otimes \shift \varphi^{25}\otimes \varphi^3\varphi^4\varphi^5}_{174}\boldsymbol{\uwave{+}}0
\end{align*}
and
\begin{align*}
&-\mu_1^{\otimes}(-\varphi^1\varphi^2\varphi^3\varphi^{41}\dots \varphi^{43}\varphi^{441}\varphi^{442}\otimes \shift \varphi^{443}\otimes \varphi^{444}\otimes \shift \varphi^{445}\otimes \shift \varphi^{45}\otimes \shift \varphi^5\otimes 1)\\
=\,&0\boldsymbol{\uwave{+}}0\boldsymbol{\uwave{+}}0\boldsymbol{\uwave{+}}\underbrace{\varphi^1\varphi^2\varphi^3\varphi^{41}\dots \varphi^{43}\varphi^{441}\varphi^{442}\hat{\omega}\varphi^{443}\varphi^{444}\otimes \shift \varphi^{445}\otimes \shift \varphi^{45}\otimes \shift \varphi^5\otimes 1}_{322}\\
&-\underbrace{\varphi^1\varphi^2\varphi^3\varphi^{41}\dots \varphi^{43}\varphi^{441}\dots \varphi^{443}\hat{\omega}\varphi^{444}\otimes \shift \varphi^{445}\otimes \shift \varphi^{45}\otimes \shift \varphi^5\otimes 1}_{323}\\
&-\underbrace{\varphi^1\varphi^2\varphi^3\varphi^{41}\dots \varphi^{43}\varphi^{441}\varphi^{442}\shift \partial(\varphi^{443})\varphi^{444}\otimes \shift \varphi^{445}\otimes \shift \varphi^{45}\otimes \shift \varphi^5\otimes 1}_y\\
&\boldsymbol{\uwave{+}}0\boldsymbol{\uwave{-}}\underbrace{\varphi^1\varphi^2\varphi^3\varphi^{41}\dots \varphi^{43}\varphi^{441}\varphi^{442}\otimes \shift \varphi^{443}\otimes \varphi^{444}\varphi^{445}\otimes \shift \varphi^{45}\otimes \shift \varphi^5\otimes 1}_{177}\\
&+\underbrace{\varphi^1\varphi^2\varphi^3\varphi^{41}\dots \varphi^{43}\varphi^{441}\varphi^{442}\otimes \shift \varphi^{443}\otimes \varphi^{444}\otimes \shift (\varphi^{445}\varphi^{45})\otimes \shift \varphi^5\otimes 1}_y\\
&-\underbrace{\varphi^1\varphi^2\varphi^3\varphi^{41}\dots \varphi^{43}\varphi^{441}\varphi^{442}\otimes \shift \varphi^{443}\otimes \varphi^{444}\otimes \shift \varphi^{445}\otimes \shift (\varphi^{45}\varphi^5)\otimes 1}_g\\
&+\underbrace{\varphi^1\varphi^2\varphi^3\varphi^{41}\varphi^{42}\varphi^{43}\varphi^{441}\varphi^{442}\otimes \shift \varphi^{443}\otimes \varphi^{444}\otimes \shift \varphi^{445}\otimes \shift \varphi^{45}\otimes \varphi^5}_{176}
\end{align*}
and
\begin{align*}
&-\mu_1^{\otimes}(-\varphi^1(\shift \otimes \shift)\partial(\shift^{-1}\varphi^2)\otimes \shift \varphi^3\otimes \shift \varphi^{41}\otimes \varphi^{42}\dots \varphi^{44}\otimes \shift \varphi^{45}\otimes \shift \varphi^5\otimes 1)\\
=\,&0\boldsymbol{\uwave{+}}0\boldsymbol{\uwave{+}}\underbrace{\varphi^1(\shift \otimes \shift)\partial(\shift^{-1}\varphi^2)\varphi^3\otimes \shift \varphi^{41}\otimes \varphi^{42}\dots \varphi^{44}\otimes \shift \varphi^{45}\otimes \shift \varphi^5\otimes 1}_g\\
&-\underbrace{\varphi^1(\shift \otimes \shift)\partial(\shift^{-1}\varphi^2)\otimes \shift (\varphi^3\varphi^{41})\otimes \varphi^{42}\dots \varphi^{44}\otimes \shift \varphi^{45}\otimes \shift \varphi^5\otimes 1}_{353}\\
&+\underbrace{\varphi^1(\shift \otimes \shift)\partial(\shift^{-1}\varphi^2)\otimes \shift \varphi^3\otimes \varphi^{41}\dots \varphi^{44}\otimes \shift \varphi^{45}\otimes \shift\varphi^5\otimes 1}_{164}\\
&\boldsymbol{\uwave{+}}0\boldsymbol{\uwave{+}}0\boldsymbol{\uwave{+}}\underbrace{\varphi^1(\shift \otimes \shift)\partial(\shift^{-1}\varphi^2)\otimes \shift \varphi^3\otimes \shift \varphi^{41}\otimes \varphi^{42}\dots \varphi^{45}\otimes \shift \varphi^5\otimes 1}_{161}\\
&-\underbrace{\varphi^1(\shift \otimes \shift)\partial(\shift^{-1}\varphi^2)\otimes \shift \varphi^3\otimes \shift \varphi^{41}\otimes \varphi^{42}\dots \varphi^{44}\otimes \shift (\varphi^{45}\varphi^5)\otimes 1}_{162}\\
&+\underbrace{\varphi^1(\shift \otimes \shift)\partial(\shift^{-1}\varphi^2)\otimes \shift \varphi^3\otimes \shift \varphi^{41}\otimes \varphi^{42}\dots \varphi^{44}\otimes \shift \varphi^{45}\otimes \varphi^5}_{163}
\end{align*}
and
\begin{align*}
&-\mu_1^{\otimes}(-\varphi^1\varphi^{21}\dots \varphi^{24}\otimes \shift \varphi^{25}\otimes \shift \varphi^3\otimes \varphi^{41}\dots \varphi^{44}\otimes \shift \varphi^{45}\otimes \shift \varphi^5\otimes 1)\\
=\,&0\boldsymbol{\uwave{+}}0\boldsymbol{\uwave{-}}\underbrace{\varphi^1(\shift \otimes \shift)\partial(\shift^{-1}\varphi^2)\otimes \shift \varphi^3\otimes \varphi^{41}\dots \varphi^{44}\otimes \shift \varphi^{45}\otimes \shift \varphi^5\otimes 1}_{164}\\
&-\underbrace{\varphi^1\varphi^{21}\dots \varphi^{24}\otimes \shift (\varphi^{25}\varphi^3)\otimes \varphi^{41}\dots \varphi^{44}\otimes \shift \varphi^{45}\otimes \shift \varphi^5\otimes 1}_g\\
&+\underbrace{\varphi^1\varphi^{21}\dots \varphi^{24}\otimes \shift \varphi^{25}\otimes \varphi^3\varphi^{41}\dots \varphi^{44}\otimes \shift \varphi^{45}\otimes \shift \varphi^5\otimes 1}_{165}\\
&\boldsymbol{\uwave{+}}0\boldsymbol{\uwave{+}}0\boldsymbol{\uwave{-}}\underbrace{\varphi^1\varphi^{21}\dots \varphi^{24}\otimes \shift \varphi^{25}\otimes \shift \varphi^3\otimes (\shift \otimes \shift)\partial(\shift^{-1}\varphi^4)\otimes \shift \varphi^5\otimes 1}_{166}\\
&-\underbrace{\varphi^1\varphi^{21}\dots \varphi^{24}\otimes \shift \varphi^{25}\otimes \shift \varphi^3\otimes \varphi^{41}\dots \varphi^{44}\otimes \shift (\varphi^{45}\varphi^5)\otimes 1}_{167}\\
&+\underbrace{\varphi^1\varphi^{21}\dots \varphi^{24}\otimes \shift \varphi^{25}\otimes \shift \varphi^3\otimes \varphi^{41}\dots \varphi^{44}\otimes \shift \varphi^{45}\otimes \varphi^5}_{169}
\end{align*}
and
\begin{align*}
&-\mu_1^{\otimes}(-\varphi^1\varphi^{21}\varphi^{221}\dots \varphi^{224}\otimes \shift \varphi^{225}\otimes \shift \varphi^{23}\otimes \varphi^{24}\otimes \shift \varphi^{25}\otimes \shift \varphi^3\otimes \varphi^4\varphi^5)\\
=\,&0\boldsymbol{\uwave{+}}0\boldsymbol{\uwave{-}}\underbrace{\varphi^1\varphi^{21}(\shift \otimes \shift)\partial(\shift^{-1}\varphi^{22})\otimes \shift \varphi^{23}\otimes \varphi^{24}\otimes \shift \varphi^{25}\otimes \shift \varphi^3\otimes \varphi^4\varphi^5}_{170}\\
&-\underbrace{\varphi^1\varphi^{21}\varphi^{221}\dots \varphi^{224}\otimes \shift (\varphi^{225}\varphi^{23})\otimes \varphi^{24}\otimes \shift \varphi^{25}\otimes \shift \varphi^3\otimes \varphi^4\varphi^5}_y\\
&+\underbrace{\varphi^1\varphi^{21}\varphi^{221}\dots \varphi^{224}\otimes \shift \varphi^{225}\otimes \varphi^{23}\varphi^{24}\otimes \shift \varphi^{25}\otimes \shift \varphi^3\otimes \varphi^4\varphi^5}_{171}\\
&\boldsymbol{\uwave{+}}0\boldsymbol{\uwave{+}}\underbrace{\varphi^1\varphi^{21}\varphi^{221}\dots \varphi^{224}\otimes \shift \varphi^{225}\otimes \shift \varphi^{23}\otimes \varphi^{24}\varphi^{25}\otimes \shift \varphi^3\otimes \varphi^4\varphi^5}_{172}\\
&-\underbrace{\varphi^1\varphi^{21}\varphi^{221}\dots \varphi^{224}\otimes \shift \varphi^{225}\otimes \shift \varphi^{23}\otimes \varphi^{24}\otimes \shift (\varphi^{25}\varphi^3)\otimes \varphi^4\varphi^5}_g\\
&+\underbrace{\varphi^1\varphi^{21}\varphi^{221}\dots \varphi^{224}\otimes \shift \varphi^{225}\otimes \shift \varphi^{23}\otimes \varphi^{24}\otimes \shift \varphi^{25}\otimes \varphi^3\varphi^4\varphi^5}_{173}\boldsymbol{\uwave{+}}0
\end{align*}
and
\begin{align*}
&-\mu_1^\otimes(-1\otimes \shift \varphi^1\otimes \shift \varphi^{21}\otimes \varphi^{22}\dots \varphi^{25}\varphi^3\varphi^{41}\dots \varphi^{44}\otimes \shift \varphi^{45}\otimes \shift \varphi^5\otimes 1)\\
=\,&+\underbrace{\varphi^1\otimes \shift \varphi^{21}\otimes \varphi^{22}\dots \varphi^{25}\varphi^3\varphi^{41}\dots \varphi^{44}\otimes \shift \varphi^{45}\otimes \shift \varphi^5\otimes 1}_{180}\\
&-\underbrace{1\otimes \shift (\varphi^1\varphi^{21})\otimes \varphi^{22}\dots \varphi^{25}\varphi^3\varphi^{41}\dots \varphi^{44}\otimes \shift \varphi^{45}\otimes \shift \varphi^5\otimes 1}_g\\
&-\underbrace{1\otimes \shift \varphi^1\otimes (\shift \otimes \shift)\partial(\shift^{-1}\varphi^2)\varphi^3\varphi^{41}\dots \varphi^{44}\otimes \shift \varphi^{45}\otimes \shift \varphi^5\otimes 1}_{179}\\
&\boldsymbol{\uwave{+}}0\boldsymbol{\uwave{+}}0\boldsymbol{\uwave{+}}0\boldsymbol{\uwave{+}}0\boldsymbol{\uwave{-}}\underbrace{1\otimes \shift \varphi^1\otimes \shift \varphi^{21}\otimes \varphi^{22}\dots \varphi^{25}\varphi^3(\shift \otimes \shift)\partial(\shift^{-1}\varphi^4)\otimes \shift \varphi^5\otimes 1}_{181}\\
&-\underbrace{1\otimes \shift \varphi^1\otimes \shift \varphi^{21}\otimes \varphi^{22}\dots \varphi^{25}\varphi^3\varphi^{41}\dots \varphi^{44}\otimes \shift (\varphi^{45}\varphi^5)\otimes 1}_g\\
&+\underbrace{1\otimes \shift \varphi^1\otimes \shift \varphi^{21}\otimes \varphi^{22}\dots \varphi^{25}\varphi^3\varphi^{41}\dots \varphi^{44}\otimes \shift \varphi^{45}\otimes \varphi^5}_{182}
\end{align*}
and
\begin{align*}
&-\mu_1^{\otimes}(-\varphi^1\varphi^{21}\dots \varphi^{24}\otimes \shift \varphi^{25}\otimes \shift \varphi^3\otimes \shift \varphi^{41}\otimes \shift \varphi^{421}\otimes \varphi^{422}\dots \varphi^{425}\varphi^{43}\dots \varphi^{45}\varphi^5)\\
=\,&0\boldsymbol{\uwave{+}}0\boldsymbol{\uwave{-}}\underbrace{\varphi^1(\shift \otimes \shift)\partial(\shift^{-1}\varphi^2) \otimes \shift \varphi^3\otimes \shift \varphi^{41}\otimes \shift \varphi^{421}\otimes \varphi^{422}\dots \varphi^{425}\varphi^{43}\dots \varphi^{45}\varphi^5}_{183}\\
&+\underbrace{\varphi^1\varphi^{21}\dots \varphi^{24}\otimes \shift (\varphi^{25}\varphi^3)\otimes \shift \varphi^{41}\otimes \shift \varphi^{421}\otimes \varphi^{422}\dots \varphi^{425}\varphi^{43}\dots \varphi^{45}\varphi^5}_{gg}\\
&+\underbrace{\varphi^1\varphi^{21}\dots \varphi^{24}\otimes \shift \varphi^{25}\otimes \shift (\varphi^3\varphi^{41})\otimes \shift \varphi^{421}\otimes \varphi^{422}\dots \varphi^{425}\varphi^{43}\dots \varphi^{45}\varphi^5}_{gg}\\
&-\underbrace{\varphi^1\varphi^{21}\dots \varphi^{24}\otimes \shift \varphi^{25}\otimes \shift \varphi^3\otimes \shift (\varphi^{41}\varphi^{421})\otimes \varphi^{422}\dots \varphi^{425}\varphi^{43}\dots \varphi^{45}\varphi^5}_y\\
&-\underbrace{\varphi^1\varphi^{21}\dots \varphi^{24}\otimes \shift \varphi^{25}\otimes \shift \varphi^3\otimes \shift \varphi^{41}\otimes (\shift \otimes \shift)\partial(\shift^{-1}\varphi^{42})\varphi^{43}\dots \varphi^{45}\varphi^5}_{184}\boldsymbol{\uwave{+}}0\boldsymbol{\uwave{+}}0\boldsymbol{\uwave{+}}0
\end{align*}
and
\begin{align*}
&-\mu_1^{\otimes}(-1\otimes \shift \varphi^1\otimes \shift \varphi^{21}\otimes \varphi^{22}\dots \varphi^{25}\otimes \shift \varphi^3\otimes \shift \varphi^{41}\otimes \varphi^{42}\dots \varphi^{45}\varphi^5)\\
=\,&+\underbrace{\varphi^1\otimes \shift \varphi^{21}\otimes \varphi^{22}\dots \varphi^{25}\otimes \shift \varphi^3\otimes \shift \varphi^{41}\otimes \varphi^{42}\dots \varphi^{45}\varphi^5}_{187}\\
&-\underbrace{1\otimes \shift (\varphi^1\varphi^{21})\otimes \varphi^{22}\dots \varphi^{25}\otimes \shift \varphi^3\otimes \shift \varphi^{41}\otimes \varphi^{42}\dots \varphi^{45}\varphi^5}_{186}\\
&-\underbrace{1\otimes \shift \varphi^1\otimes (\shift \otimes \shift)\partial(\shift^{-1}\varphi^2)\otimes \shift \varphi^3\otimes \shift \varphi^{41}\otimes \varphi^{42}\dots \varphi^{45}\varphi^5}_{185}\\
&\boldsymbol{\uwave{+}}0\boldsymbol{\uwave{+}}0\boldsymbol{\uwave{+}}\underbrace{1\otimes \shift \varphi^1\otimes \shift \varphi^{21}\otimes \varphi^{22}\dots \varphi^{25}\varphi^3\otimes \shift \varphi^{41}\otimes \varphi^{42}\dots \varphi^{45}\varphi^5}_{193}\\
&-\underbrace{1\otimes \shift \varphi^1\otimes \shift \varphi^{21}\otimes \varphi^{22}\dots \varphi^{25}\otimes \shift (\varphi^3\varphi^{41})\otimes \varphi^{42}\dots \varphi^{45}\varphi^5}_g\\
&-\underbrace{1\otimes \shift \varphi^1\otimes \shift \varphi^{21}\otimes \varphi^{22}\dots \varphi^{25}\otimes \shift \varphi^3\otimes (\shift \otimes \shift)\partial(\shift^{-1}\varphi^4)\varphi^5}_{192}\boldsymbol{\uwave{+}}0\boldsymbol{\uwave{+}}0
\end{align*}
and
\begin{align*}
&-\mu_1^{\otimes}(-1\otimes \shift \varphi^1\otimes \shift \varphi^{21}\otimes \varphi^{22}\dots \varphi^{24}\otimes \shift \varphi^{25}\otimes \shift \varphi^3\otimes (\shift \otimes \shift)\partial(\shift^{-1}\varphi^4)\varphi^5)\\
=\,&+\underbrace{\varphi^1\otimes \shift \varphi^{21}\otimes \varphi^{22}\dots \varphi^{24}\otimes \shift \varphi^{25}\otimes \shift \varphi^3\otimes (\shift \otimes \shift)\partial(\shift^{-1}\varphi^4)\varphi^5}_{190}\\
&-\underbrace{1\otimes \shift (\varphi^1\varphi^{21})\otimes \varphi^{22}\dots \varphi^{24}\otimes \shift \varphi^{25}\otimes \shift \varphi^3\otimes (\shift \otimes \shift)\partial(\shift^{-1}\varphi^4)\varphi^5}_{191}\\
&+\underbrace{1\otimes \shift \varphi^1\otimes \varphi^{21}\dots \varphi^{24}\otimes \shift \varphi^{25}\otimes \shift \varphi^3\otimes (\shift \otimes \shift)\partial(\shift^{-1}\varphi^4)\varphi^5}_{189}\\
&\boldsymbol{\uwave{+}}0\boldsymbol{\uwave{+}}0\boldsymbol{\uwave{+}}\underbrace{1\otimes \shift \varphi^1\otimes \shift \varphi^{21}\otimes \varphi^{22}\dots \varphi^{25}\otimes \shift \varphi^3\otimes (\shift \otimes \shift)\partial(\shift^{-1}\varphi^4)\varphi^5}_{192}\\
&-\underbrace{1\otimes \shift \varphi^1\otimes \shift \varphi^{21}\otimes \varphi^{22}\dots \varphi^{24}\otimes \shift (\varphi^{25}\varphi^3)\otimes (\shift \otimes \shift)\partial(\shift^{-1}\varphi^4)\varphi^5}_{188}\\
&+\underbrace{1\otimes \shift \varphi^1\otimes \shift \varphi^{21}\otimes \varphi^{22}\dots \varphi^{24}\otimes \shift \varphi^{25}\otimes \varphi^3(\shift \otimes \shift)\partial(\shift^{-1}\varphi^4)\varphi^5}_g\boldsymbol{\uwave{+}}0\boldsymbol{\uwave{+}}0
\end{align*}
and
\begin{align*}
&-\mu_1^{\otimes}(-\varphi^1\varphi^2\otimes \shift \varphi^3\otimes \shift \varphi^{41}\otimes (\shift \otimes \shift)\partial(\shift^{-1}\varphi^{42})\varphi^{43}\varphi^{44}\otimes \shift \varphi^{45}\otimes \shift \varphi^5\otimes 1)\\
=\,&0\boldsymbol{\uwave{+}}\underbrace{\varphi^1\varphi^2\varphi^3\otimes \shift \varphi^{41}\otimes (\shift \otimes \shift)\partial(\shift^{-1}\varphi^{42})\varphi^{43}\varphi^{44}\otimes \shift \varphi^{45}\otimes \shift \varphi^5\otimes 1}_{242}\\
&-\underbrace{\varphi^1\varphi^2\otimes \shift (\varphi^3\varphi^{41})\otimes (\shift \otimes \shift)\partial(\shift^{-1}\varphi^{42})\varphi^{43}\varphi^{44}\otimes \shift \varphi^{45}\otimes \shift \varphi^5\otimes 1}_{233}\\
&+\underbrace{\varphi^1\varphi^2\otimes \shift \varphi^3\otimes \varphi^{41}(\shift \otimes \shift)\partial(\shift^{-1}\varphi^{42})\varphi^{43}\varphi^{44}\otimes \shift \varphi^{45}\otimes \shift \varphi^5\otimes 1}_y\\
&\boldsymbol{\uwave{+}}0\boldsymbol{\uwave{+}}0\boldsymbol{\uwave{+}}0\boldsymbol{\uwave{+}}\underbrace{\varphi^1\varphi^2\otimes \shift \varphi^3\otimes \shift \varphi^{41}\otimes (\shift \otimes \shift)\partial(\shift^{-1}\varphi^{42})\varphi^{43}\dots \varphi^{45}\otimes \shift \varphi^5\otimes 1}_{234}\\
&-\underbrace{\varphi^1\varphi^2\otimes \shift \varphi^3\otimes \shift \varphi^{41}\otimes (\shift \otimes \shift)\partial(\shift^{-1}\varphi^{42})\varphi^{43}\varphi^{44}\otimes \shift (\varphi^{45}\varphi^5)\otimes 1}_{231}\\
&+\underbrace{\varphi^1\varphi^2\otimes \shift \varphi^3\otimes \shift \varphi^{41}\otimes (\shift \otimes \shift)\partial(\shift^{-1}\varphi^{42})\varphi^{43}\varphi^{44}\otimes \shift \varphi^{45}\otimes \varphi^5}_{232}
\end{align*}
and
\begin{align*}
&-\mu_1^{\otimes}(-\varphi^1\otimes \shift \varphi^{21}\otimes \shift \varphi^{221}\otimes \varphi^{222}\dots \varphi^{225}\varphi^{23}\varphi^{24}\otimes \shift \varphi^{25}\otimes \shift \varphi^3\otimes \varphi^4\varphi^5)\\
=\,&+\underbrace{\varphi^1\varphi^{21}\otimes \shift \varphi^{221}\otimes \varphi^{222}\dots \varphi^{225}\varphi^{23}\varphi^{24}\otimes \shift \varphi^{25}\otimes \shift \varphi^3\otimes \varphi^4\varphi^5}_{240}\\
&-\underbrace{\varphi^1\otimes \shift (\varphi^{21}\varphi^{221})\otimes \varphi^{222}\dots \varphi^{225}\varphi^{23}\varphi^{24}\otimes \shift \varphi^{25}\otimes \shift \varphi^3\otimes \varphi^4\varphi^5}_y\\
&-\underbrace{\varphi^1\otimes \shift \varphi^{21}\otimes (\shift \otimes \shift)\partial(\shift^{-1}\varphi^{22})\varphi^{23}\varphi^{24}\otimes \shift \varphi^{25}\otimes \shift \varphi^3\otimes \varphi^4\varphi^5}_{237}\\
&\boldsymbol{\uwave{+}}0\boldsymbol{\uwave{+}}0\boldsymbol{\uwave{+}}0\boldsymbol{\uwave{+}}\underbrace{\varphi^1\otimes \shift \varphi^{21}\otimes \shift \varphi^{221}\otimes \varphi^{222}\dots \varphi^{225}\varphi^{23}\dots \varphi^{25}\otimes \shift \varphi^3\otimes \varphi^4\varphi^5}_{239}\\
&-\underbrace{\varphi^1\otimes \shift \varphi^{21}\otimes \shift \varphi^{221}\otimes \varphi^{222}\dots \varphi^{225}\varphi^{23}\varphi^{24}\otimes \shift (\varphi^{25}\varphi^3)\otimes \varphi^4\varphi^5}_{238}\\
&+\underbrace{\varphi^1\otimes \shift \varphi^{21}\otimes \shift \varphi^{221}\otimes \varphi^{222}\dots \varphi^{225}\varphi^{23}\varphi^{24}\otimes \shift \varphi^{25}\otimes \varphi^3\varphi^4\varphi^5}_{241}\boldsymbol{\uwave{+}}0
\end{align*}
and
\begin{align*}
&-\mu_1^{\otimes}(-\varphi^1\varphi^2\varphi^3\otimes \shift \varphi^{41}\otimes \shift \varphi^{421}\otimes \varphi^{422}\dots \varphi^{425}\varphi^{43}\varphi^{44}\otimes \shift \varphi^{45}\otimes \shift \varphi^5\otimes 1)\\
=\,&0\boldsymbol{\uwave{+}}\underbrace{\varphi^1\varphi^2\varphi^3\varphi^{41}\otimes \shift \varphi^{421}\otimes \varphi^{422}\dots \varphi^{425}\varphi^{43}\varphi^{44}\otimes \shift \varphi^{45}\otimes \shift \varphi^5\otimes 1}_{245}\\
&-\underbrace{\varphi^1\varphi^2\varphi^3\otimes \shift (\varphi^{41}\varphi^{421})\otimes \varphi^{422}\dots \varphi^{425}\varphi^{43}\varphi^{44}\otimes \shift \varphi^{45}\otimes \shift \varphi^5\otimes 1}_y\\
&-\underbrace{\varphi^1\varphi^2\varphi^3\otimes \shift \varphi^{41}\otimes (\shift \otimes \shift)\partial(\shift^{-1}\varphi^{42})\varphi^{43}\varphi^{44}\otimes \shift \varphi^{45}\otimes \shift \varphi^5\otimes 1}_{242}\\
&\boldsymbol{\uwave{+}}0\boldsymbol{\uwave{+}}0\boldsymbol{\uwave{+}}0\boldsymbol{\uwave{+}}\underbrace{\varphi^1\varphi^2\varphi^3\otimes \shift \varphi^{41}\otimes \shift \varphi^{421}\otimes \varphi^{422}\dots \varphi^{425}\varphi^{43}\dots \varphi^{45}\otimes \shift \varphi^5\otimes 1}_{246}\\
&-\underbrace{\varphi^1\varphi^2\varphi^3\otimes \shift \varphi^{41}\otimes \shift \varphi^{421}\otimes \varphi^{422}\dots \varphi^{425}\varphi^{43}\varphi^{44}\otimes \shift (\varphi^{45}\varphi^5)\otimes 1}_{243}\\
&+\underbrace{\varphi^1\varphi^2\varphi^3\otimes \shift \varphi^{41}\otimes \shift \varphi^{421}\otimes \varphi^{422}\dots \varphi^{425}\varphi^{43}\varphi^{44}\otimes \shift \varphi^{45}\otimes \varphi^5}_{244}
\end{align*}
and
\begin{align*}
&-\mu_1^{\otimes}(-1\otimes \shift \varphi^1\otimes \shift \varphi^{21}\otimes \varphi^{221}\dots \varphi^{224}\otimes \shift \varphi^{225}\otimes \shift \varphi^{23}\otimes \varphi^{24}\varphi^{25}\varphi^3\varphi^4\varphi^5)\\
=\,&+\underbrace{\varphi^1\otimes \shift \varphi^{21}\otimes \varphi^{221}\dots \varphi^{224}\otimes \shift \varphi^{225}\otimes \shift \varphi^{23}\otimes \varphi^{24}\varphi^{25}\varphi^3\varphi^4\varphi^5}_{251}\\
&-\underbrace{1\otimes \shift (\varphi^1\varphi^{21})\otimes \varphi^{221}\dots \varphi^{224}\otimes \shift \varphi^{225}\otimes \shift \varphi^{23}\otimes \varphi^{24}\varphi^{25}\varphi^3\varphi^4\varphi^5}_g\\
&+\underbrace{1\otimes \shift \varphi^1\otimes \varphi^{21}\varphi^{221}\dots \varphi^{224}\otimes \shift \varphi^{225}\otimes \shift \varphi^{23}\otimes \varphi^{24}\varphi^{25}\varphi^3\varphi^4\varphi^5}_{250}\\
&\boldsymbol{\uwave{+}}0\boldsymbol{\uwave{+}}0\boldsymbol{\uwave{-}}\underbrace{1\otimes \shift \varphi^1\otimes \shift \varphi^{21}\otimes (\shift \otimes \shift)\partial(\shift^{-1}\varphi^{22})\otimes \shift \varphi^{23}\otimes \varphi^{24}\varphi^{25}\varphi^3\varphi^4\varphi^5}_{247}\\
&-\underbrace{1\otimes \shift \varphi^1\otimes \shift \varphi^{21}\otimes \varphi^{221}\dots \varphi^{224}\otimes \shift (\varphi^{225}\varphi^{23})\otimes \varphi^{24}\varphi^{25}\varphi^3\varphi^4\varphi^5}_{248}\\
&+\underbrace{1\otimes \shift \varphi^1\otimes \shift \varphi^{21}\otimes \varphi^{221}\dots \varphi^{224}\otimes \shift \varphi^{225}\otimes \varphi^{23}\dots \varphi^{25}\varphi^3\varphi^4\varphi^5}_{249}\boldsymbol{\uwave{+}}0\boldsymbol{\uwave{+}}0
\end{align*}
and
\begin{align*}
&-\mu_1^{\otimes}(-1\otimes \shift \varphi^1\otimes \shift \varphi^{21}\otimes \varphi^{22}\varphi^{23}(\shift \otimes \shift)\partial(\shift^{-1}\varphi^{24})\otimes \shift \varphi^{25}\otimes \shift \varphi^3\otimes \varphi^4\varphi^5)\\
=\,&+\underbrace{\varphi^1\otimes \shift \varphi^{21}\otimes \varphi^{22}\varphi^{23}(\shift \otimes \shift)\partial(\shift^{-1}\varphi^{24})\otimes \shift \varphi^{25}\otimes \shift\varphi^3\otimes \varphi^4\varphi^5}_{254}\\
&-\underbrace{1\otimes \shift (\varphi^1\varphi^{21})\otimes \varphi^{22}\varphi^{23}(\shift \otimes \shift)\partial(\shift^{-1}\varphi^{24})\otimes \shift \varphi^{25}\otimes \shift \varphi^3\otimes \varphi^4\varphi^5}_{252}\\
&+\underbrace{1\otimes \shift \varphi^1\otimes \varphi^{21}\varphi^{22}\varphi^{23}(\shift \otimes \shift)\partial(\shift^{-1}\varphi^{24})\otimes \shift \varphi^{25}\otimes \shift \varphi^3\otimes \varphi^4\varphi^5}_{253}\\
&\boldsymbol{\uwave{+}}0\boldsymbol{\uwave{+}}0\boldsymbol{\uwave{+}}0\boldsymbol{\uwave{+}}\underbrace{1\otimes \shift \varphi^1\otimes \shift \varphi^{21}\otimes \varphi^{22}\varphi^{23}(\shift \otimes \shift)\partial(\shift^{-1}\varphi^{24})\varphi^{25}\otimes \shift \varphi^3\otimes \varphi^4\varphi^5}_y\\
&-\underbrace{1\otimes \shift \varphi^1\otimes \shift \varphi^{21}\otimes \varphi^{22}\varphi^{23}(\shift \otimes \shift)\partial(\shift^{-1}\varphi^{24})\otimes \shift (\varphi^{25}\varphi^3)\otimes \varphi^4\varphi^5}_{255}\\
&+\underbrace{1\otimes \shift \varphi^1\otimes \shift \varphi^{21}\otimes \varphi^{22}\varphi^{23}(\shift \otimes \shift)\partial(\shift^{-1}\varphi^{24})\otimes \shift \varphi^{25}\otimes \varphi^3\varphi^4\varphi^5}_{256}\boldsymbol{\uwave{+}}0
\end{align*}
and
\begin{align*}
&-\mu_1^{\otimes}(-1\otimes \shift \varphi^1\otimes \varphi^{21}\varphi^{22}\otimes \shift \varphi^{23}\otimes \varphi^{24}\otimes \shift \varphi^{25}\otimes \shift \varphi^3\otimes \varphi^4\varphi^5)\\
=\,&+\underbrace{\hat{\omega}\varphi^1\varphi^{21}\varphi^{22}\otimes \shift \varphi^{23}\otimes \varphi^{24}\otimes \shift \varphi^{25}\otimes \shift \varphi^3\otimes \varphi^4\varphi^5}_{128}\\
&-\underbrace{\varphi^1\hat{\omega}\varphi^{21}\varphi^{22}\otimes \shift \varphi^{23}\otimes \varphi^{24}\otimes \shift \varphi^{25}\otimes \shift \varphi^3\otimes \varphi^4\varphi^5}_{130}\\
&-\underbrace{\shift \partial(\varphi^1)\varphi^{21}\varphi^{22}\otimes \shift \varphi^{23}\otimes \varphi^{24}\otimes \shift \varphi^{25}\otimes \shift \varphi^3\otimes \varphi^4\varphi^5}_g\\
&\boldsymbol{\uwave{+}}0\boldsymbol{\uwave{-}}\underbrace{1\otimes \shift \varphi^1\otimes \varphi^{21}\varphi^{22}\hat{\omega}\varphi^{23}\varphi^{24}\otimes \shift \varphi^{25}\otimes \shift \varphi^3\otimes \varphi^4\varphi^5}_{126}\\
&+\underbrace{1\otimes \shift \varphi^1\otimes \varphi^{21}\varphi^{22}\varphi^{23}\hat{\omega}\varphi^{24}\otimes \shift \varphi^{25}\otimes \shift \varphi^3\otimes \varphi^4\varphi^5}_{127}\\
&+\underbrace{1\otimes \shift \varphi^1\otimes \varphi^{21}\varphi^{22}\shift \partial(\varphi^{23})\varphi^{24}\otimes \shift \varphi^{25}\otimes \shift \varphi^3\otimes \varphi^4\varphi^5}_y\\
&\boldsymbol{\uwave{+}}0\boldsymbol{\uwave{+}}\underbrace{1\otimes \shift \varphi^1\otimes \varphi^{21}\varphi^{22}\otimes \shift \varphi^{23}\otimes \varphi^{24}\varphi^{25}\otimes \shift \varphi^3\otimes \varphi^4\varphi^5}_{145}\\
&-\underbrace{1\otimes \shift \varphi^1\otimes \varphi^{21}\varphi^{22}\otimes \shift \varphi^{23}\otimes \varphi^{24}\otimes \shift (\varphi^{25}\varphi^3)\otimes \varphi^4\varphi^5}_{224}\\
&+\underbrace{1\otimes \shift \varphi^1\otimes \varphi^{21}\varphi^{22}\otimes \shift \varphi^{23}\otimes \varphi^{24}\otimes \shift \varphi^{25}\otimes \varphi^3\varphi^4\varphi^5}_{222}\boldsymbol{\uwave{+}}0
\end{align*}
and
\begin{align*}
&-\mu_1^{\otimes}(1\otimes \shift \varphi^1\otimes \varphi^2\varphi^3\varphi^{41}\varphi^{42}\otimes \shift \varphi^{43}\otimes \varphi^{44}\otimes \shift \varphi^{45}\otimes \shift \varphi^5\otimes 1)\\
=\,&-\underbrace{\hat{\omega}\varphi^1\varphi^2\varphi^3\varphi^{41}\varphi^{42}\otimes \shift \varphi^{43}\otimes \varphi^{44}\otimes \shift \varphi^{45}\otimes \shift \varphi^5\otimes 1}_{82}\\
&+\underbrace{\varphi^1\hat{\omega}\varphi^2\varphi^3\varphi^{41}\varphi^{42}\otimes \shift \varphi^{43}\otimes \varphi^{44}\otimes \shift \varphi^{45}\otimes \shift \varphi^5\otimes 1}_{88}\\
&+\underbrace{\shift \partial(\varphi^1)\varphi^2\varphi^3\varphi^{41}\varphi^{42}\otimes \shift \varphi^{43}\otimes \varphi^{44}\otimes \shift \varphi^{45}\otimes \shift \varphi^5\otimes 1}_g\\
&\boldsymbol{\uwave{+}}0\boldsymbol{\uwave{+}}0\boldsymbol{\uwave{+}}\underbrace{1\otimes \shift \varphi^1\otimes \varphi^2\varphi^3\varphi^{41}\varphi^{42}\hat{\omega}\varphi^{43}\varphi^{44}\otimes \shift \varphi^{45}\otimes \shift \varphi^5\otimes 1}_{87}\\
&-\underbrace{1\otimes \shift \varphi^1\otimes \varphi^2\varphi^3\varphi^{41}\varphi^{42}\varphi^{43}\hat{\omega}\varphi^{44}\otimes \shift \varphi^{45}\otimes \shift \varphi^5\otimes 1}_{89}\\
&-\underbrace{1\otimes \shift \varphi^1\otimes \varphi^2\varphi^3\varphi^{41}\varphi^{42}\shift \partial(\varphi^{43})\varphi^{44}\otimes \shift \varphi^{45}\otimes \shift \varphi^5\otimes 1}_y\\
&\boldsymbol{\uwave{+}}0\boldsymbol{\uwave{-}}\underbrace{1\otimes \shift \varphi^1\otimes \varphi^2\varphi^3\varphi^{41}\varphi^{42}\otimes \shift \varphi^{43}\otimes \varphi^{44}\varphi^{45}\otimes \shift \varphi^5\otimes 1}_{212}\\
&+\underbrace{1\otimes \shift \varphi^1\otimes \varphi^2\varphi^3\varphi^{41}\varphi^{42}\otimes \shift \varphi^{43}\otimes \varphi^{44}\otimes \shift (\varphi^{45}\varphi^5)\otimes 1}_y\\
&-\underbrace{1\otimes \shift \varphi^1\otimes \varphi^2\varphi^3\varphi^{41}\varphi^{42}\otimes \shift \varphi^{43}\otimes \varphi^{44}\otimes \shift \varphi^{45}\otimes \varphi^5}_{211}
\end{align*}
and
\begin{align*}
&-\mu_1^{\otimes}(1\otimes \shift \varphi^1\otimes \varphi^2\otimes \shift \varphi^3\otimes \shift \varphi^{41}\otimes \varphi^{42}\dots \varphi^{45}\otimes \shift \varphi^5\otimes 1)\\
=\,&-\underbrace{\hat{\omega}\varphi^1\varphi^2\otimes \shift \varphi^3\otimes \shift \varphi^{41}\otimes \varphi^{42}\dots \varphi^{45}\otimes \shift \varphi^5\otimes 1}_{75}\\
&+\underbrace{\varphi^1\hat{\omega}\varphi^2\otimes \shift \varphi^3\otimes \shift \varphi^{41}\otimes \varphi^{42}\dots \varphi^{45}\otimes \shift \varphi^5\otimes 1}_{76}\\
&+\underbrace{\shift \partial(\varphi^1)\varphi^2\otimes \shift \varphi^3\otimes \shift \varphi^{41}\otimes \varphi^{42}\dots \varphi^{45}\otimes \shift \varphi^5\otimes 1}_g\\
&\boldsymbol{\uwave{+}}0\boldsymbol{\uwave{+}}\underbrace{1\otimes \shift \varphi^1\otimes \varphi^2\varphi^3\otimes \shift \varphi^{41}\otimes \varphi^{42}\dots \varphi^{45}\otimes \shift \varphi^5\otimes 1}_{210}\\
&-\underbrace{1\otimes \shift \varphi^1\otimes \varphi^2\otimes \shift (\varphi^3\varphi^{41})\otimes \varphi^{42}\dots \varphi^{45}\otimes \shift \varphi^5\otimes 1}_{195}\\
&-\underbrace{1\otimes \shift \varphi^1\otimes \varphi^2\otimes \shift \varphi^3\otimes (\shift \otimes \shift)\partial(\shift^{-1}\varphi^4)\otimes \shift \varphi^5\otimes 1}_{148}\\
&\boldsymbol{\uwave{+}}0\boldsymbol{\uwave{+}}0\boldsymbol{\uwave{+}}\underbrace{1\otimes \shift \varphi^1\otimes \varphi^2\otimes \shift \varphi^3\otimes \shift \varphi^{41}\otimes \varphi^{42}\dots \varphi^{45}\hat{\omega}\varphi^5}_{90}\\
&-\underbrace{1\otimes \shift \varphi^1\otimes \varphi^2\otimes \shift \varphi^3\otimes \shift \varphi^{41}\otimes \varphi^{42}\dots \varphi^{45}\varphi^5\hat{\omega}}_{78}\\
&-\underbrace{1\otimes \shift \varphi^1\otimes \varphi^2\otimes \shift \varphi^3\otimes \shift \varphi^{41}\otimes \varphi^{42}\dots \varphi^{45}\shift \partial(\varphi^5)}_g
\end{align*}
and
\begin{align*}
&-\mu_1^{\otimes}(1\otimes \shift (\varphi^1\varphi^{21})\otimes \varphi^{22}\dots \varphi^{24}\otimes \shift \varphi^{25}\otimes \shift \varphi^3\otimes \shift \varphi^{41}\otimes \varphi^{42}\dots \varphi^{45}\varphi^5)\\
=\,&-\underbrace{\hat{\omega}\varphi^1\varphi^{21}\dots \varphi^{24}\otimes \shift \varphi^{25}\otimes \shift \varphi^3\otimes \shift \varphi^{41}\otimes \varphi^{42}\dots \varphi^{45}\varphi^5}_{123}\\
&+\underbrace{\varphi^1\varphi^{21}\hat{\omega}\varphi^{22}\dots \varphi^{24}\otimes \shift \varphi^{25}\otimes \shift \varphi^3\otimes \shift \varphi^{41}\otimes \varphi^{42}\dots \varphi^{45}\varphi^5}_{124}\\
&+\underbrace{\shift \partial(\varphi^1)\varphi^{21}\dots \varphi^{24}\otimes \shift \varphi^{25}\otimes \shift \varphi^3\otimes \shift \varphi^{41}\otimes \varphi^{42}\dots \varphi^{45}\varphi^5}_{gg}\\
&\boldsymbol{\uwave{+}}0\boldsymbol{\uwave{+}}0\boldsymbol{\uwave{+}}\underbrace{\varphi^1\shift \partial(\varphi^{21})\varphi^{22}\dots \varphi^{24}\otimes \shift \varphi^{25}\otimes \shift \varphi^3\otimes \shift \varphi^{41}\otimes \varphi^{42}\dots \varphi^{45}\varphi^5}_y\\
&+\underbrace{1\otimes \shift (\varphi^1\varphi^{21})\otimes \varphi^{22}\dots \varphi^{25}\otimes \shift \varphi^3\otimes \shift \varphi^{41}\otimes \varphi^{42}\dots \varphi^{45}\varphi^5}_{186}\\
&-\underbrace{1\otimes \shift (\varphi^1\varphi^{21})\otimes \varphi^{22}\dots \varphi^{24}\otimes \shift (\varphi^{25}\varphi^3)\otimes \shift \varphi^{41}\otimes \varphi^{42}\dots \varphi^{45}\varphi^5}_g\\
&+\underbrace{1\otimes \shift (\varphi^1\varphi^{21})\otimes \varphi^{22}\dots \varphi^{24}\otimes \shift \varphi^{25}\otimes \shift (\varphi^3\varphi^{41})\otimes \varphi^{42}\dots \varphi^{45}\varphi^5}_g\\
&+\underbrace{1\otimes \shift (\varphi^1\varphi^{21})\otimes \varphi^{22}\dots \varphi^{24}\otimes \shift \varphi^{25}\otimes \shift \varphi^3\otimes (\shift \otimes \shift)\partial(\shift^{-1}\varphi^4)\varphi^5}_{191}\boldsymbol{\uwave{+}}0\boldsymbol{\uwave{+}}0
\end{align*}
and
\begin{align*}
&-\mu_1^{\otimes}(-\varphi^1\varphi^2\otimes \shift \varphi^3\otimes \shift \varphi^{41}\otimes \varphi^{42}\otimes \shift \varphi^{43}\otimes \varphi^{44}\varphi^{45}\otimes \shift \varphi^5\otimes 1)\\
=\,&0\boldsymbol{\uwave{+}}\underbrace{\varphi^1\varphi^2\varphi^3\otimes \shift \varphi^{41}\otimes \varphi^{42}\otimes \shift \varphi^{43}\otimes \varphi^{44}\varphi^{45}\otimes \shift \varphi^5\otimes 1}_{215}\\
&-\underbrace{\varphi^1\varphi^2\otimes \shift (\varphi^3\varphi^{41})\otimes \varphi^{42}\otimes \shift \varphi^{43}\otimes \varphi^{44}\varphi^{45}\otimes \shift \varphi^5\otimes 1}_{209}\\
&+\underbrace{\varphi^1\varphi^2\otimes \shift \varphi^3\otimes \varphi^{41}\varphi^{42}\otimes \shift \varphi^{43}\otimes \varphi^{44}\varphi^{45}\otimes \shift \varphi^5\otimes 1}_{223}\\
&\boldsymbol{\uwave{+}}0\boldsymbol{\uwave{+}}\underbrace{\varphi^1\varphi^2\otimes \shift \varphi^3\otimes \shift \varphi^{41}\otimes \varphi^{42}\hat{\omega}\varphi^{43}\varphi^{44}\varphi^{45}\otimes \shift \varphi^5\otimes 1}_{113}\\
&-\underbrace{\varphi^1\varphi^2\otimes \shift \varphi^3\otimes \shift \varphi^{41}\otimes \varphi^{42}\varphi^{43}\hat{\omega}\varphi^{44}\varphi^{45}\otimes \shift \varphi^5\otimes 1}_{114}\\
&-\underbrace{\varphi^1\varphi^2\otimes \shift \varphi^3\otimes \shift \varphi^{41}\otimes \varphi^{42}\shift \partial(\varphi^{43})\varphi^{44}\varphi^{45}\otimes \shift \varphi^5\otimes 1}_{g}\\
&\boldsymbol{\uwave{+}}0\boldsymbol{\uwave{-}}\underbrace{\varphi^1\varphi^2\otimes \shift \varphi^3\otimes \shift \varphi^{41}\otimes \varphi^{42}\otimes \shift \varphi^{43}\otimes \varphi^{44}\varphi^{45}\hat{\omega}\varphi^5}_{112}\\
&+\underbrace{\varphi^1\varphi^2\otimes \shift \varphi^3\otimes \shift \varphi^{41}\otimes \varphi^{42}\otimes \shift \varphi^{43}\otimes \varphi^{44}\varphi^{45}\varphi^5\hat{\omega}}_{110}\\
&+\underbrace{\varphi^1\varphi^2\otimes \shift \varphi^3\otimes \shift \varphi^{41}\otimes \varphi^{42}\otimes \shift \varphi^{43}\otimes \varphi^{44}\varphi^{45}\shift \partial(\varphi^5)}_g
\end{align*}
and
\begin{align*}
&-\mu_1^{\otimes}(-\varphi^1\varphi^{21}\varphi^{22}\otimes \shift \varphi^{23}\otimes \varphi^{24}\otimes \shift \varphi^{25}\otimes \shift \varphi^3\otimes \shift \varphi^{41}\otimes \varphi^{42}\dots \varphi^{45}\varphi^5)\\
=\,&0\boldsymbol{\uwave{+}}\underbrace{\varphi^1\varphi^{21}\varphi^{22}\hat{\omega}\varphi^{23}\varphi^{24}\otimes \shift \varphi^{25}\otimes \shift \varphi^3\otimes \shift \varphi^{41}\otimes \varphi^{42}\dots \varphi^{45}\varphi^5}_{133}\\
&-\underbrace{\varphi^1\varphi^{21}\varphi^{22}\varphi^{23}\hat{\omega}\varphi^{24}\otimes \shift \varphi^{25}\otimes \shift \varphi^3\otimes \shift \varphi^{41}\otimes \varphi^{42}\dots \varphi^{45}\varphi^5}_{134}\\
&-\underbrace{\varphi^1\varphi^{21}\varphi^{22}\shift \partial(\varphi^{23})\varphi^{24}\otimes \shift \varphi^{25}\otimes \shift \varphi^3\otimes \shift \varphi^{41}\otimes \varphi^{42}\dots \varphi^{45}\varphi^5}_y\\
&\boldsymbol{\uwave{+}}0\boldsymbol{\uwave{-}}\underbrace{\varphi^1\varphi^{21}\varphi^{22}\otimes \shift \varphi^{23}\otimes \varphi^{24}\varphi^{25}\otimes \shift \varphi^3\otimes \shift \varphi^{41}\otimes \varphi^{42}\dots \varphi^{45}\varphi^5}_{216}\\
&+\underbrace{\varphi^1\varphi^{21}\varphi^{22}\otimes \shift \varphi^{23}\otimes \varphi^{24}\otimes \shift (\varphi^{25}\varphi^3)\otimes \shift \varphi^{41}\otimes \varphi^{42}\dots \varphi^{45}\varphi^5}_y\\
&-\underbrace{\varphi^1\varphi^{21}\varphi^{22}\otimes \shift \varphi^{23}\otimes \varphi^{24}\otimes \shift \varphi^{25}\otimes \shift (\varphi^3\varphi^{41})\otimes \varphi^{42}\dots \varphi^{45}\varphi^5}_g\\
&-\underbrace{\varphi^1\varphi^{21}\varphi^{22}\otimes \shift \varphi^{23}\otimes \varphi^{24}\otimes \shift \varphi^{25}\otimes \shift \varphi^3\otimes (\shift \otimes \shift)\partial(\shift^{-1}\varphi^4)\varphi^5}_{155}\boldsymbol{\uwave{+}}0\boldsymbol{\uwave{+}}0
\end{align*}
and
\begin{align*}
&-\mu_1^{\otimes}(-\varphi^1\varphi^{21}\dots \varphi^{24}\otimes \shift \varphi^{25}\otimes \shift \varphi^3\otimes \shift \varphi^{41}\otimes \varphi^{42}\otimes \shift \varphi^{43}\otimes \varphi^{44}\varphi^{45}\varphi^5)\\
=\,&0\boldsymbol{\uwave{+}}0\boldsymbol{\uwave{-}}\underbrace{\varphi^1(\shift \otimes \shift)\partial(\shift^{-1}\varphi^2)\otimes \shift \varphi^3\otimes \shift \varphi^{41}\otimes \varphi^{42}\otimes \shift \varphi^{43}\otimes \varphi^{44}\varphi^{45}\varphi^5}_{345}\\
&-\underbrace{\varphi^1\varphi^{21}\dots \varphi^{24}\otimes \shift (\varphi^{25}\varphi^3)\otimes \shift \varphi^{41}\otimes \varphi^{42}\otimes \shift \varphi^{43}\otimes \varphi^{44}\varphi^{45}\varphi^5}_g\\
&+\underbrace{\varphi^1\varphi^{21}\dots \varphi^{24}\otimes \shift \varphi^{25}\otimes \shift (\varphi^3\varphi^{41})\otimes \varphi^{42}\otimes \shift \varphi^{43}\otimes \varphi^{44}\varphi^{45}\varphi^5}_g\\
&-\underbrace{\varphi^1\varphi^{21}\dots \varphi^{24}\otimes \shift \varphi^{25}\otimes \shift \varphi^3\otimes \varphi^{41}\varphi^{42}\otimes \shift \varphi^{43}\otimes \varphi^{44}\varphi^{45}\varphi^5}_{168}\\
&\boldsymbol{\uwave{+}}0\boldsymbol{\uwave{-}}\underbrace{\varphi^1\varphi^{21}\dots \varphi^{24}\otimes \shift \varphi^{25}\otimes \shift \varphi^3\otimes \shift \varphi^{41}\otimes \varphi^{42}\hat{\omega}\varphi^{43}\varphi^{44}\varphi^{45}\varphi^5}_{140}\\
&+\underbrace{\varphi^1\varphi^{21}\dots \varphi^{24}\otimes \shift \varphi^{25}\otimes \shift \varphi^3\otimes \shift \varphi^{41}\otimes \varphi^{42}\varphi^{43}\hat{\omega}\varphi^{44}\varphi^{45}\varphi^5}_{139}\\
&+\underbrace{\varphi^1\varphi^{21}\dots \varphi^{24}\otimes \shift \varphi^{25}\otimes \shift \varphi^3\otimes \shift \varphi^{41}\otimes \varphi^{42}\shift \partial(\varphi^{43})\varphi^{44}\varphi^{45}\varphi^5}_y\boldsymbol{\uwave{+}}0
\end{align*}
and
\begin{align*}
&-\mu_1^{\otimes}(\varphi^1\varphi^{21}\dots \varphi^{24}\otimes \shift \varphi^{25}\otimes \shift \varphi^3\otimes \shift \varphi^{41}\otimes \varphi^{42}\dots \varphi^{44}\otimes \shift (\varphi^{45}\varphi^5)\otimes 1)\\
=\,&0\boldsymbol{\uwave{+}}0\boldsymbol{\uwave{+}}\underbrace{\varphi^1(\shift \otimes \shift)\partial(\shift^{-1}\varphi^2)\otimes \shift \varphi^3\otimes \shift \varphi^{41}\otimes \varphi^{42}\dots \varphi^{44}\otimes \shift (\varphi^{45}\varphi^5)\otimes 1}_{162}\\
&+\underbrace{\varphi^1\varphi^{21}\dots \varphi^{24}\otimes \shift (\varphi^{25}\varphi^3)\otimes \shift \varphi^{41}\otimes \varphi^{42}\dots \varphi^{44}\otimes \shift (\varphi^{45}\varphi^5)\otimes 1}_g\\
&-\underbrace{\varphi^1\varphi^{21}\dots \varphi^{24}\otimes \shift \varphi^{25}\otimes \shift (\varphi^3\varphi^{41})\otimes \varphi^{42}\dots \varphi^{44}\otimes \shift (\varphi^{45}\varphi^5)\otimes 1}_g\\
&+\underbrace{\varphi^1\varphi^{21}\dots \varphi^{24}\otimes \shift \varphi^{25}\otimes \shift \varphi^3\otimes \varphi^{41}\dots \varphi^{44}\otimes \shift (\varphi^{45}\varphi^5)\otimes 1}_{167}\\
&\boldsymbol{\uwave{+}}0\boldsymbol{\uwave{+}}0\boldsymbol{\uwave{+}}\underbrace{\varphi^1\varphi^{21}\dots \varphi^{24}\otimes \shift \varphi^{25}\otimes \shift \varphi^3\otimes \shift \varphi^{41}\otimes \varphi^{42}\varphi^{43}\varphi^{44}\hat{\omega}\varphi^{45}\varphi^5}_{122}\\
&-\underbrace{\varphi^1\varphi^{21}\dots \varphi^{24}\otimes \shift \varphi^{25}\otimes \shift \varphi^3\otimes \shift \varphi^{41}\otimes \varphi^{42}\dots \varphi^{45}\varphi^5\hat{\omega}}_{121}\\
&-\underbrace{\varphi^1\varphi^{21}\dots \varphi^{24}\otimes \shift \varphi^{25}\otimes \shift \varphi^3\otimes \shift \varphi^{41}\otimes \varphi^{42}\dots \varphi^{44}\shift \partial(\varphi^{45})\varphi^5}_y\\
&-\underbrace{\varphi^1\varphi^{21}\dots \varphi^{24}\otimes \shift \varphi^{25}\otimes \shift \varphi^3\otimes \shift \varphi^{41}\otimes \varphi^{42}\dots \varphi^{45}\shift \partial(\varphi^5)}_{gg}
\end{align*}
and
\begin{align*}
&-\mu_1^{\otimes}(-\varphi^1\varphi^2\otimes \shift (\varphi^3\varphi^{41})\otimes \varphi^{42}\varphi^{43}\varphi^{441}\dots \varphi^{444}\otimes \shift \varphi^{445}\otimes \shift \varphi^{45}\otimes \shift \varphi^5\otimes 1)\\
=\,&0\boldsymbol{\uwave{+}}\underbrace{\varphi^1\varphi^2\hat{\omega}\varphi^3\varphi^{41}\varphi^{42}\varphi^{43}\varphi^{441}\dots \varphi^{444}\otimes \shift \varphi^{445}\otimes \shift \varphi^{45}\otimes \shift \varphi^5\otimes 1}_{310}\\
&-\underbrace{\varphi^1\varphi^2\varphi^3\varphi^{41}\hat{\omega}\varphi^{42}\varphi^{43}\varphi^{441}\dots \varphi^{444}\otimes \shift \varphi^{445}\otimes \shift \varphi^{45}\otimes \shift \varphi^5\otimes 1}_{311}\\
&-\underbrace{\varphi^1\varphi^2\shift \partial(\varphi^3)\varphi^{41}\dots \varphi^{43}\varphi^{441}\dots \varphi^{444}\otimes \shift \varphi^{445}\otimes \shift \varphi^{45}\otimes \shift \varphi^5\otimes 1}_{gg}\\
&-\underbrace{\varphi^1\varphi^2\varphi^3\shift \partial(\varphi^{41})\varphi^{42}\varphi^{43}\varphi^{441}\dots \varphi^{444}\otimes \shift \varphi^{445}\otimes \shift \varphi^{45}\otimes \shift \varphi^5\otimes 1}_g\\
&\boldsymbol{\uwave{+}}0\boldsymbol{\uwave{+}}0\boldsymbol{\uwave{+}}0\boldsymbol{\uwave{+}}\underbrace{\varphi^1\varphi^2\otimes \shift (\varphi^3\varphi^{41})\otimes \varphi^{42}\varphi^{43}(\shift \otimes \shift)\partial(\shift^{-1}\varphi^{44})\otimes \shift \varphi^{45}\otimes \shift \varphi^5\otimes 1}_{235}\\
&+\underbrace{\varphi^1\varphi^2\otimes \shift (\varphi^3\varphi^{41})\otimes \varphi^{42}\varphi^{43}\varphi^{441}\dots \varphi^{444}\otimes \shift (\varphi^{445}\varphi^{45})\otimes \shift \varphi^5\otimes 1}_y\\
&-\underbrace{\varphi^1\varphi^2\otimes \shift (\varphi^3\varphi^{41})\otimes \varphi^{42}\varphi^{43}\varphi^{441}\dots \varphi^{444}\otimes \shift \varphi^{445}\otimes \shift (\varphi^{45}\varphi^5) \otimes 1}_g\\
&+\underbrace{\varphi^1\varphi^2\otimes \shift (\varphi^3\varphi^{41})\otimes \varphi^{42}\varphi^{43}\varphi^{441}\dots \varphi^{444}\otimes \shift \varphi^{445}\otimes \shift \varphi^{45}\otimes \varphi^5}_{284}
\end{align*}
and
\begin{align*}
&-\mu_1^{\otimes}(\varphi^1\varphi^2\varphi^3\varphi^{41}\varphi^{42}\otimes \shift (\varphi^{43}\varphi^{441})\otimes \varphi^{442}\dots \varphi^{444}\otimes \shift \varphi^{445}\otimes \shift \varphi^{45}\otimes \shift \varphi^5\otimes 1)\\
=\,&0\boldsymbol{\uwave{+}}0\boldsymbol{\uwave{-}}\underbrace{\varphi^1\varphi^2\varphi^3\varphi^{41}\varphi^{42}\hat{\omega}\varphi^{43}\varphi^{441}\dots \varphi^{444}\otimes \shift \varphi^{445}\otimes \shift \varphi^{45}\otimes \shift \varphi^5\otimes 1}_{315}\\
&+\underbrace{\varphi^1\varphi^2\varphi^3\varphi^{41}\dots \varphi^{43}\varphi^{441}\hat{\omega}\varphi^{442}\dots \varphi^{444}\otimes \shift \varphi^{445}\otimes \shift \varphi^{45}\otimes \shift \varphi^5\otimes 1}_{316}\\
&+\underbrace{\varphi^1\varphi^2\varphi^3\varphi^{41}\varphi^{42}\shift \partial(\varphi^{43})\varphi^{441}\dots \varphi^{444}\otimes \shift \varphi^{445}\otimes \shift \varphi^{45}\otimes \shift \varphi^5\otimes 1}_g\\
&+\underbrace{\varphi^1\varphi^2\varphi^3\varphi^{41}\dots \varphi^{43}\shift \partial(\varphi^{441})\varphi^{442}\dots \varphi^{444}\otimes \shift \varphi^{445}\otimes \shift \varphi^{45}\otimes \shift \varphi^5\otimes 1}_y\\
&\boldsymbol{\uwave{+}}0\boldsymbol{\uwave{+}}0\boldsymbol{\uwave{+}}\underbrace{\varphi^1\varphi^2\varphi^3\varphi^{41}\varphi^{42}\otimes \shift (\varphi^{43}\varphi^{441})\otimes \varphi^{442}\dots \varphi^{445}\otimes \shift \varphi^{45}\otimes \shift \varphi^5\otimes 1}_{288}\\
&-\underbrace{\varphi^1\varphi^2\varphi^3\varphi^{41}\varphi^{42}\otimes \shift (\varphi^{43}\varphi^{441})\otimes \varphi^{442}\dots \varphi^{444}\otimes \shift (\varphi^{445}\varphi^{45})\otimes \shift \varphi^5\otimes 1}_y\\
&+\underbrace{\varphi^1\varphi^2\varphi^3\varphi^{41}\varphi^{42}\otimes \shift (\varphi^{43}\varphi^{441})\otimes \varphi^{442}\dots \varphi^{444}\otimes \shift \varphi^{445}\otimes \shift \varphi^{45}\varphi^5\otimes 1}_g\\
&-\underbrace{\varphi^1\varphi^2\varphi^3\varphi^{41}\varphi^{42}\otimes \shift (\varphi^{43}\varphi^{441})\otimes \varphi^{442}\dots \varphi^{444}\otimes \shift \varphi^{445}\otimes \shift \varphi^{45}\otimes \varphi^5}_{217}
\end{align*}
and
\begin{align*}
&-\mu_1^{\otimes}(-1\otimes \shift \varphi^1\otimes \shift \varphi^{21}\otimes \shift \varphi^{221}\otimes \varphi^{222}\otimes \shift \varphi^{223}\otimes \varphi^{224}\varphi^{225}\varphi^{23}\dots \varphi^{25}\varphi^3\varphi^4\varphi^5)\\
=\,&+\underbrace{\varphi^1\otimes \shift \varphi^{21}\otimes \shift \varphi^{221}\otimes \varphi^{222}\otimes \shift \varphi^{223}\otimes \varphi^{224}\varphi^{225}\varphi^{23}\varphi^{24}\varphi^{25}\varphi^3\varphi^4\varphi^5}_{218}\\
&-\underbrace{1\otimes \shift (\varphi^1\varphi^{21})\otimes \shift \varphi^{221}\otimes \varphi^{222}\otimes \shift \varphi^{223}\otimes \varphi^{224}\varphi^{225}\varphi^{23}\dots \varphi^{25}\varphi^3\varphi^4\varphi^5}_g\\
&+\underbrace{1\otimes \shift \varphi^1\otimes \shift (\varphi^{21}\varphi^{221})\otimes \varphi^{222}\otimes \shift \varphi^{223}\otimes \varphi^{224}\varphi^{225}\varphi^{23}\dots \varphi^{25}\varphi^3\varphi^4\varphi^5}_y\\
&-\underbrace{1\otimes \shift \varphi^1\otimes \shift \varphi^{21}\otimes \varphi^{221}\varphi^{222}\otimes \shift \varphi^{223}\otimes \varphi^{224}\varphi^{225}\varphi^{23}\dots \varphi^{25}\varphi^3\varphi^4\varphi^5}_{219}\\
&\boldsymbol{\uwave{+}}0\boldsymbol{\uwave{-}}\underbrace{1\otimes \shift \varphi^1\otimes \shift \varphi^{21}\otimes \shift \varphi^{221}\otimes \varphi^{222}\hat{\omega}\varphi^{223}\dots \varphi^{225}\varphi^{23}\dots \varphi^{25}\varphi^3\varphi^4\varphi^5}_{317}\\
&+\underbrace{1\otimes \shift \varphi^1\otimes \shift \varphi^{21}\otimes \shift \varphi^{221}\otimes \varphi^{222}\varphi^{223}\hat{\omega}\varphi^{224}\varphi^{225}\varphi^{23}\dots \varphi^{25}\varphi^3\varphi^4\varphi^5}_{318}\\
&+\underbrace{1\otimes \shift \varphi^1\otimes \shift \varphi^{21}\otimes \shift \varphi^{221}\otimes \varphi^{222}\shift \partial(\varphi^{223})\varphi^{224}\varphi^{225}\varphi^{23}\dots \varphi^{25}\varphi^3\varphi^4\varphi^5}_y\boldsymbol{\uwave{+}}0\boldsymbol{\uwave{+}}0\boldsymbol{\uwave{+}}0
\end{align*}
and
\begin{align*}
&-\mu_1^{\otimes}(-1\otimes \shift \varphi^1\otimes \shift \varphi^{21}\otimes \shift \varphi^{221}\otimes \varphi^{222}\dots \varphi^{225}\varphi^{23}\varphi^{24}\otimes \shift (\varphi^{25}\varphi^3)\otimes \varphi^4\varphi^5)\\
=\,&+\underbrace{\varphi^1\otimes \shift \varphi^{21}\otimes \shift \varphi^{221}\otimes \varphi^{222}\dots \varphi^{225}\varphi^{23}\varphi^{24}\otimes \shift \varphi^{25}\varphi^3\otimes \varphi^4\varphi^5}_{238}\\
&-\underbrace{1\otimes \shift (\varphi^1\varphi^{21})\otimes \shift \varphi^{221}\otimes \varphi^{222}\dots \varphi^{225}\varphi^{23}\varphi^{24}\otimes \shift (\varphi^{25}\varphi^3)\otimes \varphi^4\varphi^5}_g\\
&+\underbrace{1\otimes \shift \varphi^1\otimes \shift (\varphi^{21}\varphi^{221})\otimes \varphi^{222}\dots \varphi^{225}\varphi^{23}\varphi^{24}\otimes \shift (\varphi^{25}\varphi^3)\otimes \varphi^4\varphi^5}_y\\
&+\underbrace{1\otimes \shift \varphi^1\otimes \shift \varphi^{21}\otimes (\shift \otimes \shift)\partial(\shift^{-1}\varphi^{22})\varphi^{23}\varphi^{24}\otimes \shift (\varphi^{25}\varphi^3)\otimes \varphi^4\varphi^5}_{225}\\
&\boldsymbol{\uwave{+}}0\boldsymbol{\uwave{+}}0\boldsymbol{\uwave{+}}0\boldsymbol{\uwave{-}}\underbrace{1\otimes \shift \varphi^1\otimes \shift \varphi^{21}\otimes \shift \varphi^{221}\otimes \varphi^{222}\dots \varphi^{225}\varphi^{23}\varphi^{24}\hat{\omega}\varphi^{25}\varphi^3\varphi^4\varphi^5}_{328}\\
&+\underbrace{1\otimes \shift \varphi^1\otimes \shift \varphi^{21}\otimes \shift \varphi^{221}\otimes \varphi^{222}\dots \varphi^{225}\varphi^{23}\varphi^{24} \varphi^{25}\varphi^3\hat{\omega}\varphi^4\varphi^5}_{329}\\
&+\underbrace{1\otimes \shift \varphi^1\otimes \shift \varphi^{21}\otimes \shift \varphi^{221}\otimes \varphi^{222}\dots \varphi^{225}\varphi^{23}\varphi^{24}\shift \partial(\varphi^{25})\varphi^3\varphi^4\varphi^5}_g\\
&+\underbrace{1\otimes \shift \varphi^1\otimes \shift \varphi^{21}\otimes \shift \varphi^{221}\otimes \varphi^{222}\dots \varphi^{225}\varphi^{23}\dots \varphi^{25}\shift \partial(\varphi^3)\varphi^4\varphi^5}_{gg}\boldsymbol{\uwave{+}}0
\end{align*}
and
\begin{align*}
&-\mu_1^{\otimes}(1\otimes \shift \varphi^1\otimes \shift \varphi^{21}\otimes \shift \varphi^{221}\otimes \varphi^{222}\dots \varphi^{224}\otimes \shift (\varphi^{225}\varphi^{23})\otimes \varphi^{24}\varphi^{25}\varphi^3\varphi^4\varphi^5)\\
=\,&-\underbrace{\varphi^1\otimes \shift \varphi^{21}\otimes \shift \varphi^{221}\otimes \varphi^{222}\dots \varphi^{224}\otimes \shift (\varphi^{225}\varphi^{23})\otimes \varphi^{24}\varphi^{25}\varphi^3\varphi^4\varphi^5}_{236}\\
&+\underbrace{1\otimes \shift (\varphi^1\varphi^{21})\otimes \shift \varphi^{221}\otimes \varphi^{222}\dots \varphi^{224}\otimes \shift (\varphi^{225}\varphi^{23})\otimes \varphi^{24}\varphi^{25}\varphi^3\varphi^4\varphi^5}_g\\
&-\underbrace{1\otimes \shift \varphi^1\otimes \shift (\varphi^{21}\varphi^{221})\otimes \varphi^{222}\dots \varphi^{224}\otimes \shift (\varphi^{225}\varphi^{23})\otimes \varphi^{24}\varphi^{25}\varphi^3\varphi^4\varphi^5}_y\\
&+\underbrace{1\otimes \shift \varphi^1\otimes \shift \varphi^{21}\otimes \varphi^{221}\dots \varphi^{224}\otimes \shift (\varphi^{225}\varphi^{23})\otimes \varphi^{24}\varphi^{25}\varphi^3\varphi^4\varphi^5}_{248}\\
&\boldsymbol{\uwave{+}}0\boldsymbol{\uwave{+}}0\boldsymbol{\uwave{+}}\underbrace{1\otimes \shift \varphi^1\otimes \shift \varphi^{21}\otimes \shift \varphi^{221}\otimes \varphi^{222}\dots \varphi^{224}\hat{\omega}\varphi^{225}\varphi^{23}\dots \varphi^{25}\varphi^3\varphi^4\varphi^5}_{324}\\
&-\underbrace{1\otimes \shift \varphi^1\otimes \shift \varphi^{21}\otimes \shift \varphi^{221}\otimes \varphi^{222}\dots \varphi^{225}\varphi^{23}\hat{\omega}\varphi^{24}\varphi^{25}\varphi^3\varphi^4\varphi^5}_{325}\\
&-\underbrace{1\otimes \shift \varphi^1\otimes \shift \varphi^{21}\otimes \shift \varphi^{221}\otimes \varphi^{222}\dots \varphi^{224}\shift \partial(\varphi^{225})\varphi^{23}\dots \varphi^{25}\varphi^3\varphi^4\varphi^5}_y\\
&-\underbrace{1\otimes \shift \varphi^1\otimes \shift \varphi^{21}\otimes \shift \varphi^{221}\otimes \varphi^{222}\dots \varphi^{225}\shift \partial(\varphi^{23})\varphi^{24}\varphi^{25}\varphi^3\varphi^4\varphi^5}_g\boldsymbol{\uwave{+}}0\boldsymbol{\uwave{+}}0
\end{align*}
and
\begin{align*}
&-\mu_1^{\otimes}(\varphi^1\varphi^{21}\varphi^{22}\otimes \shift (\varphi^{23}\varphi^{241})\otimes \varphi^{242}\dots \varphi^{244}\otimes \shift \varphi^{245}\otimes \shift \varphi^{25}\otimes \shift \varphi^3\otimes \varphi^4\varphi^5)\\
=\,&0\boldsymbol{\uwave{-}}\underbrace{\varphi^1\varphi^{21}\varphi^{22}\hat{\omega}\varphi^{23}\varphi^{241}\dots \varphi^{244}\otimes \shift \varphi^{245}\otimes \shift \varphi^{25}\otimes \shift \varphi^3\otimes \varphi^4\varphi^5}_{312}\\
&+\underbrace{\varphi^1\varphi^{21}\dots \varphi^{23}\varphi^{241}\hat{\omega}\varphi^{242}\dots \varphi^{244}\otimes \shift \varphi^{245}\otimes \shift \varphi^{25}\otimes \shift \varphi^3\otimes \varphi^4\varphi^5}_{313}\\
&+\underbrace{\varphi^1\varphi^{21}\varphi^{22}\shift \partial(\varphi^{23})\varphi^{241}\dots \varphi^{244}\otimes \shift \varphi^{245}\otimes \shift \varphi^{25}\otimes \shift \varphi^3\otimes \varphi^4\varphi^5}_g\\
&+\underbrace{\varphi^1\varphi^{21}\dots \varphi^{23}\shift \partial(\varphi^{241})\varphi^{242}\dots \varphi^{244}\otimes \shift \varphi^{245}\otimes \shift \varphi^{25}\otimes \shift \varphi^3\otimes \varphi^4\varphi^5}_y\\
&\boldsymbol{\uwave{+}}0\boldsymbol{\uwave{+}}0\boldsymbol{\uwave{+}}\underbrace{\varphi^1\varphi^{21}\varphi^{22}\otimes \shift (\varphi^{23}\varphi^{241})\otimes \varphi^{242}\dots \varphi^{245}\otimes \shift \varphi^{25}\otimes \shift \varphi^3\otimes \varphi^4\varphi^5}_{268}\\
&-\underbrace{\varphi^1\varphi^{21}\varphi^{22}\otimes \shift (\varphi^{23}\varphi^{241})\otimes \varphi^{242}\dots \varphi^{244}\otimes \shift (\varphi^{245}\varphi^{25})\otimes \shift \varphi^3\otimes \varphi^4\varphi^5}_y\\
&+\underbrace{\varphi^1\varphi^{21}\varphi^{22}\otimes \shift (\varphi^{23}\varphi^{241})\otimes \varphi^{242}\dots \varphi^{244}\otimes \shift \varphi^{245}\otimes \shift (\varphi^{25}\varphi^3)\otimes \varphi^4\varphi^5}_g\\
&-\underbrace{\varphi^1\varphi^{21}\varphi^{22}\otimes \shift (\varphi^{23}\varphi^{241})\otimes \varphi^{242}\dots \varphi^{244}\otimes \shift \varphi^{245}\otimes \shift \varphi^{25}\otimes \varphi^3\varphi^4\varphi^5}_{290}\boldsymbol{\uwave{+}}0
\end{align*}
and
\begin{align*}
&-\mu_1^{\otimes}(-\varphi^1\varphi^{21}\varphi^{22}\varphi^{23}\varphi^{241}\dots \varphi^{244}\otimes \shift \varphi^{245}\otimes \shift \varphi^{25}\otimes \shift \varphi^3\otimes \varphi^{4}\otimes \shift \varphi^{5}\otimes 1)\\
=&0\boldsymbol{\uwave{+}}0\boldsymbol{\uwave{+}}0\boldsymbol{\uwave{-}}\underbrace{\varphi^{1}\varphi^{21}\varphi^{22}\varphi^{23}(\shift \otimes \shift)\partial(\shift^{-1}\varphi^{24})\otimes \shift \varphi^{25}\otimes \shift \varphi^{3}\otimes \varphi^{4}\otimes \shift \varphi^{5}\otimes 1}_{154}\\
&-\underbrace{\varphi^{1}\varphi^{21}\varphi^{22}\varphi^{23}\varphi^{241}\dots \varphi^{244}\otimes \shift (\varphi^{245}\varphi^{25})\otimes \shift \varphi^{3}\otimes \varphi^{4}\otimes \shift \varphi^{5}\otimes 1}_y\\
&+\underbrace{\varphi^{1}\varphi^{21}\varphi^{22}\varphi^{23}\varphi^{241}\dots \varphi^{244}\otimes \shift \varphi^{245}\otimes \shift (\varphi^{25}\varphi^{3})\otimes \varphi^{4}\otimes \shift \varphi^{5}\otimes 1}_g\\
&-\underbrace{\varphi^1\varphi^{21}\varphi^{22}\varphi^{23}\varphi^{241}\dots \varphi^{244}\otimes \shift \varphi^{245}\otimes \shift \varphi^{25}\otimes \varphi^3\varphi^4\otimes \shift \varphi^5\otimes 1}_{270}\\
&\boldsymbol{\uwave{+}}0\boldsymbol{\uwave{-}}\underbrace{\varphi^1\varphi^{21}\varphi^{22}\varphi^{23}\varphi^{241}\dots \varphi^{244}\otimes \shift \varphi^{245}\otimes \shift \varphi^{25}\otimes \shift \varphi^3\otimes \varphi^4\hat{\omega} \varphi^5}_{146}\\
&+\underbrace{\varphi^1\varphi^{21}\varphi^{22}\varphi^{23}\varphi^{241}\dots \varphi^{244}\otimes \shift \varphi^{245}\otimes \shift \varphi^{25}\otimes \shift \varphi^{3}\otimes \varphi^4\varphi^5\hat{\omega}}_{147}\\
&+\underbrace{\varphi^1\varphi^{21}\varphi^{22}\varphi^{23}\varphi^{241}\dots \varphi^{244}\otimes \shift \varphi^{245}\otimes \shift \varphi^{25}\otimes \shift \varphi^3\otimes \varphi^4\shift \partial(\varphi^5)}_{g}
\end{align*}
and
\begin{align*}
&-\mu_1^{\otimes}(-\varphi^1\varphi^2\otimes \shift \varphi^3\otimes \shift \varphi^{41}\otimes \shift \varphi^{421}\otimes \varphi^{422}\otimes \shift \varphi^{423}\otimes \varphi^{424}\varphi^{425}\varphi^{43}\varphi^{44}\varphi^{45}\varphi^5)\\
=&0\boldsymbol{\uwave{+}}\underbrace{\varphi^1\varphi^2\varphi^3\otimes \shift \varphi^{41}\otimes \shift \varphi^{421}\otimes \varphi^{422}\otimes \shift \varphi^{423}\otimes \varphi^{424}\varphi^{425}\varphi^{43}\varphi^{44}\varphi^{45}\varphi^5}_{220}\\
&-\underbrace{\varphi^1\varphi^2\otimes \shift (\varphi^3\varphi^{41})\otimes \shift \varphi^{421}\otimes \varphi^{422}\otimes \shift \varphi^{423}\otimes \varphi^{424}\varphi^{425}\varphi^{43}\varphi^{44}\varphi^{45}\varphi^5}_g\\
&+\underbrace{\varphi^1\varphi^2\otimes \shift \varphi^3\otimes \shift (\varphi^{41}\varphi^{421})\otimes \varphi^{422}\otimes \shift \varphi^{423}\otimes \varphi^{424}\varphi^{425}\varphi^{43}\varphi^{44}\varphi^{45}\varphi^5}_y\\
&-\underbrace{\varphi^1\varphi^2\otimes \shift \varphi^3\otimes \shift \varphi^{41}\otimes \varphi^{421}\varphi^{422}\otimes \shift \varphi^{423}\otimes \varphi^{424}\varphi^{425}\varphi^{43}\varphi^{44}\varphi^{45}\varphi^5}_{156}\\
&\boldsymbol{\uwave{+}}0\boldsymbol{\uwave{-}}\underbrace{\varphi^1\varphi^2\otimes \shift \varphi^3\otimes \shift \varphi^{41}\otimes \shift \varphi^{421}\otimes \varphi^{422}\hat{\omega} \varphi^{423}\varphi^{424}\varphi^{425}\varphi^{43}\varphi^{44}\varphi^{45}\varphi^5}_{326}\\
&+\underbrace{\varphi^1\varphi^2\otimes \shift \varphi^3\otimes \shift \varphi^{41}\otimes \shift \varphi^{421}\otimes \varphi^{422}\varphi^{423}\hat{\omega}\varphi^{424}\varphi^{425}\varphi^{43}\varphi^{44}\varphi^{45}\varphi^5}_{327}\\
&+\underbrace{\varphi^1\varphi^2\otimes \shift \varphi^3\otimes \shift \varphi^{41}\otimes \shift \varphi^{421}\otimes \varphi^{422}\shift \partial(\varphi^{423})\varphi^{424}\varphi^{425}\varphi^{43}\varphi^{44}\varphi^{45}\varphi^5}_{g}\boldsymbol{\uwave{+}}0\boldsymbol{\uwave{+}}0
\end{align*}
and
\begin{align*}
&-\mu_1^{\otimes}(\varphi^1\varphi^2\otimes \shift \varphi^3\otimes \shift \varphi^{41}\otimes \shift \varphi^{421}\otimes \varphi^{422}\varphi^{423}\varphi^{424}\otimes \shift (\varphi^{425}\varphi^{43})\otimes \varphi^{44}\varphi^{45}\varphi^5)\\
=&0\boldsymbol{\uwave{-}}\underbrace{\varphi^1\varphi^2\varphi^3\otimes \shift \varphi^{41}\otimes \shift \varphi^{421}\otimes \varphi^{422}\varphi^{423}\varphi^{424}\otimes \shift (\varphi^{425}\varphi^{43})\otimes \varphi^{44}\varphi^{45}\varphi^5}_{289}\\
&+\underbrace{\varphi^1\varphi^2\otimes \shift (\varphi^3\varphi^{41})\otimes \shift \varphi^{421}\otimes \varphi^{422}\varphi^{423}\varphi^{424}\otimes \shift (\varphi^{425}\varphi^{43})\otimes \varphi^{44}\varphi^{45}\varphi^5}_g\\
&-\underbrace{\varphi^1\varphi^2\otimes \shift \varphi^3\otimes \shift (\varphi^{41}\varphi^{421})\otimes \varphi^{422}\varphi^{423}\varphi^{424}\otimes \shift (\varphi^{425}\varphi^{43})\otimes \varphi^{44}\varphi^{45}\varphi^5}_g\\
&+\underbrace{\varphi^1\varphi^2\otimes \shift \varphi^3\otimes \shift \varphi^{41}\otimes \varphi^{421}\dots \varphi^{424}\otimes \shift (\varphi^{425}\varphi^{43})\otimes \varphi^{44}\varphi^{45}\varphi^5}_{159}\\
&\boldsymbol{\uwave{+}}0\boldsymbol{\uwave{+}}0\boldsymbol{\uwave{+}}\underbrace{\varphi^1\varphi^2\otimes \shift \varphi^3\otimes \shift \varphi^{41}\otimes \shift \varphi^{421}\otimes \varphi^{422}\varphi^{423}\varphi^{424}\hat{\omega} \varphi^{425}\varphi^{43}\varphi^{44}\varphi^{45}\varphi^5}_{330}\\
&-\underbrace{\varphi^1\varphi^2\otimes \shift \varphi^3\otimes \shift \varphi^{41}\otimes \shift \varphi^{421}\otimes \varphi^{422}\dots \varphi^{425}\varphi^{43}\hat{\omega} \varphi^{44}\varphi^{45}\varphi^5}_{331}\\
&-\underbrace{\varphi^1\varphi^2\otimes \shift \varphi^3\otimes \shift \varphi^{41}\otimes \shift \varphi^{421}\otimes \varphi^{422}\varphi^{423}\varphi^{424}\shift \partial(\varphi^{425})\varphi^{43}\varphi^{44}\varphi^{45}\varphi^5}_{g}\\
&-\underbrace{\varphi^1\varphi^2\otimes \shift \varphi^3\otimes \shift \varphi^{41}\otimes \shift \varphi^{421}\otimes \varphi^{422}\varphi^{423}\varphi^{424}\varphi^{425}\shift \partial(\varphi^{43})\varphi^{44}\varphi^{45}\varphi^5}_{g}\boldsymbol{\uwave{+}}0
\end{align*}
and
\begin{align*}
&-\mu_1^{\otimes}(1\otimes \shift \varphi^1\otimes \varphi^{21}\dots \varphi^{24}\otimes \shift \varphi^{25}\otimes \shift \varphi^3\otimes \varphi^4\otimes \shift \varphi^5\otimes 1)\\
=&-\underbrace{\hat{\omega}\varphi^1\varphi^{21}\dots \varphi^{24}\otimes \shift \varphi^{25}\otimes \shift \varphi^3\otimes \varphi^4\otimes \shift \varphi^5\otimes 1}_{54}\\
&+\underbrace{\varphi^1\hat{\omega}\varphi^{21}\dots \varphi^{24}\otimes \shift \varphi^{25}\otimes \shift \varphi^3\otimes \varphi^4\otimes \shift \varphi^5\otimes 1}_{72}\\
&+\underbrace{\shift \partial(\varphi^1)\varphi^{21}\dots \varphi^{24}\otimes \shift \varphi^{25}\otimes \shift \varphi^3\otimes \varphi^4\otimes \shift \varphi^5\otimes 1}_g\\
&\boldsymbol{\uwave{+}}0\boldsymbol{\uwave{+}}0\boldsymbol{\uwave{-}}\underbrace{1\otimes \shift \varphi^1\otimes (\shift \otimes \shift)\partial(\shift^{-1}\varphi^2)\otimes \shift \varphi^3\otimes \varphi^4\otimes \shift \varphi^5\otimes 1}_{55}\\
&-\underbrace{1\otimes \shift \varphi^1\otimes \varphi^{21}\dots \varphi^{24}\otimes \shift (\varphi^{25}\varphi^3)\otimes \varphi^4\otimes \shift \varphi^5\otimes 1}_{149}\\
&+\underbrace{1\otimes \shift \varphi^1\otimes \varphi^{21}\dots \varphi^{24}\otimes \shift \varphi^{25}\otimes \varphi^3\varphi^4\otimes \shift \varphi^5\otimes 1}_{230}\\
&\boldsymbol{\uwave{+}}0\boldsymbol{\uwave{+}}\underbrace{1\otimes \shift \varphi^1\otimes \varphi^{21}\dots \varphi^{24}\otimes \shift \varphi^{25}\otimes \shift \varphi^3\otimes \varphi^4\hat{\omega} \varphi^5}_{53}\\
&-\underbrace{1\otimes \shift \varphi^1\otimes \varphi^{21}\dots \varphi^{24}\otimes \shift \varphi^{25}\otimes \shift \varphi^3\otimes \varphi^4\varphi^5\hat{\omega}}_{52}\\
&-\underbrace{1\otimes \shift \varphi^1\otimes \varphi^{21}\dots \varphi^{24}\otimes \shift \varphi^{25}\otimes \shift \varphi^3\otimes \varphi^4\shift \partial(\varphi^5)}_{g}
\end{align*}
and
\begin{align*}
&-\mu_1^{\otimes}(-1\otimes \shift \varphi^1\otimes \varphi^{2}\otimes \shift (\varphi^3\varphi^{41})\otimes \varphi^{42}\varphi^{43}\varphi^{44}\otimes \shift \varphi^{45}\otimes \shift \varphi^5\otimes 1)\\
=&\underbrace{\hat{\omega} \varphi^1\varphi^2\otimes \shift (\varphi^3\varphi^{41})\otimes \varphi^{42}\varphi^{43}\varphi^{44}\otimes \shift \varphi^{45}\otimes \shift \varphi^5\otimes 1}_{83}\\
-&\underbrace{\varphi^1\hat{\omega}\varphi^2\otimes \shift (\varphi^3\varphi^{41})\otimes \varphi^{42}\varphi^{43}\varphi^{44}\otimes \shift \varphi^{45}\otimes \shift \varphi^5\otimes 1}_{84}\\
&-\underbrace{\shift \partial(\varphi^1)\varphi^2\otimes \shift (\varphi^3\varphi^{41})\otimes \varphi^{42}\varphi^{43}\varphi^{44}\otimes \shift \varphi^{45}\otimes \shift \varphi^5\otimes 1}_g\\
&\boldsymbol{\uwave{+}}0\boldsymbol{\uwave{-}}\underbrace{1\otimes \shift \varphi^1\otimes \varphi^2\hat{\omega} \varphi^3\varphi^{41}\dots \varphi^{44}\otimes \shift \varphi^{45}\otimes \shift \varphi^5\otimes 1}_{85}\\
&+\underbrace{1\otimes \shift \varphi^1\otimes \varphi^2\varphi^3\varphi^{41}\hat{\omega}\varphi^{42}\varphi^{43}\varphi^{44}\otimes \shift \varphi^{45}\otimes \shift \varphi^5\otimes 1}_{86}\\
&+\underbrace{1\otimes \shift \varphi^1\otimes \varphi^2\shift \partial(\varphi^3)\varphi^{41}\varphi^{42}\varphi^{43}\varphi^{44}\otimes \shift \varphi^{45}\otimes \shift \varphi^5\otimes 1}_{g}\\
&+\underbrace{1\otimes \shift \varphi^1\otimes \varphi^2\varphi^3\shift \partial(\varphi^{41})\varphi^{42}\varphi^{43}\varphi^{44}\otimes \shift \varphi^{45}\otimes \shift \varphi^5\otimes 1}_{g}\\
&\boldsymbol{\uwave{+}}0\boldsymbol{\uwave{+}}0\boldsymbol{\uwave{+}}\underbrace{1\otimes \shift \varphi^1\otimes \varphi^2\otimes \shift (\varphi^3\varphi^{41})\otimes \varphi^{42}\dots \varphi^{45}\otimes \shift \varphi^5\otimes 1}_{195}\\
&-\underbrace{1\otimes \shift \varphi^1\otimes \varphi^2\otimes \shift (\varphi^3\varphi^{41})\otimes \varphi^{42}\varphi^{43}\varphi^{44}\otimes \shift (\varphi^{45}\varphi^5)\otimes 1}_{y}\\
&+\underbrace{1\otimes \shift \varphi^1\otimes \varphi^2\otimes \shift (\varphi^3\varphi^{41})\otimes \varphi^{42}\varphi^{43}\varphi^{44}\otimes \shift \varphi^{45}\otimes \varphi^5}_{194}
\end{align*}
and
\begin{align*}
&-\mu_1^{\otimes}(1\otimes \shift \varphi^1\otimes 1\otimes \shift \varphi^{21}\otimes \varphi^{22}\varphi^{23}\varphi^{24}\otimes \shift \varphi^{25}\otimes \shift \varphi^3\otimes \varphi^4\varphi^5)\\
=&-\underbrace{\hat{\omega} \varphi^1\otimes \shift \varphi^{21}\otimes \varphi^{22}\varphi^{23}\varphi^{24}\otimes \shift \varphi^{25}\otimes \shift \varphi^3\otimes \varphi^4\varphi^5}_{152}\\
&+\underbrace{\varphi^1\hat{\omega}\otimes \shift \varphi^{21}\otimes \varphi^{22}\varphi^{23}\varphi^{24}\otimes \shift \varphi^{25}\otimes \shift \varphi^3\otimes \varphi^4\varphi^5}_{94}\\
&+\underbrace{\shift \partial(\varphi^1)\otimes \shift \varphi^{21}\otimes \varphi^{22}\varphi^{23}\varphi^{24}\otimes \shift \varphi^{25}\otimes \shift \varphi^3\otimes \varphi^4\varphi^5}_{g}\\
&\boldsymbol{\uwave{+}}\underbrace{1\otimes \shift \varphi^1\otimes \hat{\omega} \varphi^{21}\dots \varphi^{24}\otimes \shift \varphi^{25}\otimes \shift \varphi^3\otimes \varphi^4\varphi^5}_{96}\\
&-\underbrace{1\otimes \shift \varphi^1\otimes \varphi^{21}\hat{\omega}\varphi^{22}\varphi^{23}\varphi^{24}\otimes \shift \varphi^{25}\otimes \shift \varphi^3\otimes \varphi^4\varphi^5}_{97}\\
&-\underbrace{1\otimes \shift \varphi^1\otimes \shift \partial(\varphi^{21})\varphi^{22}\varphi^{23}\varphi^{24}\otimes \shift \varphi^{25}\otimes \shift \varphi^3\otimes \varphi^4\varphi^5}_y\\
&\boldsymbol{\uwave{+}}0\boldsymbol{\uwave{+}}0\boldsymbol{\uwave{-}}\underbrace{1\otimes \shift \varphi^1\otimes 1\otimes \shift \varphi^{21}\otimes \varphi^{22}\dots \varphi^{25}\otimes \shift \varphi^3\otimes \varphi^4\varphi^5}_{355}\\
&+\underbrace{1\otimes \shift \varphi^1\otimes 1\otimes \shift \varphi^{21}\otimes \varphi^{22}\varphi^{23}\varphi^{24}\otimes \shift (\varphi^{25}\varphi^3)\otimes \varphi^4\varphi^5}_{95}\\
&-\underbrace{1\otimes \shift \varphi^1\otimes 1\otimes \shift \varphi^{21}\otimes \varphi^{22}\varphi^{23}\varphi^{24}\otimes \shift \varphi^{25}\otimes \varphi^3\varphi^4\varphi^5}_{356}\boldsymbol{\uwave{+}}0
\end{align*}
and
\begin{align*}
&-\mu_1^{\otimes}(1\otimes \shift \varphi^1\otimes \shift \varphi^{21}\otimes \varphi^{22}\otimes \shift \varphi^{23}\otimes \varphi^{24}\varphi^{25}\varphi^3\varphi^4\otimes \shift \varphi^5\otimes 1)\\
=&-\underbrace{\varphi^1\otimes \shift \varphi^{21}\otimes \varphi^{22}\otimes \shift \varphi^{23}\otimes \varphi^{24}\varphi^{25}\varphi^3\varphi^4\otimes \shift \varphi^5\otimes 1}_{196}\\
&+\underbrace{1\otimes \shift (\varphi^1\varphi^{21})\otimes \varphi^{22}\otimes \shift \varphi^{23}\otimes \varphi^{24}\varphi^{25}\varphi^3\varphi^4\otimes \shift \varphi^5\otimes 1}_{y}\\
&-\underbrace{1\otimes \shift \varphi^1\otimes \varphi^{21}\varphi^{22}\otimes \shift \varphi^{23}\otimes \varphi^{24}\varphi^{25}\varphi^3\varphi^4\otimes \shift \varphi^5\otimes 1}_{197}\\
&\boldsymbol{\uwave{+}}0\boldsymbol{\uwave{-}}\underbrace{1\otimes \shift \varphi^1\otimes \shift \varphi^{21}\otimes \varphi^{22}\hat{\omega} \varphi^{23}\varphi^{24}\varphi^{25}\varphi^3\varphi^4\otimes \shift \varphi^5\otimes 1}_{108}\\
&+\underbrace{1\otimes \shift \varphi^1\otimes \shift \varphi^{21}\otimes \varphi^{22}\varphi^{23}\hat{\omega}\varphi^{24}\varphi^{25}\varphi^3\varphi^4\otimes \shift \varphi^5\otimes 1}_{109}\\
&+\underbrace{1\otimes \shift \varphi^1\otimes \shift \varphi^{21}\otimes \varphi^{22}\shift \partial(\varphi^{23})\varphi^{24}\varphi^{25}\varphi^3\varphi^4\otimes \shift \varphi^5\otimes 1}_y\\
&\boldsymbol{\uwave{+}}0\boldsymbol{\uwave{+}}0\boldsymbol{\uwave{+}}\underbrace{1\otimes \shift \varphi^1\otimes \shift \varphi^{21}\otimes \varphi^{22}\otimes \shift \varphi^{23}\otimes \varphi^{24}\varphi^{25}\varphi^3\varphi^4\hat{\omega} \varphi^5}_{107}\\
&-\underbrace{1\otimes \shift \varphi^1\otimes \shift \varphi^{21}\otimes \varphi^{22}\otimes \shift \varphi^{23}\otimes \varphi^{24}\varphi^{25}\varphi^3\varphi^4\varphi^5\hat{\omega}}_{106}\\
&-\underbrace{1\otimes \shift \varphi^1\otimes \shift \varphi^{21}\otimes \varphi^{22}\otimes \shift \varphi^{23}\otimes \varphi^{24}\varphi^{25}\varphi^3\varphi^4\shift \partial(\varphi^5)}_g
\end{align*}
and
\begin{align*}
&-\mu_1^{\otimes}(-1\otimes \shift \varphi^1\otimes \shift \varphi^{21}\otimes \varphi^{22}\varphi^{23}\varphi^{24}\otimes \shift (\varphi^{25}\varphi^3)\otimes \varphi^4\otimes \shift \varphi^5\otimes 1)\\
=&+\underbrace{\varphi^1\otimes \shift \varphi^{21}\otimes \varphi^{22}\varphi^{23}\varphi^{24}\otimes \shift (\varphi^{25}\varphi^3)\otimes \varphi^4\otimes \shift \varphi^5\otimes 1}_{153}\\
&-\underbrace{1\otimes \shift (\varphi^1\varphi^{21})\otimes \varphi^{22}\varphi^{23}\varphi^{24}\otimes \shift (\varphi^{25}\varphi^3)\otimes \varphi^4\otimes \shift \varphi^5\otimes 1}_y\\
&+\underbrace{1\otimes \shift \varphi^1\otimes \varphi^{21}\dots \varphi^{24}\otimes \shift (\varphi^{25}\varphi^3)\otimes \varphi^4\otimes \shift \varphi^5\otimes 1}_{149}\\
&\boldsymbol{\uwave{+}}0\boldsymbol{\uwave{+}}0\boldsymbol{\uwave{+}}\underbrace{1\otimes \shift \varphi^1\otimes \shift \varphi^{21}\otimes \varphi^{22}\varphi^{23}\varphi^{24}\hat{\omega}\varphi^{25}\varphi^3\varphi^4\otimes \shift \varphi^5\otimes 1}_{99}\\
&-\underbrace{1\otimes \shift \varphi^1\otimes \shift \varphi^{21}\otimes \varphi^{22}\dots \varphi^{25}\varphi^3\hat{\omega} \varphi^4\otimes \shift \varphi^5\otimes 1}_{100}\\
&-\underbrace{1\otimes \shift \varphi^1\otimes \shift \varphi^{21}\otimes \varphi^{22}\varphi^{23}\varphi^{24}\shift \partial(\varphi^{25})\varphi^3\varphi^4\otimes \shift \varphi^5\otimes 1}_y\\
&-\underbrace{1\otimes \shift \varphi^1\otimes \shift \varphi^{21}\otimes \varphi^{22}\dots \varphi^{25}\shift \partial(\varphi^3)\otimes \varphi^4\otimes \shift \varphi^5\otimes 1}_g\\
&\boldsymbol{\uwave{+}}0\boldsymbol{\uwave{-}}\underbrace{1\otimes \shift \varphi^1\otimes \shift \varphi^{21}\otimes \varphi^{22}\varphi^{23}\varphi^{24}\otimes \shift (\varphi^{25}\varphi^3)\otimes \varphi^4\hat{\omega}\varphi^5}_{101}\\
&+\underbrace{1\otimes \shift \varphi^1\otimes \shift \varphi^{21}\otimes \varphi^{22}\varphi^{23}\varphi^{24}\otimes \shift (\varphi^{25}\varphi^3)\otimes \varphi^4\varphi^5\hat{\omega}}_{98}\\
&+\underbrace{1\otimes \shift \varphi^1\otimes \shift \varphi^{21}\otimes \varphi^{22}\varphi^{23}\varphi^{24}\otimes \shift (\varphi^{25}\varphi^3)\otimes \varphi^4\shift \partial(\varphi^5)}_{g}
\end{align*}
and
\begin{align*}
&-\mu_1^{\otimes}(-\varphi^1\varphi^2\otimes \shift \varphi^3\otimes \shift \varphi^{41}\otimes \varphi^{42}\otimes \shift \varphi^{43}\otimes \shift \varphi^{441}\otimes \varphi^{442}\dots \varphi^{445}\varphi^{45}\varphi^5)\\
=&0\boldsymbol{\uwave{+}}\underbrace{\varphi^1\varphi^2\varphi^3\otimes \shift \varphi^{41}\otimes \varphi^{42}\otimes \shift \varphi^{43}\otimes \shift \varphi^{441}\otimes \varphi^{442}\dots \varphi^{445}\varphi^{45}\varphi^5}_{294}\\
&-\underbrace{\varphi^1\varphi^2\otimes \shift (\varphi^3\varphi^{41})\otimes \varphi^{42}\otimes \shift \varphi^{43}\otimes \shift \varphi^{441}\otimes \varphi^{442}\dots \varphi^{445}\varphi^{45}\varphi^5}_g\\
&+\underbrace{\varphi^1\varphi^2\otimes \shift \varphi^3\otimes \varphi^{41}\varphi^{42}\otimes \shift \varphi^{43}\otimes \shift \varphi^{441}\otimes \varphi^{442}\dots \varphi^{445}\varphi^{45}\varphi^5}_{295}\\
&\boldsymbol{\uwave{+}}0\boldsymbol{\uwave{+}}\underbrace{\varphi^1\varphi^2\otimes \shift \varphi^3\otimes \shift \varphi^{41}\otimes \varphi^{42}\varphi^{43}\otimes \shift \varphi^{441}\otimes \varphi^{442}\dots \varphi^{445}\varphi^{45}\varphi^5}_{293}\\
&-\underbrace{\varphi^1\varphi^2\otimes \shift \varphi^3\otimes \shift \varphi^{41}\otimes \varphi^{42}\otimes \shift (\varphi^{43}\varphi^{441})\otimes \varphi^{442}\dots \varphi^{445}\varphi^{45}\varphi^5}_y\\
&-\underbrace{\varphi^1\varphi^2\otimes \shift \varphi^3\otimes \shift \varphi^{41}\otimes \varphi^{42}\otimes \shift \varphi^{43}\otimes (\shift \otimes \shift)\partial(\shift^{-1}\varphi^{44})\varphi^{45}\varphi^5}_{292}\boldsymbol{\uwave{+}}0\boldsymbol{\uwave{+}}0
\end{align*}
and
\begin{align*}
&-\mu_1^{\otimes}(\varphi^1\varphi^2\otimes \shift \varphi^3\otimes \shift \varphi^{41}\otimes \varphi^{42}\varphi^{43}\varphi^{44}\otimes \shift \varphi^{45}\otimes 1\otimes \shift \varphi^5\otimes 1)\\
=&0\boldsymbol{\uwave{-}}\underbrace{\varphi^1\varphi^2\varphi^3\otimes \shift \varphi^{41}\otimes \varphi^{42}\varphi^{43}\varphi^{44}\otimes \shift \varphi^{45}\otimes 1\otimes \shift \varphi^5\otimes 1}_{150}\\
&+\underbrace{\varphi^1\varphi^2\otimes \shift (\varphi^3\varphi^{41})\otimes \varphi^{42}\varphi^{43}\varphi^{44}\otimes \shift \varphi^{45}\otimes 1\otimes \shift \varphi^5\otimes 1}_{354}\\
&-\underbrace{\varphi^1\varphi^2\otimes \shift \varphi^3\otimes \varphi^{41}\dots \varphi^{44}\otimes \shift \varphi^{45}\otimes 1\otimes \shift \varphi^5\otimes 1}_{357}\\
&\boldsymbol{\uwave{+}}0\boldsymbol{\uwave{+}}0\boldsymbol{\uwave{-}}\underbrace{\varphi^1\varphi^2\otimes \shift \varphi^3\otimes \shift \varphi^{41}\otimes \varphi^{42}\varphi^{43}\varphi^{44}\hat{\omega}\varphi^{45}\otimes \shift \varphi^5\otimes 1}_{119}\\
&+\underbrace{\varphi^1\varphi^2\otimes \shift \varphi^3\otimes \shift \varphi^{41}\otimes \varphi^{42}\dots \varphi^{45}\hat{\omega}\otimes \shift \varphi^5\otimes 1}_{120}\\
&+\underbrace{\varphi^1\varphi^2\otimes \shift \varphi^3\otimes \shift \varphi^{41}\otimes \varphi^{42}\dots \shift \partial(\varphi^{45})\otimes \shift \varphi^5\otimes 1}_{g}\\
&\boldsymbol{\uwave{+}}\underbrace{\varphi^1\varphi^2\otimes \shift \varphi^3\otimes \shift \varphi^{41}\otimes \varphi^{42}\varphi^{43}\varphi^{44}\otimes \shift \varphi^{45}\otimes \hat{\omega}\varphi^5}_{118}\\
&-\underbrace{\varphi^1\varphi^2\otimes \shift \varphi^3\otimes \shift \varphi^{41}\otimes \varphi^{42}\varphi^{43}\varphi^{44}\otimes \shift \varphi^{45}\otimes \varphi^5\hat{\omega}}_{117}\\
&-\underbrace{\varphi^1\varphi^2\otimes \shift \varphi^3\otimes \shift \varphi^{41}\otimes \varphi^{42}\varphi^{43}\varphi^{44}\otimes \shift \varphi^{45}\otimes \shift \partial(\varphi^5)}_g
\end{align*}
and
\begin{align*}
&-\mu_1^{\otimes}(1\otimes \shift \varphi^1\otimes \shift \varphi^{21}\otimes \shift \varphi^{221}\otimes \varphi^{222}\dots \varphi^{225}\varphi^{23}\varphi^{24}\varphi^{25}\varphi^3\varphi^4\otimes \shift \varphi^5\otimes 1)\\
=&-\underbrace{\varphi^1\otimes \shift \varphi^{21}\otimes \shift \varphi^{221}\otimes \varphi^{222}\dots \varphi^{225}\varphi^{23}\varphi^{24}\varphi^{25}\varphi^3\varphi^4\otimes \shift \varphi^5\otimes 1}_{226}\\
&+\underbrace{1\otimes \shift (\varphi^1\varphi^{21})\otimes \shift \varphi^{221}\otimes \varphi^{222}\dots \varphi^{225}\varphi^{23}\varphi^{24}\varphi^{25}\varphi^3\varphi^4\otimes \shift \varphi^5\otimes 1}_g\\
&-\underbrace{1\otimes \shift \varphi^1\otimes \shift (\varphi^{21}\varphi^{221})\otimes \varphi^{222}\dots \varphi^{225}\varphi^{23}\varphi^{24}\varphi^{25}\varphi^3\varphi^4\otimes \shift \varphi^5\otimes 1}_y\\
&-\underbrace{1\otimes \shift \varphi^1\otimes \shift \varphi^{21}\otimes (\shift \otimes \shift)\partial(\shift^{-1}\varphi^{22})\varphi^{23}\varphi^{24}\varphi^{25}\varphi^3\varphi^4\otimes \shift \varphi^5\otimes 1}_{178}\\
&\boldsymbol{\uwave{+}}0\boldsymbol{\uwave{+}}0\boldsymbol{\uwave{+}}0\boldsymbol{\uwave{+}}0\boldsymbol{\uwave{+}}\underbrace{1\otimes \shift \varphi^1\otimes \shift \varphi^{21}\otimes \shift \varphi^{221}\otimes \varphi^{222}\dots \varphi^{225}\varphi^{23}\varphi^{24}\varphi^{25}\varphi^3\varphi^4\hat{\omega} \varphi^5}_{305}\\
&-\underbrace{1\otimes \shift \varphi^1\otimes \shift \varphi^{21}\otimes \shift \varphi^{221}\otimes \varphi^{222}\dots \varphi^{225}\varphi^{23}\varphi^{24}\varphi^{25}\varphi^3\varphi^4\varphi^5\hat{\omega}}_{300}\\
&-\underbrace{1\otimes \shift \varphi^1\otimes \shift \varphi^{21}\otimes \shift \varphi^{221}\otimes \varphi^{222}\dots \varphi^{225}\varphi^{23}\varphi^{24}\varphi^{25}\varphi^3\varphi^4\shift \partial(\varphi^5)}_{gg}
\end{align*}
and
\begin{align*}
&-\mu_1^{\otimes}(1\otimes \shift \varphi^1\otimes \varphi^2\varphi^3\varphi^{41}\varphi^{42}\varphi^{43}\varphi^{441}\dots \varphi^{444}\otimes \shift \varphi^{445}\otimes \shift \varphi^{45}\otimes \shift \varphi^5\otimes 1)\\
=&-\underbrace{\hat{\omega}\varphi^1\varphi^2\varphi^3\varphi^{41}\varphi^{42}\varphi^{43}\varphi^{441}\dots \varphi^{444}\otimes \shift \varphi^{445}\otimes \shift \varphi^{45}\otimes \shift \varphi^{5}\otimes 1}_{301}\\
&+\underbrace{\varphi^1\hat{\omega}\varphi^2\varphi^3\varphi^{41}\varphi^{42}\varphi^{43}\varphi^{441}\dots \varphi^{444}\otimes \shift \varphi^{445}\otimes \shift \varphi^{45}\otimes \shift \varphi^5\otimes 1}_{304}\\
&+\underbrace{\shift \partial(\varphi^1)\varphi^2\varphi^3\varphi^{41}\varphi^{42}\varphi^{43}\varphi^{441}\dots \varphi^{444}\otimes \shift \varphi^{445}\otimes \shift \varphi^{45}\otimes \shift \varphi^5\otimes 1}_{gg}\\
&\boldsymbol{\uwave{+}}0\boldsymbol{\uwave{+}}0\boldsymbol{\uwave{+}}0\boldsymbol{\uwave{+}}0\boldsymbol{\uwave{-}}\underbrace{1\otimes \shift \varphi^1\otimes \varphi^2\varphi^3\varphi^{41}\varphi^{42}\varphi^{43}(\shift \otimes \shift)\partial(\shift^{-1}\varphi^{44})\otimes \shift \varphi^{45}\otimes \shift \varphi^5\otimes 1}_{227}\\
&-\underbrace{1\otimes \shift \varphi^1\otimes \varphi^2\varphi^3\varphi^{41}\varphi^{42}\varphi^{43}\varphi^{441}\dots \varphi^{444}\otimes \shift (\varphi^{445}\varphi^{45})\otimes \shift \varphi^5\otimes 1}_{g}\\
&+\underbrace{1\otimes \shift \varphi^1\otimes \varphi^2\varphi^3\varphi^{41}\varphi^{42}\varphi^{43}\varphi^{441}\dots \varphi^{444}\otimes \shift \varphi^{445}\otimes \shift (\varphi^{45}\varphi^5)\otimes 1}_{g}\\
&-\underbrace{1\otimes \shift \varphi^1\otimes \varphi^2\varphi^3\varphi^{41}\varphi^{42}\varphi^{43}\varphi^{441}\dots \varphi^{444}\otimes \shift \varphi^{445}\otimes \shift \varphi^{45}\otimes \varphi^5}_{228}
\end{align*}
and
\begin{align*}
&-\mu_1^{\otimes}(-1\otimes \shift \varphi^1\otimes \varphi^2\otimes \shift \varphi^3\otimes \shift \varphi^{41}\otimes \shift \varphi^{421}\otimes \varphi^{422}\dots \varphi^{425}\varphi^{43}\varphi^{44}\varphi^{45}\varphi^5)\\
=&+\underbrace{\hat{\omega}\varphi^1\varphi^2\otimes \shift \varphi^3\otimes \shift \varphi^{41}\otimes \shift \varphi^{421}\otimes \varphi^{422}\dots \varphi^{425}\varphi^{43}\varphi^{44}\varphi^{45}\varphi^5}_{302}\\
&-\underbrace{\varphi^1\hat{\omega}\varphi^2\otimes \shift \varphi^3\otimes \shift \varphi^{41}\otimes \shift \varphi^{421}\otimes \varphi^{422}\dots \varphi^{425}\varphi^{43}\varphi^{44}\varphi^{45}\varphi^5}_{303}\\
&-\underbrace{\shift \partial(\varphi^1)\varphi^2\otimes \shift \varphi^3\otimes \shift \varphi^{41}\otimes \shift \varphi^{421}\otimes \varphi^{422}\dots \varphi^{425}\varphi^{43}\varphi^{44}\varphi^{45}\varphi^5}_{g}\\
&\boldsymbol{\uwave{+}}0\boldsymbol{\uwave{-}}\underbrace{1\otimes \shift \varphi^1\otimes \varphi^2\varphi^3\otimes \shift \varphi^{41}\otimes \shift \varphi^{421}\otimes \varphi^{422}\dots \varphi^{425}\varphi^{43}\varphi^{44}\varphi^{45}\varphi^5}_{157}\\
&+\underbrace{1\otimes \shift \varphi^1\otimes \varphi^2\otimes \shift (\varphi^3\varphi^{41})\otimes \shift \varphi^{421}\otimes \varphi^{422}\dots \varphi^{425}\varphi^{43}\varphi^{44}\varphi^{45}\varphi^5}_{g}\\
&-\underbrace{1\otimes \shift \varphi^1\otimes \varphi^2\otimes \shift \varphi^3\otimes \shift (\varphi^{41}\varphi^{421})\otimes \varphi^{422}\dots \varphi^{425}\varphi^{43}\varphi^{44}\varphi^{45}\varphi^5}_{y}\\
&-\underbrace{1\otimes \shift \varphi^1\otimes \varphi^2\otimes \shift \varphi^3\otimes \shift \varphi^{41}\otimes (\shift \otimes \shift)\partial(\shift^{-1}\varphi^{42})\varphi^{43}\varphi^{44}\varphi^{45}\varphi^5}_{229}\boldsymbol{\uwave{+}}0\boldsymbol{\uwave{+}}0\boldsymbol{\uwave{+}}0
\end{align*}
and
\begin{align*}
&-\mu_1^{\otimes}(-1\otimes \shift (\varphi^1\varphi^{21})\otimes \varphi^{22}\varphi^{23}\varphi^{241}\dots \varphi^{244}\otimes \shift \varphi^{245}\otimes \shift \varphi^{25}\otimes \shift \varphi^3\otimes \varphi^4\varphi^5)\\
=&+\underbrace{\hat{\omega} \varphi^1\varphi^{21}\varphi^{22}\varphi^{23}\varphi^{241}\dots \varphi^{244}\otimes \shift \varphi^{245}\otimes \shift \varphi^{25}\otimes \shift \varphi^3\otimes \varphi^4\varphi^5}_{306}\\
&-\underbrace{\varphi^1\varphi^{21}\hat{\omega} \varphi^{22}\varphi^{23}\varphi^{241}\dots \varphi^{244}\otimes \shift \varphi^{245}\otimes \shift \varphi^{25}\otimes \shift \varphi^3\otimes \varphi^4\varphi^5}_{307}\\
&-\underbrace{\shift \partial(\varphi^1)\varphi^{21}\varphi^{22}\varphi^{23}\varphi^{241}\dots \varphi^{244}\otimes \shift \varphi^{245}\otimes \shift \varphi^{25}\otimes \shift \varphi^3\otimes \varphi^4\varphi^5}_{gg}\\
&-\underbrace{\varphi^1\shift \partial(\varphi^{21})\varphi^{22}\varphi^{23}\varphi^{241}\dots \varphi^{244}\otimes \shift \varphi^{245}\otimes \shift \varphi^{25}\otimes \shift \varphi^3\otimes \varphi^4\varphi^5}_{g}\\
&\boldsymbol{\uwave{+}}0\boldsymbol{\uwave{+}}0\boldsymbol{\uwave{+}}0\boldsymbol{\uwave{+}}\underbrace{1\otimes \shift (\varphi^1\varphi^{21})\otimes \varphi^{22}\varphi^{23}(\shift \otimes \shift)\partial(\shift^{-1}\varphi^{24})\otimes \shift \varphi^{25}\otimes \shift \varphi^3\otimes \varphi^4\varphi^5}_{252}\\
&+\underbrace{1\otimes \shift (\varphi^1\varphi^{21})\otimes \varphi^{22}\varphi^{23}\varphi^{241}\dots \varphi^{244}\otimes \shift (\varphi^{245}\varphi^{25})\otimes \shift \varphi^3\otimes \varphi^4\varphi^5}_y\\
&-\underbrace{1\otimes \shift (\varphi^1\varphi^{21})\otimes \varphi^{22}\varphi^{23}\varphi^{241}\dots \varphi^{244}\otimes \shift \varphi^{245}\otimes \shift (\varphi^{25}\varphi^3)\otimes \varphi^4\varphi^5}_g\\
&+\underbrace{1\otimes \shift (\varphi^1\varphi^{21})\otimes \varphi^{22}\varphi^{23}\varphi^{241}\dots \varphi^{244}\otimes \shift \varphi^{245}\otimes \shift \varphi^{25}\otimes \varphi^3\varphi^4\varphi^5}_{258}\boldsymbol{\uwave{+}}0
\end{align*}
and 
\begin{align*}
&-\mu_1^{\otimes}(-\varphi^1\varphi^2\otimes \shift \varphi^3\otimes \shift \varphi^{41}\otimes \shift \varphi^{421}\otimes \varphi^{422}\dots \varphi^{425}\varphi^{43}\varphi^{44}\otimes \shift (\varphi^{45}\varphi^5)\otimes 1)\\
=&0\boldsymbol{\uwave{+}}\underbrace{\varphi^1\varphi^2\varphi^3\otimes \shift \varphi^{41}\otimes \shift \varphi^{421}\otimes \varphi^{422}\dots \varphi^{425}\varphi^{43}\varphi^{44}\otimes \shift (\varphi^{45}\varphi^5)\otimes 1}_{243}\\
&-\underbrace{\varphi^1\varphi^2\otimes \shift (\varphi^3\varphi^{41})\otimes \shift \varphi^{421}\otimes \varphi^{422}\dots \varphi^{425}\varphi^{43}\varphi^{44}\otimes \shift (\varphi^{45}\varphi^5)\otimes 1}_{g}\\
&+\underbrace{\varphi^1\varphi^2\otimes \shift \varphi^3\otimes \shift (\varphi^{41}\varphi^{421})\otimes \varphi^{422}\dots \varphi^{425}\varphi^{43}\varphi^{44}\otimes \shift (\varphi^{45}\varphi^5)\otimes 1}_{g}\\
&+\underbrace{\varphi^1\varphi^2\otimes \shift \varphi^3\otimes \shift \varphi^{41}\otimes (\shift \otimes \shift)\partial(\shift^{-1}\varphi^{42})\varphi^{43}\varphi^{44}\otimes \shift (\varphi^{45}\varphi^5)\otimes 1}_{231}\\
&\boldsymbol{\uwave{+}}0\boldsymbol{\uwave{+}}0\boldsymbol{\uwave{+}}0\boldsymbol{\uwave{-}}\underbrace{\varphi^1\varphi^2\otimes \shift \varphi^3\otimes \shift \varphi^{41}\otimes \shift \varphi^{421}\otimes \varphi^{422}\dots \varphi^{425}\varphi^{43}\varphi^{44}\hat{omega} \varphi^{45}\varphi^5}_{309}\\
&+\underbrace{\varphi^1\varphi^2\otimes \shift \varphi^3\otimes \shift \varphi^{41}\otimes \shift \varphi^{421}\otimes \varphi^{422}\dots \varphi^{425}\varphi^{43}\varphi^{44}\varphi^{45}\varphi^5\hat{\omega}}_{308}\\
&+\underbrace{\varphi^1\varphi^2\otimes \shift \varphi^3\otimes \shift \varphi^{41}\otimes \shift \varphi^{421}\otimes \varphi^{422}\dots \varphi^{425}\varphi^{43}\varphi^{44}\partial(\varphi^{45})\varphi^5}_{g}\\
&+\underbrace{\varphi^1\varphi^2\otimes \shift \varphi^3\otimes \shift \varphi^{41}\otimes \shift \varphi^{421}\otimes \varphi^{422}\dots \varphi^{425}\varphi^{43}\varphi^{44}\varphi^{45}\shift \partial(\varphi^5)}_{gg}
\end{align*}
and
\begin{align*}
&-\mu_1^{\otimes}(-\varphi^1\varphi^2\otimes \shift \varphi^3\otimes \shift \varphi^{41}\otimes \shift \varphi^{421}\otimes \shift \varphi^{4221}\otimes \varphi^{4222}\dots \varphi^{4225}\varphi^{423}\varphi^{424}\varphi^{425}\varphi^{43}\varphi^{44}\varphi^{45}\varphi^5)\\
=&0\boldsymbol{\uwave{+}}\underbrace{\varphi^1\varphi^2\varphi^3\otimes \shift \varphi^{41}\otimes \shift \varphi^{421}\otimes \shift \varphi^{4221}\otimes \varphi^{4222}\dots \varphi^{4225}\varphi^{423}\varphi^{424}\varphi^{425}\varphi^{43}\varphi^{44}\varphi^{45}\varphi^5}_{298}\\
&-\underbrace{\varphi^1\varphi^2\otimes \shift (\varphi^3\varphi^{41})\otimes \shift \varphi^{421}\otimes \shift \varphi^{4221}\otimes \varphi^{4222}\dots \varphi^{4225}\varphi^{423}\varphi^{424}\varphi^{425}\varphi^{43}\varphi^{44}\varphi^{45}\varphi^5}_{gg}\\
&+\underbrace{\varphi^1\varphi^2\otimes \shift \varphi^3\otimes \shift (\varphi^{41}\varphi^{421})\otimes \shift \varphi^{4221}\otimes \varphi^{4222}\dots \varphi^{4225}\varphi^{423}\varphi^{424}\varphi^{425}\varphi^{43}\varphi^{44}\varphi^{45}\varphi^5}_{g}\\
&-\underbrace{\varphi^1\varphi^2\otimes \shift \varphi^3\otimes \shift \varphi^{41}\otimes \shift (\varphi^{421}\varphi^{4221})\otimes \varphi^{4222}\dots \varphi^{4225}\varphi^{423}\varphi^{424}\varphi^{425}\varphi^{43}\varphi^{44}\varphi^{45}\varphi^5}_{y}\\
&-\underbrace{\varphi^1\varphi^2\otimes \shift \varphi^3\otimes \shift \varphi^{41}\otimes \shift \varphi^{421}\otimes (\shift \otimes \shift)\partial(\shift^{-1}\varphi^{422})\varphi^{423}\varphi^{424}\varphi^{425}\varphi^{43}\varphi^{44}\varphi^{45}\varphi^5}_{297}\boldsymbol{\uwave{+}}0\boldsymbol{\uwave{+}}0\boldsymbol{\uwave{+}}0\boldsymbol{\uwave{+}}0
\end{align*}
and
\begin{align*}
&-\mu_1^{\otimes}(-\varphi^1\varphi^2\varphi^3\varphi^{41}\varphi^{42}\varphi^{43}\varphi^{441}\varphi^{442}\varphi^{443}\varphi^{4441}\dots \varphi^{4444}\otimes \shift \varphi^{4445}\otimes \shift \varphi^{445}\otimes \shift \varphi^{45}\otimes \shift \varphi^5\otimes 1)\\
=&0\boldsymbol{\uwave{+}}0\boldsymbol{\uwave{+}}0\boldsymbol{\uwave{+}}0\boldsymbol{\uwave{+}}0\boldsymbol{\uwave{-}}\underbrace{\varphi^1\varphi^2\varphi^3\varphi^{41}\varphi^{42}\varphi^{43}\varphi^{441}\varphi^{442}\varphi^{443}(\shift \otimes \shift)\partial(\shift^{-1}\varphi^{444})\otimes \shift \varphi^{445}\otimes \shift \varphi^{45}\otimes \shift \varphi^5\otimes 1}_{299}\\
&-\underbrace{\varphi^1\varphi^2\varphi^3\varphi^{41}\varphi^{42}\varphi^{43}\varphi^{441}\varphi^{442}\varphi^{443}\varphi^{4441}\dots \varphi^{4444}\otimes \shift (\varphi^{4445}\varphi^{445})\otimes \shift \varphi^{45}\otimes \shift \varphi^5\otimes 1}_y\\
&+\underbrace{\varphi^1\varphi^2\varphi^3\varphi^{41}\varphi^{42}\varphi^{43}\varphi^{441}\varphi^{442}\varphi^{443}\varphi^{4441}\dots \varphi^{4444}\otimes \shift \varphi^{4445}\otimes \shift (\varphi^{445}\varphi^{45})\otimes \shift \varphi^5\otimes 1}_g\\
&-\underbrace{\varphi^1\varphi^2\varphi^3\varphi^{41}\varphi^{42}\varphi^{43}\varphi^{441}\varphi^{442}\varphi^{443}\varphi^{4441}\dots \varphi^{4444}\otimes \shift \varphi^{4445}\otimes \shift \varphi^{445}\otimes \shift (\varphi^{45}\varphi^5)\otimes 1}_{gg}\\
&+\underbrace{\varphi^1\varphi^2\varphi^3\varphi^{41}\varphi^{42}\varphi^{43}\varphi^{441}\varphi^{442}\varphi^{443}\varphi^{4441}\dots \varphi^{4444}\otimes \shift \varphi^{4445}\otimes \shift \varphi^{445}\otimes \shift \varphi^{45}\otimes \varphi^5.}_{158}
\end{align*}

\end{proof}

\end{appendices}

\end{document}